%% file: MCgood_Revision_decolorized_renamed.tex
\documentclass[hidelinks,onefignum,onetabnum]{siamart220329}

\input{MCgood_shared_Revision}

\ifpdf
\hypersetup{
  pdftitle={Monte Carlo is a good sampling strategy},
  pdfauthor={B. Adcock and S. Brugiapaglia}
}
\fi

\begin{document}

{\maketitle}

% REQUIRED
\begin{abstract}
This paper concerns the approximation of smooth, high-dimensional functions from limited samples using polynomials. This task lies at the heart of many applications in computational science and engineering -- notably, {some of} those arising from parametric modelling and computational uncertainty quantification. It is common to use Monte Carlo sampling in such applications, so as not to succumb to the curse of dimensionality. However, it is well known that such a strategy is theoretically suboptimal. Specifically, there are many polynomial spaces of dimension $n$ for which the sample complexity scales log-quadratically, i.e., like $c \cdot n^2 \cdot \log(n)$ as $n \rightarrow \infty$. This well-documented phenomenon has led to a concerted effort over the last decade to design improved, and moreover, near-optimal strategies, whose sample complexities scale log-linearly, or even linearly in $n$.

{I}n this work we demonstrate that Monte Carlo is actually a perfectly good strategy in high dimensions, despite its apparent suboptimality. We first document this phenomenon empirically via a systematic set of numerical experiments. Next, we present a theoretical analysis {that rigorously justifies this fact}
in the case of holomorphic functions of infinitely-many variables. We show that there is a least-squares approximation based on $m$ Monte Carlo samples whose error decays algebraically fast in $m/\log(m)$, with a rate that is the same as that of the best $n$-term polynomial approximation. This result is non-constructive, since it assumes knowledge of a suitable polynomial subspace in which to perform the approximation. We next present a compressed sensing-based scheme that achieves the same rate, except for a larger polylogarithmic factor.  This scheme is practical, and numerically it performs as well as or better than well-known adaptive least-squares schemes.

Overall, our findings in this paper {demonstrate that} Monte Carlo sampling is {a good choice for polynomial approximation in high dimensions and shed light on why this is the case}. Therefore, the benefits of improved sampling strategies are generically limited to lower-dimensional settings. 
\end{abstract}

% REQUIRED

% REQUIRED

\section{Introduction}
\label{s:introduction}

Approximating a smooth function $f : \cU \subseteq \bbR^d \rightarrow \bbC$ from (noisy) sample values
\be{
\label{sample-values}
f(\bm{y}_i) + e_i,\quad i = 1,\ldots,m,\qquad \text{where }\bm{y}_1,\ldots \bm{y}_m \in \mathcal{U},
}
is a task of fundamental importance in computational science and engineering. This task is well understood in low dimensions. But many modern applications {\cite{smith2013uncertainty,sullivan2015introduction}} call for the approximation of functions depending on many (and potentially infinitely many) variables. Such high-dimensional approximation problems occur in many fields, for example in parametric modelling and computational Uncertainty Quantification (UQ). {Here, smooth, high-dimensional functions commonly arise as solution maps of parametric Differential Equations (DEs).} 

{M}ethods based on polynomials have proved to {be} effective tools for approximating such functions. {Both classical techniques such as Least Squares (LS) and more recent tools such as Compressed Sensing (CS) have been intensively investigated over the last several decades. See \S \ref{ss:history} for a historical discussion. Although by no means the only choice, these tools have been quite widely adopted in the aforementioned applications.}

\subsection{{Monte Carlo s}ampling and sample complexity}

This paper is about the choice of sampling strategy for polynomial approximation in high 
dimensions. In any approximation scheme, the choice of sample points $\bm{y}_1,\ldots,\bm{y}_m$ is of singular importance. Obtaining samples is often a key bottleneck in applications, since they may require time- or resource-consuming numerical simulations or physical experiments. Thus, it is vital to choose sample points in a judicious manner, so as to facilitate accurate and stable approximations from as few samples as possible. 

Care is needed when selecting a sampling strategy for a high-dimensional problem  so as not to succumb to the {\textit{curse of dimensionality}}. Taking inspiration from high-dimensional quadrature, it is common to use \textit{Monte Carlo (MC)} sampling in polynomial approximation schemes, i.e., the points $\bm{y}_1,\ldots,\bm{y}_m$ are drawn randomly and independently from some underlying probability measure. Such sampling strategies also occur naturally {in} UQ settings, where the variables are stochastic.

{MC sampling is a remarkably effective technique for high-dimensional integration. Yet, in the context of high-dimensional polynomial approximation,} MC sampling suffers from a critical limitation. The \textit{sample complexity} -- the number of samples $m$ required for a stable and accurate approximation -- usually scales poorly with the problem dimension $d$ or approximation space dimension $n$. For example, if $\cU = [-1,1]^d$ and the underlying measure is the uniform measure, then there are polynomial subspaces of dimension $n$ for which the sample complexity scales \textit{log-quadratically}, i.e., 
\be{
\label{LS-samp-comp-intro}
m = c \cdot n^2 \cdot \log(n)
}
{for some $c > 0$.} See, e.g., \cite{chkifa2015discrete,migliorati2014analysis,cohen2013stability}. {Note that $n$ also scales exponentially with $d$ for many classical approximation spaces, such as those corresponding to \emph{tensor-product}, \emph{total degree}, or \emph{hyperbolic cross} index sets (see \S \ref{ss:choosing-S}). As a result, \eqref{LS-samp-comp-intro} implies that $m$ may need to grow twice as exponentially fast in $d$ as $n$ for these choices.}

\subsection{Improved and (near-) near-optimal sampling strategies}\label{s:improved-sampling-intro}

This well-known limitation of MC sampling has led to a concerted effort in the development of sampling strategies that offer either better theoretical sample complexity, better practical performance or (ideally) both. One finds many different approaches in the recent literature, including \textit{preconditioning} \cite{rauhut2012sparse}, \textit{asymptotic sampling} \cite{hampton2015coherence}, \textit{coherence-optimal sampling} \cite{hampton2015coherence}, \textit{Christoffel sampling} \cite{narayan2017christoffel}, \textit{randomly subsampled quadratures} \cite{zhou2015weighted}, \textit{low-discrepancy points} \cite{migliorati2015analysis}, \textit{Quasi Monte Carlo sampling}, \textit{Latin hypercube sampling}, \textit{deterministically subsampled quadratures} \cite{seshadri2017effectively}, \textit{boosting techniques} \cite{haberstich2022boosted}, and methods based on \textit{optimal design of experiments} \cite{fajraoui2017sequential,diaz2018sparse}. {See \cite{guo2020constructing,hadigol2018least} for further discussion and references.} However, perhaps most notably, in the last several years random sampling strategies have been introduced that are provably \textit{near-optimal} \cite{cohen2017optimal} {(see also \cite{hampton2015coherence} for earlier work in this direction and \cite{adcock2022towards,guo2020constructing} for reviews).} Specifically, for any fixed (polynomial or, in fact, nonpolynomial) approximation space of dimension $n$, the corresponding sample complexity scales \textit{log-linearly}, i.e.,
\be{
\label{LS-samp-comp-intro-opt}
m = c \cdot n \cdot \log(n).
}

\subsection{Aims of this paper}\label{ss:aims}

\begin{figure}[t!]
\begin{center}
\begin{small}
 \begin{tabular}{@{\hspace{0pt}}c@{\hspace{\errplotsp}}c@{\hspace{\errplotsp}}c@{\hspace{0pt}}}
\includegraphics[width = \errplotimg]{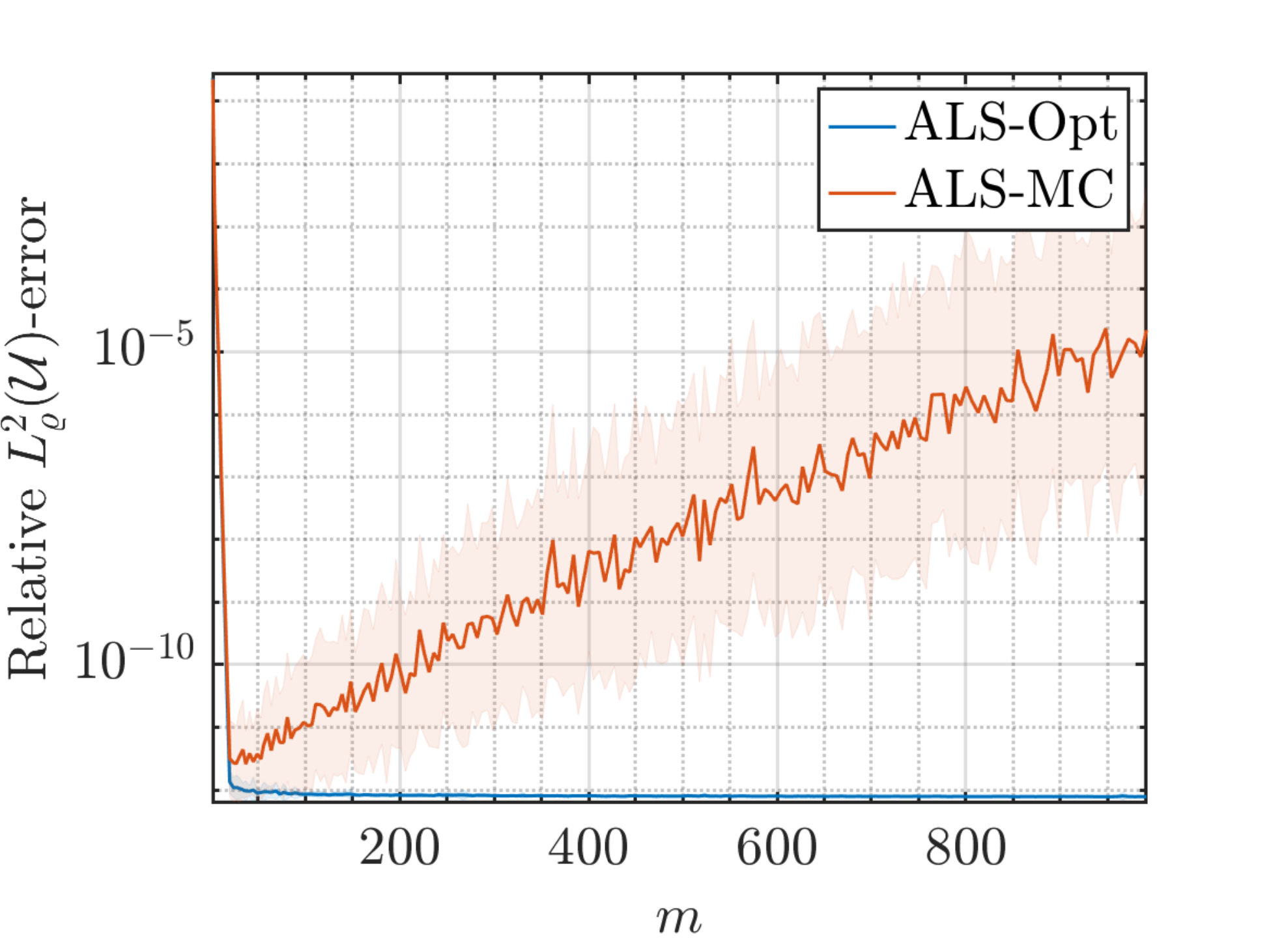}
&
\includegraphics[width = \errplotimg]{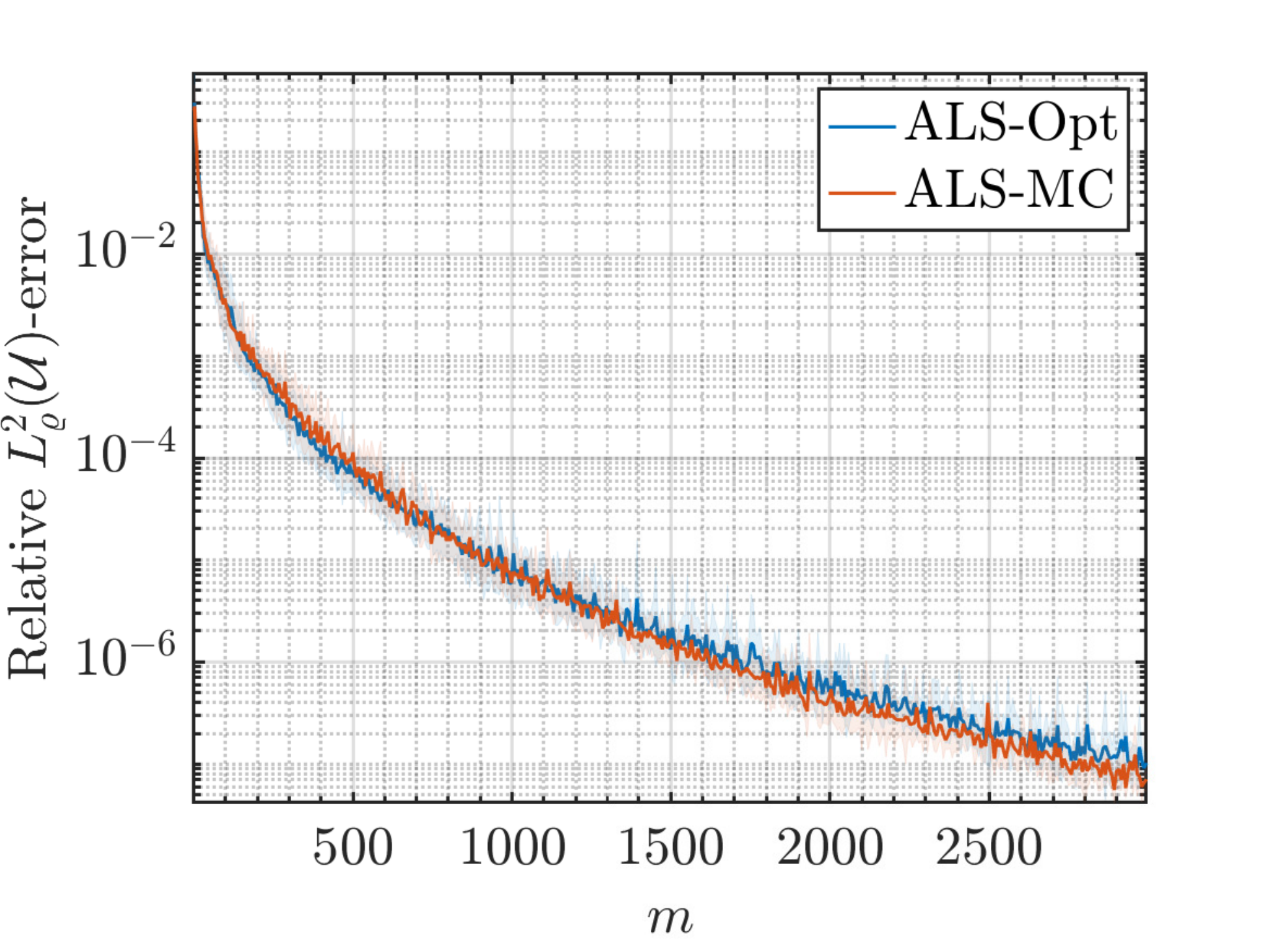}
&
\includegraphics[width = \errplotimg]{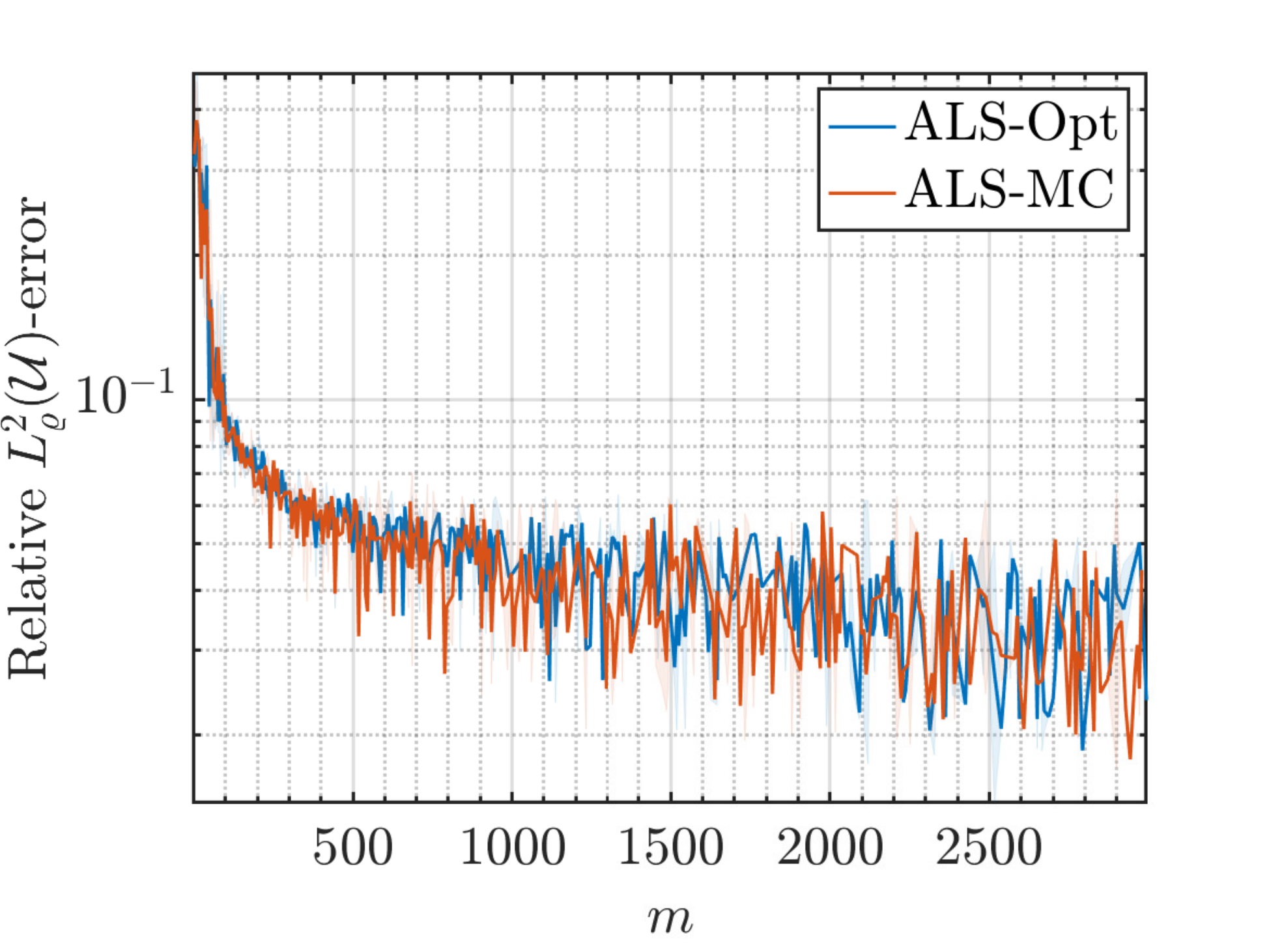}
\\[\errplotgraphsp]
\includegraphics[width = \errplotimg]{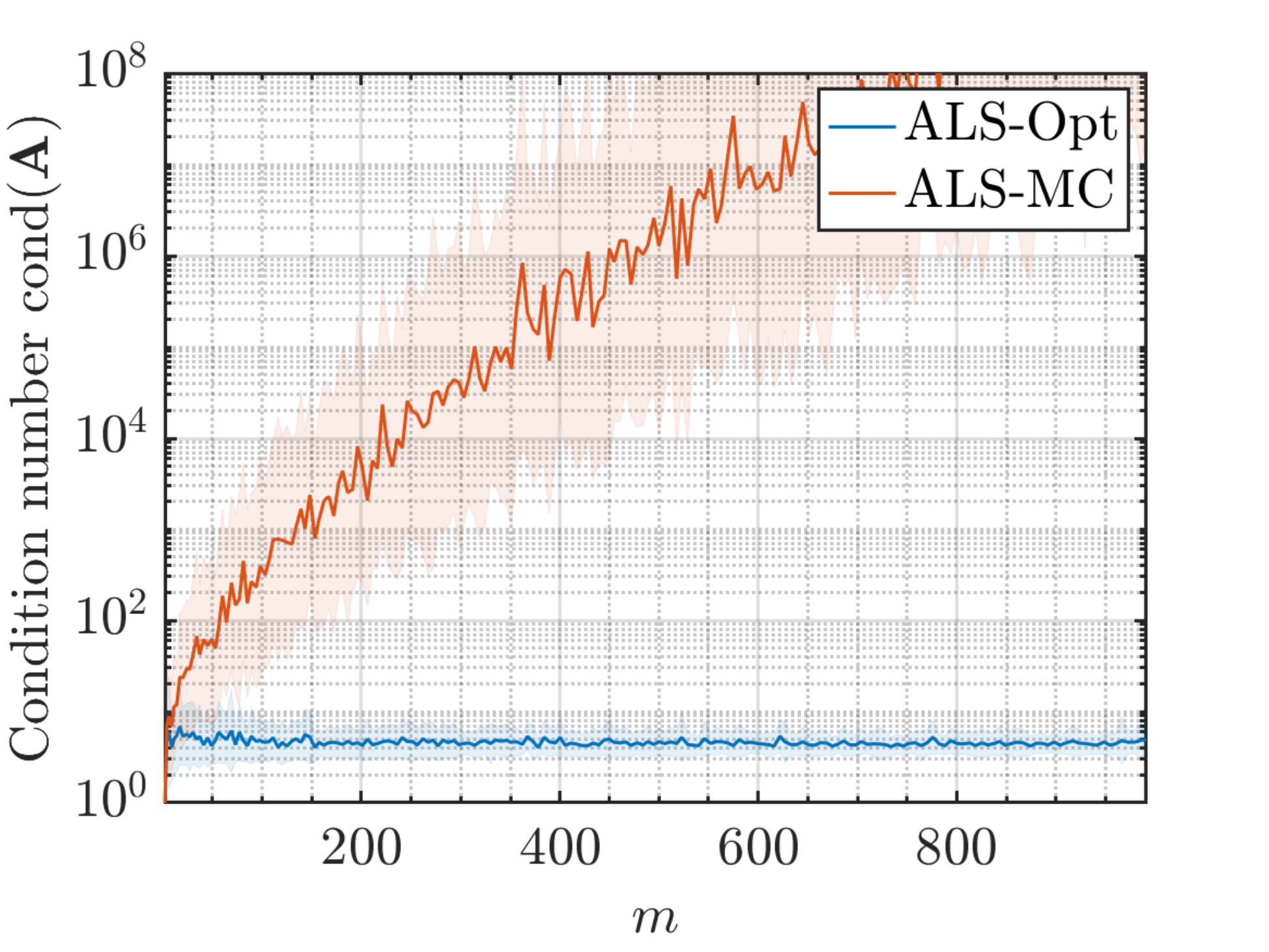}
&
\includegraphics[width = \errplotimg]{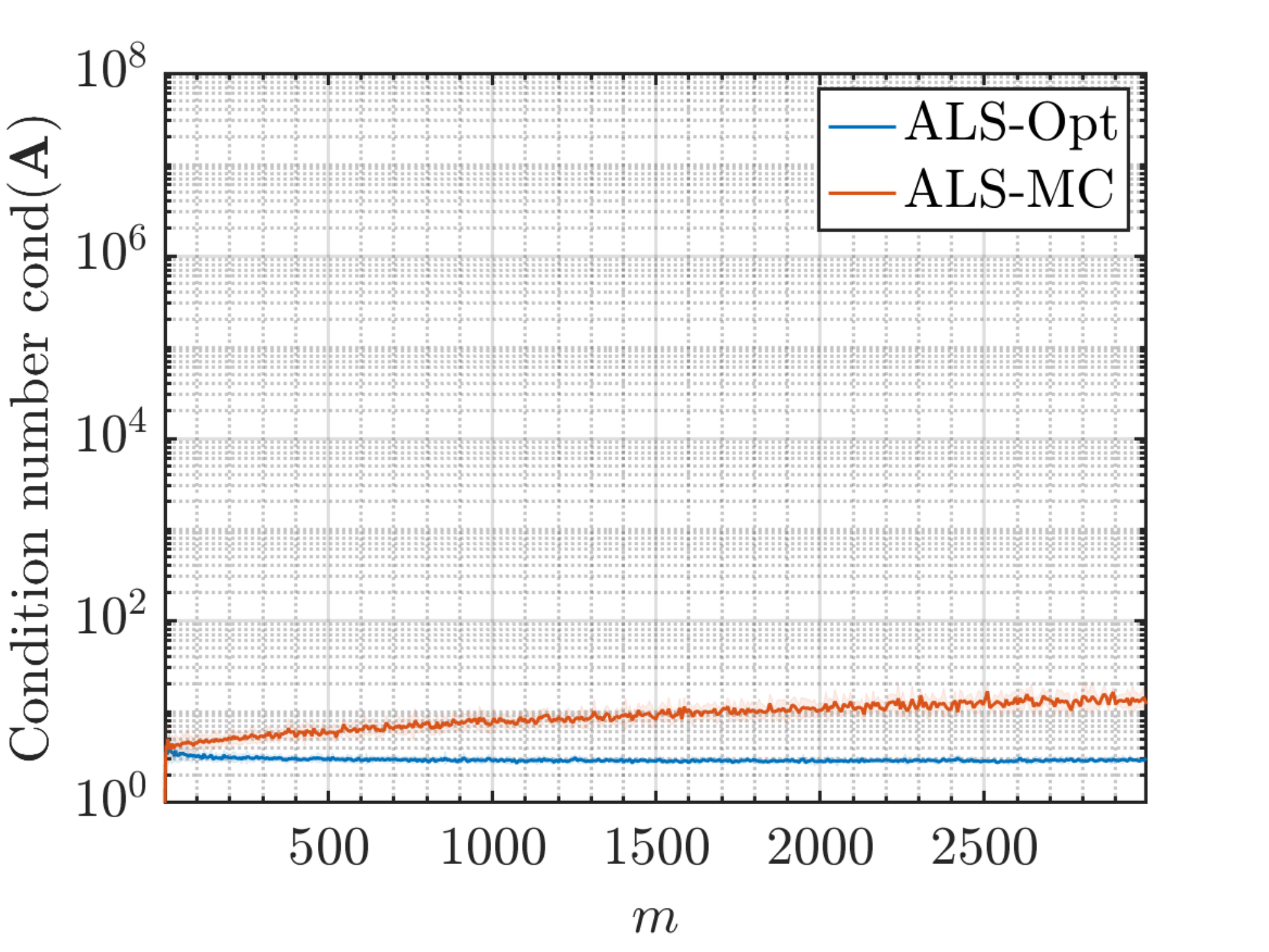}
&
\includegraphics[width = \errplotimg]{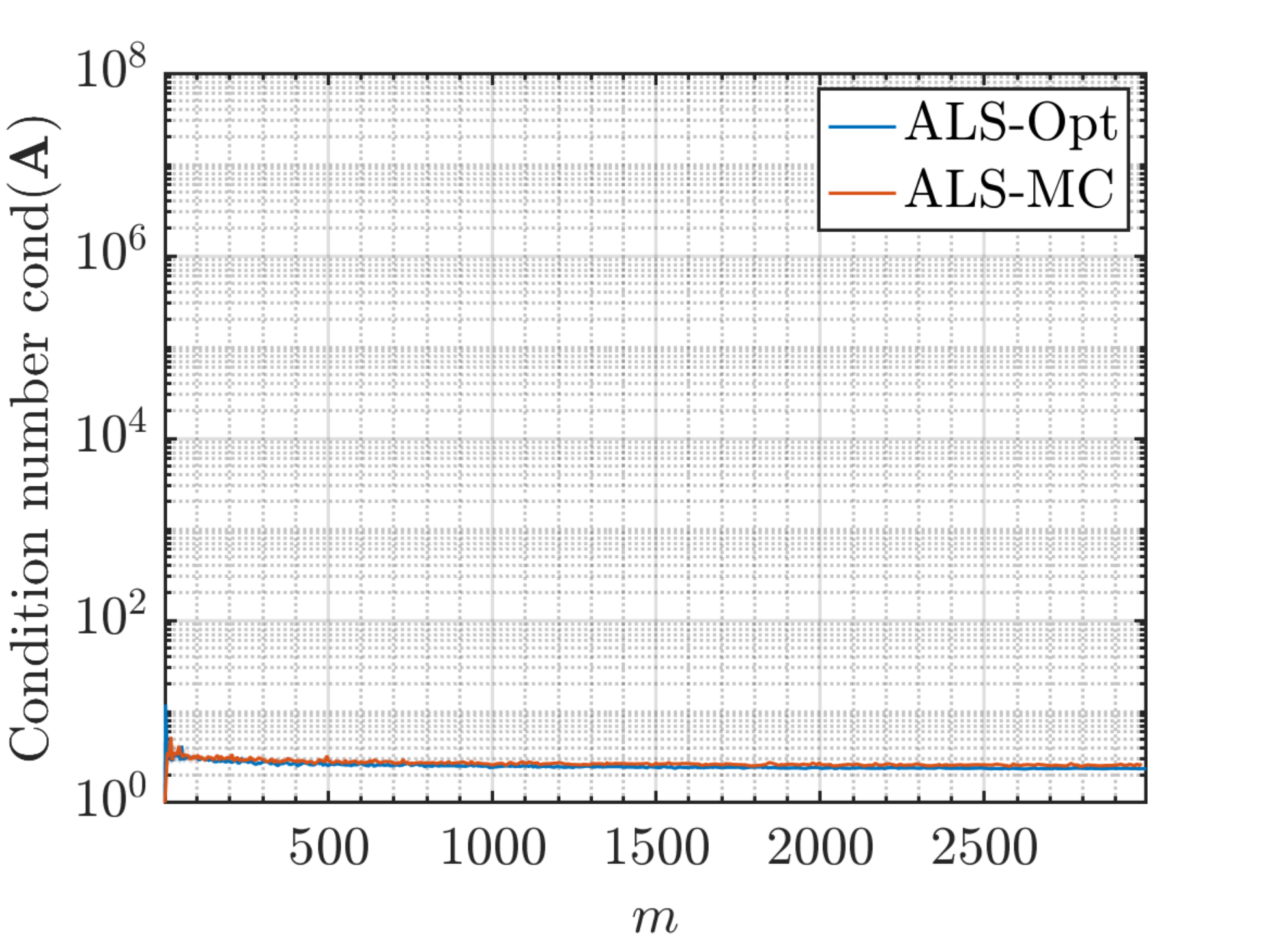}
\\[\errplottextsp]
$d = 1$ & $d = 4$ & $d = 32$
\end{tabular}
\end{small}
\end{center}
\caption{In low dimensions, MC sampling from the uniform measure on $\cU = [-1,1]^d$ is a poor choice. However, it becomes far better as the dimension increases. This figure plots the error (top row) and {($2$-norm)} condition number of the LS matrix (bottom row) versus $m$ for (adaptive) LS polynomial approximation using MC samples (ALS-MC) and a near-optimal sampling strategy (ALS-Opt) described in \S \ref{s:improved-sampling-intro}. The function considered is a quantity of interest of a parametric {DE} with lognormal diffusion term (see \S \ref{ss:lognormal-pDE}). The scaling $m = \max \{ n+1 , \lceil n \log(n) \rceil \}$ is used in both cases. For MC sampling, the LS problem is highly ill-conditioned when $d = 1$, since $m$ scales log-linearly with $n$ rather than log-quadratically, as in \eqref{LS-samp-comp-intro}. This causes a dramatic loss of accuracy due to round-off error. However, as $d$ increases, this LS problem becomes much better conditioned, with both its condition number and the error closely tracking those of the near-optimal strategy. }
\label{fig:fig1}
\end{figure}

{The previous discussion suggests that MC is a bad sampling strategy for polynomial approximation.} And indeed, almost all of the above strategies yield significant benefits over MC sampling in low dimensions, where the dimension $d$ may be on the order of five or less. However, such benefits have been consistently observed to lessen as the dimension increases. When $d$ is on the order of 10 or more, the performance gap between MC sampling and any improved sampling strategy is often strikingly less. This is even the case for the near-optimal sampling strategies, in spite of their theoretical optimality. 
A typical example of effect is shown in Fig.\ \ref{fig:fig1}, with further examples presented in \S \ref{s:MC-bad}. {Here LS with MC sampling is unstable and potentially divergent in low dimensions as the number of samples $m$ increases, but in higher dimensions its performance is essentially the same as that LS with a theoretically near-optimal sampling strategy.}

{The aim of this paper is to shed light on this phenomenon.} 
We first document {it} through a series of numerical experiments {on different smooth test functions in both low and high dimensions (in the range $d = 1$ to $d = 32$)}. Then we provide a theoretical explanation for why it occurs. {We summarize these theoretical results next.}

\subsection{Nontechnical summary of main theoretical results}\label{s:main-res-summary}

{Our theoretical analysis centres on polynomial approximation of holomorphic functions of infinitely-many variables. Specifically, we combine the theory of \textit{best $n$-term polynomial approximation} of such functions with the theories of LS and CS. }

\subsubsection{Polynomial approximation theory in infinite dimensions}
\label{ss:main-res-approx-theory}

The topic of polynomial approximation theory for infinite-dimensional functions has developed over the last decade or so, motivated by applications in parametric DEs and UQ. 
{See \cite{cohen2015approximation} and \cite[Chpts.\ 3 \& 4]{adcock2022sparse} for reviews}. A key result in this area is the assertion of certain \textit{algebraic} rates of convergence for best $n$-term polynomial approximations to certain classes of holomorphic, infinite-dimensional functions.

Consider the infinite hypercube $\cU = [-1,1]^{\bbN}$ equipped with the uniform measure $\varrho$ and let $L^2_{\varrho}(\cU)$ be the Lebesgue space of (complex-valued) square-integrable functions $ [-1,1]^{\bbN} \rightarrow \bbC$. (For the reader who is unfamiliar with concepts of infinite-dimensional measures and function spaces, we refer to \S \ref{s:polyapp-inf-dim}). Any function $f \in L^2_{\varrho}(\cU)$ has an infinite expansion with respect to multivariate orthogonal polynomials, which in this case are tensor products of the univariate Legendre polynomials. The \textit{best $n$-term approximation} to $f$ is the polynomial formed simply by picking the largest $n$ terms of this expansion in absolute value. In this paper, we denote this approximation as $f_n$. 

Now let $0 < p <1$, $\varepsilon > 0$ and $\bm{b} \in \ell^p(\bbN)$. Then {a key result in this area is that there} is a class $\cH(\bm{b},\varepsilon) \subset L^2_{\varrho}(\cU)$ 
of holomorphic functions within which the best $n$-term approximation converges algebraically with rate depending on $1/p$. Specifically,
\be{
\label{best-n-term-err}
\nmu{f - f_n}_{L^2_{\varrho}(\cU)} \leq C(\bm{b},\varepsilon,p) \cdot n^{\frac12-\frac1p},\quad \forall f \in \cH(\bm{b},\varepsilon),\ n \in \bbN.
}
{In other words, there are classes of functions whose approximation is free from the curse of dimensionality: despite depending on infinitely many variables, they can be approximated with algebraic rate by a finite ($n$-term)  polynomial approximation.} 

The class $\cH(\bm{b},\varepsilon)$, which is known as the class of (unit-norm) \textit{$(\bm{b},\varepsilon)$-holomorphic} functions \cite{chkifa2015breaking,schwab2019deep}, is both theoretically interesting, since it defines a class where algebraic rates are obtained, \textit{and} practically relevant. It was first introduced to describe the parametric regularity of the solutions of certain parametric DEs. It is now known that {several important families} 
of parametric DEs have parametric solution maps that belong to such classes (for suitable $\bm{b}$, $\varepsilon$ depending on the DE). {These include parametric diffusion equations, parametric heat equations, PDEs over parametrized domains, and parametric initial value problems with affine parametric dependence.} See \S \ref{s:pDE-background} for some further discussion.  

Note that $(\bm{b},\varepsilon)$-holomorphic functions are \textit{anisotropic}. As we see in \S \ref{s:polyapp-inf-dim}, the $j$th component $b_j$ of the sequence $\bm{b}$ effectively controls the degree of smoothness of functions in $\cH(\bm{b},\varepsilon)$ with respect to the $j$th variable $y_j$.

\subsubsection{First contribution: sharpness of the rate \eqref{best-n-term-err}}\label{ss:first-contrib}

{In our first contribution, Theorem \ref{t:alg-rate-sharp}, we show that the rate \eqref{best-n-term-err} is effectively sharp. In particular, there exist sequences $\bm{b} \in \ell^p(\bbN)$ and infinitely many functions $f \in \cH(\bm{b},\varepsilon)$ for which $n^{1/2-1/p}$ is the best possible algebraic rate of convergence.}

\subsubsection{Second contribution: near-optimal approximation in the case of known anisotropy}

Our {other two} theoretical results concern constructing polynomial approximations to functions in $\cH(\bm{b},\varepsilon)$ from MC samples, i.e., $\bm{y}_1,\ldots,\bm{y}_m \sim_{\mathrm{i.i.d.}} \varrho${, that attain this algebraic rate}. 
In view of {\eqref{LS-samp-comp-intro}, one would generally expect to be able to recover an $n$-term polynomial approximation from $m$ samples via LS for $n$ at most $\sqrt{2m / (c \log m)}$. Combining this upper bound on $n$ with \eqref{best-n-term-err}}, it is therefore tempting to believe that the best one could hope to achieve is
\bes{
\nmu{f - \hat{f}}_{L^2_{\varrho}(\cU)} \leq C \cdot \left ( m / \log(m) \right )^{\frac12(\frac12-\frac1p)},
} 
when $\bm{b} \in \ell^p(\bbN)$.
In other words, the algebraic rate is halved, due to the log-quadratic sample complexity \eqref{LS-samp-comp-intro}. However, this is not the case. In Theorem \ref{t:main-res} we show that there is a LS approximation for which the error behaves like
\be{
\label{MC-LS-rate}
\nmu{f - \hat{f} }_{L^2_{\varrho}(\cU)} \leq C \cdot \left ( m / \log(m) \right )^{\frac12 - \frac1p}.
}
Hence the approximation converges with the same rate as that of the best $n$-term approximation, up to constants and the log term. We conclude that MC sampling is not only a good sampling strategy in infinitely many dimensions{, but in view of the result in \S \ref{ss:first-contrib},} it is in fact near optimal for the class $\cH(\bm{b},\varepsilon)$.

\subsubsection{Second contribution: near-optimal approximation in the case of unknown anisotropy}

A limitation of any LS scheme is that it requires knowledge of a space in which the function is well approximated. The above result is no different, since it involves a judicious choice of polynomial subspace of dimension $\leq n$ (depending on the parameters $\bm{b}$ and $\varepsilon$) in which one simultaneously obtains the same algebraic rate  of convergence \eqref{best-n-term-err} as the best $n$-term approximation, while also not suffering from poor sample worst-case complexity bound \eqref{LS-samp-comp-intro-opt} for MC sampling. At the very least, constructing such a subspace requires \textit{a priori} knowledge of the parameters $\bm{b}$ and $\varepsilon$. Unfortunately, these are {generally} unknown in practice.

Because of this limitation, LS is often used as part of an adaptive approximation scheme \cite{migliorati2019adaptive}. Here, rather than a single approximation, one computes a sequence of approximations $\hat{f}^{(1)},\hat{f}^{(2)},\ldots$ and (nested) polynomial subspaces $\cP^{(1)} \subseteq \cP^{(2)} \subseteq \cdots$, in which the approximation $\hat{f}^{(l)}$ is used to construct the next subspace $\cP^{(l+1)}$, typically in a greedy manner. As we see later in this paper, so-called \textit{Adaptive Least-Squares (ALS)} approximation is often quite effective. However, it lacks theoretical guarantees. It is currently unknown whether such an approximation achieves the algebraic rates of convergence \eqref{best-n-term-err} of the best $n$-term approximation (up to constants and log factors).

With this in mind, in the final part of this paper we demonstrate that this issue can be avoided by changing the approximation procedure. By using ideas from CS, we show it is possible to compute a polynomial approximation from $m$ MC samples for which the error behaves like
\be{
\label{MC-CS-rate}
\nmu{f - \hat{f} }_{L^2_{\varrho}(\cU)} \leq C \cdot \left ( m / \log^4(m) \right )^{\frac12-\frac1p},
}
subject to the slightly stricter requirement $\bm{b} \in \ell^p_{\mathsf{M}}(\bbN)$, where $\ell^p_{\mathsf{M}}(\bbN)$ is the \textit{monotone} $\ell^p$-space. See Theorem \ref{t:main-res-2}. Thus, for a slightly restricted class of functions, it is possible to achieve the same algebraic rates of convergence as the best $n$-approximation in terms of the number of samples $m$, up to a polylogarithmic factor.

This scheme is also practical. We present a series of numerical experiments comparing it to ALS. These show that the CS scheme offers consistently as good as, or sometimes better, performance, while also having theoretical guarantees.

\subsubsection{{The gap between theory and practice}}\label{ss:theory-practice}

{Our theoretical analysis asserts that MC sampling is near-optimal for polynomial approximation in infinite dimensions. Throughout this article we present numerical experiments for functions of finitely many variables only. Our analysis does apply in finite dimensions -- since a function of $d$ variables can be viewed as a function of infinitely many variables that is constant with respect to all but the first $d$ variables -- but, as we discuss later, it may not accurately predict the convergence rate when the dimension is small. There is consequently a gap between the theoretical analysis and practice, which leads one naturally to wonder `how high is high-dimensional'? In our experiments, we consistently witness this effect in finite dimensions. Most often, it occurs when $d$ is still relatively small (e.g., $d \leq 8$), and certainly by 
$d = 32$, which is the largest dimension considered in this paper.  However, we currently have no theoretical understanding of how large $d$ should be for such behaviour to kick in. This is an interesting problem that we leave for future work.}

\subsubsection{Why bother?}

As noted, the development of improved sampling strategies for high-dimensional approximation has been an active area of interest over the last few years. The purpose of this work is to show that MC sampling is actually eminently suitable (in fact, near-optimal) for certain high-dimensional approximation tasks -- in particular, those arising from parametric DE problems, which served as the original motivation for much of this line of research.

However, the reader may justifiably be wondering why we bother. Why use MC sampling when we know how to design sampling strategies that are near-optimal regardless of the dimension? For this, we offer several arguments.
First, it is both academically interesting and practically relevant to understand the limits of what one stands to gain by changing the sampling strategy. Our results {suggest} that these gains are limited to lower-dimensional problems, at least for smooth function approximation on bounded domains. 
Second, {methods based on} MC {sampling have} the distinct advantage of {allowing one to decouple the function evaluation and function approximation phases}. Thus, in the parametric DE context, the major computational burden of evaluating the target function at the sample points can be trivially parallelized. On the other hand, when adaptive approximation schemes (e.g., the aforementioned ALS method) are combined with near-optimal sampling, the sampling strategy also becomes adaptive. 
Hence, the tasks of sampling the function (i.e., numerically solving a DE) and constructing the polynomial approximation cannot be decoupled. In practice, this might require the design of \emph{ad hoc} software interfaces or data transfer strategies.
Third and finally, MC samples are {exchangeable and} usually easy to generate. {They are} found ubiquitously in UQ and machine learning applications; in particular, applications involving legacy data, where one is not afforded the luxury to adapt the samples to the target function and/or approximation scheme. 

In summary, there are ample reasons why MC samples may be preferred in practice. Hence investigations into their theoretical and practical performance in high dimensions are particularly relevant.

\subsection{{Historical context and discussion}}\label{ss:history}

{Our primary interest is the polynomial approximation of smooth high- and infinite-dimensional functions, with specific focus on \textit{Least Squares (LS)} and  \textit{Compressed Sensing (CS)} 
techniques. LS approximation is a classical topic, with a genesis that 
can be traced back to very early work by Gauss and Legendre at the turn of the 18th and 19th centuries \cite{stigler1981gauss}. LS is utilized in countless areas of applied mathematics and statistics, with eminent examples such as function approximation (of direct interest in this paper) and statistical regression.}

{Our work stems from a vigorous research stream started in the early 2010s, characterized by renewed interest towards LS \cite{cohen2013stability,chkifa2015discrete, migliorati2014analysis} %
and the application of the CS paradigm \cite{doostan2011nonadapted,mathelin2012compressed,rauhut2012sparse}
in the context of high-dimensional function approximation. The fast development of this research area is in tandem with progress in the fields of stochastic and parametric DEs. For stochastic DEs, this is strongly related to \emph{polynomial chaos} technique. This dates back to the 1930s \cite{wiener1938homogeneous}, and has found tremendous success in computational science and engineering since its inception in the 1990s \cite{ghanem1990polynomial} and thanks to a series of influential works in the 2000s, such as \cite{xiu2002wiener,babuska2008stochastic}, which contributed to the development of the UQ field (see also \cite{smith2013uncertainty,sullivan2015introduction}). As noted, a signature result for parametric DEs was the discovery that solution maps of a large class of parametric models admit best $n$-term polynomial approximation rates that are independent of the number of parameters, and therefore free from the curse of dimensionality. For more details about parametric DEs and historical remarks we refer to \cite{cohen2015approximation} and \cite[Chpt.\ 4]{adcock2022sparse}.}

{
Although considerably younger than LS, MC sampling is also a classical technique. The happenings that led to its widespread use in computational science are intertwined with the construction of the first electronic programmable digital computer (ENIAC), pioneering simulations of complex physical systems, and -- regrettably
 -- the creation of the first nuclear weapons in the 1940s. The history of the MC method features scientists of the calibre of Fermi, von Neumann, Ulam and Metropolis \cite{metropolis1990alamos, metropolis1987beginning}.}

{One of the most popular applications of MC sampling is numerical quadrature. It is well known that the MC quadrature error scales proportionally to $m^{-1/2}$, where $m$ is the number of MC samples of the function to be integrated (see, e.g., \cite[Section 2.1]{owen2013monte}). Remarkably, this decay rate is independent of the number of function's variables, which -- at least, intuitively -- motivates the use of MC in high dimensions. We recall, though, that in this work we are concerned with high-dimensional function approximation, as opposed to quadrature. In particular, it is worth noting that best $n$-term approximation error rates of the smooth high-dimensional functions of interest in this paper decay much faster than the MC quadrature error (see \S\ref{ss:main-res-approx-theory}). 
}

{Naturally, any discussion of LS and MC sampling would be incomplete without a word on regression. We now warn the reader that this is \emph{not} a paper on statistical regression. Statistical regression refers to the problem of fitting a model to noisy data under some probabilistic assumptions on the noise (e.g., homoscedasticity). This is not the perspective adopted in this paper. In fact, our goal is the accurate and stable approximation of a `ground truth' function from pointwise samples, where the noise corrupting the samples need not be modelled as a random variable. For a more in-depth discussion of the differences between the function approximation and statistical regression settings, we refer to \cite[\S 1.1]{guo2020constructing}, of which we share the same viewpoint.}

\subsection{Outline}

In \S \ref{s:polyLS} we present an overview of polynomial approximation of multivariate functions and (weighted) LS. In \S \ref{s:theoryLS} we present relevant theory for (weighted) LS, before discussing MC sampling and the near-optimal strategy mentioned in \S \ref{s:improved-sampling-intro}.  In \S \ref{s:MC-bad} we present numerical experiments demonstrating the main phenomenon considered in this work. We then turn our attention to its theoretical explanation. In \S \ref{s:polyapp-inf-dim} we introduce and review polynomial approximation theory for $(\bm{b},\varepsilon)$-holomorphic functions. We present our first result in \S \ref{s:main-res}, i.e., the existence of a LS approximation achieving the bound \eqref{MC-LS-rate}. Finally, we consider CS schemes in \S \ref{s:polyappCS}, including both the error bound \eqref{MC-CS-rate} and a numerical comparison between ALS and CS.

This paper also has supplementary materials. These contain {background on infinite-dimensional measures (\S \ref{s:inf_dim_measures}),} proofs of various results (\S \ref{s:proofs}), additional information on the experimental setup (\S \ref{s:additional-info}), background on parametric DEs (\S \ref{s:pDE-background}) and further experiments (\S \ref{s:further-experiments}). MATLAB code reproducing all the experiments is available at \url{https://github.com/benadcock/is-MC-bad}.

\section{Polynomials and least-squares polynomial approximation}\label{s:polyLS}

In this section, we describe multivariate orthogonal polynomials and polynomial approximation via LS. In order to keep the technical level reasonable, we consider the finite-dimensional case only. The infinite-dimensional case is introduced in \S \ref{s:polyapp-inf-dim}.

\subsection{Univariate notation}\label{ss:univariate}

Let $\varrho$ be a probability measure on $[-1,1]$ and write $L^2_{\varrho}([-1,1])$ for the corresponding Lebesgue space of square-integrable functions $f: [-1,1] \rightarrow \bbC$. We assume that $\varrho$ generates a unique sequence of orthonormal polynomials $\{ \psi_{\nu} \}_{\nu \in \bbN_0} \subset L^2_{\varrho}([-1,1])$. In other words,
\bes{
\ip{\psi_{\nu}}{\psi_{\nu'}} = \delta_{\nu,\nu'},\quad \forall \nu , \nu' \in \bbN_0,\quad \text{and}\quad 
\spn \{ \psi_{0},\ldots,\psi_{n} \} = \bbP_{n},\quad \forall n \in \bbN_0,
}
where $\bbP_n$ is the space of polynomials of degree at most $n$.
Note that this is a mild assumption (see, e.g., \cite[\S 2.1]{narayan2018computation}). Two particular cases we focus on in this paper are the uniform and Chebyshev (arcsine) measures
\bes{
\D \varrho(y) = 2^{-1} \D y,\qquad \D \varrho(y) = (\pi \sqrt{1-y^2})^{-1} \D y,
}
which generate the Legendre and (first kind) Chebyshev polynomials, respectively.

\subsection{Multivariate polynomial approximation}\label{ss:multivar-poly-app}

Let $d \in \bbN$ and consider the symmetric hypercube $\cU =[-1,1]^d$. We define a probability measure over $\cU$ via tensor products. Abusing notation, we write $\varrho = \varrho \times \cdots \times \varrho$ for this measure,
where on the right-hand side $\varrho$ denotes the probability measure on $[-1,1]$. 

Let  $L^2_{\varrho}(\cU)$ be the Lebesgue space of square-integrable functions $f : \cU \rightarrow \bbC$.
We construct an orthonormal polynomial basis for this space via tensor products. Writing $\bm{\nu} = (\nu_1,\ldots,\nu_d) \in \bbN^d_0$
for an arbitrary multi-index, we define 
\bes{
\Psi_{\bm{\nu}}(\bm{y}) = \psi_{\nu_1}(y_1) \cdots \psi_{\nu_d}(y_d),\qquad \bm{y} = (y_i)^{d}_{i=1} \in \cU.
}
The set $\{ \Psi_{\bm{\nu}} \}_{\bm{\nu} \in \bbN^d_0} \subset L^2_{\varrho}(\cU)$ forms an orthonormal basis. Hence any function $f \in L^2_{\varrho}(\cU)$ has the convergent expansion
\be{
\label{f-exp}
f = \sum_{\bm{\nu} \in \bbN^d_0} c_{\bm{\nu}} \Psi_{\bm{\nu}},\qquad \text{where }
c_{\bm{\nu}} = \ip{f}{\Psi_{\bm{\nu}}}_{L^2_{\varrho}(\cU)} = \int_{\cU} f(\bm{y}) \overline{\Psi_{\bm{\nu}}(\bm{y})} \D \varrho(\bm{y}).
}
We are interested in $n$-term approximations to such functions. An \textit{$n$-term approximation} to $f$ based on an multi-index set $S \subset \bbN^d_0$, $|S| = n$, has the form
\be{
\label{f-S-def}
f_{S} = \sum_{\bm{\nu} \in S} c_{\bm{\nu}} \Psi_{\bm{\nu}}.
}
Due to Parseval's identity, the error of such an approximation is determined by the size of the coefficients $c_{\bm{\nu}}$ not included in $S$. Specifically,
\be{
\label{Parseval}
\nmu{f - f_S}^2_{L^2_{\varrho}(\cU)} = \sum_{\bm{\nu} \notin S} | c_{\bm{\nu}} |^2.
}

\subsection{{Choosing $S$}}\label{ss:choosing-S}

{
In standard multivariate polynomial approximation it is common to make an  \emph{a priori} choice of $S$ with some prescribed structure (see, e.g., \cite{guo2020constructing}). Given $\ell \in \mathbb{N}_0$, classical examples are the \emph{tensor-product} index set of order $\ell$,
$$
S^{\mathsf{TP}}_\ell = 
\left\{\bm{\nu} = (\nu_k)_{k=1}^d \in \mathbb{N}_0^d : \max  \{ \nu_1,\ldots,\nu_d \} \leq \ell \right\},
$$
or the \emph{total degree} index set of order $\ell$, 
$$
S^{\mathsf{TD}}_\ell = 
\left\{\bm{\nu} = (\nu_k)_{k=1}^d \in \mathbb{N}_0^d : \nu_1 + \cdots + \nu_d \leq \ell \right\}.
$$
However, the cardinality of these index sets scales exponentially with the dimension $d$, making them poorly suited in all but low-dimensional problems. Specifically, $|S^{\mathsf{TP}}_\ell| = (\ell+1)^d$ and $|S^{\mathsf{TD}}_\ell| = {\ell+d \choose d} $ for any $\ell \in \mathbb{N}_0$ and $d \in \mathbb{N}$. Applying Stirling's formula, one deduces that $|S^{\mathsf{TD}}_\ell| \sim \ell^d/d!$ as $\ell \to \infty$ for fixed $d$. This effect can be partially mitigated by working with \emph{hyperbolic cross} index sets (see, e.g., \cite{guo2020constructing}). Yet, the cardinalities of these sets also grow exponentially with dimension.
}

{Alternatively, one can look to replace such \textit{isotropic} sets with anisotropic variants. High-dimensional functions may be highly \textit{anisotropic}, i.e., they may depend more strongly on some variables than others. The aforementioned index sets fail to capture this behaviour, thus potentially leading to poor approximations. However, selecting a good anisotropic set \textit{a priori} requires detailed knowledge about the function $f$, which is typically unavailable in practice.}

{A theoretical alternative is provided by the concept of best $n$-term approximation. Here, one does away with the classical notion of polynomial degree (as described by the parameter $\ell$ in $S^{\mathsf{TP}}_{\ell}$ and $S^{\mathsf{TD}}_{\ell}$) and instead seeks an optimal index set in a function-dependent manner by minimizing the error \eqref{Parseval}} over all possible sets. A \textit{best $n$-term approximation} $f_n$ of $f$ is thus defined as
\be{
\label{best-n-term-inf}
f_n = f_{S^*},\qquad \text{where }S^* \in \argmin{} \{ \nmu{f - f_S}_{L^2_{\varrho}(\cU)} : S \subset \bbN^d_0,\ |S| = n \}.
}
It follows straightforwardly from \eqref{Parseval} that $S^*$ consists of those multi-indices corresponding to the $n$ largest coefficients $c_{\bm{\nu}}$ of $f$ in absolute value; that is to say, $S^* = \{ \bm{\nu}_1,\bm{\nu}_2,\ldots,\bm{\nu}_n \}$, where $\bm{\nu}_1,\bm{\nu}_2,\ldots$ are such that $|c_{\bm{\nu}_1} | \geq | c_{\bm{\nu}_2} | \geq \cdots $.

{Unfortunately,} the best $n$-term approximation is a theoretical benchmark. It is generally impossible to construct, since doing so would generically involve computing and sorting infinitely-many coefficients -- something that clearly cannot be done from the finite data \eqref{sample-values}. {A more practical approach is therefore to} construct a set $S$ -- or more precisely, a nested sequence of set $S^{(1)} \subseteq S^{(2)} \subseteq \cdots $ -- {in a function dependent, adaptive manner}. This is typically done {via a greedy procedure}, with multi-indices being added according to some importance criterion. The aforementioned ALS method is a procedure of this type. We discuss it further in \S \ref{s:ALS} and \S \ref{ss:adaptive-LS}. 

Regardless of how the index set (or sets) is constructed, however, it is generally useful to restrict one's attention sets with certain structure. A \textit{lower} set (also known as \textit{monotone} or \textit{downward closed} -- see, e.g., {\cite{cohen2015approximation}}) is a set $S \subseteq \bbN^d_0$ for which
\bes{
(\bm{\nu} \in S\text{ and }\bm{\mu} \leq \bm{\nu})\ \Rightarrow\ \bm{\mu} \in S,
}
(here the inequality $\bm{\mu} \leq \bm{\nu}$ is understood componentwise). {Lower sets are ubiquitous in multivariate polynomial approximation, with (isotropic or anisotropic)} tensor-product, total degree and hyperbolic cross index sets all being examples. Lower sets are also commonly employed in adaptive strategies such as ALS (see \S \ref{ss:adaptive-LS}), so as to make the greedy selection procedure tractable. Indeed, given a lower set $S^{(l)}$ there are only finitely many multi-indices $\bm{\nu} \in \bbN^d_0 \backslash S^{(l)}$ for which $S^{(l)} \cup \{ \bm{\nu} \}$ is also lower.

\subsection{Weighted least-squares polynomial approximation}

Fix a set multi-indices $S \subset \bbN^d_0$, $|S| = n$, and consider $m \geq n$ noisy samples \eqref{sample-values} 
 of a function $f \in L^2_{\varrho}(\cU)$ at \textit{sample points} $\bm{y}_1,\ldots,\bm{y}_m \in \cU$. We wish to compute an approximation to $f$ from the polynomial space defined by $S$: namely, the subspace
\bes{
\cP_S = \spn \{ \Psi_{\bm{\nu}} : \bm{\nu} \in S \} \subset L^2_{\varrho}(\cU).
}
Given a positive and almost everywhere finite  weight function $w : \cU \rightarrow \bbR$, we define a \textit{weighted least-squares approximation} to $f$ as
\be{
\label{wLS-def}
\hat{f} \in \argmin{p \in \cP_S} {\frac1m} \sum^{m}_{i=1} w(\bm{y}_i) | f(\bm{y}_i) + e_i - p(\bm{y}_i) |^2.
}
This is readily computed by solving an algebraic LS problem for the coefficients of $\hat{f}$. Indeed, let $\bm{\nu}_1,\ldots,\bm{\nu}_n$ be an enumeration of the indices in $S$. Then
\bes{
\hat{f} = \sum^{n}_{i=1} \hat{c}_i \Psi_{\bm{\nu}_i},\quad \text{where }\hat{\bm{c}} = (\hat{c}_i)^{n}_{i=1}  \in \argmin{\bm{z} \in \bbC^n} \nm{\bm{A} \bm{z} - \bm{f}}_2^{2},
}
and the LS matrix $\bm{A} \in \bbC^{m \times n}$ and vector $\bm{f} \in \bbC^m$ are given by
\be{
\label{LS-meas-mat-vec}
\bm{A} = \left ( \sqrt{w(\bm{y}_i) / m} \Psi_{\bm{\nu}_j}(\bm{y}_i) \right )^{m,n}_{i,j=1},\qquad \bm{f} = \left ( \sqrt{w(\bm{y}_i)/m} (f(\bm{y}_i) + e_i) \right )^{m}_{i=1}.
}

\section{Theory of (weighted) least-squares approximation}\label{s:theoryLS}

In this section, we present some elementary theory for weighted LS approximation. For the sake of generality, in the majority of this section we consider an arbitrary $n$-dimensional subspace $\cP \subset L^2_{\varrho}(\cU)$ (i.e., not necessarily a polynomial subspace of the form $\cP = \cP_S$). We write $\{ \Psi_i \}^{n}_{i=1}$ for an orthonormal basis for $\cP$.
However, for convenience, we make the mild assumption that the constant function $p(\bm{y}) = 1$, $\forall \bm{y} \in \cU$,  is an element of $\cP$. Note that this always holds when $\cP = \cP_S$ and $S$ is a lower set.

\subsection{Accuracy and stability}

Given sample points $\bm{y}_1,\ldots,\bm{y}_m \in \cU$ and a weight function $w : \cU \rightarrow \bbR$, define the discrete semi-norm
\bes{
\nmu{f}^2_{\mathsf{disc},w} = \frac1m \sum^{m}_{i=1} w(\bm{y}_i) | f(\bm{y}_i) |^2,\quad \forall f \in L^2_{\varrho}(\cU) \cap C(\cU),
}
(here $C(\cU)$ is the set of continuous functions on $\cU$)
and the \textit{discrete stability} constants
\eas{
\alpha_{w} &= \inf \left \{ \nmu{p}_{\mathsf{disc},w} : p \in \cP,\ \nm{p}_{L^2_{\varrho}(\cU)} = 1 \right \},
\\
 \beta_w &=  \sup \left \{ \nmu{p}_{\mathsf{disc},w} : p \in \cP,\ \nm{p}_{L^2_{\varrho}(\cU)} = 1 \right \},
}
Notice that $\alpha_{w} = \sigma_{\min}(\bm{A})$ and $\beta_{w} = \sigma_{\max}(\bm{A})$, where is the LS matrix \eqref{LS-meas-mat-vec}.

\lem{
[Accuracy and stability of weighted LS]
\label{l:wLS-acc-stab}
Let $\cP \subset L^2_{\varrho}(\cU) \cap C(\cU)$ with $1 \in \cP$, $f \in  L^2_{\varrho}(\cU) \cap C(\cU)$, $\bm{e} \in \bbC^m$, $\bm{y}_1,\ldots,\bm{y}_m \in \cU$ and $w : \cU \rightarrow \bbR$ be such that $w(\bm{y}_i) > 0$, $\forall i \in [m] : = \{1,\ldots,m\}$. If $\alpha_w > 0$ then the problem
\be{
\label{wLS-general}
\min_{p \in \cP} {\frac1m} \sum^{m}_{i=1} w(\bm{y}_i) | f(\bm{y}_i) + e_i - p(\bm{y}_i) |^2
}
has a unique solution $\hat{f}$. This solution satisfies
\be{
\label{LSerrbd}
\nmu{f - \hat{f}}_{L^2_{\varrho}(\cU)} \leq \inf_{p \in \cP} \left \{ \nmu{f - p}_{L^{2}_{\varrho}(\cU)} + \frac{1}{\alpha_w} \nmu{f - p}_{\mathsf{disc},w} \right \} + \frac{\beta_w}{\alpha_w} \nm{\bm{e}}_{\infty}.
}
Also, the {($2$-norm)} condition number of the LS matrix satisfies $\mathrm{cond}(\bm{A}) = \beta_w / \alpha_w$.
}

This result is standard {(see, e.g., \cite[Prop.\ 1]{migliorati2014analysis})}. We include a short proof in \S \ref{s:proofs} for completeness. Note that the last statement is immediate, since $\mathrm{cond}(\bm{A}) = \sigma_{\max}(\bm{A}) / \sigma_{\min}(\bm{A})$ by definition. 
The main takeaway from this lemma is that the accuracy and stability of $\hat{f}$ are determined by the size of the constants $\alpha_w$ and $\beta_w$. In the next subsection, we discuss how to control these constants in the case of random sampling. In doing so, we also obtain an estimate for the sample complexity.

\subsection{Sample complexity}

We now consider the sample complexity of weighted LS approximation in the case of random sampling. Specifically, we now assume that $\bm{y}_1,\ldots,\bm{y}_m \sim_{\mathrm{i.i.d.}} \mu$, where $\mu$ is some probability measure with support in $\cU$.

The analysis of this type of sampling strategy involves the (reciprocal) \textit{Christoffel function} of the subspace $\cP$ \cite{cohen2013stability}:
\bes{
\cK(\cP)(\bm{y}) = \sup \left \{ | p(\bm{y}) |^2 : p \in \cP,\ \nm{p}_{L^2_{\varrho}(\cU)} = 1 \right \},\quad \forall \bm{y} \in \cU.
}
Observe that $\cK(\cP)(\bm{y}) \geq 1$, $\forall \bm{y} \in \cU$, since, by assumption, the function $1 \in \cP$.
It is also a short argument to show that $\cK(\cP)$ has the equivalent expression 
\be{
\label{K-alternative-def}
\cK(\cP)(\bm{y}) = \sum^{n}_{i=1} | \Psi_i(\bm{y}) |^2,\quad \forall \bm{y} \in \cU,
}
where $\{ \Psi_i \}^{n}_{i=1}$ is any orthonormal basis for $\cP$.
Given $\cK$ and a weight function $w : \cU \rightarrow (0,\infty)$, we now also define
\be{
\label{kappa-def}
\kappa(\cP ; w) = \sup_{\bm{y} \in \cU} w(\bm{y}) \cK(\cP)(\bm{y}).
}

\thm{
[Sample complexity of weighted LS with random sampling]
\label{t:wLS-samp-comp}
Let $\cP \subset L^2_{\varrho}(\cU) \cap C(\cU)$ with $\dim(\cP) = n$ and $1 \in \cP$, $0 < \epsilon < 1$ and $\mu$ be a probability measure on $\cU$ such that
\be{
\label{mu-rho-w}
\D \mu(\bm{y}) = (w(\bm{y}))^{-1} \D \varrho(\bm{y})
}
for some strictly positive and finite almost everywhere weight function $w : \cU \rightarrow \bbR$. Let $\bm{y}_1,\ldots,\bm{y}_m \sim_{\mathrm{i.i.d.}} \mu$, where $m$ satisfies
\be{
\label{wLS-samp-comp}
m \geq 7 \cdot \kappa(\cP ; w) \cdot \log(2n/\epsilon).
}
Then the following holds with probability at least $1-\epsilon$. For any $f \in L^2_{\varrho}(\cU) \cap C(\cU)$, the solution $\hat{f}$ of the weighted LS problem \eqref{wLS-general} is unique and satisfies
\be{
\label{wLS-err-bd}
\nmu{f - \hat{f}}_{L^2_{\varrho}(\cU)} \leq \inf_{p \in \cP} \left \{ \nmu{f - p}_{L^2_{\varrho}(\cU)} + 2 \nmu{f-p}_{\mathsf{disc},w} \right \} + 2 \nmu{\bm{e}}_{\infty}.
} 
Moreover, the condition number satisfies $\mathrm{cond}(\bm{A}) \leq 2$.
}

{This result {is based on standard ideas, in particular, the use of the matrix Chernoff inequality (see, e.g., \cite[Thm.\ 2.1]{cohen2017optimal})}. We include a short proof in \S \ref{s:proofs} for completeness.}
Note that the number $7$ in \eqref{wLS-samp-comp} is somewhat arbitrary. As can be seen from the proof, one can replace it with a smaller constant at the expense of larger constants in the error and condition number bounds. 

The most important aspect of this result is the sample complexity bound \eqref{wLS-samp-comp}. Since $\mu$ is a probability measure, \eqref{mu-rho-w} implies that
\be{
\label{w-int-one}
\int_{\cU} (w(\bm{y}))^{-1} \D \varrho(\bm{y}) = 1,
}
and orthonormality and the alternative expression \eqref{K-alternative-def} for $\cK(\cP)$ imply that 
\be{
\label{K-int-n}
\int_{\cU} \cK(\cP)(\bm{y}) \D \varrho(\bm{y}) = n 
}
and
\bes{
\kappa(\cP ; w) = \int_{\cU}\kappa(\cP ; w) /w(\bm{y}) \D \varrho(\bm{y}) \geq \int_{\cU} \cK(\cP)(\bm{y}) \D \varrho(\bm{y}) = n.
}
Hence $\kappa(\cP ; w) \geq n $ for any $\cP$ and $w$. As a result, the sample complexity bound \eqref{wLS-samp-comp} is always at least log-linear in $n$. In the next subsection, we show that it is generally superlinear in $n$. However, in \S \ref{ss:opt-samp-scheme} we show that it is always possible to choose the sampling measure $\mu$ so as to achieve the optimal value of the right-hand side of \eqref{wLS-samp-comp}.

\subsection{The case of Monte Carlo sampling}\label{ss:MC-samp-scheme}

MC sampling corresponds to the choice $\mu = \varrho$, in which case \eqref{mu-rho-w} holds with $w = 1$. The approximation \eqref{wLS-general} is correspondingly an \textit{unweighted} LS approximation.
In this case, the sample complexity estimate \eqref{wLS-samp-comp} takes the form
\be{
\label{MC-samp-comp-general}
m \geq c \cdot \kappa(\cP) \cdot \log(n/\epsilon),\quad \text{where }\kappa(\cP) = \kappa(\cP ; 1) = \nm{\cK(\cP)}_{L^{\infty}(\cU)}.
}
Recall that $\kappa(\cP) \geq n = \dim(\cP)$. Unfortunately, $\kappa(\cP)$ can be arbitrarily large, even in cases such as Chebyshev and Legendre polynomial approximation.
We now state two standard results, {which are based on \cite{chkifa2015discrete,chkifa2014high,chkifa2018polynomial} (see also \cite[Props.\ 5.13 \& 5.17]{adcock2022sparse}).}

\prop{
\label{prop:chebyshev-Christoffel}
Let {$d, n \in \mathbb{N}$}, $\varrho$ be the Chebyshev (arcsine) measure over $\cU = [-1,1]^d$ and $\{ \Psi_{\bm{\nu}} \}_{\bm{\nu} \in \bbN^d_0} \subset L^2_{\varrho}(\cU)$ be the orthonormal Chebyshev polynomial basis. Then
\be{
\label{eq:K_CC}
\max \left \{ \kappa(\cP_{S}) : S \subset \bbN^d_0, | S | \leq n \right \} = 2^{d} n.
}
However, for lower sets, if $1 \leq n \leq 2^{d+1}$ then
\be{
\label{eq:K_CC_lower}
n^{\log(3)/\log(2)}/3 
\leq \max \left \{ \kappa(\cP_{S}) : S \subset \bbN^d_0, | S | \leq n,\ \textnormal{$S$ lower} \right \} \leq n^{\log(3)/\log(2)}.
}
Moreover, the upper bound holds for any $n \geq 1$.
}

\prop{
\label{prop:legendre-Christoffel}
Let {$d, n \in \mathbb{N}$}, $\varrho$ be the uniform measure over $\cU = [-1,1]^d$ and $\{ \Psi_{\bm{\nu}} \} \subset L^2_{\varrho}(\cU)$ be the orthonormal Legendre polynomial basis. Then
\be{
\label{eq:K_LU}
\max \left \{ \kappa(\cP_{S}) : S \subset \bbN^d_0,\ | S | \leq n,\ \textnormal{$S$ lower} \right \} = n^2.
}
However, $\kappa(\cP_{S})$ is unbounded for arbitrary sets $S \subset \bbN^d_0$ of size $|S| \leq n$. That is, for every $c >0$ there exists a set $S$ of size $|S| \leq n$ for which $\kappa(\cP_{S}) \geq c$.
}

These two results suggest that MC sampling generally suffers from a poor sample complexity in the case of Chebyshev or Legendre polynomials. For non-lower sets the situation can be arbitrarily bad in the case of Legendre polynomial approximation, and in the case of Chebyshev polynomial approximation, highly susceptible to the curse of dimensionality. Even for lower sets, these results in combination with Theorem \ref{t:wLS-samp-comp} suggest the superlinear sample complexity bound $m \geq 7 \cdot n^{\gamma} \cdot \log(2n/\epsilon)$,
where $\gamma = \log(3)/\log(2) $ (Chebyshev) or $\gamma = 2$ (Legendre).

\subsection{Near-optimal sampling}\label{ss:opt-samp-scheme}

We now describe a random sampling scheme that obtains provably log-linear sample complexity. This was introduced in \cite{cohen2017optimal}; see also \cite{hampton2015coherence} for earlier work in this direction. 

The idea is to choose the 
weight function $w$, and therefore, via \eqref{mu-rho-w}, the sampling measure $\mu$, so as to minimize the constant $\kappa(\cP ; w)$ appearing in the sample complexity bound \eqref{wLS-samp-comp}.
It is immediate from \eqref{kappa-def} that $\kappa(\cP ; w)$ is minimized whenever the weight function $w \propto (\cK(\cP))^{-1}$. Recall that $w$ and $\cK(\cP)$ satisfy \eqref{w-int-one} and \eqref{K-int-n}, respectively. Hence, the choice of $w$ that minimizes $\kappa(\cP ; w)$ is precisely
\be{
\label{opt-w}
w(\bm{y}) = \left ( n^{-1} \cK(\cP)(\bm{y}) \right )^{-1},
}
and the corresponding sampling measure is
\be{
\label{opt-samp-meas}
\D \mu(\bm{y}) =  n^{-1} \cK(\cP)(\bm{y}) \D \varrho(\bm{y}) = n^{-1} \sum^{n}_{i=1} | \Psi_i(\bm{y}) |^2  \D \varrho(\bm{y}).
}
Since $\kappa(\cP ; w) = n$ in this case, the sample complexity bound \eqref{wLS-samp-comp} reads as
\be{
\label{opt-samp-comp}
m \geq 7 \cdot n \cdot \log(2n / \epsilon).
}
In other words, it is optimal up to the constant $7$ and the logarithmic factor $\log(2n/\epsilon)$.

Notice that this property holds for \textit{any} subspace $\cP$, regardless of whether it is a polynomial space or not. However, when applied to Chebyshev or Legendre polynomial approximation, it reduces the sample complexity from the superlinear rates asserted in Propositions \ref{prop:chebyshev-Christoffel} and \ref{prop:legendre-Christoffel} to log-linear in $n$.

{
\rem{
[Optimal sampling]
Recently, a series of works have strived to remove the logarithmic dependence in the sample complexity \eqref{opt-samp-comp} by carefully modifying the sampling strategy and weight function {(see \cite{bartel2023constructive,krieg2021functionII,dolbeault2022optimal,temlyakov2021optimal} and references therein)}. These works use nontrivial constructions based on \cite{batson2014twice,marcus2015interlacing}. We do not consider such approaches in this paper. As discussed in \cite[\S 5]{dolbeault2022optimal}, of the existing approaches, those which are computationally feasible (i.e., implementable in polynomial time) are accompanied by error bounds that involve the $L^{\infty}$-norm error $\inf_{p \in \cP} \nmu{f - p}_{L^{\infty}(\cU)}$, this being an upper bound for the term $\nmu{f - p}_{\mathsf{disc},w}$ in \eqref{wLS-err-bd}. If applied to the approximation of $(\bm{b},\varepsilon)$-holomorphic functions, this would lead to suboptimal algebraic rates of the form $(m/c)^{1-\frac1p}$ (recall \S \ref{s:main-res-summary}). Later, when we prove our main result, we use Bernstein's inequality to estimate the term $\nmu{f - p}_{\mathsf{disc},w}$ in a more careful way to obtain near-optimal rates for MC sampling. See Theorem \ref{t:main-res} and \S \ref{s:proofs}.
}
}

\section{{Monte Carlo sampling is good in high dimensions}}\label{s:MC-bad}

The previous discussion suggests that MC sampling is a poor strategy, especially in the case of Legendre polynomials, and that significant improvements may be realized with the near-optimal sampling scheme of \S \ref{ss:opt-samp-scheme}. However, as we noted in \S \ref{ss:aims} it has often been observed that MC sampling performs relatively well in high dimensions. The purpose of this section is to demonstrate this phenomenon via numerical experiments.

{This} phenomenon was briefly investigated in \cite{cohen2017optimal}. Here, phase transition plots were produced showing the empirical probability $\bbP ( \mathrm{cond}(\bm{A}) \leq 3 )$ for randomly generated sequences of lower sets. As observed therein, in low dimensions MC sampling is significantly worse than the near-optimal scheme, but in higher dimensions the difference in performance is greatly reduced.

\subsection{Adaptive (weighted) least-squares approximation}\label{s:ALS}

Since the focus of this paper is on function approximation, in our experiments we compare MC and near-optimal sampling on several different smooth function approximation tasks.

As mentioned previously, we shall use an \textit{Adaptive LS (ALS)} approximation scheme \cite{migliorati2019adaptive}. The procedure is described in full detail in \S \ref{ss:adaptive-LS}. However, in short, it proceeds as follows. Starting from an initial index set $S^{(1)} = \{\bm{0}\}$, at step $l$ this scheme first computes a (weighted) LS approximation $\hat{f}^{(l)} \in \cP_{S^{(l)}}$ using the index set $S^{(l)}$, and then uses the set $S^{(l)}$ and the approximation $\hat{f}^{(l)}$ to construct the next index set $S^{(l+1)} \supseteq S^{(l)}$ in a greedy manner. It does this by using $\hat{f}^{(l)}$ to estimate the coefficients with indices {belonging} to the so-called \textit{reduced margin} of $S^{(l)}$, and then chooses those coefficients which are largest in magnitude.
The result of this procedure is a sequence of approximations $\hat{f}_1,\hat{f}_2,\ldots$ and nested multi-index sets $S^{(1)} \subseteq S^{(2)} \subseteq \cdots$ of sizes {$1 = n^{(1)} \leq n^{(2)} \leq \cdots$}, where $n^{(l)} = |S^{(l)}|$. We remark in passing that the generated multi-index sets are always lower sets, due to the use of the reduced margin.

\subsection{Experimental setup}

We now describe the salient aspects of the experiments. See \S \ref{s:additional-info} for further details.
In these experiments, we choose the number of samples at step $l$ to be log-linear in $n^{(l)}$, i.e.,
\be{
\label{m-scaling-LS}
m = \max  \{ n^{(l)} +1, \lceil n^{(l)} \cdot \log(n^{(l)}) \rceil  \}
}
(the use of the maximum simply ensures that $m > n^{(l)}$ for all $l$, so the LS problem is overdetermined). Thus, in view of Theorem \ref{t:wLS-samp-comp} and the discussion in \S  \ref{ss:opt-samp-scheme} we expect the near-optimal sampling scheme to be stable.

In this and other experiments, we consider the (relative) error of the LS approximation $\hat{f}^{(l)}$ and the condition number $\mathrm{cond}(\bm{A}^{(l)})$, where $\bm{A}^{(l)}$ is the $l$th LS matrix. For the former, we compute the discrete $L^2_{\varrho}$-norm error over a grid of $100,000$ MC points. To ensure a valid comparison, we also use this grid to generate the sample points. See \S \ref{ss:fine-grids-measures} for further discussion on this point.

Since all experiments involve random sampling, we perform $T = 50$ trials. In these experiments, each trial corresponds to a sequence of approximations $\hat{f}^{(1)},\hat{f}^{(2)},\ldots$ and LS matrices $\bm{A}^{(1)},\bm{A}^{(2)},\ldots$. Thus, each figure shows the statistics relating to the computed quantities. As discussed in \S \ref{ss:statistical-sims}, we plot a main curve showing the geometric mean over the trials, and a shaded region showing one (geometric) standard deviation. {These shaded regions play a role similar to \emph{standard errors} in the context of MC quadrature \cite{owen2013monte} and they are sometimes too small to be visible. When this occurs, it is due to an extreme concentration of the visualized data around its mean.}

\subsection{{Test functions}}

{In addition to the parametric DE considered in Fig.\ \ref{fig:fig1}, we also consider various other test functions. First, we consider
\be{
f_1(\bm{y})  =\exp \left ( \sum^{d}_{i=1} \frac{y_i }{2i} \right ),\quad \forall \bm{y} \in [-1,1]^d \label{f1}.
}
This function is entire and anisotropic with respect to its variables. For larger $d$ it is sharply peaked near its maximal and minimal values at $\bm{y} = \pm \bm{1}$, with $f_1(+\bm{1}) \sim 1.34 \sqrt{d}$, $f_{1}(-\bm{1}) \sim (1.34 \sqrt{d})^{-1}$ as $d \rightarrow \infty$, while $f_1(\bm{0}) = 1$. Next, we consider
\be{
f_2(\bm{y}) = \left ( 1 + \frac{1}{2d} \sum^{d}_{i=1} q_i y_i \right )^{-1},\quad \forall \bm{y} \in [-1,1]^d, \qquad \text{where }q_i = 10^{-\frac{3(i-1)}{d-1}}\label{f2}
}
(see, e.g., \cite{migliorati2019adaptive}).
This function is anisotropic and holomorphic, but not entire. Unlike $f_1$ it varies between finite maximum and minimum values, with $0.93 + \ord{1/d} \leq f_2(\bm{y}) \leq 1.08 + \ord{1/d}$ as $d \rightarrow \infty$.}

 {Next, we consider several examples from the family of test functions
\be{
\label{f3}
f_3(\bm{y}) = \prod^{d}_{i=1} \frac{\sqrt{2\delta_i+\delta^2_i}}{y_i+1+\delta_i},\quad \forall \bm{y} \in [-1,1]^d,
}
with positive parameters $\delta_i > 0 $ (the factor in the numerator ensures this function has unit norm, and therefore avoids scale effects for large $d$ due to the $d$-fold product). These functions are holomorphic in $[-1,1]^d$ with singularities at any $\bm{y}$ for which $y_i = -1 - \delta_i$ for some $i$. They can be either isotropic or anisotropic depending on the choice of $\delta_i$: a larger $\delta_i$ implies a smoother dependence on the variable $y_i$, and smaller $\delta_i$ implies a less smooth dependence. As we discuss later in Remark \ref{rem:f3-theory}, our later theoretical analysis is readily applied to such functions.}

{Finally, we consider a test function $f_{\mathsf{bor}}$ from the \textit{Virtual Library of Simulation Experiments} \cite{surjanovic2013virtual}, which serves as a physical model for water flow through a borehole:
\be{
\label{borehole}
f_4(\bm{y}) = f_{\mathsf{bor}}(\bm{y}),\quad \forall \bm{y} \in [-1,1]^d.
}
The original function in \cite{surjanovic2013virtual} depends on $d_{*} = 8$ parameters that vary in finite intervals. We scale these parameters to $[-1,1]$ and then, for $d < d_{*}$, set the last $d_{*}-d$ parameters equal to their maximum value to allow tests to be carried out for different values of $d$. See \S \ref{ss:further-test-suite} for further information and additional examples from this test suite.}

\subsection{Numerical experiments}

\begin{figure}[t!]
\begin{center}
\begin{small}
 \begin{tabular}{@{\hspace{0pt}}c@{\hspace{\errplotsp}}c@{\hspace{\errplotsp}}c@{\hspace{0pt}}}
\includegraphics[width = \errplotimg]{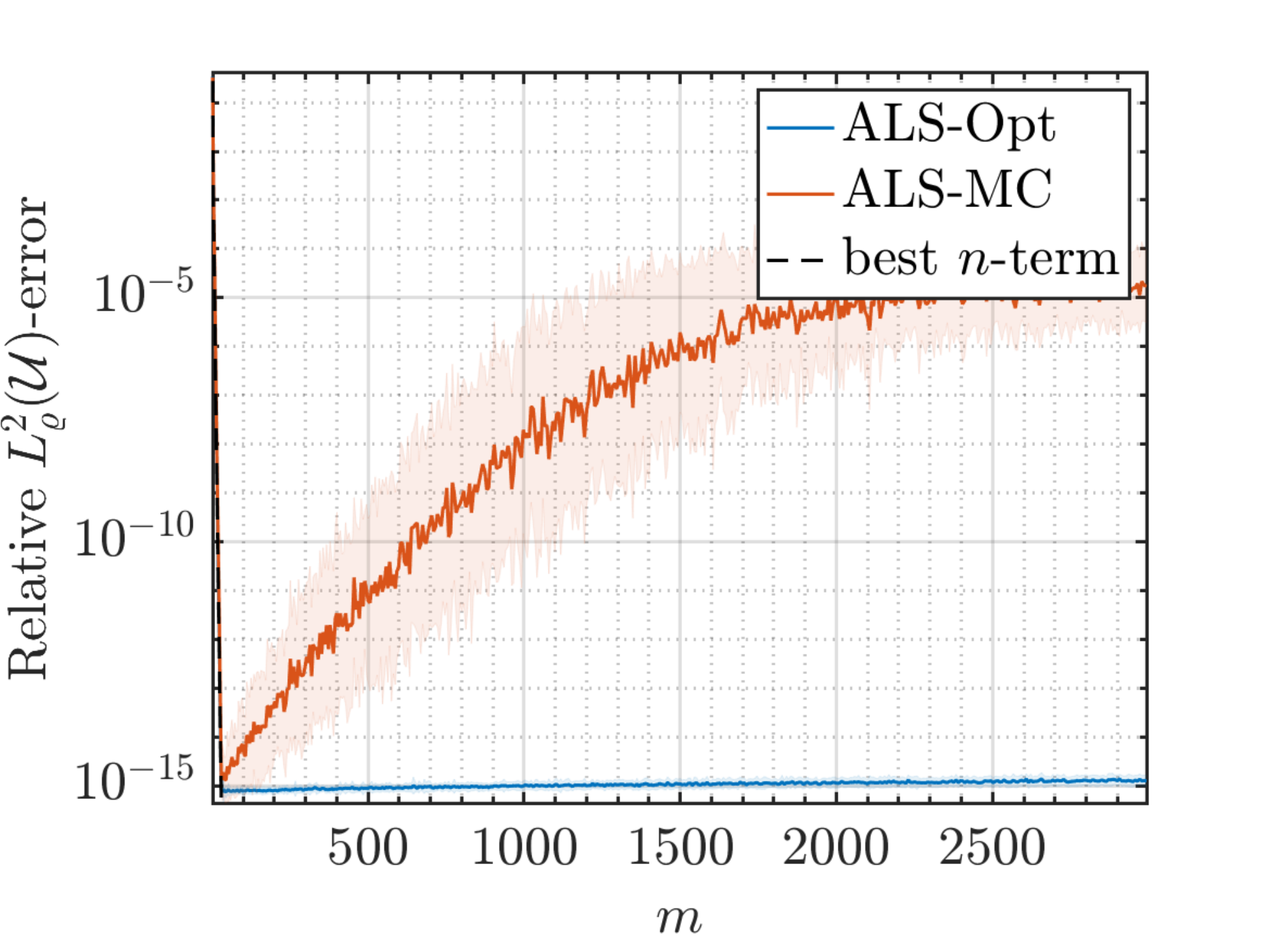}
&
\includegraphics[width = \errplotimg]{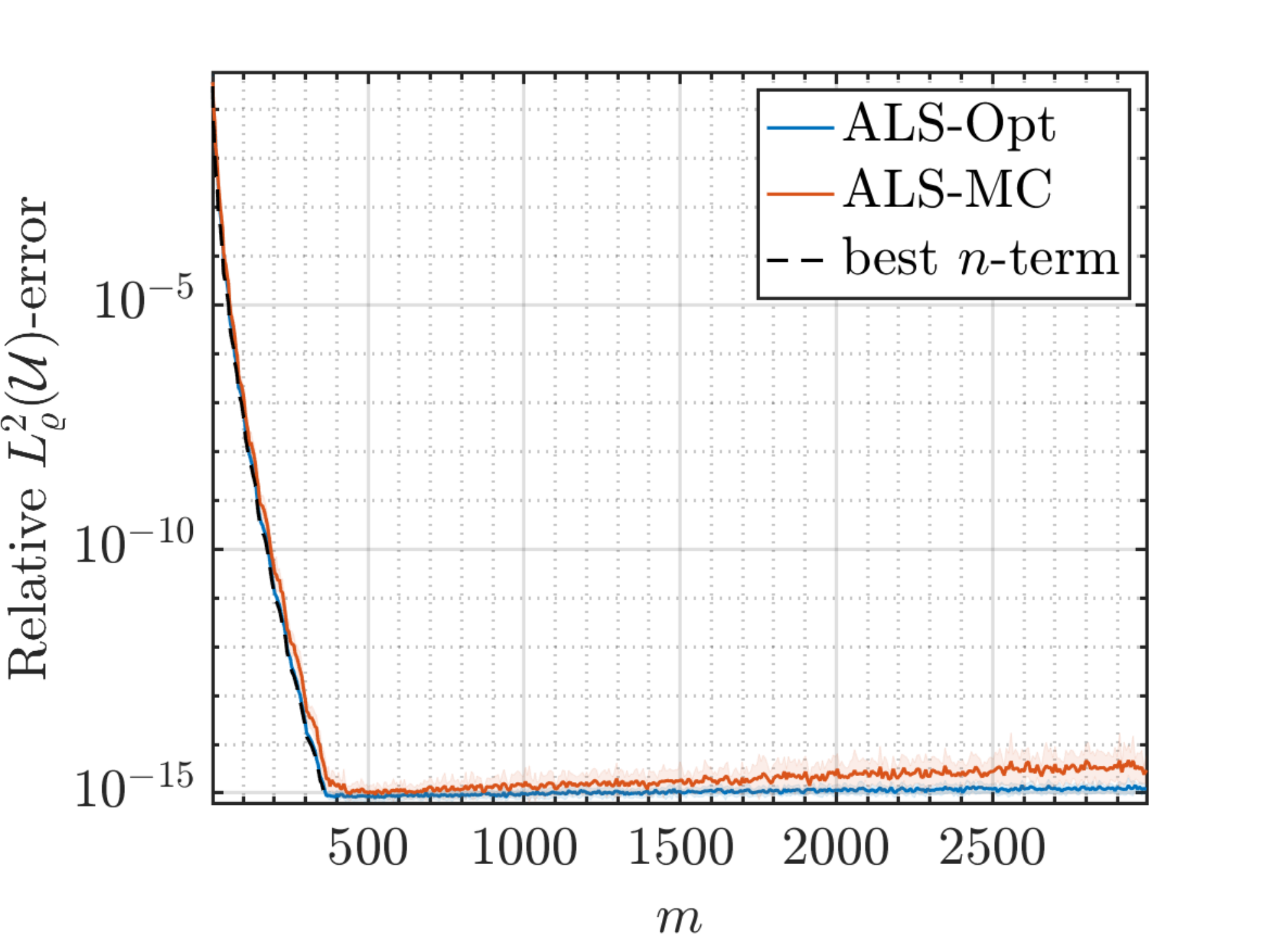}
&
\includegraphics[width = \errplotimg]{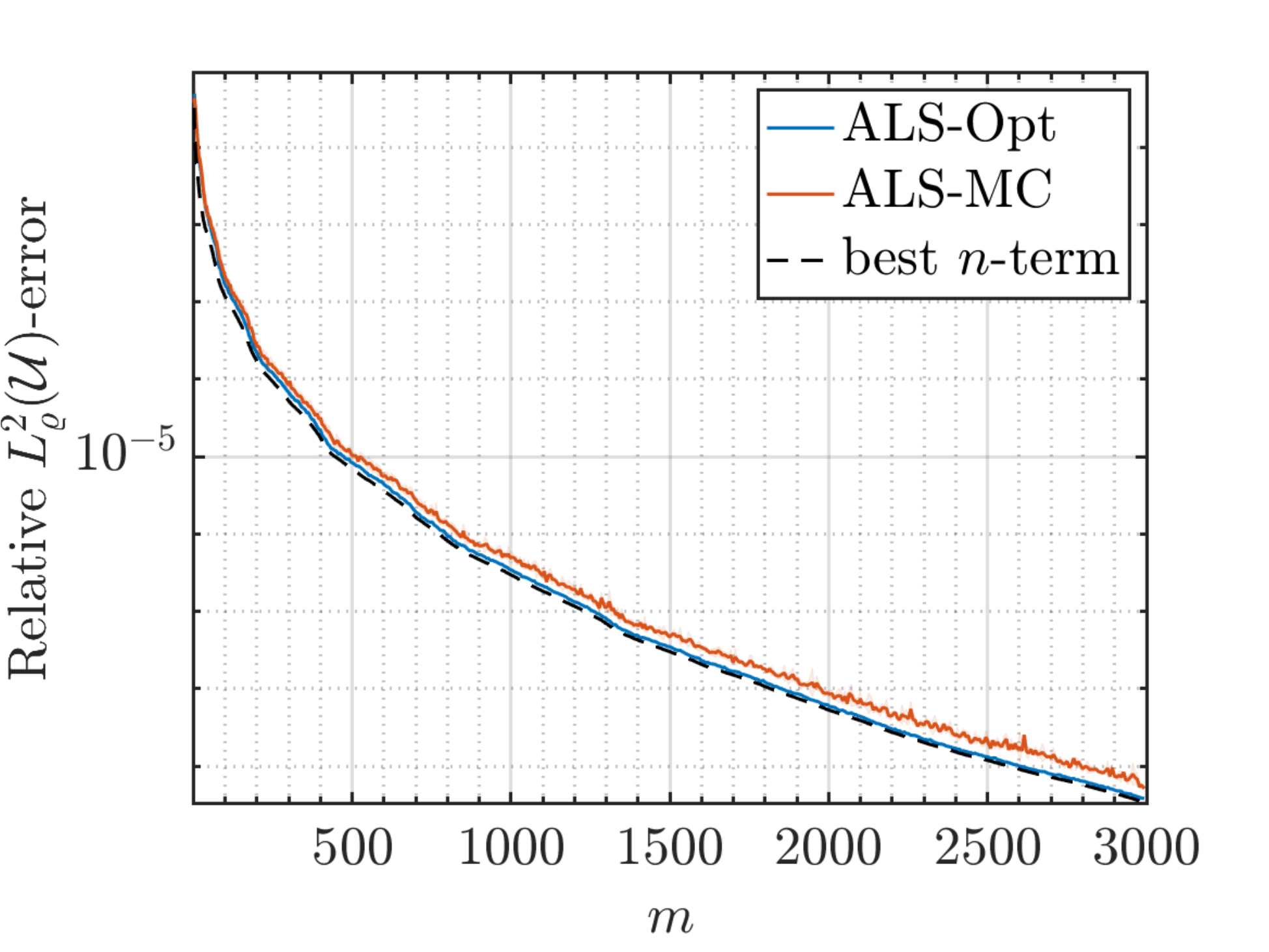}
\\[\errplotgraphsp]
\includegraphics[width = \errplotimg]{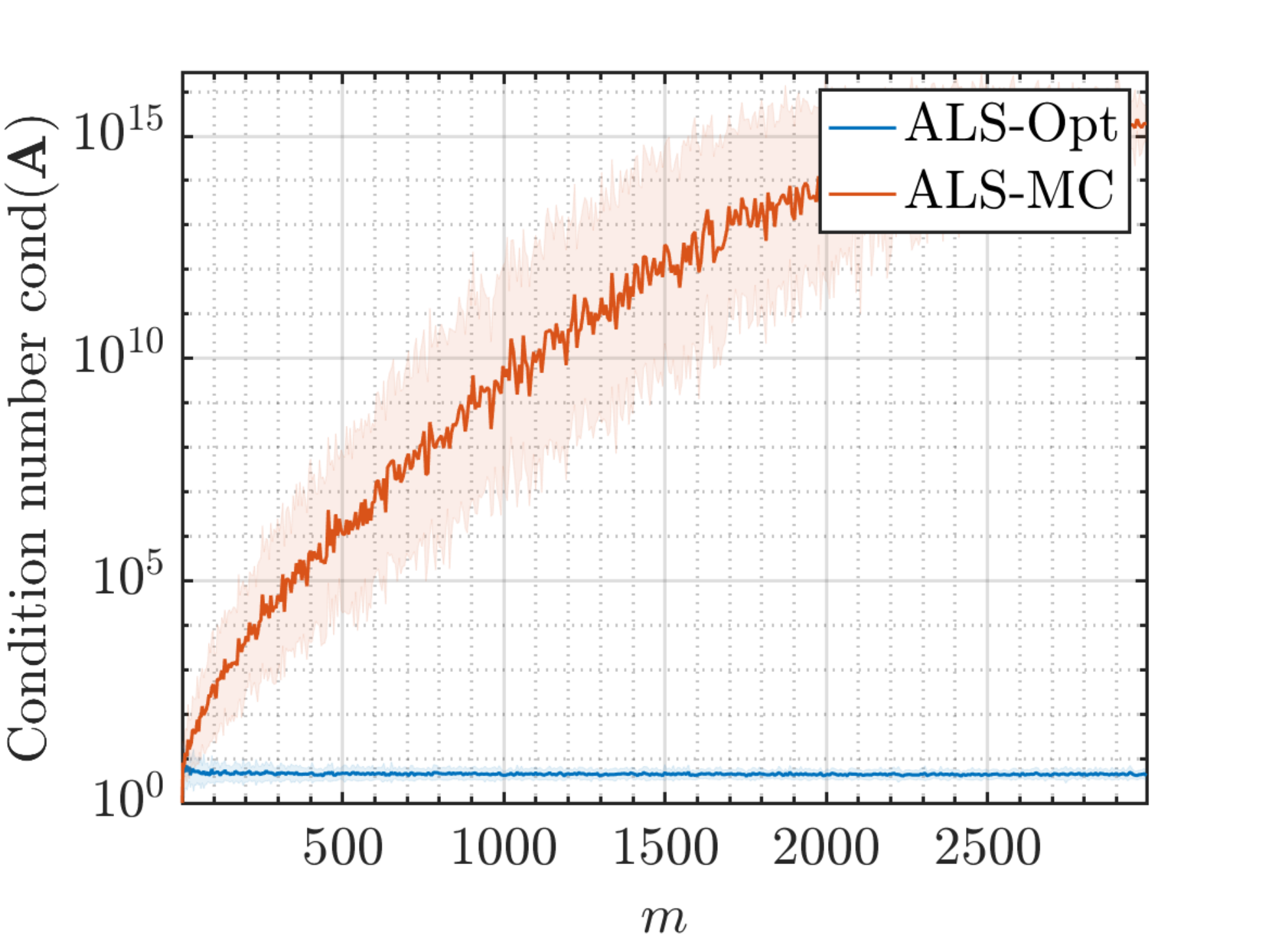}
&
\includegraphics[width = \errplotimg]{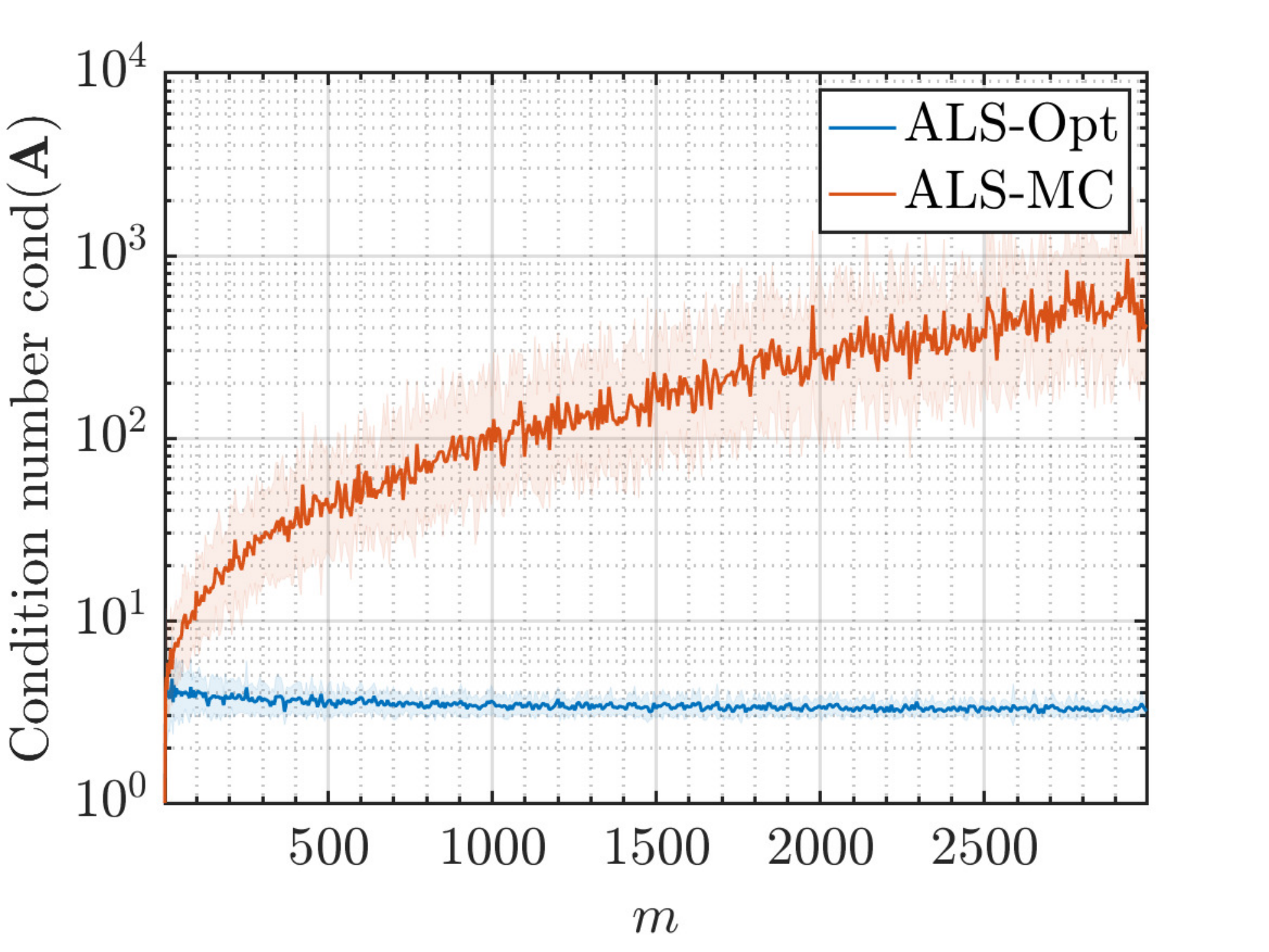}
&
\includegraphics[width = \errplotimg]{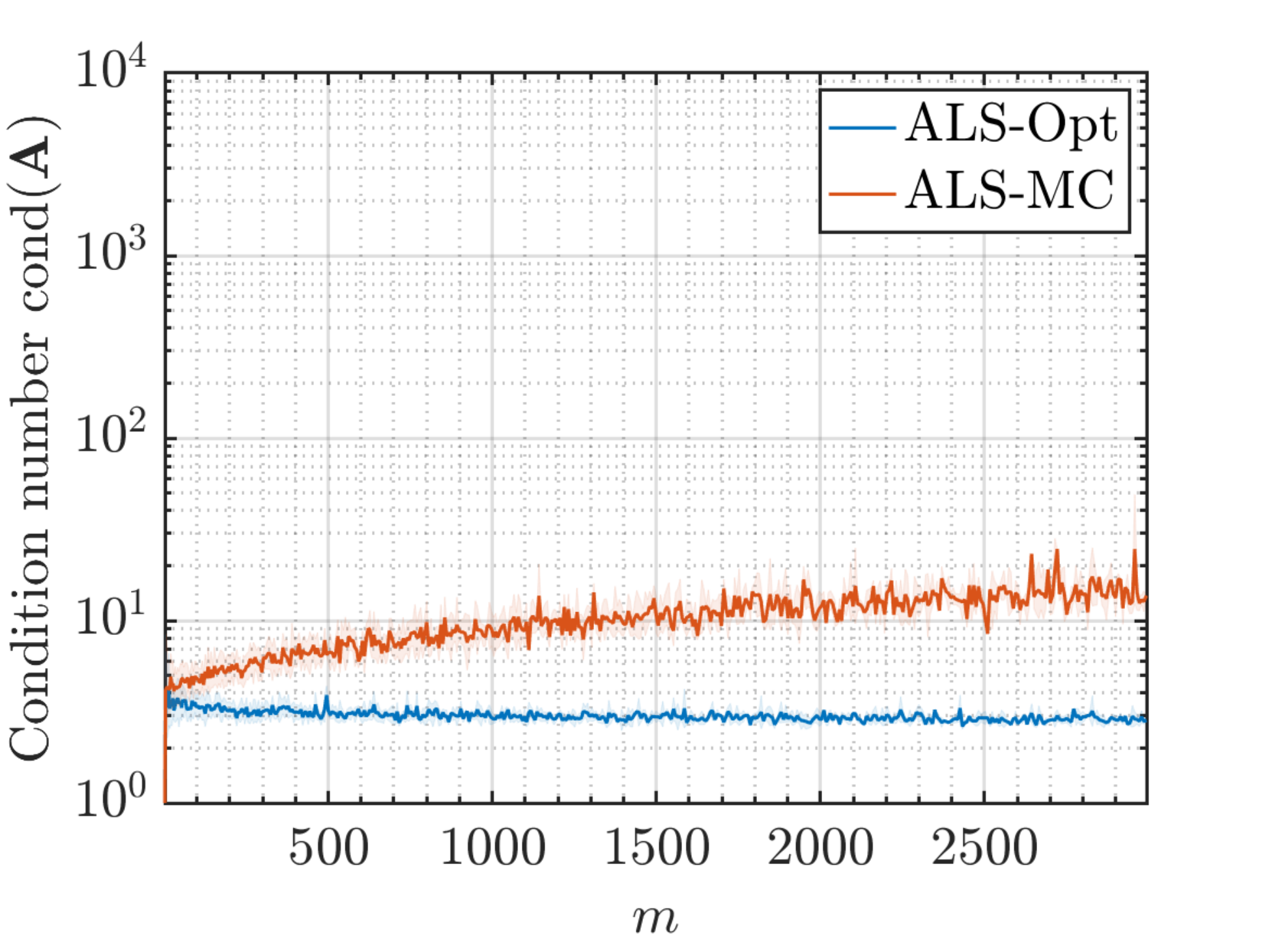}
\\[\errplottextsp]
$d = 1$ & $d = 2$ & $d = 4$
\\[\errplottextsp]
\includegraphics[width = \errplotimg]{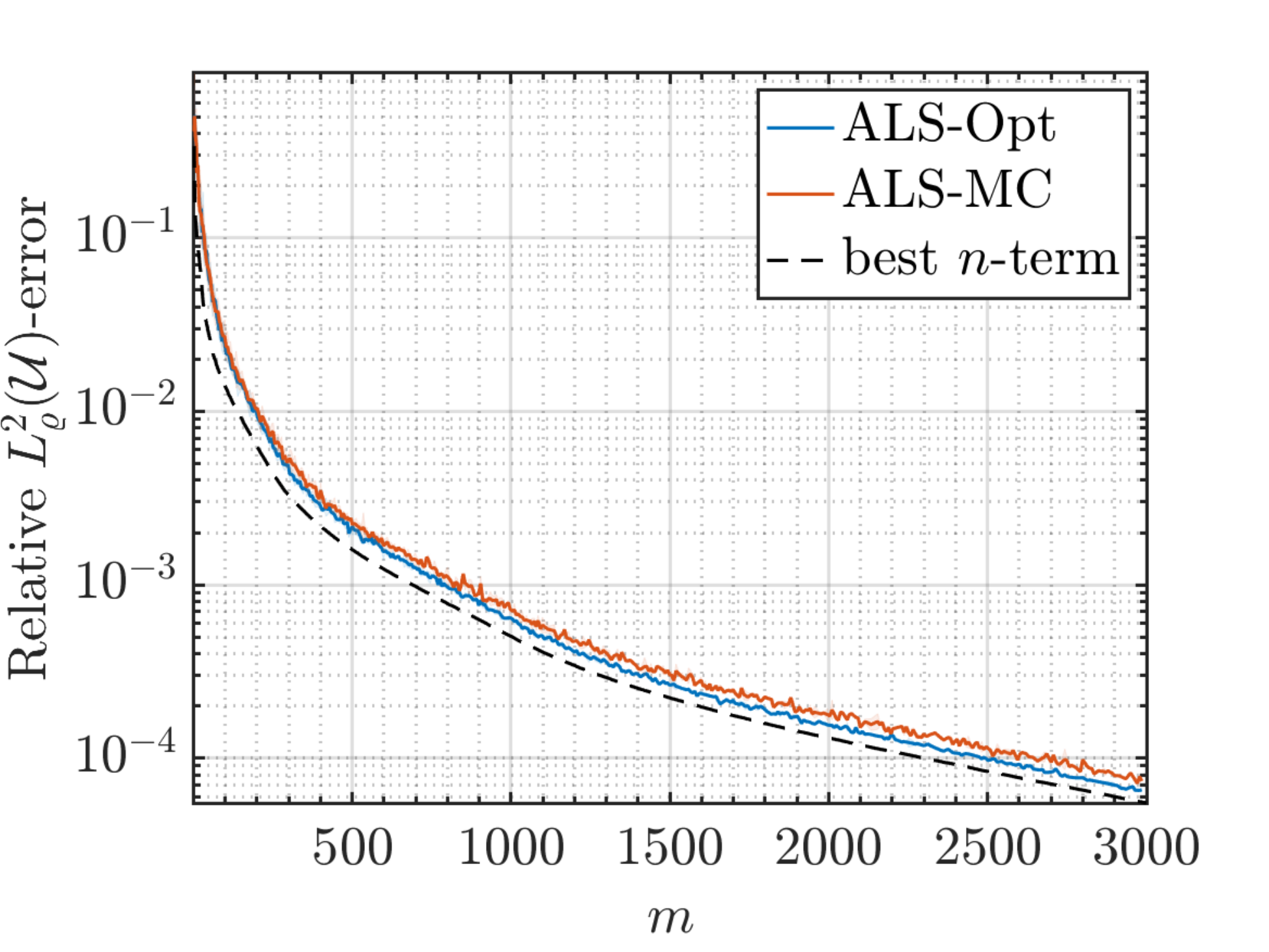}
&
\includegraphics[width = \errplotimg]{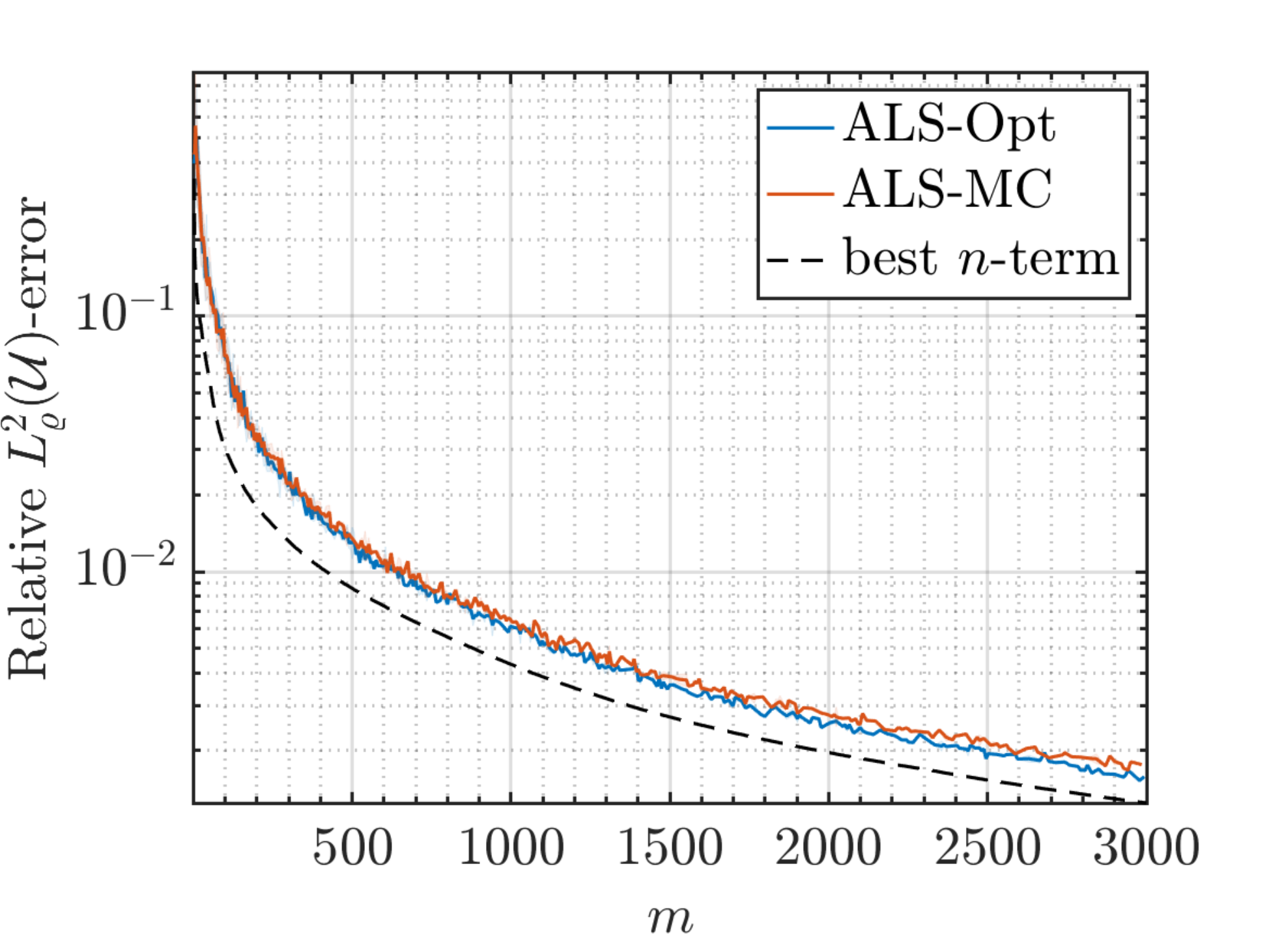}
&
\includegraphics[width = \errplotimg]{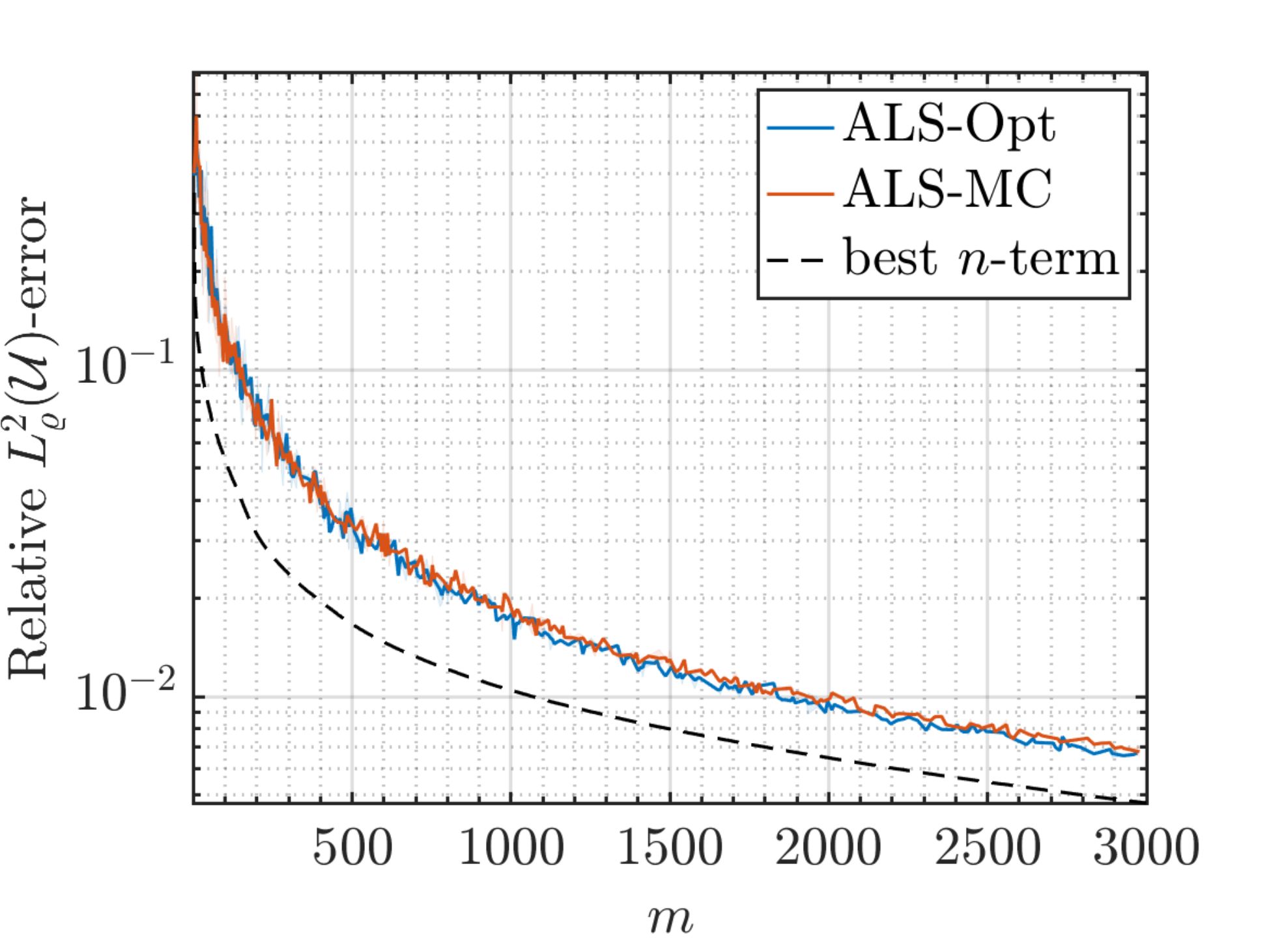}
\\[\errplotgraphsp]
\includegraphics[width = \errplotimg]{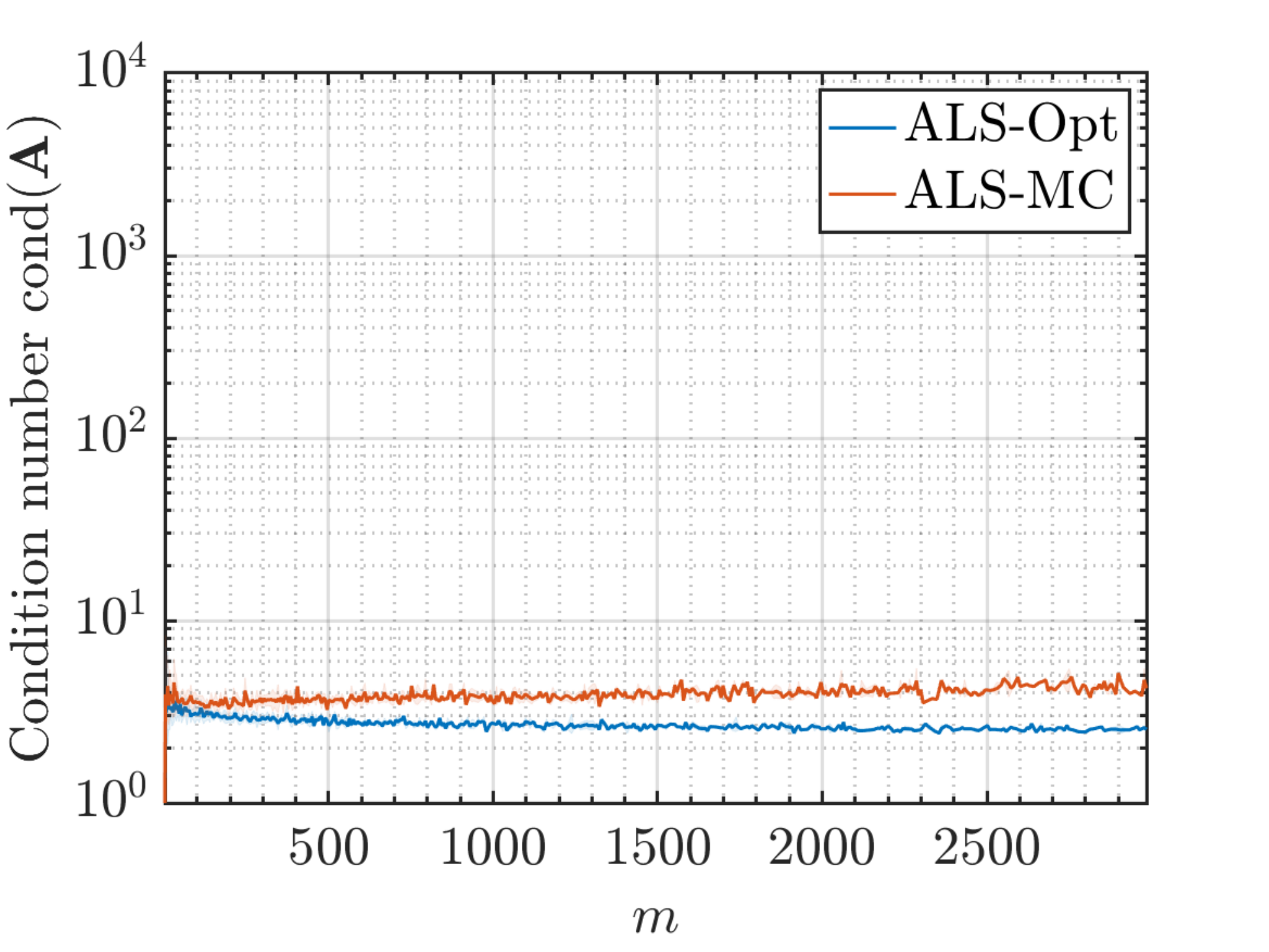}
&
\includegraphics[width = \errplotimg]{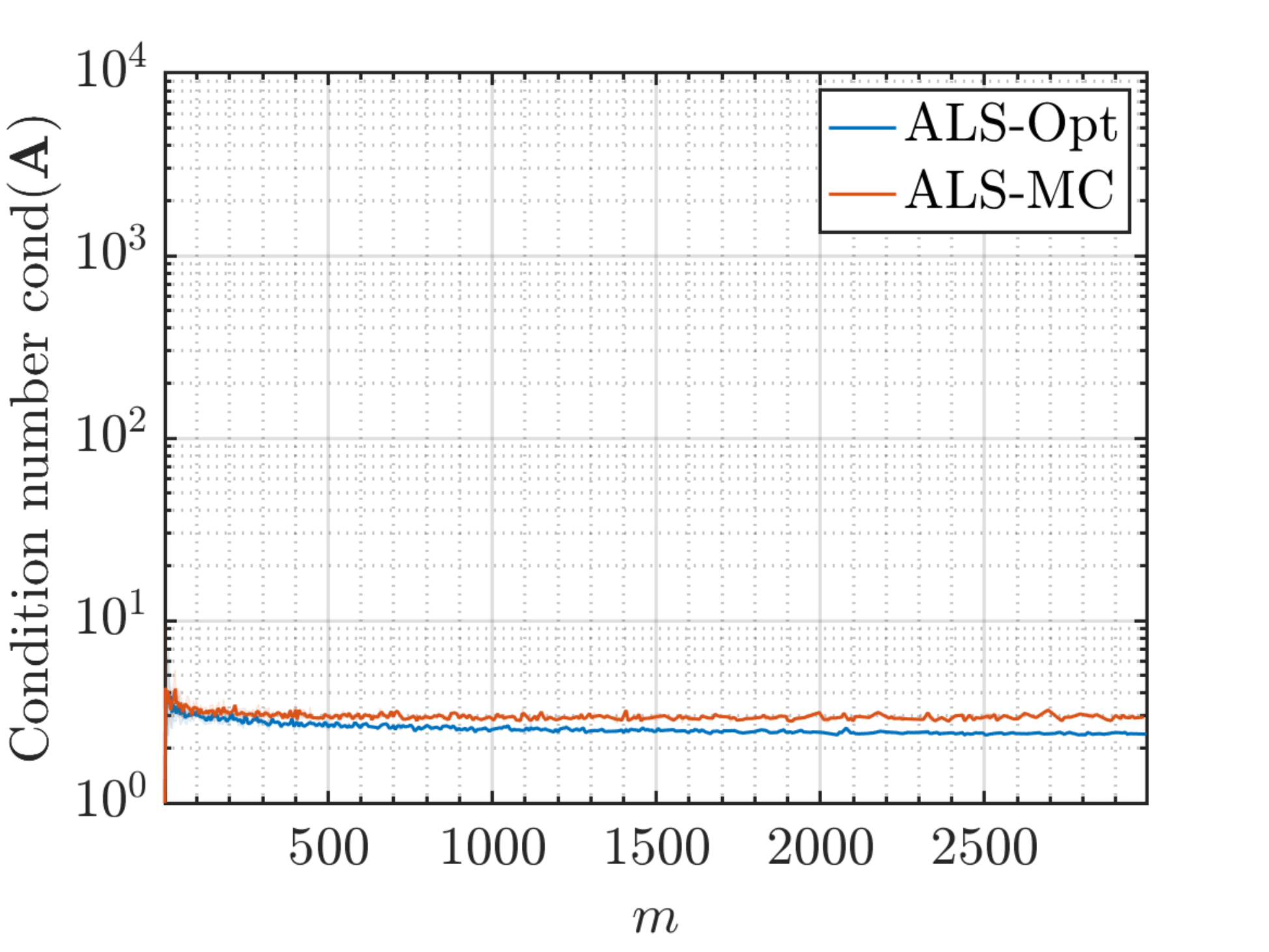}
&
\includegraphics[width = \errplotimg]{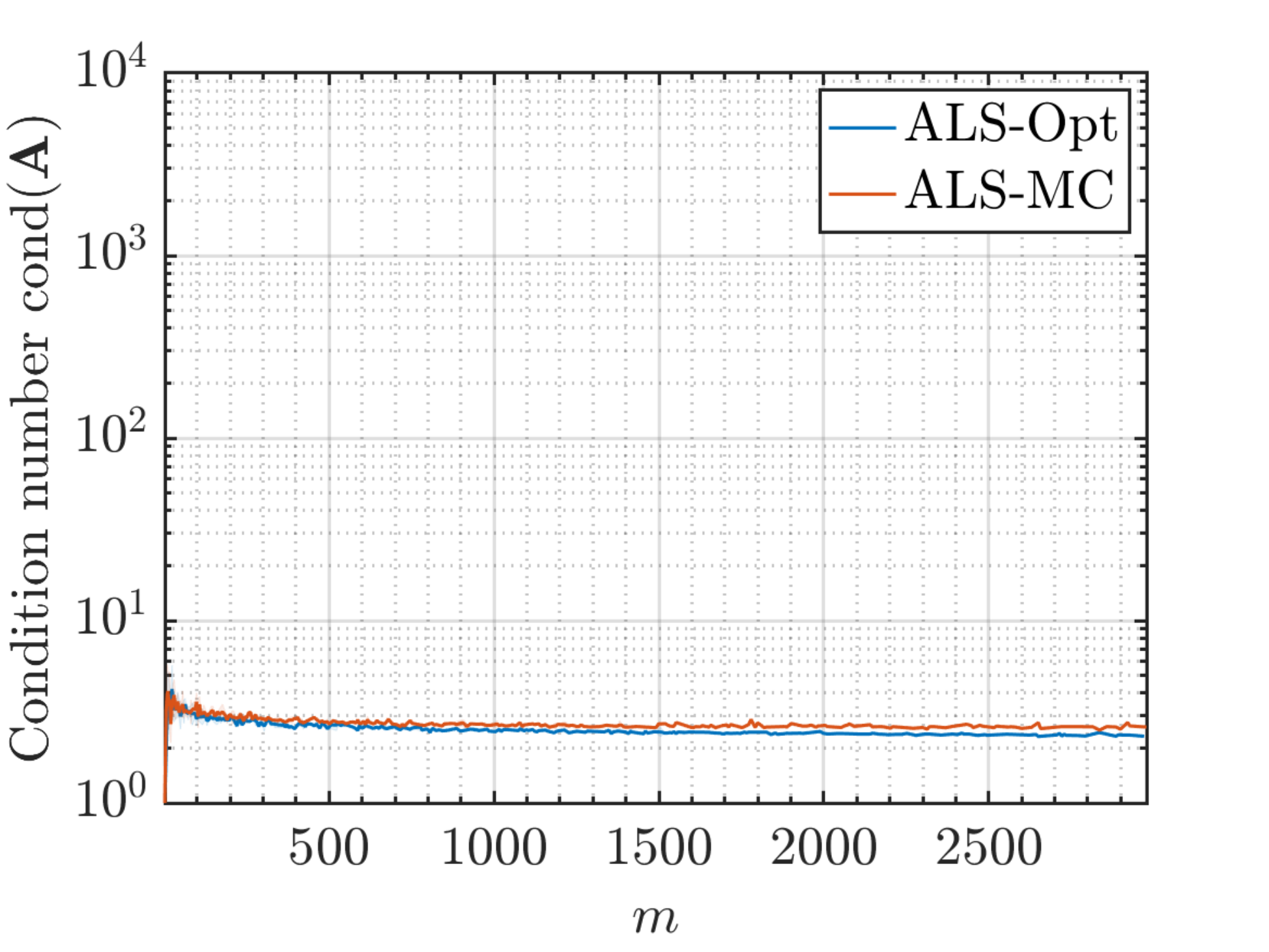}
\\[\errplottextsp]
$d = 8$ & $d = 16$ & $d = 32$
\end{tabular}
\end{small}
\end{center}
\caption{ALS approximation of the function $f = f_1$ using Legendre polynomials and either near-optimal sampling (\S \ref{ss:opt-samp-scheme}) or MC sampling (\S \ref{ss:MC-samp-scheme}). This figure shows the relative $L^2_{\varrho}(\cU)$-norm error (rows 1 and 3), computed over a grid of 100,000 MC points, and condition number $\kappa(\bm{A})$ (rows 2 and 4) versus the number of samples $m$, which is given by \eqref{m-scaling-LS}. {The dashed line shows the best $n$-term approximation error. Here, for each $m$, $n$ is chosen as in \eqref{m-scaling-LS} so that $m = \max \{ n + 1 , \lceil n \cdot \log(n) \rceil \}$.}} 
\label{fig:fig2}
\end{figure}

\begin{figure}[t!]
\begin{center}
\begin{small}
 \begin{tabular}{@{\hspace{0pt}}c@{\hspace{\errplotsp}}c@{\hspace{\errplotsp}}c@{\hspace{0pt}}}
\includegraphics[width = \errplotimg]{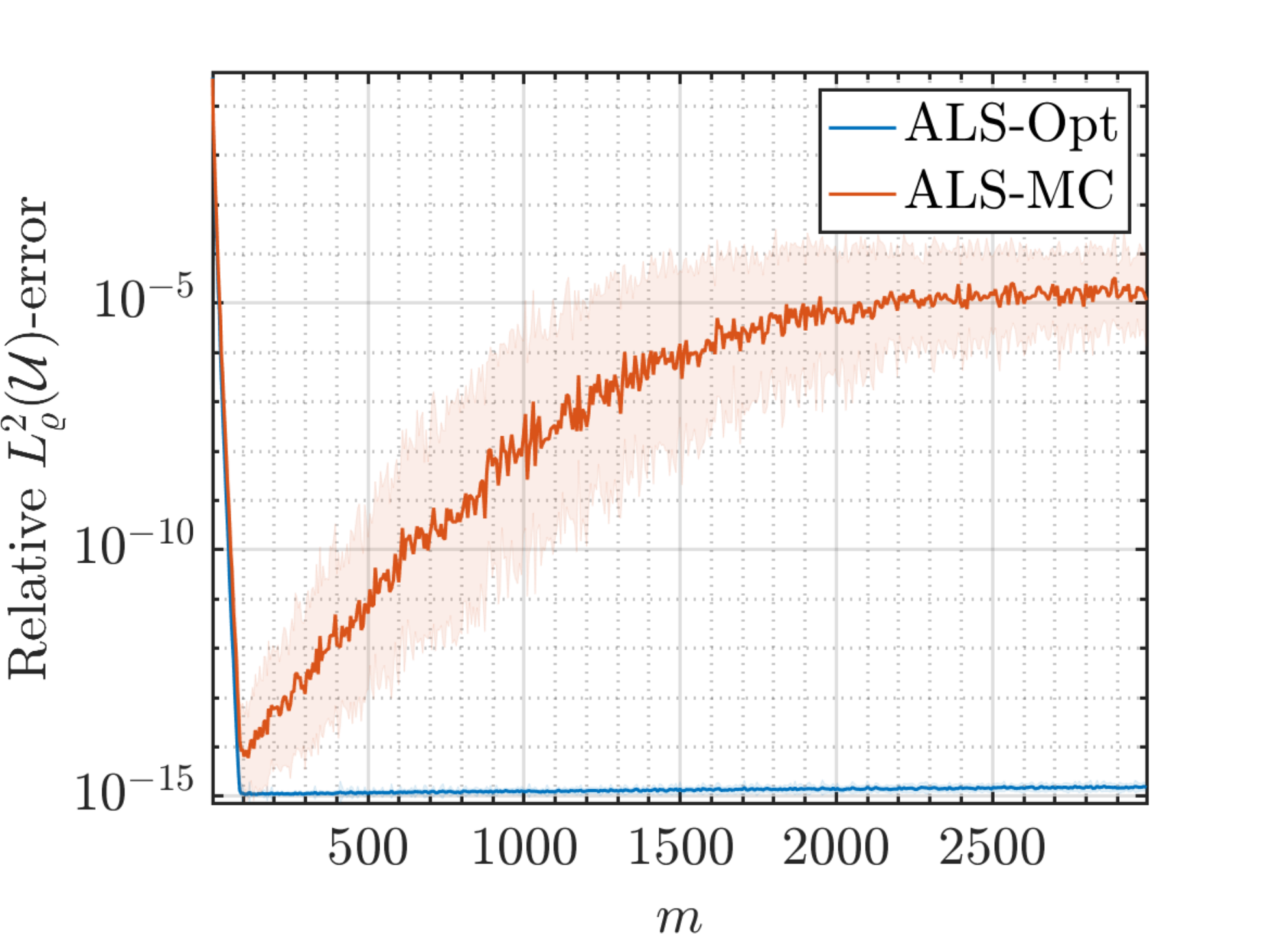}
&
\includegraphics[width = \errplotimg]{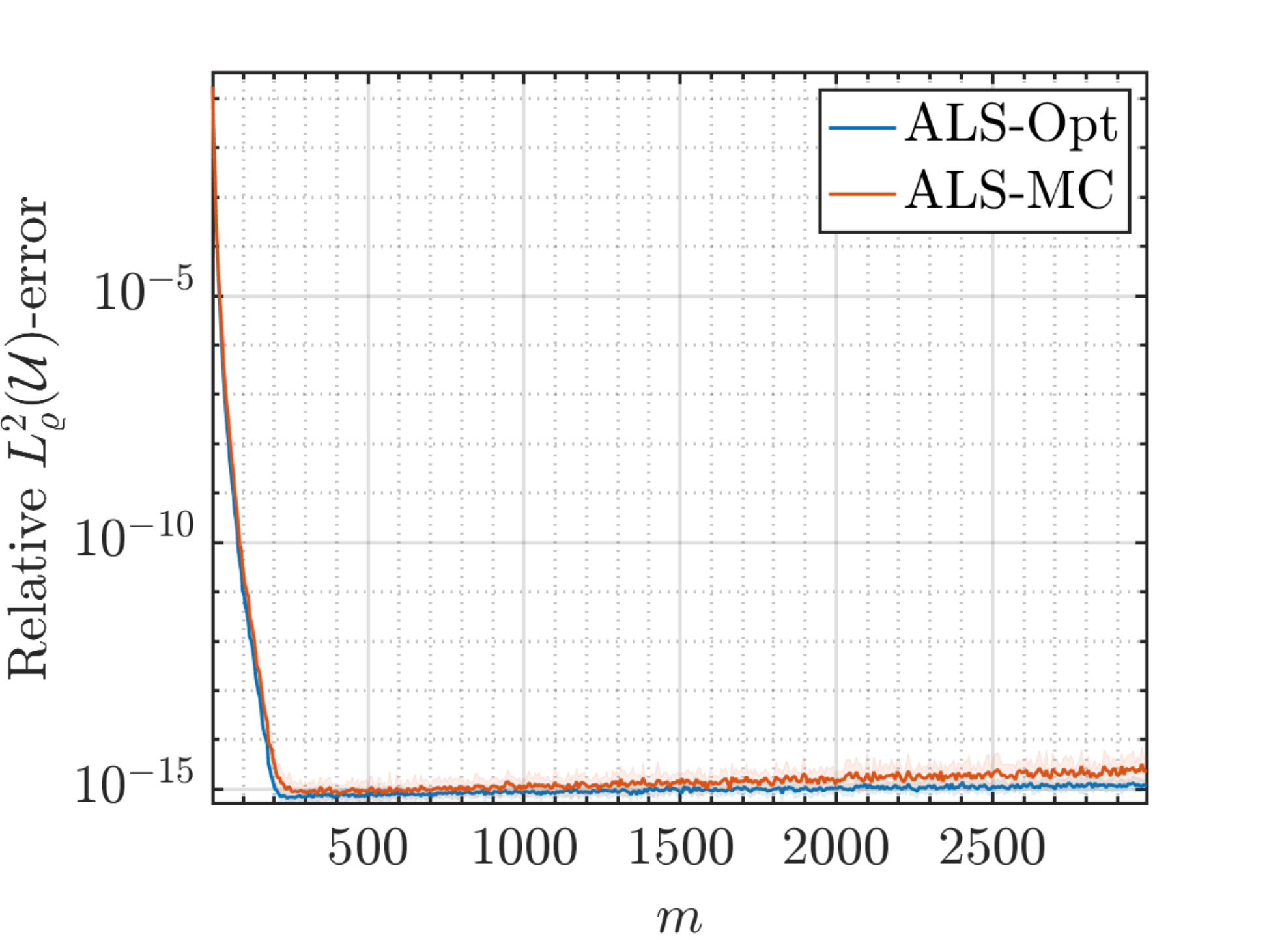}
&
\includegraphics[width = \errplotimg]{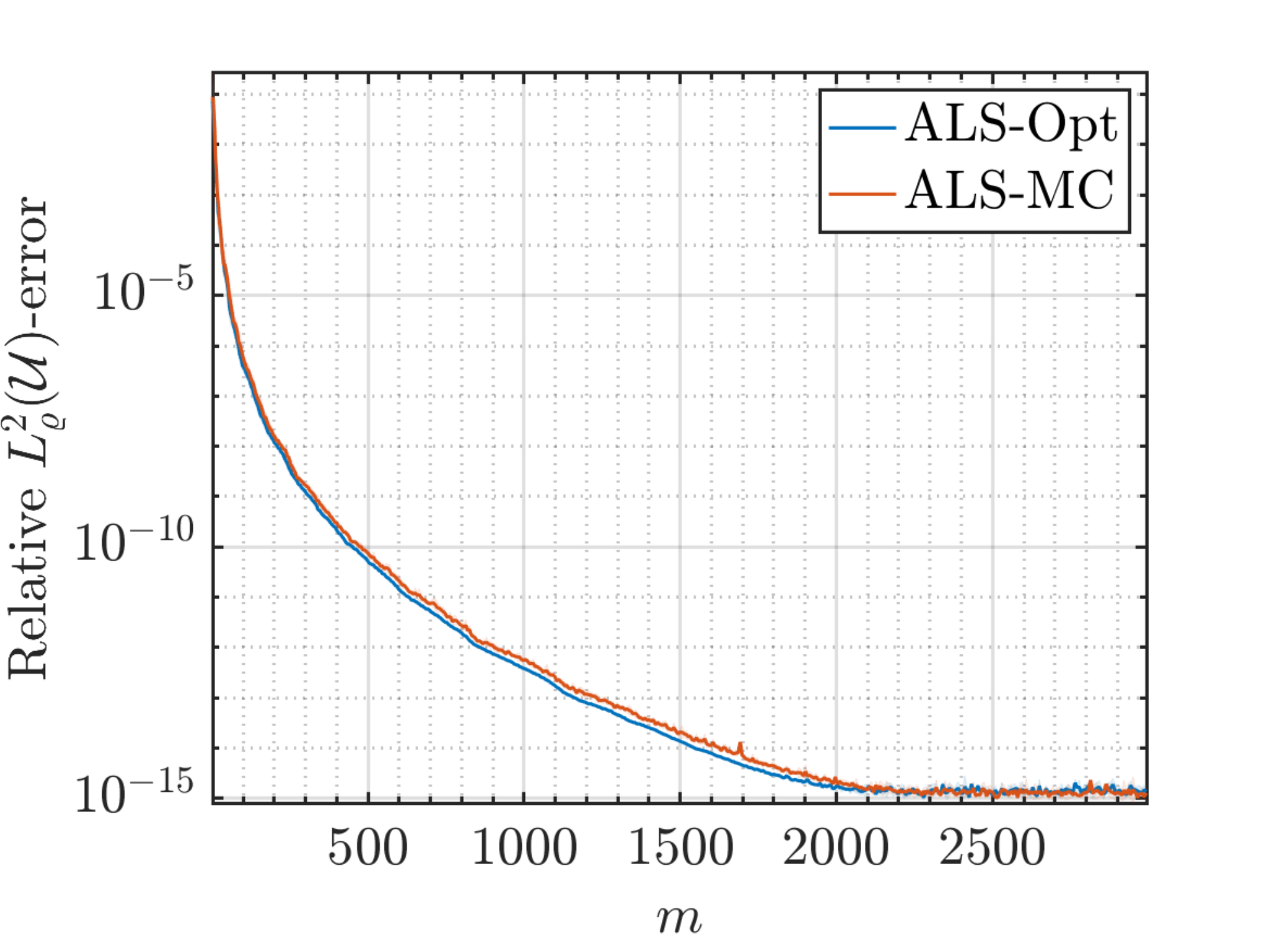}
\\[\errplotgraphsp]
\includegraphics[width = \errplotimg]{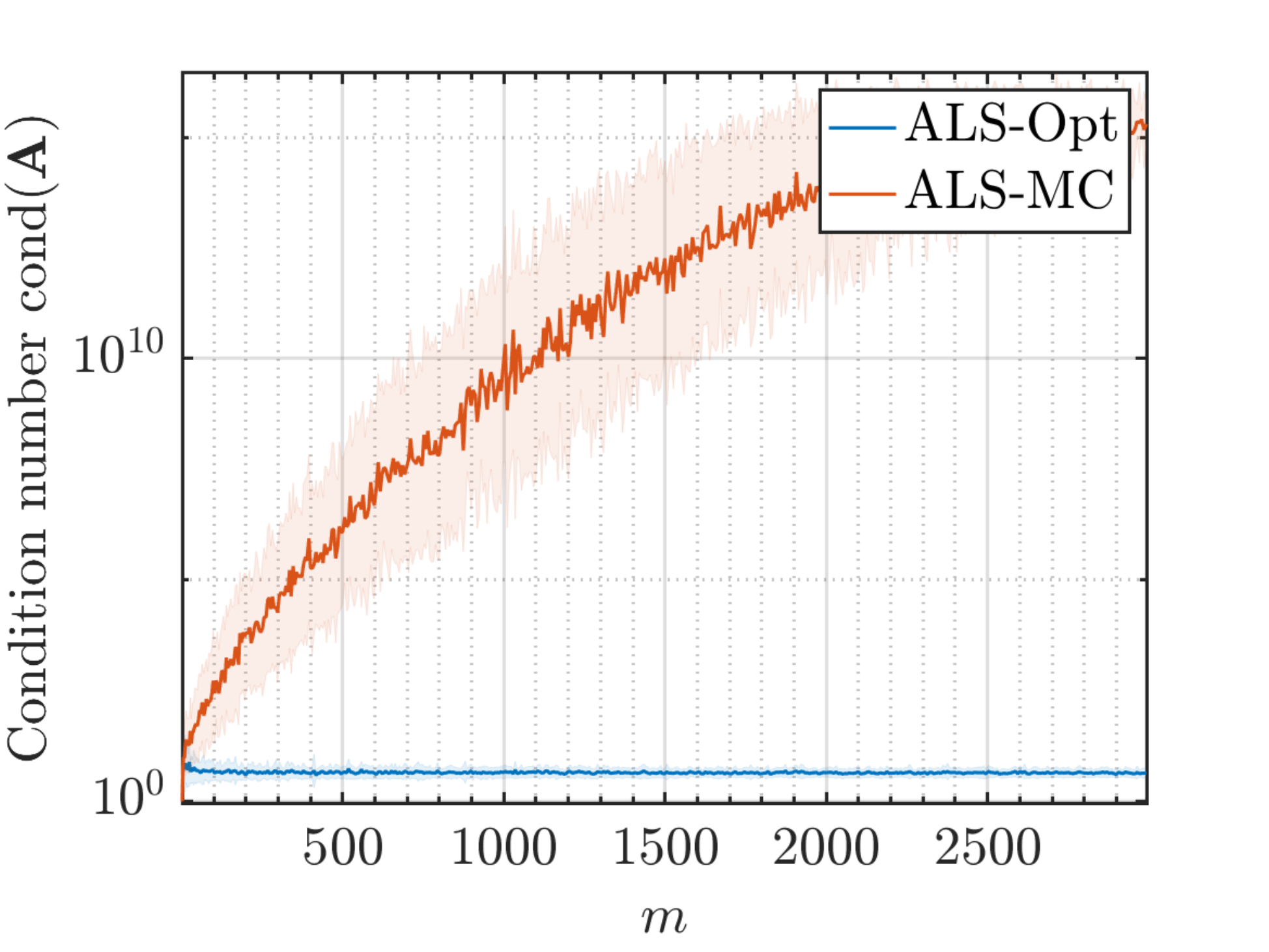}
&
\includegraphics[width = \errplotimg]{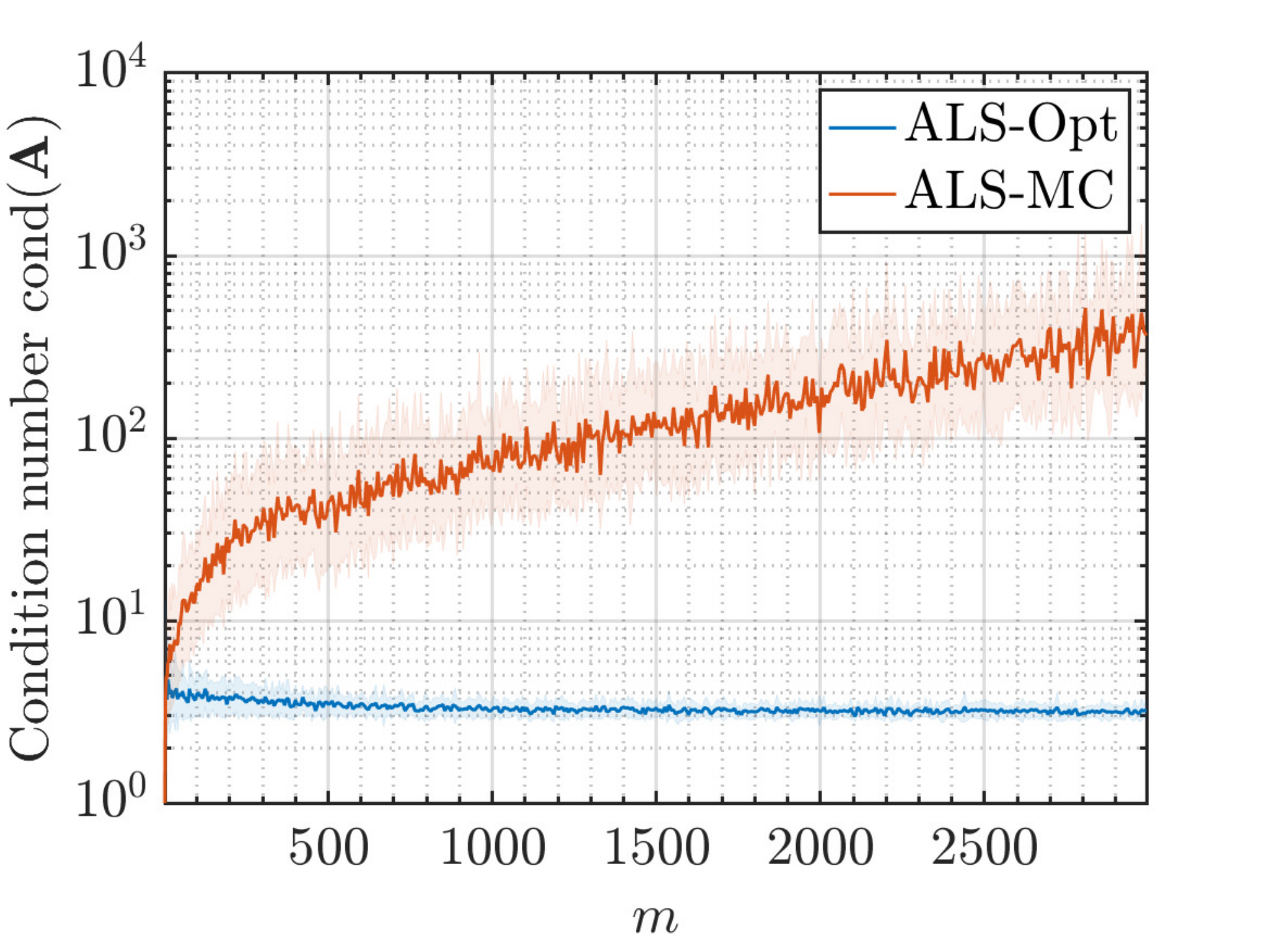}
&
\includegraphics[width = \errplotimg]{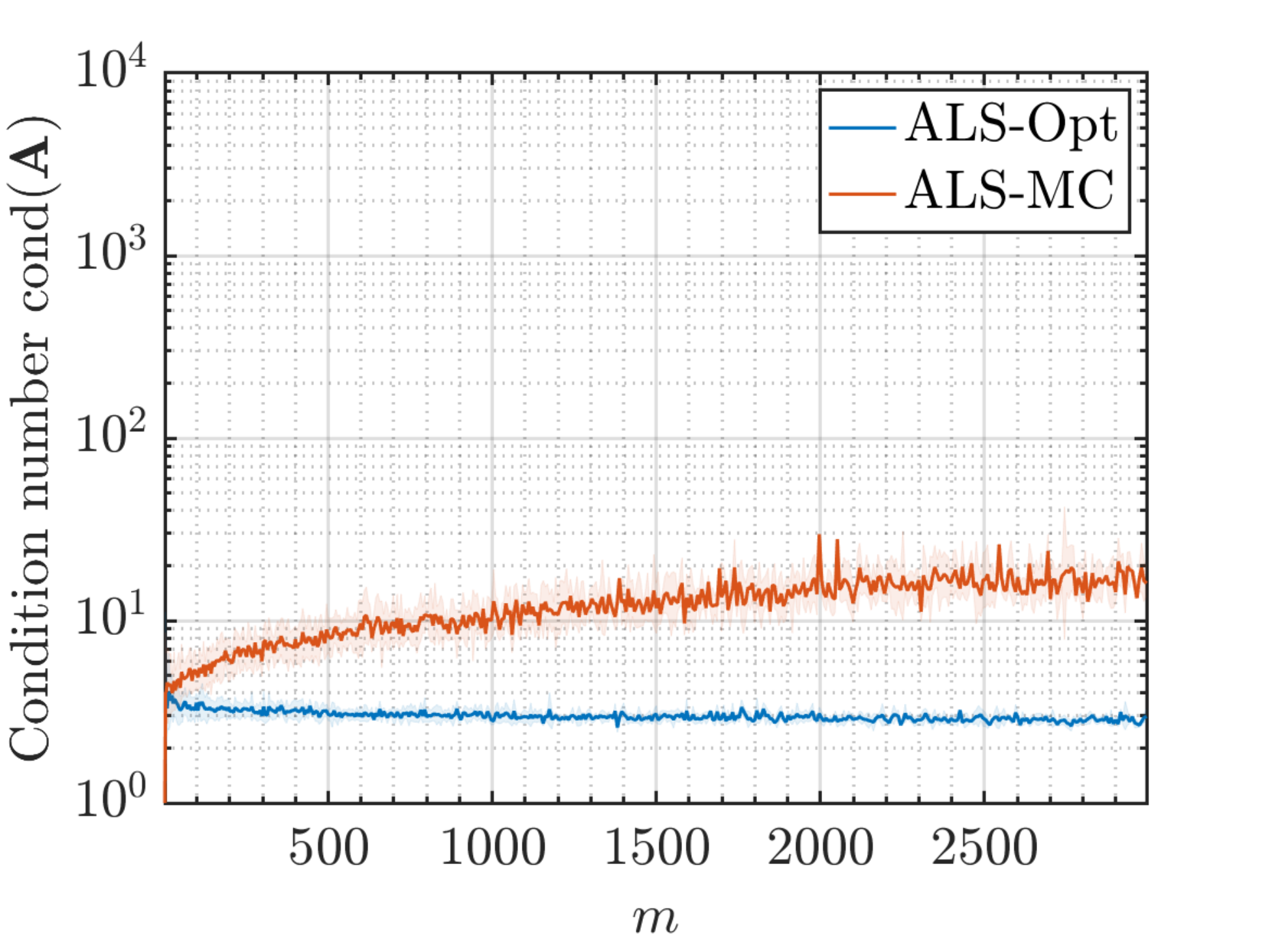}
\\[\errplottextsp]
$d = 1$ & $d = 2$ & $d = 4$
\\[\errplottextsp]
\includegraphics[width = \errplotimg]{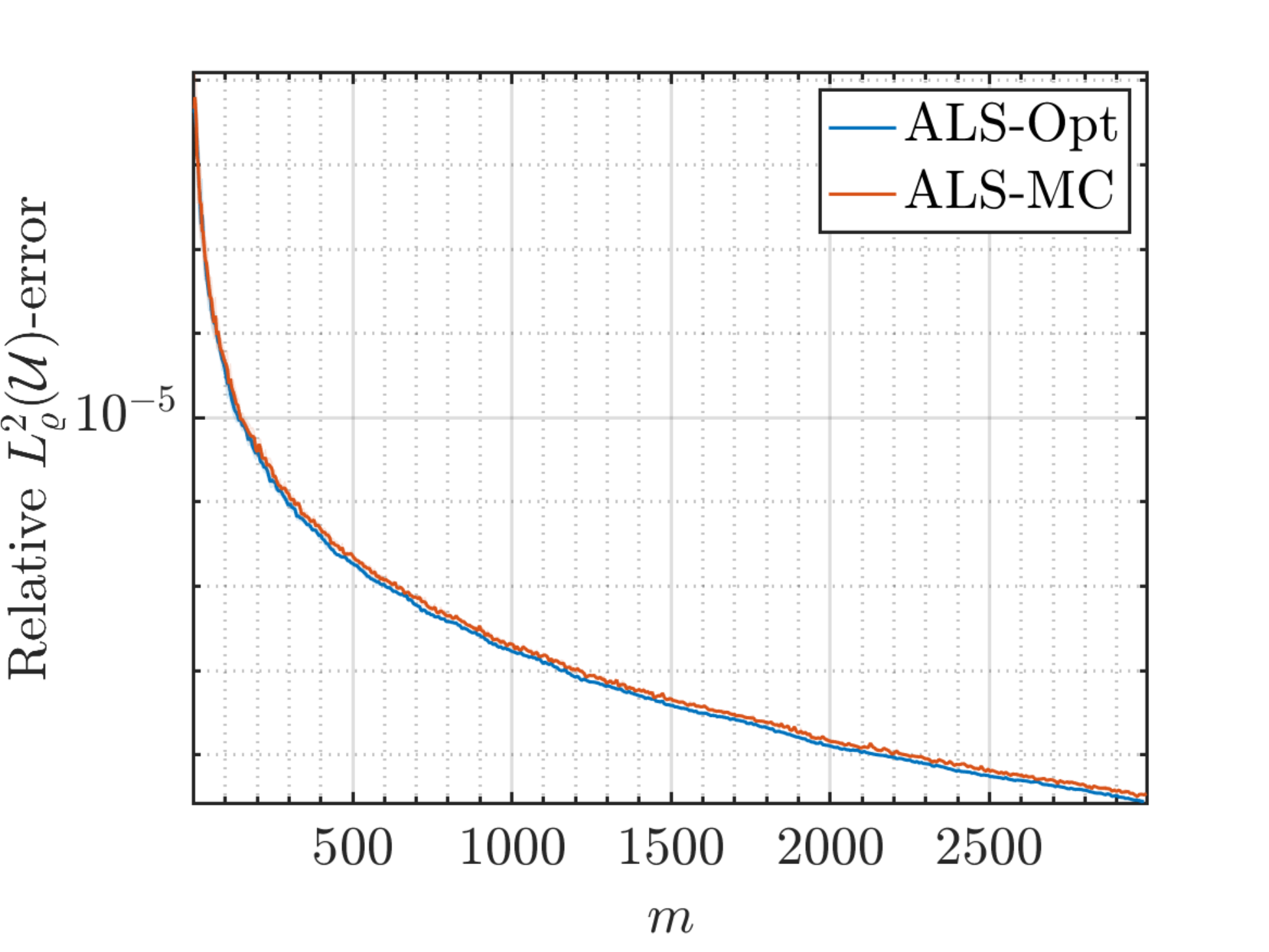}
&
\includegraphics[width = \errplotimg]{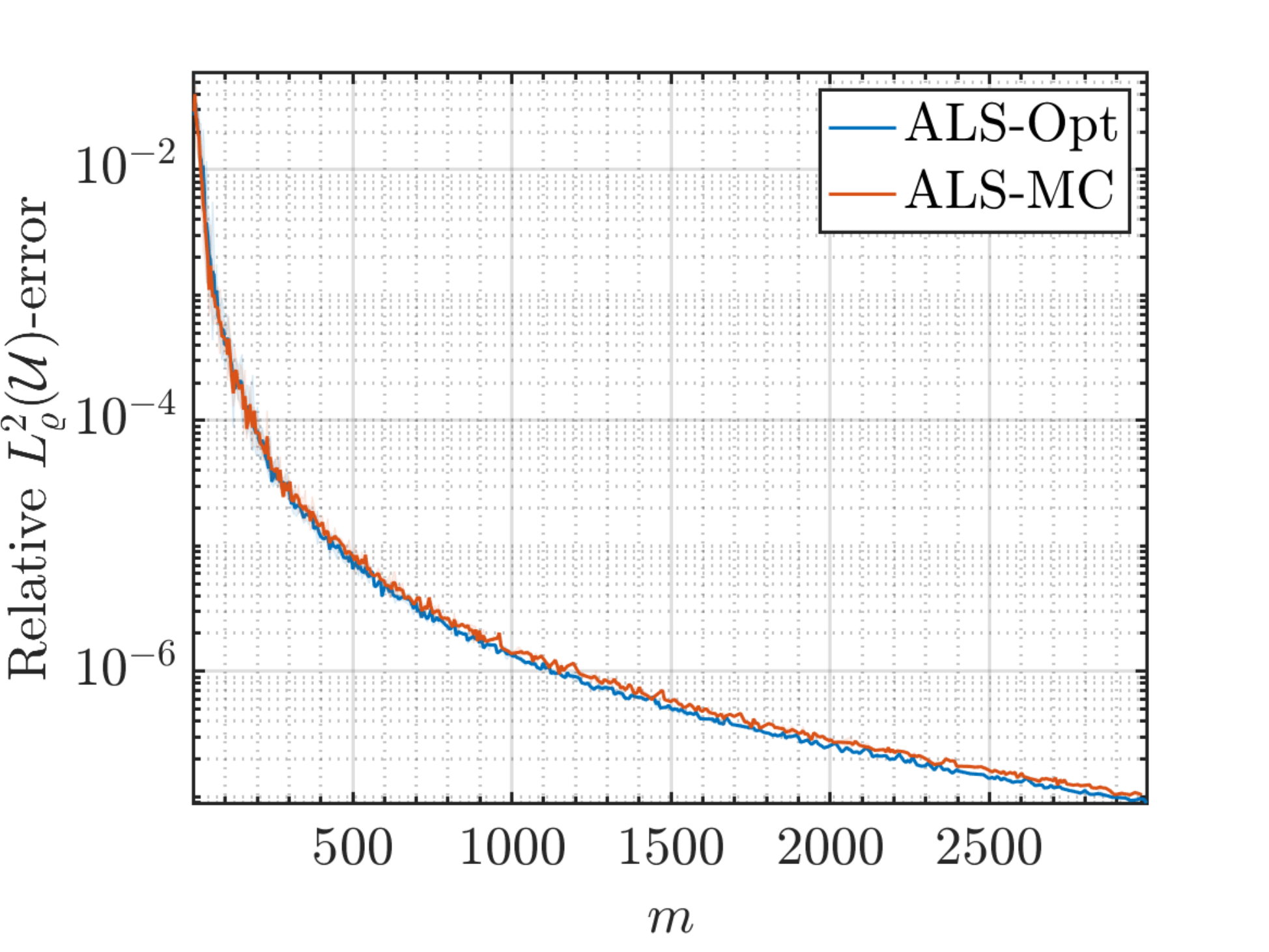}
&
\includegraphics[width = \errplotimg]{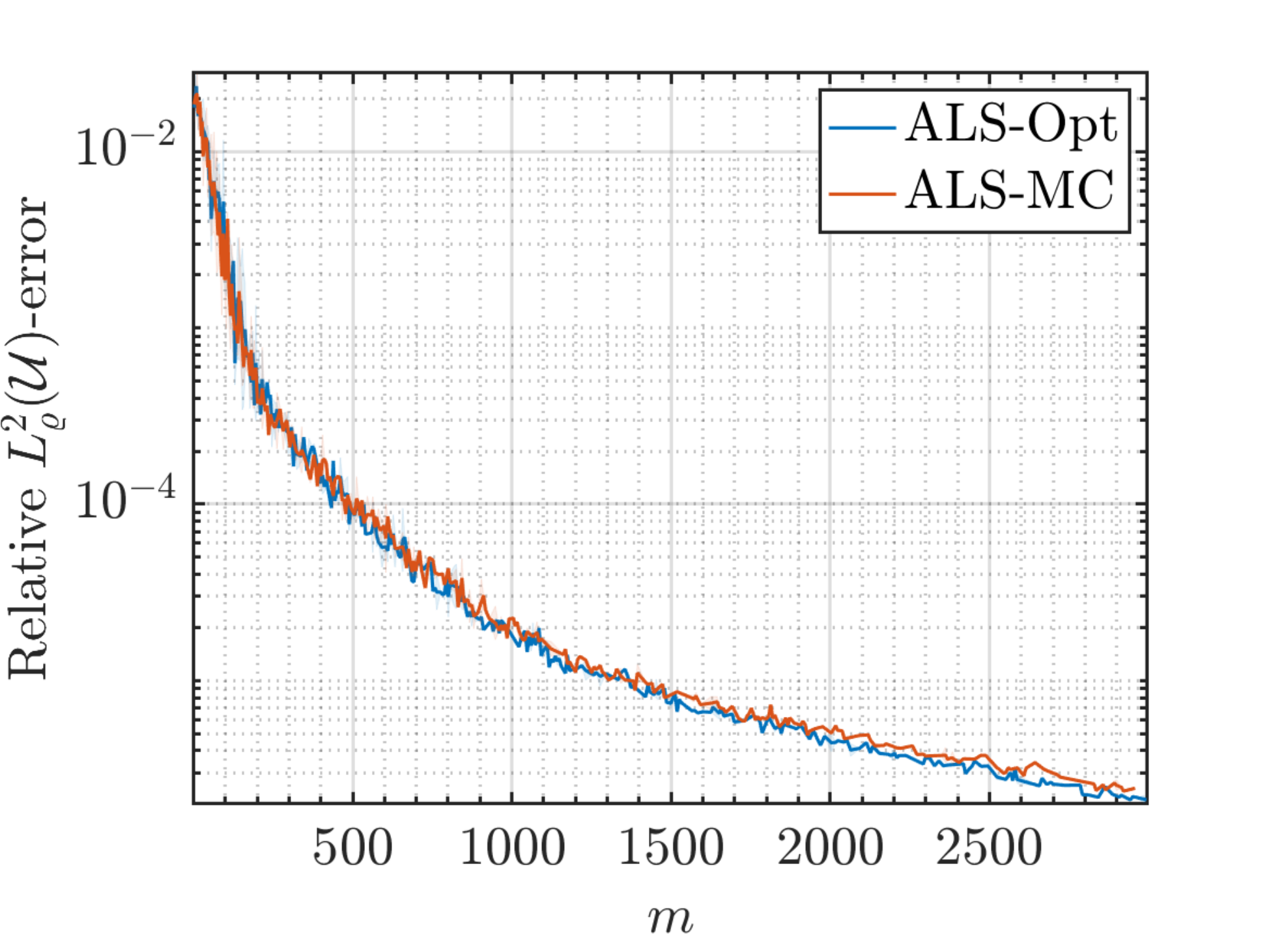}
\\[\errplotgraphsp]
\includegraphics[width = \errplotimg]{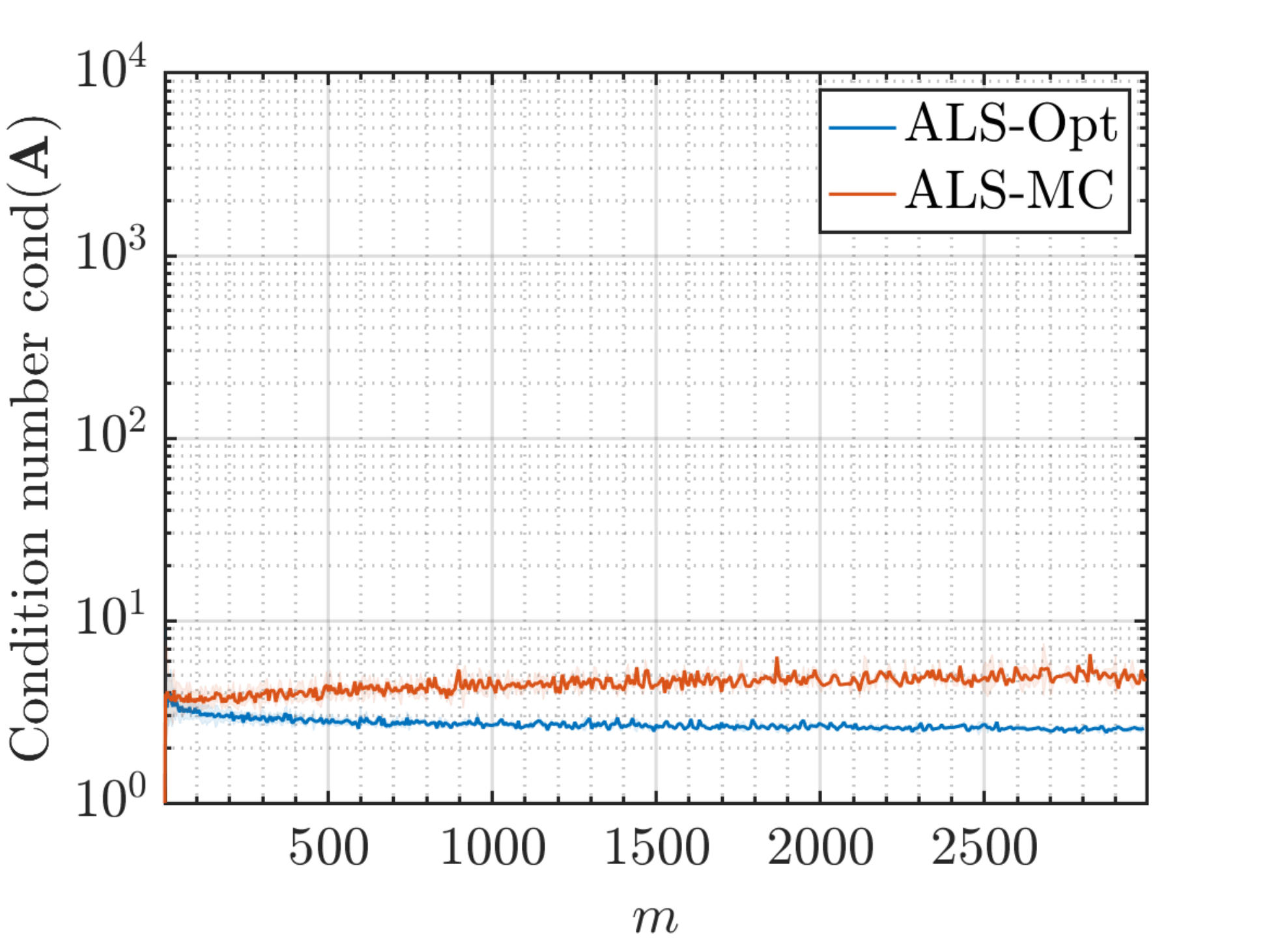}
&
\includegraphics[width = \errplotimg]{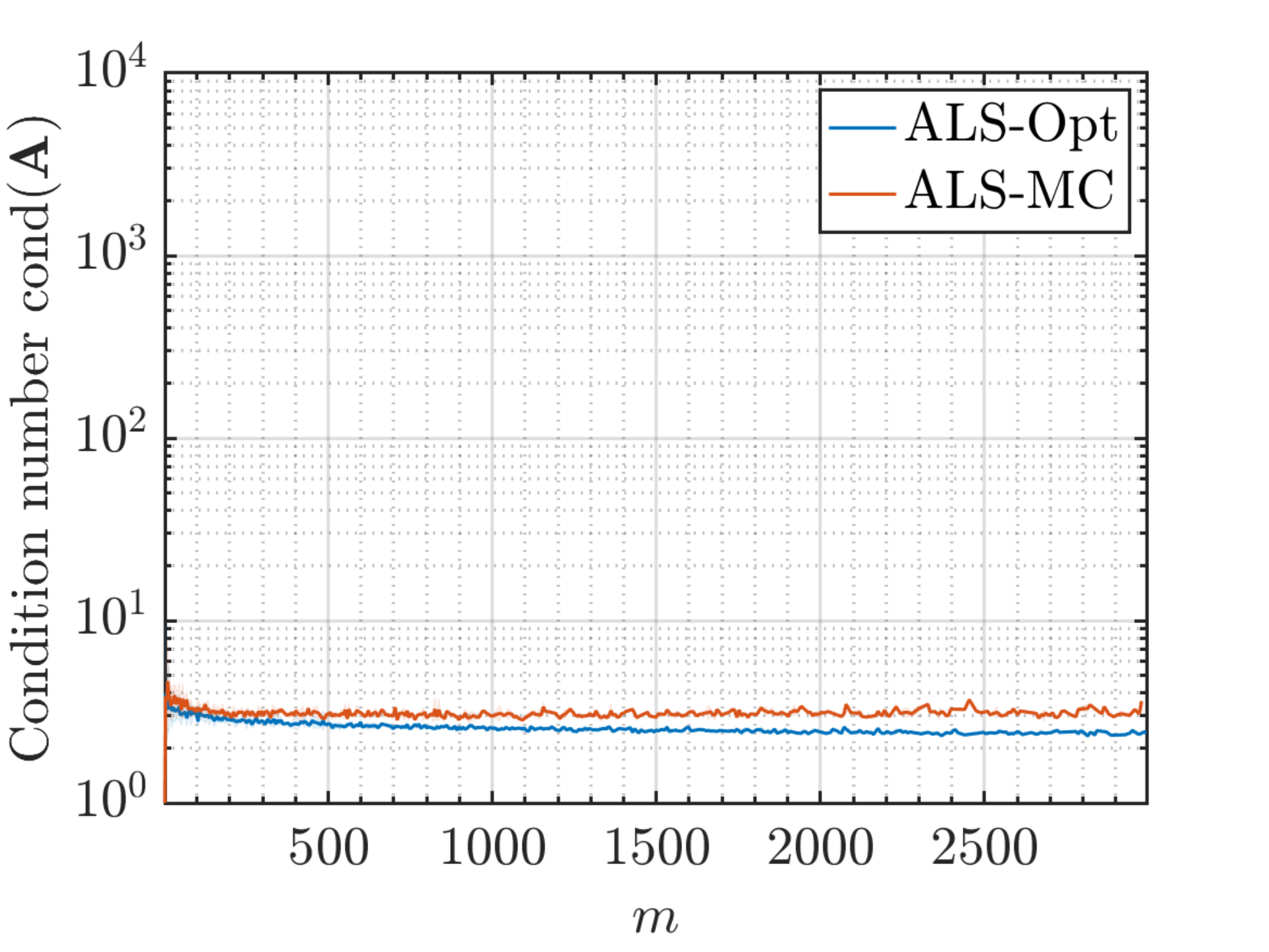}
&
\includegraphics[width = \errplotimg]{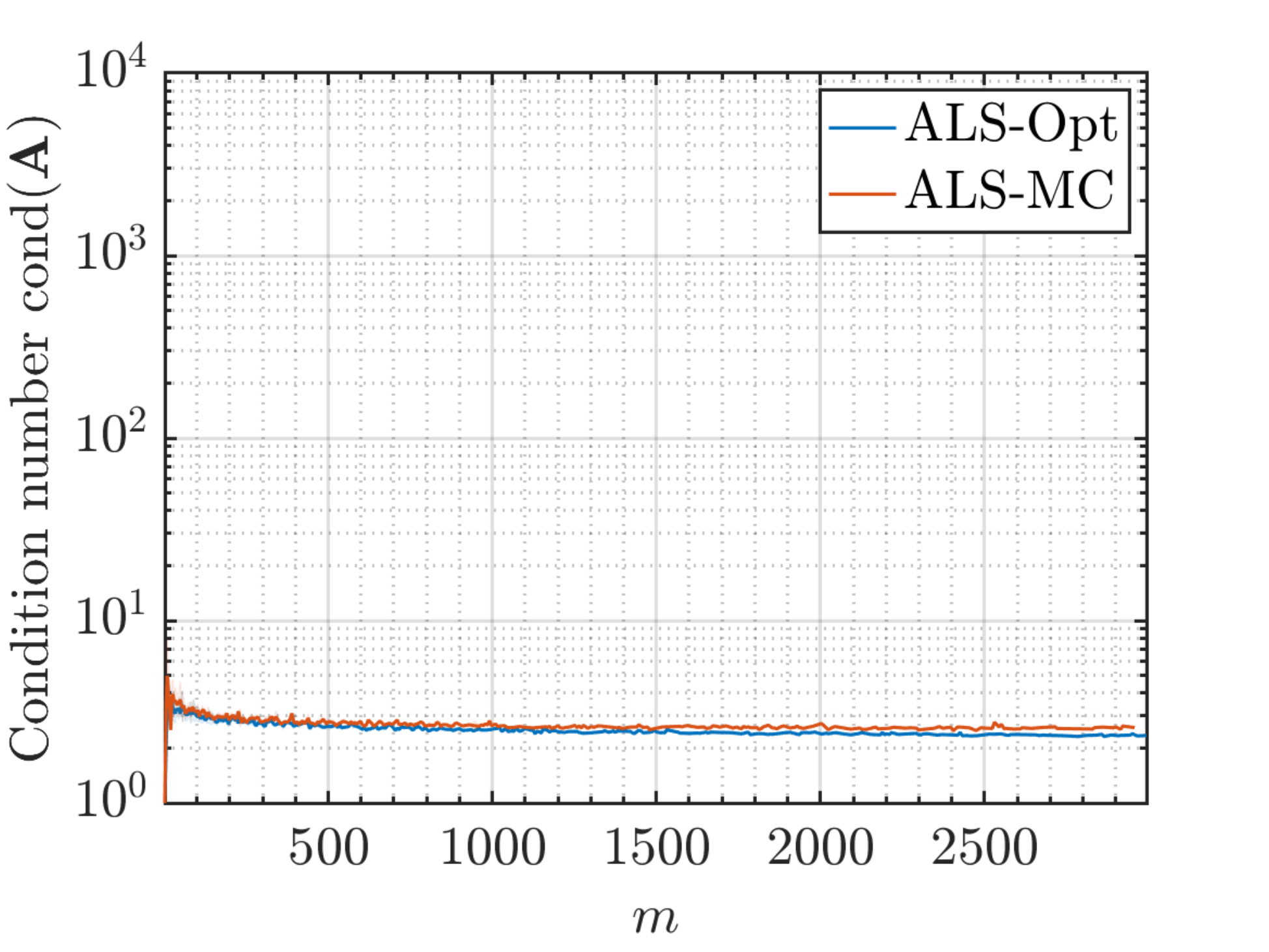}
\\[\errplottextsp]
$d = 8$ & $d = 16$ & $d = 32$
\end{tabular}
\end{small}
\end{center}
\caption{The same as Fig.\ \ref{fig:fig2} but with $f = f_2$.} 
\label{fig:fig3}
\end{figure}

\begin{figure}[t!]
\begin{center}
\begin{small}
 \begin{tabular}{@{\hspace{0pt}}c@{\hspace{\errplotsp}}c@{\hspace{\errplotsp}}c@{\hspace{0pt}}}
\includegraphics[width = \errplotimg]{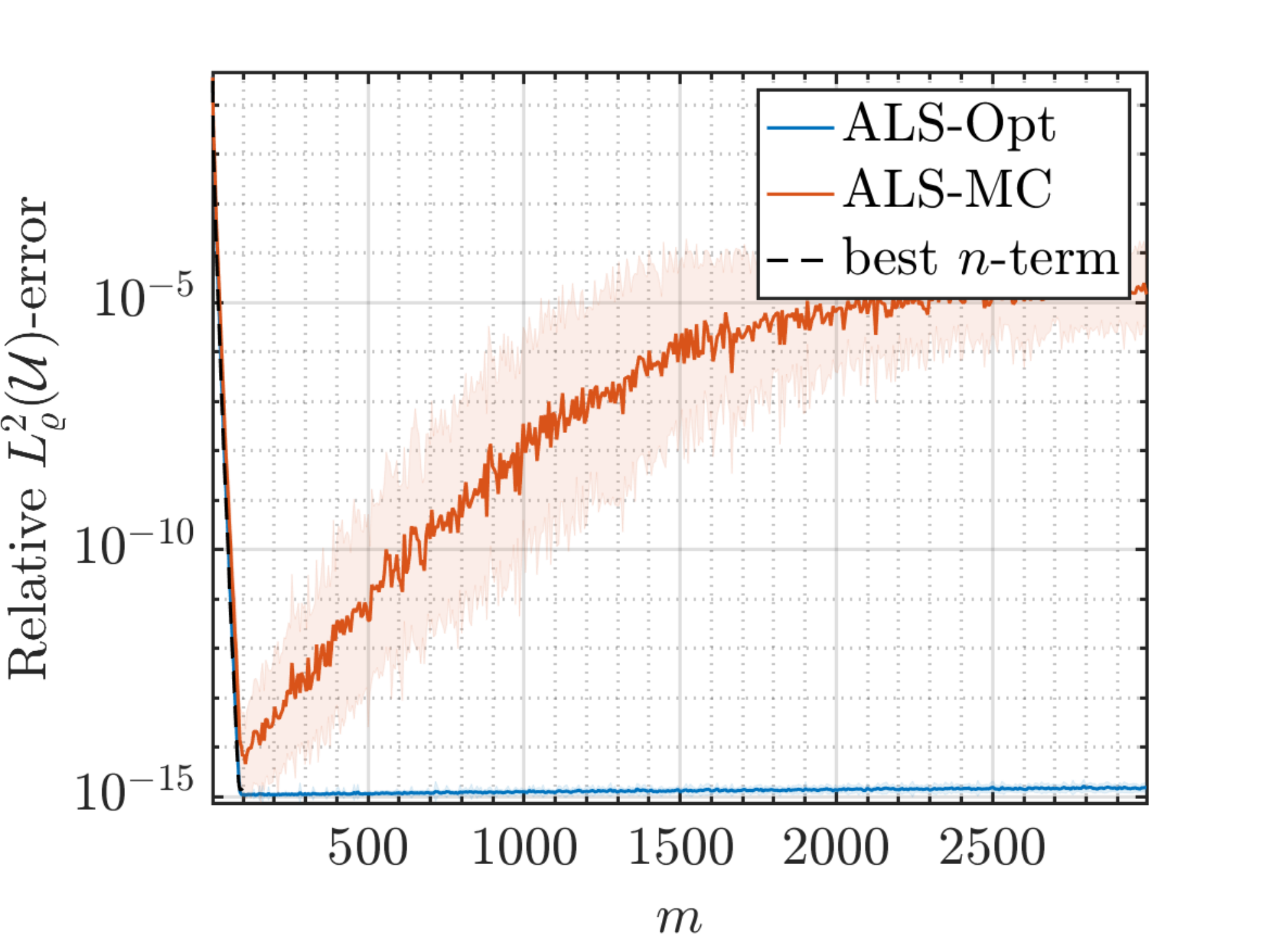}
&
\includegraphics[width = \errplotimg]{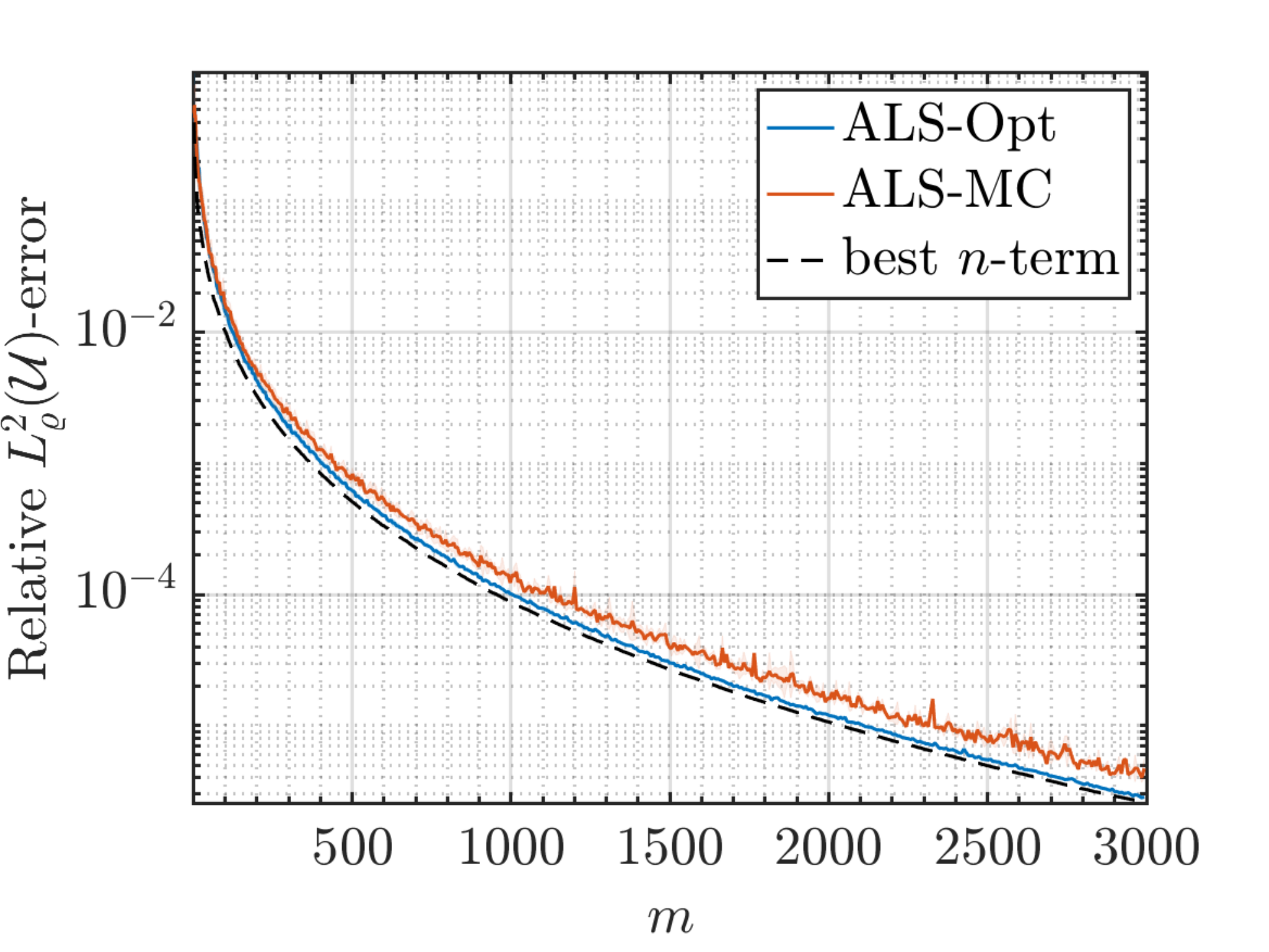}
&
\includegraphics[width = \errplotimg]{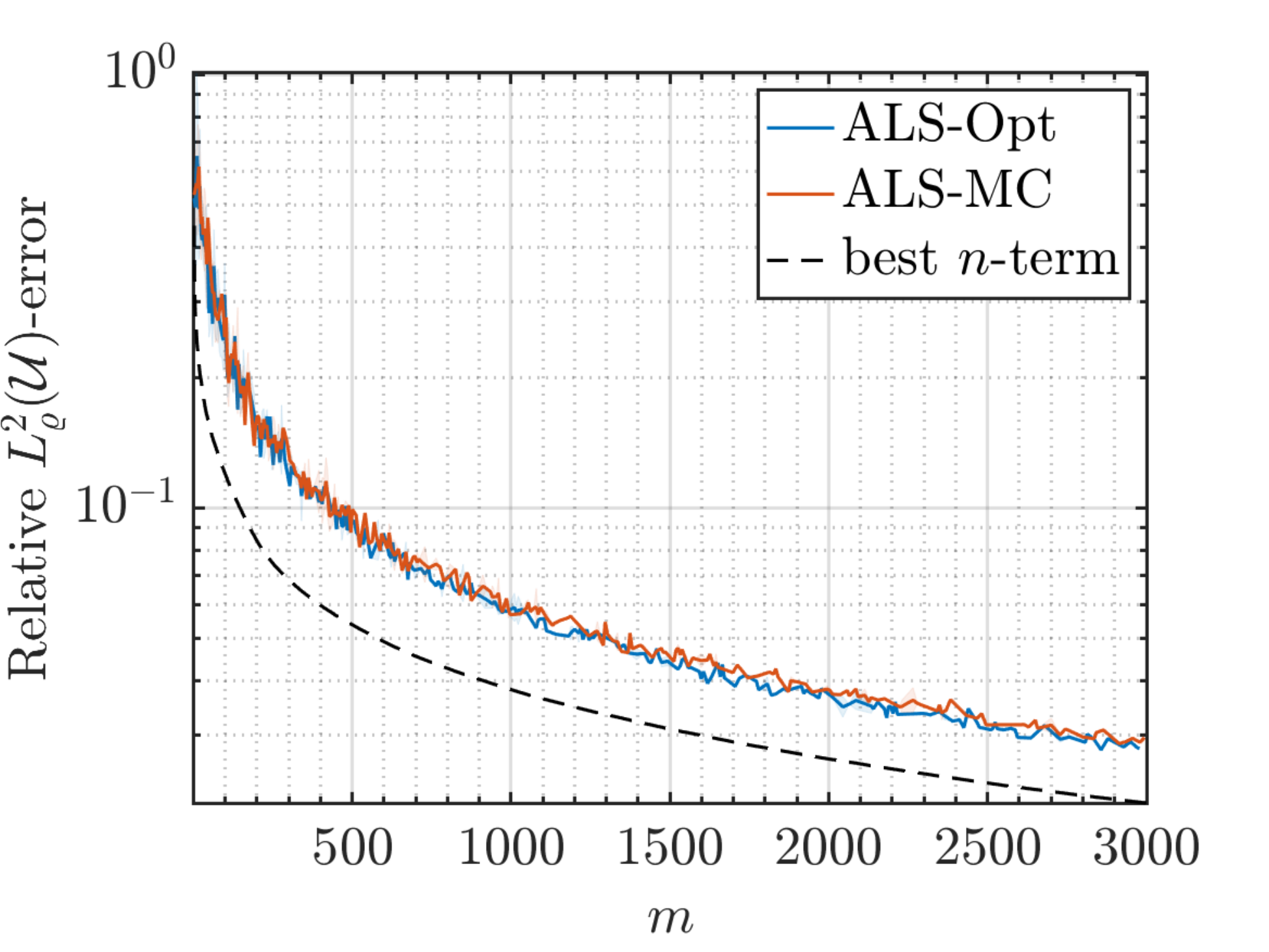}
\\[\errplotgraphsp]
\includegraphics[width = \errplotimg]{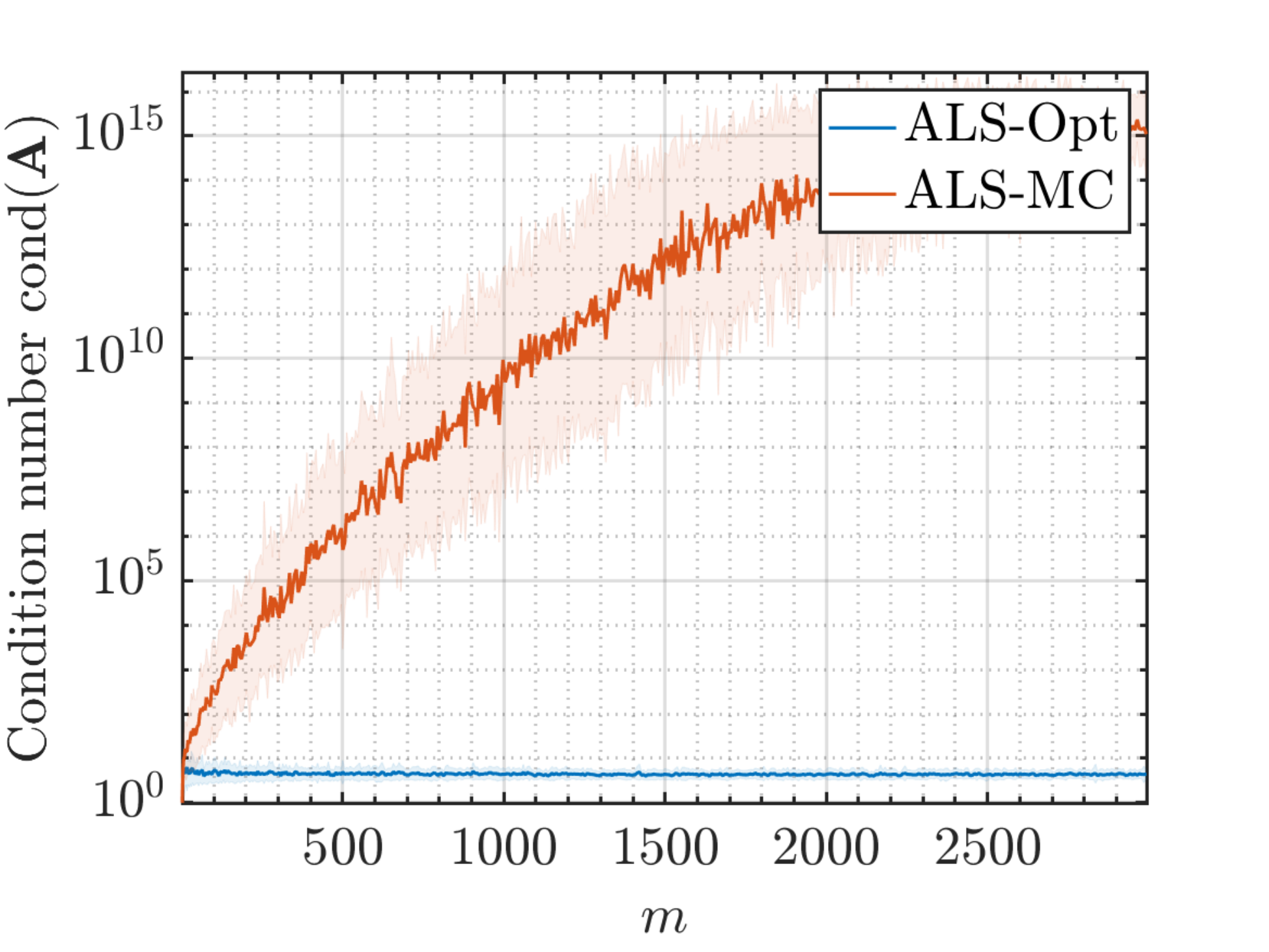}
&
\includegraphics[width = \errplotimg]{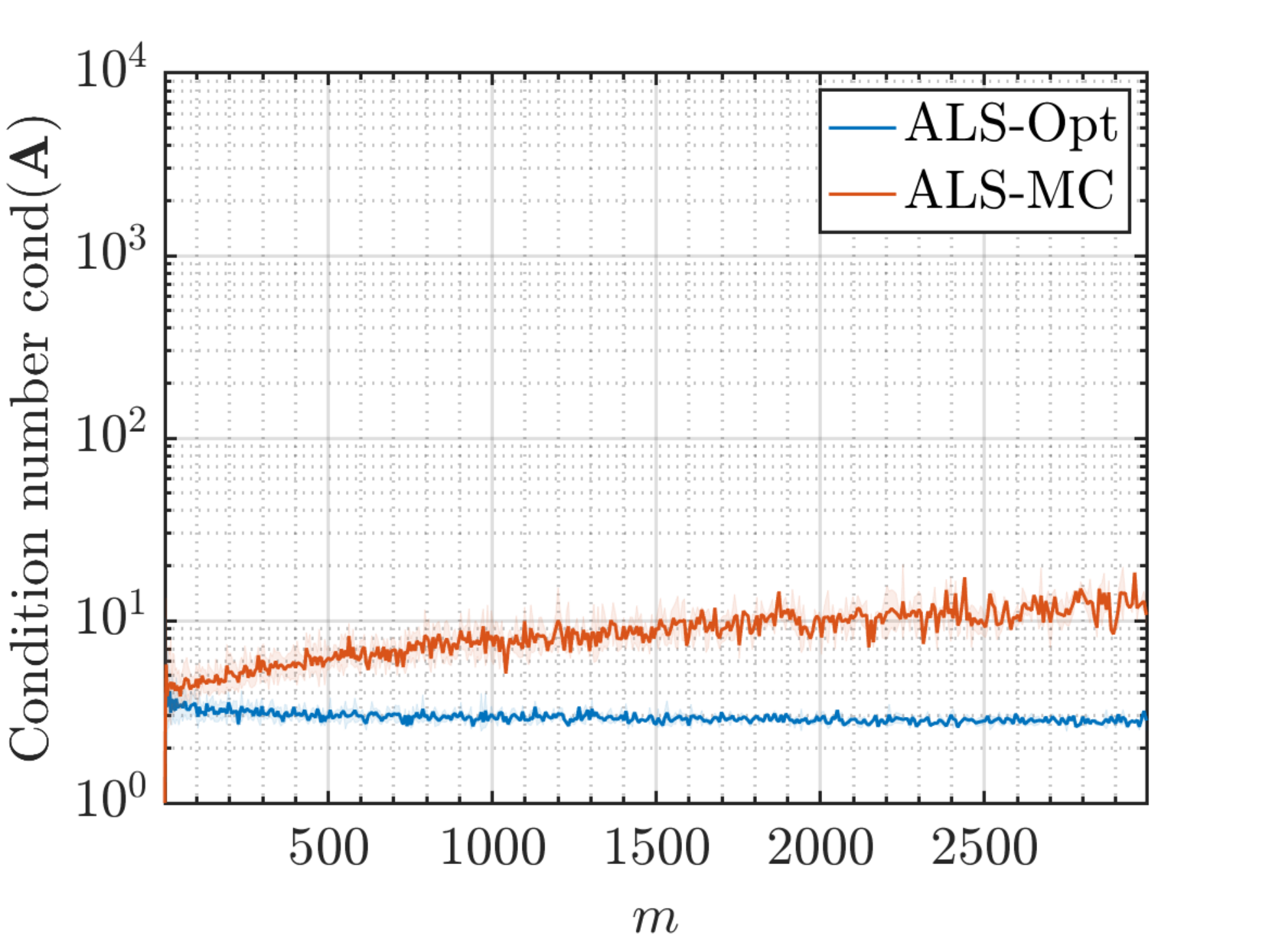}
&
\includegraphics[width = \errplotimg]{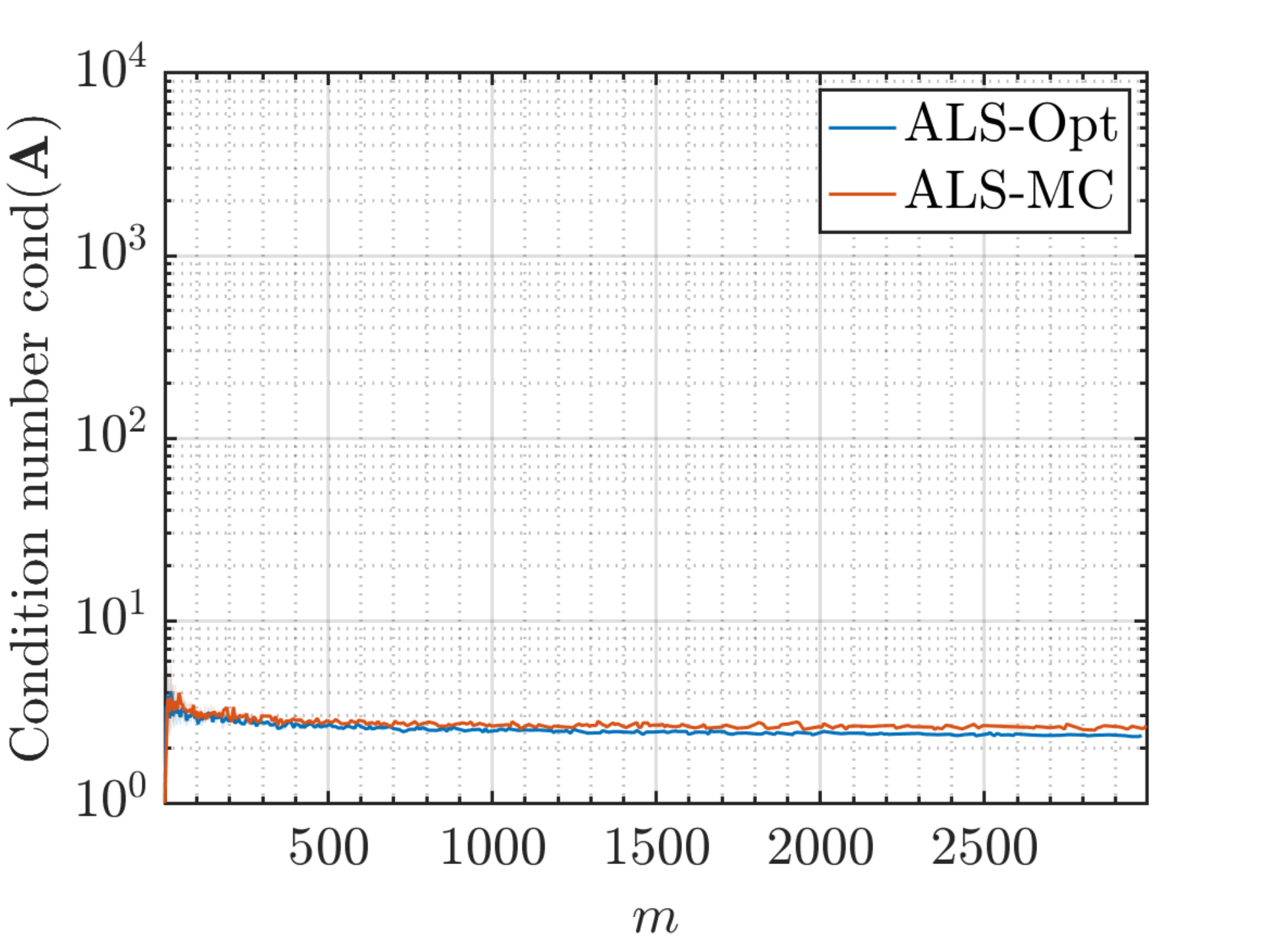}
\\[\errplotgraphsp]
\includegraphics[width = \errplotimg]{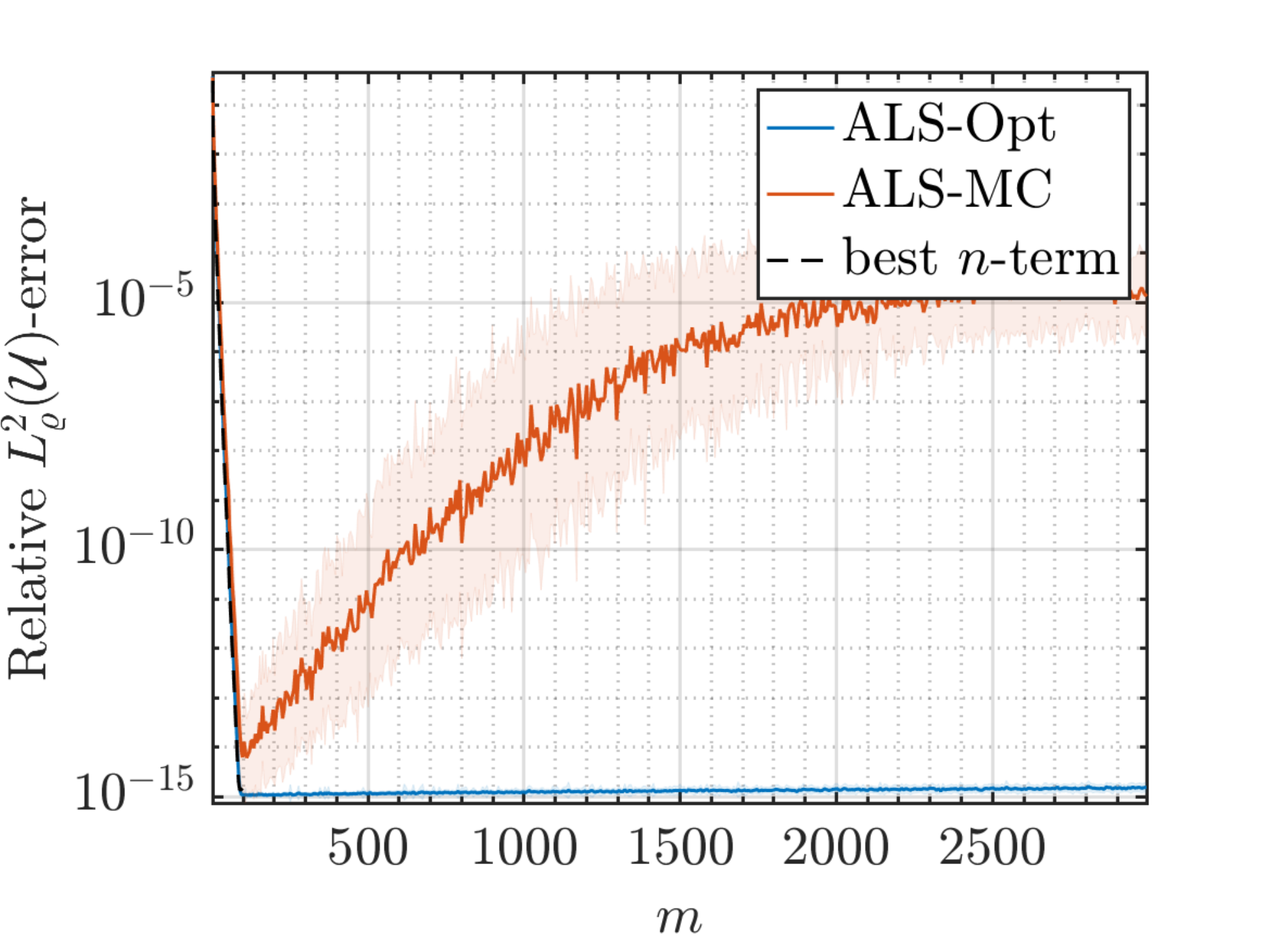}
&
\includegraphics[width = \errplotimg]{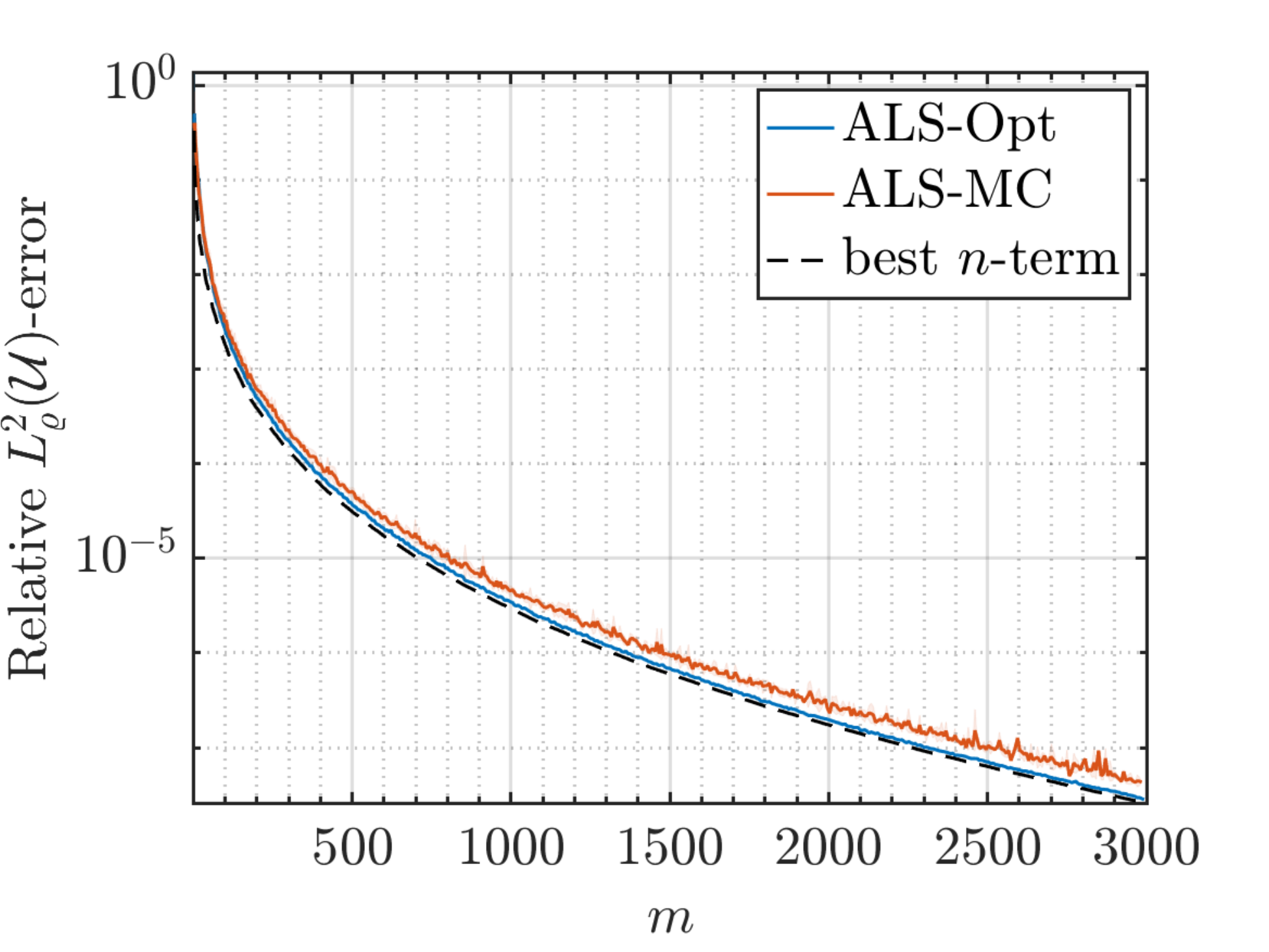}
&
\includegraphics[width = \errplotimg]{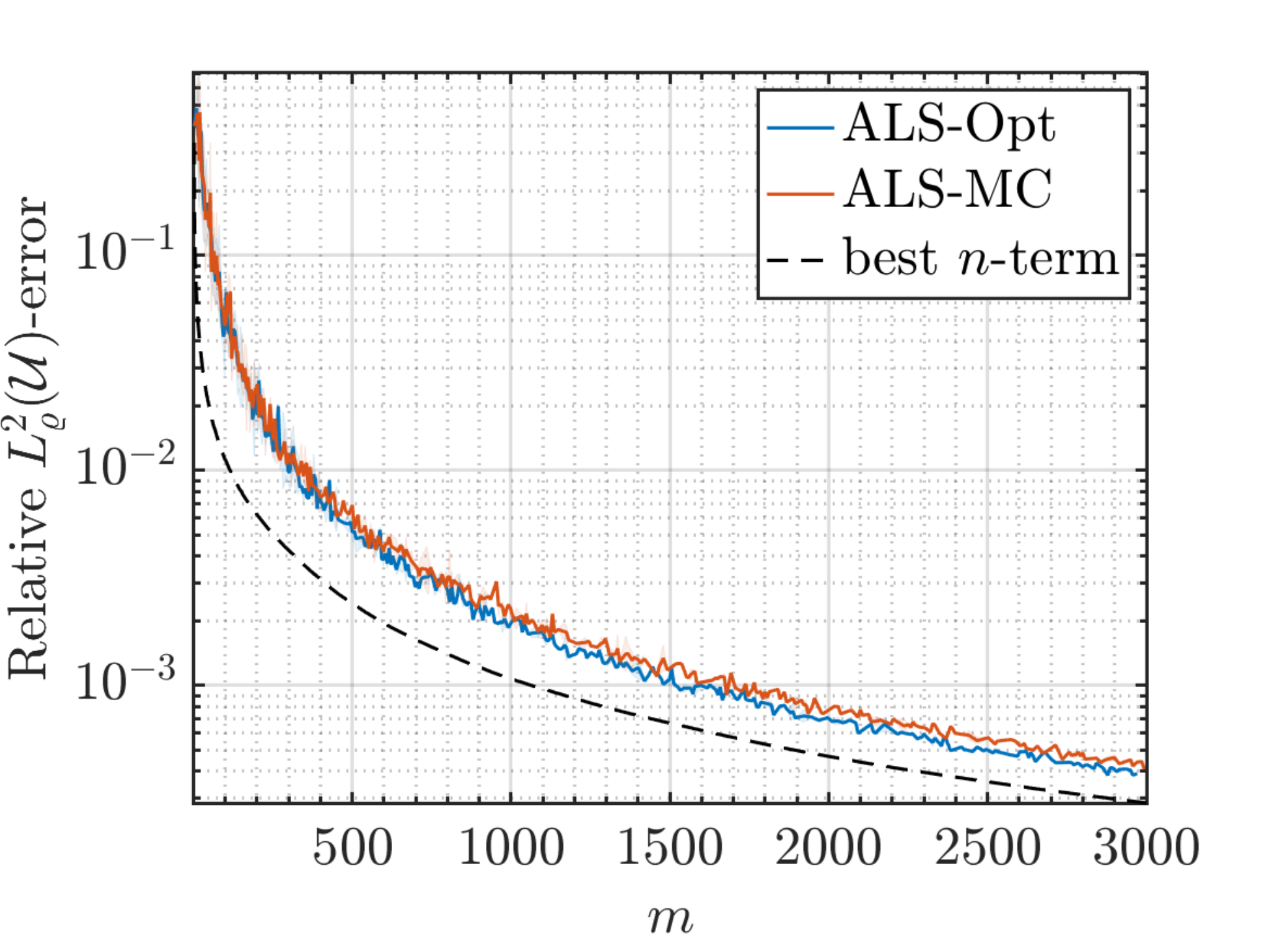}
\\[\errplotgraphsp]
\includegraphics[width = \errplotimg]{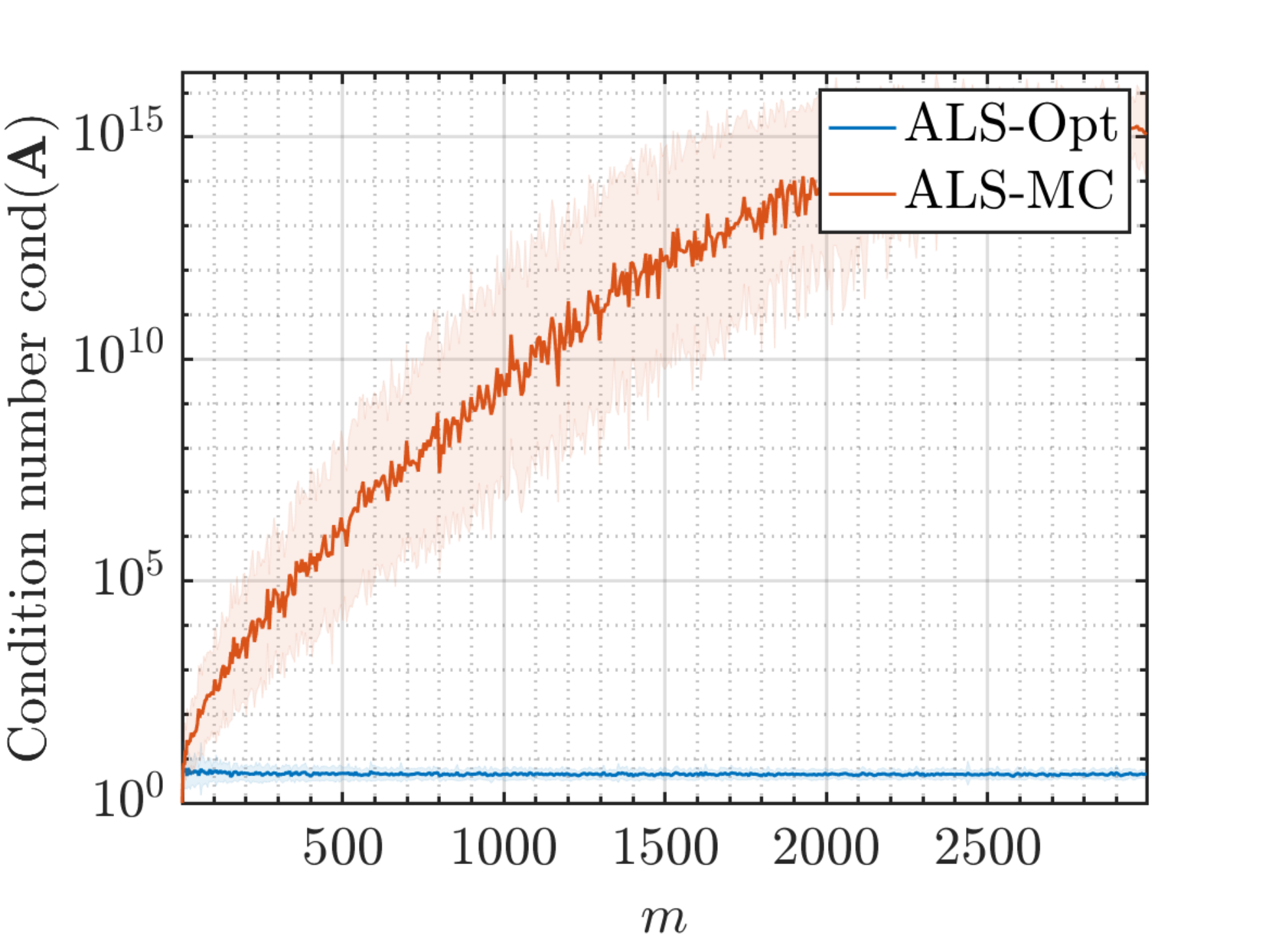}
&
\includegraphics[width = \errplotimg]{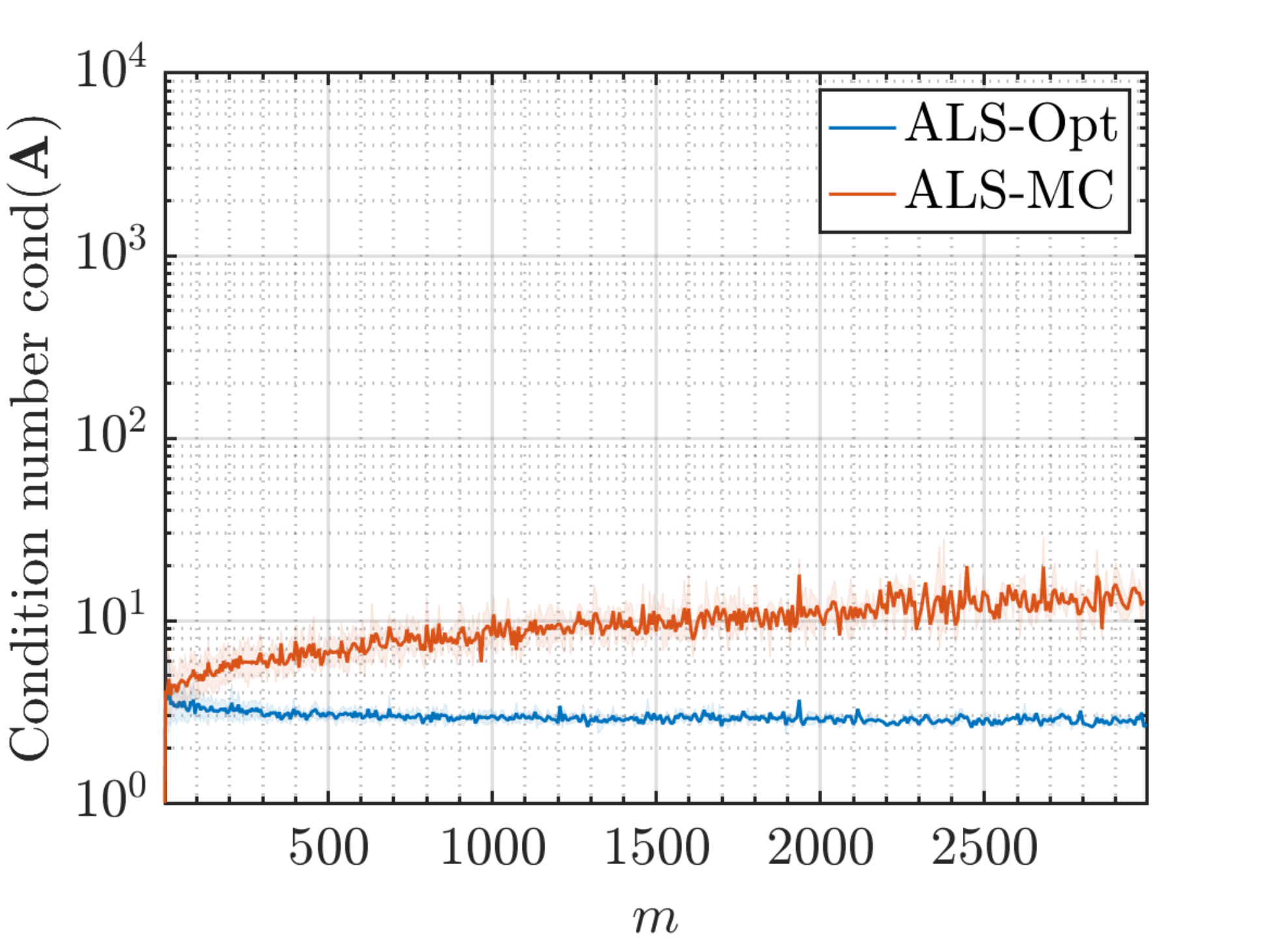}
&
\includegraphics[width = \errplotimg]{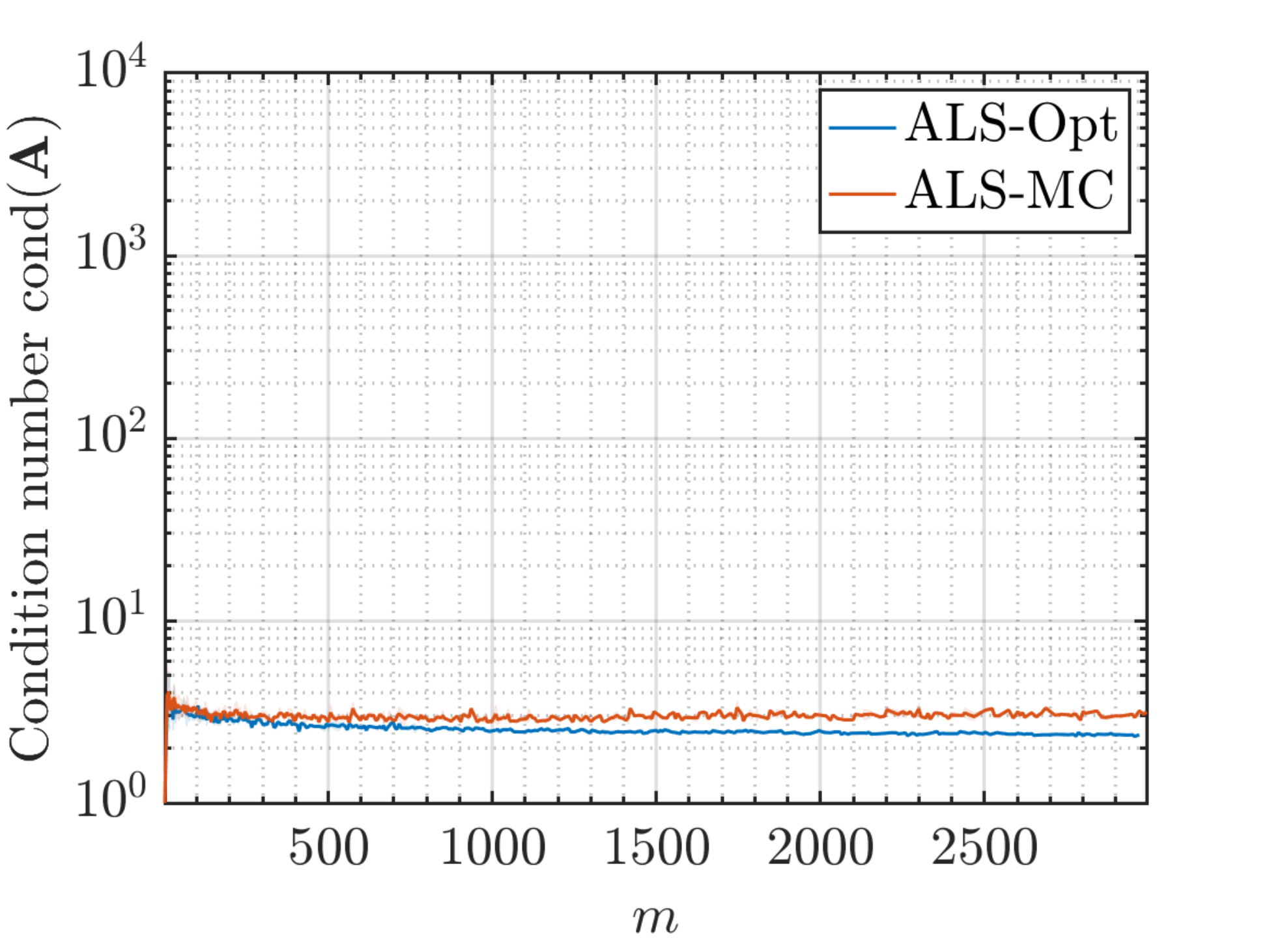}
\\[\errplottextsp]
$d = 1$ & $d = 4$ & $d = 32$
\end{tabular}
\end{small}
\end{center}
\caption{{The same as Fig.\ \ref{fig:fig2} but with $f = f_3$ with $\delta_i = i$ (top two rows) and $\delta_i = i^2$ (bottom two rows). For succinctness, we consider only the values $d = 1,4,32$.}} 
\label{fig:fig4}
\end{figure}

\begin{figure}[t!]
\begin{center}
\begin{small}
 \begin{tabular}{@{\hspace{0pt}}c@{\hspace{\errplotsp}}c@{\hspace{\errplotsp}}c@{\hspace{0pt}}}
\includegraphics[width = \errplotimg]{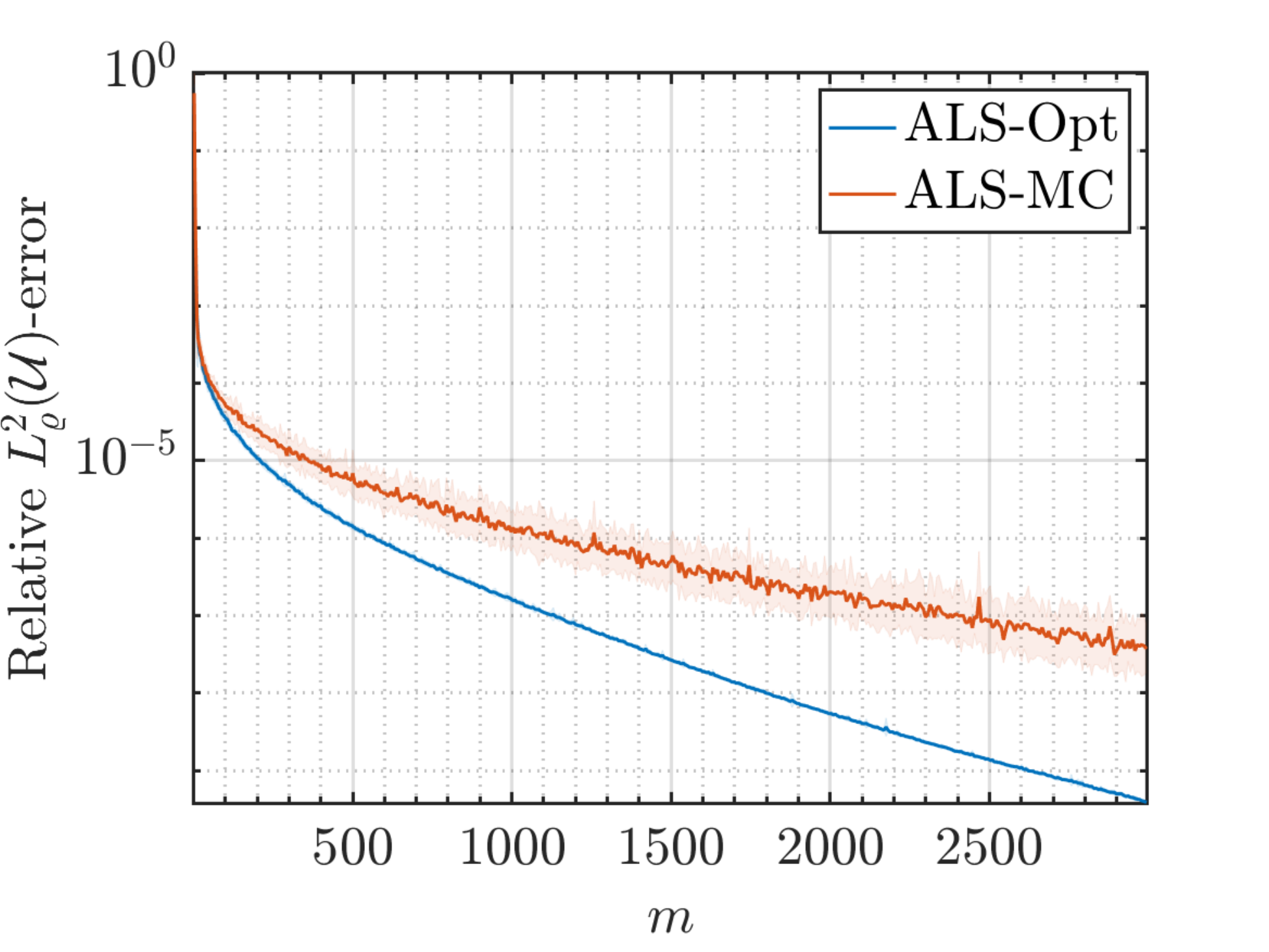}
&
\includegraphics[width = \errplotimg]{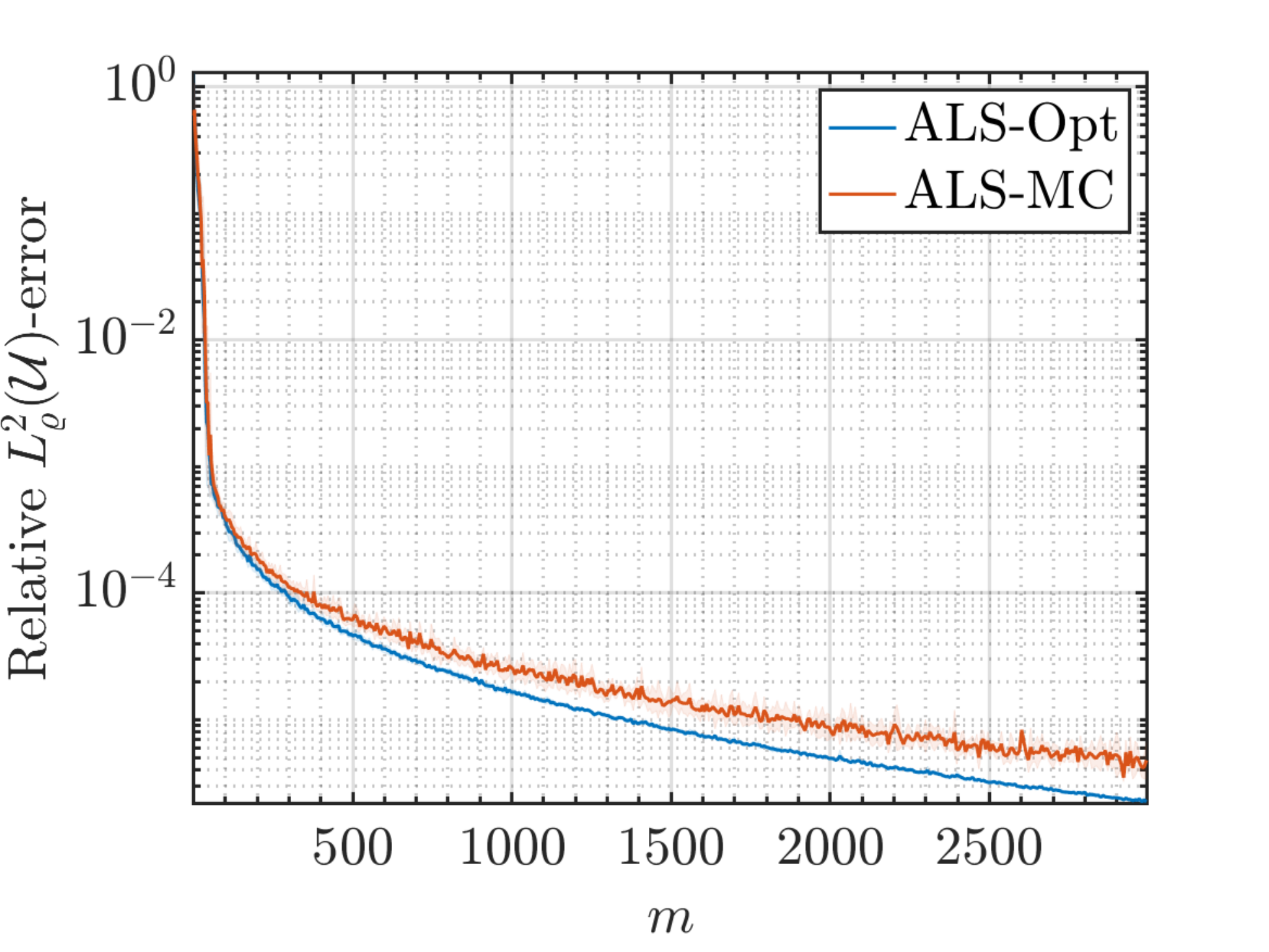}
&
\includegraphics[width = \errplotimg]{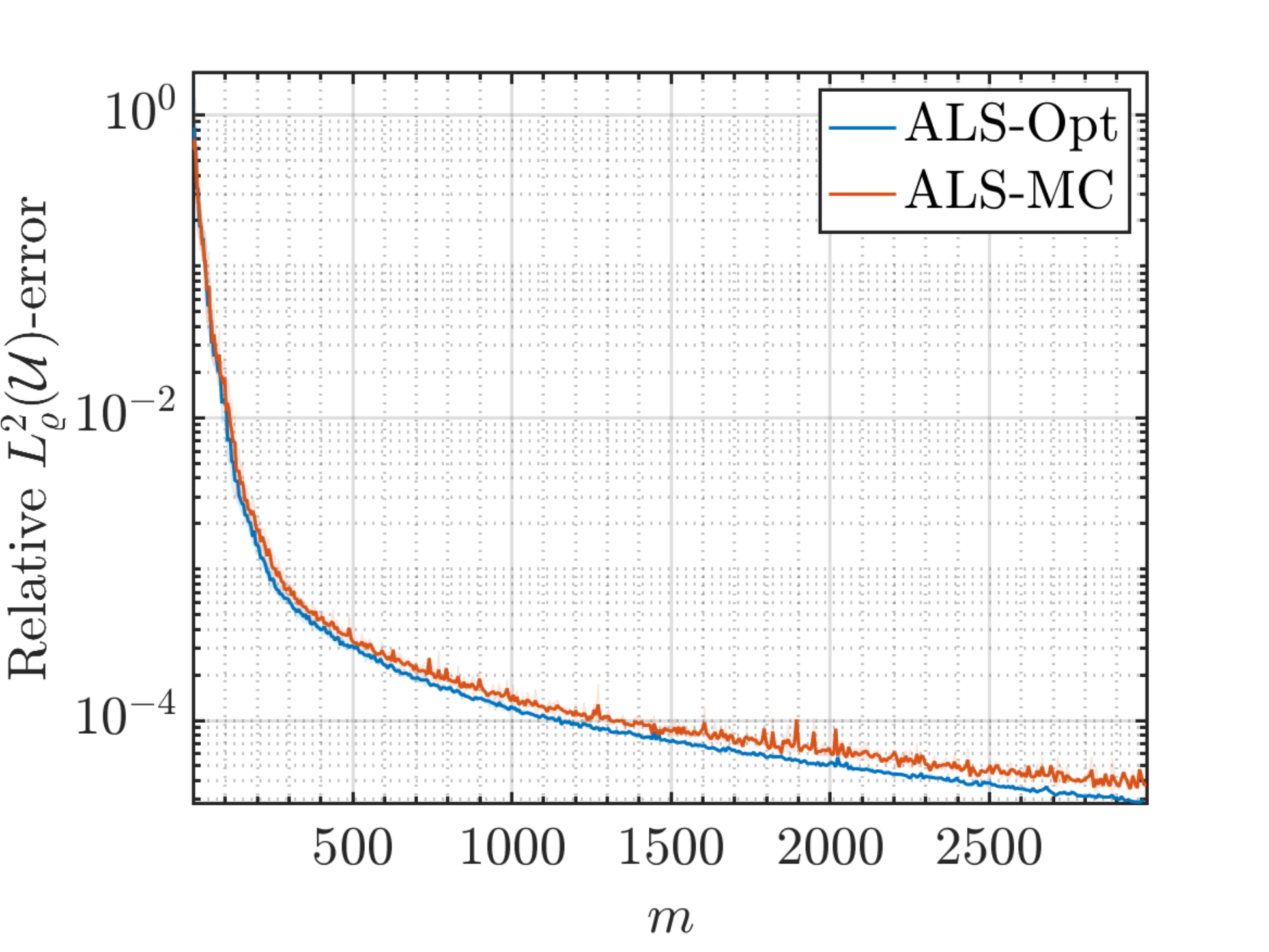}
\\[\errplotgraphsp]
\includegraphics[width = \errplotimg]{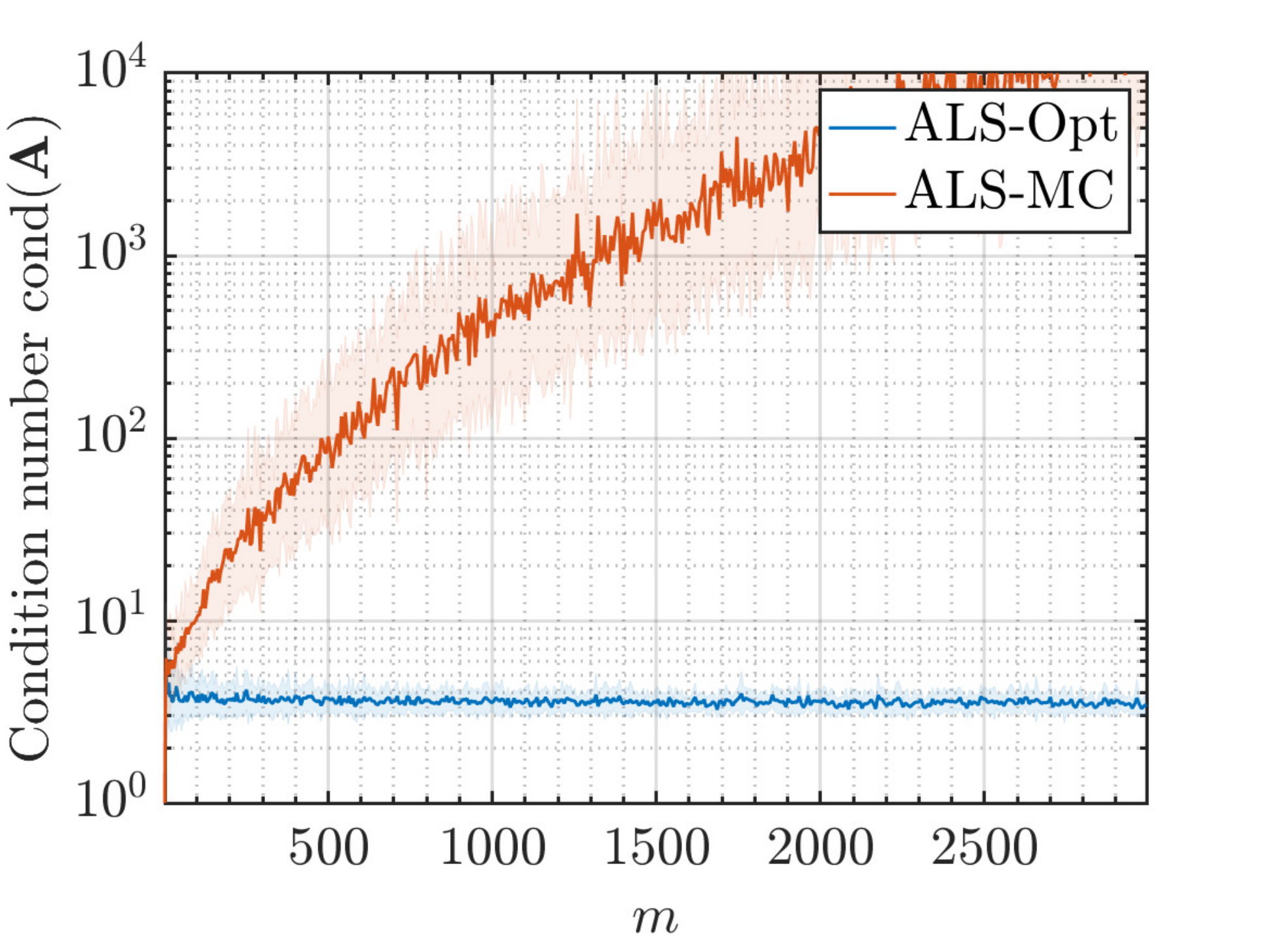}
&
\includegraphics[width = \errplotimg]{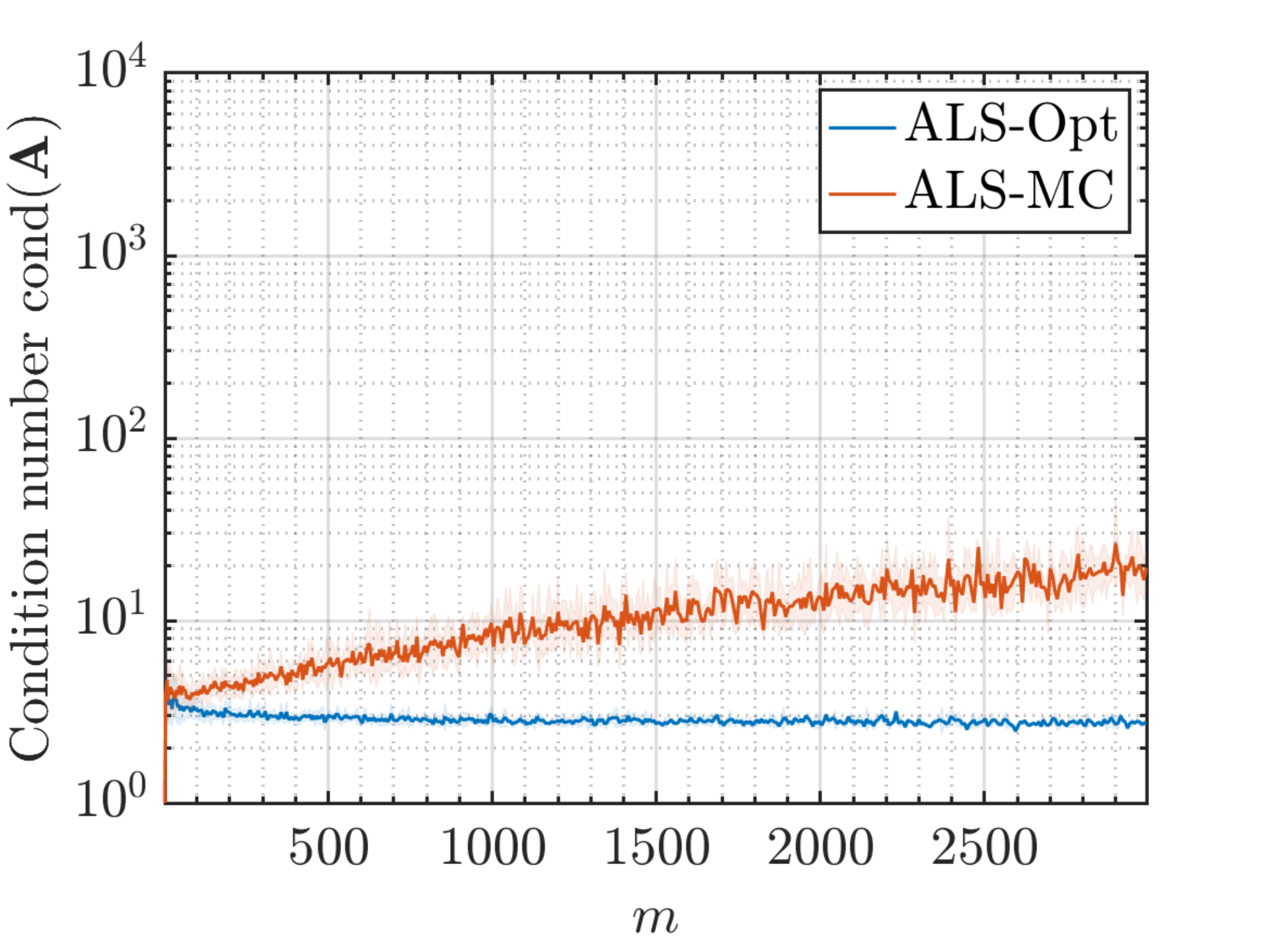}
&
\includegraphics[width = \errplotimg]{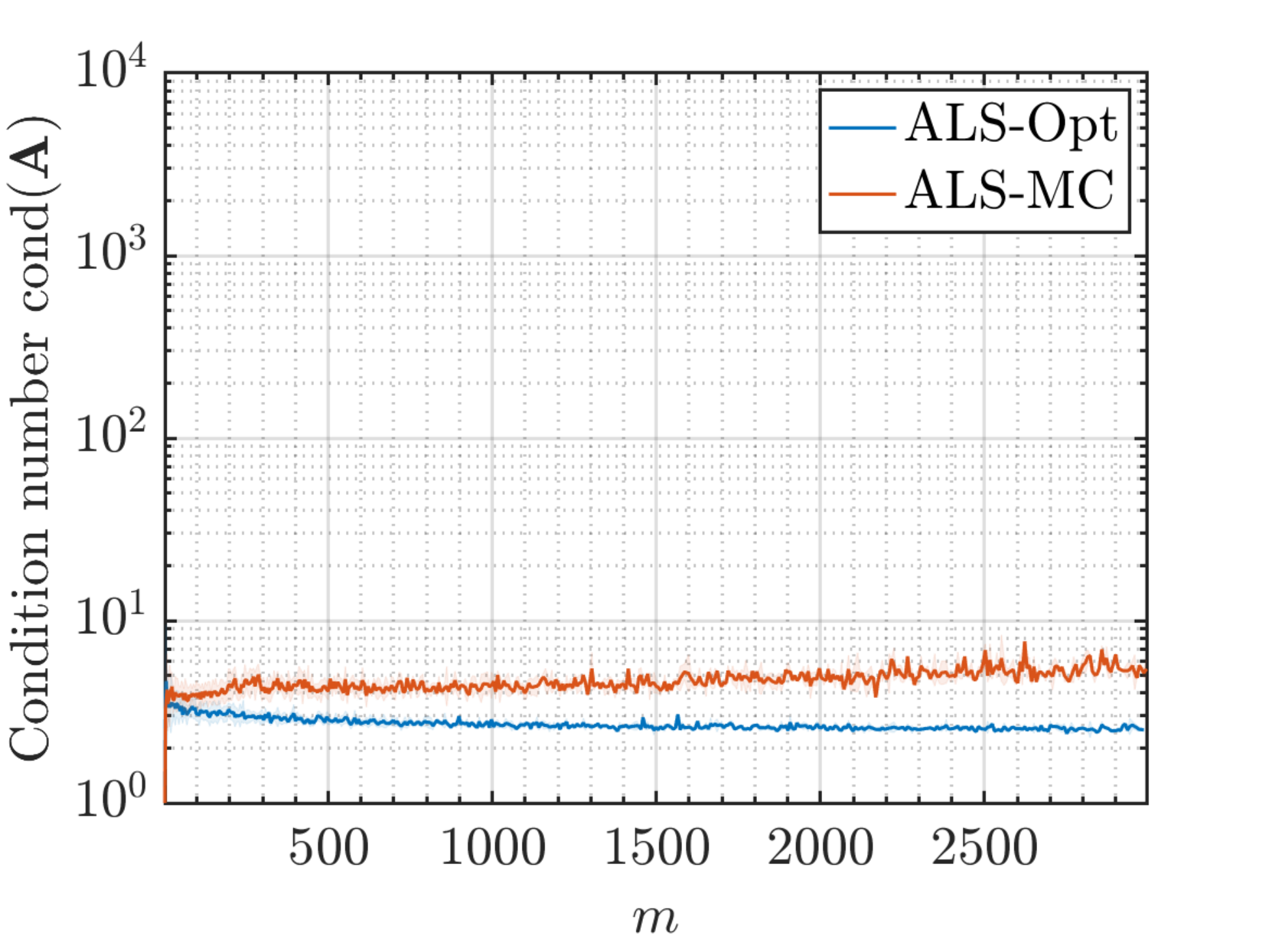}
\\[\errplottextsp]
$d = 2$ & $d = 6$ & $d = 8$
\end{tabular}
\end{small}
\end{center}
\caption{{The same as Fig.\ \ref{fig:fig2} but with $f = f_4$. }} 
\label{fig:fig5}
\end{figure}

{Figs.\ \ref{fig:fig2}--\ref{fig:fig5} show numerical results for the test functions $f_i$, $i = 1,2,3,4$.}
In all cases, there is a substantial difference between the two sampling strategies when $d = 1$. MC sampling leads to a highly ill-conditioned LS problem, with condition number increasing exponentially fast in $m$. Concurrently, the approximation error, while initially decreasing rapidly, eventually begins to increase exponentially, due to the effect of round-off error when solving the ill-conditioned LS problem. {The standard deviation (which is represented by the shaded region) is also high -- an effect that is symptomatic of an ill-conditioned problem.}

This behaviour is well known (see, e.g., \cite{migliorati2014analysis}). The reason {stems from the choice of} log-linear scaling in \eqref{m-scaling-LS}. {This scaling is} asymptotically much smaller than the log-quadratic scaling which suffices for accuracy and stability with MC sampling{. The result is ill-conditioning of the LS problem, and poor stability and accuracy of the corresponding approximation.}

On the other hand, the near-optimal sampling strategy is perfectly well conditioned. The error decreases rapidly to roughly machine epsilon, and remains at this level even as $m$ increases. The {standard deviation} of both the approximation error and condition number are also much lower than in the case of MC sampling.

{Note that functions $f_1$ and $f_3$ are products of univariate functions, which means that their best $n$-term approximation \eqref{best-n-term-inf} can be computed. We include this as a benchmark in Figs.\ \ref{fig:fig2} and \ref{fig:fig4}. It is notable that the error from ALS with near-optimal sampling closely tracks that of the best $n$-term approximation in low dimensions. Yet as $d$ increases there is a widening gap in performance -- in other words, ALS may struggle to identify good index sets in high dimensions. }

{This aside,} the situation in $d = 1$ dimensions insofar as sampling goes is unambiguous: near-optimal sampling leads to a significant improvement over MC sampling. However, the picture begins to change as the dimension increases. For $d = 2$ and $d = 4$ dimensions, the latter still has a growing condition number. However, the rate of growth is much slower than in the one-dimensional case, and as a result, there is {a lower standard deviation and} far less drift in the error as $m \rightarrow \infty$. Moreover, when the dimension is increased further, such effect dramatically lessens. In $d = 16$ or $d = 32$ dimensions, the condition number of MC sampling is virtually the same as that of the near-optimal sampling strategy, as is the error.
In other words, while MC sampling (with log-linear scaling) is a poor strategy in low dimensions, in higher dimensions, {its performance in these examples is very close to that of} the near-optimal strategy. 

\begin{figure}[t!]
\begin{center}
\begin{small}
 \begin{tabular}{@{\hspace{0pt}}c@{\hspace{\errplotsp}}c@{\hspace{\errplotsp}}c@{\hspace{0pt}}}
\includegraphics[width = \errplotimg]{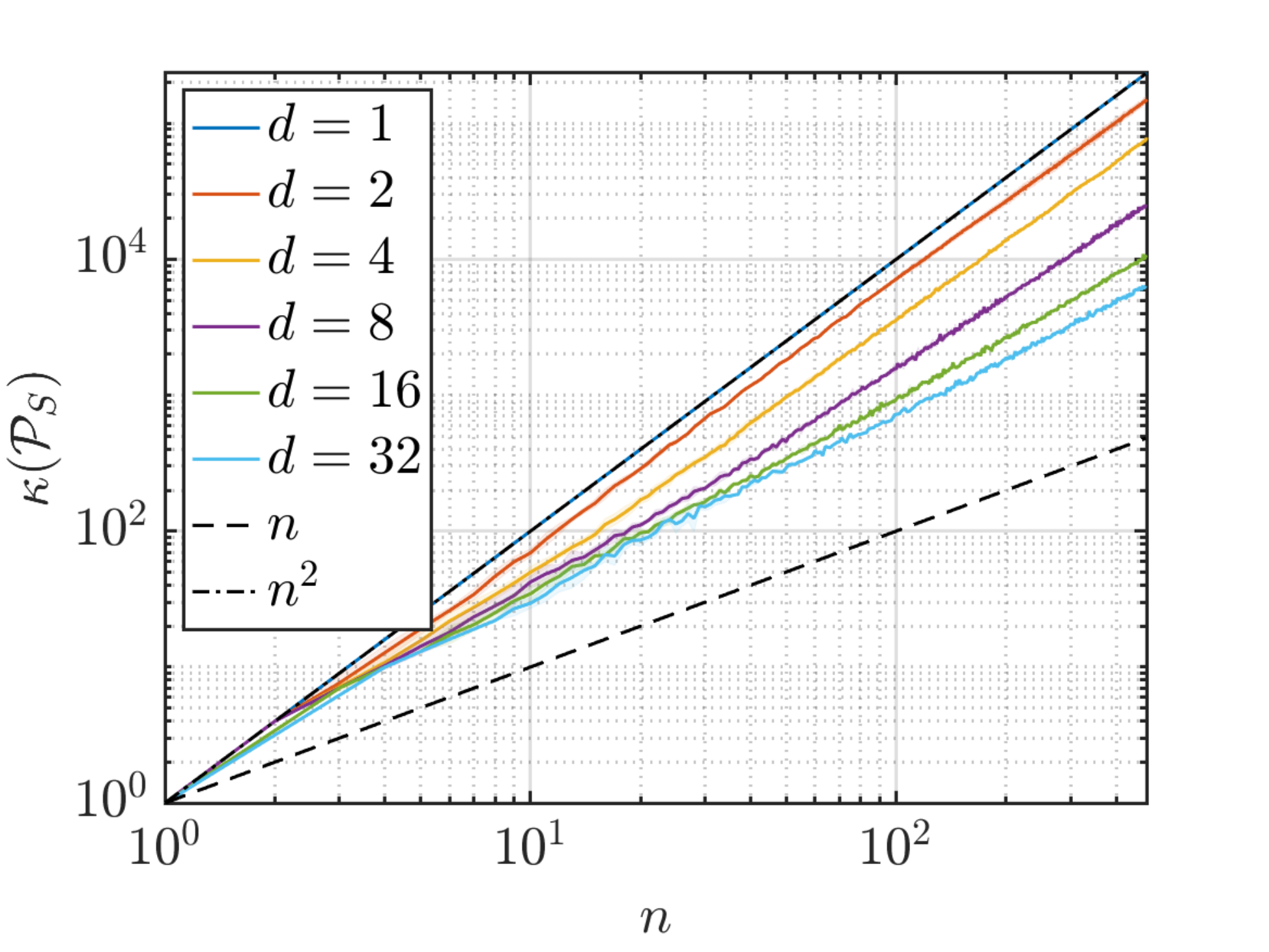}
&
\includegraphics[width = \errplotimg]{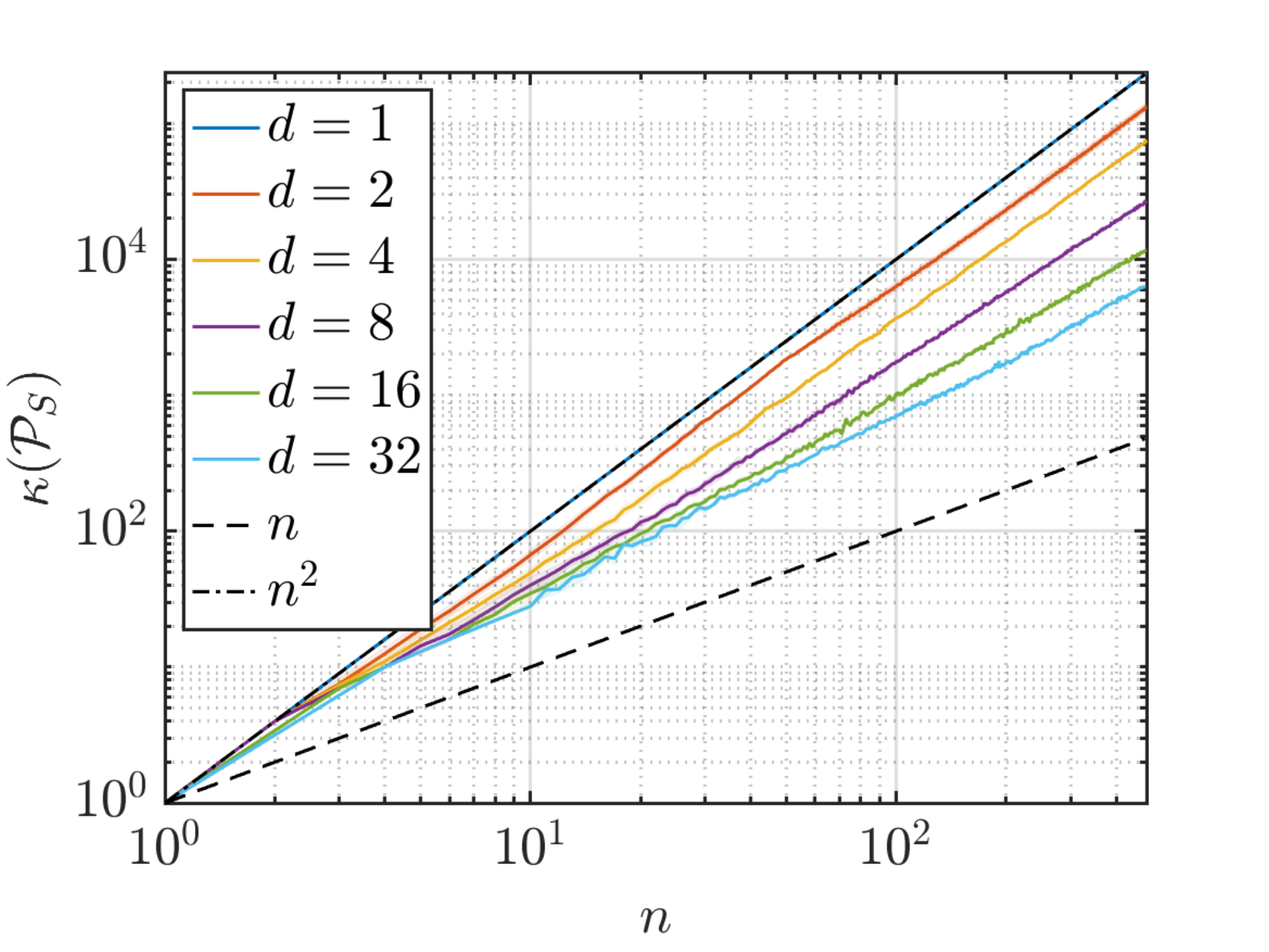}
&
\includegraphics[width = \errplotimg]{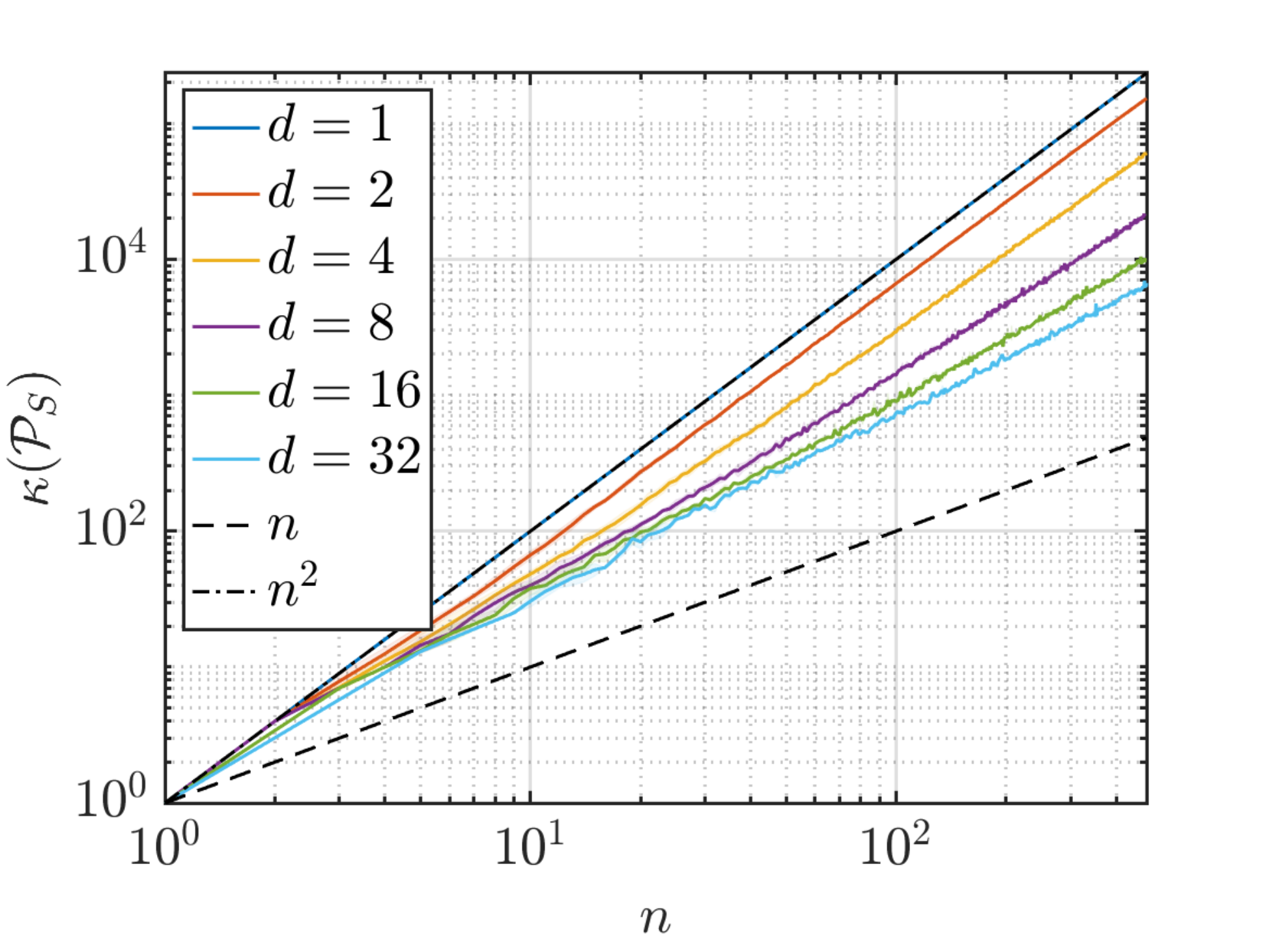}
\\[\errplottextsp]
$f = f_1$ & $f = f_2$ & $f = f_3$, $\delta_i = i$
\end{tabular}
\end{small}
\end{center}
\caption{The quantity $\kappa(\cP_S)$ versus $n$ for the index sets generated by the adaptive LS approximation using MC sampling different dimensions. The dashed lines represent the scalings $n$ (lower) and $n^2$ (upper), respectively.} 
\label{fig:fig6}
\end{figure}

To gain some further insight, in Fig.\ \ref{fig:fig6} we plot the function $\kappa(\cP_S)$ for the index sets $S^{(1)},S^{(2)},\ldots$ produced by the ALS scheme with MC sampling. In one dimension, this function behaves like $n^2$, where $n = n^{(l)}$ is the size of the index set. This is exactly as we expect. The adaptive procedure must produce a lower set, and in one dimension there is only lower set of size $n$: namely, $S = \{0,\ldots,n-1\}$. Hence, we must have $\kappa(\cP_S) = n^2$ (see Proposition \ref{prop:legendre-Christoffel}).
On the other hand, as the dimension increases, we see that $\kappa(\cP_S)$ scales more slowly with $n$. Although this scaling is still faster than the optimal scaling $\kappa(\cP_S) = n$, it again highlights the fact that MC sampling becomes a progressively less bad a sampling strategy as the dimension increases.

\subsection{Illustrative case studies}

Recall that the index sets $S^{(1)},S^{(2)},\ldots$ generated by the adaptive LS approximation scheme are function dependent. As these experiments make clear, the worst-case quadratic scaling of $\kappa(\cP_S)$ identified in Proposition \ref{prop:legendre-Christoffel} is not realized for any of the three functions when $d \geq 2$. To examine this situation more concretely, it is informative to consider some specific cases.

\subsubsection{Additively separable functions}

\begin{figure}[t!]
\begin{center}
\begin{small}
 \begin{tabular}{@{\hspace{0pt}}c@{\hspace{\errplotsp}}c@{\hspace{\errplotsp}}c@{\hspace{0pt}}}
\includegraphics[width = \errplotimg]{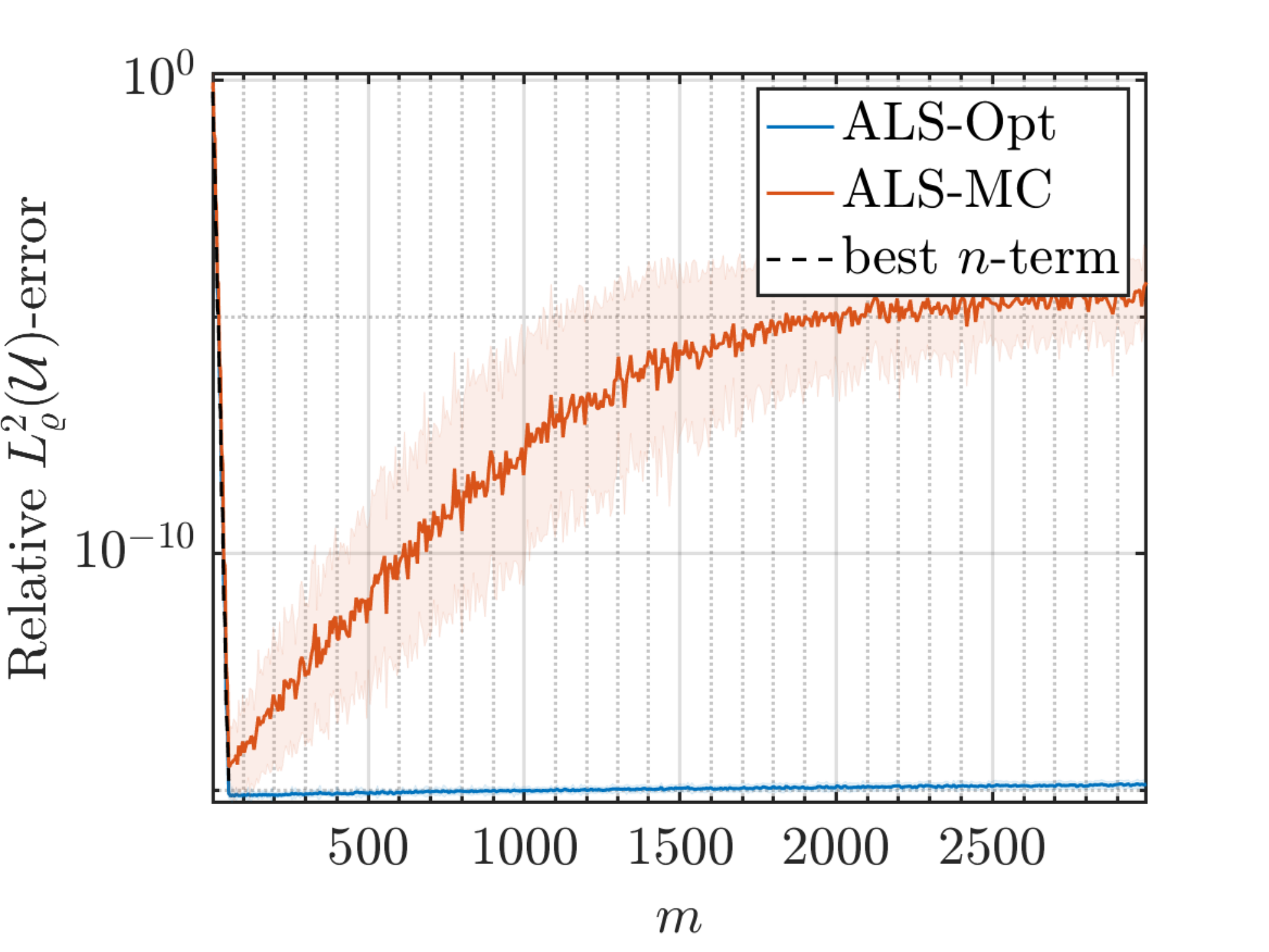}
&
\includegraphics[width = \errplotimg]{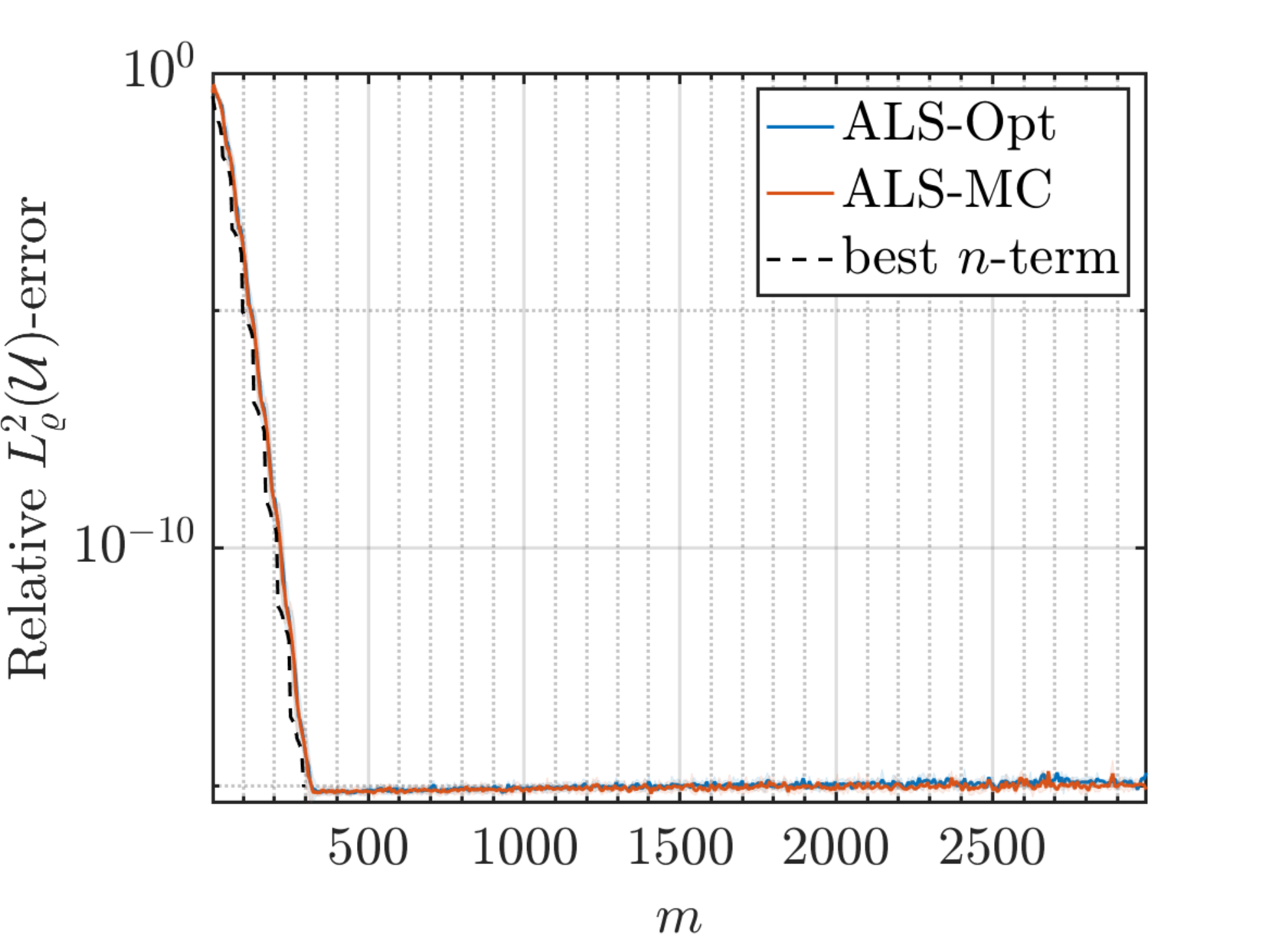}
&
\includegraphics[width = \errplotimg]{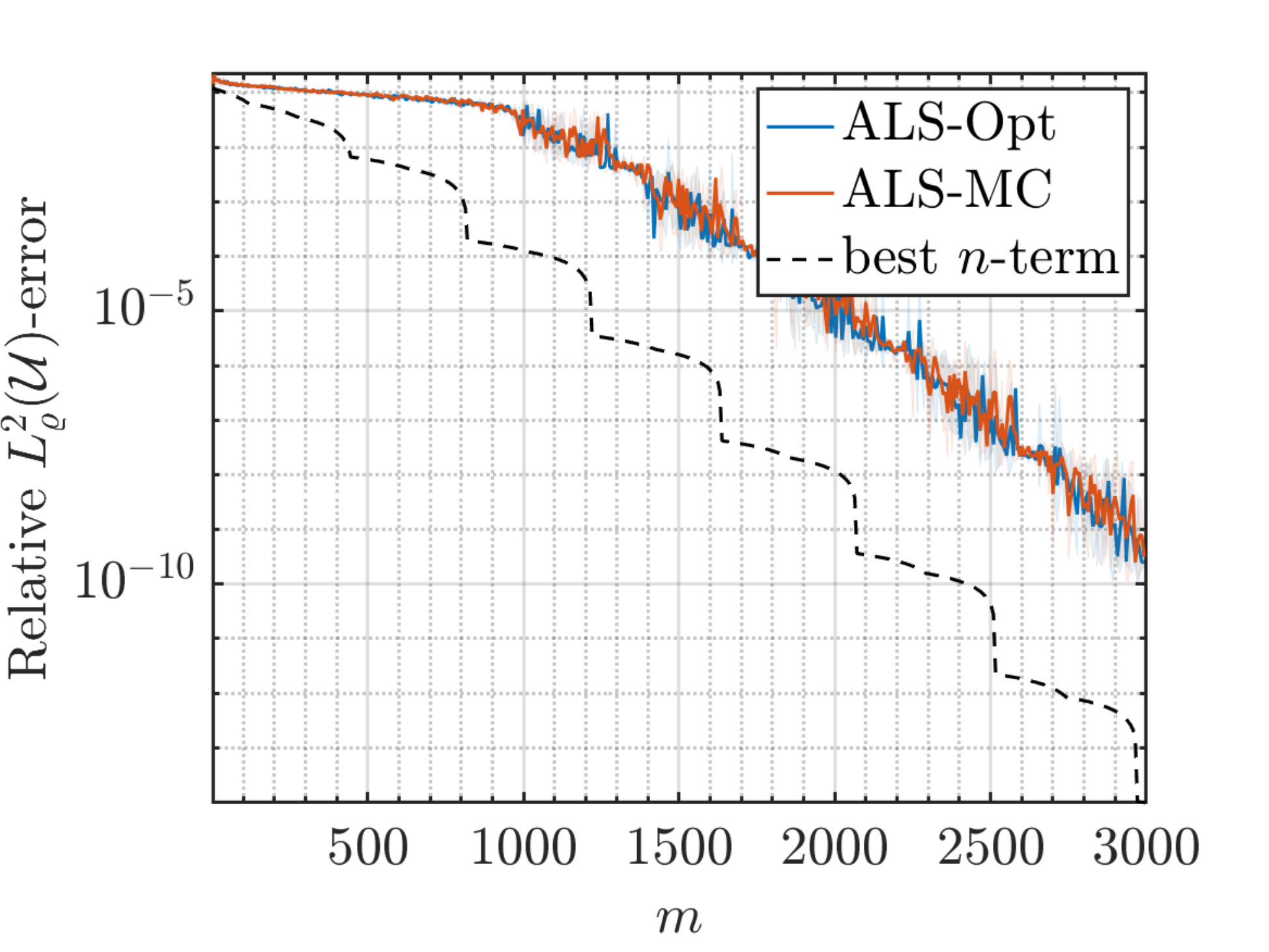}
\\[\errplotgraphsp]
\includegraphics[width = \errplotimg]{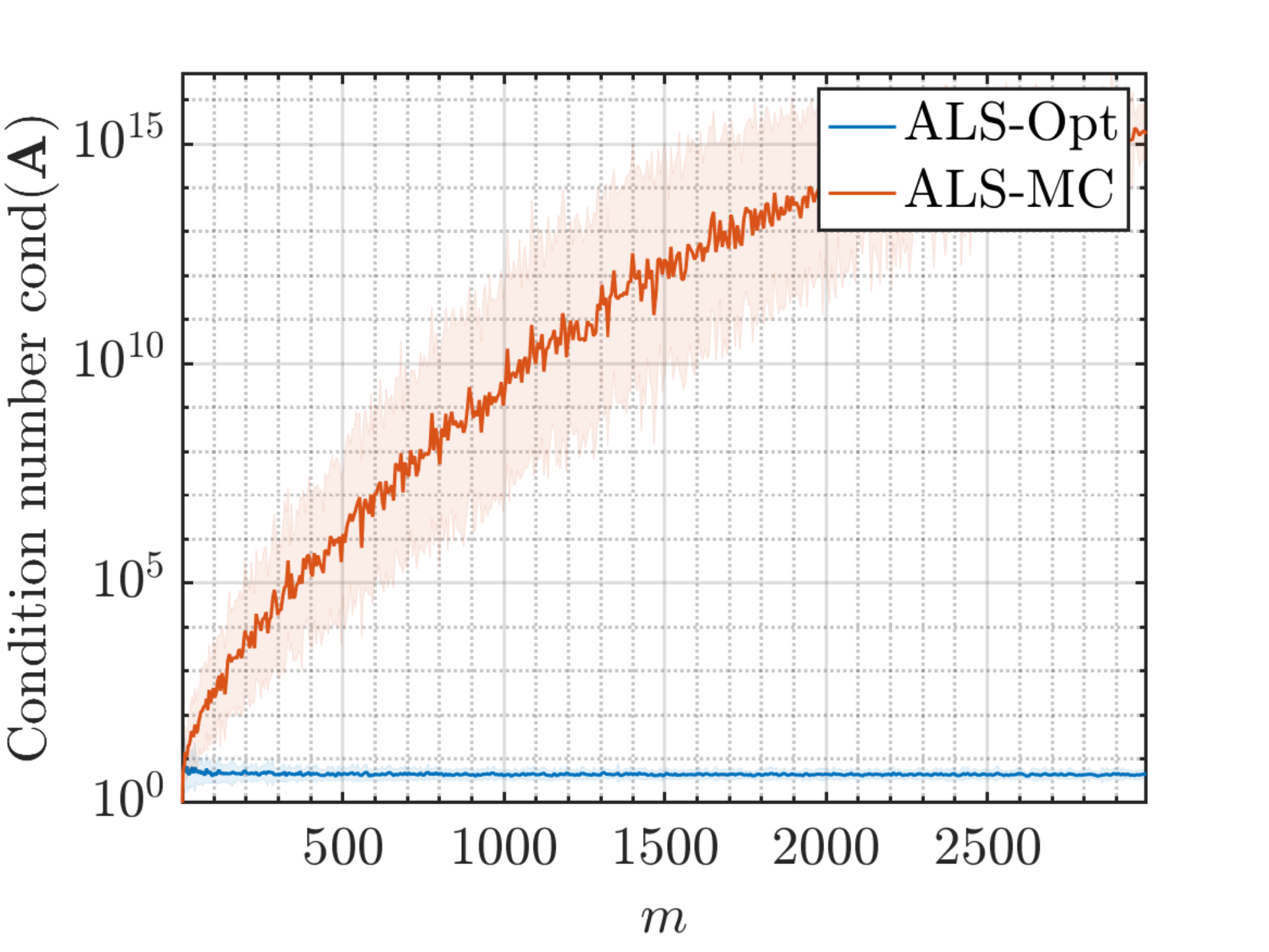}
&
\includegraphics[width = \errplotimg]{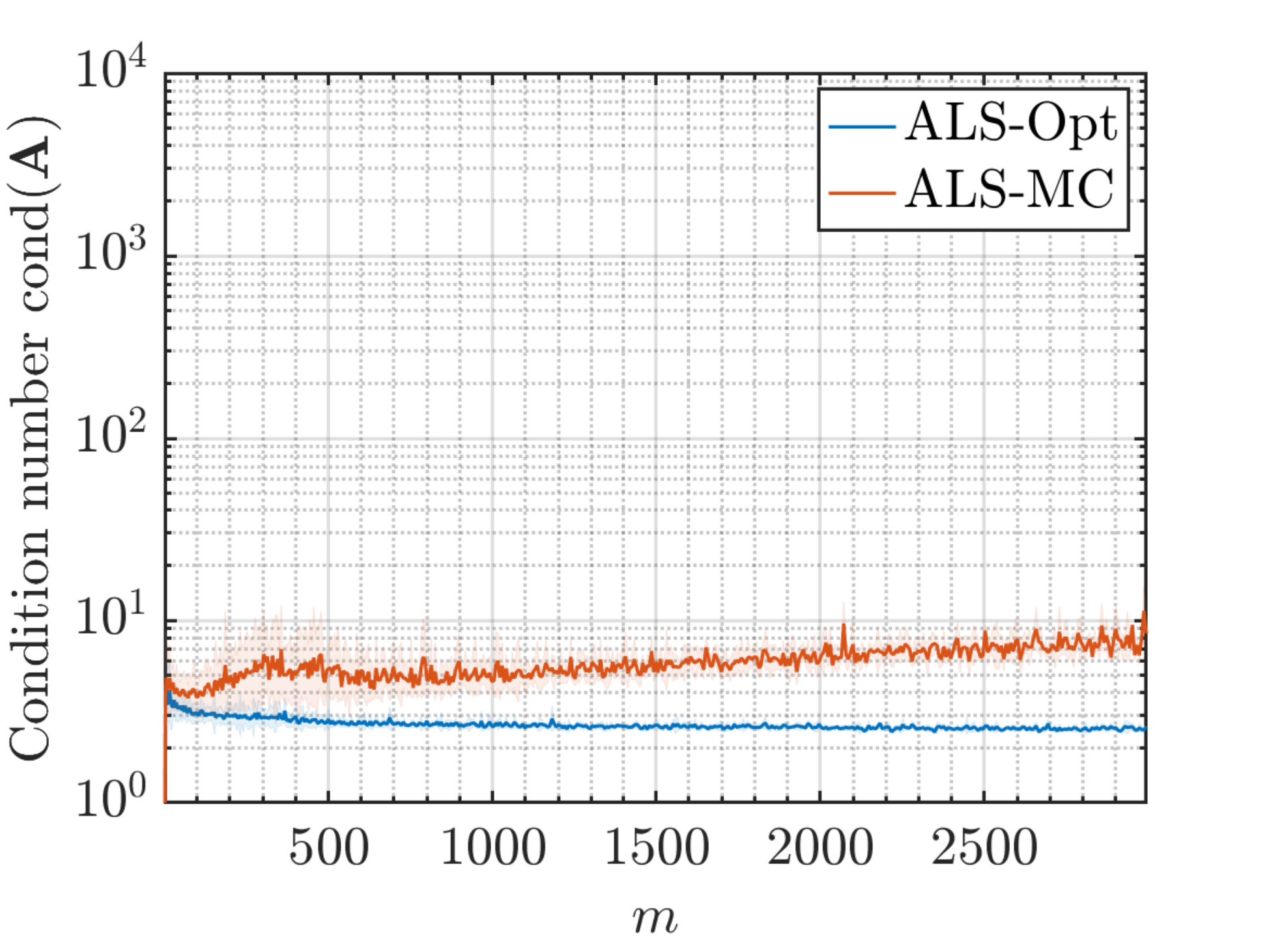}
&
\includegraphics[width = \errplotimg]{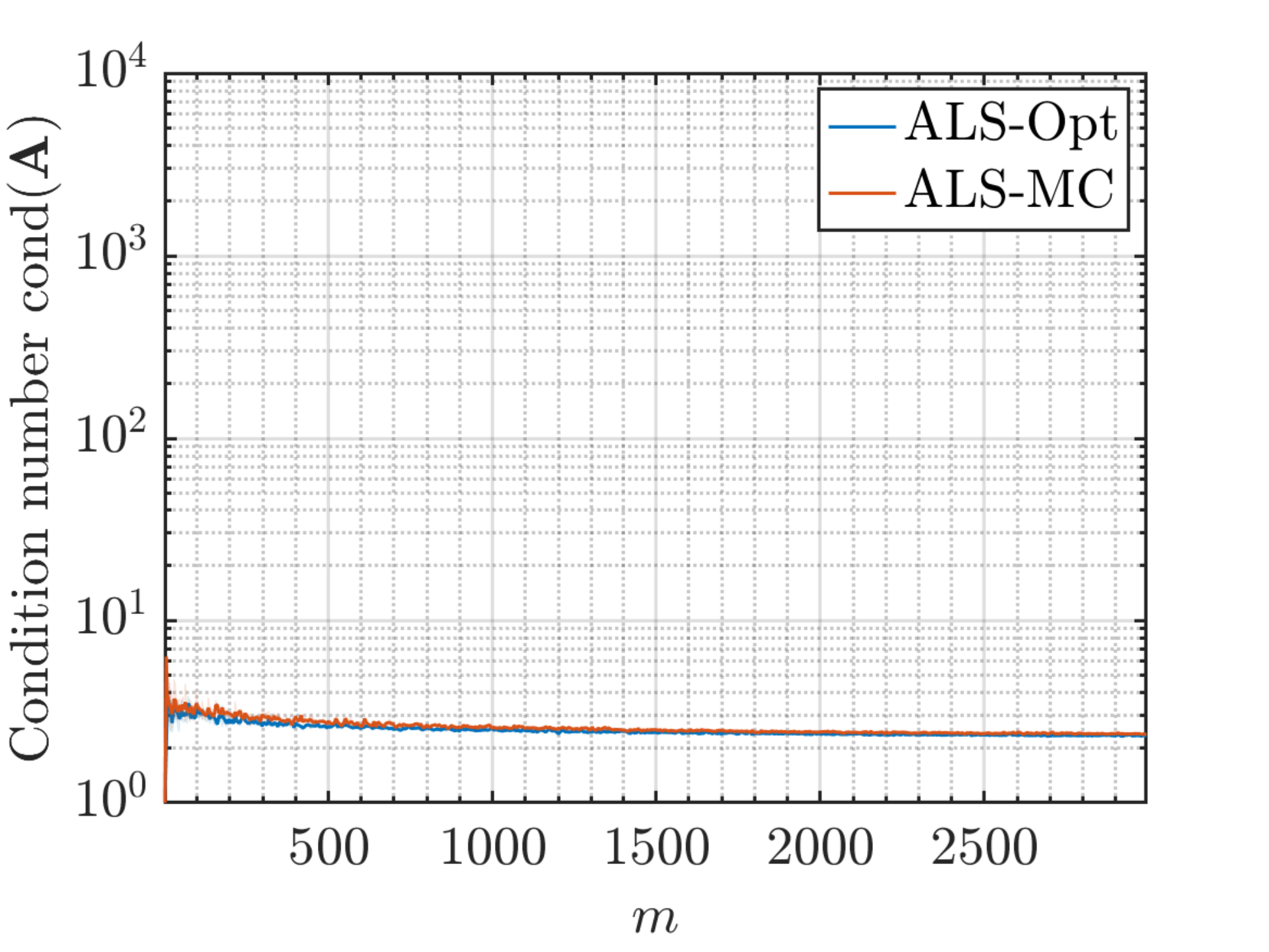}
\\[\errplottextsp]
$d = 1$ & $d = 4$ & $d = 32$
\end{tabular}
\end{small}
\end{center}
\caption{{The same as Fig.\ \ref{fig:fig2} but for the additively separable function $f(\bm{y}) = \sum^{d}_{i=1} 0.3 + \sin(16/15 y_i - 0.7) + \sin^2(16/15 y_i - 0.7)$ from \cite{guo2017gradient}}.} 
\label{fig:fig7}
\end{figure}

First, consider a function of the form $f(\bm{y}) = g(y_1)+\cdots + g(y_d)$,
where $g : [-1,1] \rightarrow \bbC$. 
Write $g = \sum^{\infty}_{\nu=0} d_{\nu} \psi_{\nu}$ and suppose without loss of generality that the coefficients are in nonincreasing order, i.e., $|d_0| \geq |d_1| \geq \cdots $.
Let $\bm{c} = (c_{\bm{\nu}})_{\bm{\nu} \in \bbN^d_0}$ be the coefficients of $f$. Then $c_{\bm{0}} = d \cdot d_0$, $c_{\bm{\nu}} = d_{\nu}$ if $\bm{\nu} = \nu \bm{e}_j$ for some $j \in [d]$ and $\nu \in \bbN$, and $c_{\bm{\nu}} = \bm{0}$ otherwise.
Here $\bm{e}_j \in \bbR^d$ is the multi-index with $1$ in its $j$th component and zero otherwise. 
Now consider the best $n$-term approximation $f_n$ to $f$. Suppose that {$n = d p+1$} for some $p \in \bbN$. Then it is clear that the set $S^*$ defined in \eqref{best-n-term-inf} is precisely
\bes{
S^* = \left \{ k \bm{e}_j : k = 0,\ldots,p,\ j = 1,\ldots,d \right \}.
}
It is then a short argument based on \eqref{K-alternative-def} to show that
\bes{
\kappa(\cP_{S^*}) = d \cdot \kappa(\cQ_p) {+1-d},\qquad \cQ_p : = \spn \{ \psi_0,\ldots,{\psi_{p}} \} \subset L^2_{\varrho}([-1,1]).
}
In particular, {in} the case of Legendre polynomials, Proposition \ref{prop:legendre-Christoffel} implies that 
\bes{
\kappa(\cP_{S^*}) = d \cdot {\kappa(\bbP_{p}) + 1-d = (n-1)^2/d + 2n - 1.}
}
Hence, $\kappa(\cP_{S^*})$ is quadratic in $n$, but with a constant of $1/d$. In higher dimensions, one therefore expects less oversampling to be needed to ensure stability. {In Fig.\ \ref{fig:fig7} we examine the performance of ALS for such a function. As expected, MC sampling gives nearly as good performance in higher dimensions as the near-optimal scheme.}

\subsubsection{Low-dimensional functions}\label{ss:low-dim-fns}

On the other hand, there are functions for which this phenomenon may not occur. As an extreme example, consider a function of the form $f(\bm{y}) = g(y_1)$,
where we assume once more that $g = \sum^{\infty}_{\nu=0} d_{\nu} \psi_{\nu}$ with coefficients in nonincreasing order. Then the set of the best $n$-term approximation is
\be{
\label{Sstar-oned-fun}
S^* = \{ k \bm{e}_1 : 0 \leq k < n \},
}
for which one has $\kappa(\cP_{S^*}) = n^2$
in the case of Legendre polynomials for any $d$. Thus, if MC sampling is used with this index set in combination with log-linear oversampling one expects instability and potential divergence of the LS approximation, regardless of the dimension. More generally, one expects a similar effect to occur for functions of the form $f(\bm{y}) = g(y_{i_1},\ldots,y_{i_r}) + h(\bm{y})$,
where $i_1,\ldots,i_r \in [d]$, $r \ll d$ and $\nmu{h}_{L^2_{\varrho}(\cU)} \ll \nmu{f}_{L^2_{\varrho}(\cU)}$, i.e., functions that are, up to a small perturbation, low dimensional.

\begin{figure}[t!]
\begin{center}
\begin{small}
 \begin{tabular}{@{\hspace{0pt}}c@{\hspace{\errplotsp}}c@{\hspace{\errplotsp}}c@{\hspace{0pt}}}
\includegraphics[width = \errplotimg]{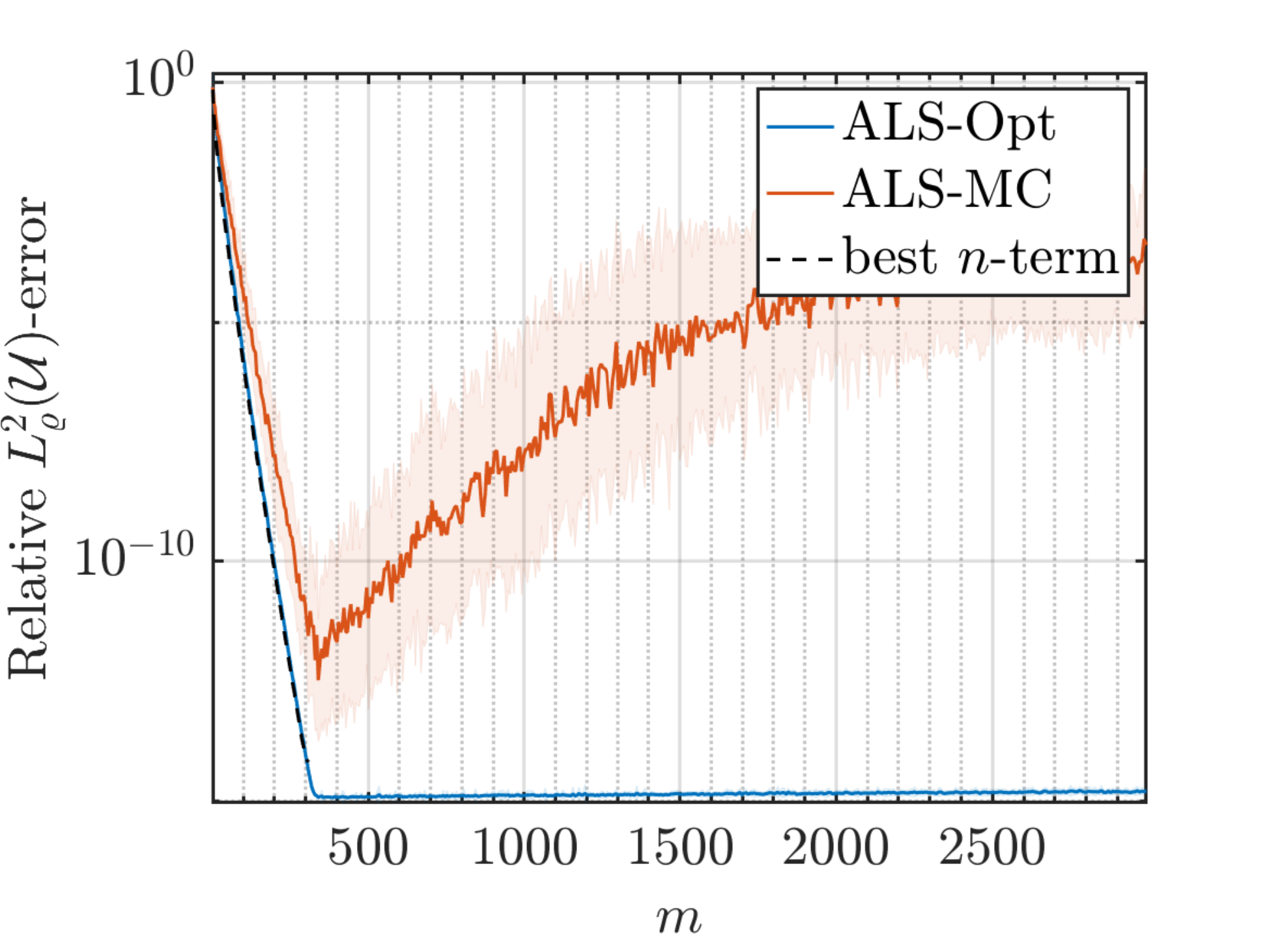}
&
\includegraphics[width = \errplotimg]{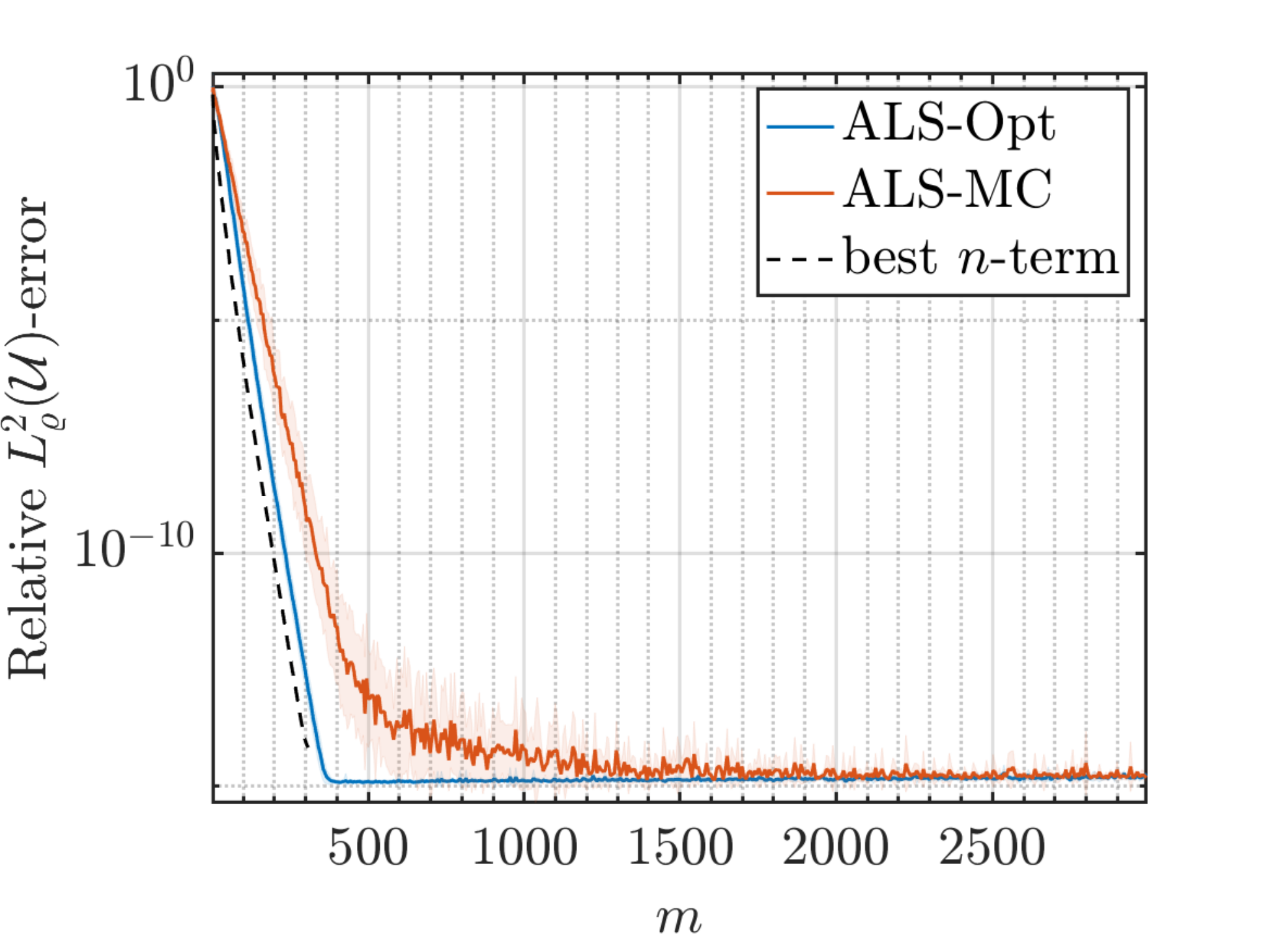}
&
\includegraphics[width = \errplotimg]{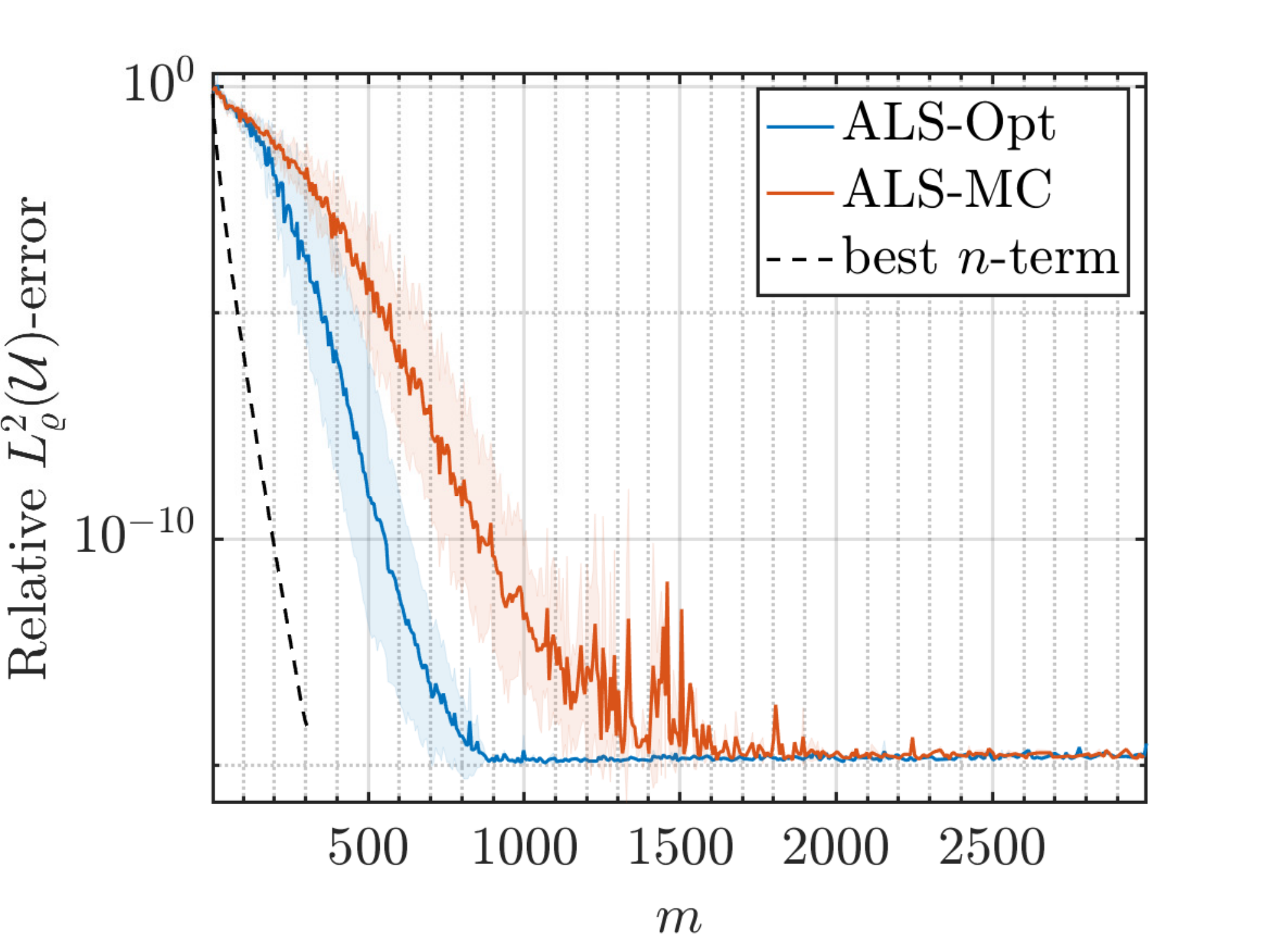}
\\[\errplotgraphsp]
\includegraphics[width = \errplotimg]{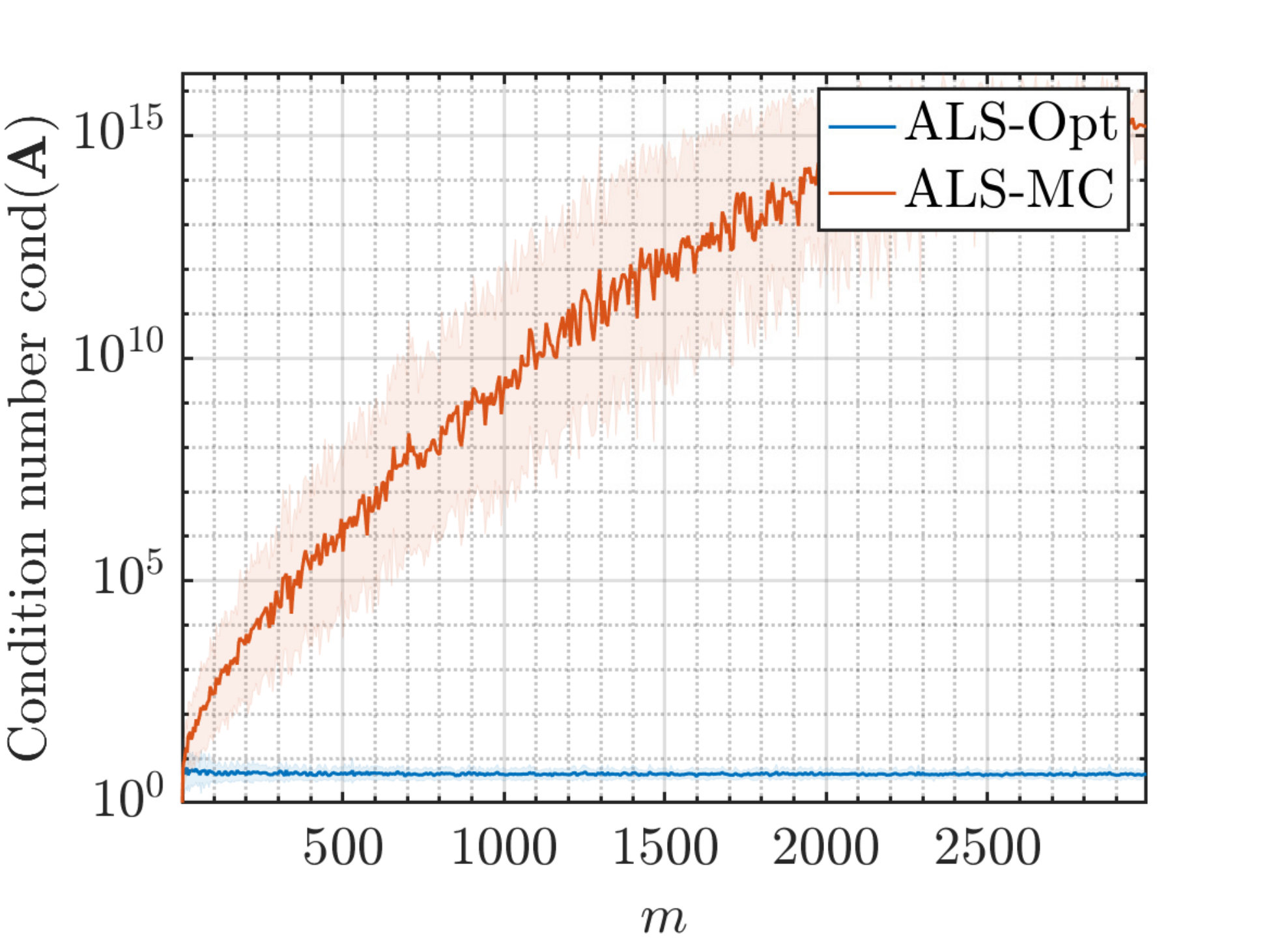}
&
\includegraphics[width = \errplotimg]{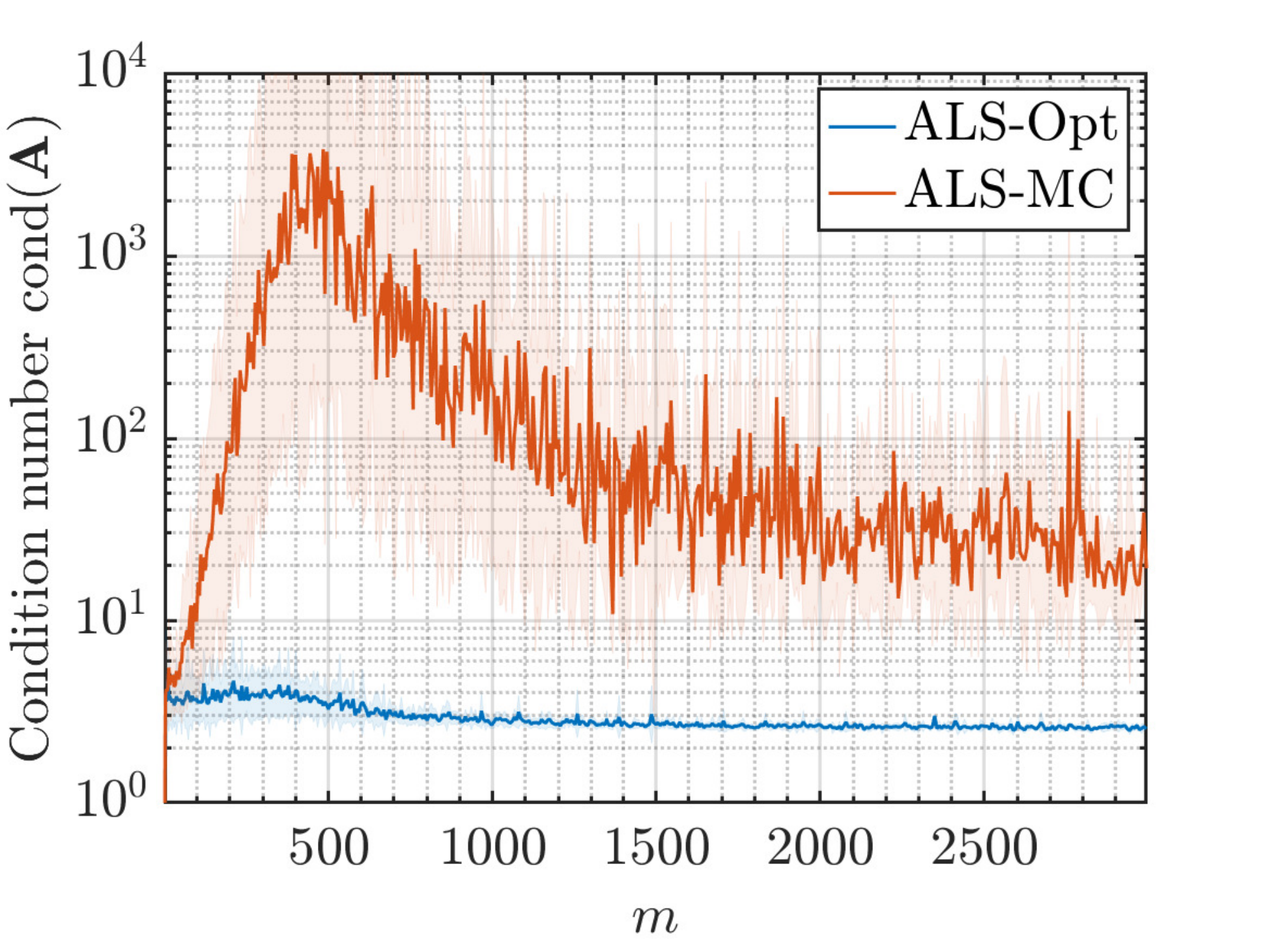}
&
\includegraphics[width = \errplotimg]{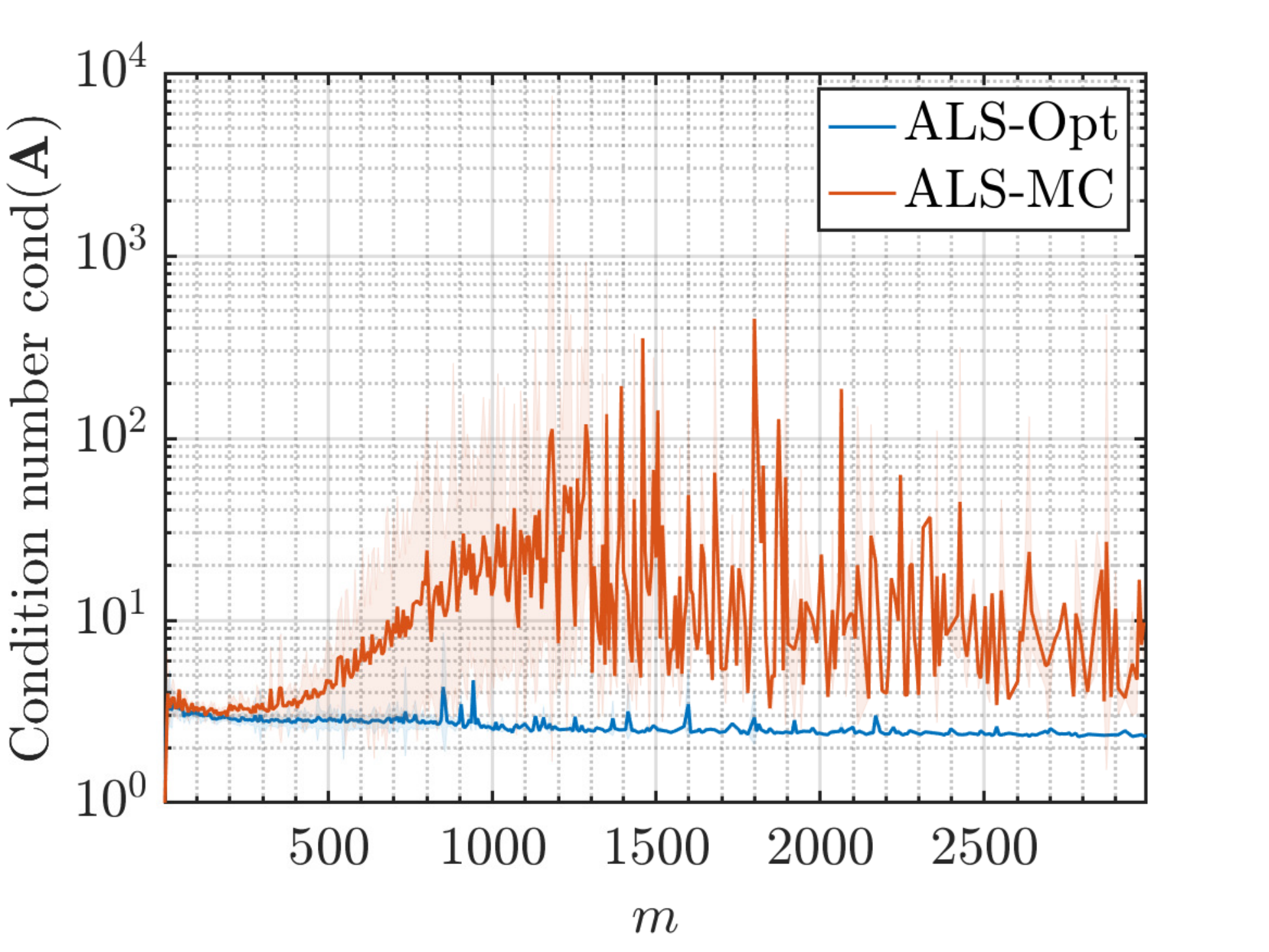}
\\[\errplottextsp]
$d = 1$ & $d = 4$ & $d = 32$
\end{tabular}
\end{small}
\end{center}
\caption{The same as Fig.\ \ref{fig:fig2} but for the function of one variable $f(\bm{y}) = (10-9 y_1)^{-1}$. For succinctness, we consider only the values $d = 1,4,32$.} 
\label{fig:fig8}
\end{figure}

We consider an example of such a function in Fig.\ \ref{fig:fig8}. Interestingly, the results deviate from our expectation: as the dimension increases, MC sampling becomes better conditioned, and the error approaches that of optimal sampling. The reason behind this seeming contradiction stems from the fact that the ALS scheme may not choose the index set of the best $n$-term approximation, i.e., \eqref{Sstar-oned-fun}, in $d \geq 2$ dimensions. Indeed, it uses estimates for the polynomial coefficients based on the current iterate to construct each new index set. By computing `suboptimal' index sets, the LS approximation based on MC sampling actually performs relatively better in higher dimensions. However, it is notable that its performance in comparison to the near-optimal scheme is still worse than for the functions considered in Figs.\ \ref{fig:fig2}--\ref{fig:fig5}.

\section{Polynomial approximation theory in infinite dimensions}\label{s:polyapp-inf-dim}

We now turn our attention to theoretical analysis of the phenomenon seen in the previous section. As mentioned in \S \ref{s:main-res-summary}, we shall do this by considering polynomial approximation of infinite-dimensional holomorphic functions. In this section, we {describe the standard setup for this approach}, before presenting our theoretical result in \S \ref{s:main-res}. For further information on the material presented in this section, see, e.g., \cite[Chpt.\ 3]{adcock2022sparse} and \cite{cohen2015approximation,chkifa2015breaking}.

\subsection{Infinite-dimensional setup}
\label{ss:inf_dim_setup}

We now consider scalar-valued functions of the form $f : \cU \rightarrow \bbC$, where $\cU = [-1,1]^{\bbN}$. As in \S \ref{ss:univariate}, we let $\varrho$ be a probability measure on $[-1,1]$. The Kolmogorov extension theorem  {(see \S \ref{s:inf_dim_measures})} guarantees the existence of a probability measure on $\cU$ formed as the infinite tensor-product of this measure. Abusing notation, we denote this measure as $\varrho = \varrho \times \varrho \times \cdots$ and let $L^2_{\varrho}(\cU)$ be the corresponding Lebesgue space of square-integrable functions.

We construct an orthonormal polynomial basis of $L^2_{\varrho}(\cU)$ by tensor products. Let $\bm{\nu} = (\nu_1,\nu_2,\ldots) \in \bbN^{\bbN}_0$ be an infinite multi-index and define the set
\be{
\label{F-set-def}
\cF = \{ \bm{\nu} \in \bbN^{\bbN}_{0} : \nm{\bm{\nu}}_0 < \infty \} \subset \bbN^{\bbN}_0.
} 
Here $\nm{\bm{\nu}}_0 = | \{ i : \nu_i \neq 0 \} |$ is the $\ell^0$-`norm'. Thus, $\cF$ is the set of infinite multi-indices with at most finitely-many nonzero terms.
For any $\bm{\nu} = (\nu_i)^{\infty}_{i=1} \in \cF$, we define the corresponding tensor-product orthonormal polynomial $\Psi_{\bm{\nu}}$ as
\bes{
\Psi_{\bm{\nu}}(\bm{y}) = \prod_{i \in \mathrm{supp}(\bm{\nu})} \psi_{\nu_i}(y_i),\quad \bm{y} = (y_i)^{\infty}_{i=1} \in \cU,
}
where $\mathrm{supp}(\bm{\nu}) = \{ i : \nu_i \neq 0 \}$ is the \textit{support} of $\bm{\nu}$. 
It can be shown that the set of such functions $\{ \Psi_{\bm{\nu}} : \bm{\nu} \in \cF \} \subset L^2_{\varrho}(\cU)$
forms an orthonormal basis for $L^2_{\varrho}(\cU)$. Hence, much as in the finite-dimensional case (see \S \ref{ss:multivar-poly-app}), any function $f \in L^2_{\varrho}(\cU)$ has an expansion
\be{
\label{f-exp-inf}
f = \sum_{\bm{\nu} \in \cF} c_{\bm{\nu}} \Psi_{\bm{\nu}},\qquad \text{where }
c_{\bm{\nu}} = \ip{f}{\Psi_{\bm{\nu}}}_{L^2_{\varrho}(\cU)} = \int_{\cU} f(\bm{y}) \overline{\Psi_{\bm{\nu}}(\bm{y})} \D \varrho(\bm{y}).
}

\subsection{\boldmath Smoothness and the class of $(b,\varepsilon)$-holomorphic functions}

We now introduce the class of holomorphic functions considered. In the univariate setting, it is well known that the convergence rate of a polynomial approximation to $f$ is dictated by the largest \textit{Bernstein ellipse} within which $f$ is holomorphic {(see, e.g., \cite[Chpt.\ 8]{trefethen2013approximation})}. Recall that the Bernstein ellipse of parameter $\rho > 1$ is defined by
\bes{
\cE_{\rho} = \{ (z + z^{-1})/2 : z \in \bbC,\ 1 \leq |z| \leq \rho  \} \subset \bbC.
}
For convenience, we set $\cE_{1} = [-1,1]$. In the infinite-dimensional setting, we consider Cartesian products of Bernstein ellipses. Given a parameter $\bm{\rho} = (\rho_1,\rho_2,\ldots) \in [1,\infty)^{\bbN}$ we define the \textit{Bernstein polyellipse} of parameter $\bm{\rho}$ as
\bes{
\cE_{\bm{\rho}} = \cE_{\rho_1} \times \cE_{\rho_2} \times \cdots \subset \bbC^{\bbN}.
}

\defn{[$(\bm{b},\varepsilon)$-holomorphy]
\label{d:holomorphy}
Let $\bm{b} \in [0,\infty)^{\bbN}$ and $\varepsilon > 0$. A function $f : \cU \rightarrow \bbC$ is \textit{$(\bm{b},\varepsilon)$-holomorphic} if it is holomorphic in every Bernstein polyellipse $\cE_{\bm{\rho}}$ with parameter $\bm{\rho} = (\rho_i)^{\infty}_{i=1} \in [1,\infty)^{\bbN}$ satisfying
\be{
\label{b-eps-holo}
\sum^{\infty}_{i=1} \left ( \frac{\rho_i + \rho^{-1}_i}{2} - 1 \right ) b_i \leq \varepsilon.
}
}
See, e.g., \cite{chkifa2015breaking,schwab2019deep}.
For convenience we denote the corresponding region as
\bes{
\cR_{\bm{b},\varepsilon} = \bigcup \{ \cE_{\bm{\rho}} : \bm{\rho} \in [1,\infty)^{\bbN},\ \text{$\bm{\rho}$ satisfies \eqref{b-eps-holo}} \} \subseteq \bbC^{\bbN}.
}
We also write $\cH(\bm{b},\varepsilon) =  \{ \text{$f : \cU \rightarrow \bbC$ $(\bm{b},\varepsilon)$-holomorphic},\ \nm{f}_{L^{\infty}(\cR_{\bm{b},\epsilon})} \leq 1 \}$
for the set of functions that are holomorphic in $\cR_{\bm{b},\varepsilon}$ with uniform norm at most one. 

Note that the sequence $\bm{b}$ determines the type of anisotropic behaviour of functions in $\cH(\bm{b},\varepsilon)$. Indeed, if $b_j$ is large, then \eqref{b-eps-holo} holds only for small values of $\rho_j$, meaning that $f$ is less smooth with respect to the variable $y_j$. Conversely, if $b_j$ is small (or even $b_j = 0$), then $f$ is more smooth (entire) in the variable $y_j$.

The additional parameter $\varepsilon$ in Definition \ref{d:holomorphy} is technically redundant. However, it is customary to include it because of the parametric DE context. As observed previously, parametric DEs were the original motivations for the study of this class of functions. As noted (see also \S \ref{s:pDE-background}), the parametric solution maps of many different classes of parametric DEs are $(\bm{b},\varepsilon)$-holomorphic functions.

\rem{[Finite-dimensional functions]
\label{rem:fin-dim}
Definition \ref{d:holomorphy} is somewhat complicated, in that it requires the function to have a holomorphic extension to a union of Bernstein polyellipses. This is needed in infinite dimensions to obtain algebraic rates of convergence of the best $n$-term approximation. In finite dimensions, it is enough for the function to be holomorphic in a single Bernstein polyellipse. However, any such function can be considered within this definition. Indeed, let $f : [-1,1]^d \rightarrow \bbC$ be a function of finitely-many variables that is holomorphic in the finite-dimensional Bernstein polyellipse $\cE_{\bar\rho_1} \times \cdots \times \cE_{\bar\rho_d} \subset \bbC^d$. Now let 
\eas{
b_i = \varepsilon \left ( (\bar\rho_i + \bar\rho^{-1}_{i})/2 - 1 \right )^{-1},\ i \in [d],\qquad b_i = 0,\ i \in \bbN \backslash [d].
}
Then the extension of $f$ to a function of infinitely-many variables is $(\bm{b},\varepsilon)$-holomorphic. Hence, the various results that follow also apply to finite-dimensional functions.
}

\subsection{\boldmath  Best $n$-term polynomial approximation in $\cH(b,\varepsilon)$}

Let $f \in L^2_{\varrho}(\cU)$. As in \S \ref{ss:choosing-S}, we consider $n$-term approximations to $f$, i.e., those taking the form
\bes{
f \approx f_S : = \sum_{\bm{\nu} \in S} c_{\bm{\nu}} \Psi_{\bm{\nu}},
}
where, in this case, $S \subset \cF$, $|S| = n$. Also as before, we write 
\bes{
f_n = f_{S^*},\qquad \text{where }S^* \in \argmin{} \{ \nmu{f - f_S}_{L^2_{\varrho}(\cU)} : S \subset {\mathcal{F}},\ |S| = n \}
}
for the best $n$-ter{m} approximation to $f$.
The following result is well known (see, e.g., \cite[Thm.\ 3.28]{adcock2022sparse} or \cite[\S 3.2]{cohen2015approximation}, as well as \S \ref{ss:best-n-term}). It demonstrates that the best $n$-term approximation of any function in $\cH(\bm{b},\varepsilon)$ converges with algebraic rate.

\thm{
[Algebraic convergence of the best $n$-term approximation]
\label{t:best-s-term}
Let $\varepsilon > 0$ and {$\bm{b} \in [0,\infty)^{\bbN}$ be such that $\bm{b} \in \ell^p(\bbN)$ for some $0 < p < 1$}. Then
\bes{
\nmu{f - f_n}_{L^2_{\varrho}(\cU)} \leq C = C(\bm{b},\varepsilon,p) \cdot n^{\frac12-\frac1p},\qquad \forall f \in \cH(\bm{b},\varepsilon),\ n \in \bbN.
}
}

\rem{[Differences between finite and infinite dimensions]
\label{rem:fin-inf-diff}
The appearance of algebraic rates is one way in which the infinite-dimensional setting differs from the $d$-dimensional setting. When $d$ is fixed and $n \rightarrow \infty$, the best $n$-term polynomial approximation of any $d$-dimensional function that is holomorphic in a Bernstein polyellipse $\cE_{\rho_1} \times \cdots \cE_{\rho_d}$ converges exponentially fast in $n^{\frac1d}$ {(and therefore faster than any algebraic power of $n$)} with the precise rate depending on the parameters $\rho_1,\ldots,\rho_d$ {(see \cite[\S 3.9]{cohen2015approximation} or \cite[\S 3.5-3.6]{adcock2022sparse})}.
{However, the number $n^{\frac1d}$ grows exceedingly slowly with $n$ for moderate to large $d$. Hence, this asymptotic regime of exponential convergence is almost never encountered in practice, besides low-dimensional cases. In this case, the practical convergence behaviour in finite, but large dimensions is often better described by the above algebraic rates -- see next for an illustration.}
}

\begin{figure}[t!]
\begin{center}
\begin{small}
 \begin{tabular}{@{\hspace{0pt}}c@{\hspace{\errplotsp}}c@{\hspace{\errplotsp}}c@{\hspace{0pt}}}
\includegraphics[width = \errplotimg]{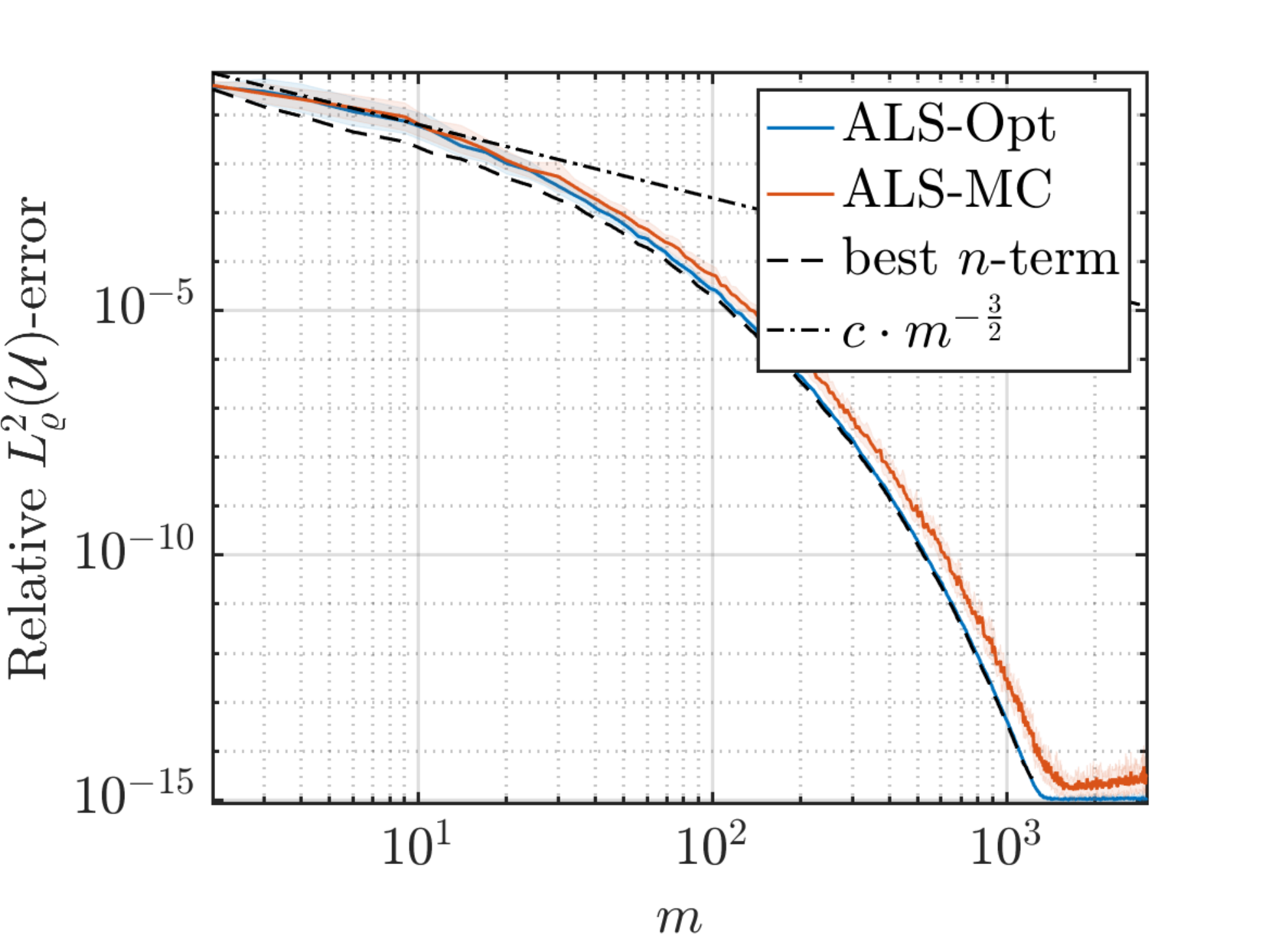}
&
\includegraphics[width = \errplotimg]{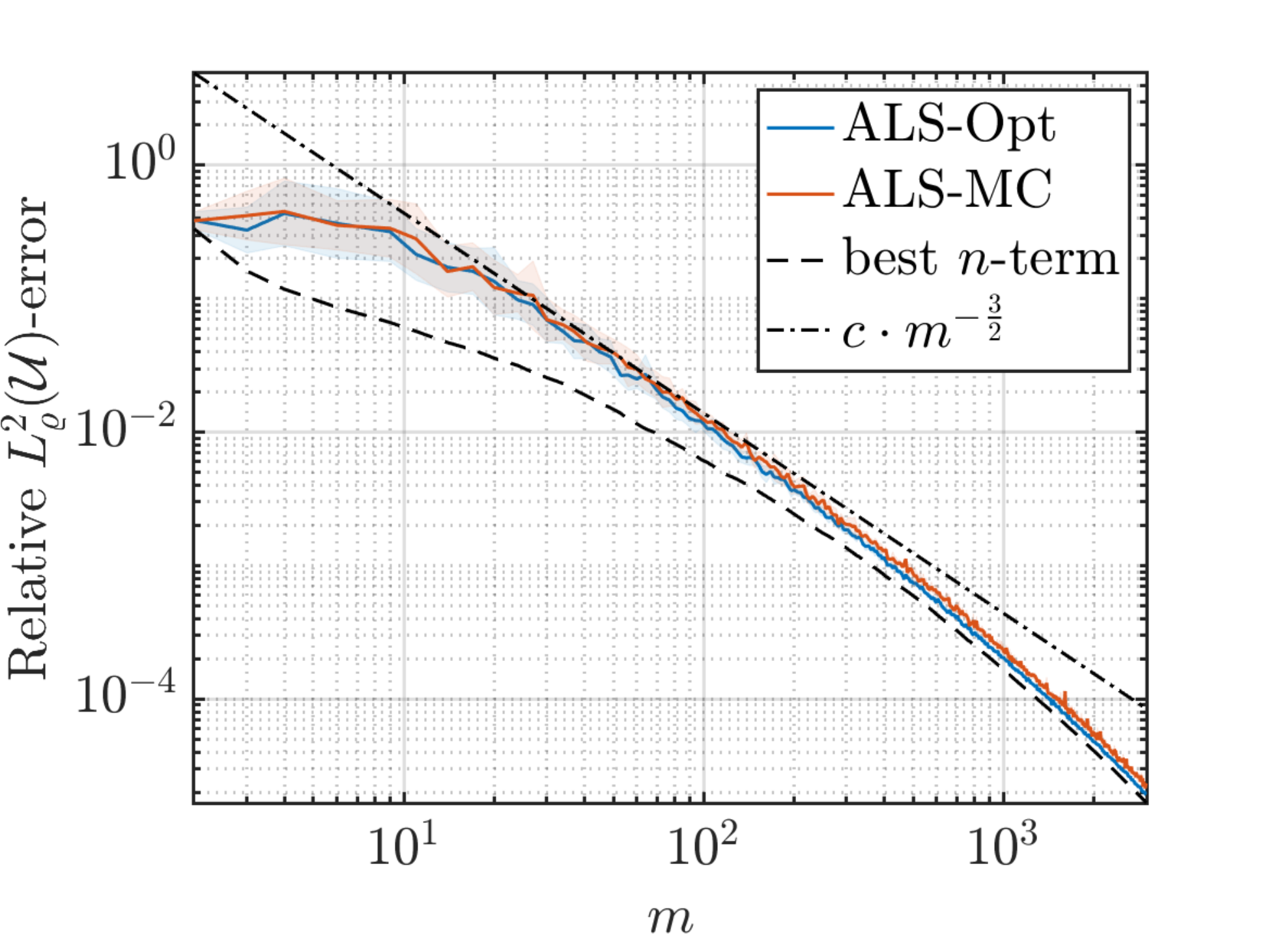}
&
\includegraphics[width = \errplotimg]{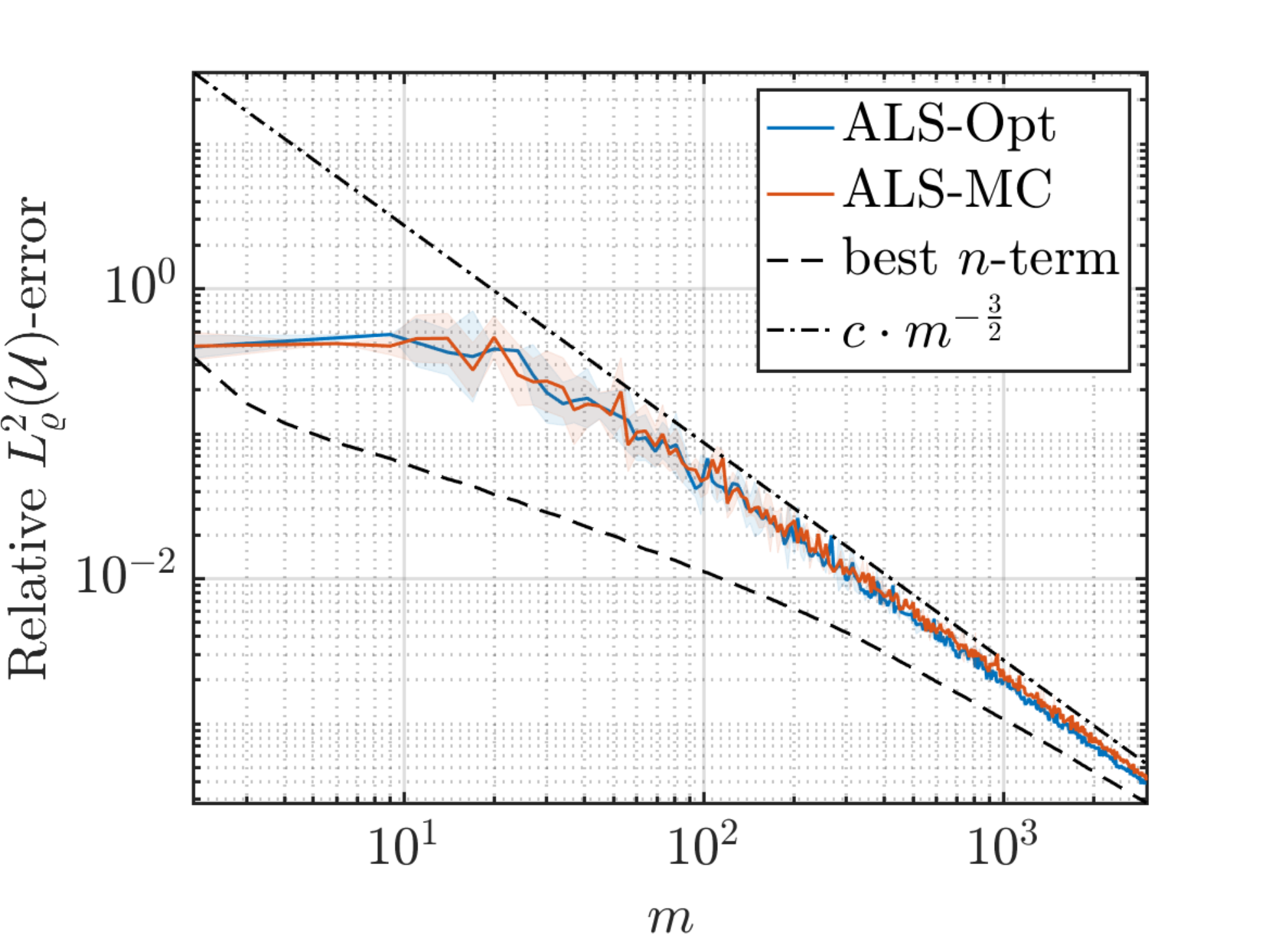}
\\[\errplottextsp]
$d = 2$ & $d = 8$ & $d = 32$
\end{tabular}
\end{small}
\end{center}
\caption{{Comparison between ALS approximation of the function $f = f_3$ with $\delta_i = i^2$ and the theoretical rate $m^{-\frac32}$.} }
\label{fig:fig9}
\end{figure}

\rem{[The function \eqref{f3}]
\label{rem:f3-theory}
{
The function \eqref{f3} has a singularity at $y_i = - 1 - \delta_i$. Hence it is holomorphic in the Bernstein polyellipse $\cE_{\rho_1} \times \cdots \times \cE_{\rho_d}$ for any $\rho_i \geq 1$ satisfying $(\rho_i + \rho_i^{-1})/2 < 1+\delta_i$ (recall the left-hand side is major semi-axis length of the Bernstein ellipse $\cE_{\rho_i}$). Due to Remark \ref{rem:fin-dim}, it is therefore $(\bm{b},\varepsilon)$-holomorphic for any $\bm{b} = (b_i)^{\infty}_{i=1}$ with $b_i > \varepsilon / \delta_i$, $i \in [d]$, and $b_i = 0$ otherwise. 
Theorem \ref{t:best-s-term} now implies the dimension-independent algebraic convergence rate of $n^{\frac12-\frac1p}$ whenever the sequence $(1/\delta_i)_{i \in \bbN} \in \ell^p(\bbN)$. In particular, when $\delta_i = i^2$ as in Fig.\ \ref{fig:fig4}, this guarantees algebraic convergence with rate arbitrarily close to $n^{-\frac32}$ in any dimensions. Fig.\ \ref{fig:fig9} compares this theoretical rate to the practical performance of ALS. In low dimensions, as noted in the previous remark, the algebraic rate is a poor predictor of the method's behaviour. However, in higher ($d = 32$) dimensions, we see a very close agreement between the theoretical rate and empirical performance.
}
}

{ In the next result,} we show that the rate $n^{\frac12-\frac1p}$ appearing in Theorem \ref{t:best-s-term} is sharp. To the best of our knowledge, this result is new. 
This result involves the so-called \textit{weak-$\ell^p$} space. Let $\bm{c} = (c_i)^{\infty}_{i=1}$ be a sequence and write $\bm{c}^* = (c^*_i)^{\infty}_{i=1}$ for a nonincreasing rearrangement of the absolute value sequence $(|c_i|)^{\infty}_{i=1}$. Then, for $p > 0$ the \textit{weak $\ell^p$-quasinorm} is defined as $\nm{\bm{c}}_{p,\infty} = \sup_{i=1,\ldots,\infty} i^{1/p} c^*_i$.
The weak $\ell^p$-space $w \ell^p(\bbN)$ consists of all sequences for which $\nm{\bm{c}}_{p,\infty} < \infty$ \cite{devore1998nonlinear}. 

\thm{
[Sharpness of the rate $n^{\frac12-\frac1p}$]
\label{t:alg-rate-sharp}
Let $\bm{b} \in [0,\infty)^{\bbN}$ and $\varepsilon > 0$. Suppose that $\bm{b} \in \ell^p(\bbN)$ for some $0 < p \leq 1$ but $\bm{b} \notin w \ell^r(\bbN)$ for some $0 < r < p$. Then there are infinitely many non-linearly dependent functions $f \in \cH(\bm{b},\varepsilon)$ such that
\bes{
\limsup_{n \rightarrow \infty} \frac{\nmu{f - f_n}_{L^2_{\varrho}(\cU)}}{n^{\frac12-\frac1r}} = +\infty.
}
}
{Note that there exist sequences with $\bm{b} \in \ell^p(\bbN)$ but $\bm{b} \notin w \ell^r(\bbN)$ for \textit{any} $0 < r < p$. For example, the sequence $b_i = (i \log^2(i))^{-1/p}$ has this property.}

\section{Near-best least-squares polynomial approximation in infinite dimensions with Monte Carlo sampling}\label{s:main-res}

We now present our first main theoretical result on polynomial approximation from MC samples. This shows that, in infinite dimensions, MC sampling is \textit{near best}, in the sense that there is a LS procedure that achieves the same algebraic rates as the best $n$-term approximation, up to a log factor.

\thm{
[MC sampling is near-best for known $\bm{b}$ and $\varepsilon$]
\label{t:main-res}
Let $0 < \epsilon < 1$, $\varrho$ be either the uniform or Chebyshev measure on $\cU = [-1,1]^{\bbN}$, {$\bm{b} \in [0,\infty)^{\bbN}$ with $\bm{b} \in \ell^p(\bbN)$ for some $0 < p <1$}, $\varepsilon > 0$, $m \geq 3$ and $\bm{y}_1,\ldots,\bm{y}_m \sim_{\mathrm{i.i.d.}} \varrho$. Then there exists a set $S \subset \cF$ (depending on $\bm{b}$ and $\varepsilon$ only) of cardinality $|S| \leq \left \lceil m / \log(m/\epsilon) \right \rceil$
such that the following holds with probability at least $1-\epsilon$ for each fixed $f \in \cH(\bm{b},\varepsilon)$. For any $\bm{e} \in \bbC^m$, the LS approximation
\bes{
\hat{f} = \argmin{p \in \cP_S} {\frac1m}\sum^{m}_{i=1} |f(\bm{y}_i) + e_i - p(\bm{y}_i)|^2
}
is unique and satisfies
\bes{
\nmu{f - \hat{f}}_{L^2_{\varrho}(\cU)} \leq C(\bm{b},\varepsilon,p) \cdot \left ( m / \log(m/\epsilon) \right )^{\frac12-\frac1p} +  2 \cdot  \nmu{\bm{e}}_{\infty} .
}
Moreover, the condition number of the LS matrix $\bm{A}$ satisfies $\mathrm{cond}(\bm{A}) \leq 2$.
}

We give the proof of this result in \S \ref{s:proofs}. Its main ingredient is showing the existence of a set $S$ of size $|S| \leq  k : = \frac{m}{\log(m/\epsilon)}$ for which: (i) the error $\nmu{f - f_S}_{L^2_{\varrho}(\cU)}$ is $\ordu{k^{\frac12-\frac1p}}$, and (ii) the constant $\kappa(\cP_{S})$ in \eqref{MC-samp-comp-general} is at most $\ord{k}$.
Having done this, the result follows from Theorem \ref{t:wLS-samp-comp} with $p = f_S$ in \eqref{wLS-err-bd}, albeit after taking some additional care to bound discrete norm $\nmu{f - p}_{\mathsf{disc},1}$.

{Note that the error bound of Theorem~\ref{t:main-res} (also Theorem \ref{t:main-res-2} later) depends on $\|\bm{e}\|_\infty$, which allows for bounded, adversarial noise corrupting the samples. The noise vector $\bm{e}$ might be either deterministic or random. In the random case, it could be, e.g., of the form $e_i = g(\bm{y}_i)$ for some function $g$ for which $\|g\|_{L^{\infty}(\mathcal{U})}$ is finite. Note that it is possible to deal with unbounded random noise models using tools from \cite{migliorati2015convergence}. }

{
\rem{[The function \eqref{f3}]
Recall Remark \ref{rem:f3-theory} and the function $f_3$ with $\delta_i = i^2$. Theorem \ref{t:main-res} implies that there exists a least-squares approximation based on MC sampling that achieves the rate $m^{-\frac32+\gamma}$ for arbitrarily small $\gamma > 0$. It does not show that the ALS scheme achieves this rate: as discussed, there is no theoretical guarantee that ALS selects `good' index sets in practice. However, Fig.\ \ref{fig:fig9} suggests that ALS does appear to achieve the theoretical rate, at least in the case of this function.
}
}

\rem{
In \S \ref{ss:low-dim-fns} we considered a function of one variable $f(\bm{y}) = g(y_1)$, wherein $\kappa(\cP_{S^*}) = n^2$, where $S^*$, $|S^*| = n$, is the index set corresponding to the best $n$-term approximation. See \eqref{Sstar-oned-fun}. Such pathological examples do not contradict Theorem \ref{t:main-res}, even though stable polynomial approximation in $\cP_{S^*}$ necessitates a quadratic scaling of $m$ with $n$ in such cases. The reason is that the best $n$-term approximation error decays exponentially fast in $n$ for such functions (Remark \ref{rem:fin-inf-diff}), and therefore root exponentially fast in $m$. But this is still faster than the algebraic convergence rate asserted in Theorem \ref{t:main-res}. A similar argument applies for functions of $d$-variables, since in such cases, the best $n$-approximation error decays exponentially-fast in $n^{1/d}$. 
In general, it is the presence of large multi-indices in $S$ that cause $\kappa(\cP_S)$ to scale quadratically in $n = |S|$. Theorem \ref{t:main-res} essentially says that algebraic rates of convergence can be obtained without using multi-indices that are too large. 
}

\subsection{Weighted $k$-term approximation}\label{ss:weighted-k-term}

As a way to motivate the approach considered in the next section, we now elaborate on the construction employed in Theorem \ref{t:main-res}. This is based on \textit{weighted $(k,\bm{u})$-term approximation} \cite{rauhut2016interpolation}. Let
\be{
\label{weights-def}
\bm{u} = (u_{\bm{\nu}})_{\bm{\nu} \in \cF},\qquad \text{where }u_{\bm{\nu}} = \nmu{\Psi_{\bm{\nu}}}_{L^{\infty}(\cU)},\ \forall \bm{\nu} \in \cF,
}
and observe that $u_{\bm{\nu}} \geq \nmu{\Psi_{\bm{\nu}}}_{L^{2}_{\varrho}(\cU)} = 1$ since $\varrho$ is a probability measure and the $\Psi_{\bm{\nu}}$ are orthonormal. Next, we define the weighted cardinality of a set $S \subset \cF$ as $|S|_{\bm{u}} = \sum_{\bm{\nu} \in S} u^2_{\bm{\nu}}$.
Then, for $k \geq 0$ (note that $k$ need not be an integer in this case) we define the \textit{weighted best $(k,\bm{u})$-term approximation} to $f$ as
\bes{
f_{k,\bm{u}} = f_{S^*},\qquad \text{where }S^* \in \argmin{} \{ \nmu{f - f_S}_{L^2_{\varrho}(\cU)} : S \subset \cF,\ |S|_{\bm{u}} \leq k \}.
}
In the proof of Theorem \ref{t:main-res}, we exploit the key fact that
\be{
\label{f-weighted-err}
\nmu{f - f_{k,\bm{u}}}_{L^2_{\varrho}(\cU)} \leq C(\bm{b},\varepsilon,p) \cdot k^{\frac12-\frac1p},\quad \forall f \in \cH(\bm{b},\varepsilon).
}
Thus, the set $S^* = S$ {provides} a suitable choice in terms of (i) above. However, it transpires that it is also suitable in the terms of (ii). Indeed, consider any set $S \subseteq \cF$ with weighted cardinality $|S|_{\bm{u}} \leq k$. Then \eqref{MC-samp-comp-general} states that the sample complexity of LS in the subspace $\cP_{S}$ with MC sampling is determined by the constant
\bes{
\kappa(\cP_S) = \nmu{\cK(\cP_S)}_{L^{\infty}(\cU)} = \sup_{\bm{y} \in \cU} \sum_{\bm{\nu} \in S} | \Psi_{\bm{\nu}}(\bm{y}) |^2.
}
However, Chebyshev and Legendre polynomials attain the maximum absolute value at the point $\bm{y} = \bm{1}$. Hence 
\be{
\label{kappa-weighted-bound}
\kappa(\cP_S) = \sum_{\bm{\nu} \in S}| \Psi_{\bm{\nu}}(\bm{1}) |^2 = \sum_{\bm{\nu} \in S} \nmu{\Psi_{\bm{\nu}}}_{L^{\infty}(\cU)} = |S|_{\bm{u}} \leq k.
}
This is the essence of the proof. Sets of fixed weighted cardinality a `good' candidates for MC sampling due to \eqref{kappa-weighted-bound}, but also good for approximating functions, due to \eqref{f-weighted-err}.

\section{Near-best polynomial approximation in infinite dimensions via compressed sensing}\label{s:polyappCS}

Theorem \ref{t:main-res} asserts the existence of a set $S$ for which the corresponding LS approximation from MC samples converges with the desired algebraic rate. It says nothing about how to construct such a set in practice. At the very least, this set depends on the anisotropy parameters $\bm{b}$ and $\varepsilon$, which are generally unknown.

We now show how it is possible to compute polynomial approximations which attain the same rates of convergence without any knowledge of these parameters, subject to a slightly stronger assumption on $\bm{b}$. 
Furthermore, we also demonstrate that these lead to practical methods, with performance comparable or sometimes better than the ALS approximation studied previously.

\subsection{Polynomial approximation via compressed sensing}

{Originating in the works \cite{doostan2011nonadapted,mathelin2012compressed,rauhut2012sparse},} polynomial approximation via \textit{Compressed Sensing (CS)} is now well established.
In this section we primarily follow the {layout of \cite[Chpt.\ 7]{adcock2022sparse}, which is based on ideas originating in \cite{rauhut2016interpolation}.}

In \S \ref{ss:weighted-k-term}, we observed that any function in $\cH(\bm{b},\varepsilon)$ can be approximated with algebraically-decay error using its weighted best $(k,\bm{u})$-term approximation. However, the set $S^*$ that yields this approximation is unknown. The idea now is to use CS techniques to promote such approximate \textit{weighted sparsity} by using a weighted $\ell^1$-norm penalty term based on the weights \eqref{weights-def}.

Before doing this, however, we first need to truncate the infinite expansion \eqref{f-exp}. Let $\Lambda \subset \cF$ be a finite multi-index set and write $f_{\Lambda} = \sum_{\bm{\nu} \in \Lambda} c_{\bm{\nu}} \Psi_{\bm{\nu}}$ for the truncated expansion of $f$. Then, given samples points $\bm{y}_1,\ldots,\bm{y}_m$, we have
\bes{
f(\bm{y}_i) + e_i = f_{\Lambda}(\bm{y}_i) + (f-f_{\Lambda})(\bm{y}_i) +e_i = \sum_{\bm{\nu} \in \Lambda} c_{\bm{\nu}} \Psi_{\bm{\nu}}(\bm{y}_i) + (f-f_{\Lambda})(\bm{y}_i) + e_i.
}
Now let $\bm{\nu}_1,\ldots,\bm{\nu}_N$ be an enumeration of $\Lambda$ and $\bm{c}_{\Lambda} = (c_{\bm{\nu}_i})^{N}_{i=1}$. Then we have
\be{
\label{approx-LS-CS}
\bm{f} = \left ( (f(\bm{y}_i) + e_i )/\sqrt{m} \right )^{m}_{i=1}= \bm{A} \bm{c}_{\Lambda} + \bm{n},
}
where
\be{
\label{Ab-CS-def}
\bm{A} = \left (  \Psi_{\bm{\nu}_j}(\bm{y}_i)/\sqrt{m} \right )^{m,N}_{i,j=1} \in \bbC^{m \times N},\quad 
\bm{n} = \left ( ((f-f_{\Lambda})(\bm{y}_i) + e_i)/\sqrt{m} \right )^{m}_{i=1} \in \bbC^m.
}
The idea now is to search for (approximate) solutions of \eqref{approx-LS-CS} that have small $\ell^1_{\bm{u}}$-norm. {W}e do this via the following \textit{weighted square-root LASSO} program: 
\bes{
\hat{f} = \sum_{\bm{\nu} \in \Lambda} \hat{c}_{\bm{\nu}} \Psi_{\bm{\nu}},\quad \text{where } \hat{\bm{c}} = (\hat{c}_{\bm{\nu}})_{\bm{\nu} \in \Lambda} \in \argmin{\bm{z} \in \bbC^{N}} \lambda \nmu{\bm{z}}_{1,\bm{u}} + \nmu{\bm{A} \bm{z} - \bm{f}}_2.
}
Here $\nmu{\bm{z}}_{1,\bm{u}} = \sum_{\bm{\nu} \in \Lambda} u_{\bm{\nu}} |z_{\bm{\nu}}| $ is the $\ell^1_{\bm{u}}$-norm and $\lambda > 0$ is a parameter.

\subsection{Truncation error and anchored sets}

It remains to choose the set $\Lambda$. This set must be sufficiently large so that the truncation error $f - f_{\Lambda}$ is no larger than the approximation error \eqref{f-weighted-err}. However, $\Lambda$ must also be a finite set. 
To ensure this property, we now introduce the concept of anchored sets. A multi-index set $S \subset \cF$ is \textit{anchored} if it is lower and if, for every $j \in \bbN$,
\bes{
\bm{e}_j \in S\ \Rightarrow\ \{ \bm{e}_1,\ldots,\bm{e}_j \} \subseteq S.
}
Here $\bm{e}_j$ denotes the sequence with one in its $j$th entry and zero elsewhere. We now need the following result (see, e.g., \cite[Thm.\ 3.33]{adcock2022sparse} or \cite[\S 3.8]{cohen2015approximation}), which states that near-best $n$-term 
can be attained in anchored sets, for certain $\bm{b}$. 

\thm{
[Algebraic convergence in anchored sets]
\label{t:best-s-term-anchored}
Let $\varepsilon > 0$ and $\bm{b} \in [0,\infty)^{\bbN}$ be monotonically nonincreasing and such that $\bm{b} \in \ell^p(\bbN)$ for some $0 < p < 1$. Then, for every $n \in \bbN$, there exists an anchored set $S \subset \cF$, $|S| = n$, such that
\bes{
\nmu{f - f_S}_{L^2_{\varrho}(\cU)} \leq C(\bm{b},\varepsilon,p) \cdot n^{\frac12-\frac1p},\qquad \forall f \in \cH(\bm{b},\varepsilon),\ n \in \bbN.
}
}

With this in hand, we now construct $\Lambda$ so that it contains all anchored sets of a given size $n$. This turns out to be a finite set (see, e.g., {\cite{cohen2017discrete} or }\cite[Prop.\ 2.18]{adcock2022sparse}):
\bes{
\bigcup \{ S \subset \cF : |S| \leq n,\text{ $S$ anchored} \} \subset \Lambda^{\mathsf{HCI}}_{n},
}
where
\be{
\label{HCI-def}
\Lambda^{\mathsf{HCI}}_{n} = \left \{ \bm{\nu} = (\nu_k)^{\infty}_{k=1} \in \cF : \prod^{n-1}_{k=1} (\nu_k+1) \leq n,\ \nu_k = 0,\ \forall k \geq n \right \},
}
is isomorphic to the $(n-1)$-dimensional hyperbolic cross index set of order $n-1$.

\subsection{\boldmath Near-best polynomial approximation for unknown $b$ and $\varepsilon$}

We now present the main result of this section. For this, we need one additional concept. Let $\bm{b} = (b_i)_{i \in \bbN}$ be a sequence. We define its \textit{minimal monotone majorant} as the sequence $\tilde{\bm{b}} = (\tilde{b}_i)_{i \in \bbN}$, where $\tilde{b}_i = \sup_{j \geq i} {| b_j|}$, $\forall i \in \bbN$.
Then, given $0 < p < \infty$, we define the \textit{monotone $\ell^p$} space $\ell^p_{\mathsf{M}}(\bbN)$ as
$
\ell^p_{\mathsf{M}}(\bbN) = \{ \bm{b} \in \ell^{\infty}(\bbN) : \nmu{\bm{b}}_{p,\mathsf{M}} : = \nmu{\tilde{\bm{b}}}_{p} < \infty  \}.
$

\thm{
[MC sampling is near-best for unknown $\bm{b}$ and $\varepsilon$]
\label{t:main-res-2}
Let $0 < \epsilon < 1$, $\varrho$ be either the uniform or Chebyshev measure on $\cU = [-1,1]^{\bbN}$, $m \geq 3$ and $\bm{y}_1,\ldots,\bm{y}_m \sim_{\mathrm{i.i.d.}} \varrho$. Let $\Lambda = \Lambda^{\mathsf{HCI}}_{n}$, where
\bes{
n = \left \lceil m / L(m,\epsilon) \right \rceil,\qquad L(m,\epsilon) : = \log(m) \cdot (\log^3(m) + \log(\epsilon^{-1})),
}
$\varepsilon > 0$ and $\bm{b} \in [0,\infty)^{\bbN}$ with $\bm{b} \in \ell^p_{\mathsf{M}}(\bbN)$ for some $0 < p < 1$. Then the following holds with probability at least $1-\epsilon$ for each fixed $f \in \cH(\bm{b},\varepsilon)$. For any $\bm{e} \in \bbC^m$, every
\be{
\label{SRLASSO-minimizer}
\hat{\bm{c}} = (\hat{c}_{\bm{\nu}})_{\bm{\nu} \in \Lambda} \in \argmin{\bm{z} \in \bbC^{N}} \lambda \nmu{\bm{z}}_{1,\bm{u}} + \nmu{\bm{A} \bm{z} - \bm{f}}_2,
}
where $\bm{A}$ and $\bm{f}$ 
are as in \eqref{approx-LS-CS}--\eqref{Ab-CS-def} and $\lambda = (4 \sqrt{m/L(m,\epsilon)})^{-1}$, yields an approximation $\hat{f} = \sum_{\bm{\nu} \in \Lambda} \hat{c}_{\bm{\nu}} \Psi_{\bm{\nu}}$ satisfying
\bes{
\nmu{f - \hat{f}}_{L^2_{\varrho}(\cU) } \leq C = C(\bm{b},\varepsilon,p) \cdot \left ( m / L(m,\epsilon) \right )^{\frac12 - \frac1p} + c \cdot \nm{\bm{e}}_{\infty}.
}
Here $c \geq 1$ is a universal constant. Furthermore, given $\hat{\bm{c}}$ it is possible to compute a set $S \subset \Lambda$ of cardinality $|S| \leq \left \lceil m / L(m,\epsilon) \right \rceil$
for which $\hat{f}_S = \sum_{\bm{\nu} \in S} \hat{c}_{\bm{\nu}} \Psi_{\bm{\nu}}$ satisfies
\bes{
\nmu{f - \hat{f}_S}_{L^2_{\varrho}(\cU) } \leq 10 \cdot C(\bm{b},\varepsilon,p) \cdot \left ( m / L(m,\epsilon) \right )^{\frac12 - \frac1p} + 3 \cdot c \cdot \nm{\bm{e}}_{\infty}.
}
}

This result is {an extension of \cite[Thm.\ 3.7]{adcock2022efficient}} -- see \S \ref{s:proofs} for further details. Comparing it with Theorem \ref{t:main-res}, we conclude the following. In the absence of knowledge about the parameters $\bm{b}$ and $\varepsilon$, it is still possible to obtain the same algebraic rates of convergence using MC sampling, up to a larger polylogarithmic factor $L(m,\epsilon)$, subject to the slightly stronger assumption $\bm{b} \in \ell^p_{\mathsf{M}}(\bbN)$. In particular, MC sampling remains near-optimal in infinite dimensions, even when the anisotropy parameters are unknown.  {Note that \eqref{SRLASSO-minimizer} yields a polynomial approximation $\hat{f}$ that generally has many more nonzero terms than the approximation in Theorem \ref{t:main-res}, since it is a polynomial defined over the large index set $\Lambda = \Lambda^{\mathsf{HCI}}_n$. The second part of Theorem \ref{t:main-res-2} states that one can compute a polynomial approximation with a comparable number of terms to that of Theorem \ref{t:main-res} that also achieves the desired rate.} 

Note that $\ell^p_{\mathsf{M}}(\bbN) \subset \ell^p(\bbN)$ and that $\nm{\bm{b}}_{p,\mathsf{M}} = \nm{\bm{b}}_{p}$ whenever $\bm{b}$ is monotonically nonincreasing. Monotonicity of $\bm{b}$ means that the variables are ordered in terms of importance. Thus, the assumption $\ell^p_{\mathsf{M}}(\bbN) $ in effect says that important variables, while not necessarily being ordered, cannot occur at arbitrarily high indices.

We remark in passing that this theorem does not give an algorithm for computing a minimizer \eqref{SRLASSO-minimizer}. However, it has been shown in \cite{adcock2022efficient} that this can be done via efficient iterative algorithms. We use such algorithms in the numerical experiments in the next subsection. See \S \ref{ss:efficient-PDI} and \S \ref{s:conclusions} for further details and discussion.

\subsection{Numerical experiments}\label{ss:CS-ALS-numerical}

We conclude with several numerical experiments demonstrating the practical performance of CS-based polynomial approximation -- see \S \ref{ss:efficient-PDI} for details on the numerical implementation. We do this by comparing it with the ALS approximation scheme considered previously. 

\begin{figure}[t!]
\begin{center}
\begin{small}
 \begin{tabular}{@{\hspace{0pt}}c@{\hspace{\errplotsp}}c@{\hspace{\errplotsp}}c@{\hspace{0pt}}}
\includegraphics[width = \errplotimg]{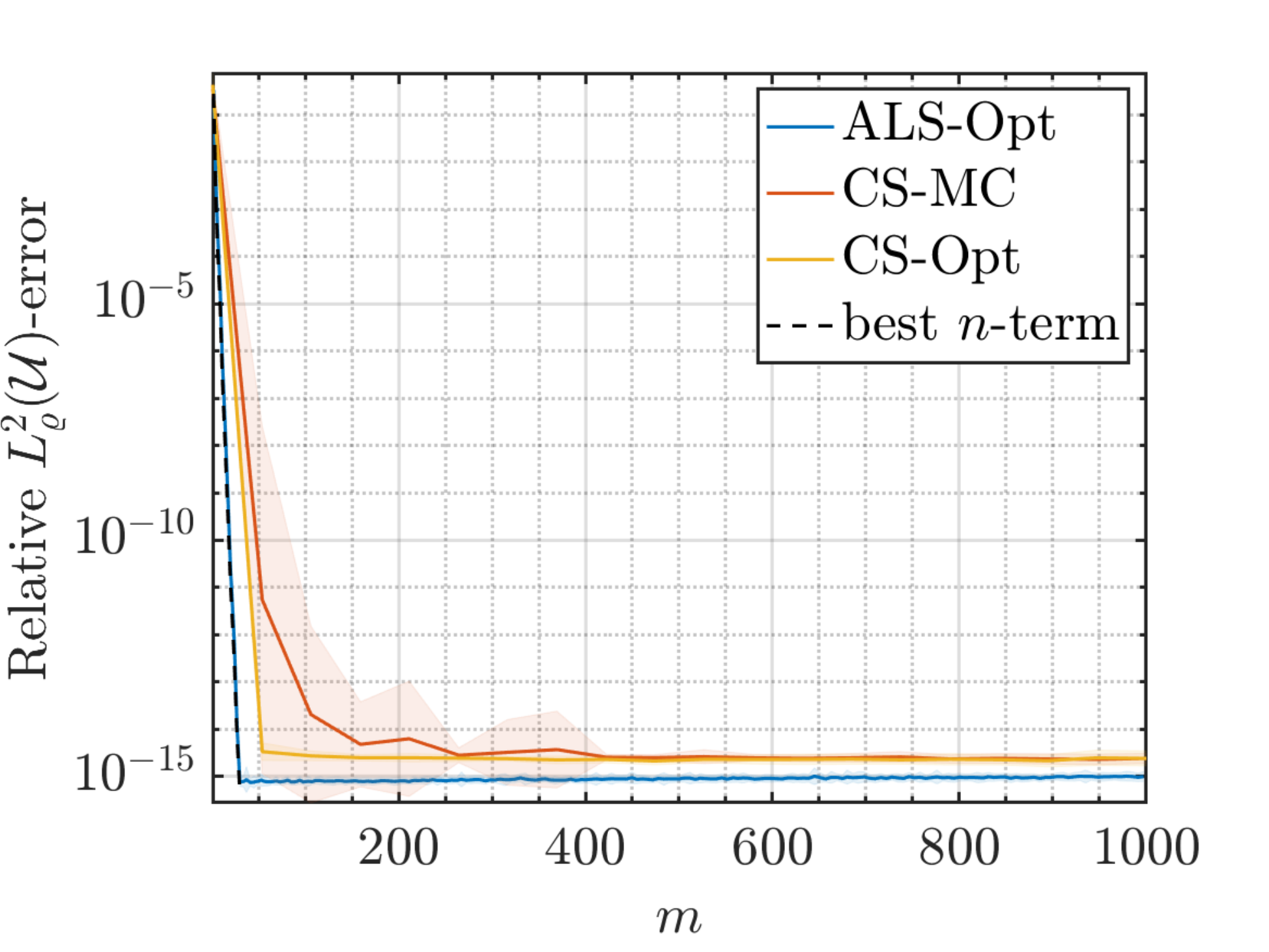}
&
\includegraphics[width = \errplotimg]{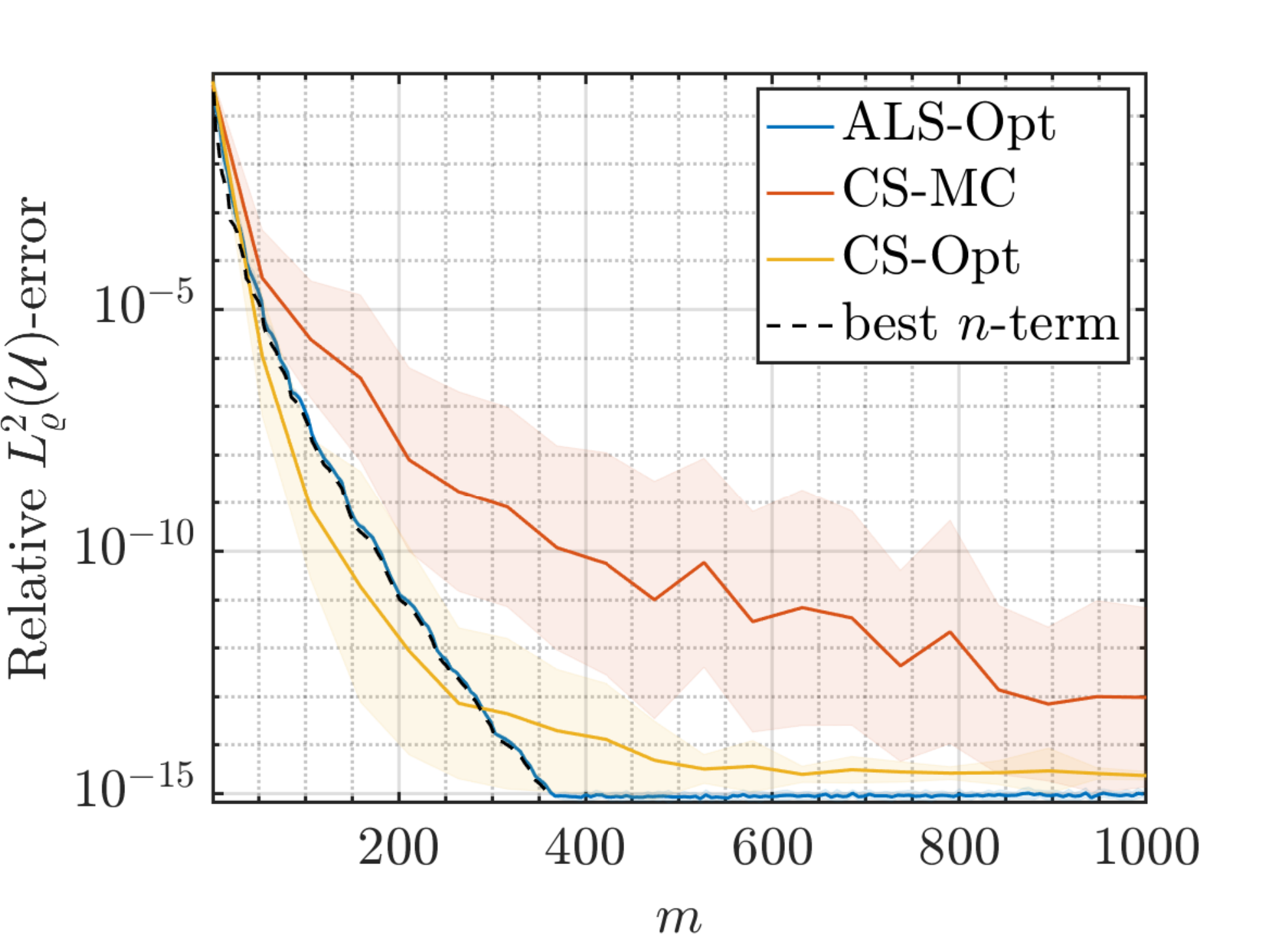}
&
\includegraphics[width = \errplotimg]{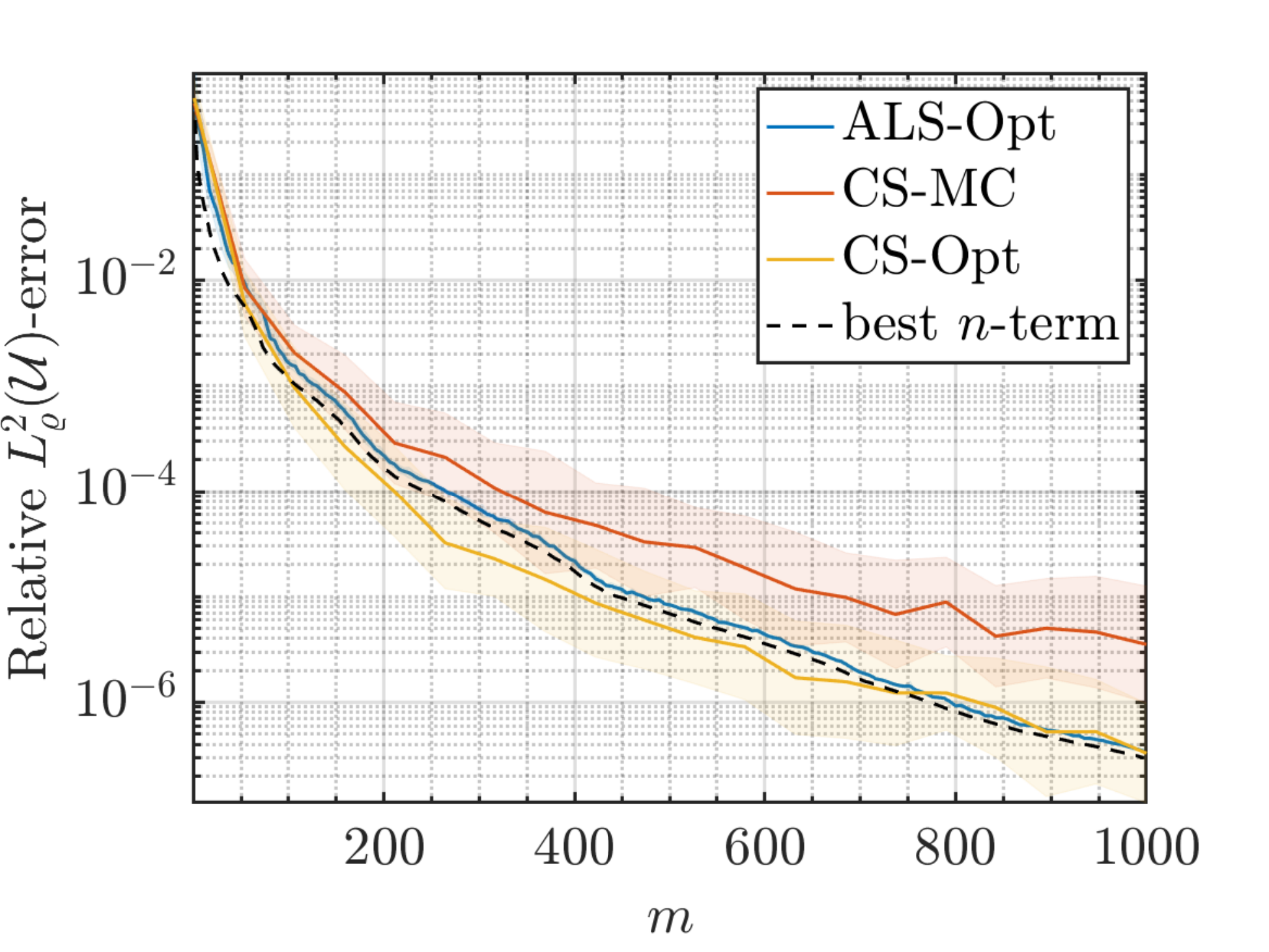}
\\[\errplottextsp]
$d = 1$ & $d = 2$ & $d = 4$
\\[\errplottextsp]
\includegraphics[width = \errplotimg]{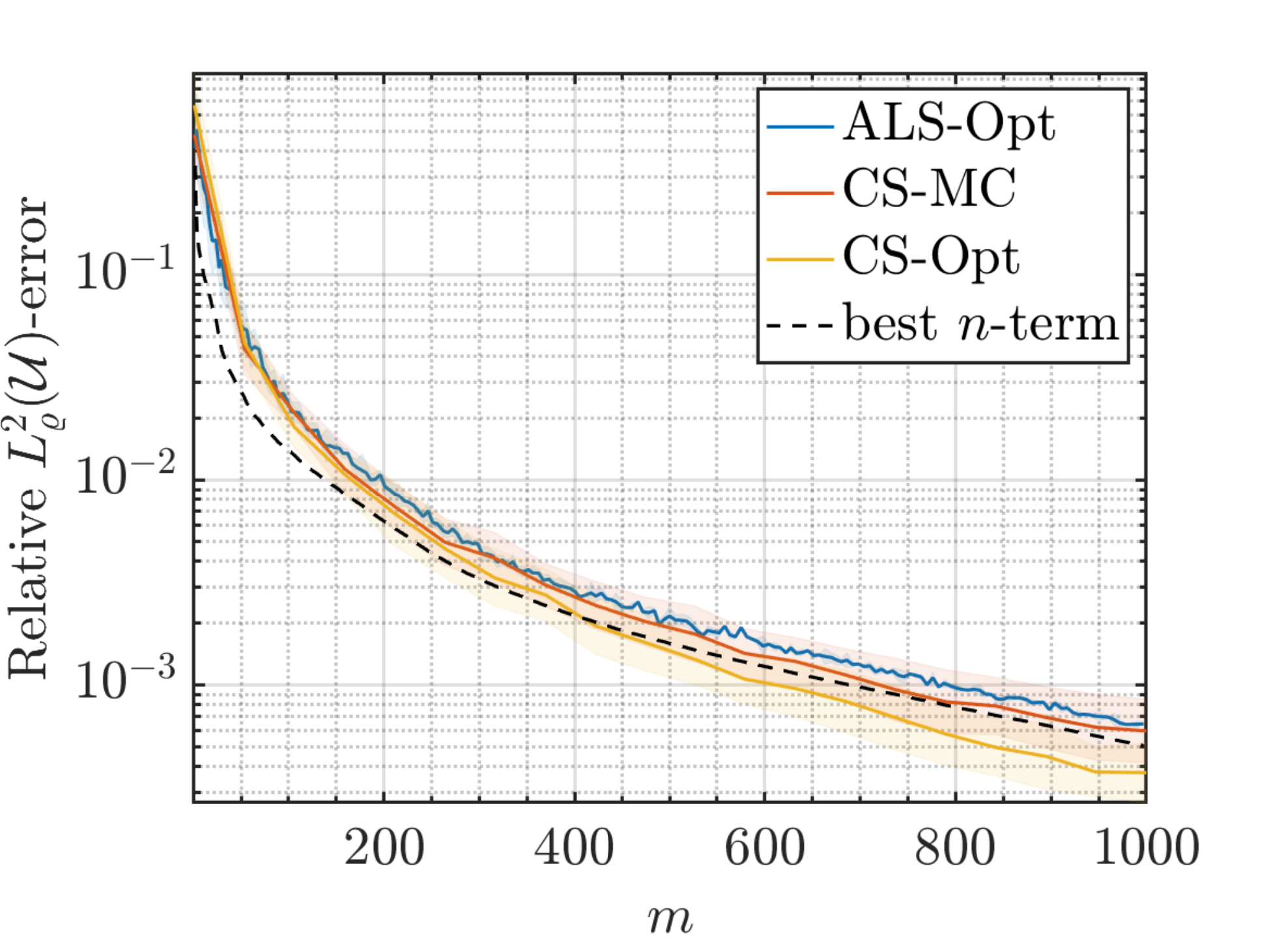}
&
\includegraphics[width = \errplotimg]{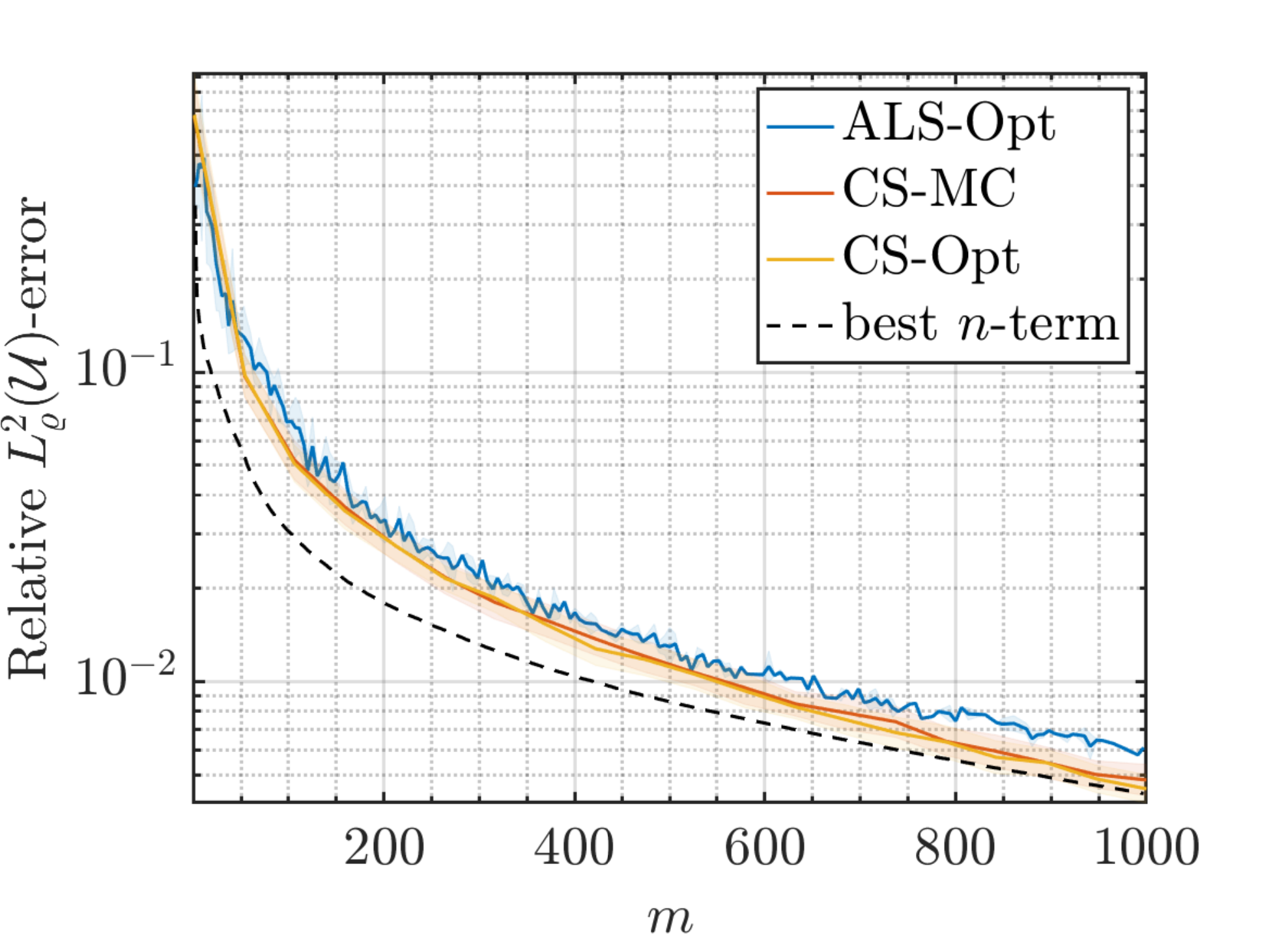}
&
\includegraphics[width = \errplotimg]{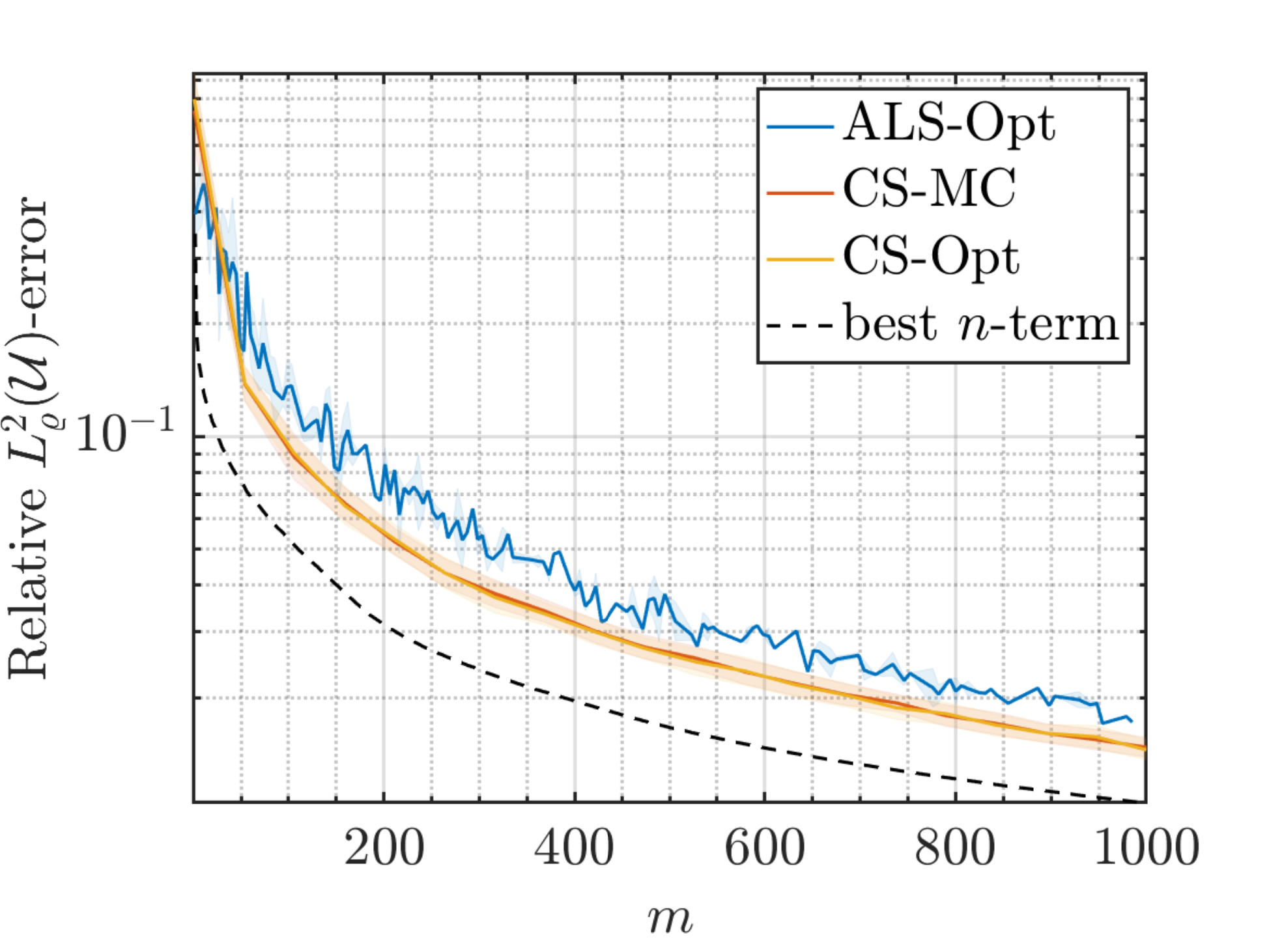}
\\[\errplottextsp]
$d = 8$ & $d = 16$ & $d = 32$
\end{tabular}
\end{small}
\end{center}
\caption{Comparison of ALS with the near-optimal sampling strategy (blue), CS with MC sampling (red) and CS with random sampling from \eqref{CS-meas-opt} (yellow) for the function $f = f_1$ using Legendre polynomials. For ALS we use the scaling \eqref{m-scaling-LS}.} 
\label{fig:fig10}
\end{figure}

\begin{figure}[t!]
\begin{center}
\begin{small}
 \begin{tabular}{@{\hspace{0pt}}c@{\hspace{\errplotsp}}c@{\hspace{\errplotsp}}c@{\hspace{0pt}}}
\includegraphics[width = \errplotimg]{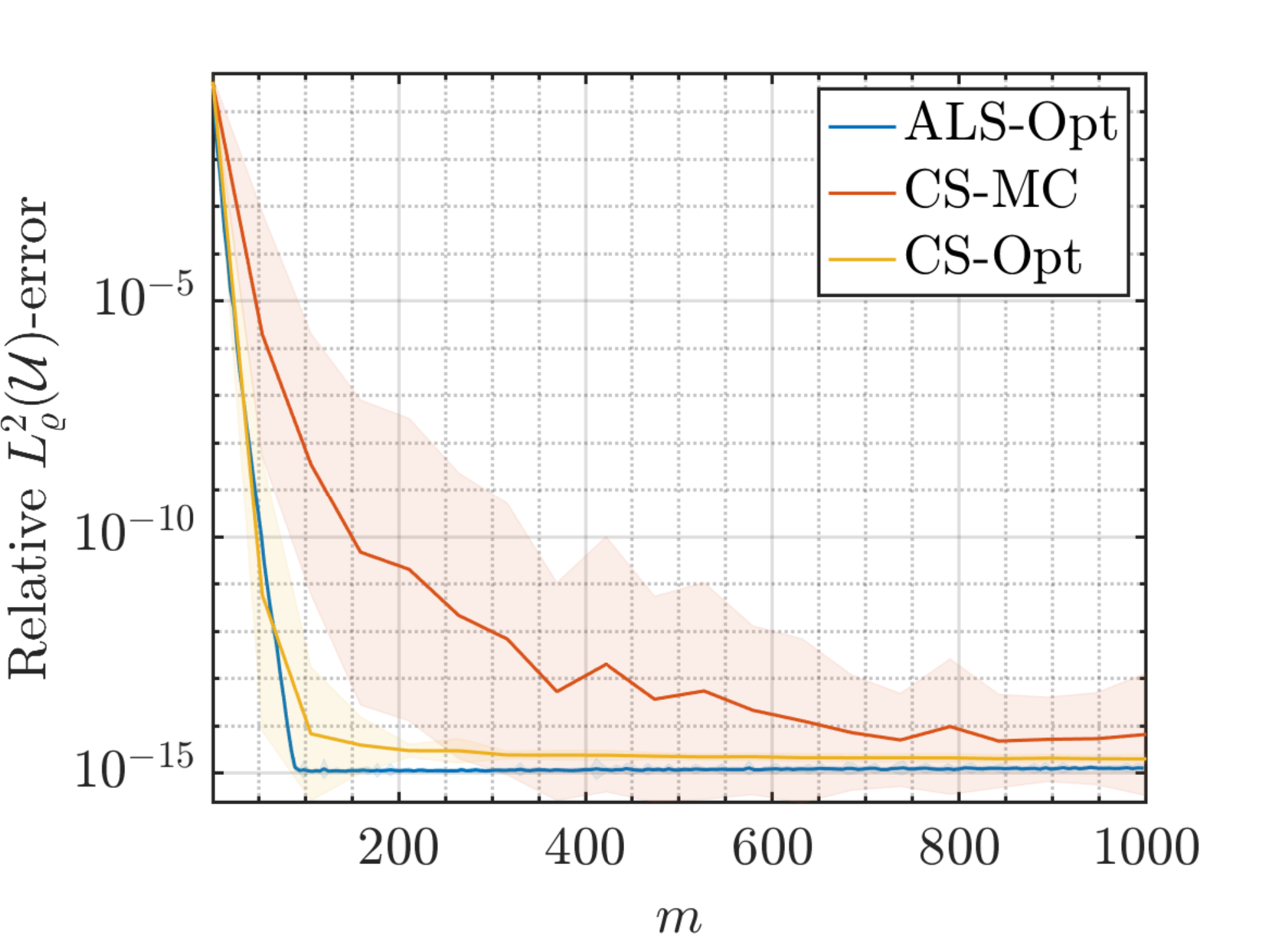}
&
\includegraphics[width = \errplotimg]{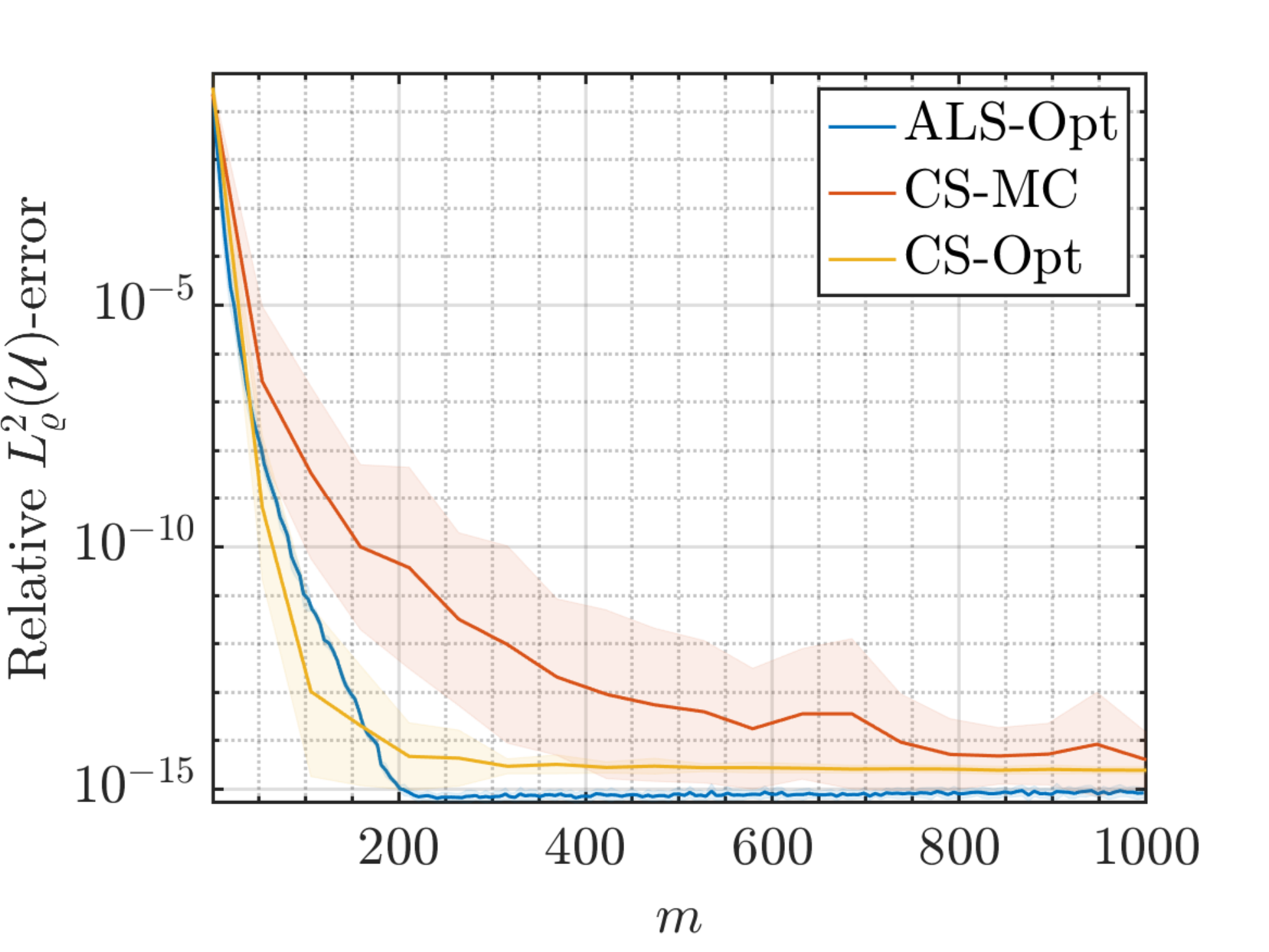}
&
\includegraphics[width = \errplotimg]{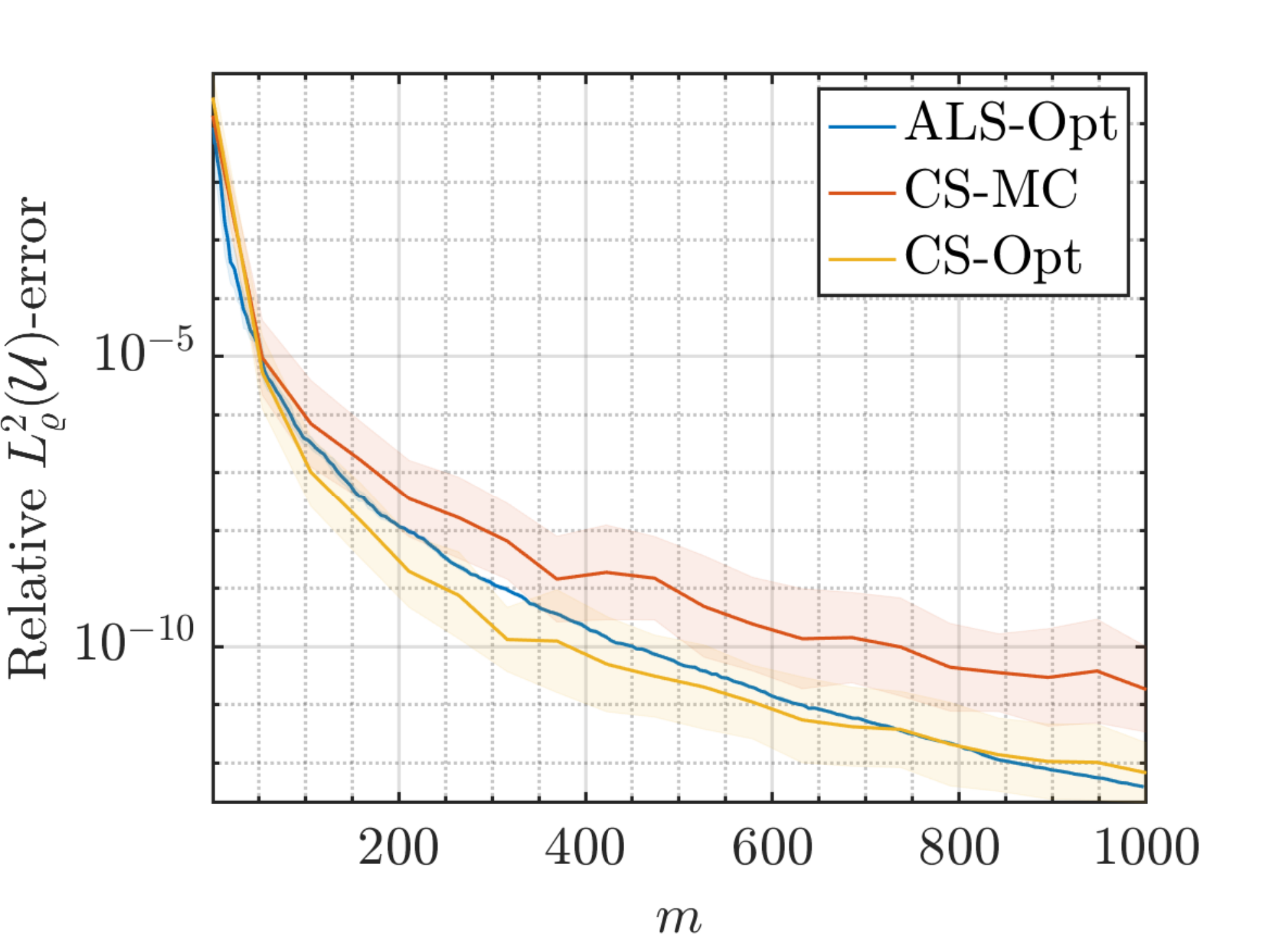}
\\[\errplottextsp]
$d = 1$ & $d = 2$ & $d = 4$
\\[\errplottextsp]
\includegraphics[width = \errplotimg]{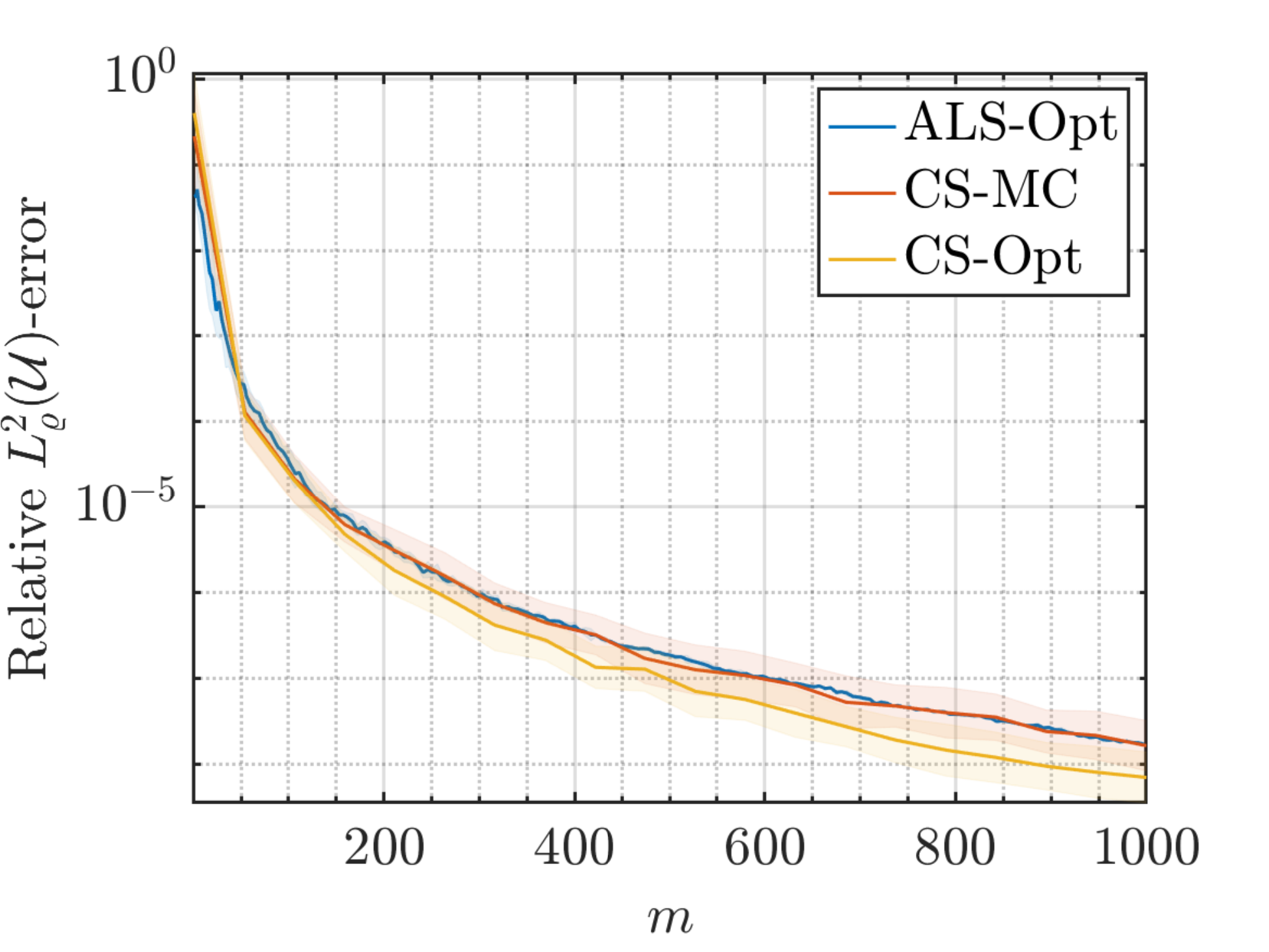}
&
\includegraphics[width = \errplotimg]{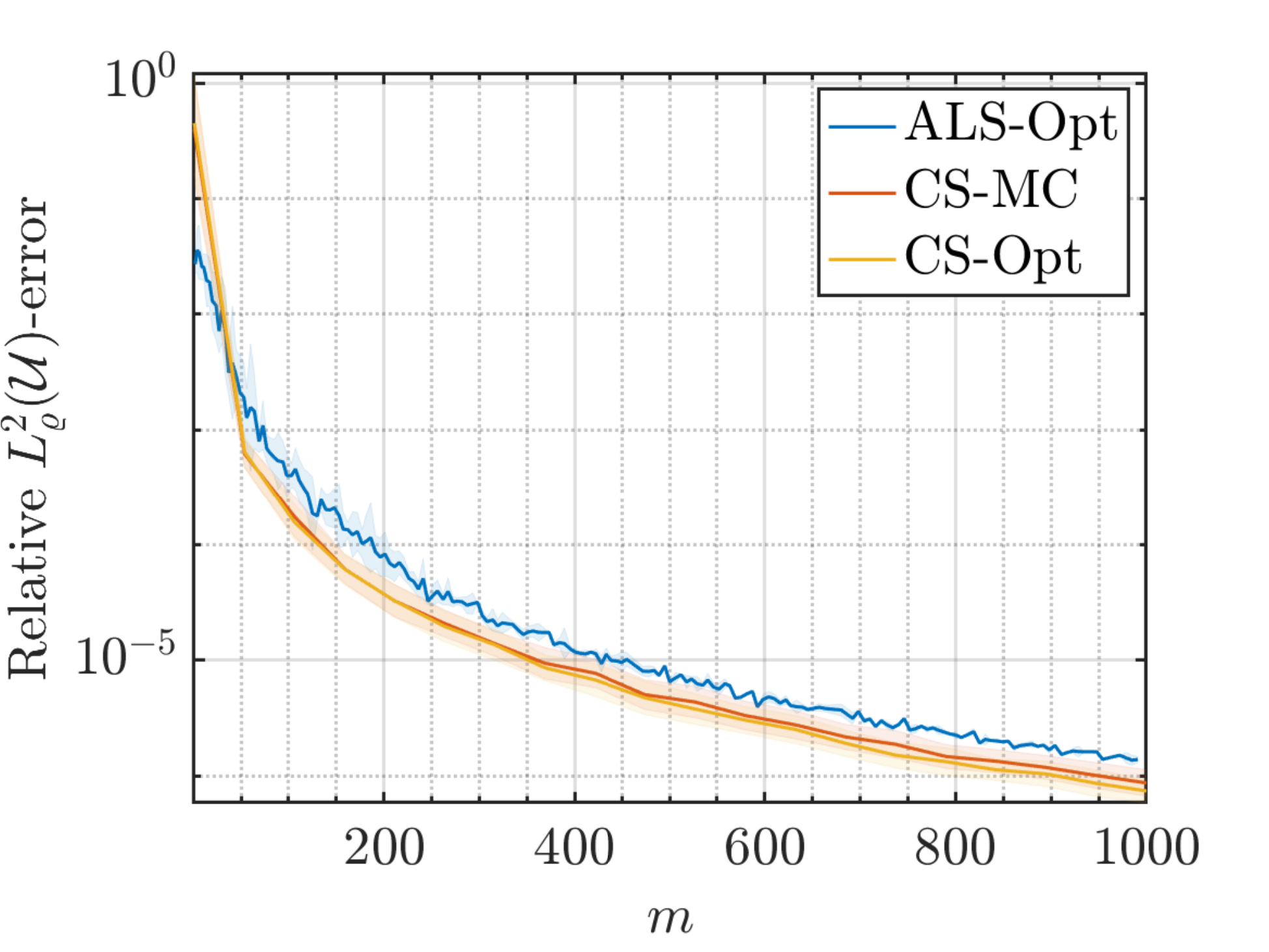}
&
\includegraphics[width = \errplotimg]{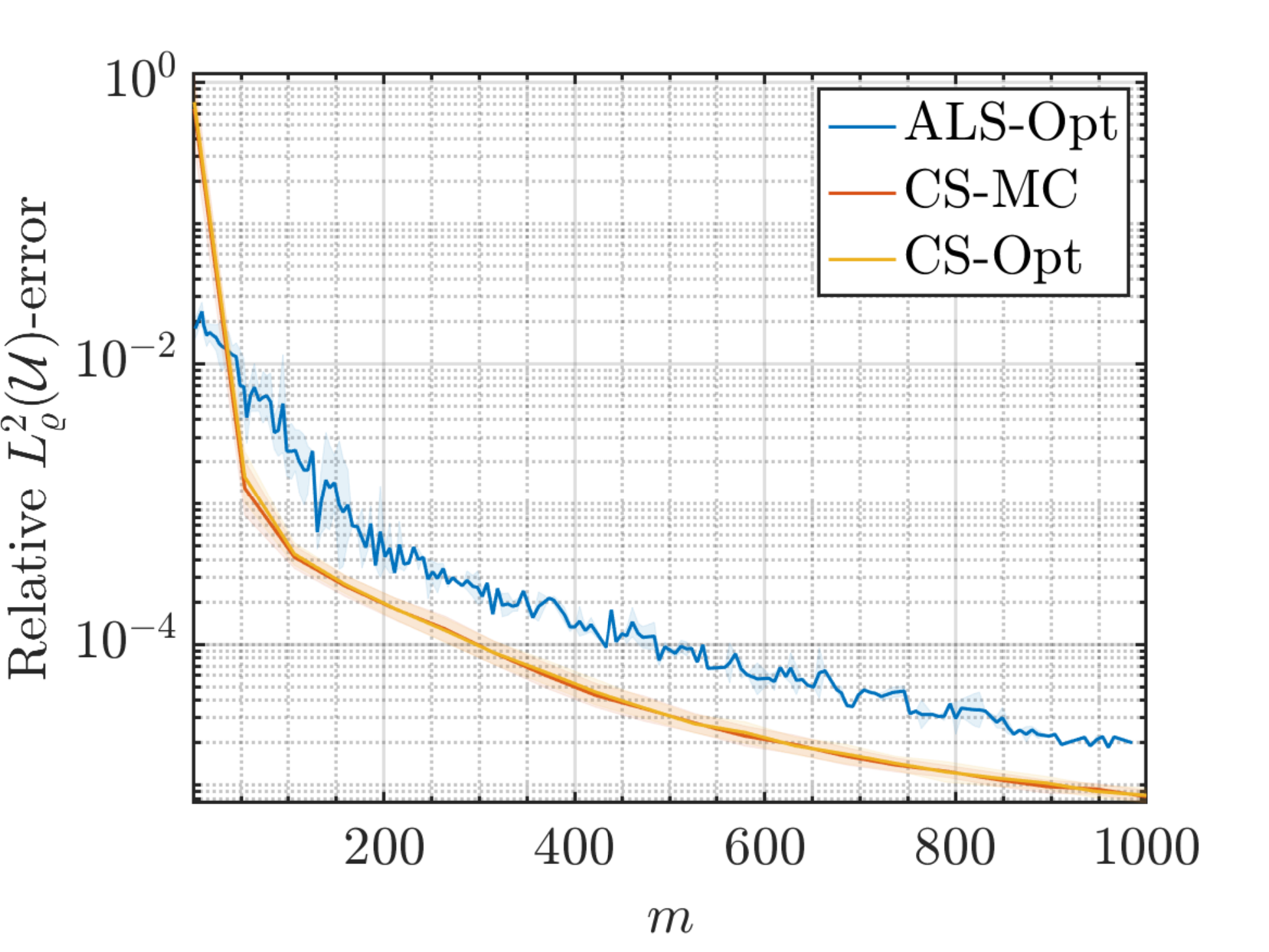}
\\[\errplottextsp]
$d = 8$ & $d = 16$ & $d = 32$
\end{tabular}
\end{small}
\end{center}
\caption{The same as Fig.\ \ref{fig:fig10}, except with $f = f_2$.} 
\label{fig:fig11}
\end{figure}

\begin{figure}[t!]
\begin{center}
\begin{small}
 \begin{tabular}{@{\hspace{0pt}}c@{\hspace{\errplotsp}}c@{\hspace{\errplotsp}}c@{\hspace{0pt}}}
\includegraphics[width = \errplotimg]{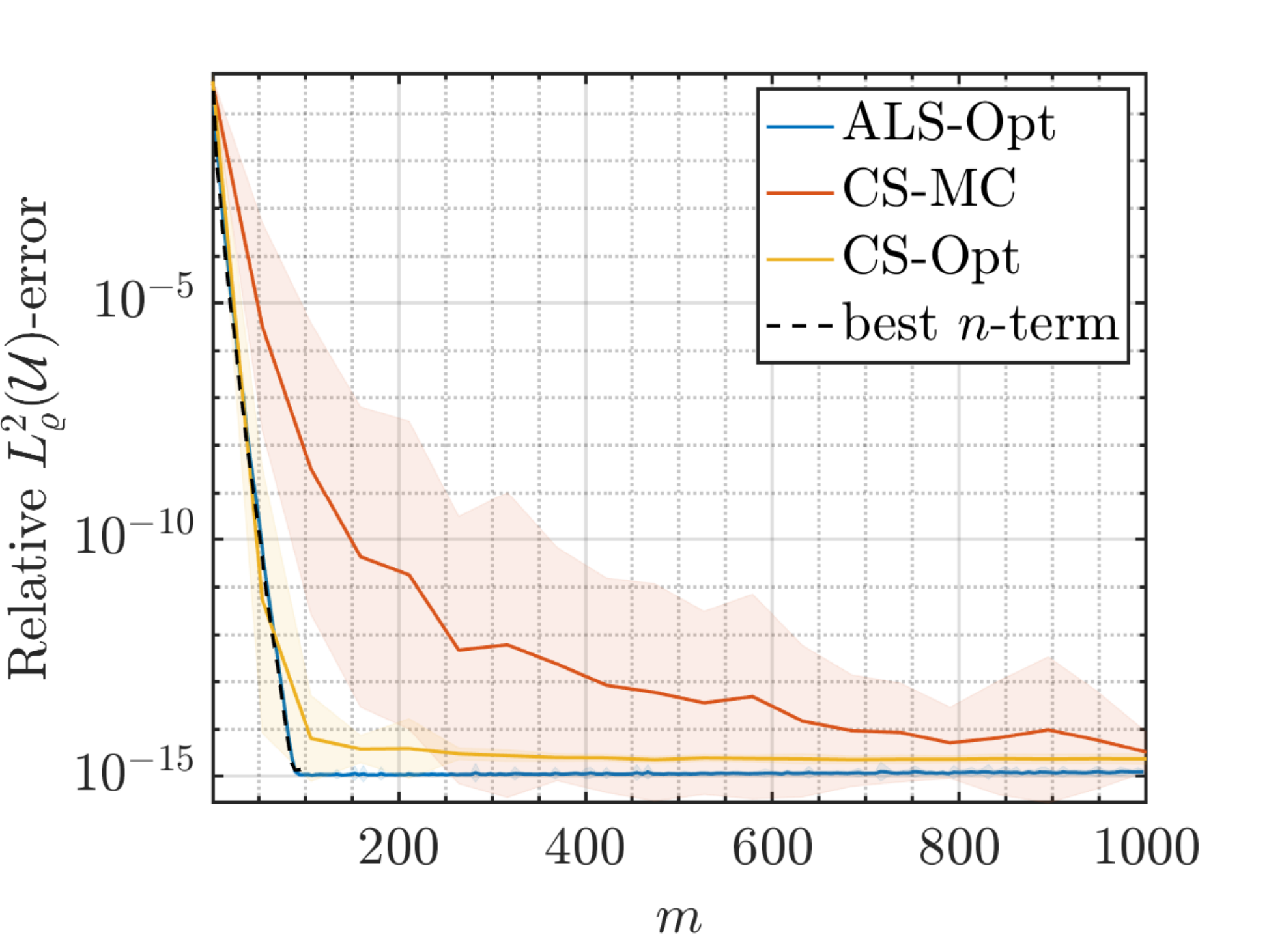}
&
\includegraphics[width = \errplotimg]{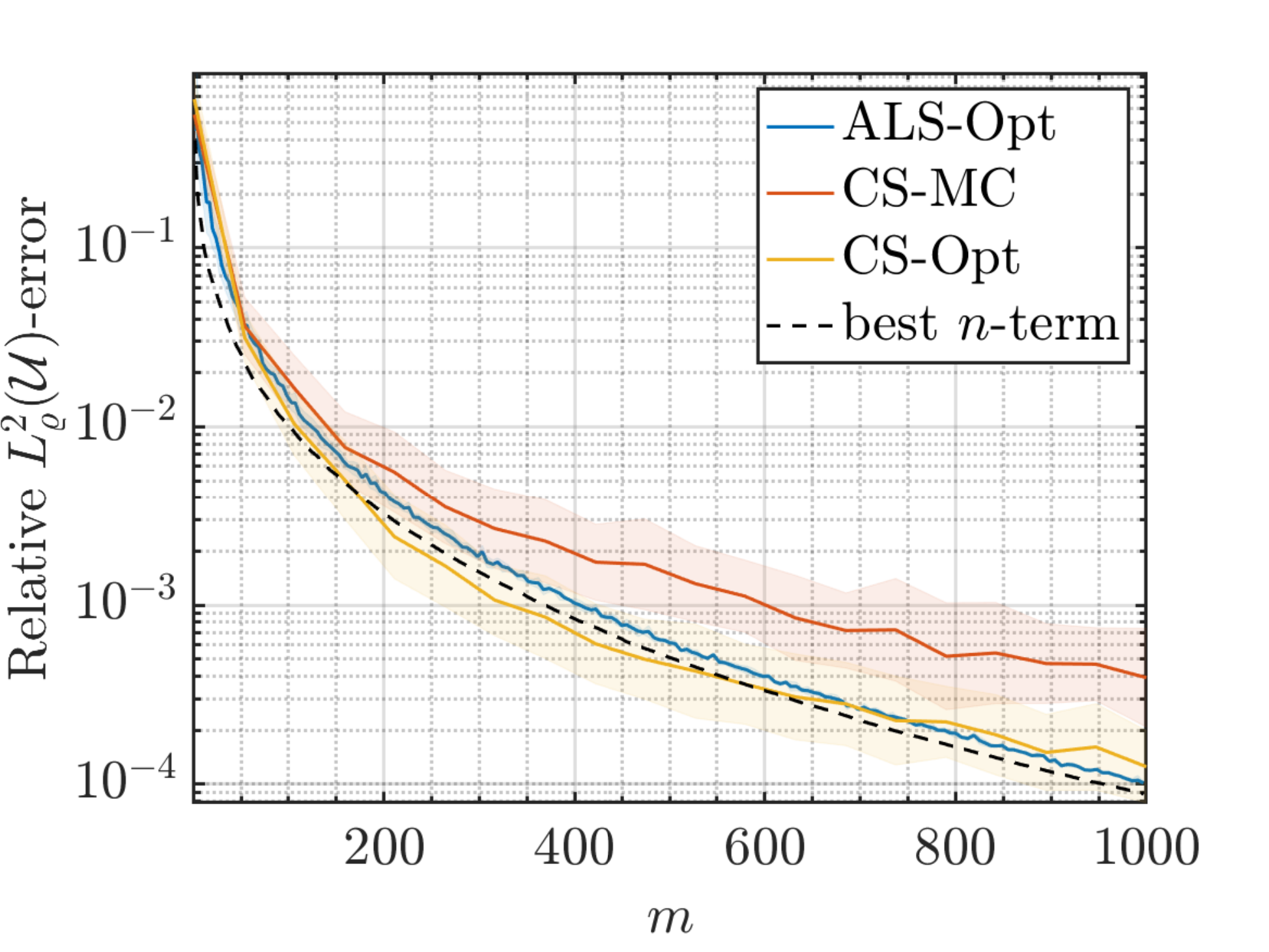}
&
\includegraphics[width = \errplotimg]{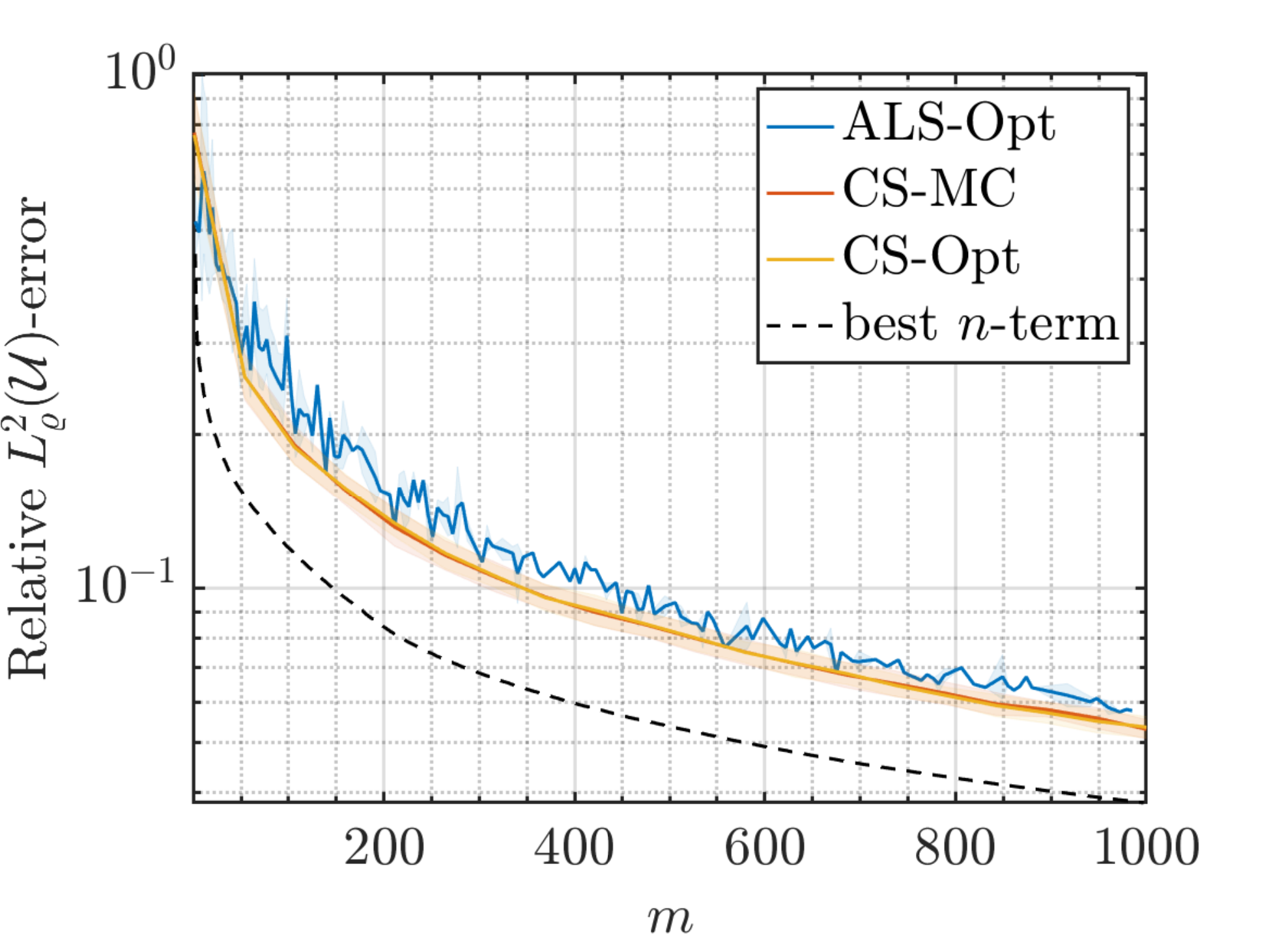}
\\[\errplottextsp]
\includegraphics[width = \errplotimg]{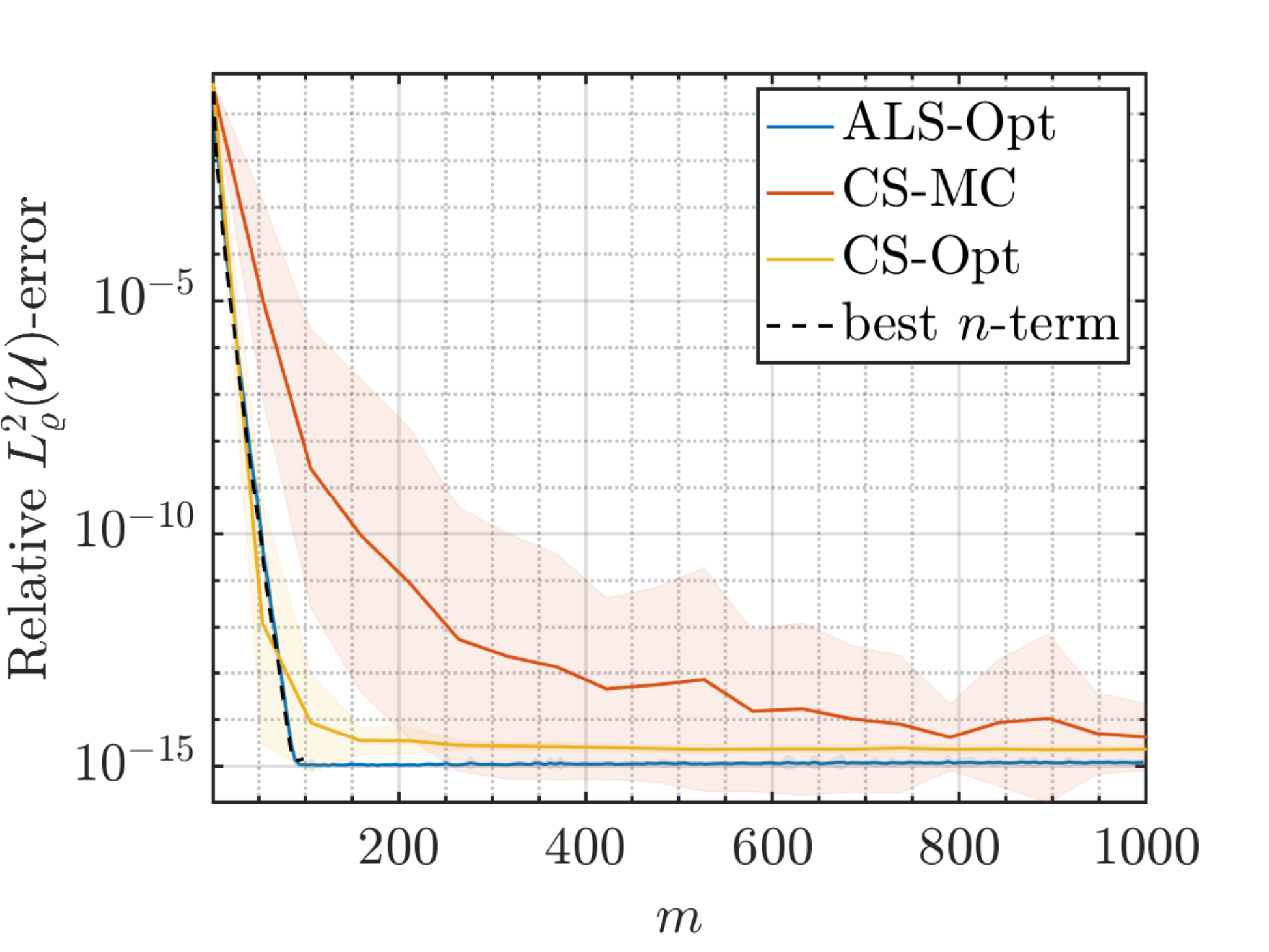}
&
\includegraphics[width = \errplotimg]{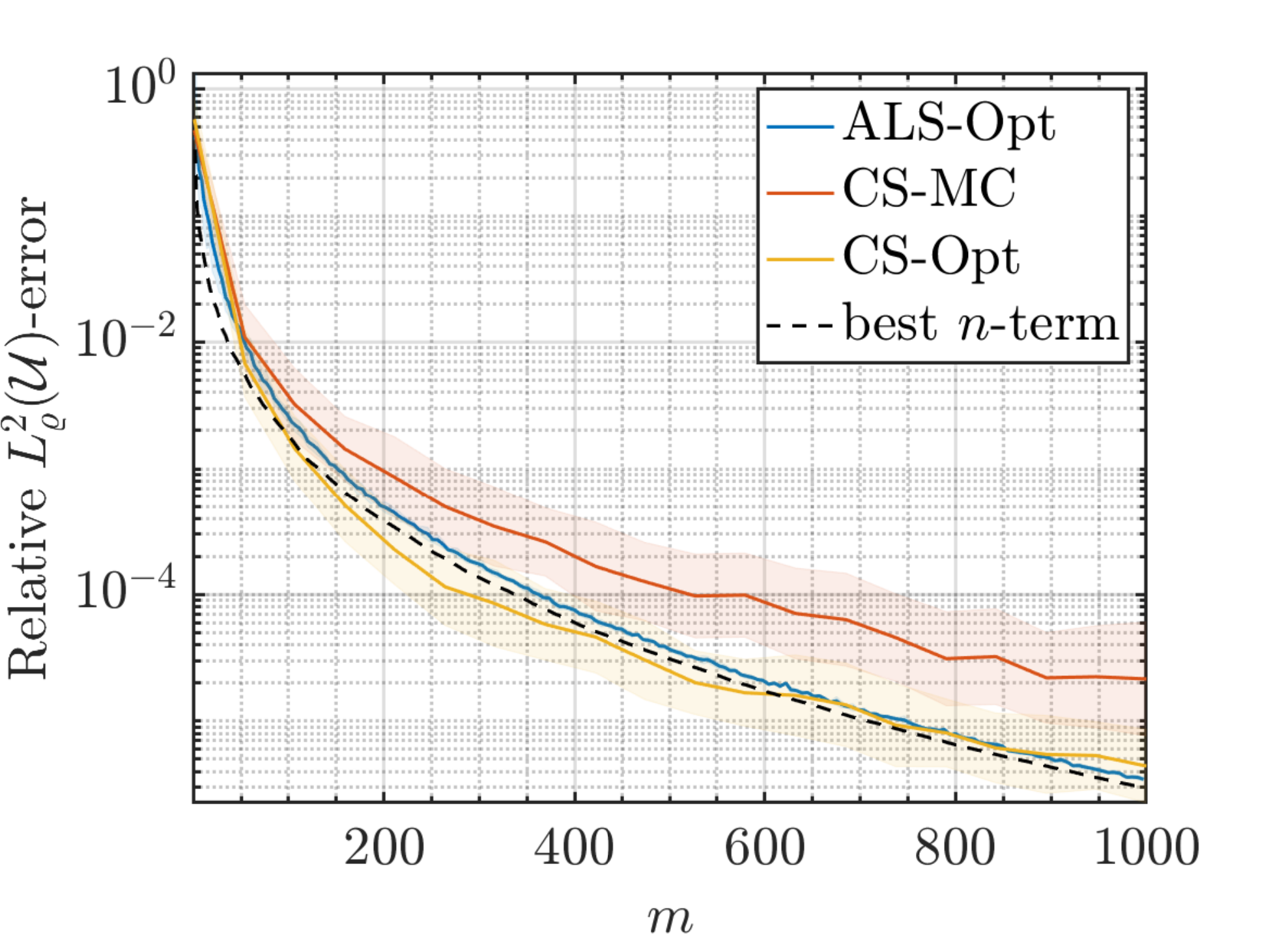}
&
\includegraphics[width = \errplotimg]{fig_12a_d32}
\\[\errplottextsp]
$d = 1$ & $d = 4$ & $d = 32$
\end{tabular}
\end{small}
\end{center}
\caption{The same as Fig.\ \ref{fig:fig10}, except with $f = f_3$ with $\delta_i = i$ (first row) and $\delta_i = i^2$ (second row).} 
\label{fig:fig12}
\end{figure}

\begin{figure}[t!]
\begin{center}
\begin{small}
 \begin{tabular}{@{\hspace{0pt}}c@{\hspace{\errplotsp}}c@{\hspace{\errplotsp}}c@{\hspace{0pt}}}
\includegraphics[width = \errplotimg]{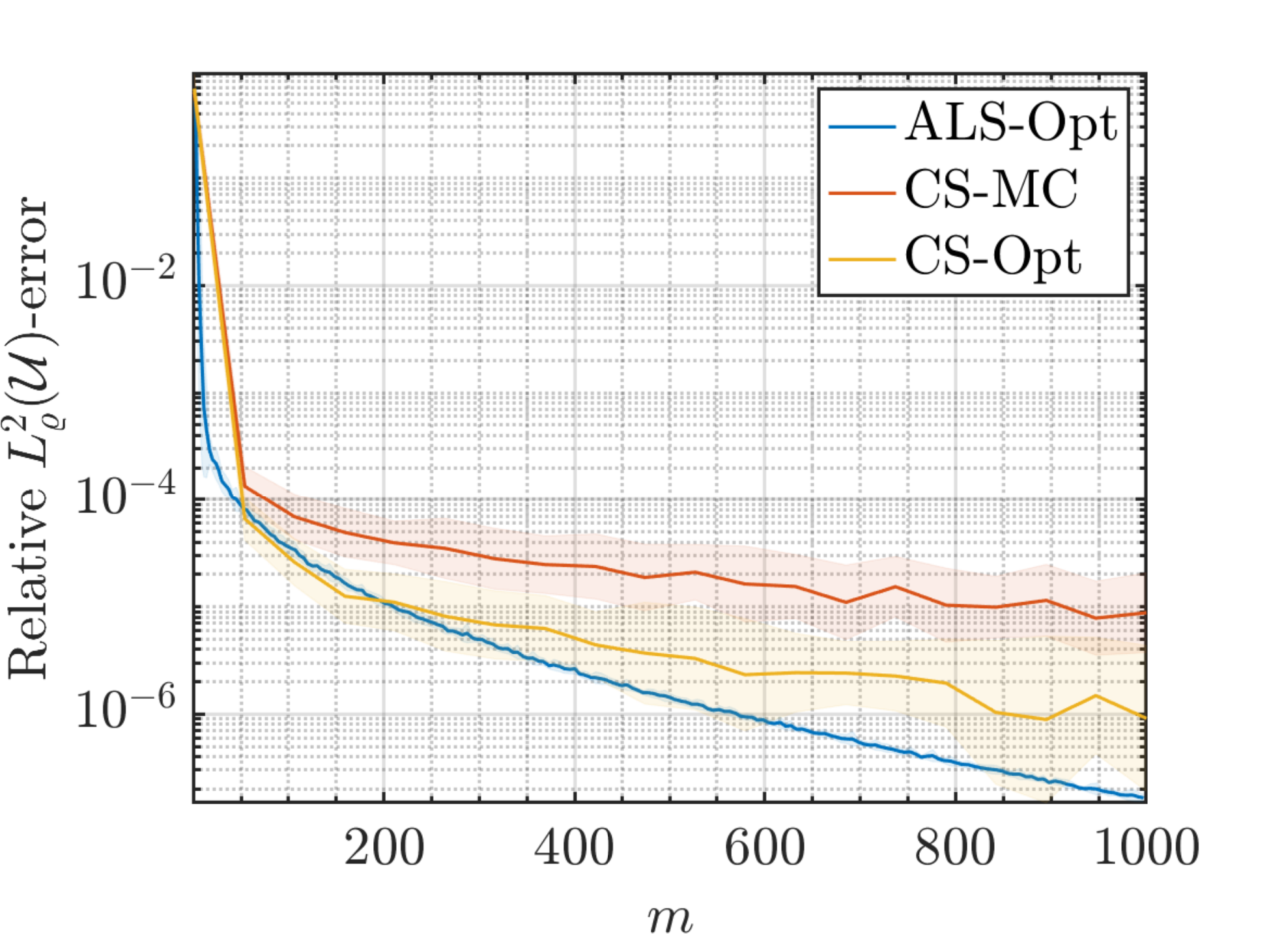}
&
\includegraphics[width = \errplotimg]{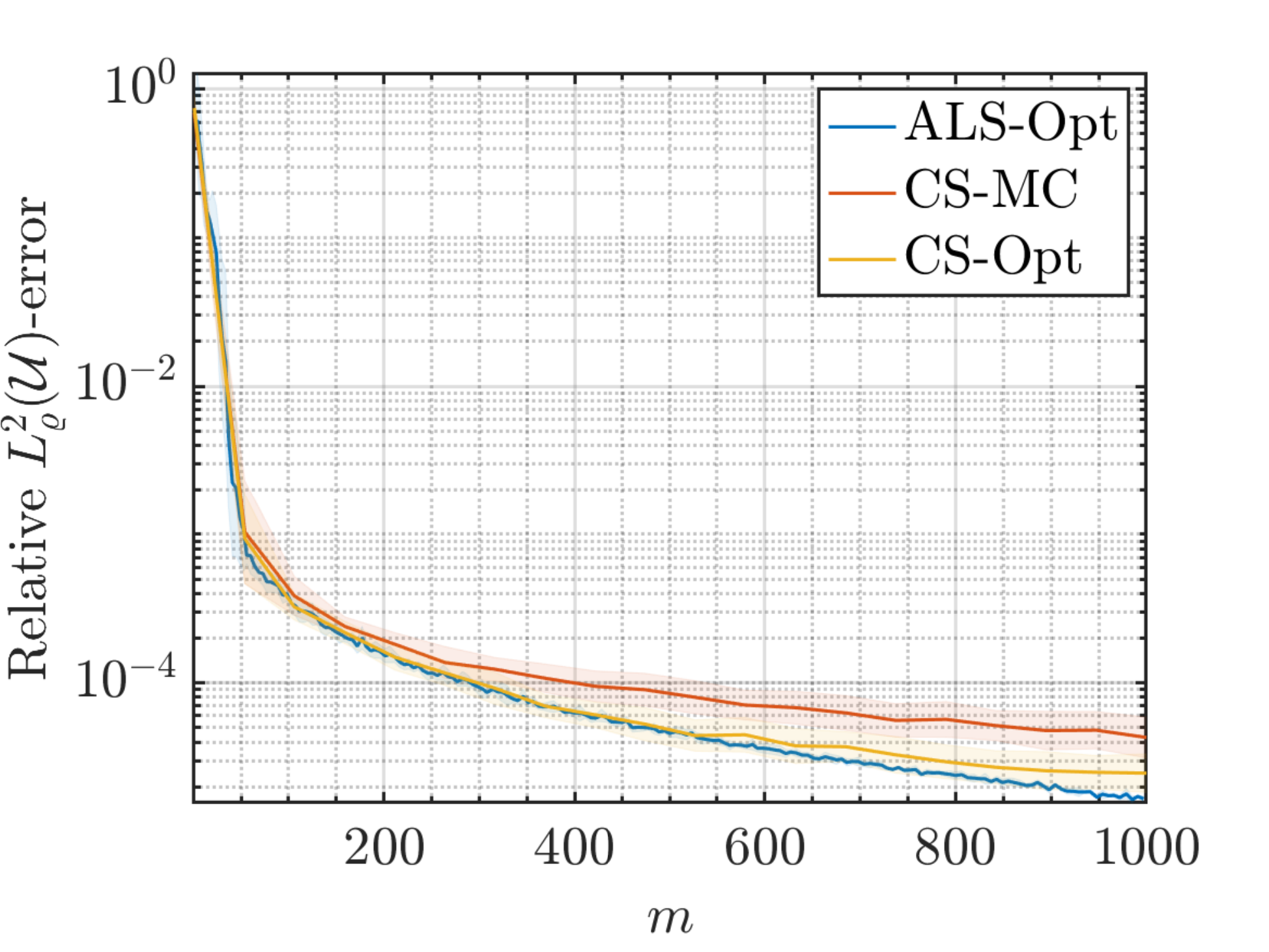}
&
\includegraphics[width = \errplotimg]{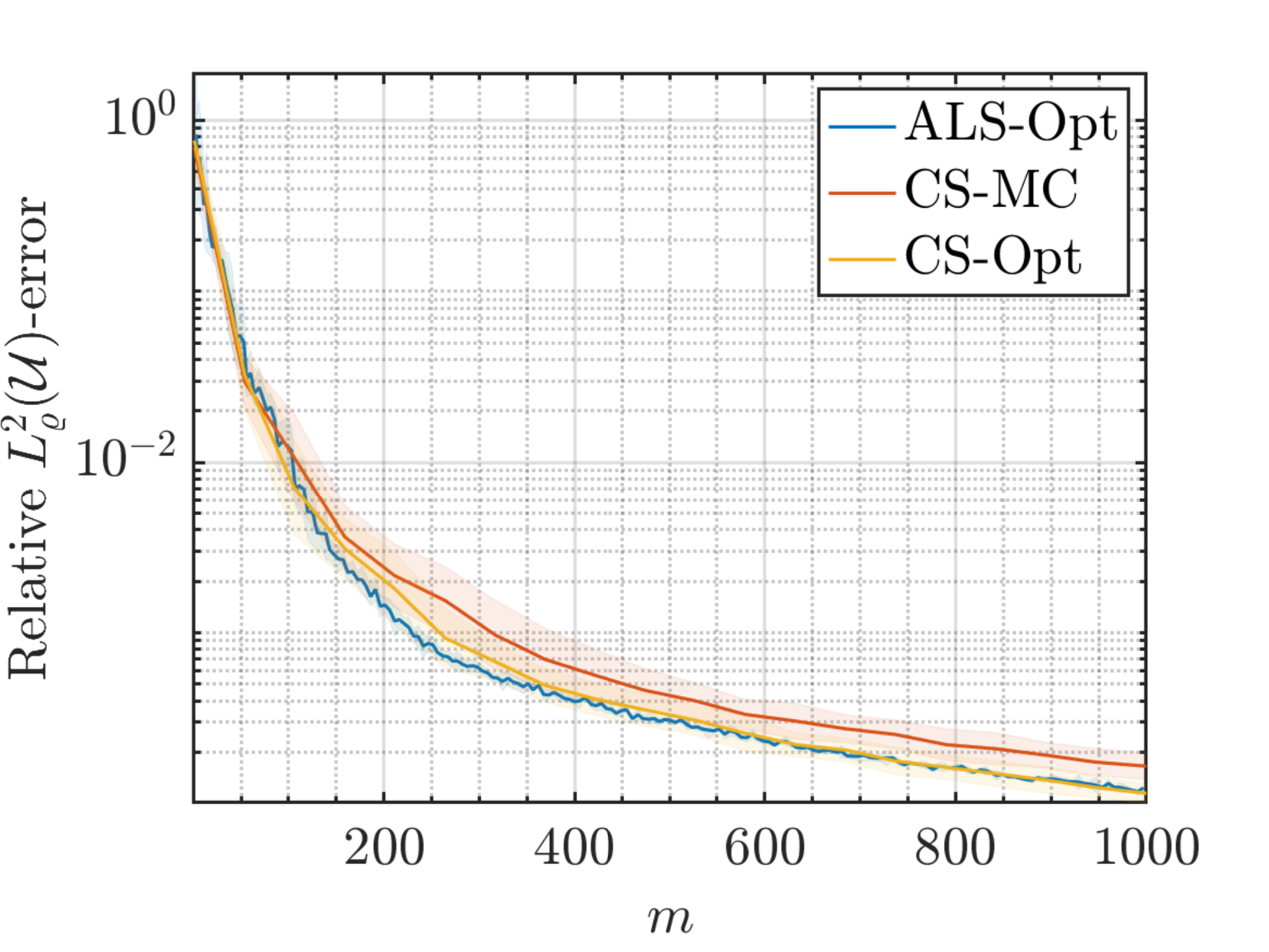}
\\[\errplottextsp]
$d = 2$ & $d = 6$ & $d = 8$
\end{tabular}
\end{small}
\end{center}
\caption{{The same as Fig.\ \ref{fig:fig10}, except with $f = f_4$}.} 
\label{fig:fig13}
\end{figure}

\begin{figure}[t!]
\begin{center}
\begin{small}
 \begin{tabular}{@{\hspace{0pt}}c@{\hspace{\errplotsp}}c@{\hspace{\errplotsp}}c@{\hspace{0pt}}}
\includegraphics[width = \errplotimg]{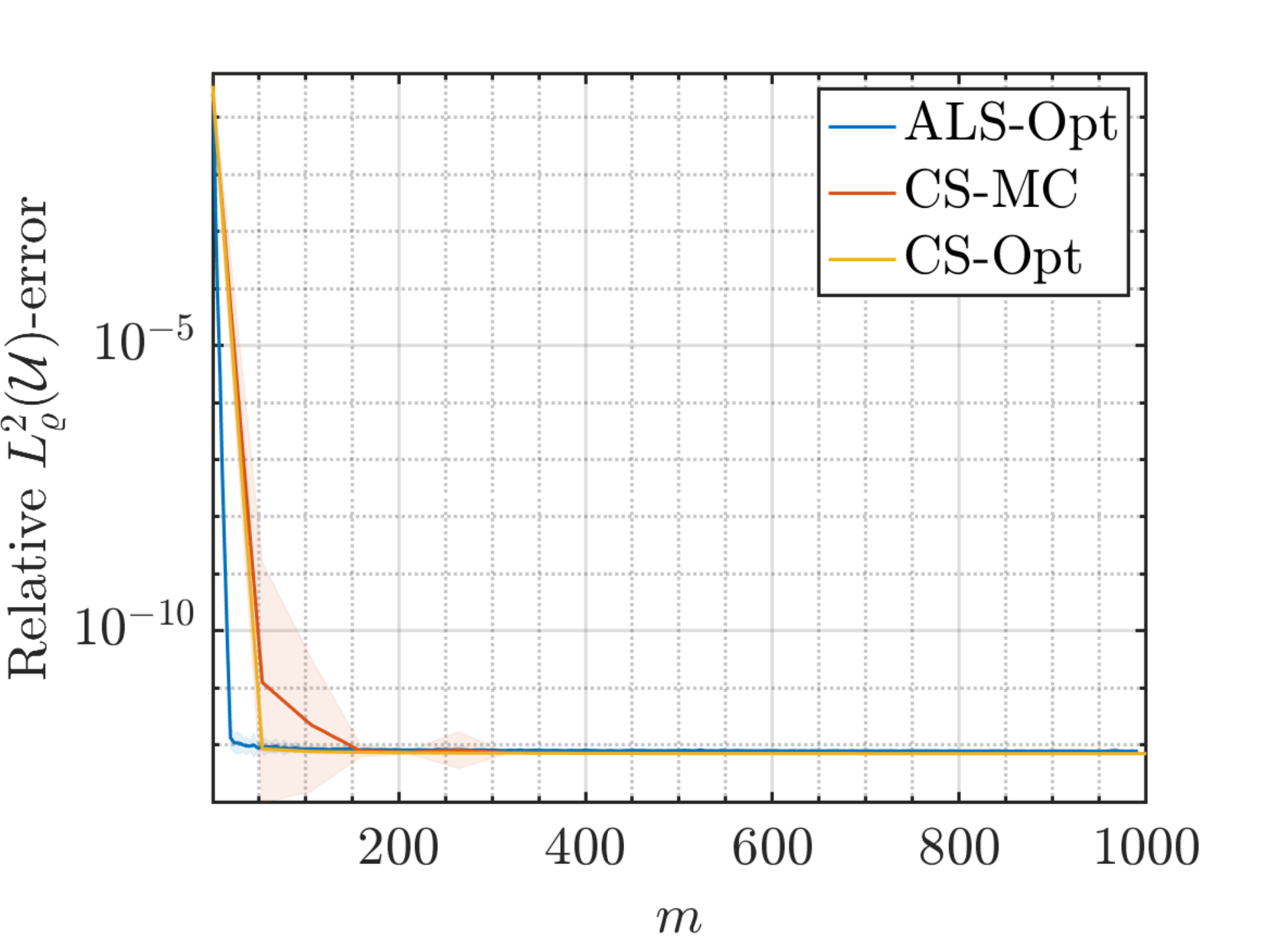}
&
\includegraphics[width = \errplotimg]{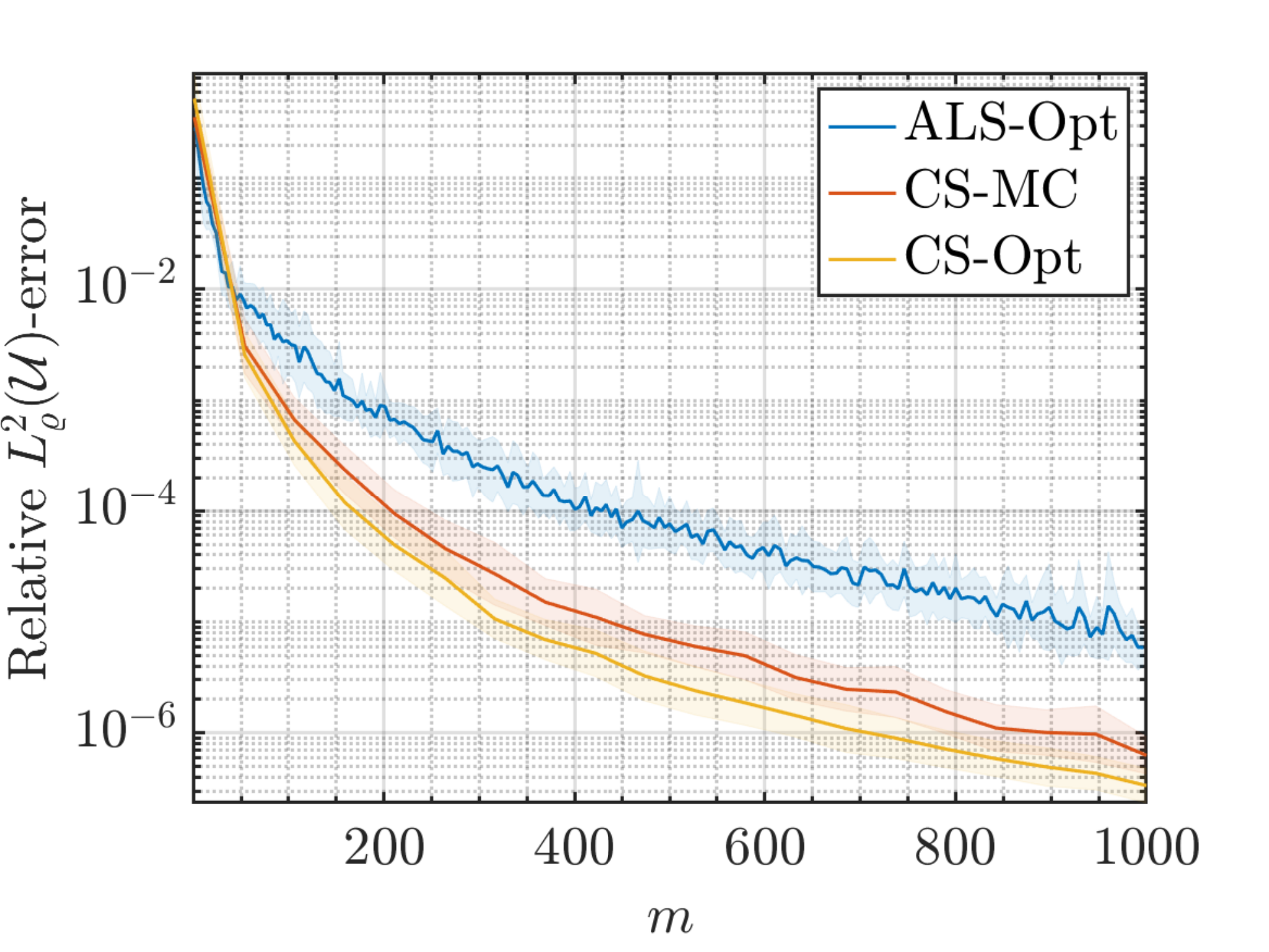}
&
\includegraphics[width = \errplotimg]{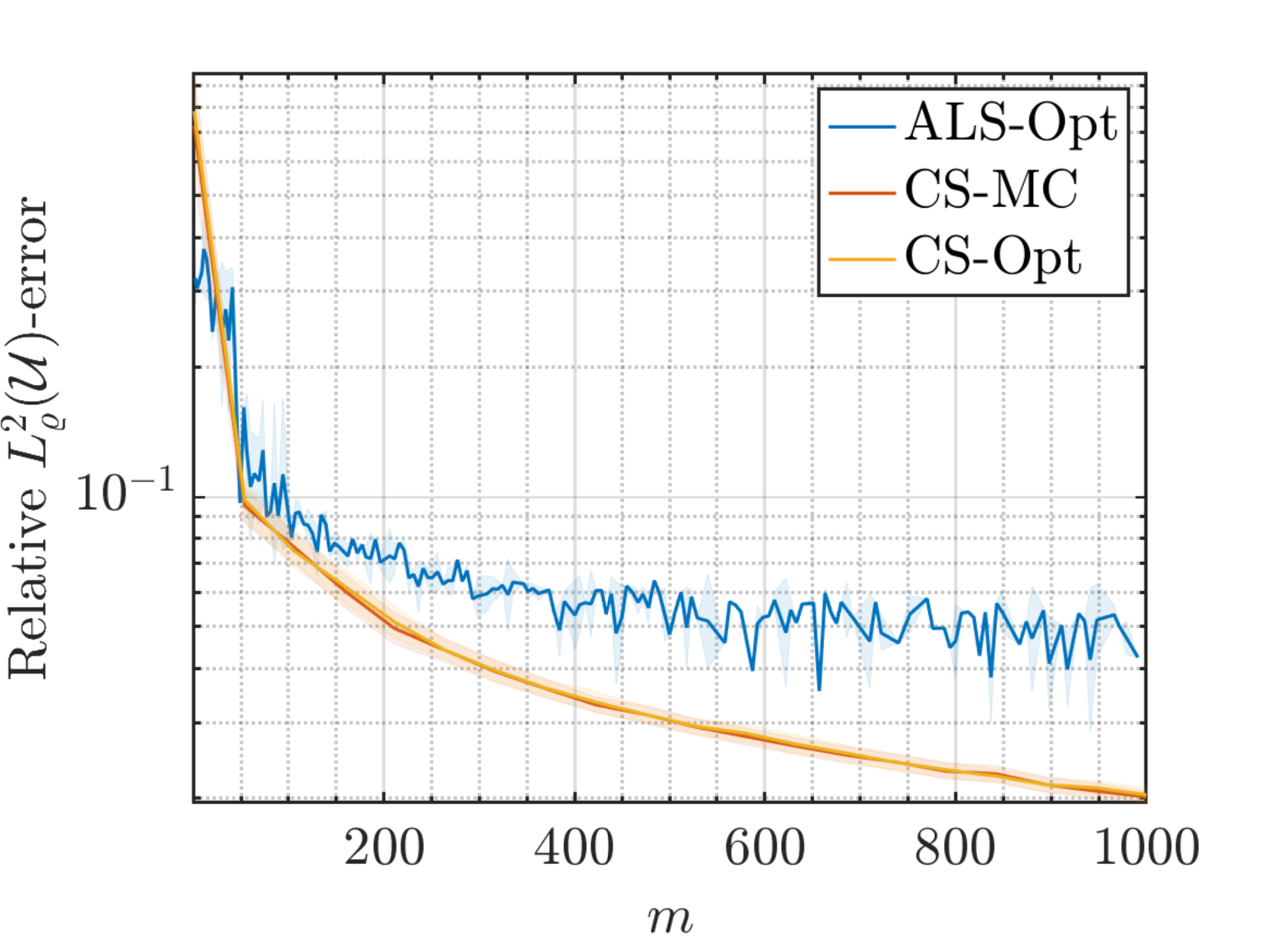}
\\[\errplottextsp]
$d = 1$ & $d = 4$ & $d = 32$
\end{tabular}
\end{small}
\end{center}
\caption{The same as Fig.\ \ref{fig:fig10}, except for the parametric DE example considered in Fig.\ \ref{fig:fig1}.} 
\label{fig:fig14}
\end{figure}

Results are shown in Figs.\ \ref{fig:fig10}--\ref{fig:fig14}. In low dimensions, ALS with the near-optimal sampling strategy typically outperforms CS with MC sampling, both in terms of the mean of the approximation error and its variance. It is notable, however, that the CS scheme performs much better than the previously-studied case of ALS with MC sampling (see Figs.\ \ref{fig:fig2}--\ref{fig:fig5}). Moreover, once the dimension increases, the CS scheme performs at least as well as, or sometimes better than the ALS scheme.

The worse performance of CS with MC sampling in low dimensions can be largely attributed to the sampling strategy. To see this, we also compute a CS approximation where the samples are drawn randomly from the measure
\be{
\label{CS-meas-opt}
\D \mu(\bm{y}) = N^{-1} \cK(\cP_{\Lambda})(\bm{y}) \D \varrho(\bm{y}) = N^{-1}  \sum_{\bm{\nu} \in \Lambda} | \Psi_{\bm{\nu}}(\bm{y}) |^2 \D \varrho(\bm{y}).
}
This is precisely the near-optimal sampling measure (if one were to use LS) for the polynomial space $\cP_{\Lambda}$ in which the CS approximation is sought. Note, however, that there are no theoretical guarantees that it is near-optimal for CS; in fact, theoretically-optimal sampling measures for CS currently do not currently exist \cite{adcock2022towards}. Nonetheless, this scheme significantly improves the performance of CS in low dimensions, rendering it competitive with ALS in this regime. However, as per the main thesis of this paper, the improvement over CS with MC sampling lessens as the dimension increases.

Finally, it is worth noting that the CS scheme also achieves high accuracy: the limiting error (when reached) is around $10^{-14}$ or $10^{-15}$ in all cases. A bane of CS-based polynomial approximations in past work was the inability of off-the-shelf solvers to achieve high accuracy. In this paper, we employ an efficient algorithm developed in \cite{adcock2022efficient} {and related to \cite{colbrook2022warpd}.}
With this scheme, an error of $\eta$ can be achieved efficiently, using a number of iterations proportional to $\log(\eta^{-1})$. See \S \ref{ss:efficient-PDI} for further information.

\section{Conclusions}\label{s:conclusions}

The purpose of this paper has been to show that MC sampling is a not only a good sampling strategy for {polynomial} approximation in high dimensions, but, in fact, a near-optimal one. Hence, efforts to design optimal sampling strategies are, for this problem at least, only effective in low dimensions.

\subsection{{Caveats}}

There are several caveats {we now mention}. First, our study is limited to polynomial approximation{. There are many other popular tools for high-dimensional approximation which may outperform such schemes in practice. We mention in passing Gaussian process regression (kriging), kernel methods, reduced-order models and, recently, deep neural networks. These approaches are outside the scope of this work. Second, our study is also limited to} bounded hypercubes, which we take to be $[-1,1]^{\bbN}$ without loss of generality. The situation in either the half space $[0,\infty)^{\bbN}$ or whole space $\bbR^{\bbN}$ is quite different. Here MC sampling performs substantially worse even in moderately large dimensions \cite{guo2020constructing,narayan2017christoffel}. Whether it is optimal in infinite dimensions {is currently unknown}. {Third,} our study is limited to the Chebyshev and uniform measures on $[-1,1]^{\bbN}$. However, we expect identical conclusions to hold for more general Jacobi measures. See \S \ref{ss:different-measures} for some further discussion on this point.  {Fourth, in this paper we have only considered errors measured in the $L^2_{\varrho}$-norm. In many applications, one also desires pointwise control of the functions being approximated, i.e., $L^{\infty}$-norm errors. Empirically, one observes the same effect if the error is measured in this norm instead (see \S\ref{ss:uniform-norm-exp} for several such experiments). Whether our theoretical analysis extends to the $L^{\infty}$-norm is currently unknown.}

{Fifth,} our analysis is also limited to classes of holomorphic functions. Clearly, other classes -- in particular, classes of piecewise smooth functions or functions with local features --  stand to benefit from changing the sampling strategy, e.g., in an adaptive, function-dependent manner. 
However, we note that MC sampling has also shown to be optimal for approximation in Sobolev spaces \cite{krieg2023random}.

{Finally, as mentioned in \S \ref{ss:theory-practice}, there is a gap between our theoretical analysis and our numerical examples. 
Fig.\ \ref{fig:fig9} and Remark \ref{rem:f3-theory} suggest that the infinite-dimensional analysis does indeed govern the empirical performance in finite (and moderate) dimensions, at least for some functions. Further investigation is needed in this direction.}

\subsection{{Future work}}

Our main results are \textit{nonuniform} guarantees, in the sense that a single draw of the sample points ensures the error bound for each fixed $f \in \cH(\bm{b},\varepsilon)$ with high probability. {Optimal and near-optimal sampling measures that yield \textit{uniform} guarantees for weighted LS approximation, i.e., those holding simultaneously for all functions in a certain class, can be constructed. See, e.g., \cite{krieg2021functionII}}. We believe these techniques could be modified to prove uniform guarantees for LS with MC sampling in the class $\cH(\bm{b},\varepsilon)$. See also \cite[Thm.\ 7.13]{adcock2022sparse} in the case of CS.

The numerical comparison in \S \ref{ss:CS-ALS-numerical} warrants some further discussion. In parametric modelling and UQ, the substantial computational burden is usually in generating the samples. Nonetheless, it is worth noting that the CS scheme is more expensive than the ALS scheme, both in terms of time and memory. For both schemes, the cost-per-iteration of an iterative solver -- e.g., conjugate gradients in the case of LS or the scheme described in \S \ref{ss:efficient-PDI} for CS -- is proportional to $m p$ (the cost of matrix-vector multiplications) and the number of iterations needed for an accuracy of $\eta$ is proportional to $\log(\eta^{-1})$. Here $p$ is the number of columns of the matrix $\bm{A}$. For LS, we have $p = |S| \leq m$. However, in the CS it is much larger: namely, $p = |\Lambda|$, where $\Lambda$ is the truncated index set. In  \cite{adcock2022efficient} it was shown that the computational cost per iteration for achieving the rates of Theorem \ref{t:main-res-2} is \textit{subexponential} in $m$, specifically, $\ord{m^{1+(1+\alpha) \log_2(4 m)}}$, where $\alpha = \log(3) / \log(4) \approx 0.79$. This situation arises because $\Lambda = \Lambda^{\mathsf{HCI}}_n$ is an isotropic index set, even though the underlying function is anisotropic. 

Whether it is possible to achieve the same rates as Theorem \ref{t:main-res-2} with a polynomial-time algorithm is {currently unknown}. Recent work of \cite{choi2021sparse2} may provide an answer. Interestingly, the algorithms in \cite{choi2021sparse2} use structured sampling strategies to achieve polynomial-time complexity. Hence, this could lead to the interesting conclusion that MC sampling is acceptable for approximation purposes, but not for fast computations. Conversely, if one is willing to forgo theoretical guarantees with MC sampling, then one can develop significantly more efficient CS-based schemes using \textit{basis adaptivity} \cite{jakeman2015enhancing,hampton2018basis}, in which $\Lambda$ is adapted at each step to the anisotropy of the function.

\section*{Acknowledgments}
The authors would like to thank David Krieg and Mario Ullrich for useful feedback on an early draft of the manuscript {and the two anonymous referees for their insightful comments, which considerably improved its quality}.

\bibliographystyle{siamplain}
\bibliography{MCgoodrefs}
\end{document}

{\maketitle}

These supplementary materials include further details in support of the main paper. They are organized as follows. {We start by providing some technical details about the existence of infinite tensor-product measures in \S \ref{s:inf_dim_measures}.} In \S \ref{s:proofs} we give the proofs of the main results in the paper. In \S \ref{s:additional-info} we provide additional information on the numerical experiments. Next, in \S \ref{s:pDE-background} we give further background on parametric DEs. Finally, we end in \S \ref{s:further-experiments} with several additional experiments.

\section{{Existence of probability measures on $[-1,1]^{\mathbb{N}}$}}\label{s:inf_dim_measures} 
{In \S\ref{ss:inf_dim_setup} we mentioned in passing that the existence of a tensor-product measure on the infinite product set $[-1,1]^{\mathbb{N}}$ is guaranteed by the Kolmogorov extension theorem. We now provide further details. For the sake of reader's convenience, we state the Kolmogorov extension theorem below (see, e.g.,  {\cite[Theorem 2.4.4]{tao2011introduction}}). Recall that a topological set is \emph{locally compact} if its points have a compact neighbourhood and it is \emph{$\sigma$-compact} if it is the union of countably-many compact subsets.}

\begin{theorem}[Kolmogorov's extension theorem]
{Let $A$ be an arbitrary set. For each $\alpha \in A$, let $(X_\alpha, \mathcal{B}_\alpha, \mu_{\alpha})$ be a probability space in which $X_\alpha$ is a locally compact, $\sigma$-compact metric space, with $\mathcal{B}_\alpha$ being its Borel $\sigma$-algebra (i.e. the $\sigma$-algebra generated by the open sets). Then there exists a unique probability measure $\mu_A =: \prod_{\alpha \in A} \mu_\alpha$  
on $(X_A, \mathcal{B}_A) := (\prod_{\alpha \in A} X_\alpha, \prod_{\alpha \in A} \mathcal{B}_\alpha)$ with the property that
\begin{equation}
\label{eq:inf_prod_measure}
\mu_A\left(\prod_{\alpha \in A}E_\alpha\right) = \prod_{\alpha \in A} \mu_{\alpha}(E_\alpha)
\end{equation}
whenever $E_\alpha \in \mathcal{B}_\alpha$ for each $\alpha \in A$, and one has $E_\alpha = X_\alpha$ for all but finitely many of the $\alpha$.}
\end{theorem}

{This result readily applies to our setting (see also \cite[Example 2.4.7]{tao2011introduction}). In fact, it suffices to let $A = \mathbb{N}$, and $X_{\alpha}  = [-1,1]$ and $\mu_\alpha = \varrho$, for every $\alpha \in A$, where $\varrho$ is a probability measure on $[-1,1]$. Note that the interval  $[-1,1]$ is a compact Hausdorff space and it is, therefore, also a locally compact, $\sigma$-compact metric space. We also observe that, in the context of our work, property \eqref{eq:inf_prod_measure} reads\footnote{{Recall that, with a slight abuse of notation, $\varrho$ denotes both a probability measure on the interval $[-1,1]$ and the corresponding tensor-product measure $\varrho = \prod_{k=1}^\infty\varrho $ over $[-1,1]^\mathbb{N}$.}}
$$
\varrho\left(\prod_{k = 1}^\infty X_k\right) = \prod_{k = 1}^\infty \varrho (X_k),
$$
for any sequence $(X_k)_{k=1}^\infty$ of measurable subsets of $[-1,1]$ such that $X_k \equiv [-1,1]$ for all but finitely many values of $k \in \mathbb{N}$.
}

\section{Proofs}\label{s:proofs}

\subsection{Lemma \ref{l:wLS-acc-stab} and Theorem \ref{t:wLS-samp-comp}}

Both results are minor adaptations and combinations of existing results. See, for instance, \cite[Thms.\ 5.3 \& 5.19]{adcock2022sparse}. We include brief proofs for completeness:

\prf{
[Proof of Lemma \ref{l:wLS-acc-stab}]

Recall that $\alpha_{w} = \sigma_{\min}(\bm{A})$, where $\bm{A}$ is the matrix defined in \eqref{LS-meas-mat-vec} associated with the algebraic LS problem. This matrix is full rank since $\alpha_{w} > 0$, and therefore the LS problem has a unique solution. Furthermore, using the normal equations, it is a short argument to show that $\hat{f}$ is defined by
\bes{
\ip{\hat{f}}{q}_{\mathsf{disc},w} = \ip{f}{q}_{\mathsf{disc},w} + \frac1m \sum^{m}_{i=1} w(\bm{y}_i) e_i \overline{q(\bm{y}_i)},\quad \forall q \in \cP,
}
where $\ip{\cdot}{\cdot}_{\mathsf{disc},w}$ is the discrete semi-innner product associated with the seminorm $\nm{\cdot}_{\mathsf{disc},w}$.
Let $p \in \cP$ be arbitrary. Setting $q = \hat{f}-p$ and subtracting $\ip{p}{q}_{\mathsf{disc},w}$ from both sides gives
\bes{
\nmu{\hat{f} - p}^2_{\mathsf{disc},w} = \ip{f - p}{\hat{f} - p}_{\mathsf{disc},w} +  \frac1m \sum^{m}_{i=1} w(\bm{y}_i) e_i \overline{\hat{f}(\bm{y}_i) - p(\bm{y}_i)}.
}
Several applications of the Cauchy-Schwarz inequality now yield
\bes{
\nmu{\hat{f} - p}_{\mathsf{disc},w} \leq \nmu{f - p}_{\mathsf{disc},w} + \sqrt{ \frac1m \sum^{m}_{i=1} w(\bm{y}_i) |e_i|^2}. 
}
For the latter term, we recall that the constant function is an element of $\cP$ by assumption.
Thus, by definition of $\beta_{w}$, we have
\bes{
\sqrt{ \frac1m \sum^{m}_{i=1} w(\bm{y}_i) |e_i|^2} \leq \nm{\bm{e}}_{\infty} \sqrt{ \frac1m \sum^{m}_{i=1} w(\bm{y}_i) } =  \nm{\bm{e}}_{\infty} \nm{1}_{\mathsf{disc},w} \leq \beta_w \nm{\bm{e}}_{\infty} \nm{1}_{L^2_{\varrho}(\cU)} .
}
Since $\varrho$ is a probability measure we have $\nm{1}_{L^2_{\varrho}(\cU)} =1$. Hence, using the definition of $\alpha_w$ as well, we get
\bes{
\nmu{\hat{f} - p}_{L^2_{\varrho}(\cU)} \leq \frac{1}{\alpha_w} \nmu{\hat{f} - p}_{\mathsf{disc},w} \leq \frac{1}{\alpha_w} \nmu{f - p}_{\mathsf{disc},w} + \frac{\beta_w}{\alpha_w} \nm{\bm{e}}_{\infty}.
}
The result now follows by writing $\nmu{f - \hat{f}}_{L^2_{\varrho}(\cU)} \leq \nmu{f - p}_{L^2_{\varrho}(\cU)} + \nmu{\hat{f} - p}_{L^2_{\varrho}(\cU)}$.
}

\prf{[Proof of Theorem \ref{t:wLS-samp-comp}]
Recall that $\alpha_{w} = \sigma_{\min}(\bm{A}) = \sqrt{\lambda_{\min}(\bm{A}^* \bm{A})}$ and $\beta_{w} =  \sigma_{\max}(\bm{A}) = \sqrt{\lambda_{\max}(\bm{A}^* \bm{A})}$. The proof follows a standard route. First, one writes $\bm{A}^* \bm{A}$ as a sum of i.i.d.\ random matrices. Then one uses the well-known matrix Chernoff bound (see, e.g., \cite[Thm.\ 1.1]{tropp2012user}) to bound the deviation of its maximum and minimum eigenvalues from  the corresponding eigenvalues of its mean $\bbE (\bm{A}^* \bm{A}) = \bm{I}$. For any $0 < \delta < 1$, this shows that
\bes{
\sqrt{1-\delta} < \alpha_{w} \leq \beta_{w} < \sqrt{1+\delta},
}
with probability at least $1-\epsilon$, provided
\bes{
m \geq ((1+\delta) \log(1+\delta) - \delta)^{-1} \cdot \kappa(\cP ; w) \cdot \log(2 n/\epsilon).
}
See, e.g., the proof of \cite[Thms.\ 5.12 \& 5.19]{adcock2022sparse}. To obtain Theorem \ref{t:wLS-samp-comp}, we set $\delta = 3/5$ (this value is arbitrary). This gives
\bes{
\frac{1}{\alpha_w} \leq \frac{1}{\sqrt{1-\delta}} < 2,\qquad \frac{\beta_w}{\alpha_w} \leq \sqrt{\frac{1+\delta}{1-\delta}} = 2,
}
and
\bes{
((1+\delta) \log(1+\delta) - \delta)^{-1} < 7.
}
Hence, we deduce from \eqref{wLS-samp-comp} and Lemma \ref{l:wLS-acc-stab} that $\hat{f}$ is unique and satisfies
\be{
\label{partial-bound}
\nmu{f - \hat{f}}_{L^2_{\varrho}(\cU)} \leq \inf_{p \in \cP} \left \{ \nmu{f - p}_{L^{2}_{\varrho}(\cU)} +2 \nmu{f - p}_{\mathsf{disc},w} \right \} + 2\nm{\bm{e}}_{\infty},
}
with probability at least $1-\epsilon$, as required.
}

\subsection{Theorem \ref{t:best-s-term}}\label{ss:best-n-term}

As noted, Theorem \ref{t:best-s-term} is {a well-known result. See, e.g., \cite[Thm.\ 3.28]{adcock2022sparse} or \cite[\S 3.2]{cohen2015approximation}), as well as \cite{beck2012optimal,bonito2021polynomial,chkifa2015breaking,cohen2011analytic,hansen2013analytic,todor2007convergence} for further background on polynomial approximation of holomorphic functions in infinite dimensions}. Its proof involves two main steps, which we now briefly describe, since similar ideas will be used later. 

The first step employs a result commonly referred to as \textit{Stechkin's inequality} (despite seemingly never appearing in a publication of Stechkin -- see \cite[\S 7.4]{dung2018hyperbolic}). For generality, let $0 < q \leq \infty$ and $n \in \bbN_0$. Then the \textit{$\ell^q$-norm best $n$-term approximation error} of a sequence $\bm{c} = (c_{\bm{\nu}})_{\bm{\nu} \in \cF}$ is defined as
\be{
\label{sigma-n-q}
\sigma_{n}(\bm{c})_q =  \min \left \{ \nm{\bm{c} - \bm{z}}_{p} : \bm{z} \in \ell^p(\cF), | \mathrm{supp}(\bm{z}) | \leq n \right \}.
}
Here, we recall that $\mathrm{supp}(\bm{z}) = \{ \bm{\nu} : z_{\bm{\nu}} \neq 0 \}$ for $\bm{z} = (z_{\bm{\nu}})_{\bm{\nu} \in \cF}$.
This quantity is particularly relevant for our purposes in the case $q = 2$. Indeed, let $f \in L^2_{\varrho}(\cU)$ with expansion \eqref{f-exp}. Then, by Parseval's identity, we readily see that
\be{
\label{sigma-n-best-n-term}
\sigma_{n}(\bm{c})_{2} = \nmu{f - f_n}_{L^2_{\varrho}(\cU)}
}
where, as in \eqref{best-n-term-inf}, $f_n$ is the best $n$-term approximation to $f$ and $\bm{c} = (c_{\bm{\nu}})_{\bm{\nu} \in \cF}$ is the sequence of coefficients of $f$.

With \eqref{sigma-n-q} to hand, Stechkin's inequality states that
\be{
\label{stechkin}
\sigma_{n}(\bm{c})_q \leq \nmu{\bm{c}}_p (n+1)^{\frac1q-\frac1p},\quad \forall n \in \bbN_0,\ 0 < p \leq q < \infty,
}
(see, e.g., \cite[Lem.\ 3.5]{adcock2022sparse}).
Hence, taking $q = 2$ and using \eqref{sigma-n-best-n-term}, we see that Theorem \ref{t:best-s-term} holds, provided
\be{
\label{coeff-p-summable-claim}
\nm{\bm{c}}_{p} \leq C(\bm{b},\varepsilon,p),\qquad \forall f \in \cH(\bm{b},\varepsilon),
}
where $\bm{c}$ is the sequence of coefficients of $f$.

This claim is established in the second step by using certain bounds for the polynomial coefficients $c_{\bm{\nu}}$ of a $(\bm{b},\varepsilon)$-holomorphic function. Specifically, let $\bm{\rho} = (\rho_i)^{\infty}_{i=1} \in [1,\infty)^{\bbN}$ be any sequence for which \eqref{b-eps-holo} holds. Then, for any $f \in \cH(\bm{b},\varepsilon)$,
\be{
\label{coeff-bound}
{| c_{\bm{\nu}} | \leq B(\bm{\nu},\bm{\rho}) : = \begin{cases} 2^{\nm{\bm{\nu}}_0/2} \bm{\rho}^{-\bm{\nu}} & \text{(Chebyshev)} \\ \left ( \prod_{j \in \mathrm{supp}(\bm{\nu})} \xi(\rho_j) \sqrt{2 \nu_j + 1} \right ) \bm{\rho}^{-\bm{\nu}} & \text{(Legendre)} \end{cases} ,}
}
where $\bm{\rho}^{-\bm{\nu}} = \prod^{\infty}_{j=1} \rho^{-\nu_j}_j$ and $\xi(t) = \min \{ 2 t , \frac{\pi}{2} (t + t^{-1}) \} / (t-1)$ for every $t > 1$. Here we also recall that $\nm{\bm{\nu}}_0 = | \mathrm{supp}(\bm{\nu}) |$ is the number of nonzero entries in $\bm{\nu}$. The claim \eqref{coeff-p-summable-claim} is then proved using a rather technical argument involving these bounds. We omit the details, and refer to, e.g., \cite[Thm.\ 3.28]{adcock2022sparse}.

\subsection{Theorem \ref{t:alg-rate-sharp}}
\label{ss:Stechkin}

As discussed in the previous subsection, the first part of Theorem \ref{stechkin} is based on Stechkin's inequality \eqref{stechkin}. In order to prove Theorem \ref{t:alg-rate-sharp}, we first require the following result, which shows that such algebraic rates of convergence precisely characterize the weak $\ell^p$ spaces: 

\lem{
[Stechkin's inequality in $w \ell^p$]
Let $0 < p <q < \infty$. Then a sequence $\bm{c} \in w \ell^p(\bbN)$ if and only if $\sigma_n(\bm{c})_q \leq C n^{\frac1q-\frac1p}$, $\forall n \in \bbN$, for some constant $C$ depending on $\bm{c}$, $p$ and $q$ only. Specifically, if $\bm{c} \in w \ell^p(\bbN)$, then
$$
\sigma_n(\bm{c})_q \leq \frac{\|\bm{c}\|_{p,\infty}}{(q/p-1)^{1/q}} n^{\frac1q-\frac1p}, \quad \forall n \in \mathbb{N}.
$$
Conversely, if $\sigma_n(\bm{c})_q \leq C n^{\frac1q-\frac1p}$, for every $n \in \mathbb{N}$, then  
$\|\bm{c}\|_{p,\infty} \leq 2^{\frac1p+\frac1q} C$. 
}
\prf{
Suppose first that $\bm{c} \in w \ell^p(\bbN)$. We follow the proof of \cite[Prop 2.11]{foucart2013mathematical}. Let $\bm{c}^* = (c^*_i)^{\infty}_{i=1}$ be a nonincreasing rearrangement of the absolute values $(|c_i|)^{\infty}_{i=1}$.  Then
\bes{
(\sigma_n(\bm{c})_q)^q = \sum_{j > n} (c^*_j)^q 
\leq \sum_{j > n} \nm{\bm{c}}^q_{p,\infty} i^{-q/p} 
\leq  \nm{\bm{c}}^q_{p,\infty} \int_n^\infty x^{-\frac{q}{p}} \D x
= \frac{\nm{\bm{c}}^q_{p,\infty}}{q/p - 1} n^{1-\frac{q}{p}}
}
and therefore $\sigma_s(\bm{c})_q \leq C n^{\frac{1}{q}-\frac1p}$, as required.

Conversely, suppose that $\sigma_n(\bm{c})_q \leq C n^{\frac{1}{q}-\frac1p}$ for all $n \in \bbN$. Then
\bes{
n(c^*_{2n})^q \leq (c^*_{n+1})^q + \cdots + (c^*_{2n})^q \leq \sum_{i > n} (c^*_i)^q \leq C n^{1-\frac{q}{p}}.
}
Hence $c^{*}_{2n} \leq C n^{-\frac{1}{p}}$. Similarly, 
\bes{
n (c^*_{2n+1})^q \leq (c^*_{n+2})^q + \cdots + (c^{*}_{2n+1})^q \leq \sum_{i > n+1} (c^*_i)^q \leq C^q (n+1)^{1-\frac{q}{p}}.
}
Hence $c^{*}_{2n+1} \leq 2^{\frac1q} C (n+1)^{-\frac1p}$. Combining these two upper bounds, we obtain
\eas{
\|\bm{c}\|_{p,\infty}  
&= \sup_{i \in \mathbb{N}} i^{\frac1p} c^*_i
\\
&= \sup_{n \in \mathbb{N}} \left\{(2n)^{\frac1p} c^*_{2s}, (2n+1)^{\frac1p}c^*_{2n+1}\right\}
\\
&\leq C \sup_{n \in \mathbb{N}} \left\{\left(\frac{2n}{n}\right)^{\frac1p}, 2^{\frac1q}\left(\frac{2n+1}{n+1}\right)^{\frac1p}\right\}
\\
&= 2^{\frac1p + \frac1q} C.
}
The result now follows.
}

Using this, we immediately deduce the following:

\cor{
\label{c:exact-Stechkin}
Let $0 < p <q <  \infty$ and suppose that $\bm{c} \in w\ell^p(\bbN)$ but $\bm{c} \notin w \ell^r(\bbN)$ for some $r < p$. Then there exists a constant $C$ depending on $\bm{c}$, $p$ and $q$ only such that
\bes{
\sigma_{n}(\bm{c})_q \leq C n^{\frac1q-\frac1p},\quad \forall n \in \bbN.
}
However,
\bes{
\limsup_{n \rightarrow \infty} \frac{\sigma_n(\bm{c})_q}{n^{\frac1q-\frac1r}} = + \infty.
}
In other words, there does not exist a constant $C$ (depending on $p$, $r$ and $\bm{c}$) such that
\bes{
\sigma_{n}(\bm{c})_q \leq C n^{\frac1q-\frac1r},\quad  \forall n \in \bbN.
}
}

With this in mind, to prove Theorem \ref{t:alg-rate-sharp} we need only show that there are infinitely many functions $f \in \cH(\bm{b},\varepsilon)$ whose coefficients $\bm{c} \in w \ell^p(\cF)$, but $\bm{c} \notin w \ell^r(\cF)$ for $0 < r < p$. We now detail this construction.

\prf{[Proof of Theorem \ref{t:alg-rate-sharp}]
We will prove a slightly more general version of the theorem which states that for any $0 < p < q$,
\bes{
\limsup_{n \rightarrow \infty} \frac{\sigma_n(\bm{c})_q}{n^{\frac1q-\frac1r}} = +\infty,
}
where $\bm{c}$ are the coefficients of the function $f$. The result then follows as the special case $q = 2$.

Let $\bm{b} = (b_i)^{\infty}_{i=1} \in [0,\infty)^{\bbN}$, $\bm{b} \in \ell^p(\bbN)$. Then we define the function
\bes{
f(\bm{y}) = \sum^{\infty}_{j=1} b_j y_j,\qquad \bm{y} = (\bm{y}_i)^{\infty}_{i=1}.
}
We first show that this function is $(\bm{b},\varepsilon)$-holomorphic. Consider a Bernstein polyellipse $\cE_{\bm{\rho}}$ with $\bm{\rho}$ satisfying \eqref{b-eps-holo}. Recall that a Bernstein ellipse $\cE_{\rho}$ is contained in the disc $\cD_{\sigma} = \{ z \in \bbC : |z| \leq \sigma \}$ with $\sigma = (\rho + \rho^{-1})/2$. Let $\bm{y} \in \cE_{\bm{\rho}}$. Then
\bes{
\sum^{\infty}_{j=1} |b_j y_j | \leq \sum^{\infty}_{j=1} b_j \frac{\rho_j + \rho^{-1}_{j}}{2} = \sum^{\infty}_{j=1} b_j \left ( \frac{\rho_j + \rho^{-1}_{j}}{2} - 1 \right ) + \nm{b}_{1} \leq \epsilon + \nm{b}_{p} < \infty,
}
since $\bm{\rho}$ satisfies \eqref{b-eps-holo} and $\bm{b} \in \ell^p(\bbN)$ and $0 < p <1$. Hence the series converges absolutely over $\cE_{\bm{\rho}}$. By direct calculation, for $\bm{y} \in \cE_{\bm{\rho}}$, we have
\bes{
\lim_{h \rightarrow 0} \frac{f(\bm{y} + h \bm{e}_i) - f(\bm{y})}{h} = b_i.
}
Thus, $f$ is holomorphic in each variable, and therefore holomorphic at $\bm{y}$.  
Hence $f$ is holomorphic in $\cE_{\bm{\rho}}$ and, since $\bm{\rho}$ was arbitrary, in $\cR_{\bm{b},\varepsilon}$. We deduce that $f \in \cH(\bm{b},\varepsilon)$.

With this to hand, we now consider this function. Recall that the first orthonormal, univariate Legendre polynomial is precisely $\psi_1(y) = \sqrt{3} y$. Hence the nonzero Legendre polynomial coefficients of $f$ are precisely the values $b_j / \sqrt{3}$. In particular,
\bes{
\sigma_{s}(\bm{c})_q = \sigma_s(\bm{b})_q / \sqrt{3}.
}
The result now follows immediately from Corollary \ref{c:exact-Stechkin} and the fact that $\bm{b} \in \ell^p(\bbN)$, which implies in particular that $\bm{b} \in w \ell^p(\bbN)$. This shows the existence of a single function with the desired property. To show that there are infinitely many non-linearly dependent such functions, one may consider functions of the form $f(\bm{y}) = \sum^{\infty}_{j=1} b_j y_j + \gamma \Psi_{\bm{\nu}}(\bm{y})$, where $\bm{\nu} \neq \bm{e}_j$ for any $j$ and $\gamma \in \bbC$.
}

\subsection{Theorem \ref{t:main-res}}

For this proof, we require some additional concepts. First, recall the weights \eqref{weights-def}.
For $0 < q \leq 2$, we define the weighted $\ell^p$-norm of a sequence $\bm{c} = (c_{\bm{\nu}})_{\bm{\nu} \in \cF}$ as
\bes{
\nmu{\bm{c}}_{q,\bm{u}} = \left ( \sum_{\bm{\nu} \in \cF} u^{2-q}_{\bm{\nu}} | c_{\bm{\nu}} |^q \right )^{1/q}.
}
Next, as in \S \ref{ss:weighted-k-term}, we define the weighted cardinality of a set $S \subset \cF$ as
\bes{
|S|_{\bm{u}} = \sum_{\bm{\nu} \in S} u^2_{\bm{\nu}}.
}
Finally, for $k \geq 0$ and $0 < q \leq 2$ we define the \textit{weighted best $(k,\bm{u})$-term approximation error} of a sequence $\bm{c} = (c_{\bm{\nu}})_{\bm{\nu} \in \cF}$  as
\bes{
\sigma_{k}(\bm{c})_{q,\bm{u}} = \min \left \{ \nm{\bm{c} - \bm{z}}_{q,\bm{u}} : \bm{z} \in \ell^q_{\bm{u}}(\cF), | \mathrm{supp}(\bm{z}) |_{\bm{u}} \leq k \right \}.
}
Similar to the unweighted case \eqref{stechkin}, one also has a version of Stechkin's inequality for $\sigma_{k}(\bm{c})_{q,\bm{u}}$. See, e.g., \cite[Lem.\ 3.12]{adcock2022sparse}. This reads
\be{
\label{weighted-Stechkin}
\sigma_{k}(\bm{c})_{q,\bm{u}} \leq \nm{\bm{c}}_{p,\bm{u}} k^{\frac1q-\frac1p},\quad \forall k > 0, 0 < p \leq q \leq 2.
}

\prf{
[Proof of Theorem \ref{t:main-res}]
Let $f \in \cH(\bm{b},\varepsilon)$ with coefficients $\bm{c} = (c_{\bm{\nu}})_{\bm{\nu} \in \cF}$ as in \eqref{f-exp-inf}. {
Define the sequence $\tilde{\bm{c}} = (\tilde{c}_{\bm{\nu}})_{\bm{\nu} \in \cF}$ by
\bes{
\tilde{c}_{\bm{\nu}} = \inf \left \{ B(\bm{\nu},\bm{\rho}) : \text{$\bm{\rho}$ satisfies \eqref{b-eps-holo}} \right \},
}
where $B(\bm{\nu},\bm{\rho})$ is as in \eqref{coeff-bound}.}
Notice that $|c_{\bm{\nu}}| \leq \tilde{c}_{\bm{\nu}}$, $\forall \bm{\nu} \in \cF$, and that $\tilde{\bm{c}}$ depends only on $\bm{b}$ and $\varepsilon$.
Next, \cite[Lem.\ 7.23]{adcock2022sparse} and the assumption on $\bm{b}$ {give that $\tilde{\bm{c}} \in \ell^p_{\bm{u}}(\cF)$ with
\bes{
\nmu{\tilde{\bm{c}}}_{p,\bm{u}} \leq C,
}
for some constant $C = C(\bm{b},\varepsilon,p)$ depending on $\bm{b}$, $\varepsilon$ and $p$ only. Note that the statement of \cite[Lem.\ 7.23]{adcock2022sparse} asserts this bound for the true coefficient sequence $\bm{c}$. However, its proof uses the coefficient bound \eqref{coeff-bound}, and therefore the result holds for the sequence $\tilde{\bm{c}}$.}
 Let
\be{
\label{k-def}
k = \frac{m}{14 \log(m/\epsilon)}.
}
We now consider the weighted Stechkin inequality \eqref{weighted-Stechkin} with $q=1$ and $q = 2$. This implies that there are sets $S_1,S_2 \subset \cF$ with $|S_q|_{\bm{u}} \leq k$, $q = 1,2$, such that
\be{
{\nmu{\tilde{\bm{c}}-\tilde{\bm{c}}_{S_q}}_{q,\bm{u}} = \sigma_k(\tilde{\bm{c}})_{q,\bm{u}} \leq C k^{\frac1q-\frac1p}}.
}
Here $\tilde{\bm{c}}_{S_q}$ is the sequence with $\bm{\nu}$th term equal to $\tilde{c}_{\bm{\nu}}$ if $\bm{\nu} \in S_q$ and zero otherwise. Let $S = S_1 \cup S_2$. {Observe that $S_1$ and $S_2$, and therefore $S$, are independent of $f$ and depend only on $\bm{b}$, $\varepsilon$, since the sequence $\tilde{\bm{c}}$ depends only on these terms.} Parseval's identity{, the inequality $|c_{\bm{\nu}}| \leq \tilde{c}_{\bm{\nu}}$} and the fact that $\nm{\cdot}_{2,\bm{u}} = \nm{\cdot}_2$ give that
\be{
\label{fS-L2-err}
\nmu{f - f_S}_{L^2_{\varrho}(\cU)} {\leq \nmu{\tilde{\bm{c}} - \tilde{\bm{c}}_S}_{2,\bm{u}} \leq \nm{\tilde{\bm{c}} - \tilde{\bm{c}}_{S_2}}_{2,\bm{u}} } \leq C k^{\frac12-\frac1p}.
}
On the other hand, the definition \eqref{weights-def} of the weights $\bm{u}$ gives
\be{
\label{fS-Linf-err}
\nmu{f - f_S}_{L^{\infty}(\cU)} \leq \sum_{\bm{\nu} \notin S} | c_{\bm{\nu}} | \nm{\Psi_{\bm{\nu}}}_{L^{\infty}(\cU)} {\leq \nm{\tilde{\bm{c}} - \tilde{\bm{c}}_S}_{1,\bm{u}} \leq \nm{\tilde{\bm{c}} - \tilde{\bm{c}}_{S_1}}_{1,\bm{u}} }\leq C k^{1-\frac1p}.
}
Now consider the LS approximation in $\cP_{S}$. We seek to apply Theorem \ref{t:wLS-samp-comp} with $\epsilon$ replaced by $\epsilon/2$. {We consider MC samples, which corresponds to the case $\mu = \varrho$. Therefore, the weight function $w(\bm{y}) \equiv 1$ due to \eqref{mu-rho-w}.} In this case, \eqref{wLS-samp-comp} is equivalent to the condition
\be{
\label{main-samp-cond-1}
m \geq 7 \cdot \kappa(\cP_S) \cdot \log(4 |S|/\epsilon).
}
It follows from the definitions of $\cP_S$ and $\kappa(\cP_S)$ and \eqref{K-alternative-def} that
\bes{
\kappa(\cP_S) = \sup_{\bm{y} \in \cU} \sum_{\bm{\nu} \in S} | \Psi_{\bm{\nu}}(\bm{y}) |^2.
}
However, $\Psi_{\bm{\nu}}$ is a tensor-product Chebyshev or Legendre polynomials. Such polynomials achieve their absolute maximum value at the corner point $\bm{y} = \bm{1}$. Hence
\bes{
\kappa(\cP_S) = \sum_{\bm{\nu} \in S} | \Psi_{\bm{\nu}}(\bm{1}) |^2 = \sum_{\bm{\nu} \in S} \nm{\Psi_{\bm{\nu}}}^2_{L^{\infty}(\cU)} = \sum_{\bm{\nu} \in S} u^2_{\bm{\nu}} = |S|_{\bm{u}} \leq |S_1|_{\bm{u}} + |S_2|_{\bm{u}} \leq 2k.
}
Here, in the third equality we recall the definition \eqref{weights-def} of the weights $\bm{u}$. Hence \eqref{main-samp-cond-1} is implied by
\be{
\label{main-samp-cond-2}
m \geq 14 \cdot k \cdot \log(8 k/ \epsilon).
}
Here, we also used the observation that $|S| \leq |S|_{\bm{u}} \leq 2k$, since the weights $u_{\bm{\nu}} \geq 1$. Now, since $m \geq 3$ and $\epsilon < 1$, \eqref{k-def} gives that
\bes{
k \leq \frac{m}{14 \log(3)} \leq \frac{m}{8}.
}
Thus, \eqref{main-samp-cond-2}, and therefore \eqref{main-samp-cond-1} as well, is implied by
\bes{
m \geq 14 \cdot k \cdot \log(m/\epsilon).
}
However, this holds because of the choice of $k$ in \eqref{k-def}.

We are now ready to apply Theorem \ref{t:wLS-samp-comp}. This implies that, with probability at least $1-\epsilon/2$, the LS approximation $\hat{f}$ is unique and satisfies
\be{
\label{main-err-bd-1}
\nmu{f - \hat{f}}_{L^2_{\varrho}(\cU)} \leq \inf_{p \in \cP_S} \left \{ \nmu{f - p}_{L^2_{\varrho}(\cU)} + 2 \nmu{f-p}_{\mathsf{disc},1} \right \} + 2 \nmu{\bm{e}}_{\infty},
}
and moreover, the condition number $\mathrm{cond}(\bm{A}) \leq 2$. 

It remains to show that \eqref{main-err-bd-1} implies the desired error bound. Let $p = f_S \in \cP_S$. Then \eqref{main-err-bd-1} gives
\bes{
\nmu{f - \hat{f}}_{L^2_{\varrho}(\cU)} \leq \nmu{f - f_S}_{L^2_{\varrho}(\cU)} + 2 \nmu{f-f_S}_{\mathsf{disc},1}  + 2 \nmu{\bm{e}}_{\infty}.
}
The main challenge now is to bound the discrete error term. We do this by considering it as a sum of independent random variables, and then by using Bernstein's inequality. This argument is essentially the same as that given in \cite[Lem.\ 7.11]{adcock2022sparse}. We omit the details. 
 We conclude that
\be{
\label{main-disc-err-bd}
\nmu{f-f_S}_{\mathsf{disc},1} \leq \sqrt{2} \left ( \frac{\nmu{f-f_S}_{L^{\infty}(\cU)} }{\sqrt{k}} + \nmu{f-f_S}_{L^2_{\varrho}(\cU)} \right ),
}
with probability at least $1-\epsilon/2$, provided
\bes{
m \geq 2 k \cdot \log(4/\epsilon).
}
However, this holds by the definition \eqref{k-def} of $k$. Indeed, using \eqref{k-def} and the fact that $m \geq 3$, we see that
\eas{
2 k \cdot \log(4/\epsilon) = \frac{m (\log(4) + \log(\epsilon^{-1}) )}{7 ( \log(m) + \log(\epsilon^{-1}))} \leq m,
}
as required. 

Therefore, since \eqref{main-err-bd-1} holds with probability at least $1-\epsilon/2$ and \eqref{main-disc-err-bd} holds with probability at least $1-\epsilon/2$, the union bound now implies that
\bes{
\nmu{f - \hat{f}}_{L^2_{\varrho}(\cU)} \leq \nmu{f - f_S}_{L^2_{\varrho}(\cU)} + 2 \sqrt{2} \left ( \frac{\nmu{f-f_S}_{L^{\infty}(\cU)} }{\sqrt{k}} + \nmu{f-f_S}_{L^2_{\varrho}(\cU)} \right )  + 2 \nmu{\bm{e}}_{\infty},
}
holds with probability at least $1-\epsilon$. Finally, to complete the proof we use \eqref{k-def}, \eqref{fS-L2-err} and \eqref{fS-Linf-err}.
}

\subsection{Theorem \ref{t:main-res-2}}

{Theorem \ref{t:main-res-2} extends \cite[Thm.\ 3.7]{adcock2022efficient} in two main ways.} First, the latter assumes that $\bm{b} \in \ell^p(\bbN)$ is monotonically nonincreasing. {Theorem \ref{t:main-res-2} makes the weaker assumption that $\bm{b} \in \ell^p_{\mathsf{M}}(\bbN)$. For example, Theorem \ref{t:main-res-2} allows for a sequence of the form $\bm{b} = (1,0,1/2,1/4,\ldots)$, while \cite[Thm.\ 3.7]{adcock2022efficient} does not.} Second, \cite[Thm.\ 3.7]{adcock2022efficient} does not address the final part of Theorem \ref{t:main-res-2}, i.e., the statement that the same algebraic rate can be obtained using a polynomial approximation of smaller size. 

{Note that \cite[Thm.\ 3.7]{adcock2022efficient} considers the approximation of Hilbert-valued, holomorphic functions $f : \cU \rightarrow \cV$, where $\cV$ is a separable Hilbert space. This setting is relevant to parametric DE problems if one seeks to approximate the whole parametric solution map. To avoid unnecessary complications, in the current paper we consider only scalar-valued function approximation. However, all our main results -- and therefore, our conclusions about the efficacy of MC sampling -- extend to the Hilbert-valued setting with few additional technicalities.}

We now {establish Theorem \ref{t:main-res-2} by proving the two aforementioned extensions of \cite[Thm.\ 3.7]{adcock2022efficient}.}

\prf{[Proof of Theorem \ref{t:main-res-2}; assumption on $\bm{b}$]
Write $\tilde{\bm{b}}$ for the minimal monotone majorant of $\bm{b}$. Let $\bm{\rho} \in [1,\infty)^{\bbN}$ be such that 
\bes{
\sum^{\infty}_{i=1} \left ( \frac{\rho_i + \rho^{-1}_{i}}{2} -1 \right ) \tilde{b}_i \leq \varepsilon.
}
Then, since $\tilde{b}_i \geq b_i$ for every $i \in \bbN$, it follows that
\bes{
\sum^{\infty}_{i=1} \left ( \frac{\rho_i + \rho^{-1}_{i}}{2} -1 \right ) b_i \leq \varepsilon.
}
Hence
\eas{
\cR_{\tilde{\bm{b}},\varepsilon} & = \bigcup \left \{ \cE_{\bm{\rho}} : \bm{\rho} \in [1,\infty)^{\bbN},\ \sum^{\infty}_{i=1} \left ( \frac{\rho_i + \rho^{-1}_i}{2} - 1 \right ) \tilde{b}_i \leq \varepsilon \right \} 
\\
& \subseteq  \bigcup \left \{ \cE_{\bm{\rho}} : \bm{\rho} \in [1,\infty)^{\bbN},\ \sum^{\infty}_{i=1} \left ( \frac{\rho_i + \rho^{-1}_i}{2} - 1 \right ) b_i \leq \varepsilon \right \}  = \cR_{\bm{b},\varepsilon}
}
From this we deduce that $\cH(\bm{b},\varepsilon) \subseteq \cH(\tilde{\bm{b}},\varepsilon)$. {We are now in a position to apply \cite[Thm.\ 3.7]{adcock2022efficient} to the sequence $\tilde{\bm{b}}$, which is monotonically nonincreasing by construction and an element of $\ell^p(\bbN)$ by assumption. As noted, in this work we consider scalar-valued functions. Thus, we apply this result with $K = 1$, in which case the term the term $f - \cP_h(f)$ in \cite[Eqn.\ (3.16)]{adcock2022efficient} vanishes. Note that \cite[Thm.\ 3.7]{adcock2022efficient} only asserts the existence of a mapping with the desired error bound, but it is shown later in \cite[\S 4.2]{adcock2022efficient} that this mapping corresponds to minimizers of the weighted square-root LASSO problem with $\lambda = (4 \sqrt{m/L(m,\epsilon)})^{-1}$ (see \cite[Tab.\ 1]{adcock2022efficient}).} 
}

For the next part of the proof of Theorem \ref{t:main-res-2} we need the following lemma:

\lem{
Let $\bm{x},\bm{z} \in \bbC^N$ and $S \subseteq [N]$ be an index set containing the $n$ largest entries of $\bm{z}$ in absolute value. Then
\bes{
\nmu{\bm{x} - \bm{z}_S }_2 \leq 3 \nmu{\bm{x} - \bm{z}}_2 + 3 \sigma_n(\bm{x})_2,
}
where $\sigma_n(\bm{x})_2 = \min \{ \nmu{\bm{x} - \bm{z}}_2 : \bm{z} \in \bbC^N,\ | \mathrm{supp}(\bm{z}) | \leq n \}$.
}
\prf{
The proof of this lemma is based on the proof of \cite[Cor.\ 3.2]{needell2010signal}. Let $T \subseteq [N]$, $|T| \leq n$, be an index set containing the $n$ largest entries of $\bm{x}$ in absolute value, so that $\sigma_n(\bm{x})_2 = \nmu{\bm{x} - \bm{x}_T }_2$. Then
\eas{
\nmu{\bm{x} - \bm{z}_S}_{2} &\leq \nm{ \bm{x}_T - \bm{z}_S}_2 + \nm{\bm{x}_{T^c}}_2
\\
& \leq  \nm{( \bm{x}_T - \bm{z}_S )_S}_2  +  \nmu{\bm{x}_{T \backslash S}}_2 + \sigma_n(\bm{x})_2
 \\
 & = \nm{(\bm{x}_T - \bm{z})_S}_2 +  \nmu{\bm{x}_{T \backslash S}}_2 + \sigma_n(\bm{x})_2
 \\
 & \leq \nm{\bm{x}_T - \bm{z}}_2 + \nmu{\bm{x}_{T \backslash S}}_2 + \sigma_n(\bm{x})_2
 \\
 & \leq \nm{\bm{x} - \bm{z} }_2 + \nmu{\bm{x}_{T \backslash S}}_2 + 2\sigma_n(\bm{x})_2
}
Now consider the second term. Write
\bes{
\nmu{\bm{x}_{T \backslash S}}_2  \leq \nmu{(\bm{x}-\bm{z})_{T \backslash S}}_2 + \nmu{\bm{z}_{T \backslash S}}_2 \leq \nmu{\bm{x} - \bm{z}}_2 + \nmu{\bm{z}_{T \backslash S}}_2.
}
Now observe that $|S| = |T| = n$ and therefore $|S \backslash T| = | T \backslash S |$. Since $S$ contains the largest $n$ entries of $\bm{z}$ in absolute value, we must have
\bes{
\nmu{\bm{z}_{T \backslash S}}_2 \leq \nmu{\bm{z}_{S \backslash T}}_2.
}
We can bound this as follows:
\bes{
\nmu{\bm{z}_{S \backslash T}}_2 = \nmu{(\bm{z}-\bm{x}_T)_{S \backslash T}}_2 \leq \nmu{\bm{z} - \bm{x}_T}_{2} \leq \nmu{\bm{z} -\bm{x} }_2 + \sigma_n(\bm{x})_2.
}
Therefore
\bes{
\nmu{\bm{x}_{T \backslash S}}_2  \leq 2 \nmu{\bm{x} - \bm{z}}_2 + \sigma_n(\bm{x})_2.
}
Combining this with the earlier inequality, we deduce the result.
}

\prf{[Proof of Theorem \ref{t:main-res-2}; construction of the set $S$]
Let $n = \lceil m / L(m,\epsilon) \rceil$ and $\hat{\bm{c}}$ be as in \eqref{SRLASSO-minimizer}. Let $S \subseteq \Lambda$ be the index set of the largest $\min \{ n , N \}$ entries of $\hat{\bm{c}}$ in absolute value, where $N = |\Lambda|$. Then, by Parseval's identity and the above lemma,
\eas{
\nmu{f - \hat{f}_S}_{L^2_{\varrho}(\cU)} & \leq \nmu{f - f_{\Lambda}}_{L^2_{\varrho}(\cU)} + \nmu{\bm{c}_{\Lambda} - \hat{\bm{c}}_S }_{2}
\\
& \leq \nmu{f - f_{\Lambda}}_{L^2_{\varrho}(\cU)} + 3 \nmu{\bm{c}_{\Lambda} - \hat{\bm{c}}}_2 + 3 \sigma_n(\bm{c}_{\Lambda})_2
\\
& \leq 4 \nmu{f - f_{\Lambda}}_{L^2_{\varrho}(\cU)} + 3 \nmu{f - \hat{f}}_{L^2_{\varrho}(\cU)} + 3 \sigma_n(\bm{c})_2.
}
Due to the first part of the theorem, the second term satisfies
\bes{
\nmu{f - \hat{f}}_{L^2_{\varrho}(\cU)} \leq C \cdot \left ( \frac{m}{L(m,\epsilon)} \right )^{\frac12 - \frac1p} + c \cdot \nm{\bm{e}}_{\infty} .
}
Moreover, an inspection of the proof of \cite[Thm.\ 3.7]{adcock2022efficient} shows that the first term satisfies
\bes{
 \nmu{f - f_{\Lambda}}_{L^2_{\varrho}(\cU)} \leq C \cdot \left ( \frac{m}{L(m,\epsilon)} \right )^{\frac12 - \frac1p}
 }
(this is simply due to the fact that $\Lambda$ contains all anchored sets by construction and Theorem \ref{t:best-s-term-anchored}). Next, for the term $\sigma_n(\bm{c})_2$ we use Parseval's identity and Theorem \ref{t:best-s-term} (along with the observation that the constant $C$ in this theorem is not larger than the constant $C$ in any of the above estimates -- a fact that is implicit in their proofs) to get
\bes{
\sigma_n(\bm{c})_2 = \nmu{f - f_n}_{L^2_{\varrho}(\cU)} \leq C \cdot n^{\frac12-\frac1p} \leq C \cdot \left ( \frac{m}{L(m,\epsilon)} \right )^{\frac12-\frac1p} \leq C \cdot \left ( \frac{m}{L(m,\epsilon)} \right )^{\frac12-\frac1p} .
}
Here, in the final step, we used the fact that $c_1 \geq 1$ by assumption. Combining this with the previous bounds yields
\eas{
\nmu{f - \hat{f}_S}_{L^2_{\varrho}(\cU)} \leq 10 \cdot C \cdot \left ( \frac{m}{L(m,\epsilon)} \right )^{\frac12 - \frac1p}+ 3  \cdot  c_2 \cdot \nmu{\bm{e}}_{\infty},
}
as required.
}

\section{Additional information on experiments}\label{s:additional-info}

In this section, we provide additional information on the various numerical experiments. MATLAB code reproducing all the experiments is available at \url{https://github.com/benadcock/is-MC-bad}.

\subsection{Statistical simulations and plotting}\label{ss:statistical-sims}

In our experiments, we analyze the behaviour of a statistical quantity, typically the error of an approximation or condition number of a LS matrix, as a certain $x$-variable, e.g., the number of measurements $m$ of approximation space size $n$, varies. We do this by computing multiple trials, thus generating a set of data of the form
\bes{
 (x_i , y^{(t)}_i ),\quad i = 1,\ldots,i_{\max}, t = 1,\ldots,T.
}
Here $T$ is the number of trials, which is taken to be $T = 50$ in this work. In our figures we always use a base-10 logarithmic scale on the $y$-axis. Following \cite[\S A.1.3]{adcock2022sparse}, we visualize this data by computing its geometric mean and (corrected) geometric standard deviation. The main curve in our figures shows the mean and the shaded region shows one standard deviation. To be precise, the former is obtained by plotting
\bes{
(x_i , 10^{\mu_i} ),\quad i = 1,\ldots,i_{\max},\qquad \text{where } \mu_i = \frac{1}{T} \sum_{t =1}^{T} \log_{10}(y_i^{(t)}),
}
and the latter region is that enclosed by the upper and lower curves defined by the data
\bes{
(x_i , 10^{\mu_i \pm \sigma_i} ),\quad i = 1,\ldots,i_{\max},\qquad \text{where } \sigma_i =  \sqrt{\frac{1}{T - 1}
\sum_{t =1}^{T} (\log_{10}(y_i^{(t)}) - \mu_i)^2}.
}
See \cite[\S A.1.3]{adcock2022sparse} for a justification of this choice.

\subsection{Finite grids and discrete sampling measures}\label{ss:fine-grids-measures}

In this work, we perform computations over a finite grid of $K = 100,000$ MC points $ \{ \bm{z}_i \}^{M}_{i=1} \subset \cU$ drawn i.i.d.\ from the underlying continuous measure $\varrho$. This grid is drawn once before any subsequent computations are performed. We use this grid in two ways. First, to compute the error, and second to draw sample points according to the various sampling schemes. The latter is based on \cite{adcock2020nearoptimal,migliorati2021multivariate}. 

Mathematically, this means that we simply replace the continuous measure $\varrho$ by the discrete measure $\tau$ defined by
\bes{
\D \tau(\bm{y}) = \frac1K \sum^{K}_{i=1} \delta(\bm{y} - \bm{z}_i) \D \varrho(\bm{y}).
}
In particular, all analysis immediately applies upon replacing $\varrho$ by $\tau$ in the various theoretical results.

Concretely, we define the (relative) discrete $L^2_{\varrho}$-norm error of an approximation $\hat{f}$ to a function $f$ as
\be{
\label{f-error-fun}
E(f) : = \frac{\sqrt{\frac{1}{K} \sum^{K}_{i=1} | f(\bm{z}_i) - \hat{f}(\bm{z}_i) |^2}}{\sqrt{\frac{1}{K} \sum^{K}_{i=1} | f(\bm{z}_i) |^2}}.
}
Note that this is precisely
\bes{
E(f) = \frac{\nmu{f - \hat{f}}_{L^2_{\tau}(\cU)}}{\nmu{f}_{L^2_{\tau}(\cU)}}.
}
Throughout our experiments, we compare MC sampling and the near-optimal sampling strategy, both with respect to the measure $\tau$. The former is done straightforwardly. To draw a point $\bm{y} \sim \tau$, we simply draw an integer $i$ randomly and uniformly from the set $[K]$, and then set $\bm{y} = \bm{z}_i$.

We now describe how to sample from the near-optimal sampling measure of \S \ref{ss:opt-samp-scheme}.
Let $\cP \subset L^2_{\tau}(\cU)$ be an $n$-dimensional subspace. When $\varrho$ is replaced by $\tau$, the sampling measure \eqref{opt-samp-meas} takes the form
\bes{
\D \mu(\bm{y}) = \frac1n \cK(\cP)(\bm{y}) \D \tau(\bm{y}) = \frac{1}{n K} \sum^{K}_{i=1} \delta(\bm{y}-\bm{z}_i) \cK(\cP)(\bm{z}_i) \D \varrho(\bm{y}).
}
Hence $\mu$ is a discrete measure, supported over the grid $\{ \bm{z}_i \}^{K}_{i=1}$ with $\bm{y} \sim \mu$ if
\be{
\label{prob-opt-disc}
\bbP(\bm{y} = \bm{z}_i) = \frac{1}{n K} \cK(\cP)(\bm{z}_i),\quad i \in [K].
}
Thus, drawing samples from $\mu$ is straightforward, once the values $\cK(\cP)(\bm{z}_i)$ have been computed. We do this via \eqref{K-alternative-def}, which first requires computing an orthonormal basis of $\cP$ in $L^2_{\tau}(\cU)$. Let $\{ \Phi_1,\ldots,\Phi_n \}$ be a possibly nonorthonormal basis of $\cP$. In our experiments, when $\cP = \cP_S$, we simply let this be the corresponding set of Chebyshev or Legendre polynomials with multi-indices in $S$. Following \cite{adcock2020nearoptimal,migliorati2021multivariate}, we construct an orthonormal basis via QR decomposition. Let
\bes{
\bm{B} = \frac{1}{\sqrt{K}} \left ( \Phi_j(\bm{z}_i) \right )^{K,n}_{i,j = 1} \in \bbC^{K \times n},
}
and write $\bm{B} = \bm{Q} \bm{R}$, where $\bm{Q} = (q_{ij})^{K,n}_{i,j=1} \in \bbC^{K \times n}$ has orthonormal columns. These columns yield an orthonormal basis $\{ \Psi_1,\ldots, \Psi_n \}$ of $\cP$ with respect to $\tau$, given by
\bes{
\Psi_j(\bm{z}_i) = \sqrt{K} q_{ij},\quad \forall i \in [K],\ j \in [n].
}
Hence \eqref{K-alternative-def} gives
\be{
\label{K-disc}
\cK(\cP)(\bm{z}_i) = K \sum^{n}_{j=1} |q_{ij} |^2,\quad \forall i \in [K].
}
Thus, we now define the discrete probability distribution $\bm{\pi} = (\pi_i)^{K}_{i=1}$ over $[K]$ as
\bes{
\pi_i = \frac1n \sum^{n}_{j=1} |q_{ij} |^2,\quad \forall i \in [K].
}
It follows immediately from \eqref{prob-opt-disc} that drawing a point $\bm{y} \sim \mu$ is equivalent to drawing an integer $i \sim \bm{\pi}$ and then setting $\bm{y} = \bm{z}_i$. 

Finally, in order to construct the weighted LS approximation, we need the values of the weight function $w$ over the grid. Using \eqref{opt-w} and \eqref{K-disc}, we see that
\bes{
w(\bm{z}_i) = \left ( \frac{K}{n} \sum^{n}_{j=1} |q_{ij} |^2 \right )^{-1} = \left ( K \pi_i \right )^{-1},\qquad \forall i \in [K].
}

\subsection{Adaptive least-squares approximation}\label{ss:adaptive-LS}

\begin{figure}[t!]
\begin{center}
\begin{small}
 \begin{tabular}{@{\hspace{0pt}}c@{\hspace{4pt}}c@{\hspace{0pt}}}
\includegraphics[width = 0.45\textwidth]{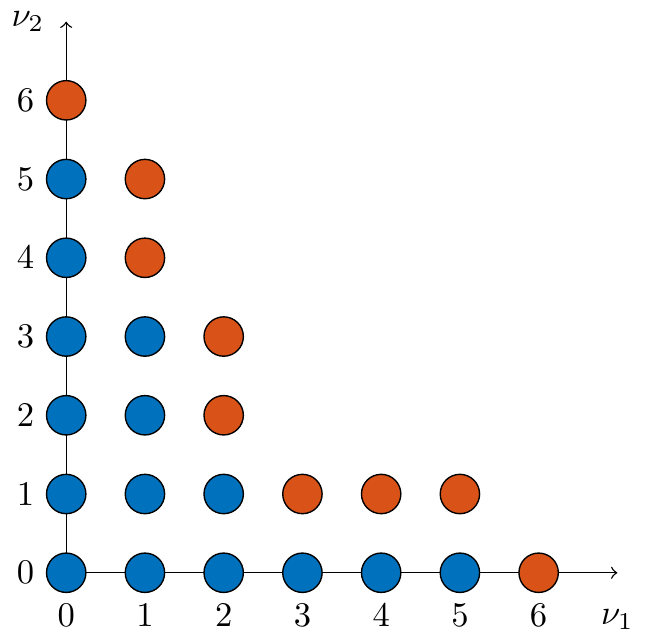}&
\includegraphics[width = 0.45\textwidth]{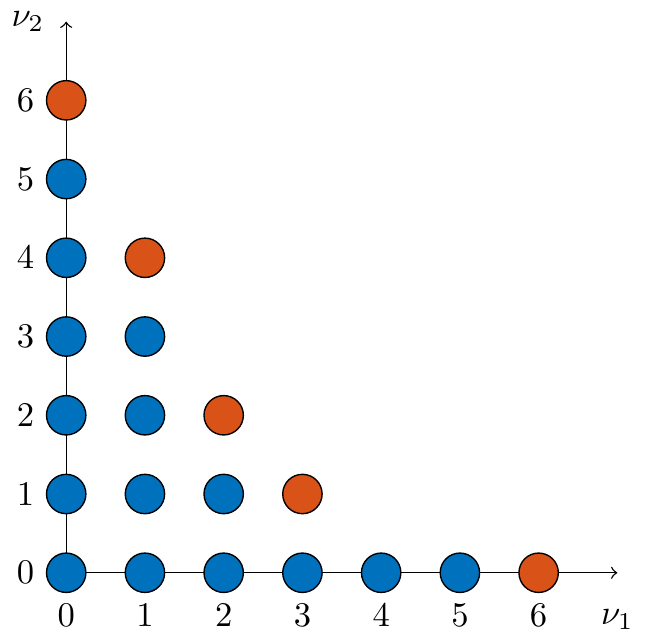}
\\[\errplottextsp]
margin & reduced margin
\end{tabular}
\end{small}
\end{center}
\caption{{The margin $\cM(S)$ (left) and reduced margin $\cR(S)$ (right) for a typical index set $S$ in two dimensions. The blue dots are the multi-indices comprising $S$ and the red dots are multi-indices in the margin or reduced margin.}} 
\label{fig:margin}
\end{figure}

We now describe the ALS approximation. This follows \cite{migliorati2019adaptive} (see also \cite{migliorati2015adaptive}), which is based on techniques for the constructing adaptive sparse grid quadratures \cite{gerstner2003dimension}. It requires several definitions. First, we define the \textit{margin} of a multi-index set $S \subseteq \bbN^d_0$ as
\bes{
\cM(S) = \left \{ \bm{\nu} \in \cF \backslash S : \exists j \in [d] : \bm{\nu} - \bm{e}_j \in S \right \}.
}
Here, as before, $\bm{e}_j \in \bbN^d_0$ is the $j$th canonical multi-index, with value $1$ in its $j$th entry and zero otherwise. Next, we define the \textit{reduced margin} of $S$ as the set
\bes{
\cR(S) = \left \{ \bm{\nu} = (\nu_j)^{d}_{j=1} \in \cM(S) : \forall j \in [d], \nu_j \neq0\ \Rightarrow\ \bm{\nu} - \bm{e}_j \in S \right \}.
}
Note that if $S$ is a lower set, then $S \cup \{ \bm{\nu} \}$ is also lower for any $\bm{\nu} \in \cR(S)$. {Fig.\ \ref{fig:margin} illustrates the concepts of margin and reduced margin for a typical multi-index set in two dimensions.}

Given a set $S \subset \bbN^d_0$, the ALS approximation constructs a new set by selecting multi-indices from its reduced margin $\cR(S)$. It does this via a \textit{bulk chasing procedure}. This procedure is defined as follows. Let $e : \cR(S) \rightarrow \bbR$ and $0 < \beta \leq 1$ be a parameter. The function $e$ serves as an estimate for the true coefficients of the function being approximated. We describe a precise choice of $e$ below. With this, the procedure $\texttt{bulk} = \texttt{bulk}(\cR(S) , e , \beta )$ computes a set $T \subseteq \cR(S)$ of minimal positive cardinality such that
\bes{
\sum_{\bm{\nu} \in T} e(\bm{\nu}) \geq \beta \sum_{\bm{\nu} \in \cR(S)} e(\bm{\nu}).
}
Following \cite{migliorati2019adaptive}, we use the value $\beta = 0.5$ in our experiments.

We are now ready to define the ALS approximation. Fix a positive and finite almost everywhere weight function $w : \cU \rightarrow [0,\infty)$. Typically, and throughout this paper, we initialize the procedure with the index set $S^{(1)} = \{ \bm{0} \}$. Now suppose at step $l$ we have computed the index set $S^{(l)}$, sample points $\{ \bm{y}^{(l)}_i \}^{m^{(l)}}_{i=1}$ and the corresponding (weighted) LS approximation $\hat{f}^{(l)} \in \cP_{S^{(l)}}$ as in \eqref{wLS-def}. Let
\bes{
\ip{f}{g}^{(l)}_{\mathsf{disc},w} = \frac1m \sum^{m}_{i=1} w(\bm{y}^{(l)}_i) f(\bm{y}^{(l)}_i) \overline{g(\bm{y}^{(l)}_i)},\quad \forall f,g \in L^2_{\varrho}(\cU) \cap C(\cU).
}
Then we define the estimator $e^{(l)} : \cR(S^{(l)}) \rightarrow \bbR$ as
\be{
\label{bulk-estimator}
e^{(l)}(\bm{\nu}) = | \ip{f - \hat{f}^{(l)}}{\Psi_{\bm{\nu}}}^{(l)}_{\mathsf{disc},w} |^2,\quad \forall \bm{\nu} \in \cR(S^{(l)}).
}
Note that computing $e^{(l)}$ involves simple linear algebra. Indeed, let $\bm{A}^{(l)}$ and $\bm{b}^{(l)}$ be the measurement matrix and vector of the LS problem \eqref{LS-meas-mat-vec}, and let $\hat{\bm{c}}^{(l)}$ be the coefficients of the LS approximation $\hat{f}^{(l)}$. Then, for every $\bm{\nu} \in \cR(S^{(l)})$,
\bes{
e^{(l)}(\bm{\nu}) = (\bm{b}^{(l)} - \bm{A}^{(l)} \hat{\bm{c}}^{(l)} )^{\top} \bm{d},\qquad \text{where }\bm{d} = \left ( \frac{\sqrt{w(\bm{y}^{(l)}_i )}}{\sqrt{m}} \Psi_{\bm{\nu}}(\bm{y}^{(l)}_i) \right )^{m^{(l)}}_{i=1}.
}
Having computed $e^{(l)}$, we then compute a set $T^{(l)} = \texttt{bulk}(\cR(S^{(l)}) , e^{(l)} , \beta)$ and define the new index set $S^{(l+1)} = S^{(l)} \cup T^{(l)}$. Note that $S^{(l+1)}$ is a lower set whenever $S^{(l)}$ is, due to the properties of the reduced margin. In particular, when the initial index set is chosen as $S^{(1)} = \{ \bm{0} \}$ (or any lower set), then the sequence $S^{(1)},S^{(2)},\ldots$ is nested and lower.

We summarize the ALS approximation in Algorithm \ref{a:ALS}. Note that in our experiments we draw a new set of $m^{(l)}$ samples at step $l$. In practice, this is wasteful, since the  function evaluations computed at step $l$ are not re-used in subsequent steps. It is often important in practice to re-use sample points, so as to not waste function evaluations. This is done by creating hierarchical sequences of sample points
\bes{
\{ \bm{y}^{(1)}_i \}^{m^{(1)}}_{i=1} \subseteq \{ \bm{y}^{(2)}_i \}^{m^{(2)}}_{i=1} \subseteq \cdots.
}
MC sampling can be easily converted to a hierarchical sampling scheme simply by drawing only $m^{(l)} - m^{(l-1)}$ new points at step $l$. It is not as straightforward to do this with the near-optimal sampling strategy. However, there are several ways to construct optimal hierarchical sampling schemes with similar near-optimal sample complexity guarantees. See \cite{arras2019sequential,migliorati2019adaptive}. In this paper we do not consider such hierarchical schemes, since we are primarily concerned with the fundamental approximation properties of the ALS approximation.

\begin{algorithm}[t]
\caption{Adaptive (weighted) least-squares approximation}
\label{a:ALS}
\begin{algorithmic}
\REQUIRE{Orthonormal basis $\{ \Psi_{\bm{\nu}} \}_{\bm{\nu} \in \bbN^d_0}$, weight function $w : \cU \rightarrow (0,\infty)$, function to approximate $f \in L^2_{\varrho}(\cU)$, bulk chasing parameter $0 < \beta \leq 1$}
\ENSURE{Sequence of lower sets $S^{(1)},S^{(2)},\ldots$ and approximations $\hat{f}^{(1)},\hat{f}^{(2)},\ldots$}
\STATE{Set $S^{(1)} = \{ \bm{0} \}$}
\FOR{$l = 1,2,\ldots$}
\STATE{Set $n^{(l)} = |S^{(l)} |$}
\STATE{Determine a number of samples $m^{(l)} \geq n^{(l)}$}
\STATE{Compute sample points $\{ \bm{y}^{(l)}_i \}^{m^{(l)}}_{i=1}$}
\STATE{Compute the (weighted) LS approximation $\hat{f} = \hat{f}^{(l)}$ via \eqref{wLS-def} with $S = S^{(l)}$ and $\bm{y}_i = \bm{y}^{(l)}_i$}
\STATE{Compute the estimator $e^{(l)}$ via \eqref{bulk-estimator}}
\STATE{Compute $T^{(l)} = \texttt{bulk}(\cR(S^{(l)}) , e^{(l)} , \beta )$}
\STATE{Set $S^{(l+1)} = S^{(l)} \cup T^{(l)}$}
\ENDFOR
\end{algorithmic}
\end{algorithm}

\subsection{Computing the CS approximation}\label{ss:efficient-PDI}

Computing the CS approximation via \eqref{SRLASSO-minimizer} involves (approximately) solving the weighted SR-LASSO problem
\be{
\label{wrs-lasso-sm}
\min_{\bm{z} \in \bbC^{N}} \cG(\bm{z}),\quad \text{where }\cG(\bm{z}) = \lambda \nmu{\bm{z}}_{1,\bm{w}} + \nmu{\bm{A} \bm{z} - \bm{f}}_2.
}
This {is a convex, nonsmooth optimization problem, widely employed in statistics since its inception \cite{belloni2011square}. It is also a powerful tool for high-dimensional polynomial approximation thanks to its stability with respect to the tuning parameter $\lambda$ when dealing with unknown noise corrupting the samples \cite{adcock2019correcting}}.
There are various algorithms and software packages for solving this problem (see, e.g., \cite[\S A.2]{adcock2022sparse}). However, a challenge is that most of these cannot deliver high accuracy in a computationally-efficient way \cite{adcock2022efficient}. In this work, we use the approach developed in \cite{adcock2022efficient}, which is based on the previous works \cite{colbrook2022warpd,colbrook2022can}. This uses a combination of the \textit{primal-dual iteration} \cite{chambolle2011first,chambolle2016ergodic} with a \textit{restarting strategy} \cite{renegar2022simple,roulet2020computational,roulet2020sharpness}.

The unrestarted primal-dual iteration for \eqref{wrs-lasso-sm} is given in Algorithm \ref{a:primal-dual}. Here, $\sgn(\cdot)$ denotes the sign of a complex number.
We omit the details of the derivation of this algorithm, and refer the reader instead to \cite[\S 4.3 \& Alg.\ 1]{adcock2022efficient}. Note that we deviate slightly from the setup described therein. First, we consider only the scalar-valued case in \eqref{wrs-lasso-sm}, whereas \cite{adcock2022efficient} considers vectors whose entries take values in a finite-dimensional subspace $\cV_h$ of a Hilbert space $\cV$. Second, we do not consider the so-called \textit{ergodic} sequences in this work.

\begin{algorithm}[t]
\caption{The unrestarted primal-dual iteration for \eqref{wrs-lasso-sm}}
\label{a:primal-dual}
\begin{algorithmic}
\REQUIRE{measurement matrix $\bm{A} \in \bbC^{m \times N}$, measurements $\bm{f} \in \bbC^m$, positive weights $\bm{w} = (w_i)^{N}_{i=1}$, parameter $\lambda > 0$, stepsizes $\tau,\sigma > 0$, number of iterations $T \geq 1$, initial values $\bm{c}^{(0)} \in \bbC^N$, $\bm{\xi}^{(0)} \in \bbC^m$}
\ENSURE{$\bm{c}^{(T)} = \text{\texttt{primal-dual-wSRLASSO}}(\bm{A},\bm{f},\bm{w},\lambda,\tau,\sigma,T,\bm{c}^{(0)},\bm{\xi}^{(0)})$, an approximate minimizer of \eqref{wrs-lasso-sm}}
%\STATE{}
\FOR{$n = 0,1,\ldots,T-1$}
\STATE{$\bm{p} = (p_i)^{N}_{i=1} = \bm{c}^{(n)} - \tau \bm{A}^* \bm{\xi}^{(n)}$}
\STATE{$\bm{c}^{(n+1)} = \left ( \max \{ | p_i | - \tau \lambda w_i , 0 \} \sgn(p_i) \right )^{N}_{i=1}$}
\STATE{$\bm{q} = \bm{\xi}^{(n)} + \sigma \bm{A} (2 \bm{c}^{(n+1)} - \bm{c}^{(n)}) - \sigma \bm{f}$}
\STATE{$\bm{\xi}^{(n+1)} = \min \{ 1 , 1/\nm{\bm{q}}_2 \} \bm{q}$}
\ENDFOR
\end{algorithmic}
\end{algorithm}

For suitable choices of the stepsize parameters $\tau,\sigma$, the values $\cG(\bm{c}^{(n)})$ converge with rate $\ordu{1/n}$ to the minimum value of the objective function of \eqref{wrs-lasso-sm}. One can obtain faster convergence using a restarting strategy. Restarting is general concept in continuous optimization \cite{renegar2022simple,roulet2020computational,roulet2020sharpness}. The idea is to run an algorithm repeatedly for a fixed, but small number of iterations, and feed the output of one instance into the input of the next instance (a `restart') while also adjusting certain parameters. In our case, the particular restarting scheme we use is described in Algorithm \ref{a:primal-dual-restart}. See  \cite[\S 4.6 \& Alg.\ 4]{adcock2022efficient} for details.

\begin{algorithm}[t]
\caption{The restarted primal-dual iteration for \eqref{wrs-lasso-sm}}
\label{a:primal-dual-restart}
\begin{algorithmic}
\REQUIRE{measurement matrix $\bm{A} \in \bbC^{m \times N}$, measurements $\bm{f} \in \bbC^m$, positive weights $\bm{w} = (w_i)^{N}_{i=1}$, parameter $\lambda > 0$, stepsizes $\tau,\sigma > 0$, number of iterations $T \geq 1$, number of restarts $R \geq 1$, tolerance $\zeta' > 0$, scale parameter $0 < r < 1$, constant $s > 0$, initial values $\bm{c}^{(0)} \in \bbC^N$, $\bm{\xi}^{(0)} \in \bbC^m$}
\ENSURE{$\tilde{\bm{c}}^{(R)} =  \text{\texttt{primal-dual-rst-wSRLASSO}}(\bm{A},\bm{f},\bm{w},\lambda,\tau,\sigma,T,R,\zeta',r,s)$, an approximate minimizer of \eqref{wrs-lasso-sm}}
\STATE{Initialize $\tilde{\bm{c}}^{(0)} = \bm{0} \in \bbC^N$, $\varepsilon_0 = \nm{\bm{f}}_2$}
\FOR{$l = 0,1,\ldots,R-1$}
\STATE{$\varepsilon_{l+1} = r (\varepsilon_l + \zeta')$}
\STATE{$a_l = s \varepsilon_{l+1}$}
\STATE{$\tilde{\bm{c}}^{(l+1)} = a_l \cdot \text{\texttt{primal-dual-wSRLASSO}}(\bm{A},\bm{f}/a_l,\bm{w},\lambda,\tau,\sigma,T,\tilde{\bm{c}}^{(l)} / a_l,\bm{0})$}
\ENDFOR
\end{algorithmic}
\end{algorithm}

For suitable choices of parameters, the values $\cG(\tilde{\bm{c}}^{(l)})$ converge with linear order to the minimum value of \eqref{wrs-lasso-sm}, up to the tolerance $\zeta'$ -- specifically, the error is proportional to
\be{
\label{primal-dual-restarts-bound}
r^{l} + \zeta'.
}
Thus, to achieve an error of size $\eta \geq \zeta'$ we need $l \propto \log(\zeta')/\log(r)$ restarts. We omit the details and refer to \cite{adcock2022efficient} for further information. Like in \cite[\S 5.2]{adcock2022efficient} we make one modification to Algorithm \ref{a:primal-dual-restart}, which is to use the stopping criterion
\be{
\label{stopping-criterion}
l \geq R\quad \text{or}\quad \nmu{\tilde{\bm{c}}^{(l)} - \tilde{\bm{c}}^{(l-1)}}_2 \leq 10 \cdot \zeta'
}
(we use the slightly more conservative value of $10$ instead of $5$ in this work since it works more reliably in the high-accuracy regime). Here $\zeta'$ is the tolerance, which determines the limiting accuracy of the iterates. In this work, since we aim for high accuracy, we set $\zeta' = 10^{-15}$.

We are now ready to describe how we setup and compute the CS approximation  \eqref{SRLASSO-minimizer} in the experiments in the paper. This is summarized in Table \ref{tab:foo}. These parameter values arise from theoretical bounds for the restarted scheme. See \cite[\S 5.1.1]{adcock2022efficient} for details. We stress at this point that with these parameter values the CS approximation is obtained robustly and reliably, and is able to achieve high accuracy, as can be seen in the various figures in this paper.

\begin{table}
\footnotesize
\caption{The parameters used in the experiments in this paper for computing \eqref{SRLASSO-minimizer}. Specifically, we set $\hat{\bm{c}} = \text{\texttt{primal-dual-rst-wSRLASSO}}(\bm{A},\bm{f},\bm{w},\lambda,\tau,\sigma,T,R,\zeta',r,s)$ with the following values.}\label{tab:foo}
\begin{center}
\begin{tabular}{|c|l|}\hline
$\Lambda$ & As in \eqref{HCI-def}, $n$ chosen as the largest integer such that $|\Lambda^{\mathsf{HCI}}_n | \leq 10,000$
\\
$\bm{A}$ & As in \eqref{Ab-CS-def} (computed via the three-term recursion relation for orthogonal polynomials)
\\
$\bm{f}$ & As in \eqref{approx-LS-CS} with $e_i = 0$, $\forall i \in [m]$
\\
$\bm{w}$ & Equal to $\bm{u}$, as in \eqref{weights-def}
\\
$\lambda$ & Equal to $(5 \sqrt{m})^{-1}$, as in \cite[\S A.2.5]{adcock2022sparse}
\\
$\tau,\sigma$ & Equal to $1/\nm{\bm{A}}_2$
\\
$T$ & Equal to $\lceil 4 \nm{\bm{A}}_2 / r \rceil$
\\
$R$ & Equal to $100$
\\
$\zeta'$ & Equal to $10^{-15}$
\\
$r$ & Equal to $\E^{-1}$
\\
$s$ & Equal to $T / (2 \nmu{\bm{A}}_2)$
\\
\hline
\end{tabular}
\end{center}
\end{table}

Finally, we note that the computational cost of Algorithm \ref{a:primal-dual-restart} is proportional to
\bes{
m \cdot N \cdot T \cdot R = m \cdot N \cdot \lceil 4 \E \nm{\bm{A}}_2  \rceil \cdot R.
}
Due to \eqref{primal-dual-restarts-bound} and the parameter choices in Table \ref{tab:foo}, the number of restarts used should be roughly $\log(\zeta') / \log(r) = -\log(\zeta') \approx 35$. In our experiments, the algorithm with stopping criterion \eqref{stopping-criterion} always halts after roughly this number of restarts.

\section{Background on parametric diffusion equations}\label{s:pDE-background}

As noted, parametric DEs were the original motivations for recent developments in polynomial approximations to high- and infinite-dimensional functions. We now provide a little more detail on {certain types of these} problems. For further information {and discussion of other holomorphy for various other families parametric DEs}, see \cite{cohen2015approximation} or \cite[Chpt.\ 4]{adcock2022sparse} and references therein.

\subsection{Parametric diffusion equations}

A typical parametric PDE is the stationary diffusion equation with parametrized diffusion coefficient and homogeneous Dirichlet boundary conditions:
\begin{align}
\begin{split}
\label{eq:fruitfly}
-\nabla_{\bm{x}} \cdot (a(\bm{x},\bm{y}) \nabla_{\bm{x}} u(\bm{x}, \bm{y})) = F(\bm{x}), &\quad  \bm{x} \in \Omega,\\
u(\bm{x},\bm{y}) = 0, &\quad \bm{x} \in \partial \Omega.
\end{split}
\end{align}
Here $\bm{x} \in \Omega \subset \mathbb{R}^k$ is the spatial variable and domain, where $k\in \mathbb{N}$, $\Omega$ has Lipschitz boundary and $\nabla_{\bm{x}}$ denotes the gradient with respect to $\bm{x}$. The function $F$ is the forcing term, and is assumed to be nonparametric.

We assume that \eqref{eq:fruitfly} satisfies a \textit{uniform ellipticity condition}: namely, 
\be{
\label{uniform-ellipticity}
 \essinf_{\bm{x} \in \Omega} a(\bm{x},\bm{y}) \geq r,\quad \forall \bm{y} \in \cU,
}
for some $r > 0$. Now let $\cV = H^1_0(\Omega)$ be the standard Sobolev space of functions with weak first-order derivatives in $L^2(\Omega)$ and traces vanishing on $\partial \Omega$. Then the problem \eqref{eq:fruitfly} has a well-defined \textit{parametric solution map}
\be{
\label{u-soln-map}
u : \cU \rightarrow \cV,\ \bm{y} \mapsto u(\cdot,\bm{y}).
}
Notice that this map is Hilbert-valued, since it takes values in the Hilbert space $\cV$. It is therefore common to study a scalar \textit{quantity of interest} of it. Given a mapping $\cQ : \cV \rightarrow \bbC$, this is defined as
\be{
\label{QoI}
f(\bm{y}) = \cQ(u(\cdot,\bm{y})).
}
A typical example involves simply evaluating $u$ at some point $\bm{x}_0 \in \Omega$ in the spatial domain, i.e., $f(\bm{y}) = u(\bm{x}_0,\bm{y})$.

\subsection{\boldmath $(b,\varepsilon)$-holomorphy of parametric elliptic DEs}\label{s:pDE-holomorphy}

Consider a parametric diffusion equation \eqref{eq:fruitfly} with affine diffusion coefficient
\be{
\label{eq:affine-diffusion}
a(\bm{x},\bm{y}) = a_0(\bm{x}) + \sum^{\infty}_{j=1} y_j \psi_j(\bm{x}),
}
where the functions $a_0, \psi_1,\psi_2,\ldots \in L^{\infty}(\Omega)$. Notice that uniform ellipticity \eqref{uniform-ellipticity} for this problem is equivalent to the condition
\bes{
\sum^{\infty}_{j=1} | \psi_j(\bm{x}) | \leq a_0(\bm{x}) - r,\quad \forall \bm{x} \in \Omega.
}
Under this assumption, the parametric solution map \eqref{u-soln-map} is well defined. Furthermore it -- and by extension, quantities of interest of it such as \eqref{QoI}  -- are holomorphic functions of $\bm{y}$. Specifically, let
\be{
\bm{b} = (b_i)^{\infty}_{i=1},\quad \text{with }b_i = \nmu{\psi_i}_{L^{\infty}(\Omega)},\ \forall i \in \bbN.
}
Then it can be shown that the parametric solution map is $(\bm{b},\varepsilon)$-holomorphic for $0 < \varepsilon < r$, where $r$ is as in \eqref{uniform-ellipticity}. See, e.g., \cite[Prop.\ 4.9]{adcock2022sparse}.

\subsection{The function considered in Fig.\ \ref{fig:fig1}}\label{ss:lognormal-pDE}

In Fig.\ \ref{fig:fig1} we consider an elliptic DE of the form \eqref{eq:fruitfly} with $k = 1$ and $\Omega = [0,1]$. Modifying an example from \cite[Eqn.\ (5.2)]{nobile2008sparse} (see also \cite{adcock2022efficient}), we consider a lognormal diffusion coefficient of the form
\bes{
\begin{split}
a(x,\bm{y}) 
& = \exp\left(1+ y_1 \left(\frac{\sqrt{\pi}\beta}{2}\right)^{\frac12} + \sum_{i=2}^d \; \zeta_i \; \vartheta_i(x) \; y_i\right) 
\\
\zeta_i 
& = (\sqrt{\pi} \beta)^{\frac12} \exp\left( \frac{-\left( \left\lfloor \frac{i}{2} \right\rfloor \pi \beta\right)^2}{8} \right), 
\quad
\vartheta_i(x) 
 = \begin{cases} 
 \sin\left( \left\lfloor\frac{i}{2}\right\rfloor \pi x/\beta_p \right) & \mbox{$i$ even}
 \\
 \cos\left(\left\lfloor\frac{i}{2}\right\rfloor \pi x/\beta_p \right) &  \mbox{$i$ odd} 
 \end{cases}.
 \\
 \beta_c &= 1/8 ,\ \beta_p = \max\{1, 2 \beta_c\},\ \beta = \beta_c/\beta_p.
 \end{split}
}
We consider the quantity of interest defined as
\bes{
f(\bm{y}) = u(0.5,\bm{y}),
}
where $u(x,\bm{y})$ is the corresponding parametric DE solution. To evaluate $f$, we use a piecewise linear finite element discretization with $1023$ degrees of freedom.

As a result of this discretization, the samples of $f$ involve a numerical error. This is equal to $f(\bm{y}_i) - f_h(\bm{y}_i)$, where $f_h$ is function computed via the finite element discretization. Since the focus of this paper is on sampling strategies, and not on numerical discretization, we ignore this error. In particular, we evaluate $f_h$ over the error grid, not $f$, and compute the error as
\bes{
E_h(f) = \frac{\sqrt{\frac{1}{K} \sum^{K}_{i=1} | f_h(\bm{z}_i) - \hat{f}(\bm{z}_i) |^2}}{\sqrt{\frac{1}{K} \sum^{K}_{i=1} | f_h(\bm{z}_i) |^2}}
}
instead of \eqref{f-error-fun}.

\section{Further experiments}\label{s:further-experiments}

We conclude this document with several additional experiments. {First, we consider several further examples from the test suite \cite{surjanovic2013virtual}. Second,} we consider what happens when the underlying measure $\varrho$ is changed. Third, we discuss how changing the oversampling amount in ALS affects its performance relative to CS. {And finally, we present additional experiments displaying the error in the $L^{\infty}$-norm.}

\subsection{{Further examples from the \textit{Virtual Library of Simulation Experiments}}}\label{ss:further-test-suite}

{
In addition to the \textit{borehole} function \eqref{borehole}, we now present several further examples from the suite of test functions \cite{surjanovic2013virtual}. Specifically, we consider functions from the \textit{Physical Models} section of the \textit{Emulation/Prediction Test Problems} test set. Each of these functions is a model for a certain physical process. For details and references, see \cite{surjanovic2013virtual}. }

{
Each function is defined has a certain input dimension $d_{*}$ and its parameters range in finite intervals that are typically not equal to $[-1,1]$. Let $\bm{\theta} = (\theta_i)^{d_{*}}_{i=1}$ be the parameters of any such function and $\bm{\theta}_{\min} = (\theta_{\min,i})^{d_{*}}_{i=1}$ and $\bm{\theta}_{\max} = (\theta_{\max,i})^{d_{*}}_{i=1}$ be the vector of parameter ranges, i.e.,
\bes{
\theta_{\min,i} \leq \theta_i \leq \theta_{\max,i},\quad i = 1,\ldots,d_{*}.
}
Write $g_{\mathsf{name}} : \prod^{d_{*}}_{i=1} [ \theta_{\min,i} , \theta_{\max,i} ] \rightarrow \bbR, \bm{\theta} \mapsto g(\bm{\theta})$ for the given function. Then we define the corresponding function $F : [-1,1]^{d_*} \rightarrow \bbR$ via an affine map
\bes{
F_{\mathsf{name}}(\bm{y}) = g_{\mathsf{name}} \left(\bar{\bm{\theta}} + \bm{\delta} \odot \bm{y} \right ),\quad \bm{y} \in [-1,1]^{d_{*}},
}
where $\bar{\bm{\theta}} = \frac12 ( \bm{\theta}_{\min} + \bm{\theta}_{\max} )$, $\bm{\delta} = \frac12 (\bm{\theta}_{\max} - \bm{\theta}_{\min})$ and $\odot$ denotes the Hadamard (or componentwise) product. Finally, for each $1 \leq d \leq d_{*}$, we define the function $f : \cU = : [-1,1]^d \rightarrow \bbR$ as
\bes{
f_{\mathsf{name}}(\bm{y}) = F_{\mathsf{name}}(y_1,\ldots,y_d,1,\ldots,1),\quad \bm{y} = (y_i)^{d}_{i=1} \in \cU
}
by setting the final $d_{*} - d$ parameters equal to their right endpoint value. In the Fig.\ \ref{fig:fig5} and the other figures in this section, we show numerical results for various values of $d$, including $d = d_{*}$.
}

\begin{figure}[t!]
\begin{center}
\begin{small}
 \begin{tabular}{@{\hspace{0pt}}c@{\hspace{\errplotsp}}c@{\hspace{\errplotsp}}c@{\hspace{0pt}}}
\includegraphics[width = \errplotimg]{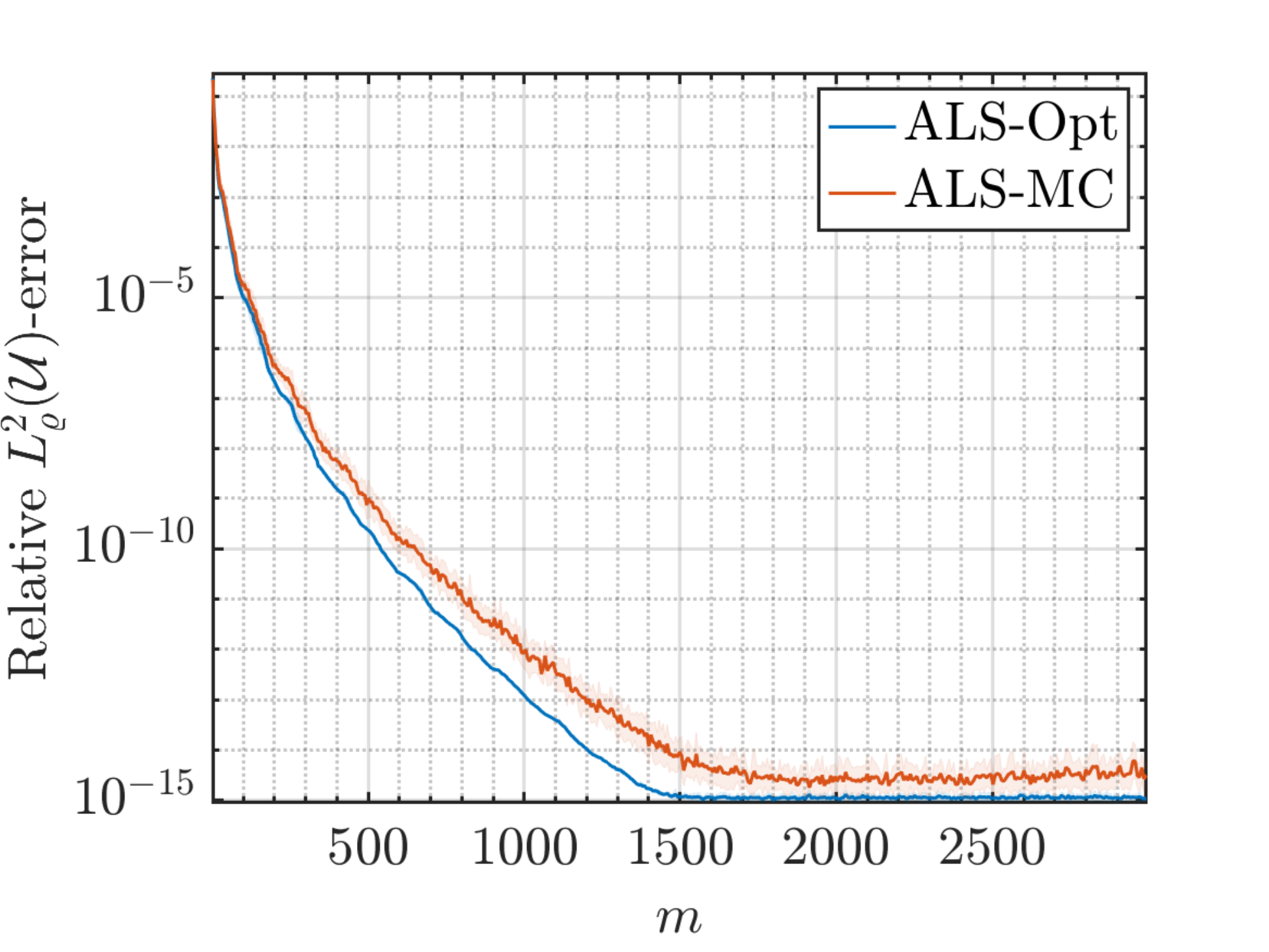}&
\includegraphics[width = \errplotimg]{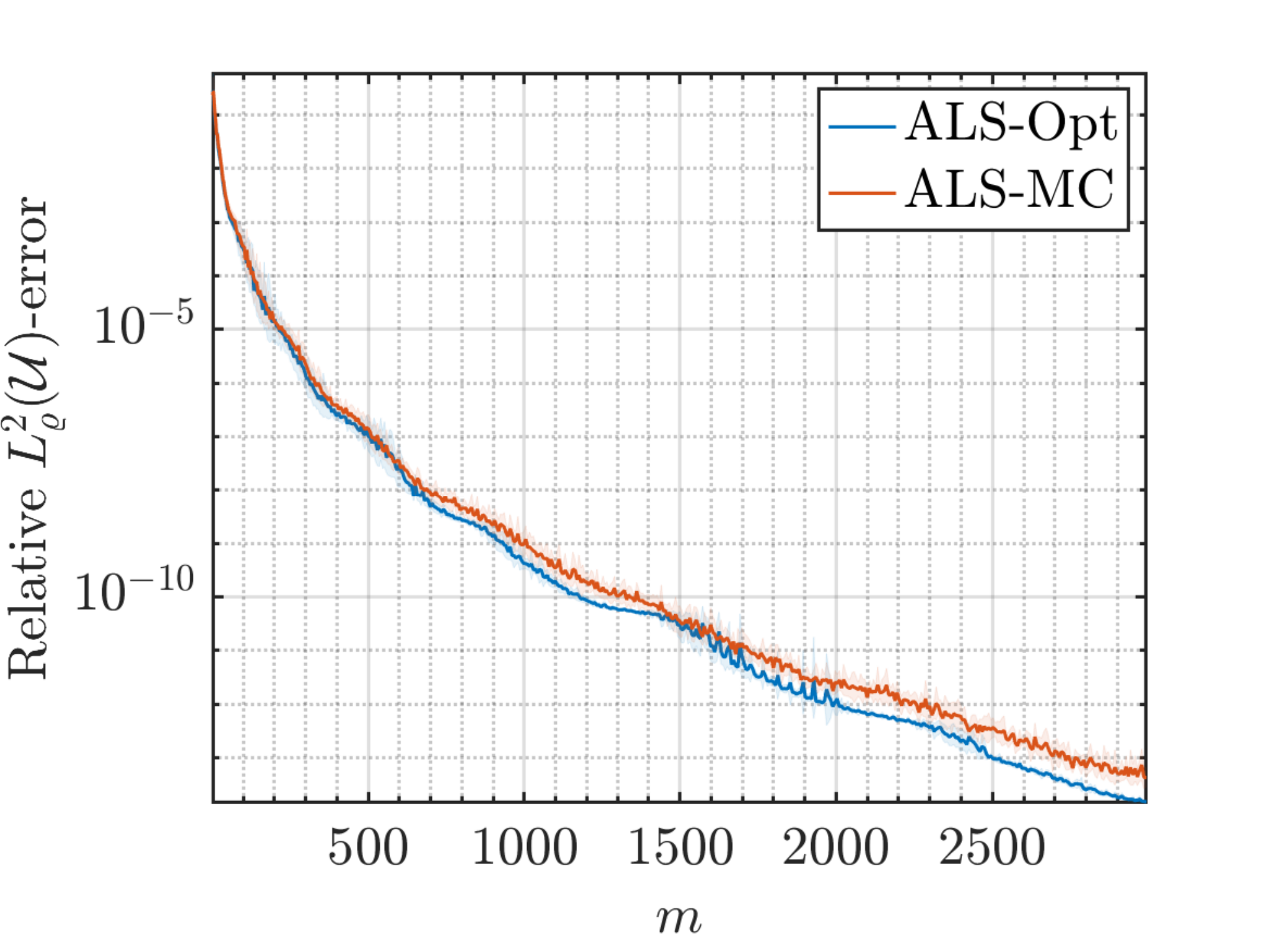}&
\includegraphics[width = \errplotimg]{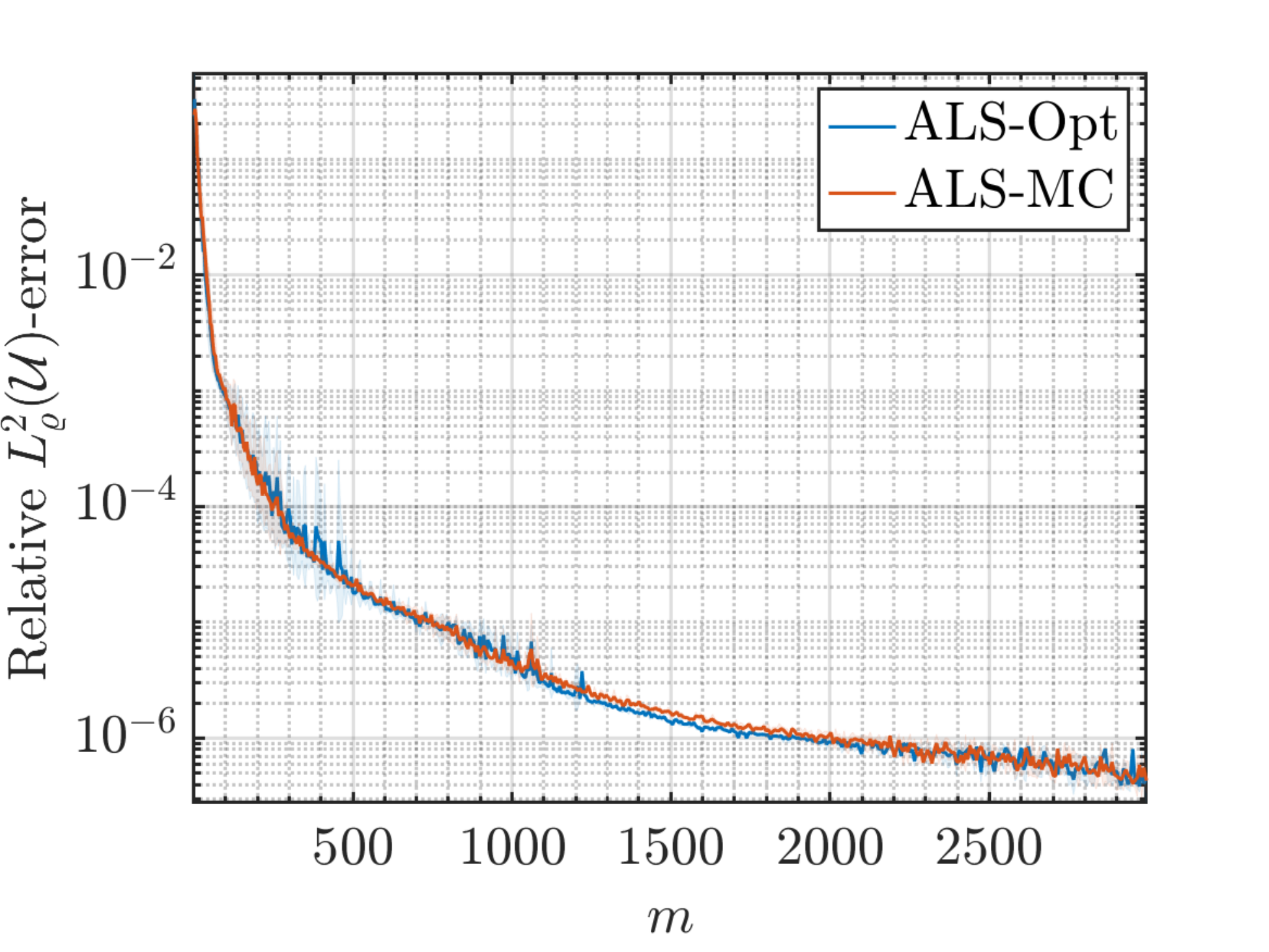}
\\[\errplotgraphsp]
\includegraphics[width = \errplotimg]{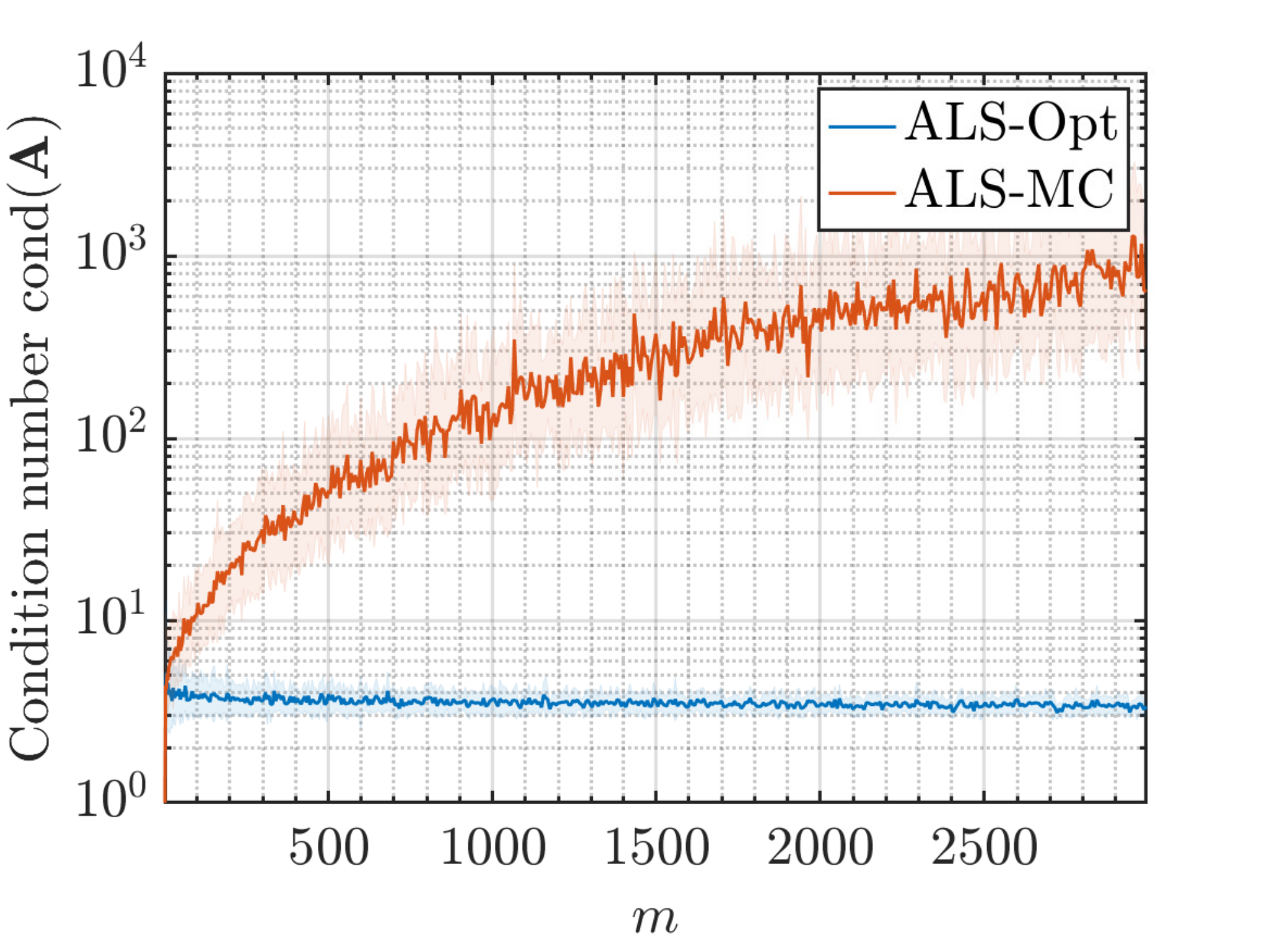}&
\includegraphics[width = \errplotimg]{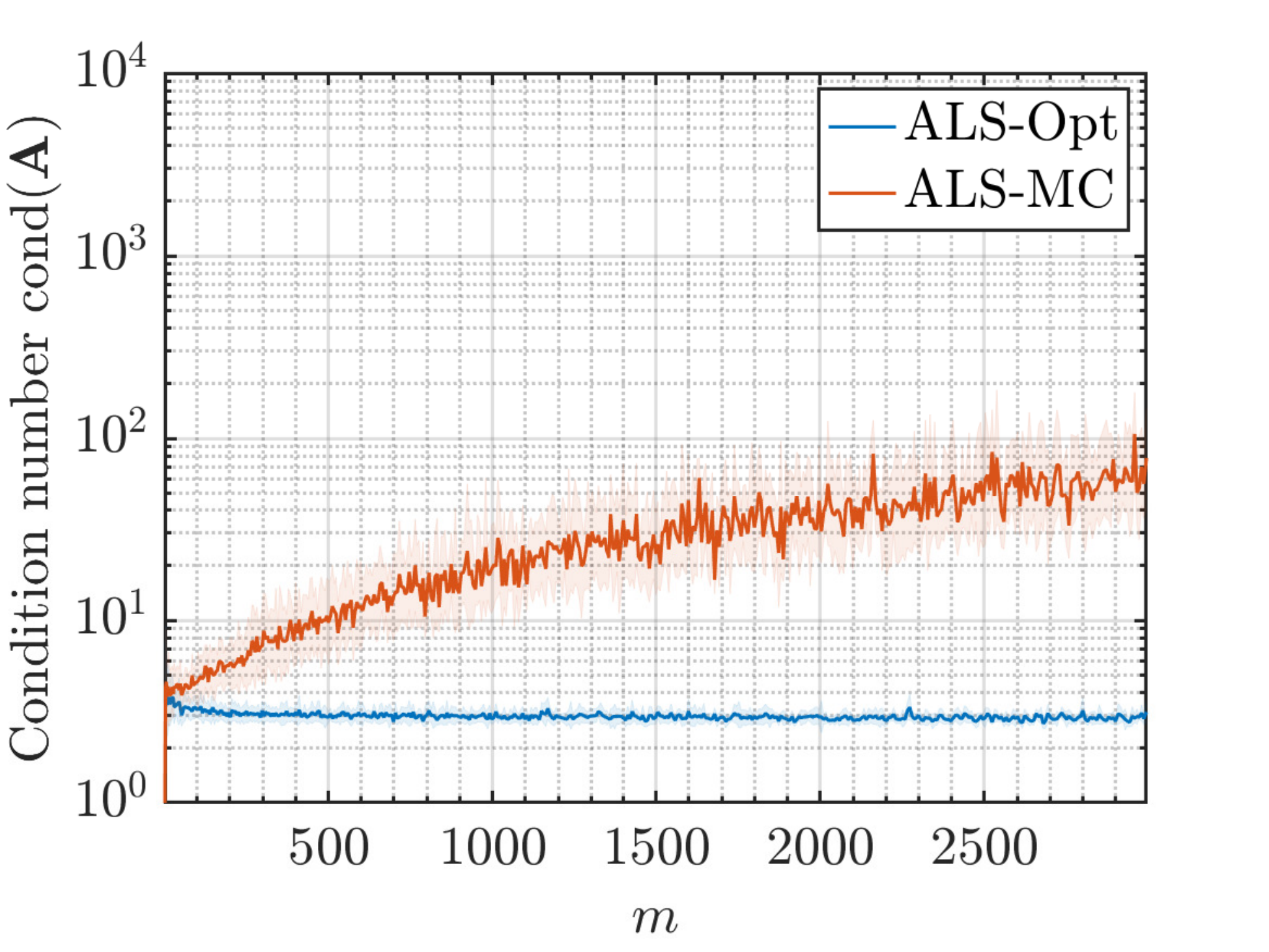}&
\includegraphics[width = \errplotimg]{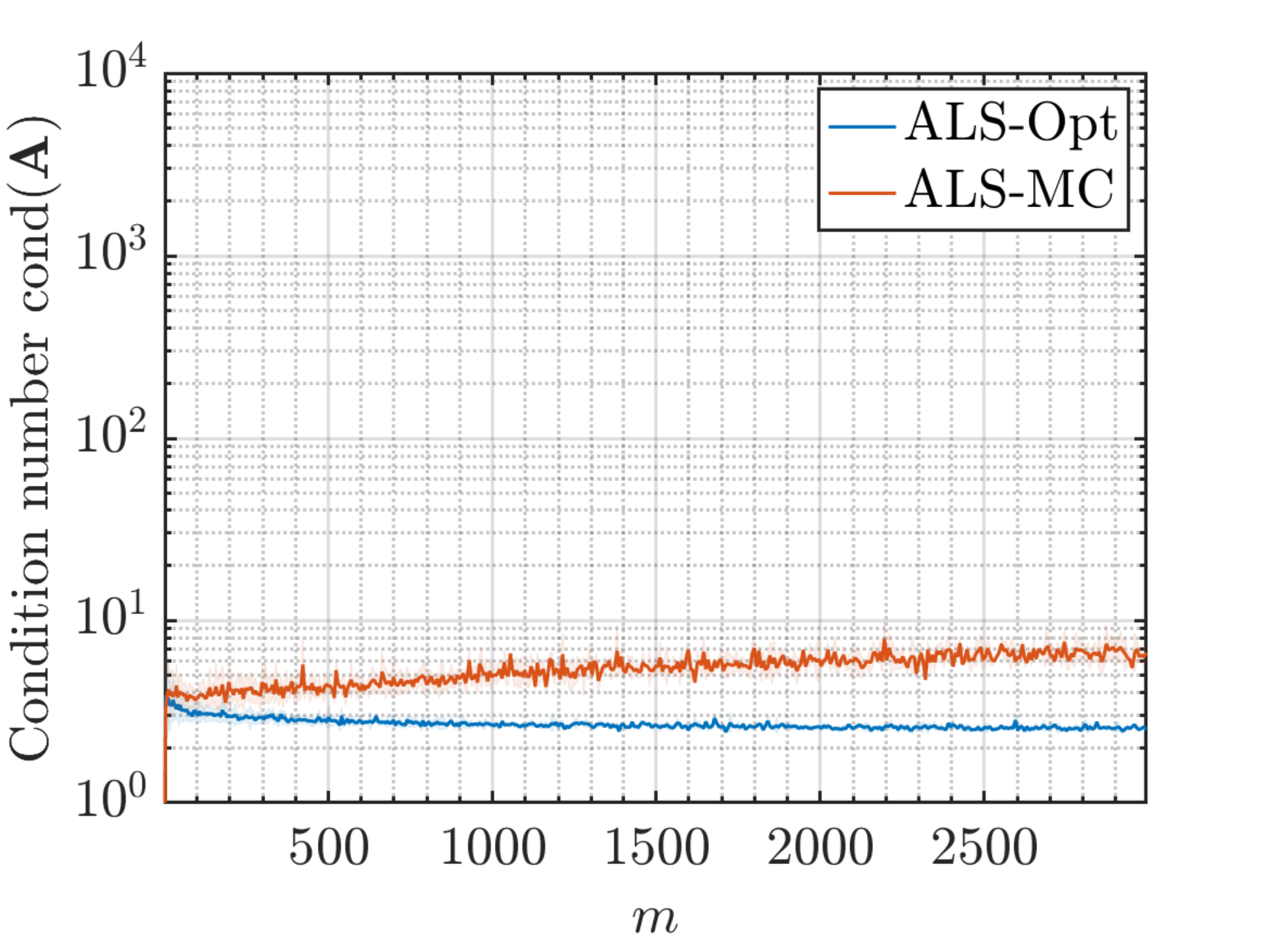}
\\[\errplottextsp]
$d = 2$ & $d = 4$ & $d = 6$
\end{tabular}
\end{small}
\end{center}
\caption{{The same as Fig.\ \ref{fig:fig2}, except for the function $f = f_{\mathsf{circuit}}$. Here $d_{*} = 6$.}}
\label{fig:fig-otlcircuit}
\end{figure}

\begin{figure}[t!]
\begin{center}
\begin{small}
 \begin{tabular}{@{\hspace{0pt}}c@{\hspace{\errplotsp}}c@{\hspace{\errplotsp}}c@{\hspace{0pt}}}
\includegraphics[width = \errplotimg]{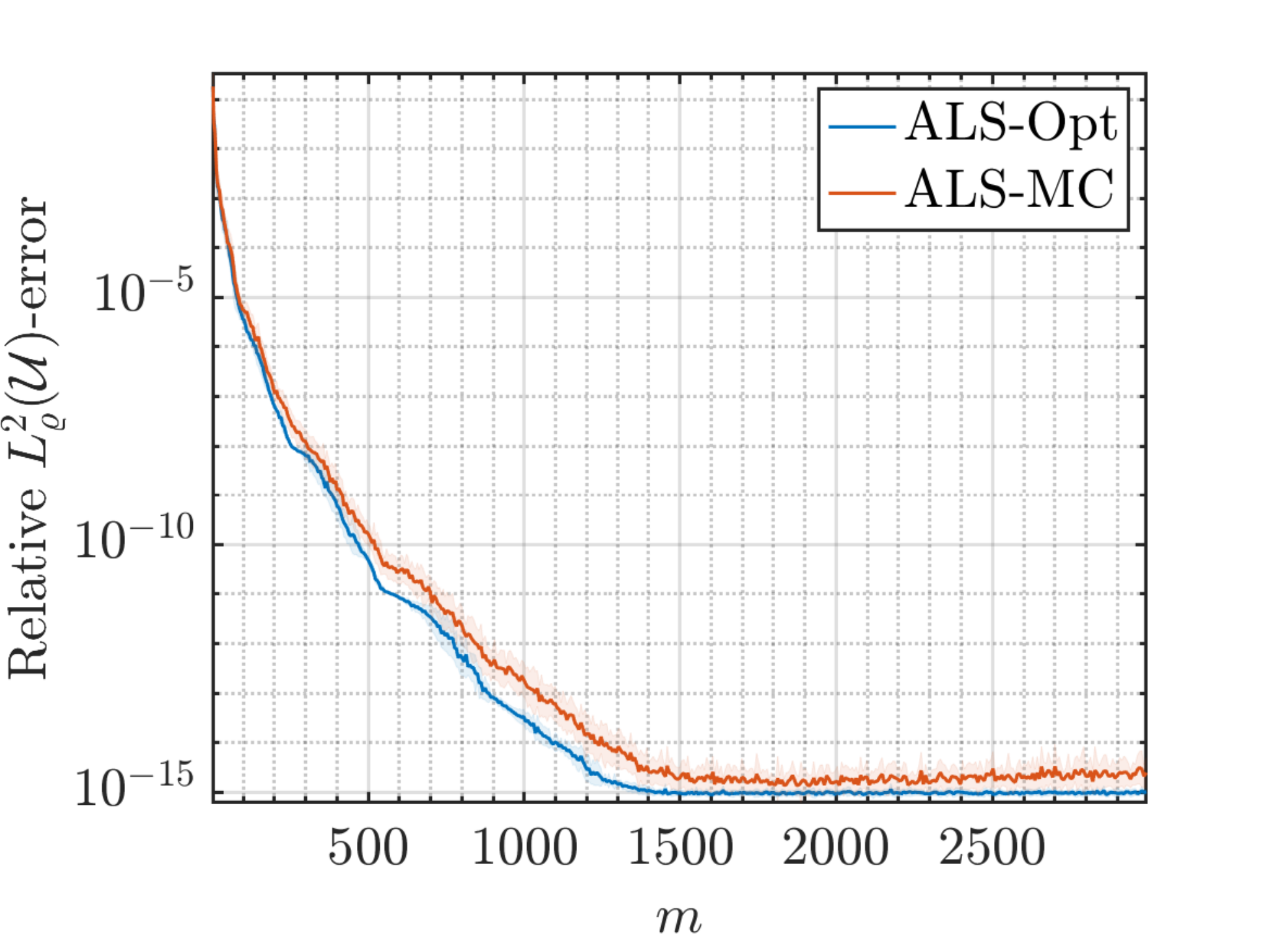}&
\includegraphics[width = \errplotimg]{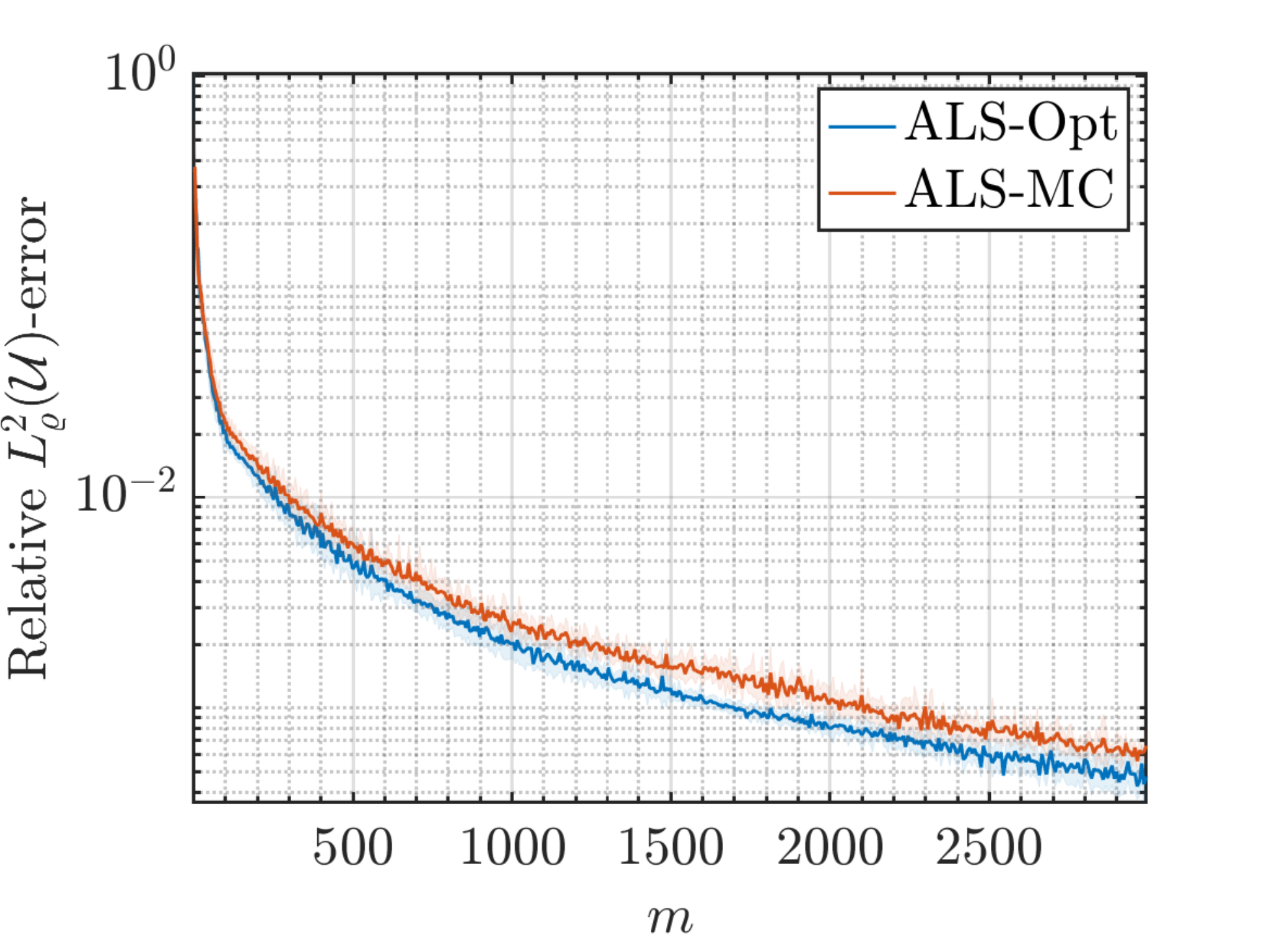}&
\includegraphics[width = \errplotimg]{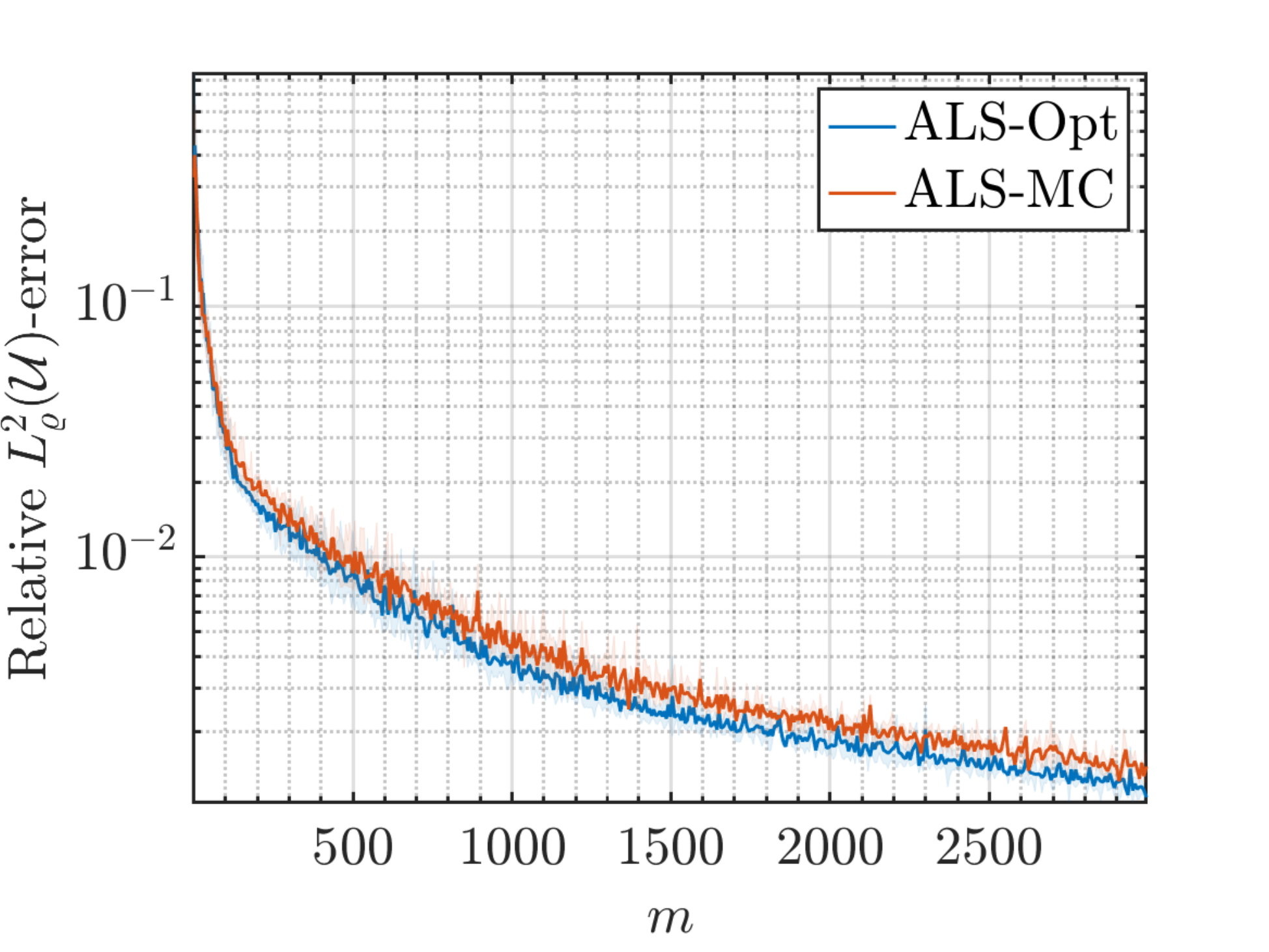}
\\[\errplotgraphsp]
\includegraphics[width = \errplotimg]{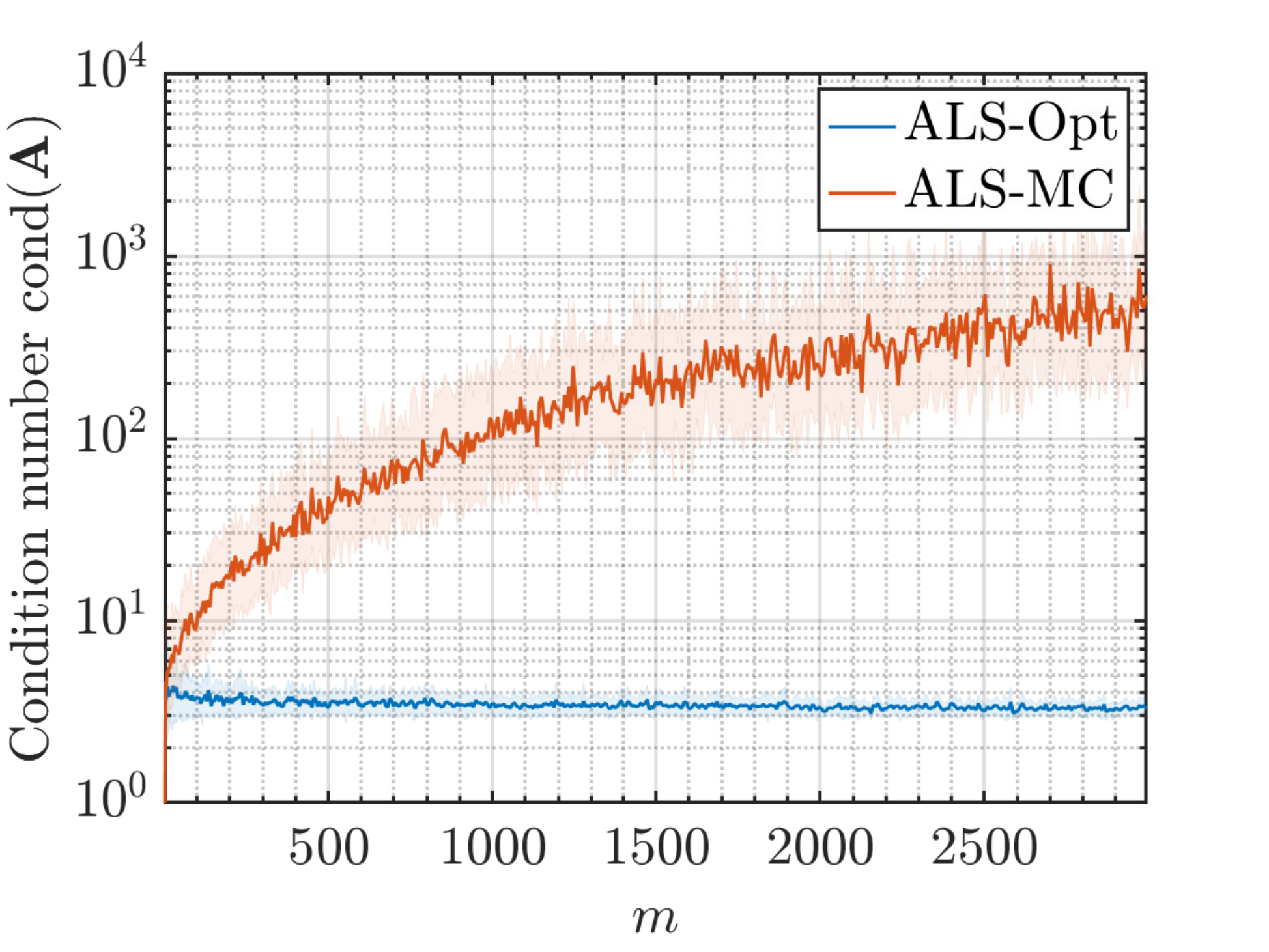}&
\includegraphics[width = \errplotimg]{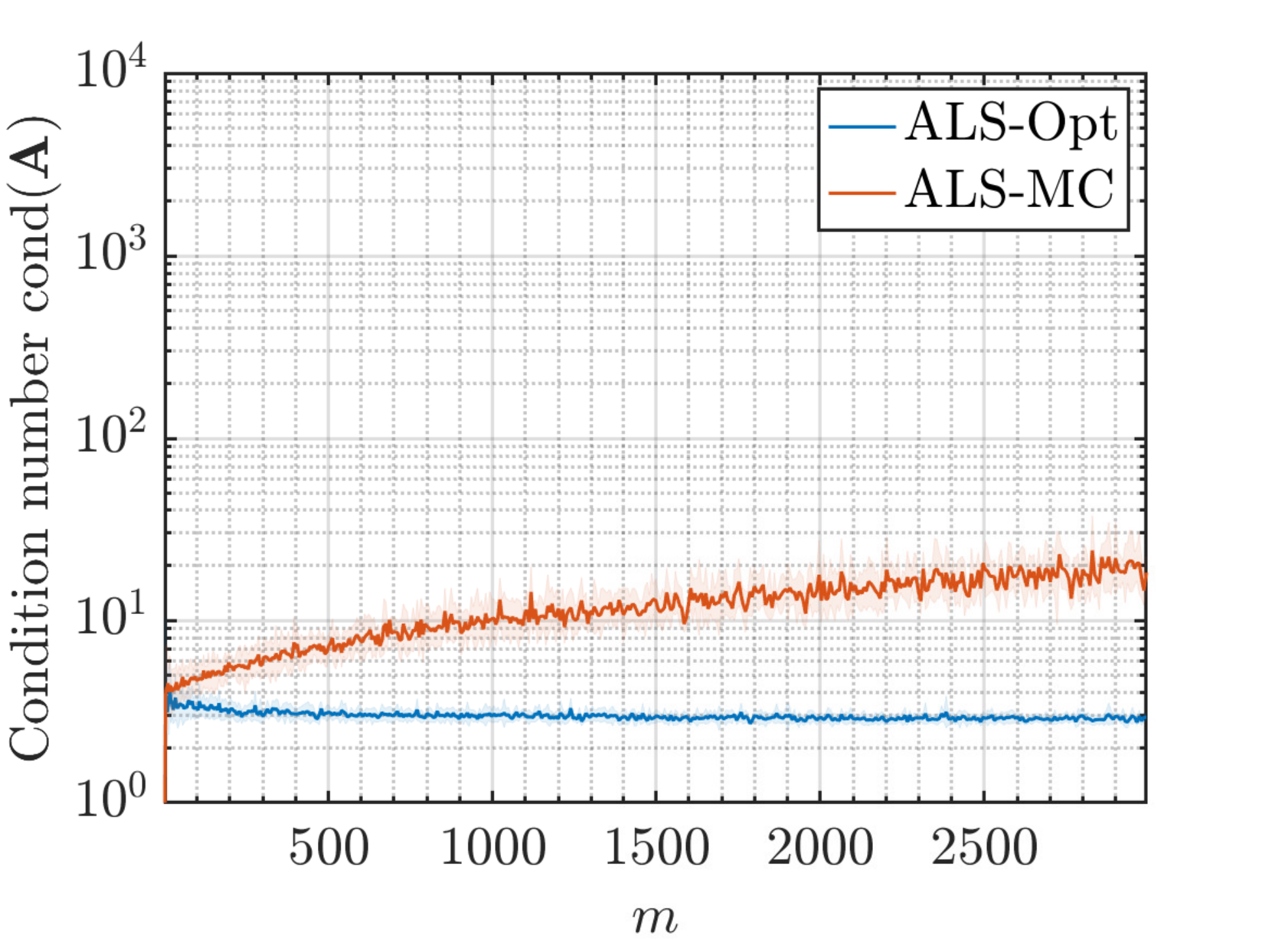}&
\includegraphics[width = \errplotimg]{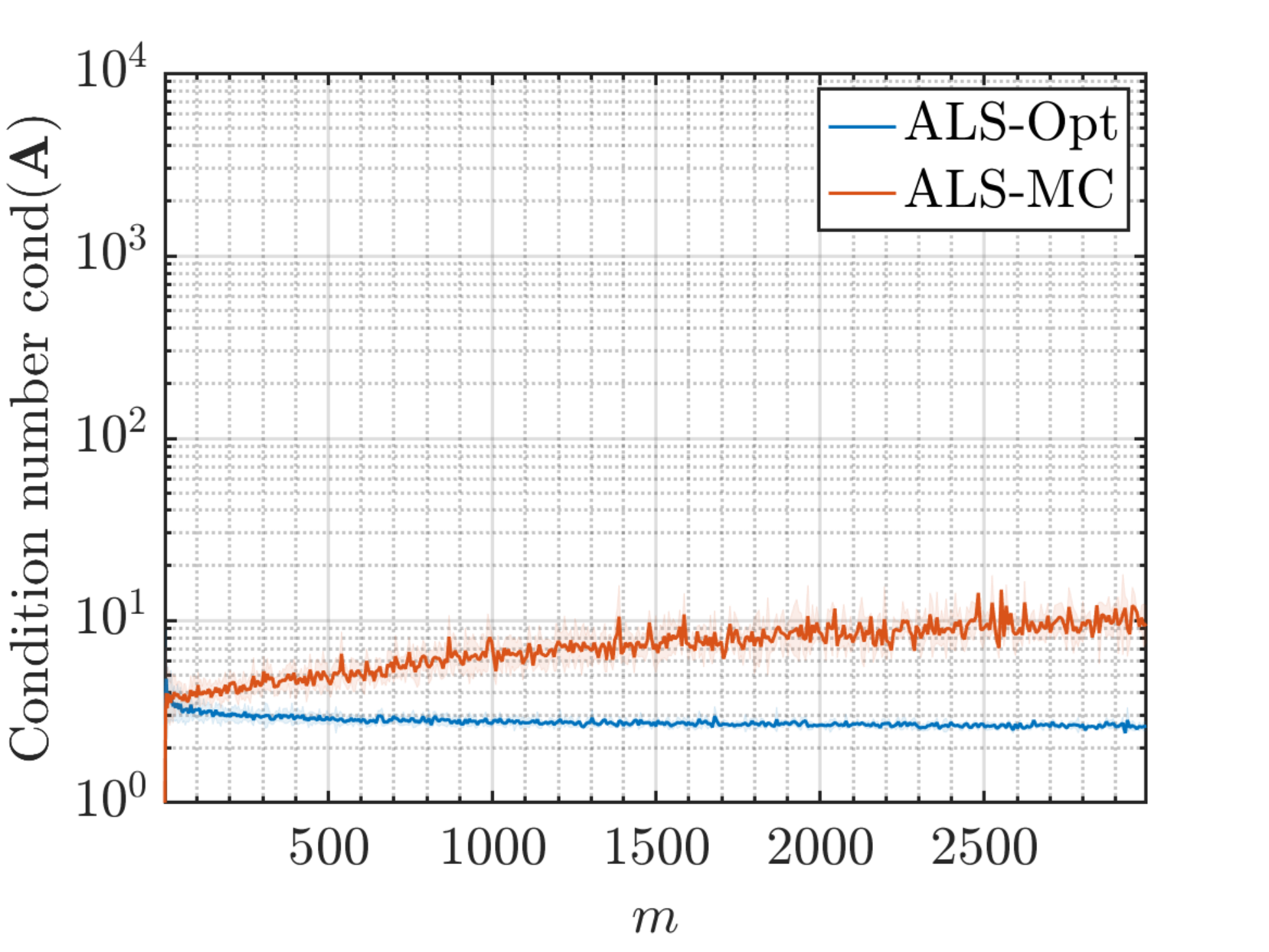}
\\[\errplottextsp]
$d = 2$ & $d = 4$ & $d = 7$
\end{tabular}
\end{small}
\end{center}
\caption{{The same as Fig.\ \ref{fig:fig2}, except for the function $f = f_{\mathsf{piston}}$. Here $d_{*} = 7$.}}
\label{fig:fig-piston}
\end{figure}

\begin{figure}[t!]
\begin{center}
\begin{small}
 \begin{tabular}{@{\hspace{0pt}}c@{\hspace{\errplotsp}}c@{\hspace{\errplotsp}}c@{\hspace{0pt}}}
\includegraphics[width = \errplotimg]{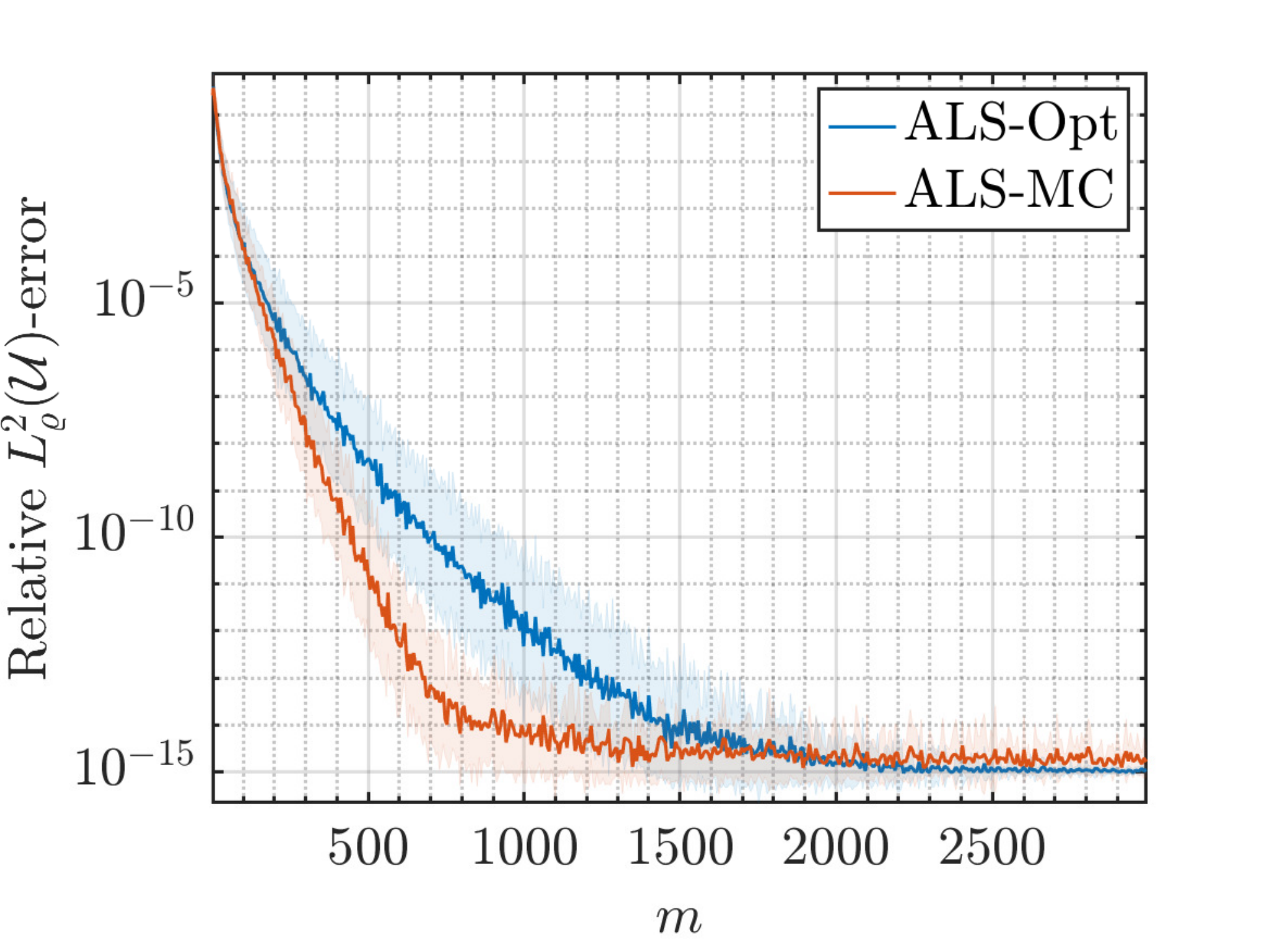}&
\includegraphics[width = \errplotimg]{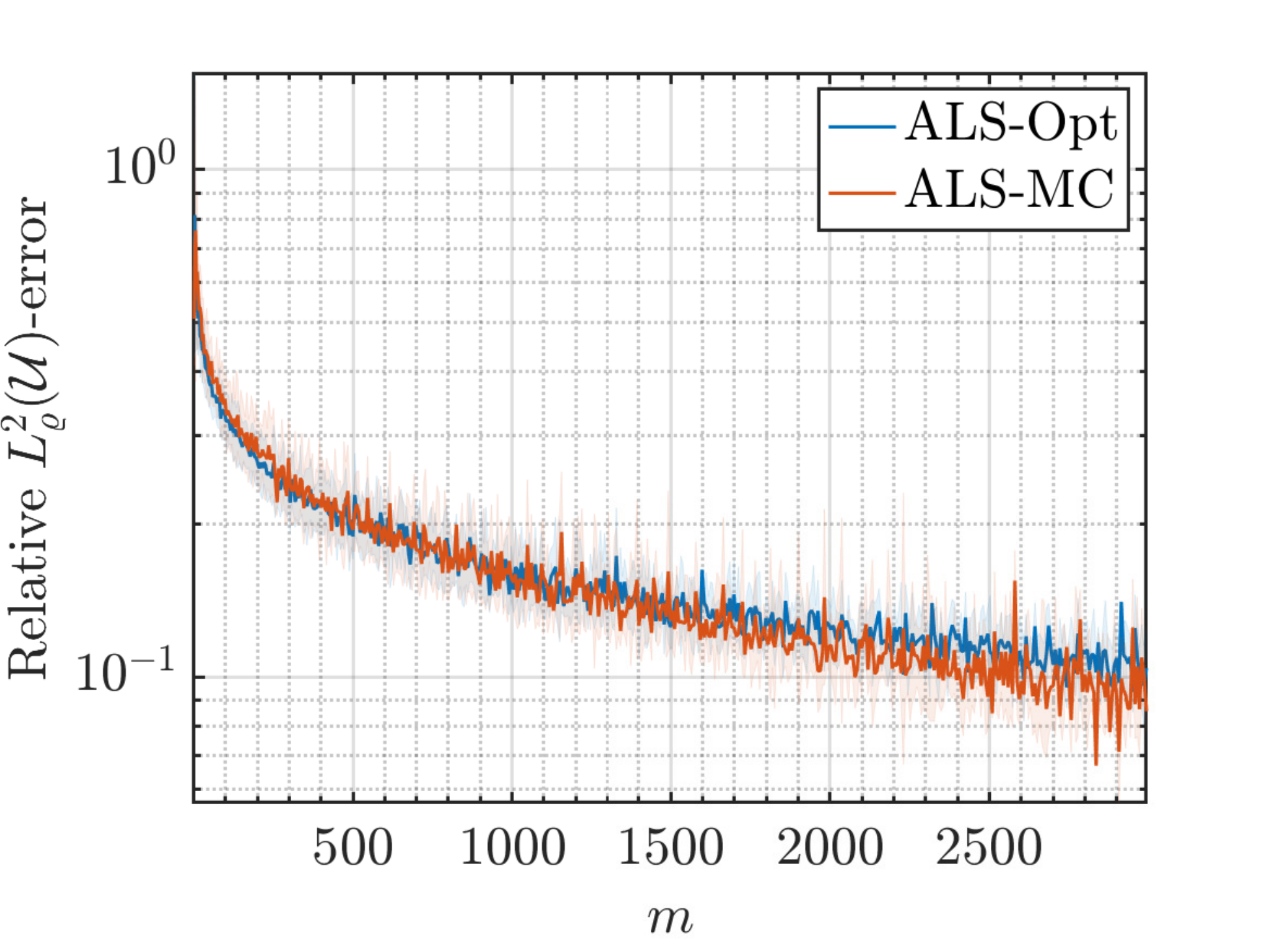}&
\includegraphics[width = \errplotimg]{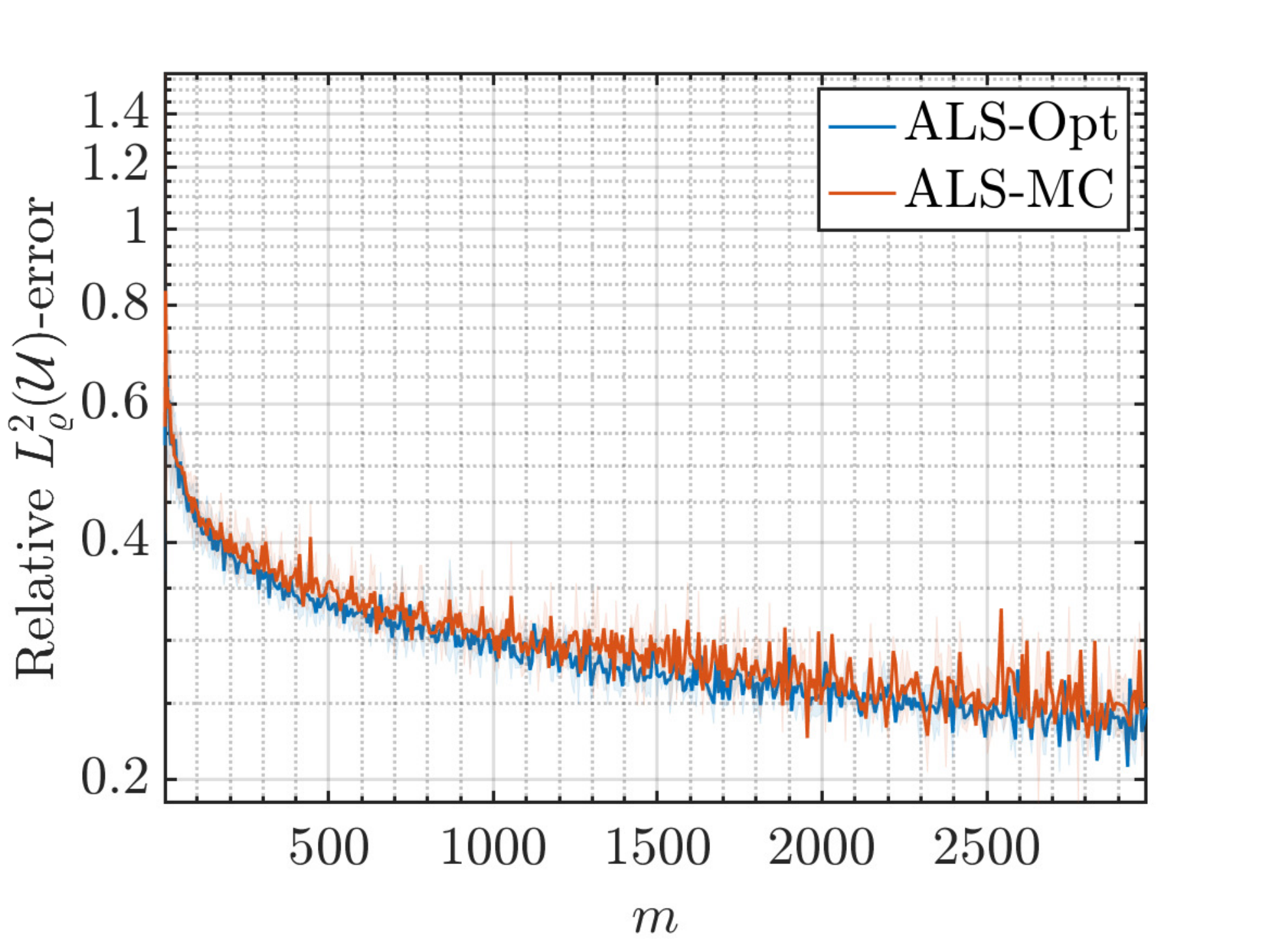}
\\[\errplotgraphsp]
\includegraphics[width = \errplotimg]{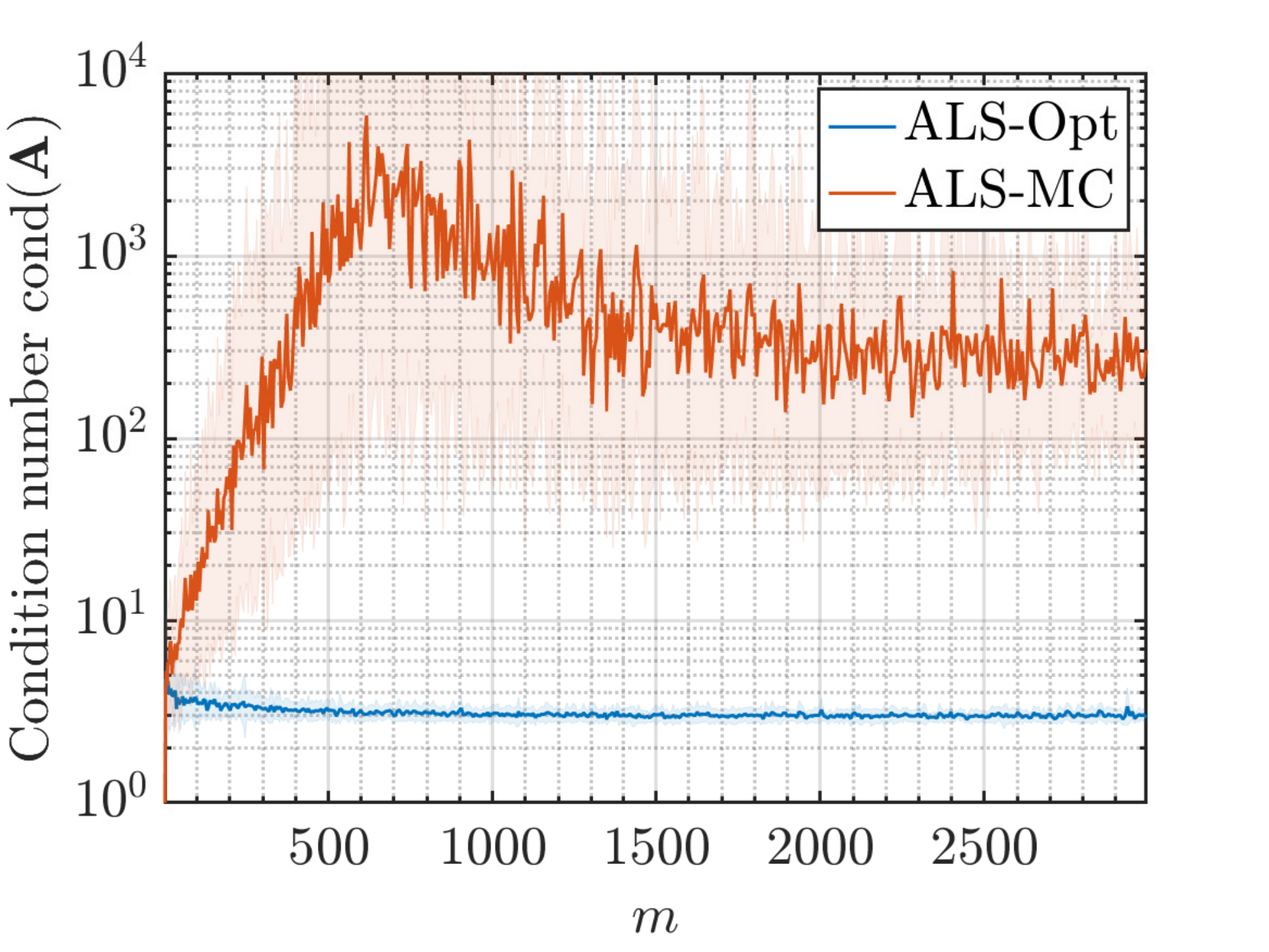}&
\includegraphics[width = \errplotimg]{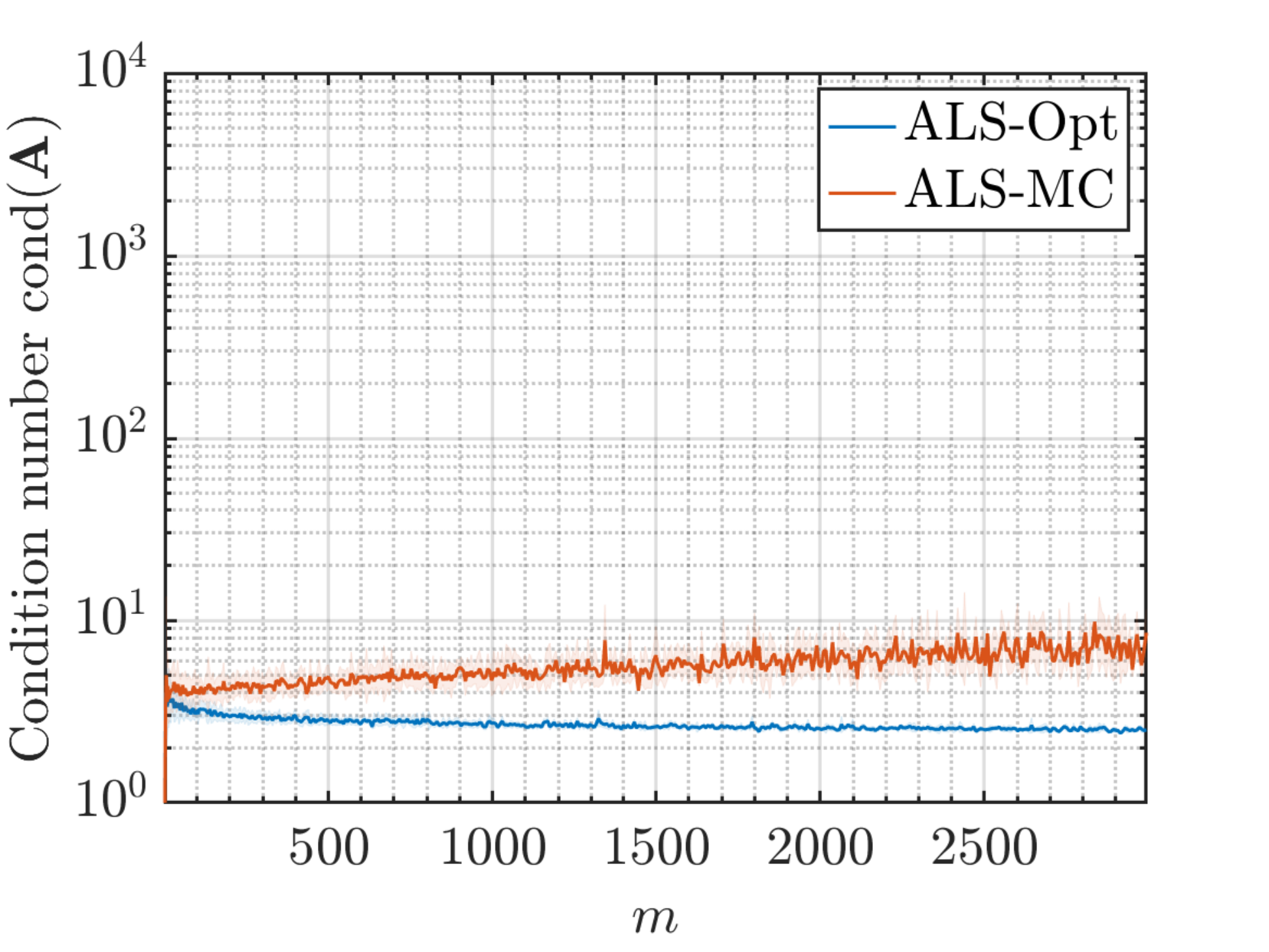}&
\includegraphics[width = \errplotimg]{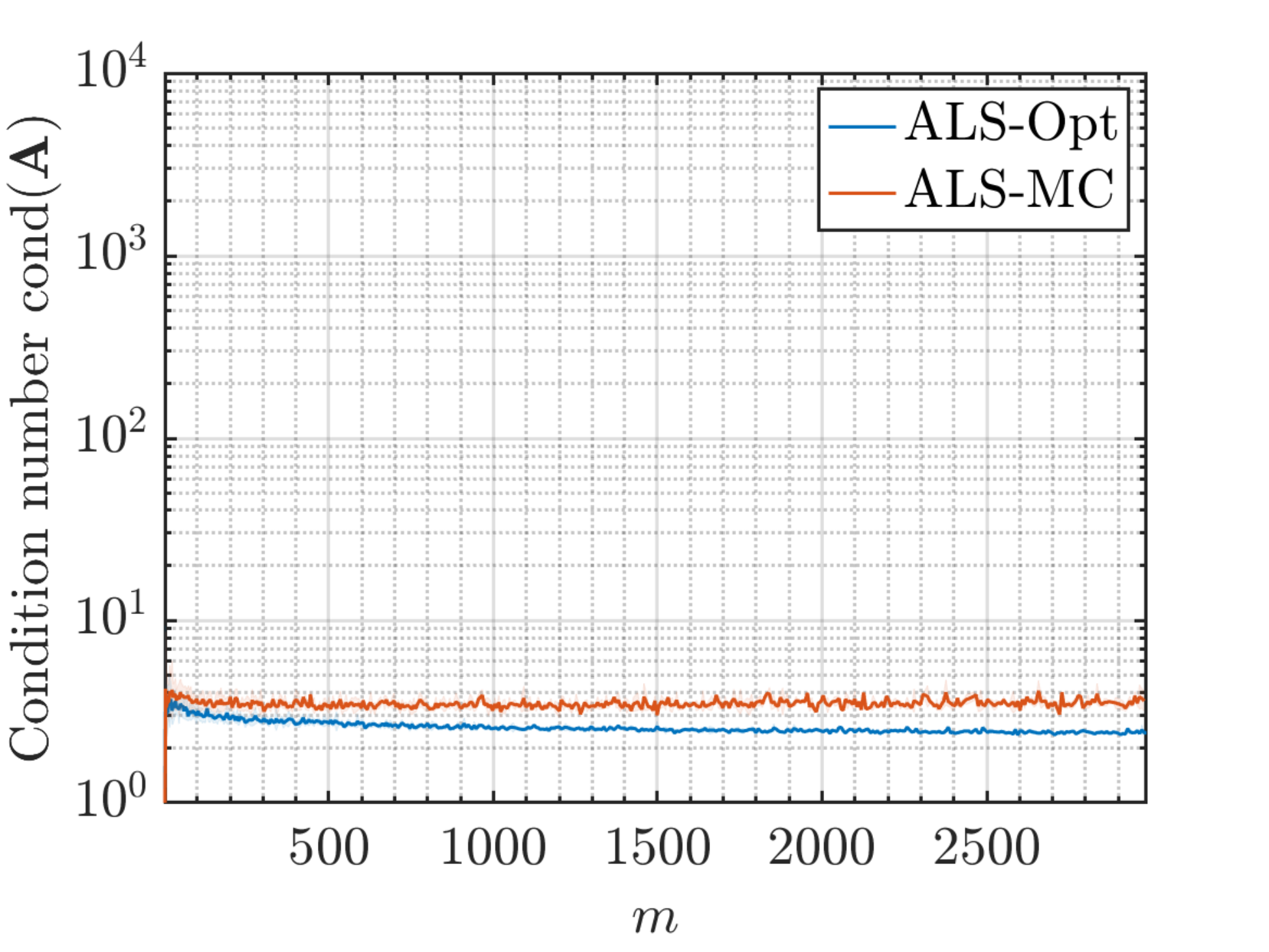}
\\[\errplottextsp]
$d = 2$ & $d = 5$ & $d = 8$
\end{tabular}
\end{small}
\end{center}
\caption{{The same as Fig.\ \ref{fig:fig2}, except for the function $f = f_{\mathsf{robot}}$. Here $d_{*} = 8$.}}
\label{fig:fig-robot}
\end{figure}

\begin{figure}[t!]
\begin{center}
\begin{small}
 \begin{tabular}{@{\hspace{0pt}}c@{\hspace{\errplotsp}}c@{\hspace{\errplotsp}}c@{\hspace{0pt}}}
\includegraphics[width = \errplotimg]{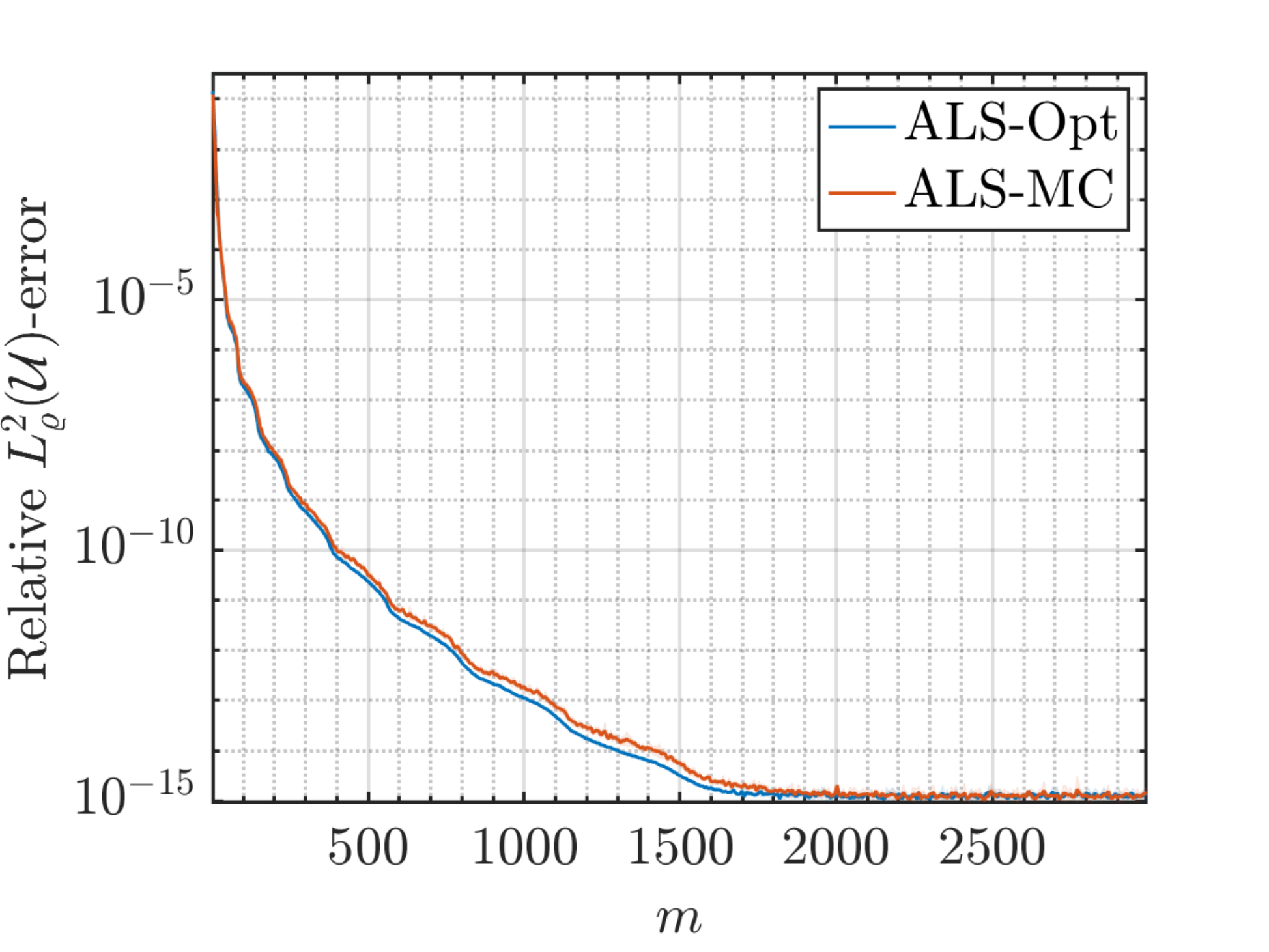}&
\includegraphics[width = \errplotimg]{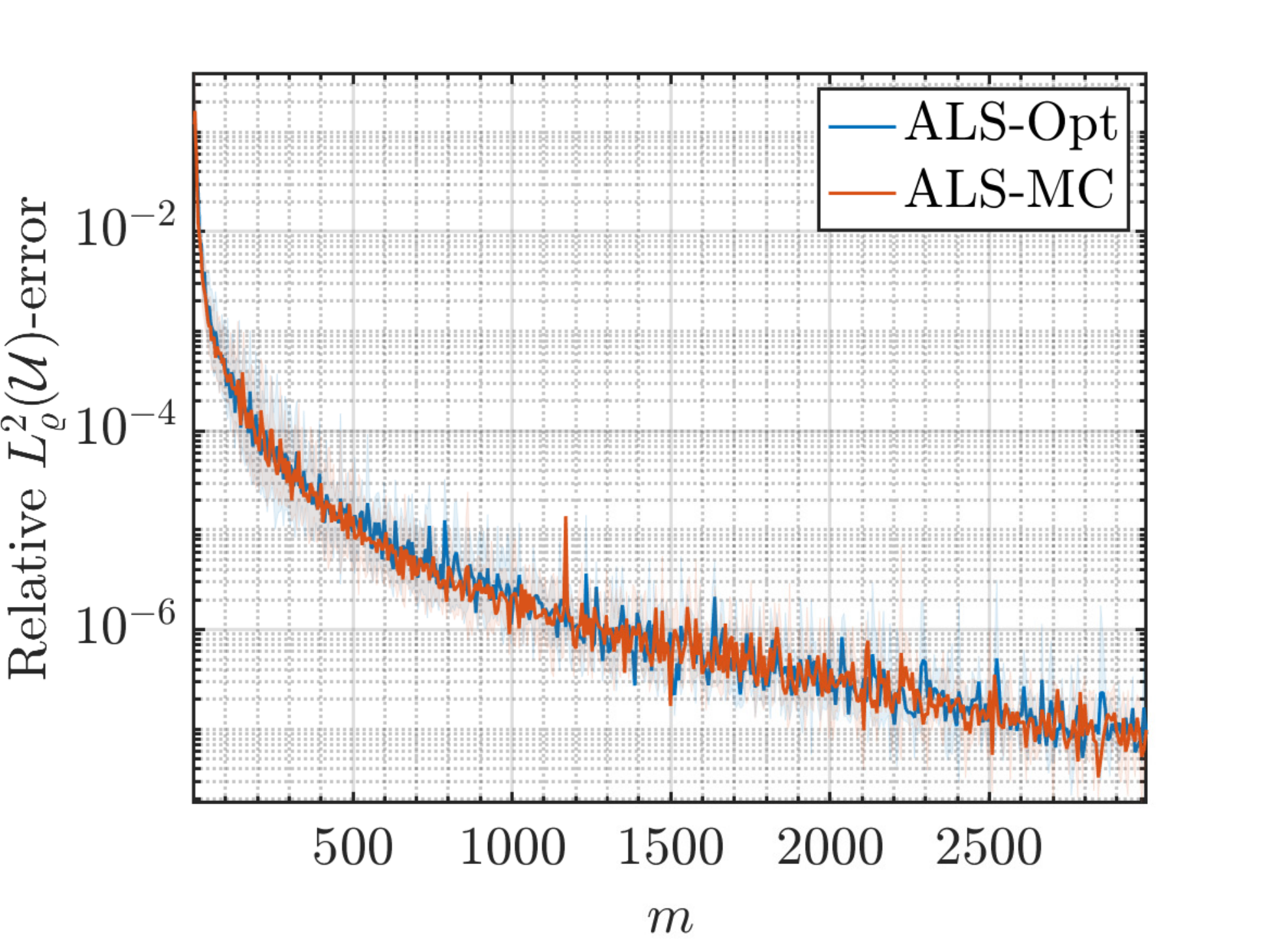}&
\includegraphics[width = \errplotimg]{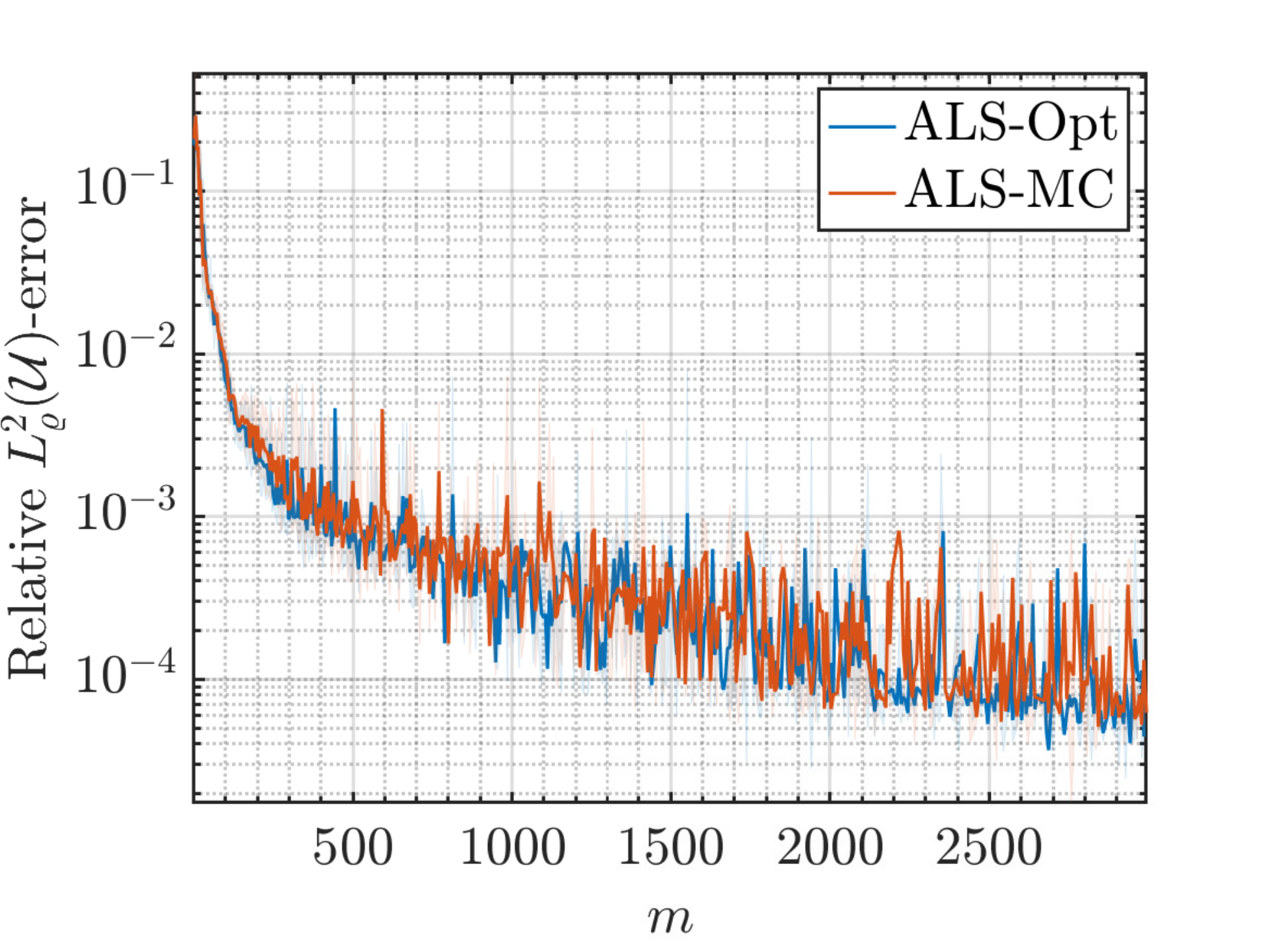}
\\[\errplotgraphsp]
\includegraphics[width = \errplotimg]{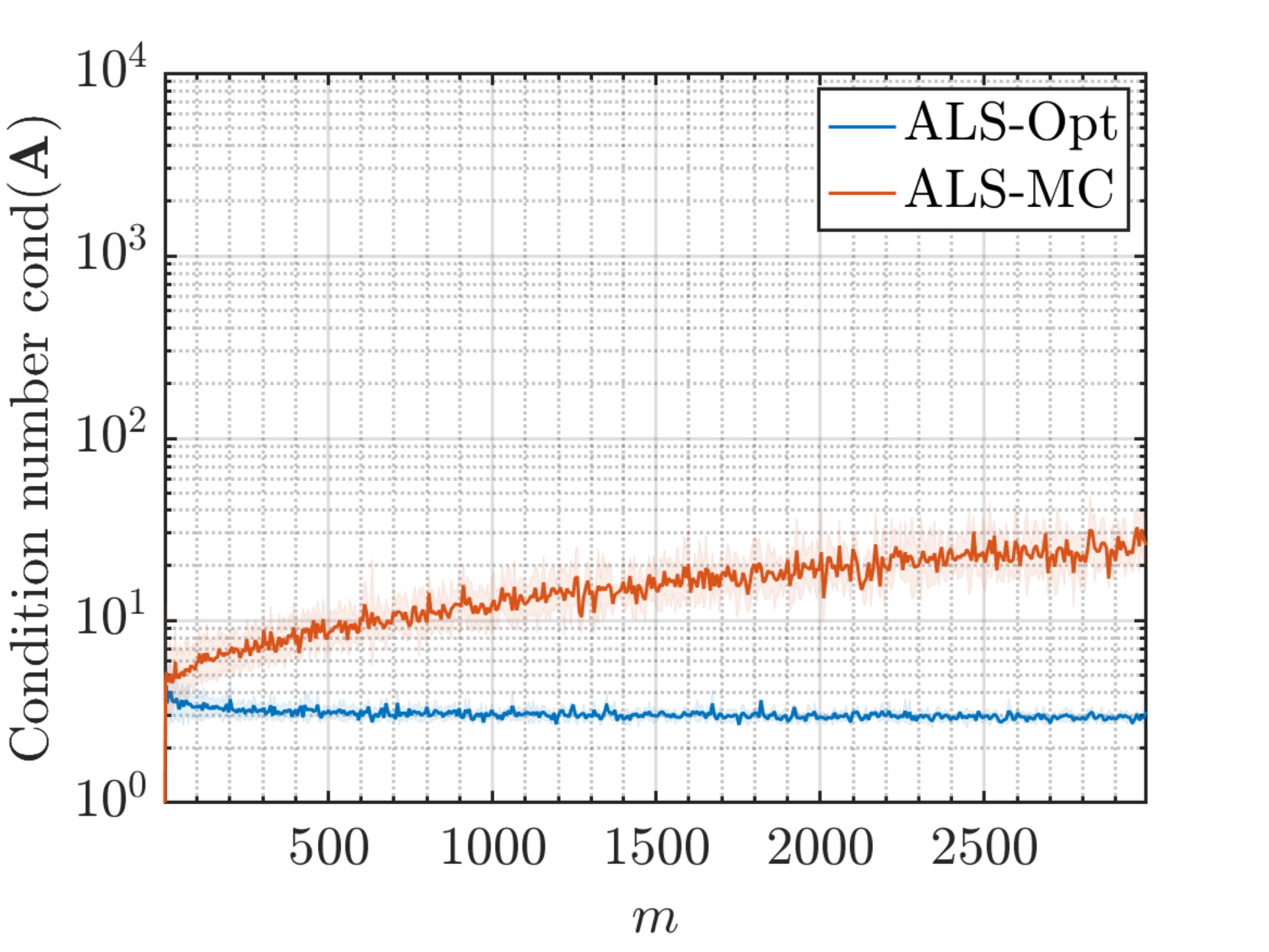}&
\includegraphics[width = \errplotimg]{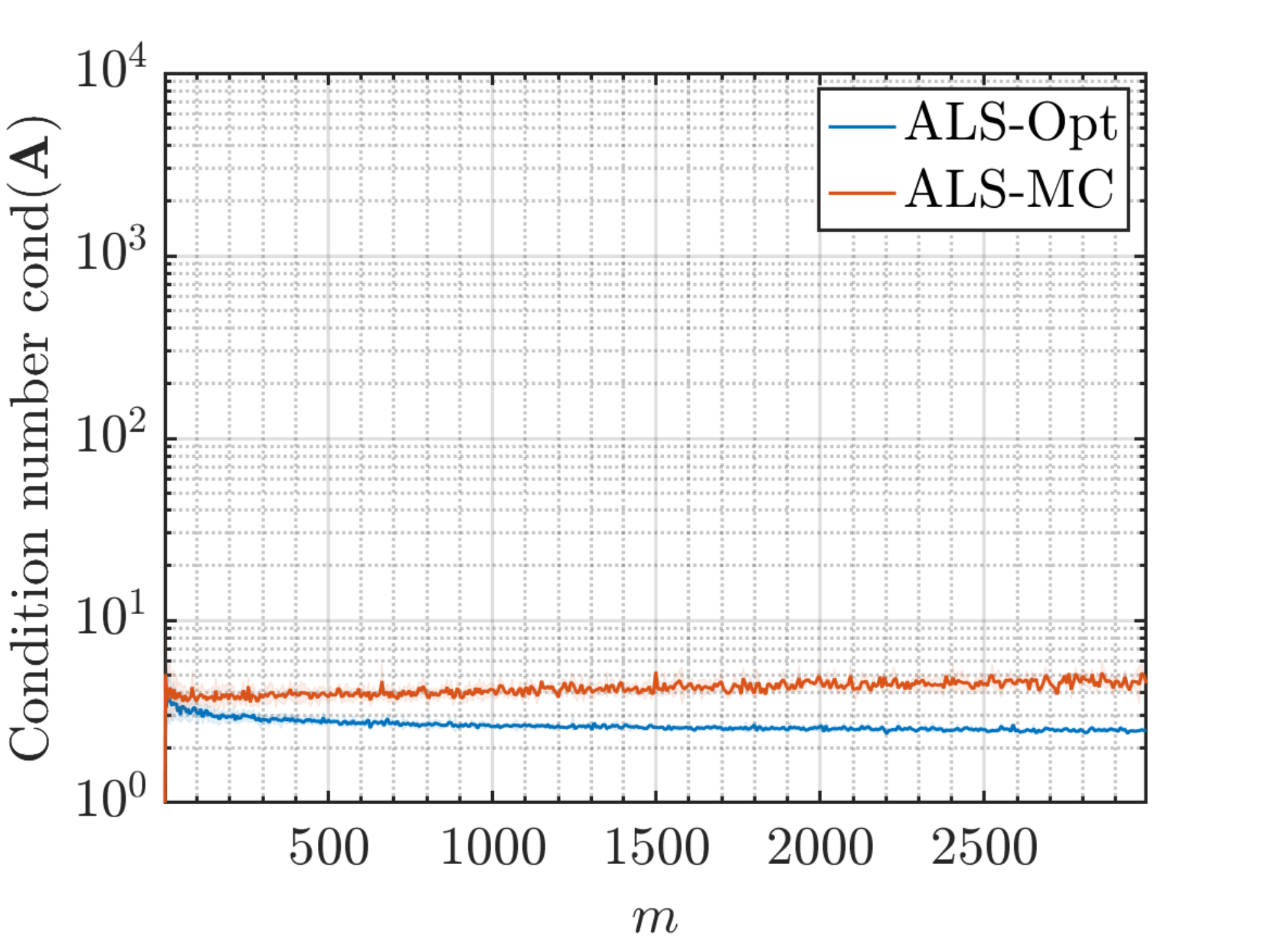}&
\includegraphics[width = \errplotimg]{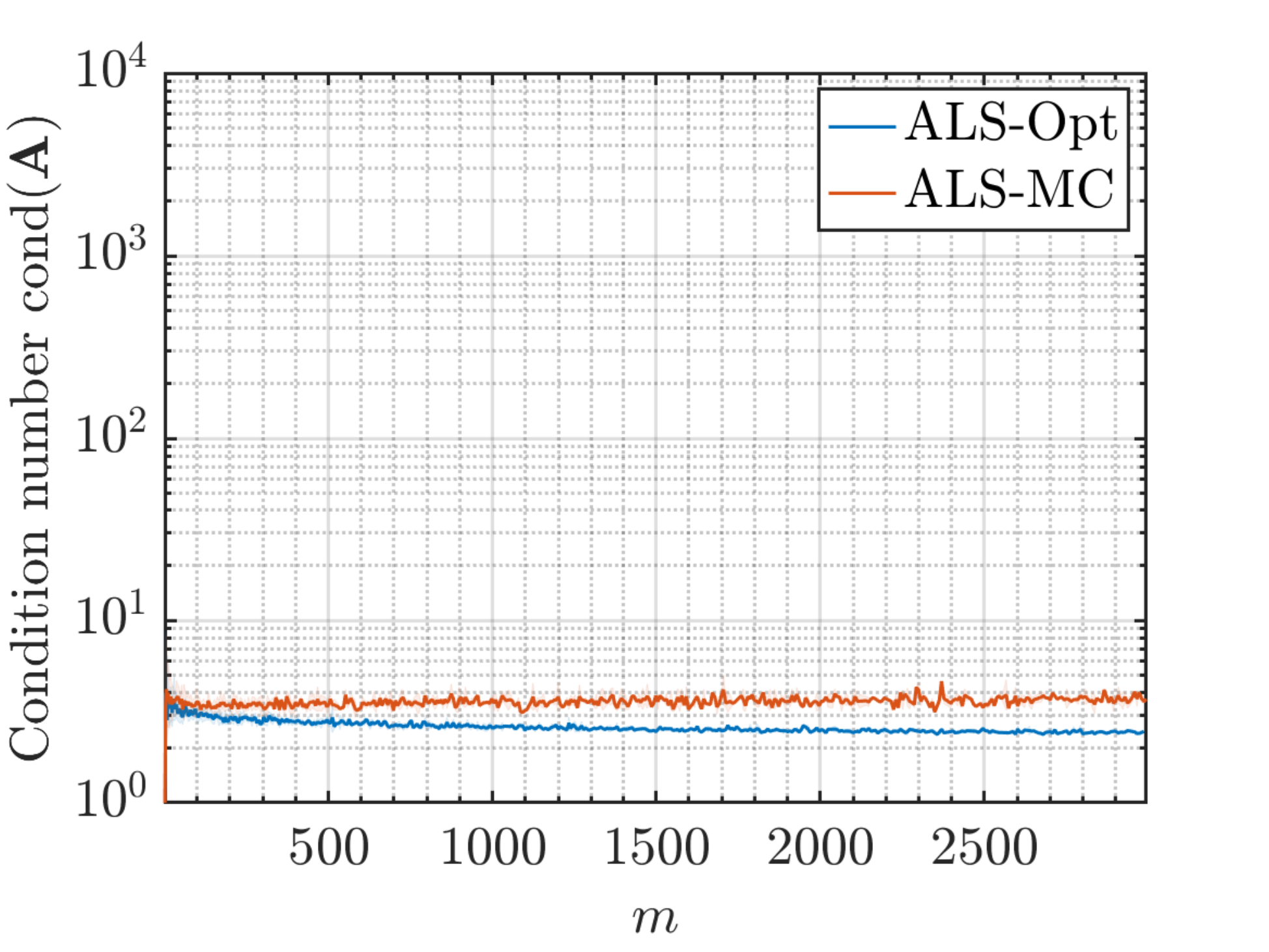}
\\[\errplottextsp]
$d = 3$ & $d = 6$ & $d = 10$
\end{tabular}
\end{small}
\end{center}
\caption{{The same as Fig.\ \ref{fig:fig2}, except for the function $f = f_{\mathsf{wing}}$. Here $d_{*} = 10$.}}
\label{fig:fig-wingweight}
\end{figure}

{In Figs.\ \ref{fig:fig-otlcircuit}--\ref{fig:fig-wingweight} we show numerical experiments for the \textit{OTL circuit} function $f = f_{\mathsf{circuit}}$, the \textit{piston simulation} function $f = f_{\mathsf{piston}}$, the \textit{robot arm} function $f = f_{\mathsf{robot}}$ and the \textit{wing weight} function $f = f_{\mathsf{wing}}$. In all cases, we witness a similar effect as in the results shown in the main paper. Namely, MC sampling is generally worse in low dimensions than the near-optimal scheme, yet in higher dimensions the two schemes perform very similarly. It is interesting that MC sampling yields a smaller error for $f = f_{\mathsf{robot}}$ in $d = 2$ dimensions than the near-optimal scheme. Notice that $f_{\mathsf{robot}}(y_1,y_2) = \sqrt{1-6 \cos( \pi y_2)}$ is constant in $y_1$ and sinusoidal in $y_2$. We believe this improved performance stems from this particular structure of the function. Note that the $f_{\mathsf{robot}}$ is nonsmooth in $d \geq 3$ dimensions, which  likely explains the significantly larger errors seen in higher dimensions. }

\subsection{Examples with different measures}\label{ss:different-measures}

\begin{figure}[t!]
\begin{center}
\begin{small}
 \begin{tabular}{@{\hspace{0pt}}c@{\hspace{\errplotsp}}c@{\hspace{\errplotsp}}c@{\hspace{0pt}}}
\includegraphics[width = \errplotimg]{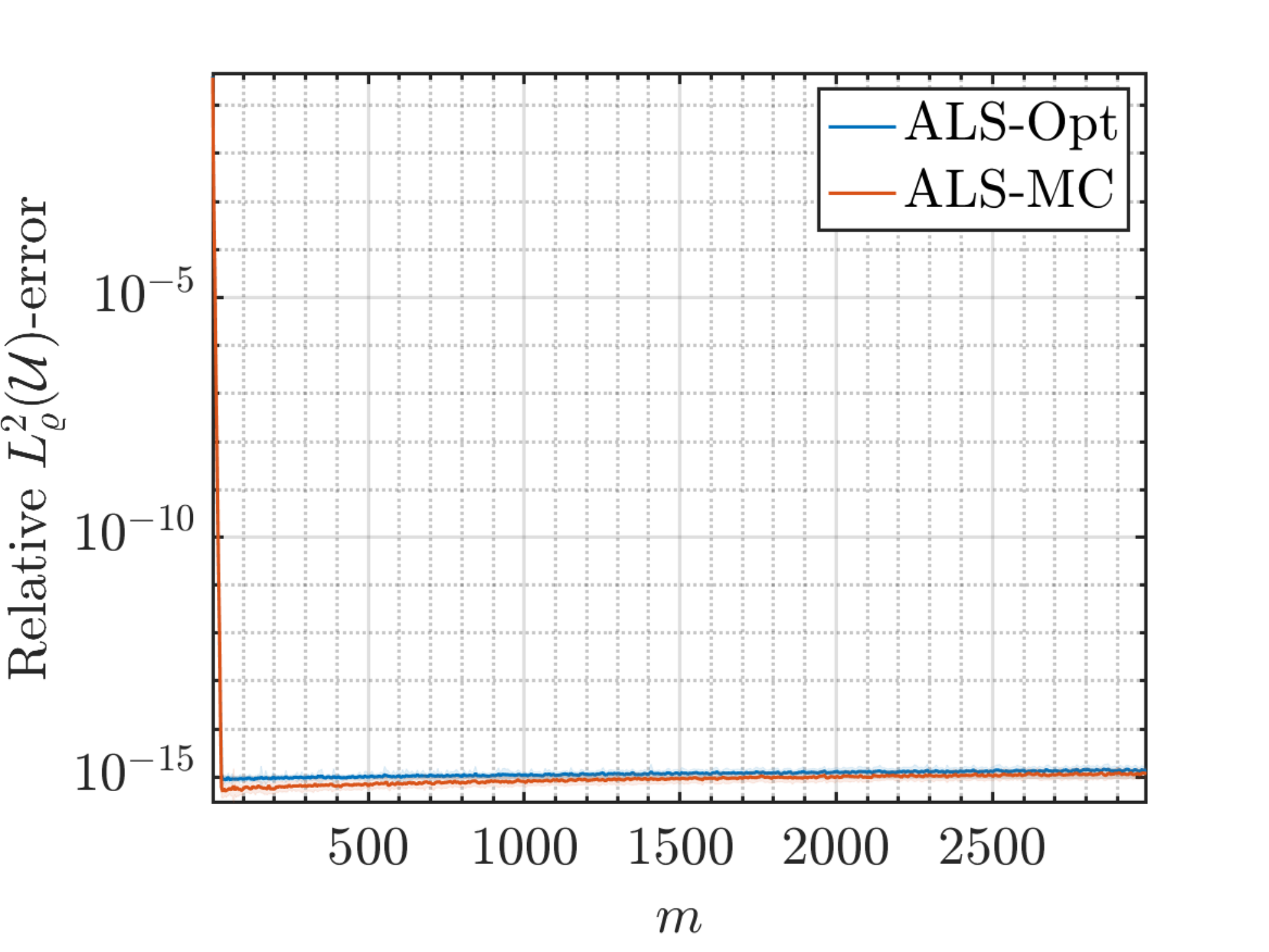}&
\includegraphics[width = \errplotimg]{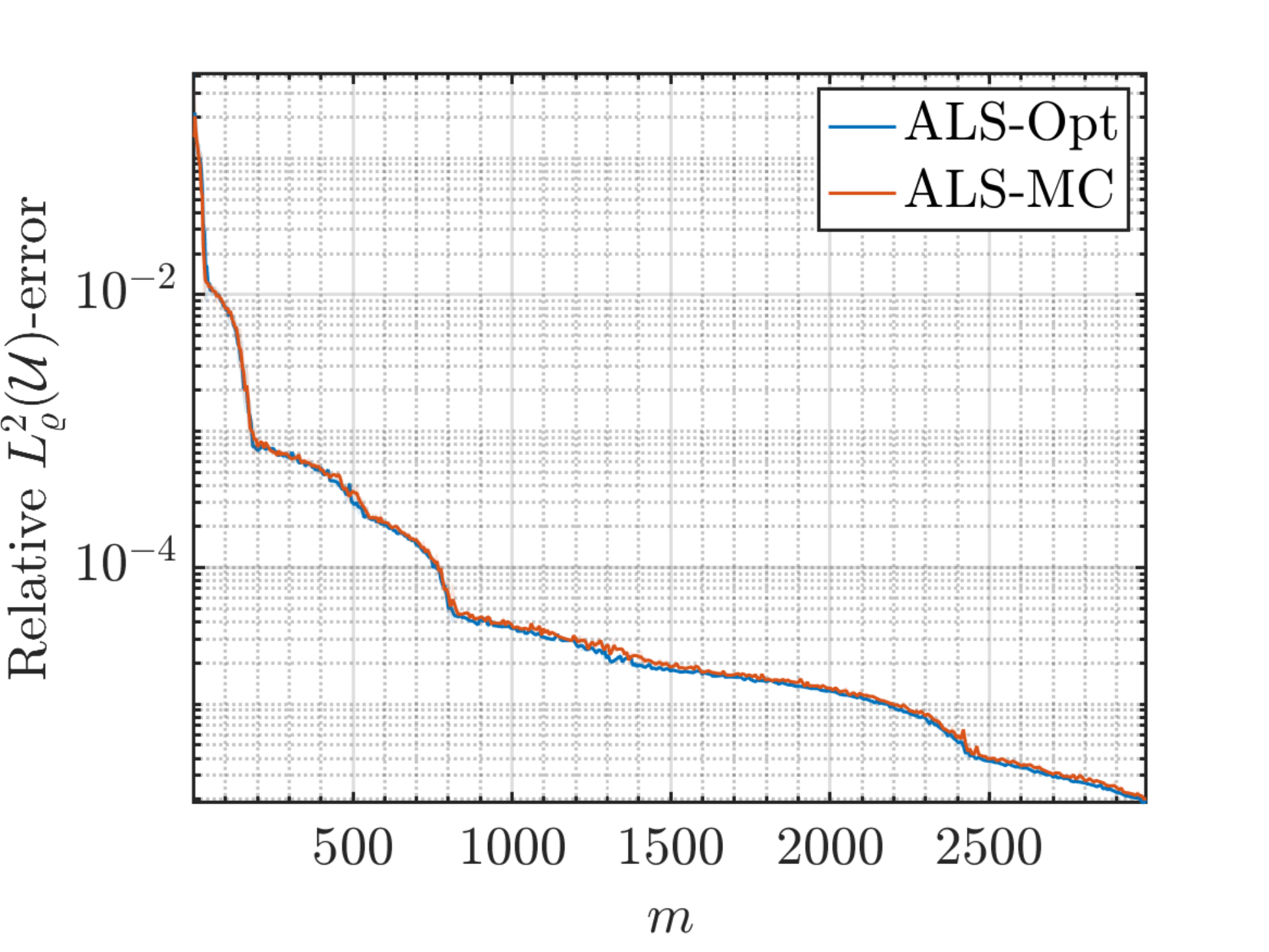}
&
\includegraphics[width = \errplotimg]{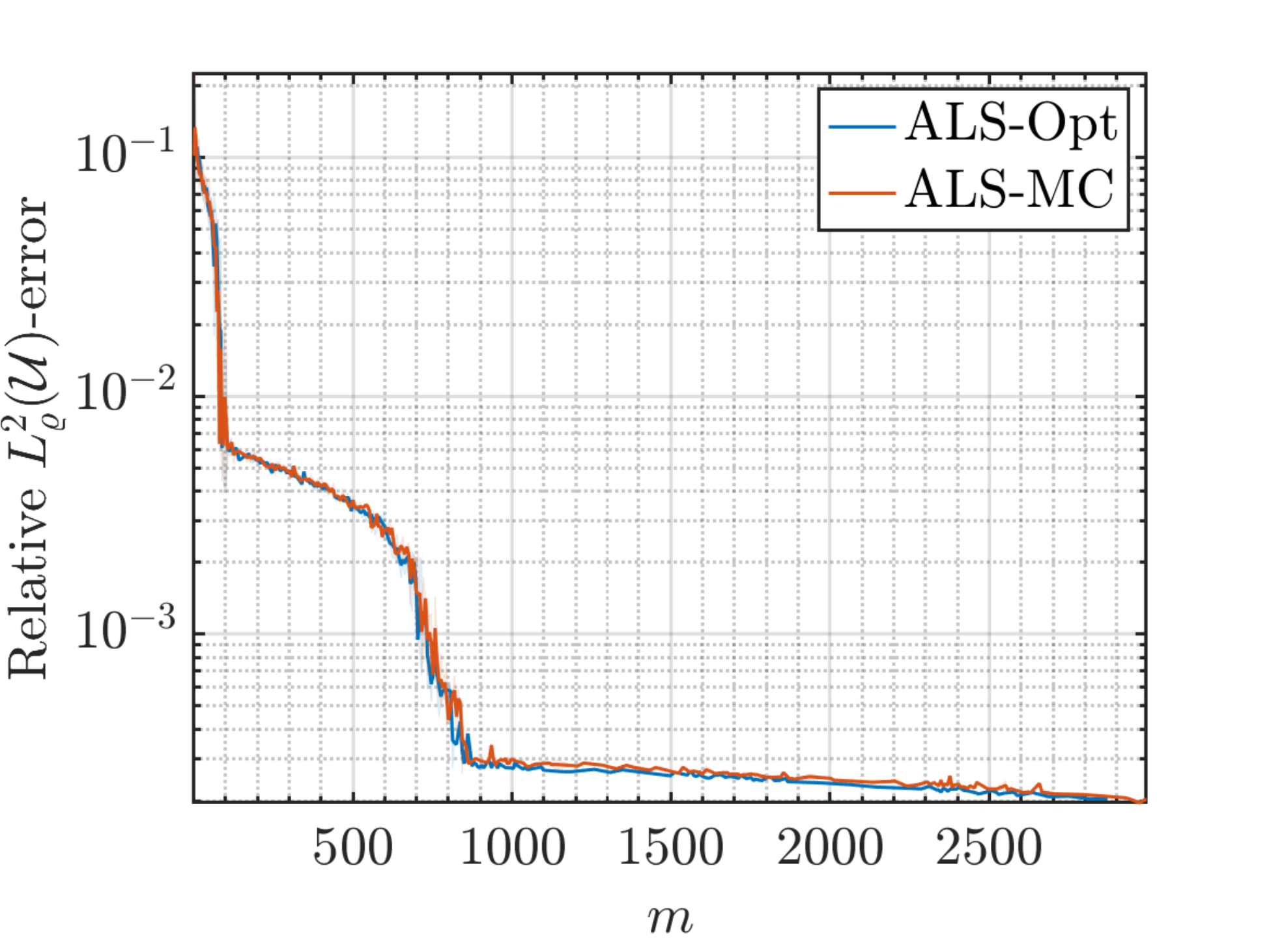}
\\[\errplotgraphsp]
\includegraphics[width = \errplotimg]{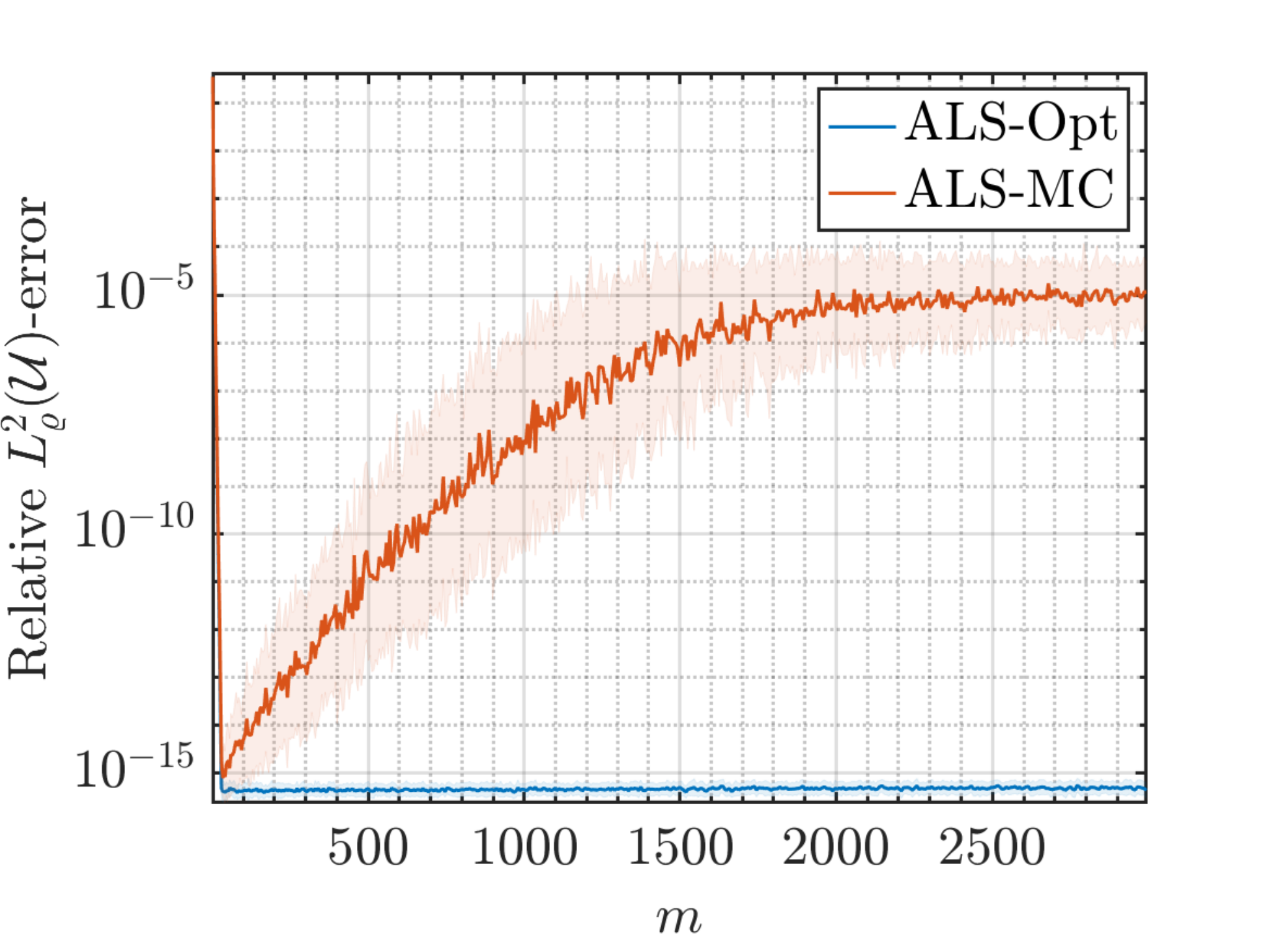}
&
\includegraphics[width = \errplotimg]{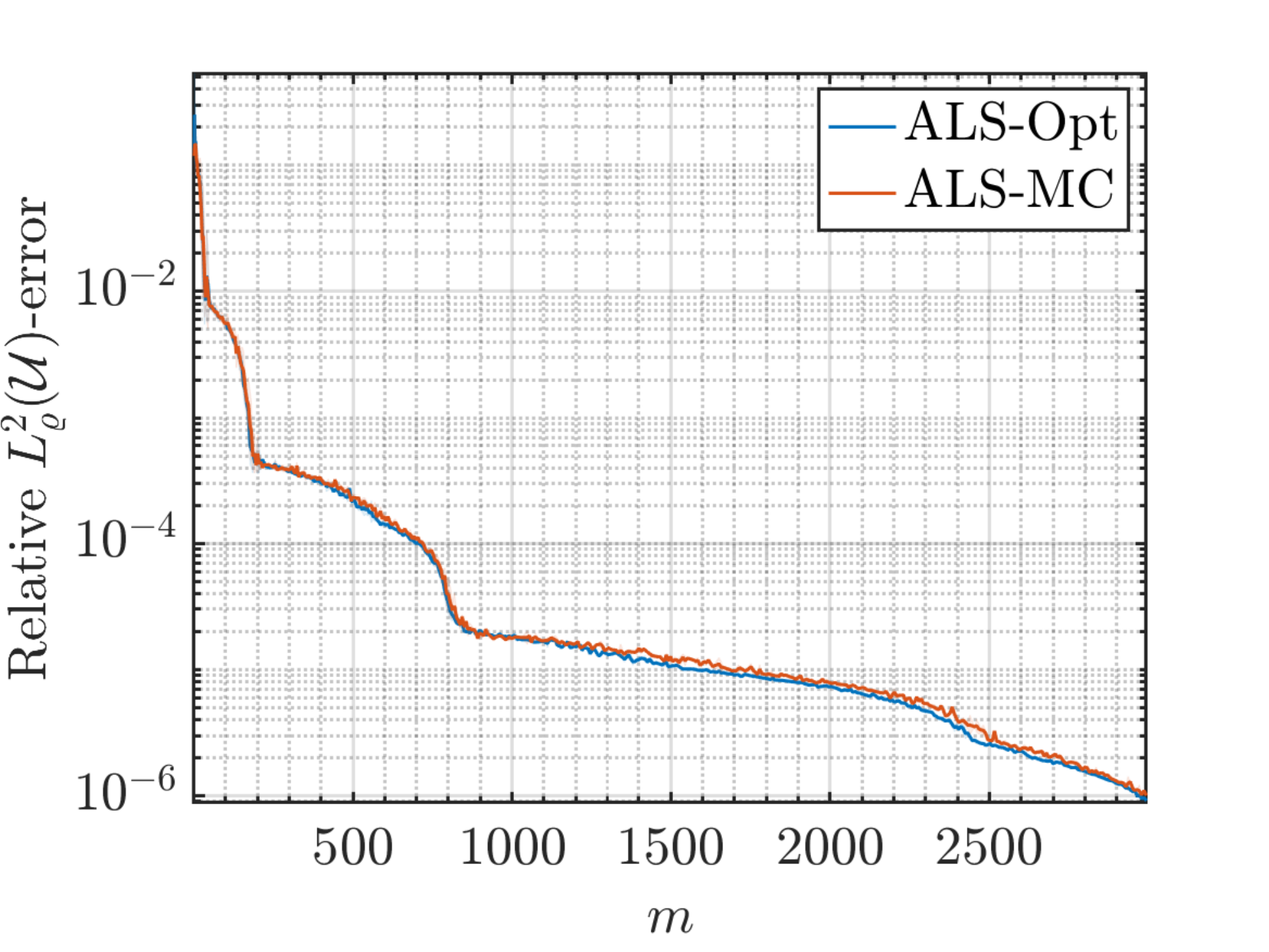}
&
\includegraphics[width = \errplotimg]{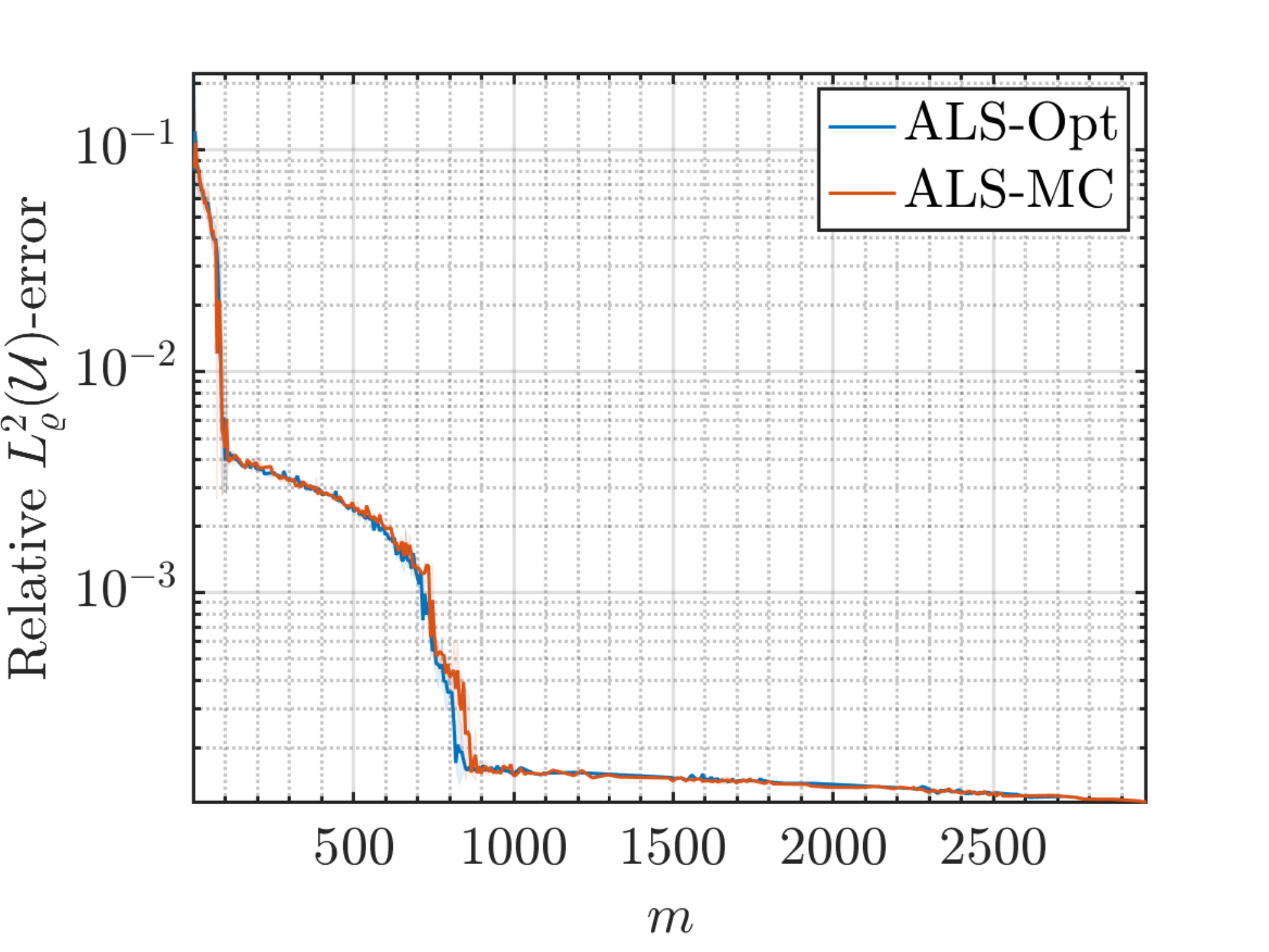}
\\[\errplotgraphsp]
\includegraphics[width = \errplotimg]{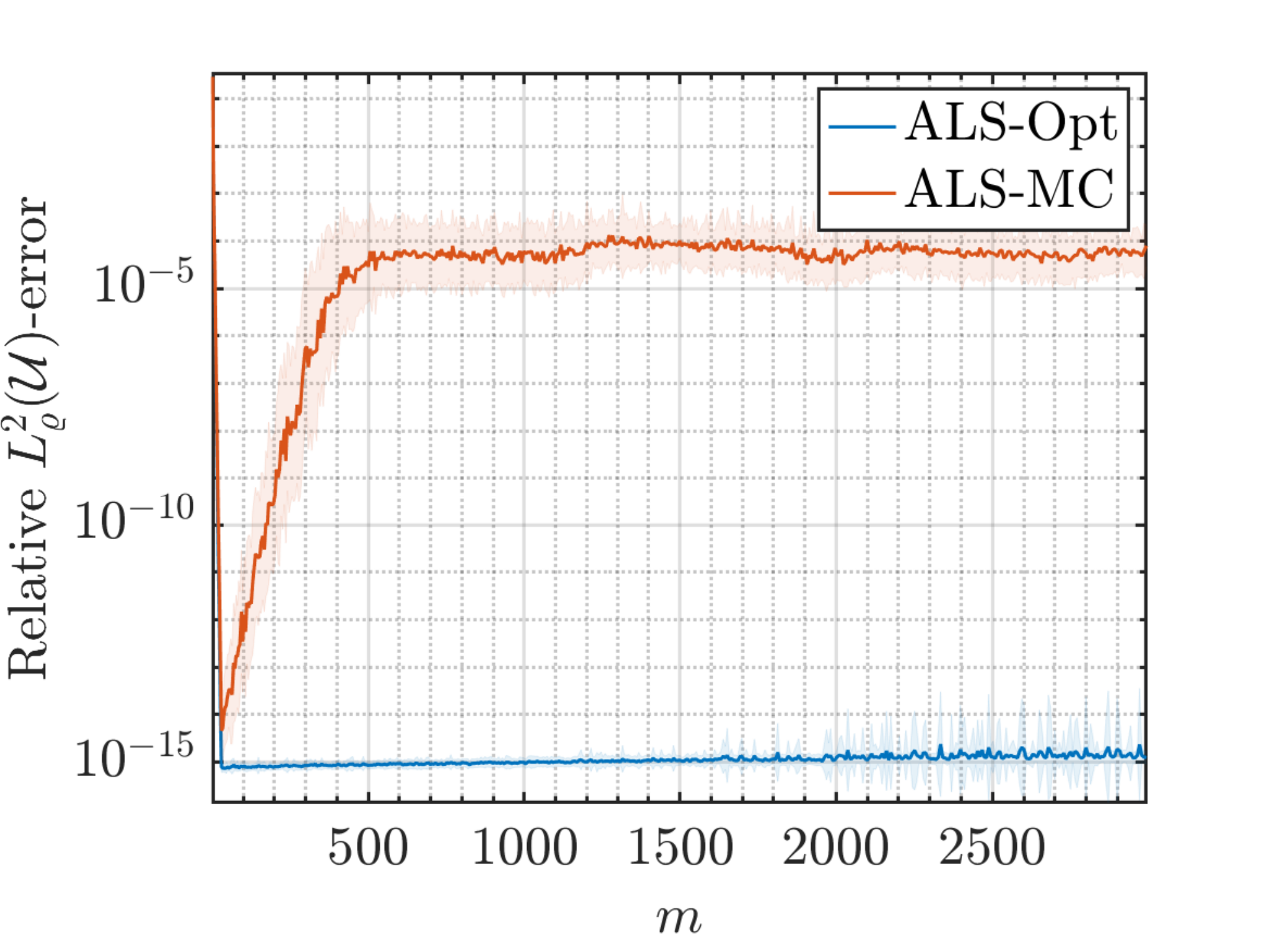}
&
\includegraphics[width = \errplotimg]{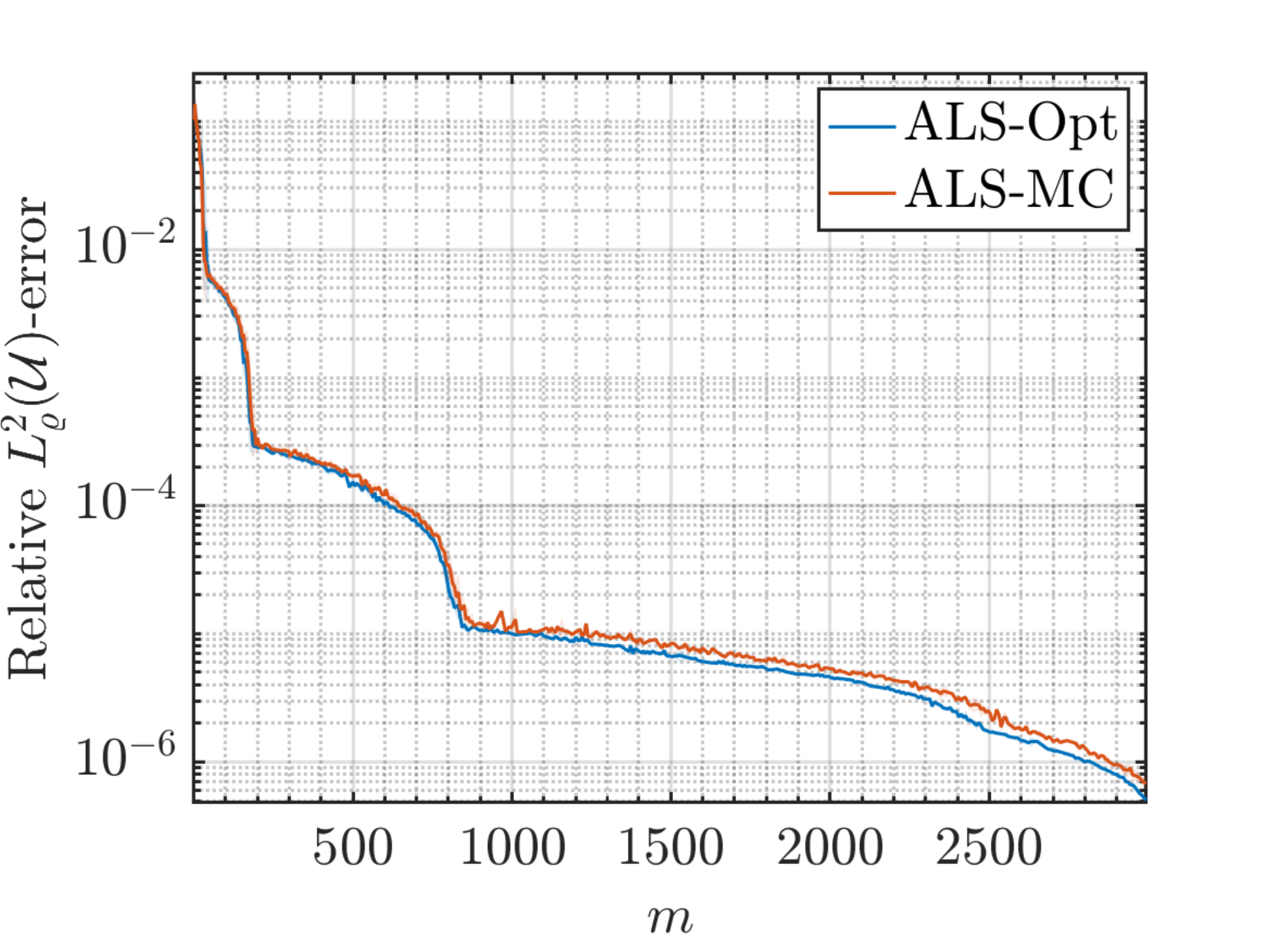}
&
\includegraphics[width = \errplotimg]{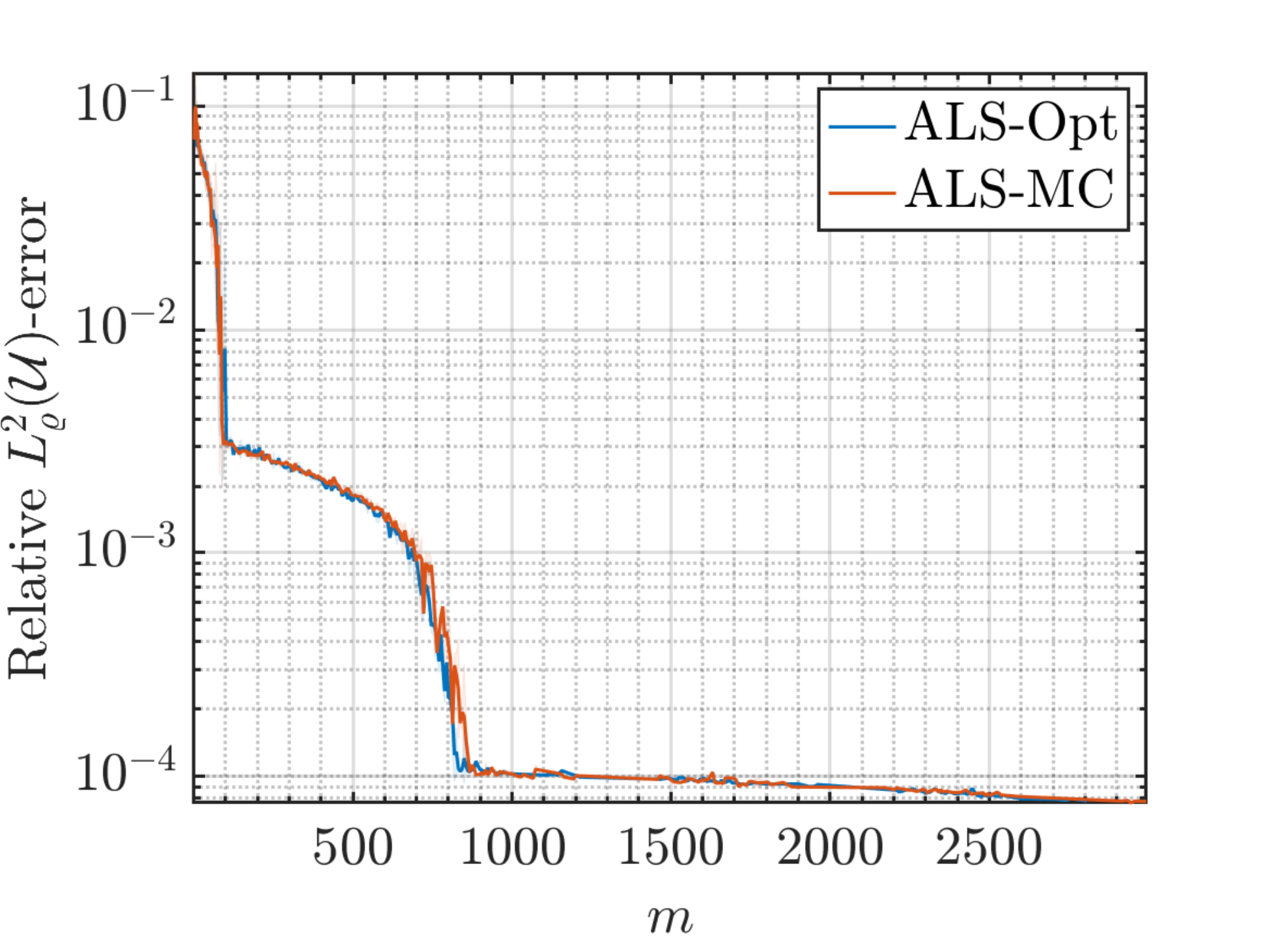}

\\[\errplottextsp]
$d = 1$ & $d = 8$ & $d = 16$
\end{tabular}
\end{small}
\end{center}
\caption{Comparison of sampling from the Chebyshev measure (top row), uniform measure (middle row) and second kind Chebyshev measure (bottom row) for the function $f = f_1$.} 
\label{fig:figSM1}
\end{figure}

\begin{figure}[t!]
\begin{center}
\begin{small}
 \begin{tabular}{@{\hspace{0pt}}c@{\hspace{\errplotsp}}c@{\hspace{\errplotsp}}c@{\hspace{0pt}}}
\includegraphics[width = \errplotimg]{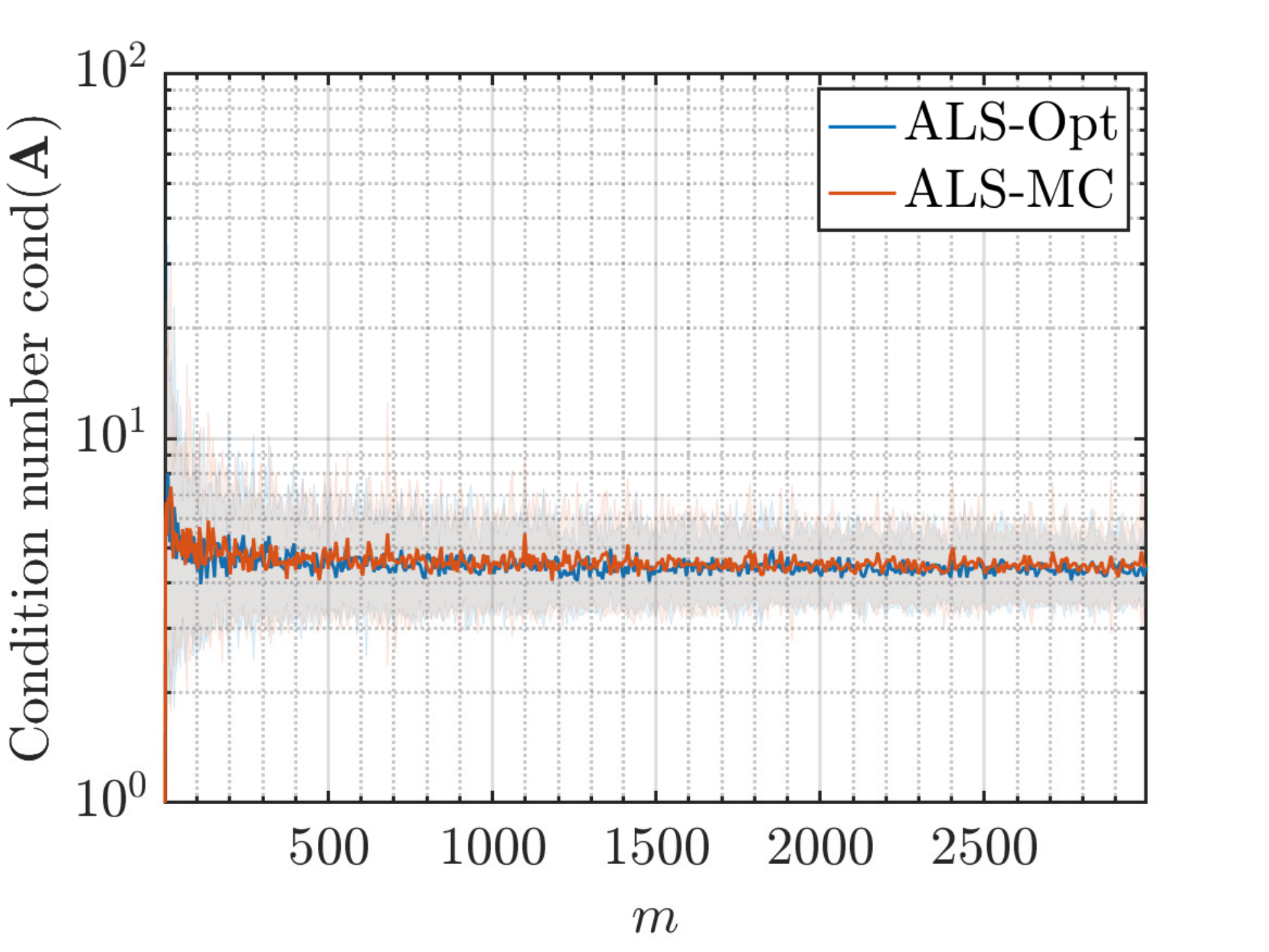}&
\includegraphics[width = \errplotimg]{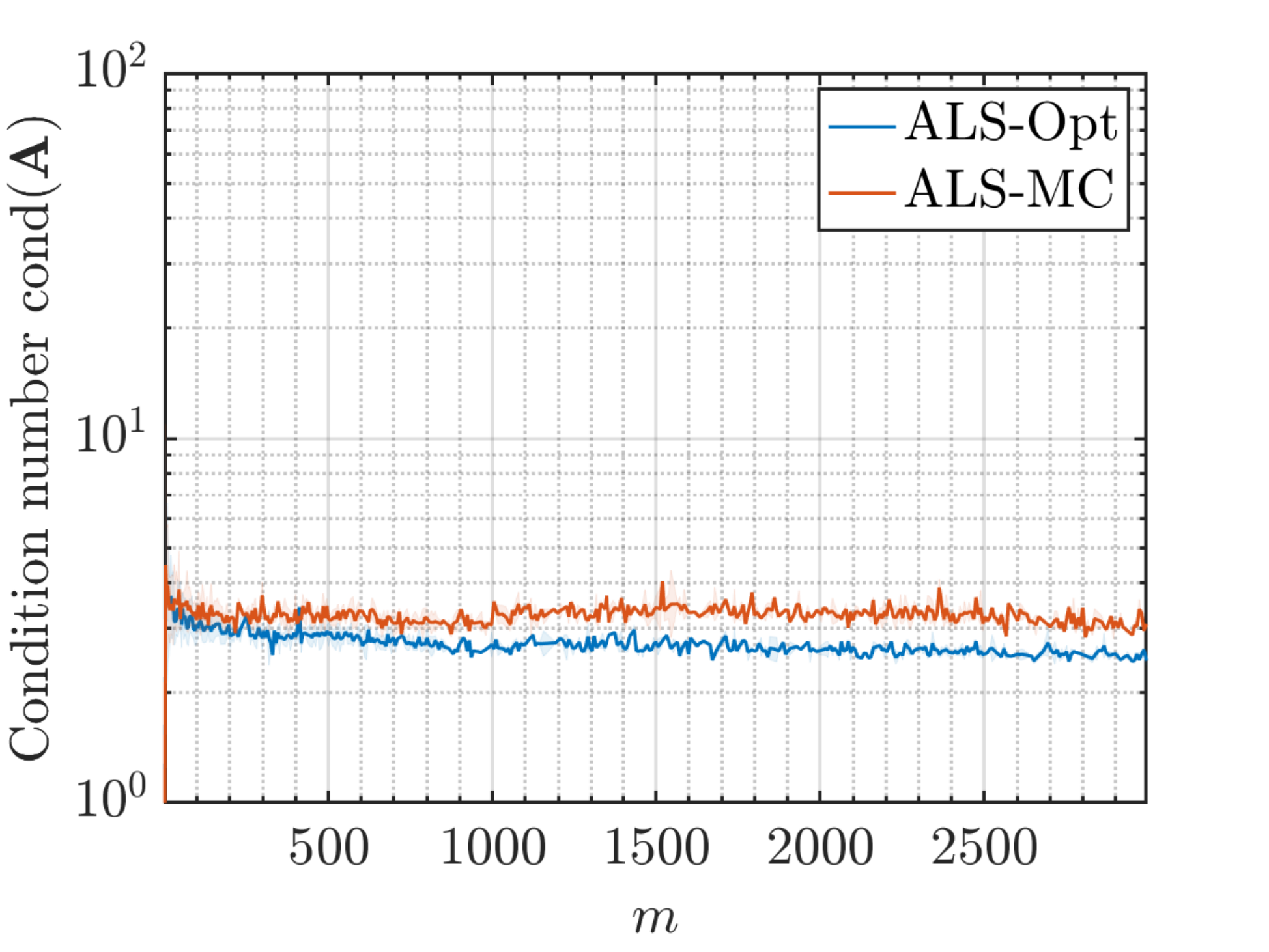}
&
\includegraphics[width = \errplotimg]{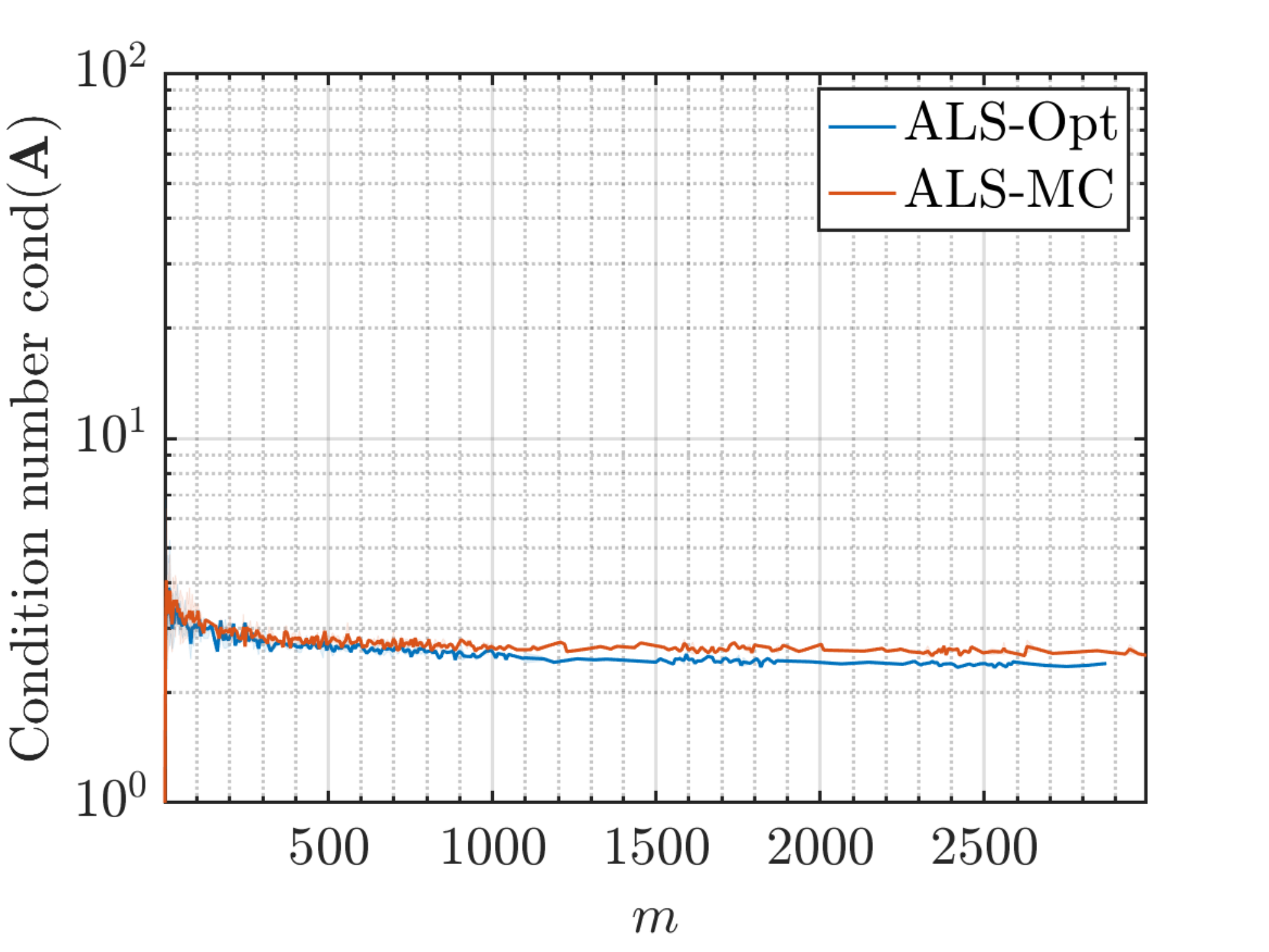}
\\[\errplotgraphsp]
\includegraphics[width = \errplotimg]{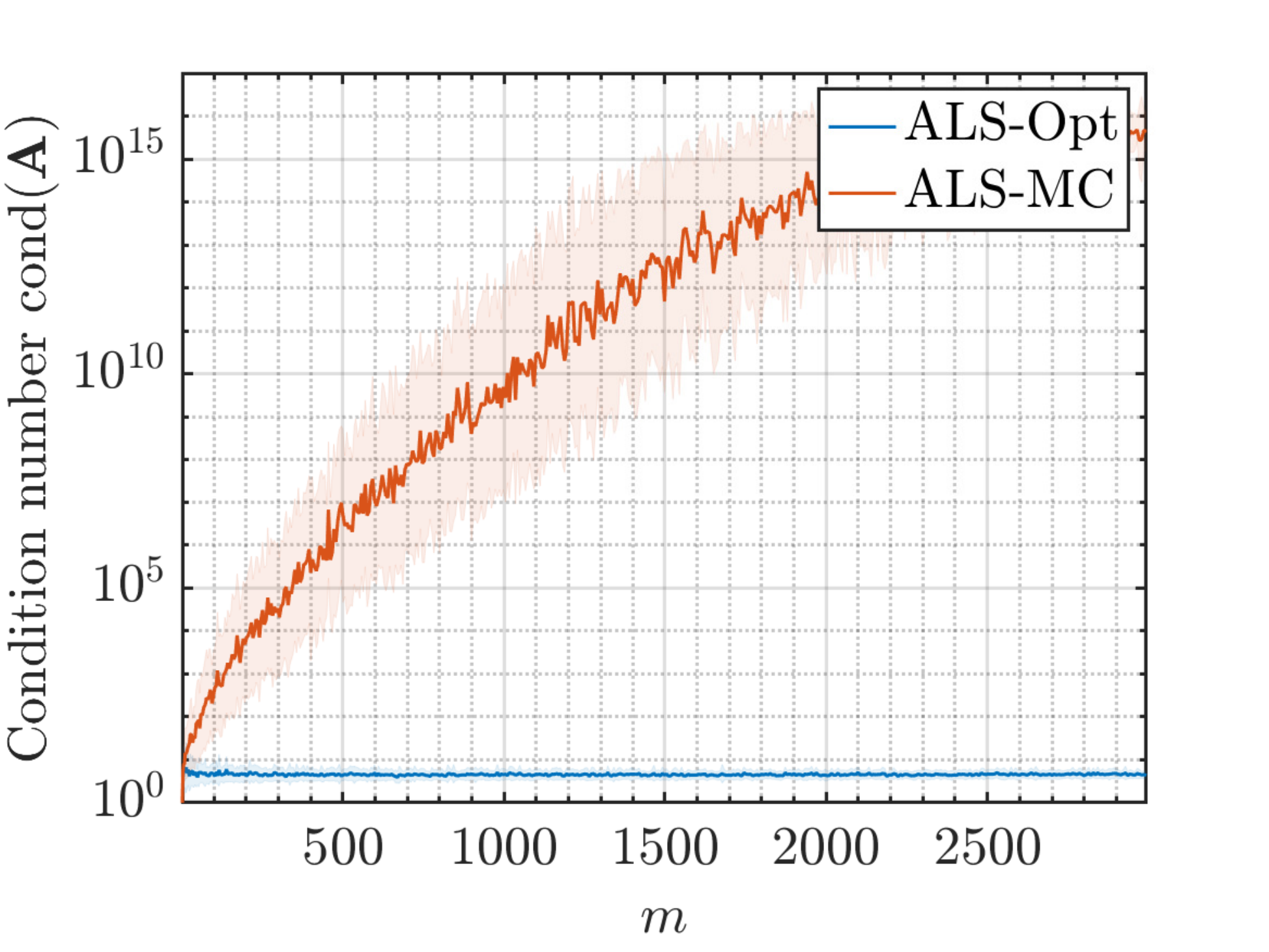}
&
\includegraphics[width = \errplotimg]{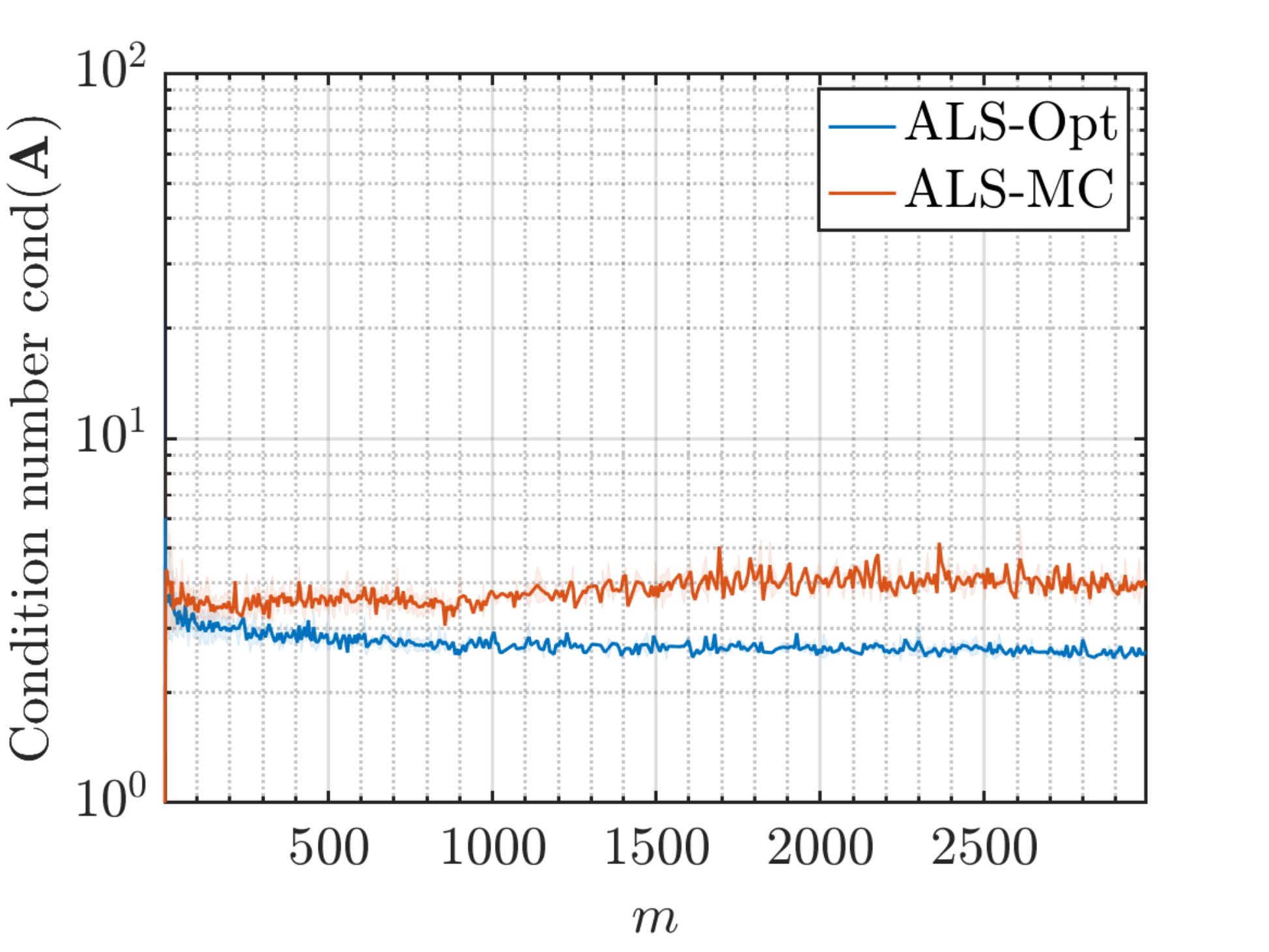}
&
\includegraphics[width = \errplotimg]{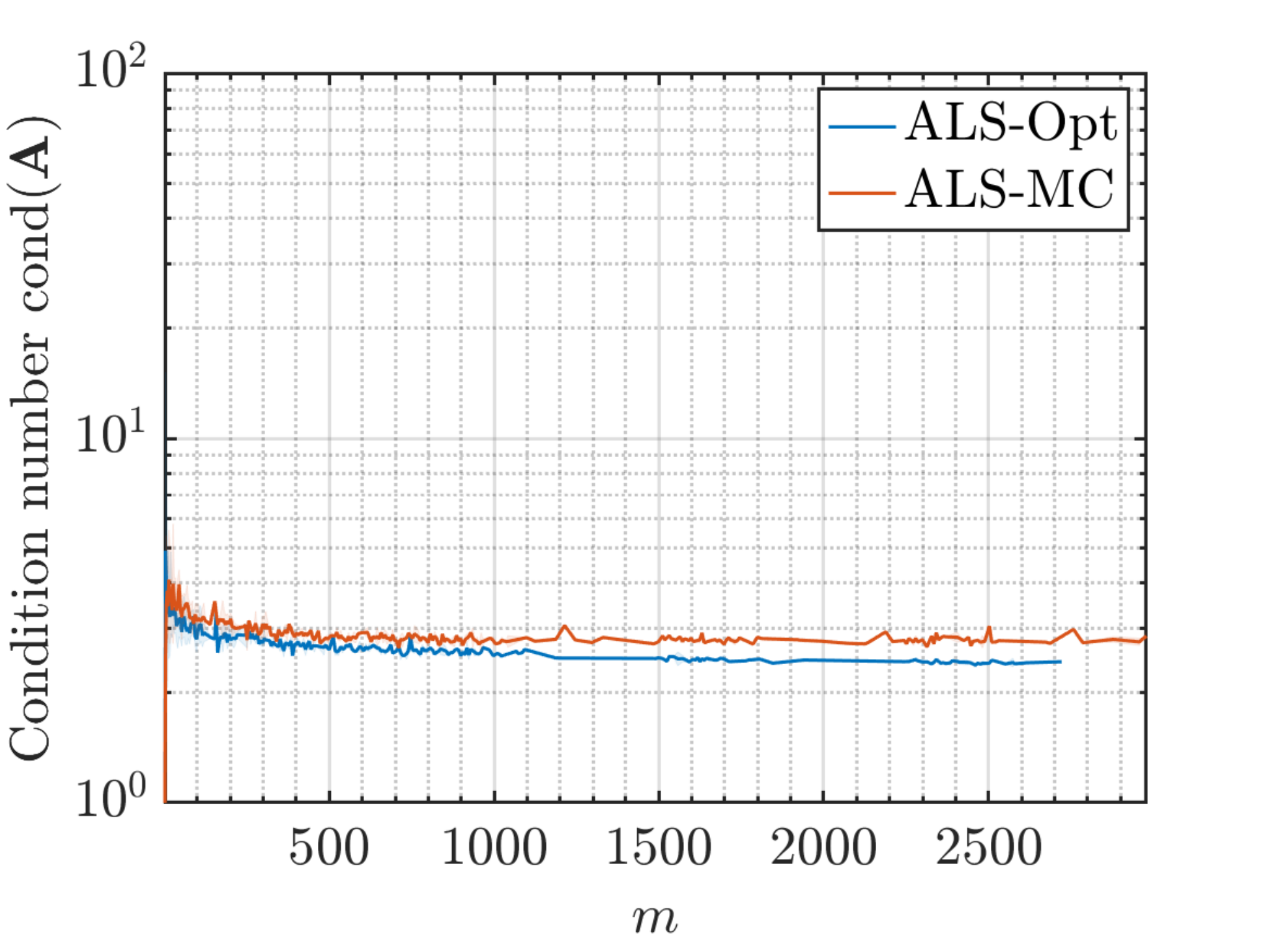}
\\[\errplotgraphsp]
\includegraphics[width = \errplotimg]{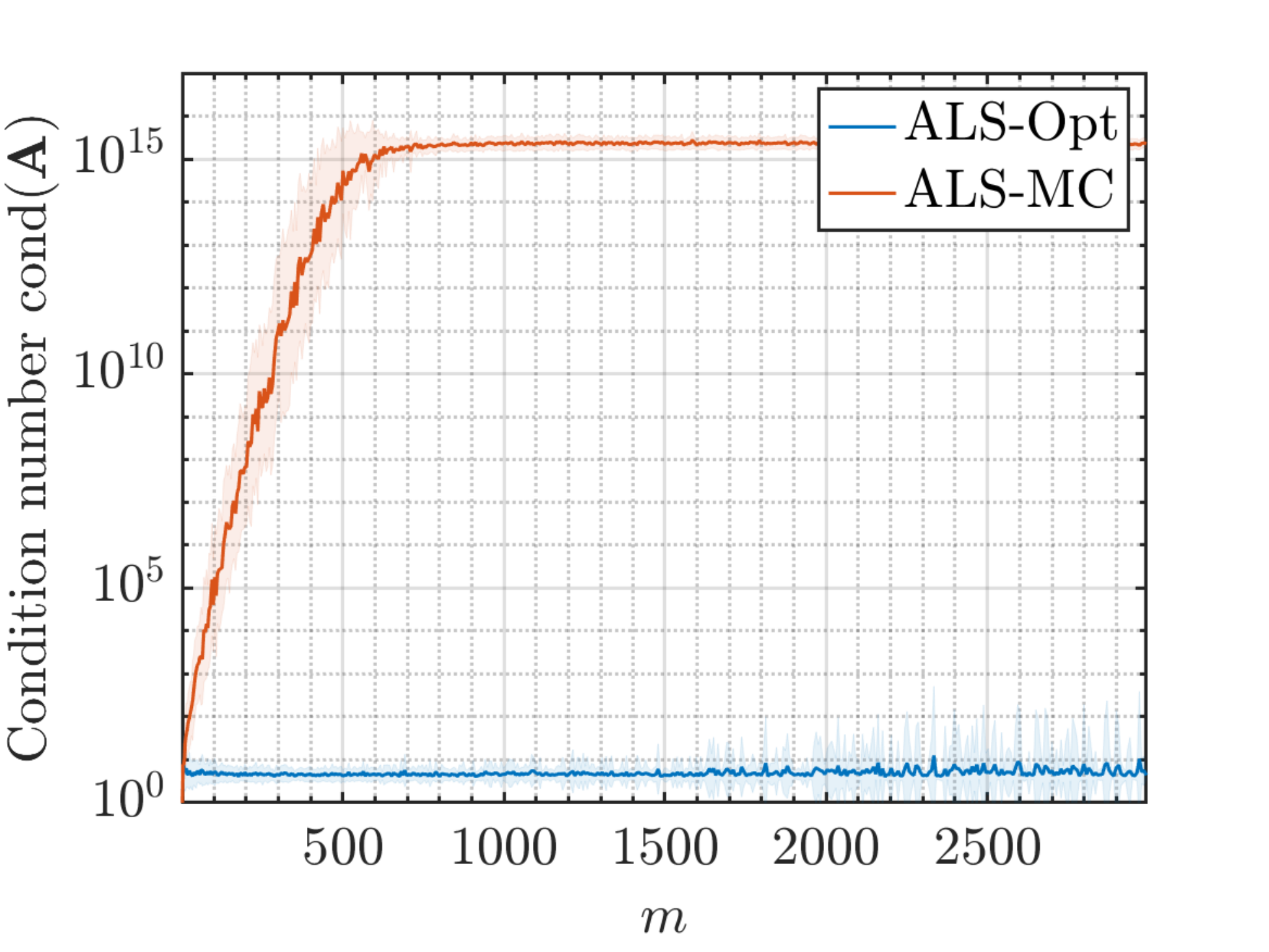}
&
\includegraphics[width = \errplotimg]{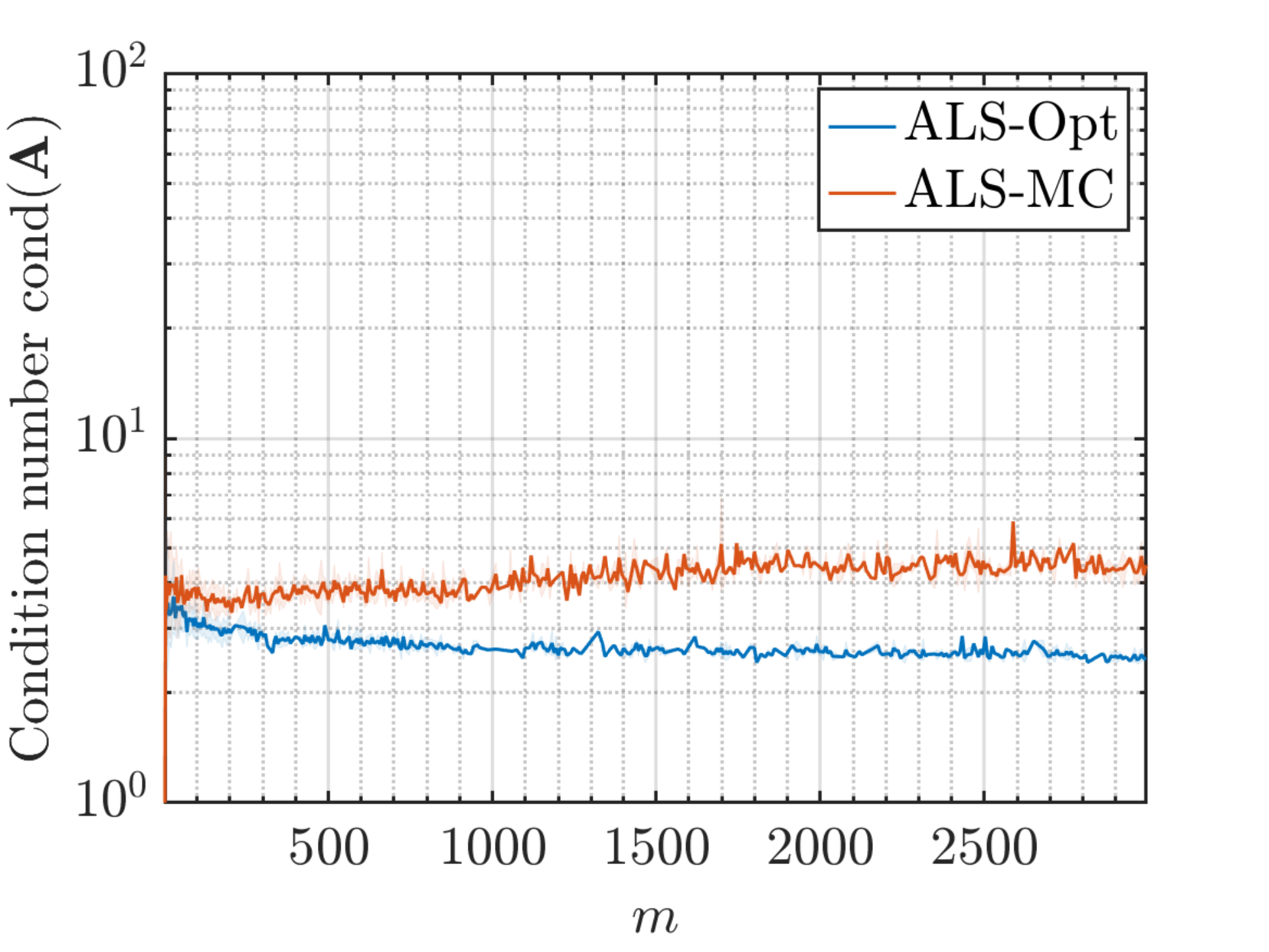}
&
\includegraphics[width = \errplotimg]{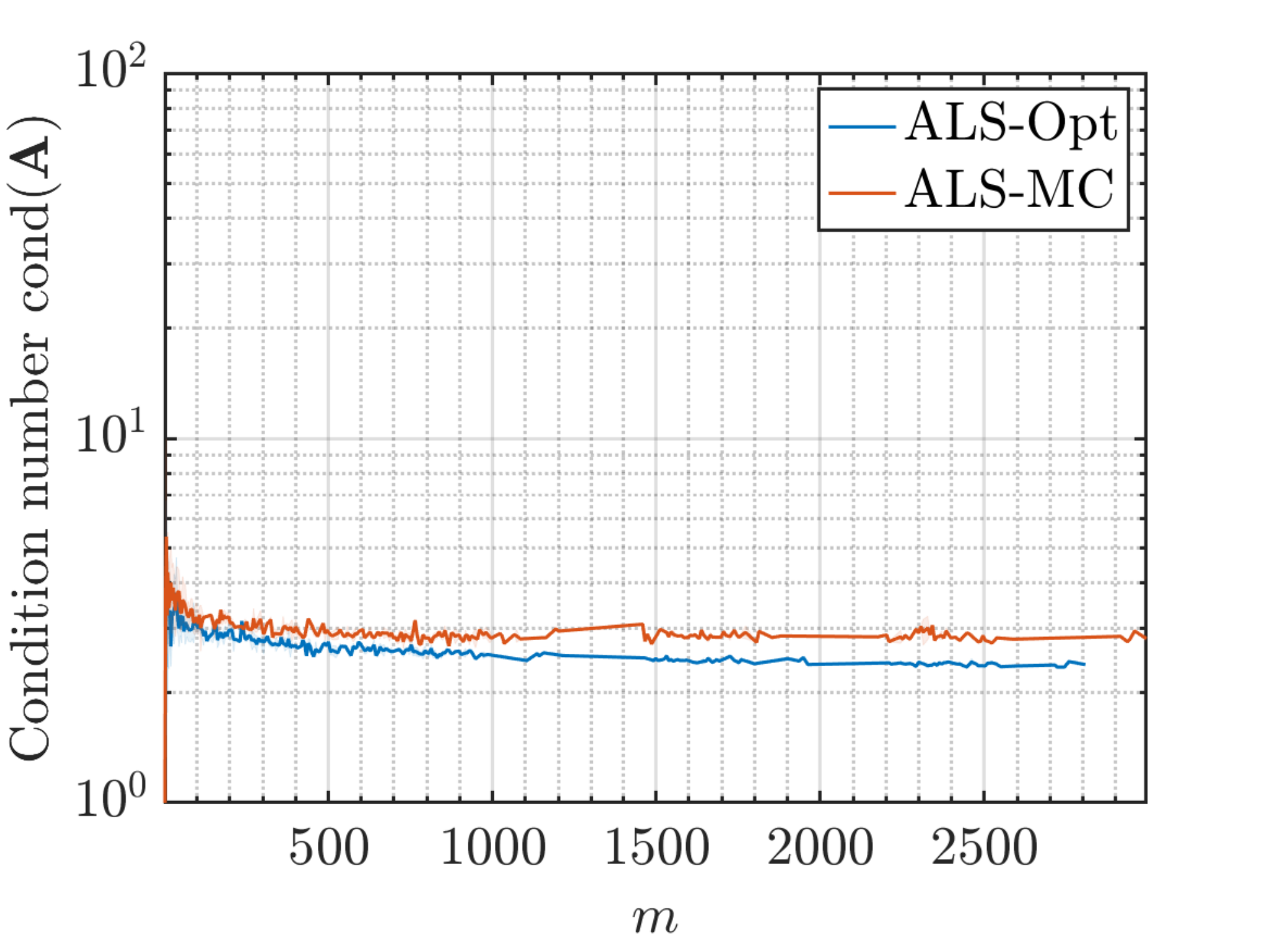}

\\[\errplottextsp]
$d = 1$ & $d = 8$ & $d = 16$
\end{tabular}
\end{small}
\end{center}
\caption{The corresponding condition numbers for the experiments in Fig.\ \ref{fig:figSM1}.} 
\label{fig:figSM2}
\end{figure}

\begin{figure}[t!]
\begin{center}
\begin{small}
 \begin{tabular}{@{\hspace{0pt}}c@{\hspace{\errplotsp}}c@{\hspace{\errplotsp}}c@{\hspace{0pt}}}
\includegraphics[width = \errplotimg]{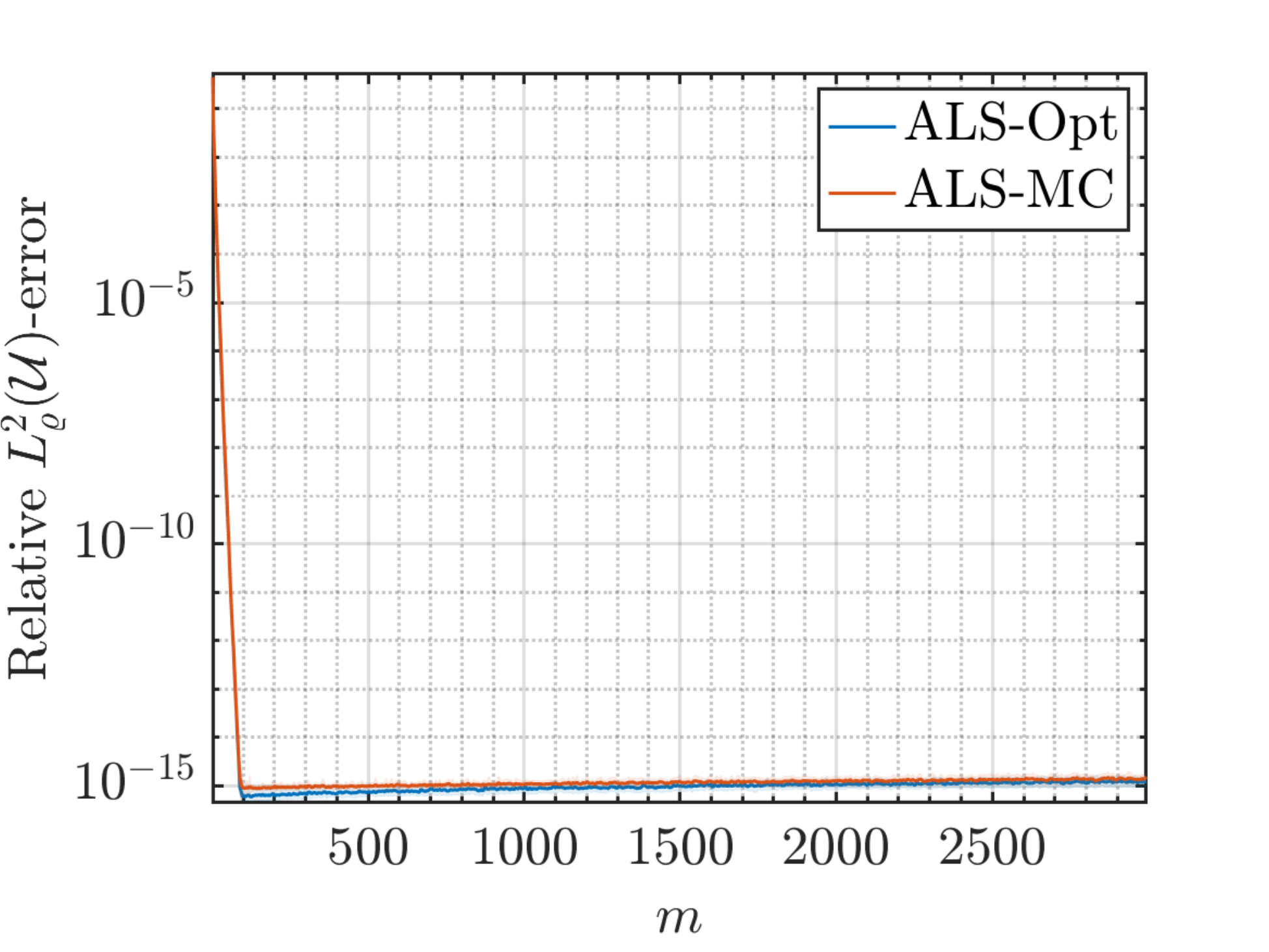}&
\includegraphics[width = \errplotimg]{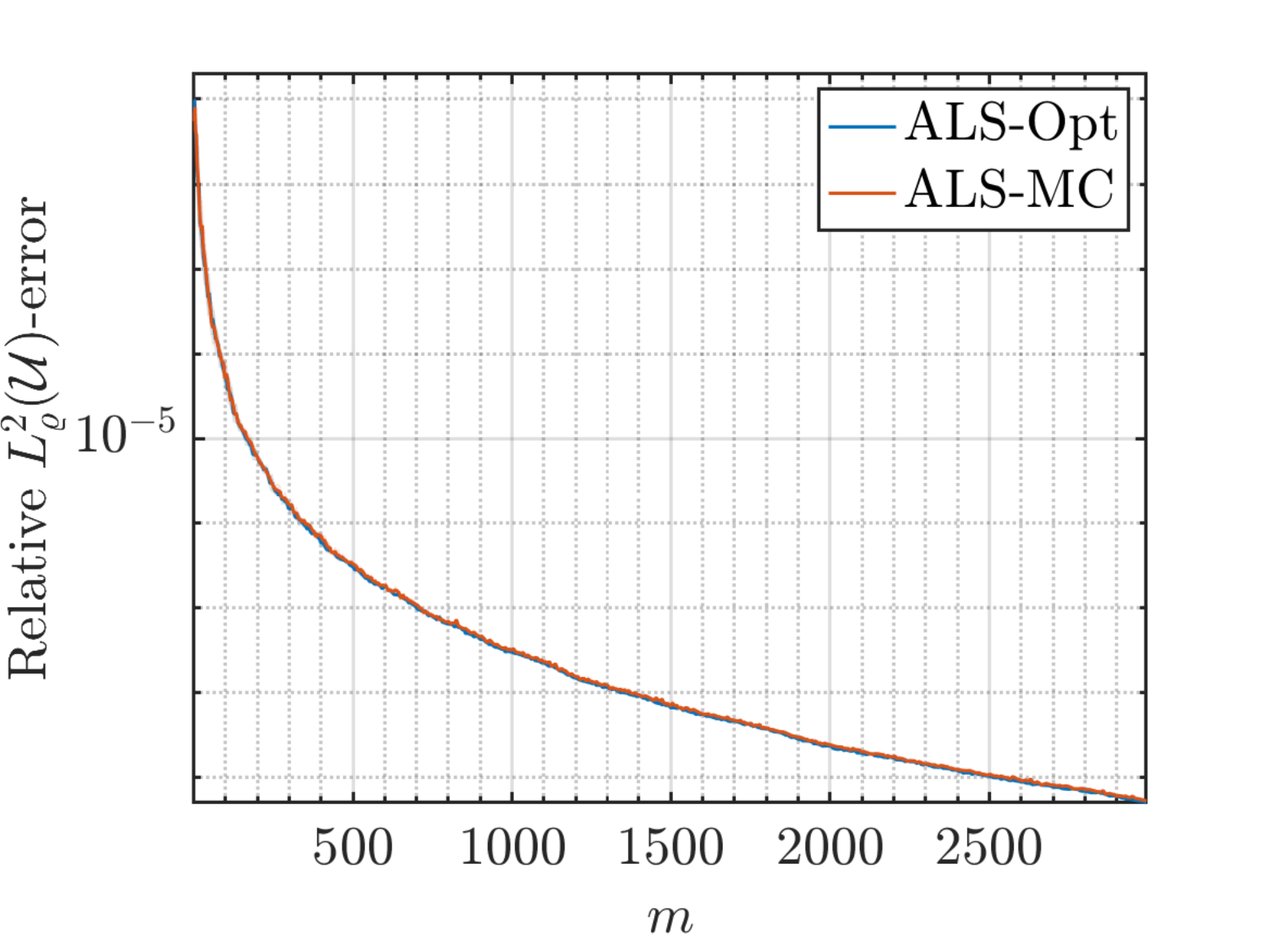}
&
\includegraphics[width = \errplotimg]{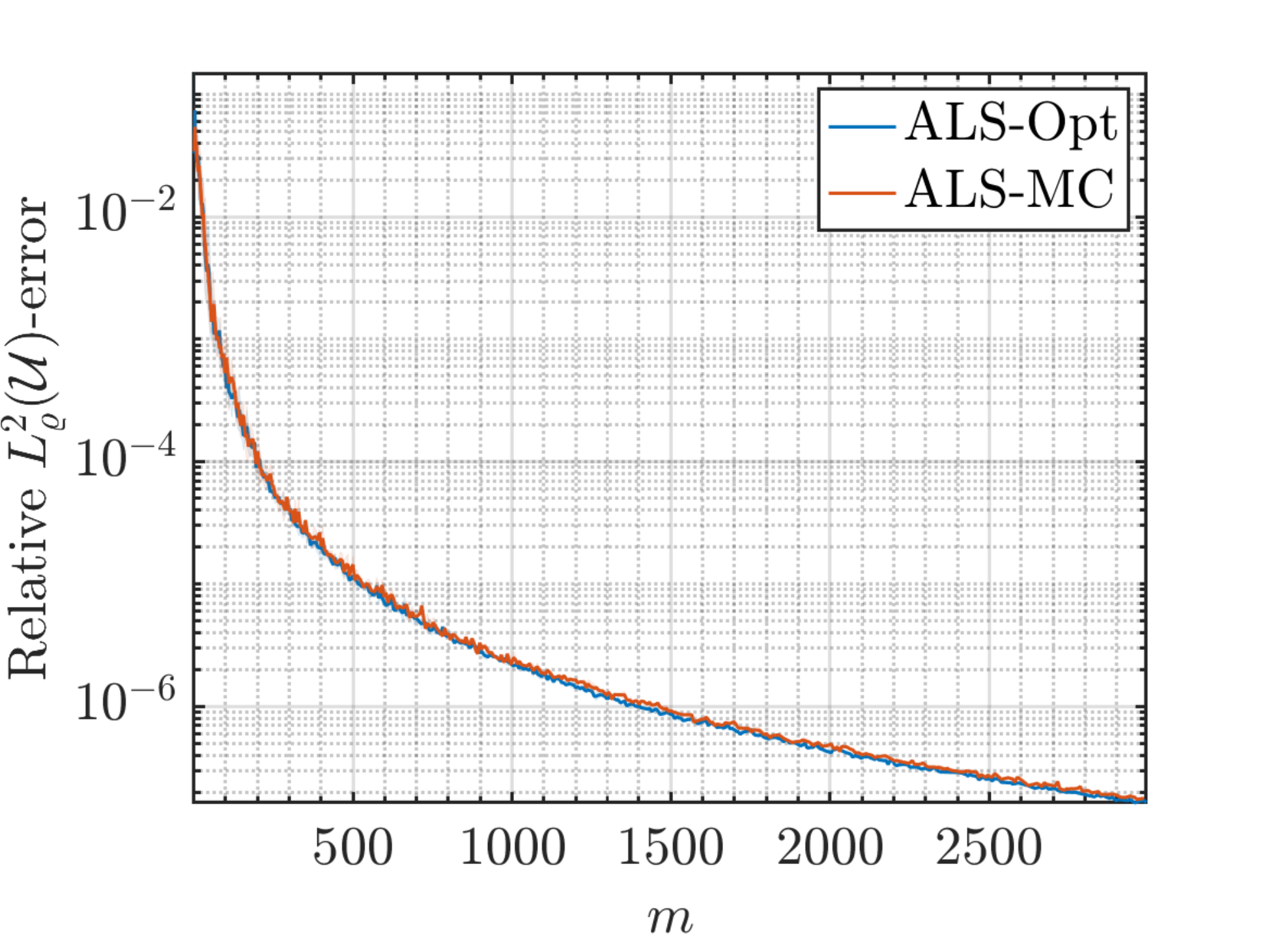}
\\[\errplotgraphsp]
\includegraphics[width = \errplotimg]{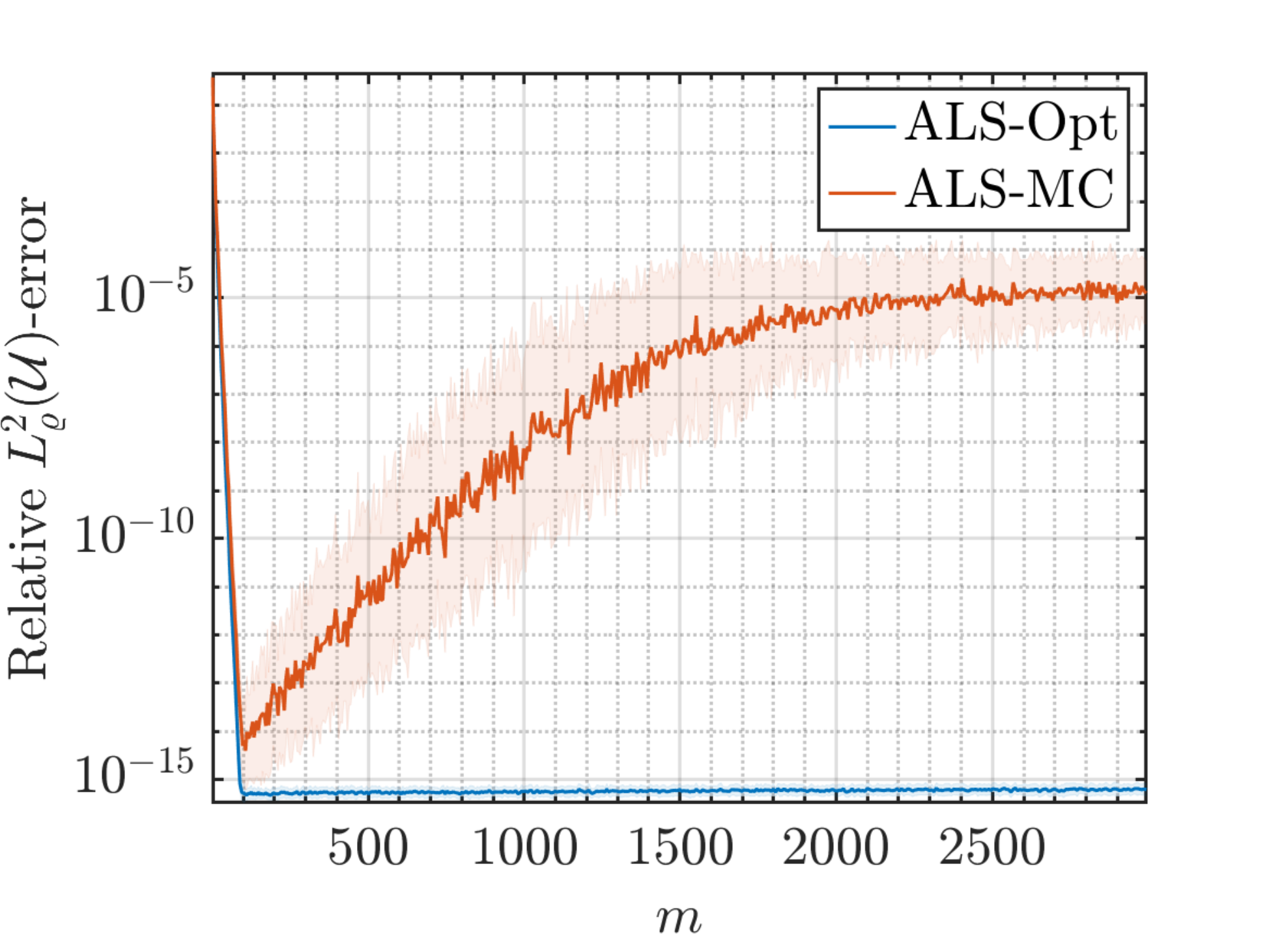}
&
\includegraphics[width = \errplotimg]{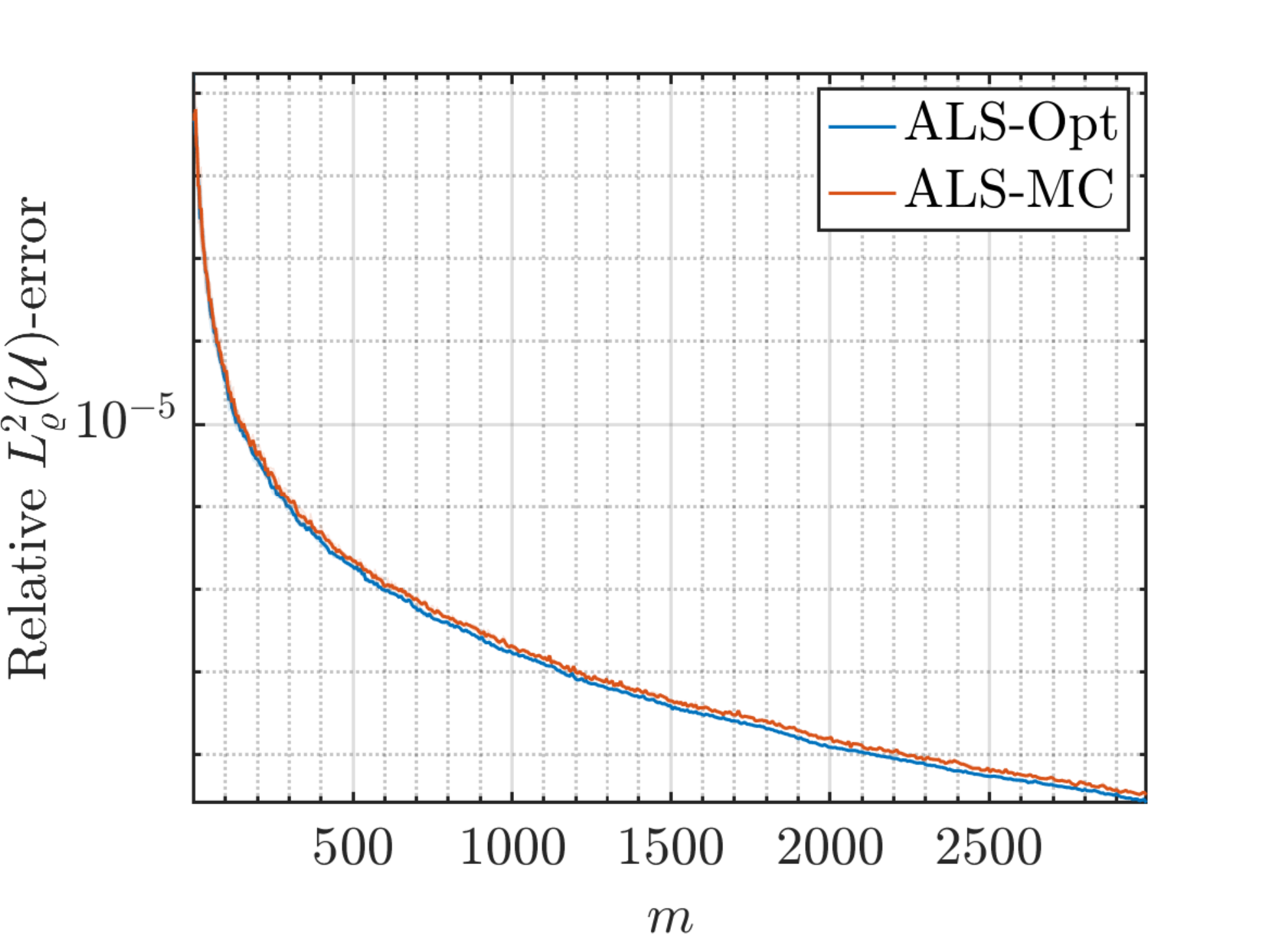}
&
\includegraphics[width = \errplotimg]{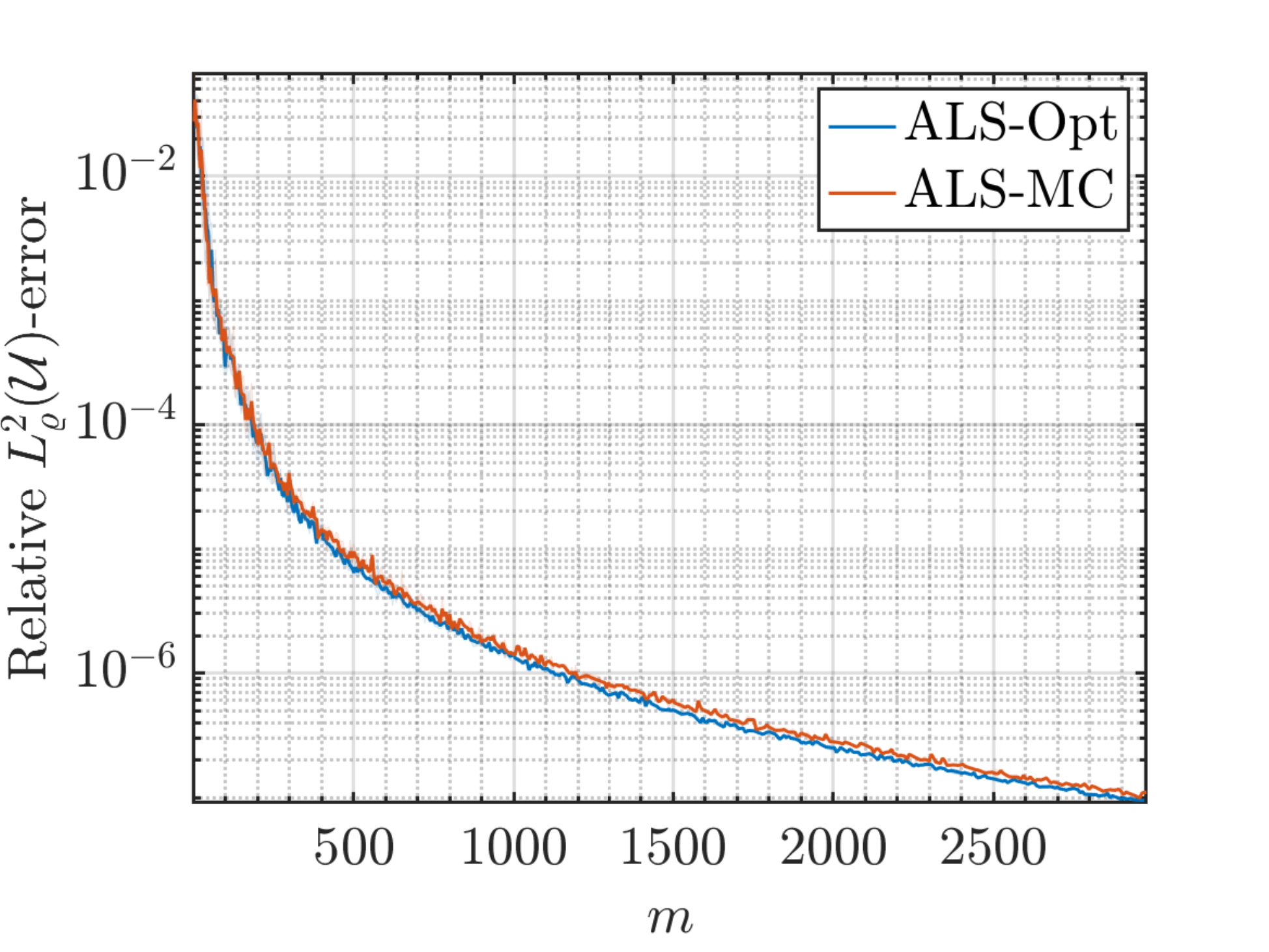}
\\[\errplotgraphsp]
\includegraphics[width = \errplotimg]{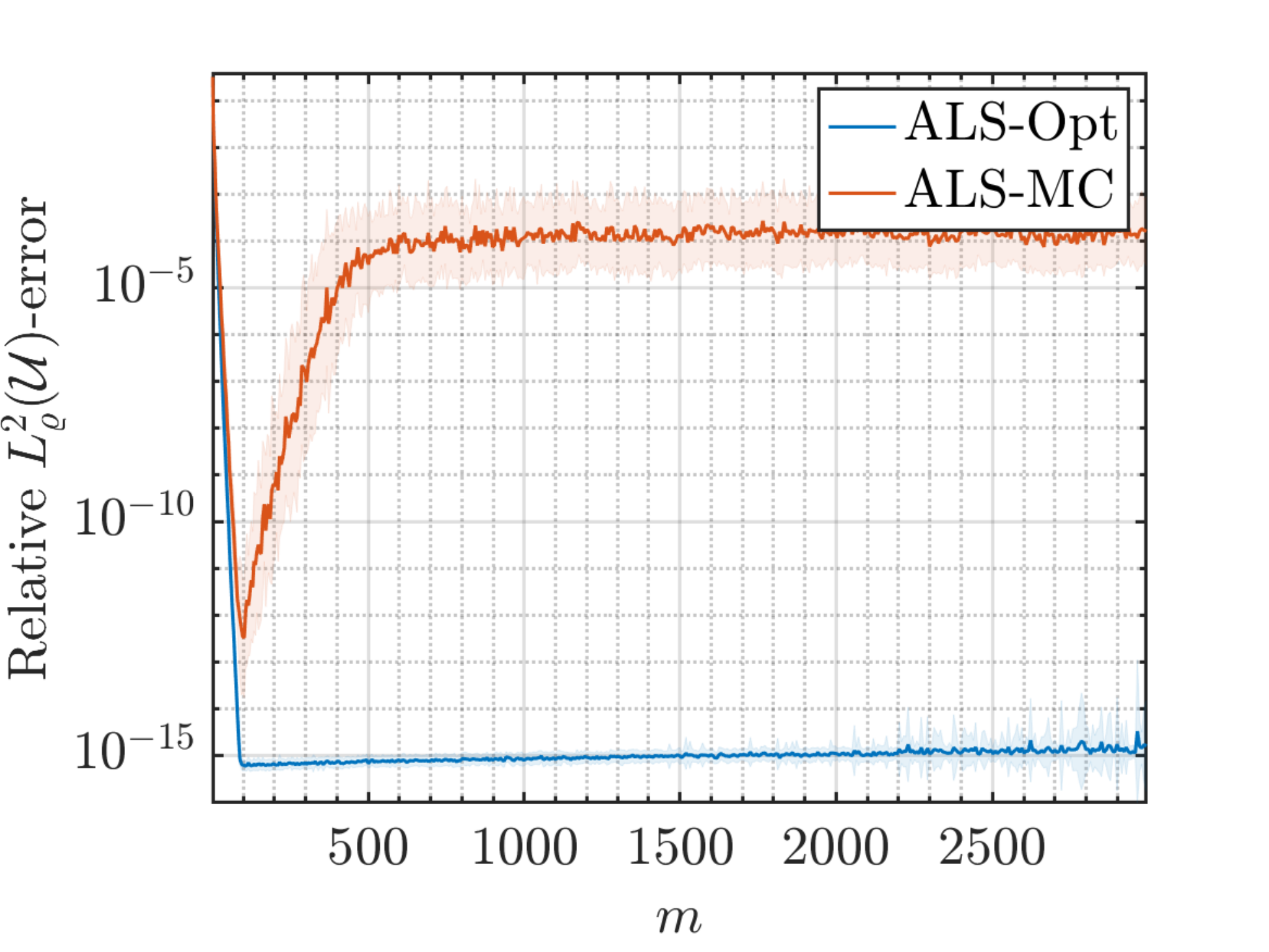}
&
\includegraphics[width = \errplotimg]{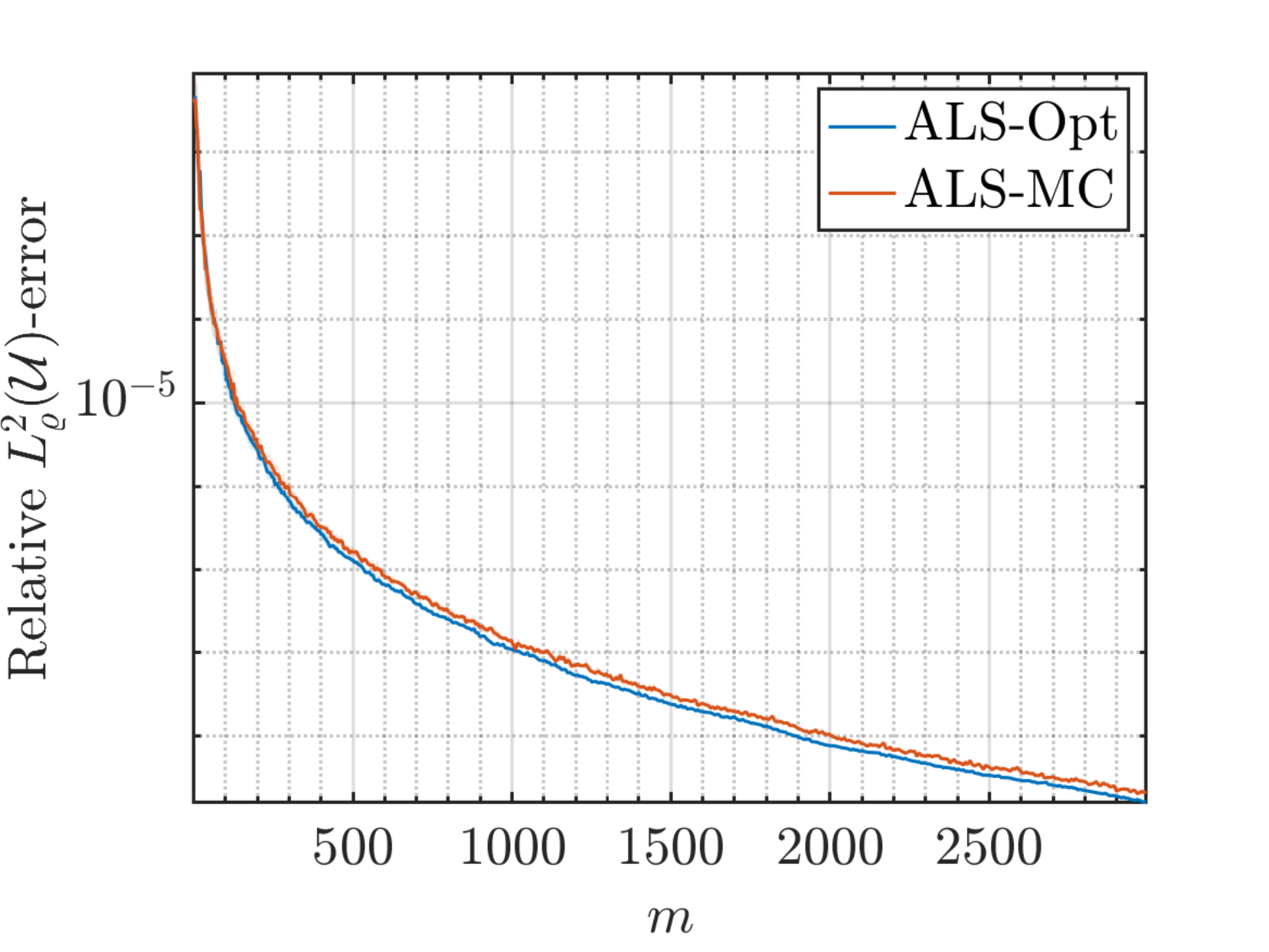}
&
\includegraphics[width = \errplotimg]{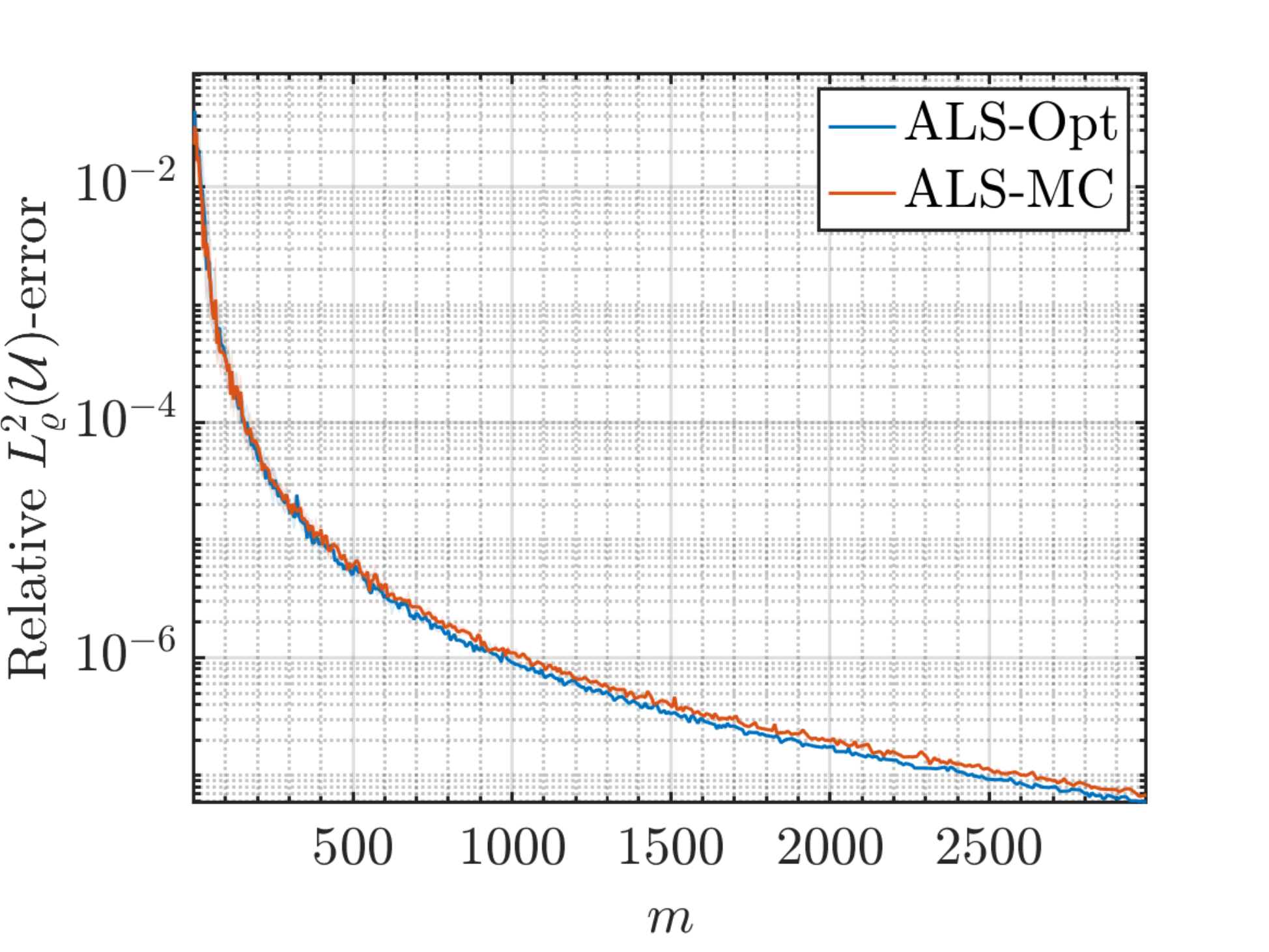}

\\[\errplottextsp]
$d = 1$ & $d = 8$ & $d = 16$
\end{tabular}
\end{small}
\end{center}
\caption{The same as Fig.\ \ref{fig:figSM1}, except with $f = f_2$.} 
\label{fig:figSM3}
\end{figure}

\begin{figure}[t!]
\begin{center}
\begin{small}
 \begin{tabular}{@{\hspace{0pt}}c@{\hspace{\errplotsp}}c@{\hspace{\errplotsp}}c@{\hspace{0pt}}}
\includegraphics[width = \errplotimg]{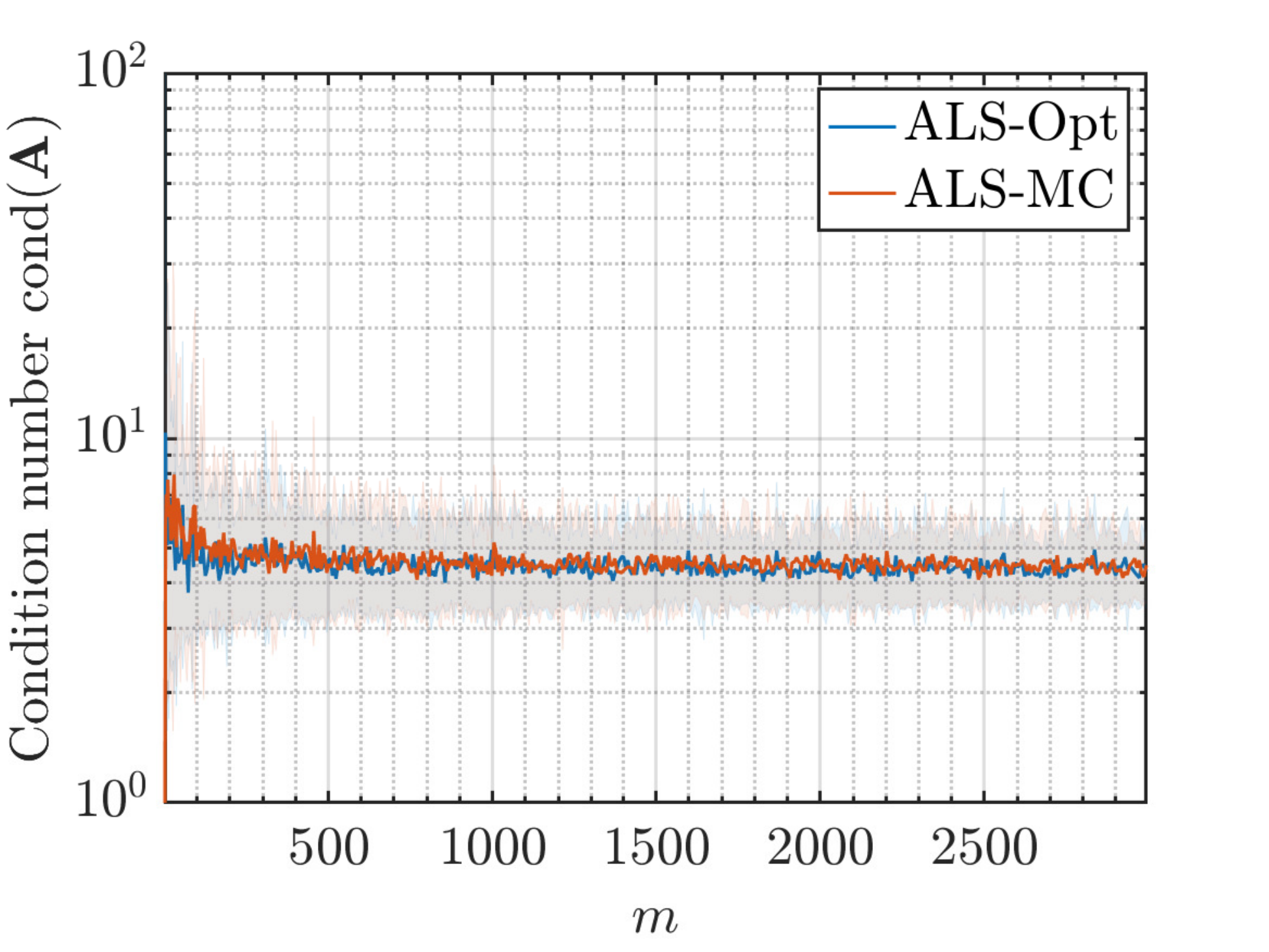}&
\includegraphics[width = \errplotimg]{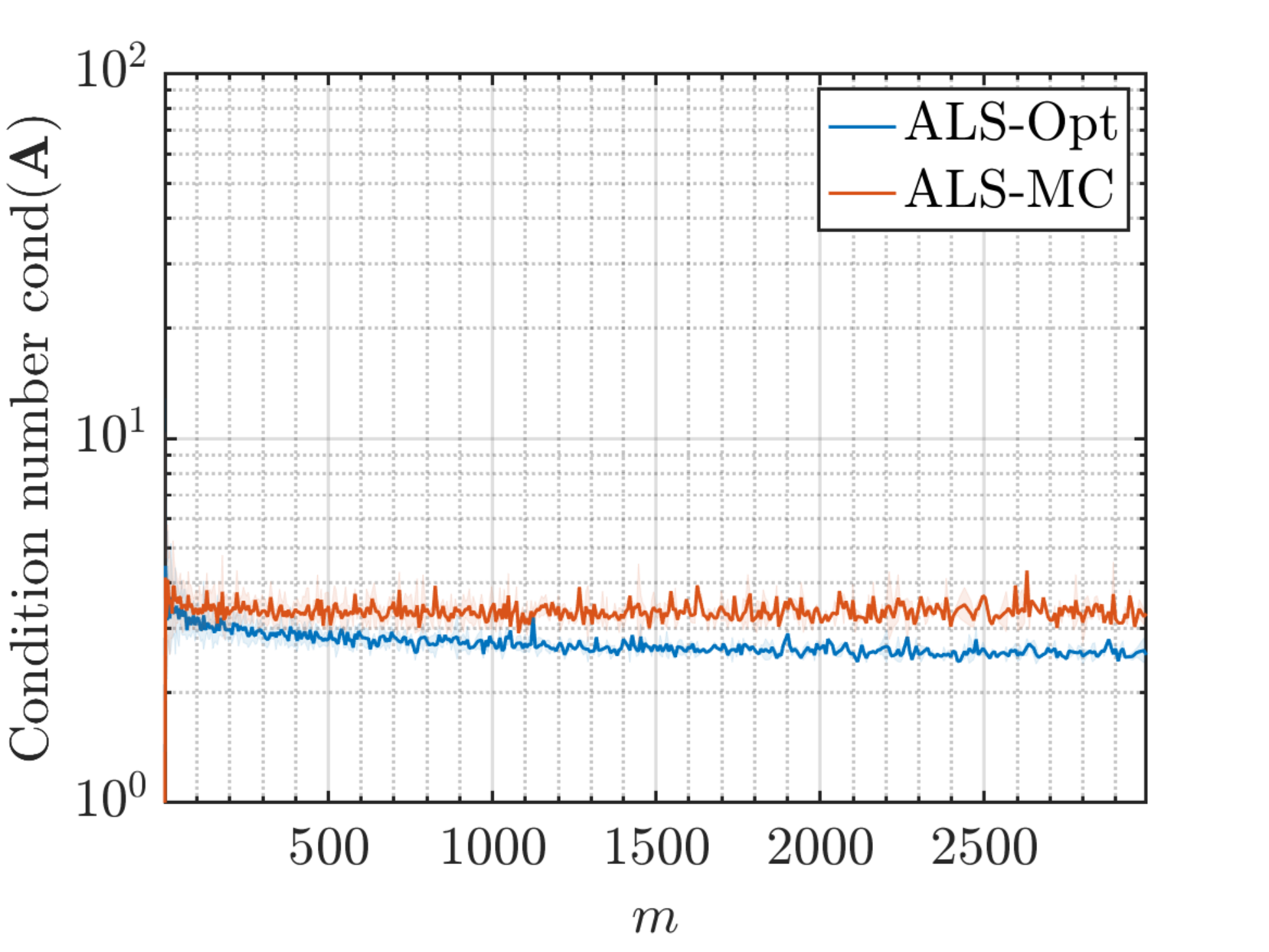}
&
\includegraphics[width = \errplotimg]{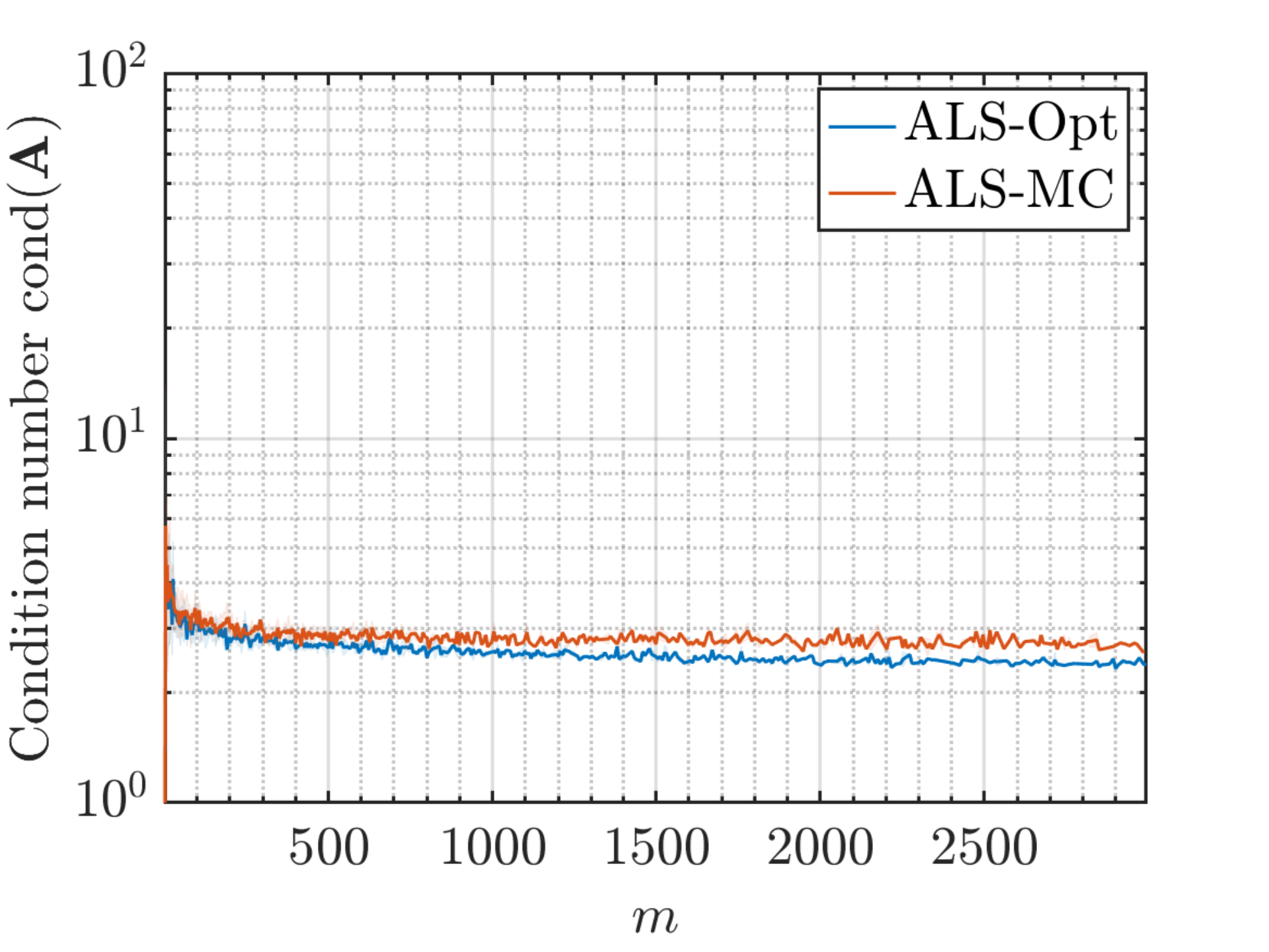}
\\[\errplotgraphsp]
\includegraphics[width = \errplotimg]{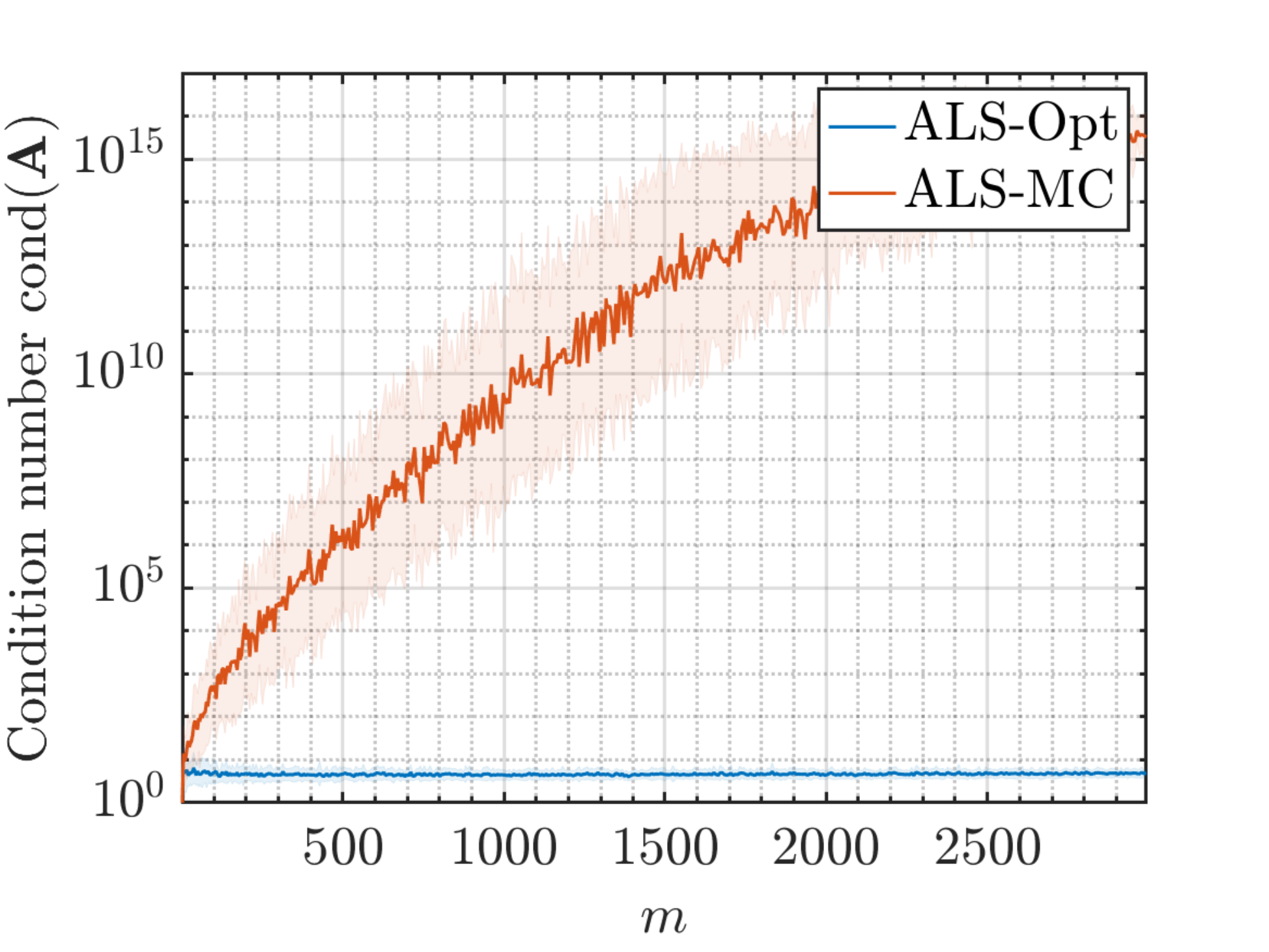}
&
\includegraphics[width = \errplotimg]{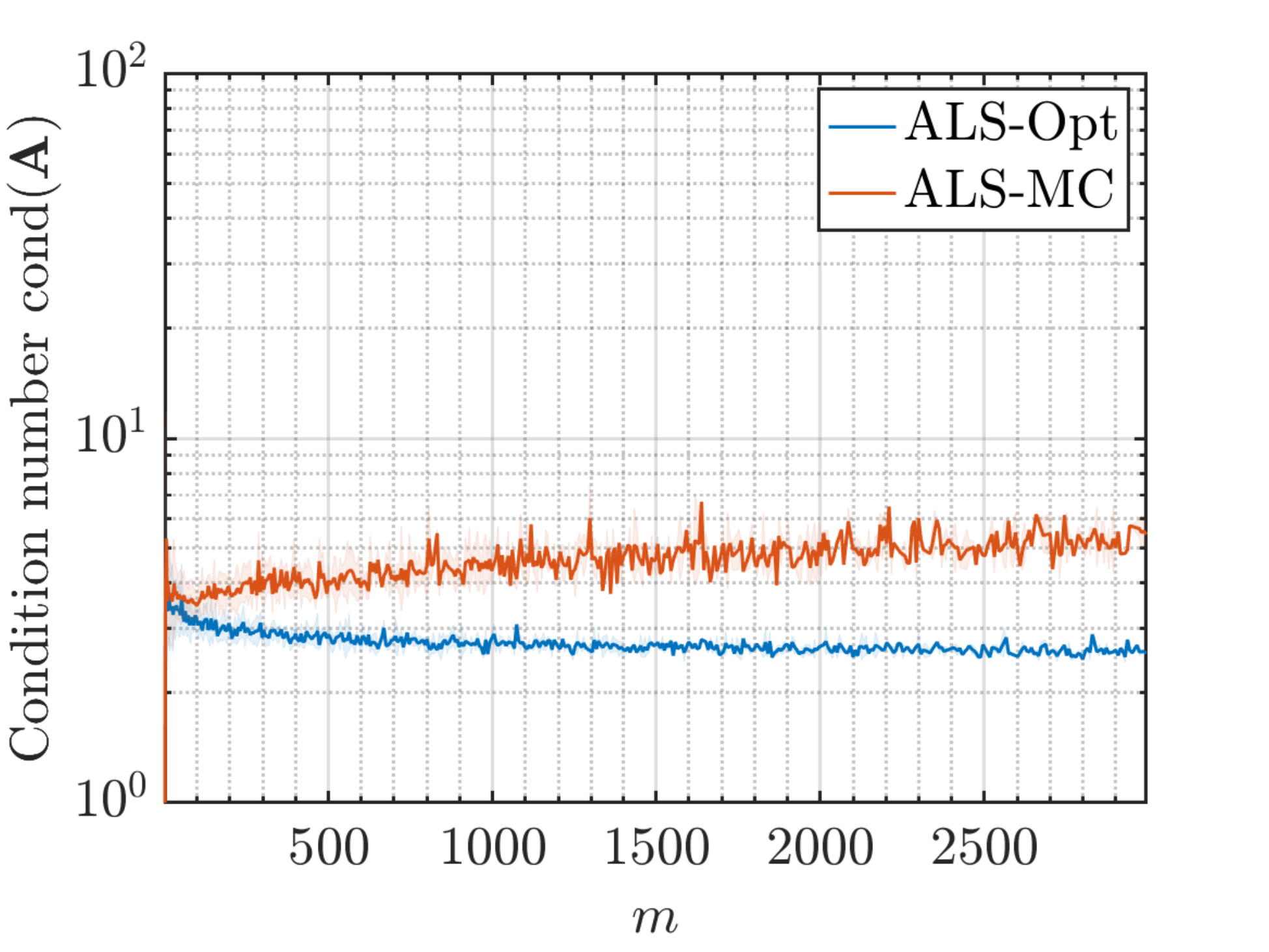}
&
\includegraphics[width = \errplotimg]{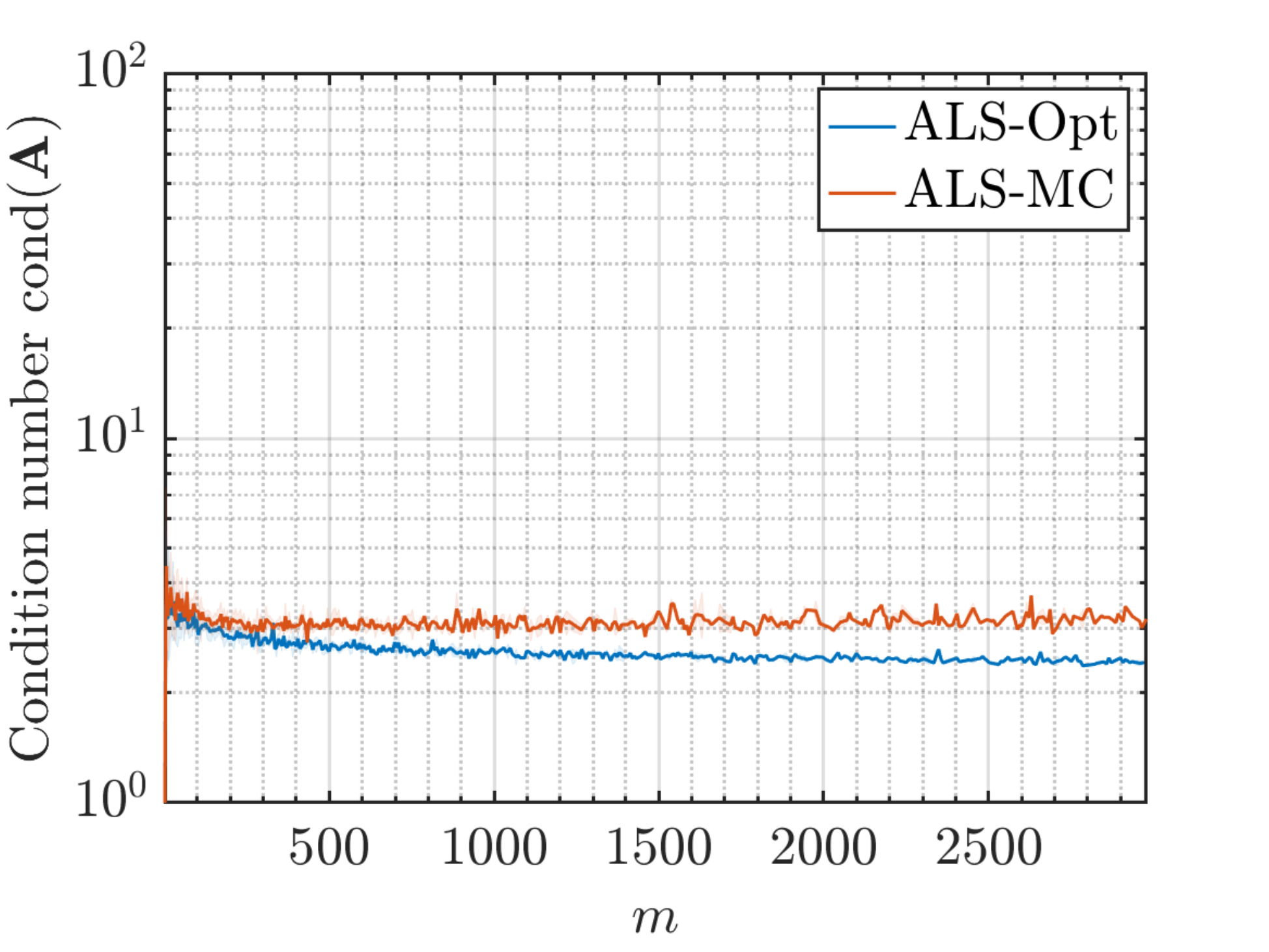}
\\[\errplotgraphsp]
\includegraphics[width = \errplotimg]{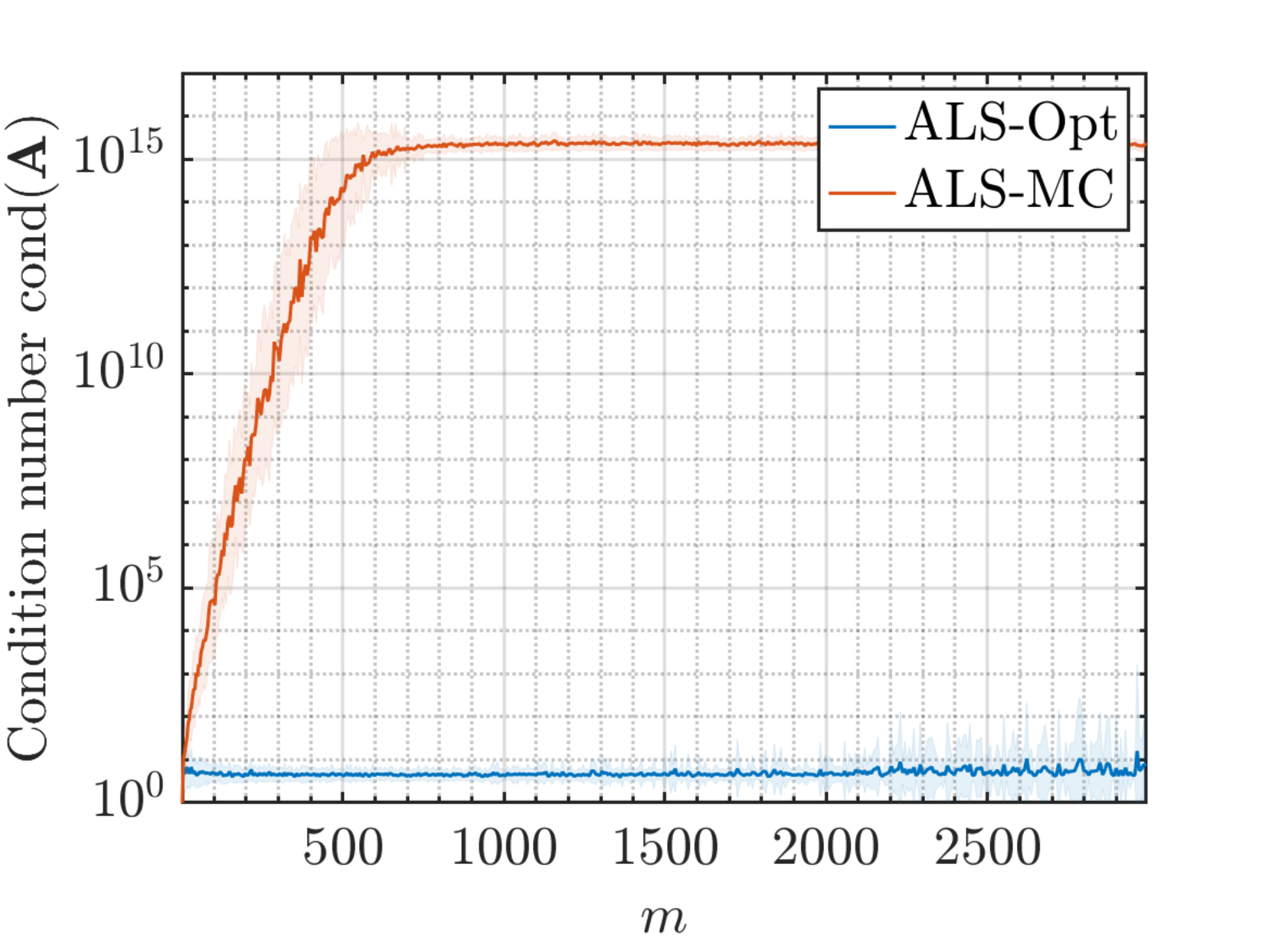}
&
\includegraphics[width = \errplotimg]{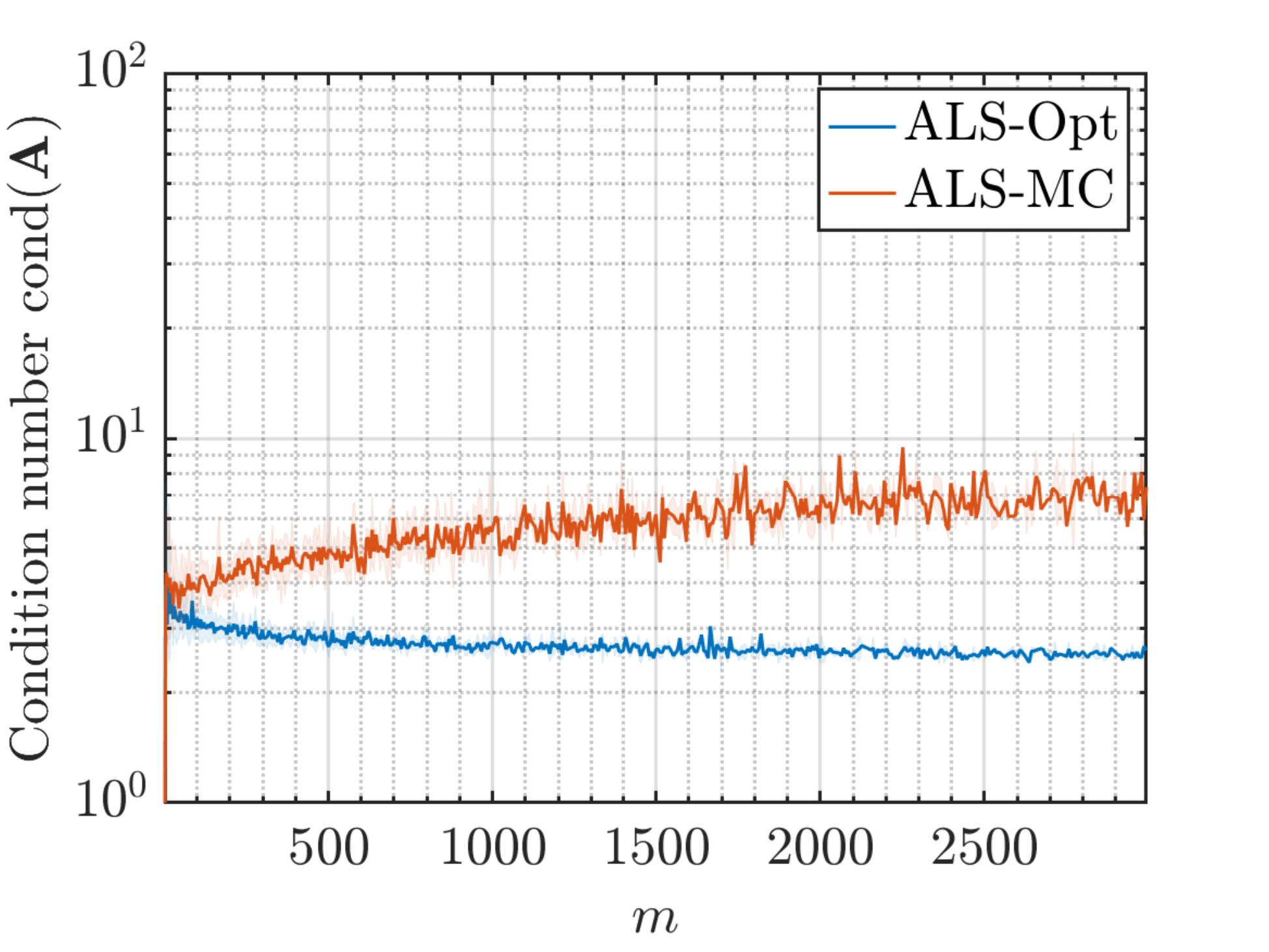}
&
\includegraphics[width = \errplotimg]{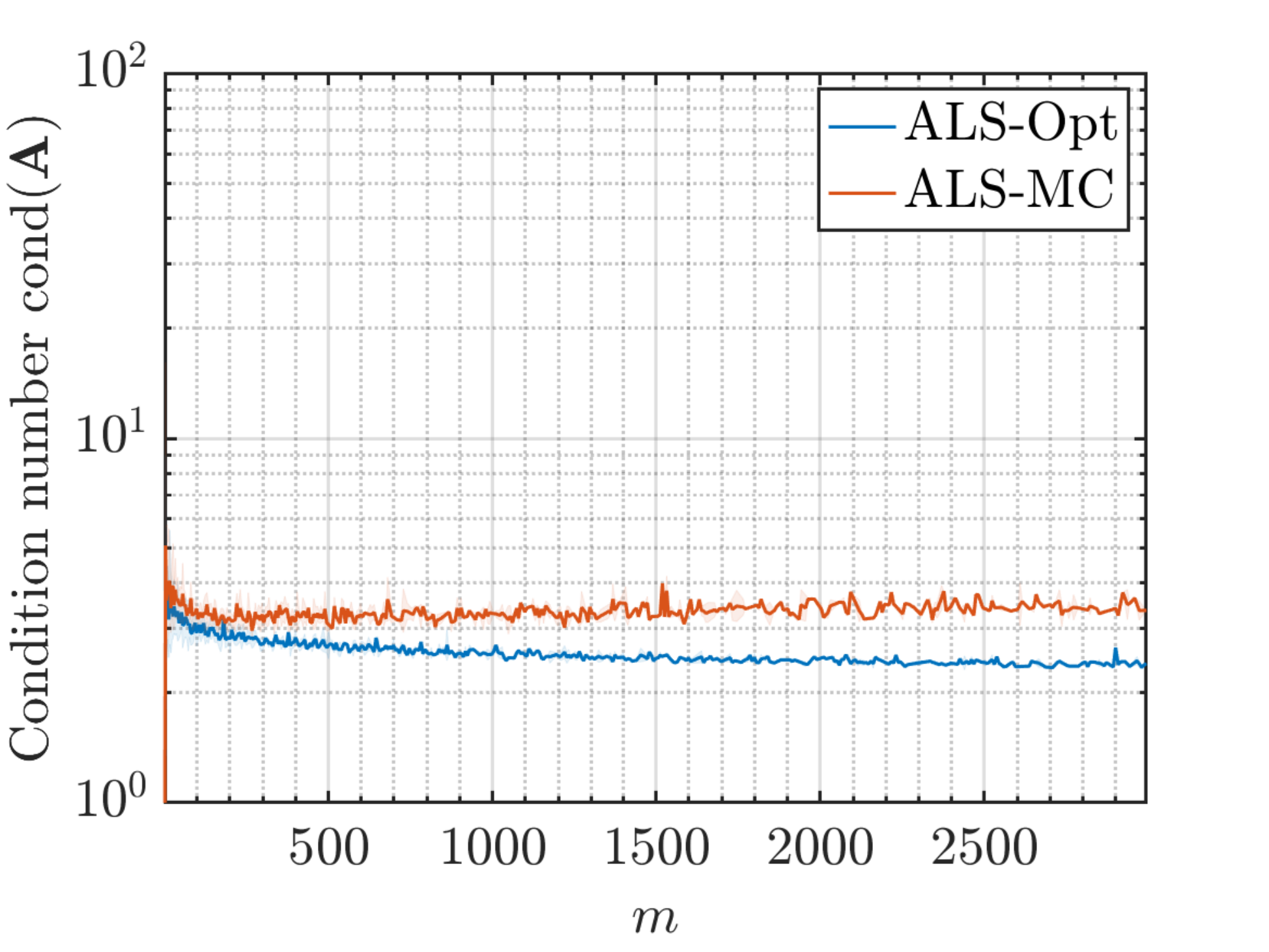}

\\[\errplottextsp]
$d = 1$ & $d = 8$ & $d = 16$
\end{tabular}
\end{small}
\end{center}
\caption{The corresponding condition numbers for the experiments in Fig.\ \ref{fig:figSM3}.} 
\label{fig:figSM4}
\end{figure}

\begin{figure}[t!]
\begin{center}
\begin{small}
 \begin{tabular}{@{\hspace{0pt}}c@{\hspace{\errplotsp}}c@{\hspace{\errplotsp}}c@{\hspace{0pt}}}
\includegraphics[width = \errplotimg]{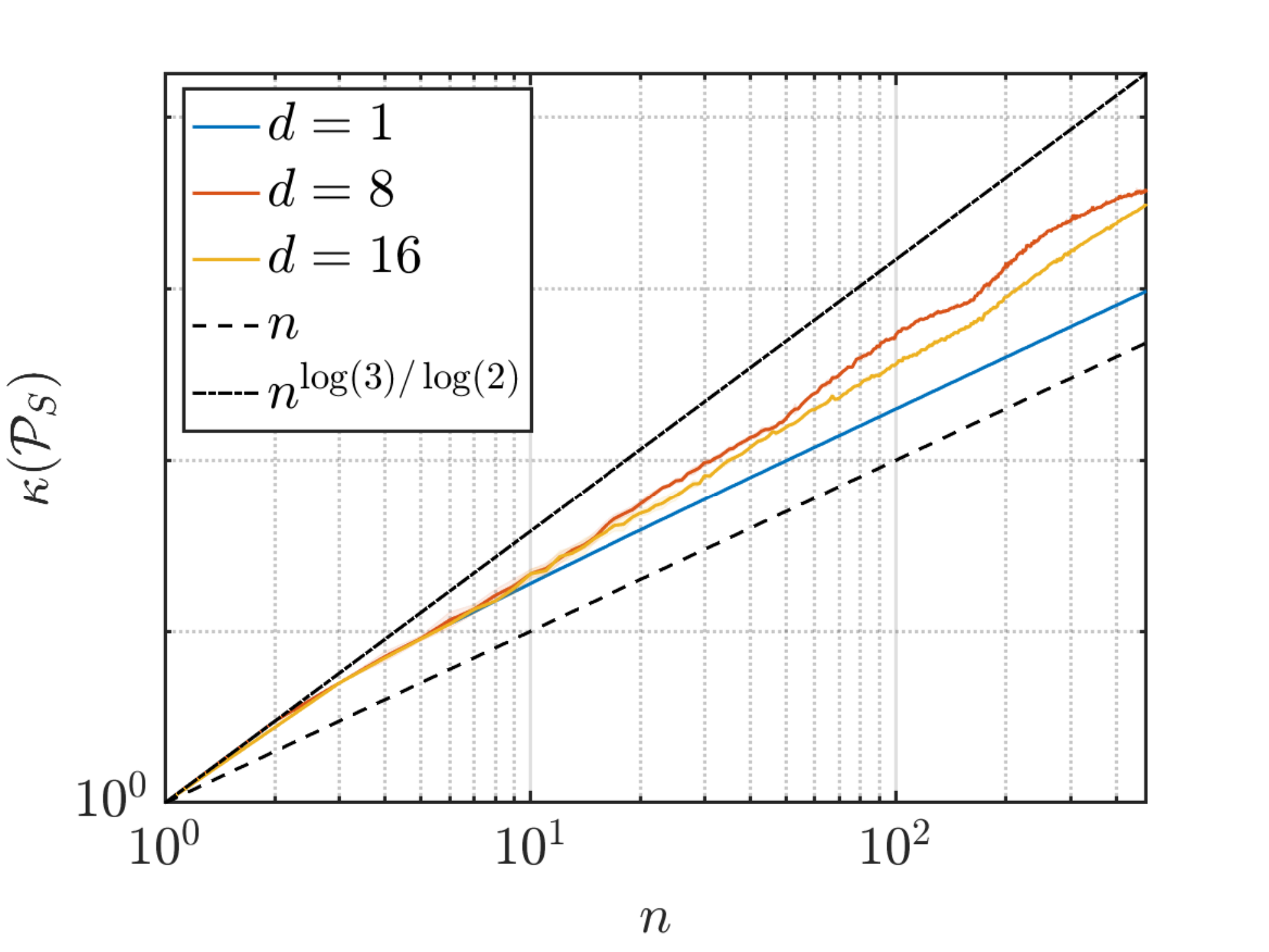}&
\includegraphics[width = \errplotimg]{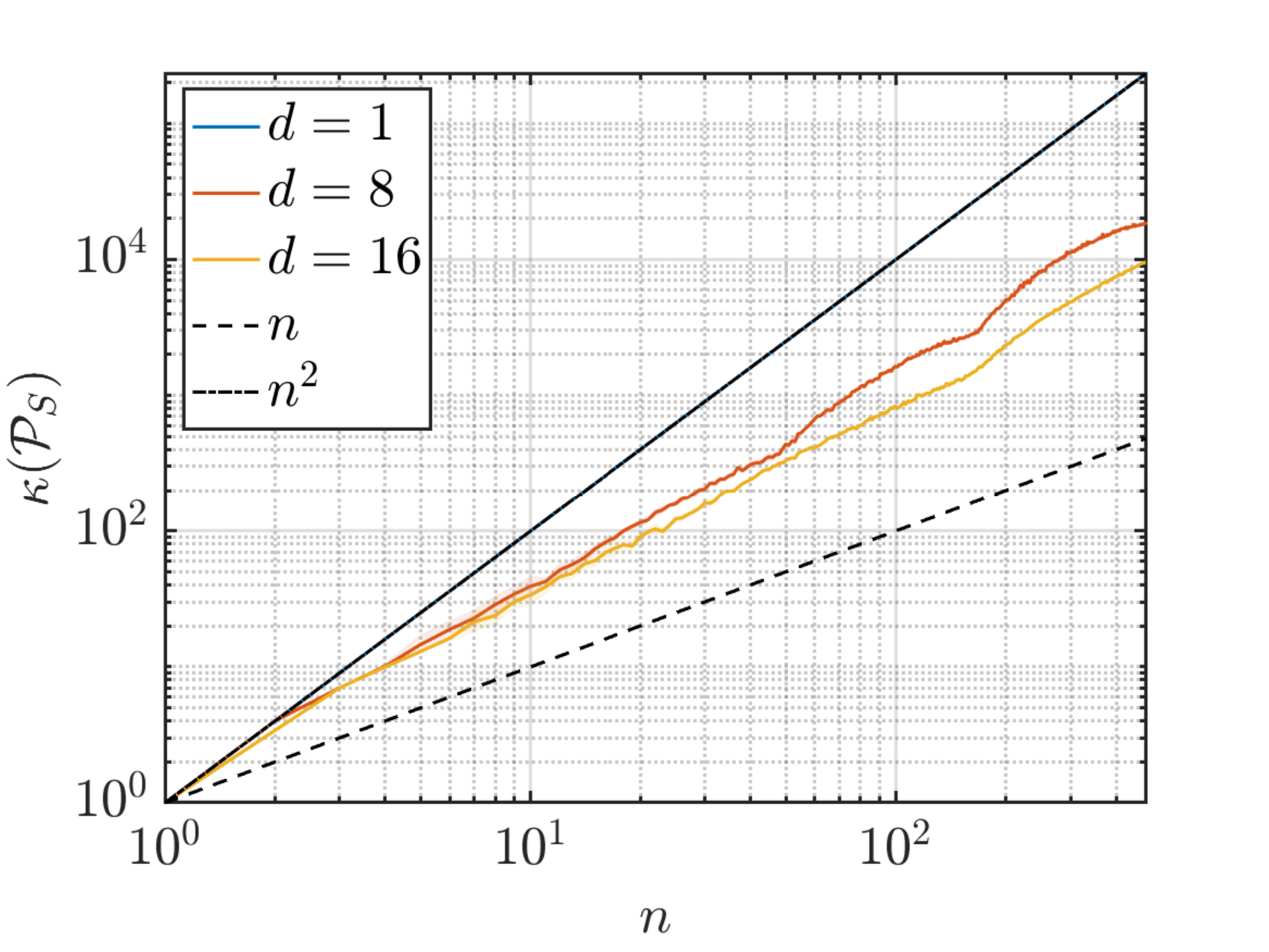}&
\includegraphics[width = \errplotimg]{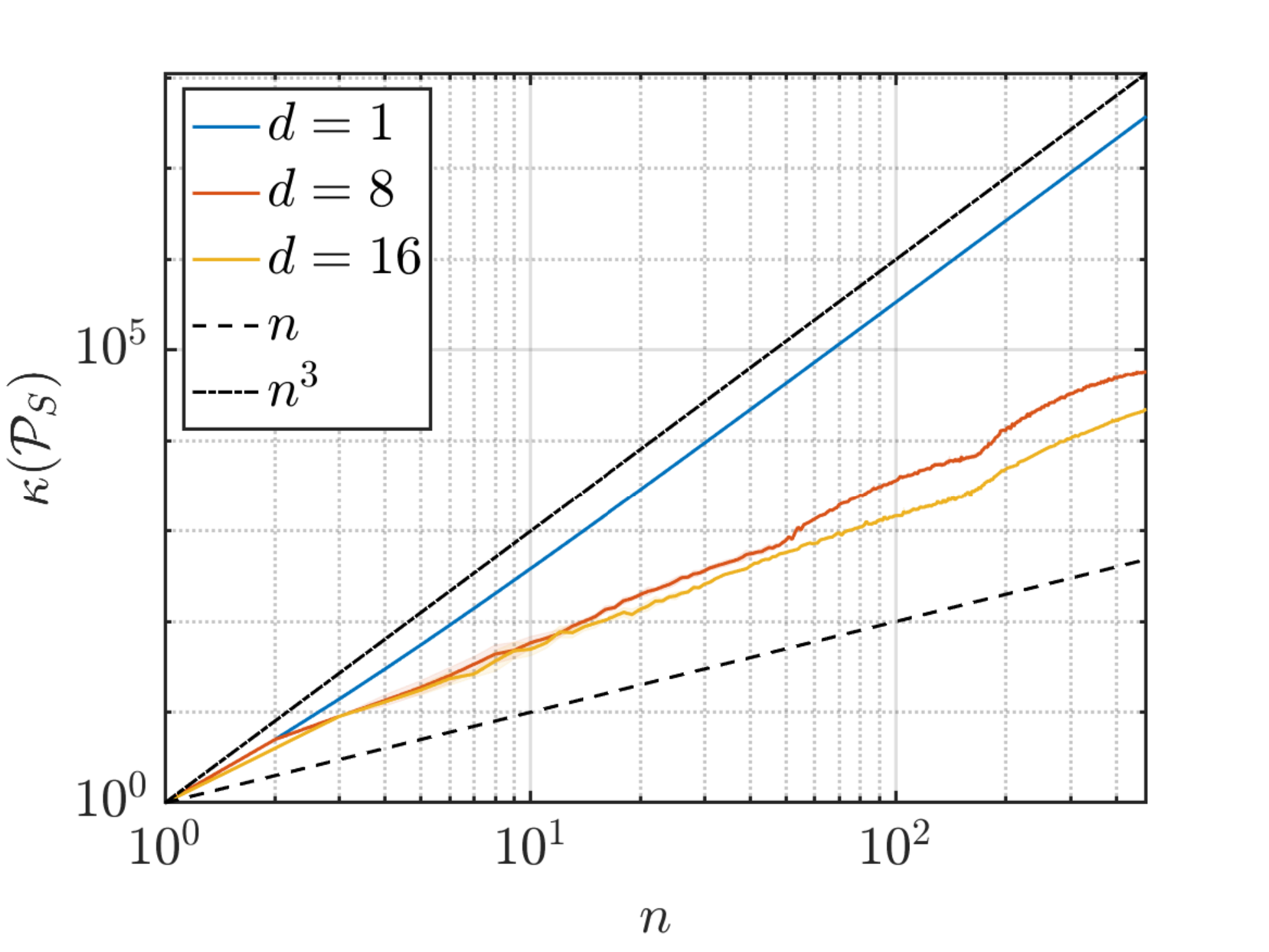}
\\[\errplotgraphsp]
\includegraphics[width = \errplotimg]{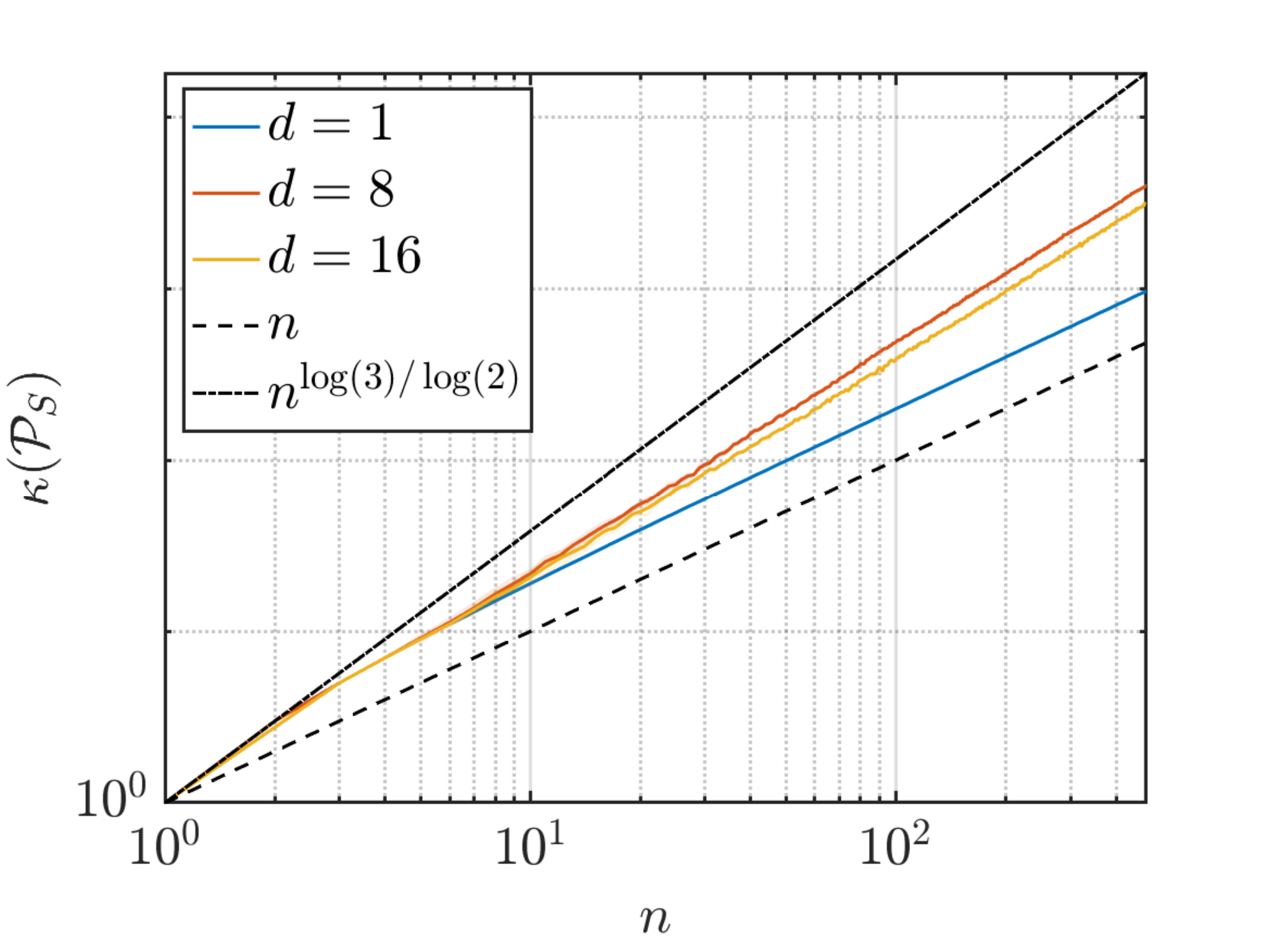}&
\includegraphics[width = \errplotimg]{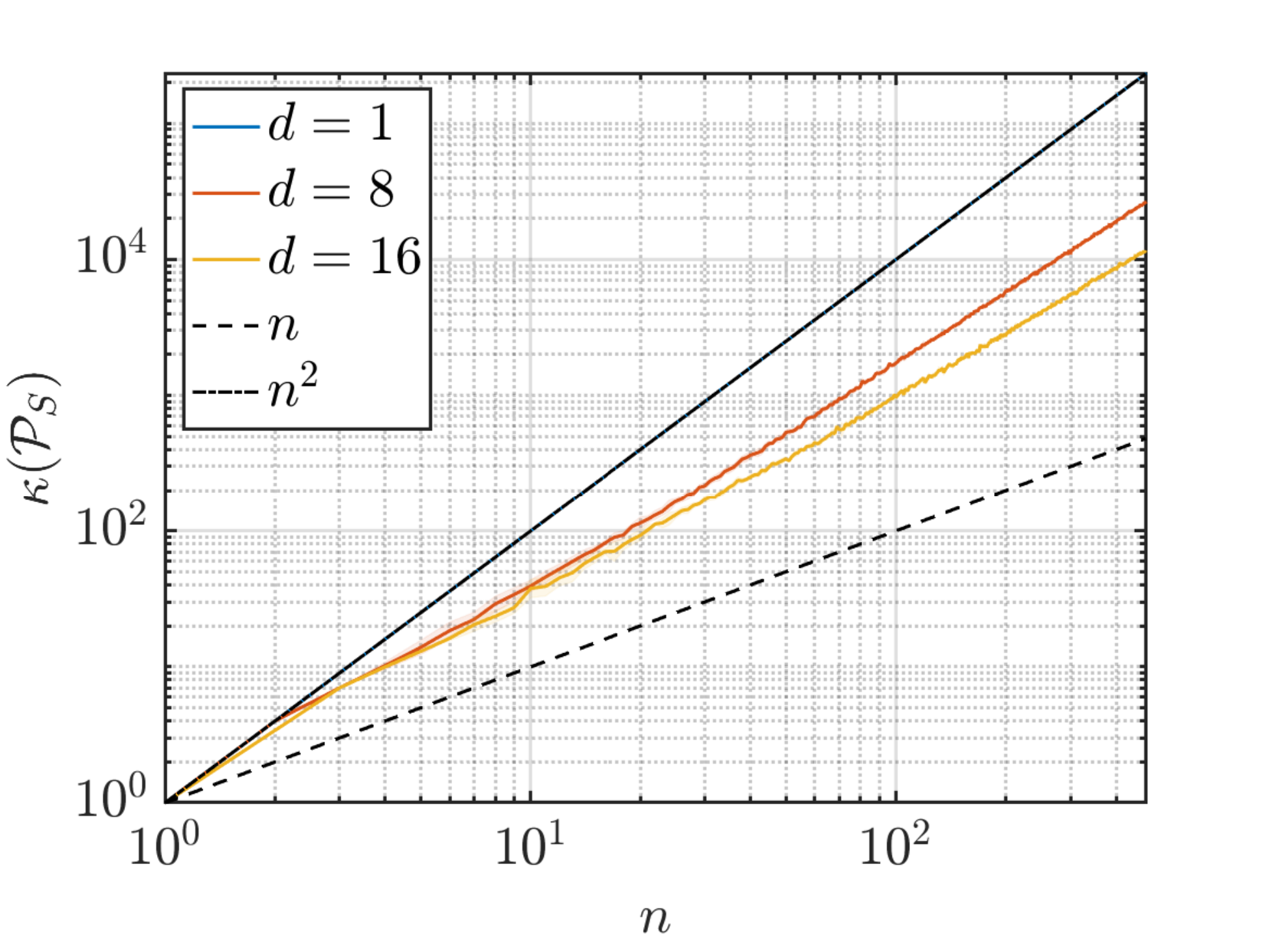}&
\includegraphics[width = \errplotimg]{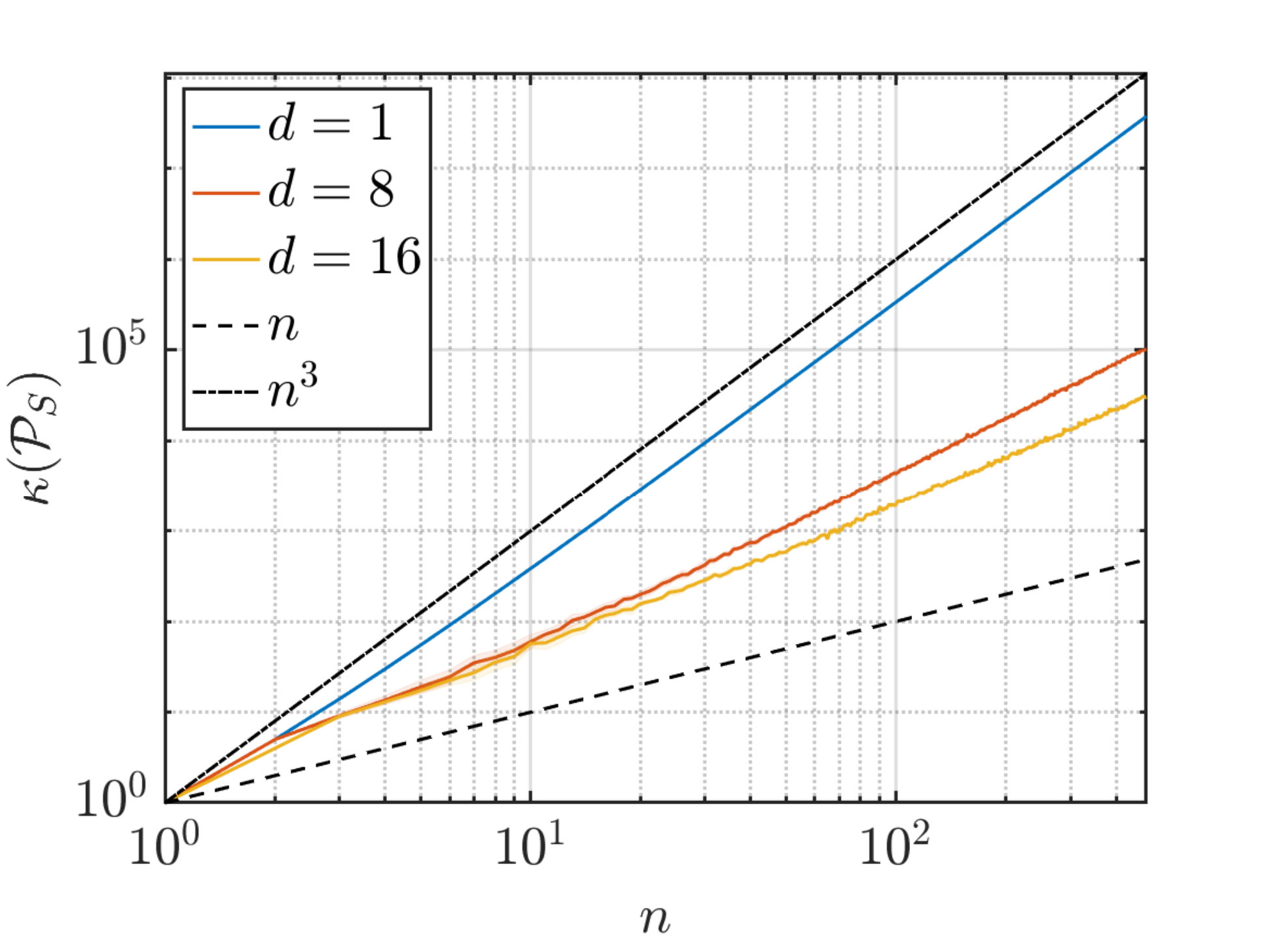}
\\[\errplottextsp]
Chebyshev & Legendre & Chebyshev second kind
\end{tabular}
\end{small}
\end{center}
\caption{The corresponding values of $\kappa(\cP_{S})$ for the experiments in Figs.\ \ref{fig:figSM1} (top row) and \ref{fig:figSM3} (bottom row).} 
\label{fig:figSM5}
\end{figure}

In the numerical experiments in the main paper we consider the uniform measure $\D \varrho(\bm{y}) = 2^{-d} \D \bm{y}$. In Figs.\ \ref{fig:figSM1}--\ref{fig:figSM5} we compare this against the first kind Chebyshev measure
\be{
\label{Cheby1st}
\D \varrho(\bm{y}) = \prod^{d}_{j=1} \frac{1}{\pi \sqrt{1-y^2_j}} \D y_j
}
and the second kind Chebyshev measure
\be{
\label{Cheby2nd}
\D \varrho(\bm{y}) = \prod^{d}_{j=1} \frac{2}{\pi} \sqrt{1-y^2_j} \D y_j.
}
The measure \eqref{Cheby1st} is covered by our main theoretical results in the paper. Hence we expect the the near-optimal sampling strategy to offer a diminishing improvement as the dimension increases over MC sampling from $\varrho$. This is clearly seen in these figures. However, we also see that the near-optimal sampling strategy offers little to no improvement in low dimensions as well. This is unsurprising. Recall from Proposition \ref{prop:chebyshev-Christoffel} that the relevant constant $\kappa(\cP_{S})$ only scales superlinearly in $n$, specifically, like $n^{\log(3)/\log(2)}$, in the high-dimenional regime  $n \leq 2^{d+1}$. In low dimensions, it scales at worst like $2^d n$. This is significantly different to the case of the uniform measure, for which the worst case scaling is quadratic in $n$, in any dimension.

While second kind Chebyshev measure \eqref{Cheby2nd} is not covered by our main theorems, Figs.\ \ref{fig:figSM1}--\ref{fig:figSM5} suggest a similar phenomenon occurs in this case. We expect our results can be extended to this and more general \textit{Jacobi} measures. It is interesting to note that MC sampling from \eqref{Cheby2nd} leads to much worse ill-conditioning in low dimensions than in the case of random sampling from the uniform measure. Indeed, in Figs.\ \ref{fig:figSM2} and \ref{fig:figSM4}, we see that the condition number is around $10^{10}$ when $m \approx 300$, whereas in the case of the uniform measure this value is only reached when $m \approx 1000$.

Insight into this behaviour can be gained by looking at the worst-case behaviour of the constant $\kappa(\cP_S)$. The univariate, orthonormal second kind Chebyshev polynomials are given by $\psi_{\nu}(y) = \sin((\nu+1) \theta) / \sin(\theta)$, $y = \cos(\theta)$, and achieve their maximum absolute value at $y = \pm 1$, taking value $| \psi_{\nu}(y) | = n+1$ there. Hence, due to \eqref{K-alternative-def},
\bes{
\kappa(\cP_S) = \sum_{\bm{\nu} \in S} \prod^{d}_{j=1} (\nu_j+1)^2.
}
One can now use this to show a similar version of Proposition \ref{prop:legendre-Christoffel} for this measure, except with \eqref{eq:K_LU} replaced by the cubic scaling
\bes{
\frac13 n^3 \leq \max \left \{ \kappa(\cP_{S}) : S \subset \bbN^d_0,\ | S | \leq n,\ \textnormal{$S$ lower} \right \} \leq n^3.
}
The upper bound follows from \cite{migliorati2015multivariate} (which also considers more general Jacobi measures). The lower bound follows by setting $S = \{ k \bm{e}_1 : k = 0,\ldots,n-1 \}$. The fact the worst-case scaling is cubic, instead of quadratic as in \eqref{eq:K_LU}, explains why log-linear oversampling in Figs.\ \ref{fig:figSM1}--\ref{fig:figSM5} is much more poorly conditioned.

\subsection{Comparing CS to ALS with different scalings}\label{ss:compare-scalings}

In Figs.\ \ref{fig:fig10}--\ref{fig:fig14} we compared CS against ALS with the log-linear oversampling \eqref{m-scaling-LS}. The reader may be wondering if this represents a fair comparison, since the scaling between $m$ and $n$ in ALS is a parameter can potentially be tuned to improve its performance. In this section we present several examples comparing CS to ALS with several linear scalings. 

\begin{figure}[t!]
\begin{center}
\begin{small}
 \begin{tabular}{@{\hspace{0pt}}c@{\hspace{\errplotsp}}c@{\hspace{\errplotsp}}c@{\hspace{0pt}}}
\includegraphics[width = \errplotimg]{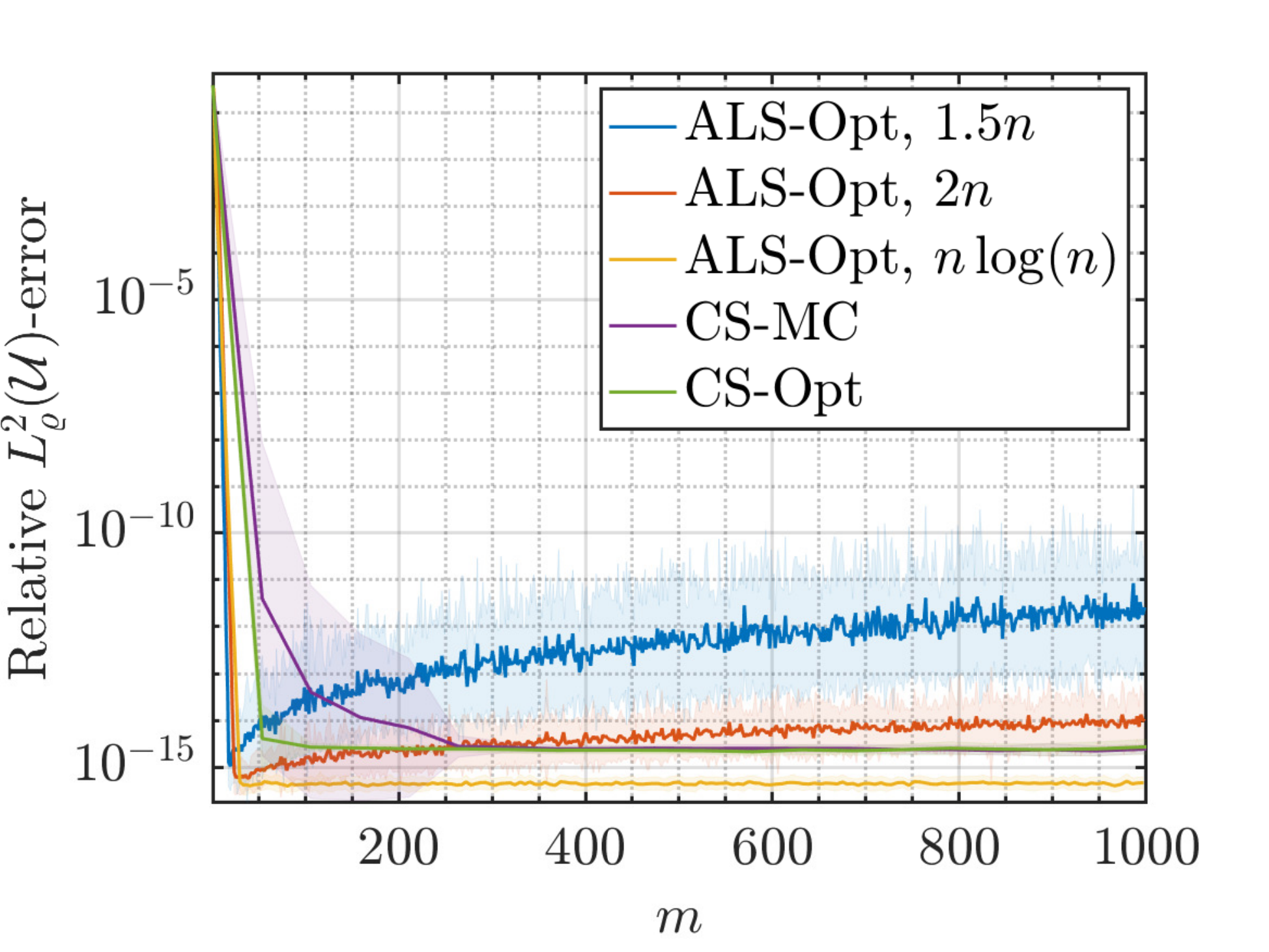}&
\includegraphics[width = \errplotimg]{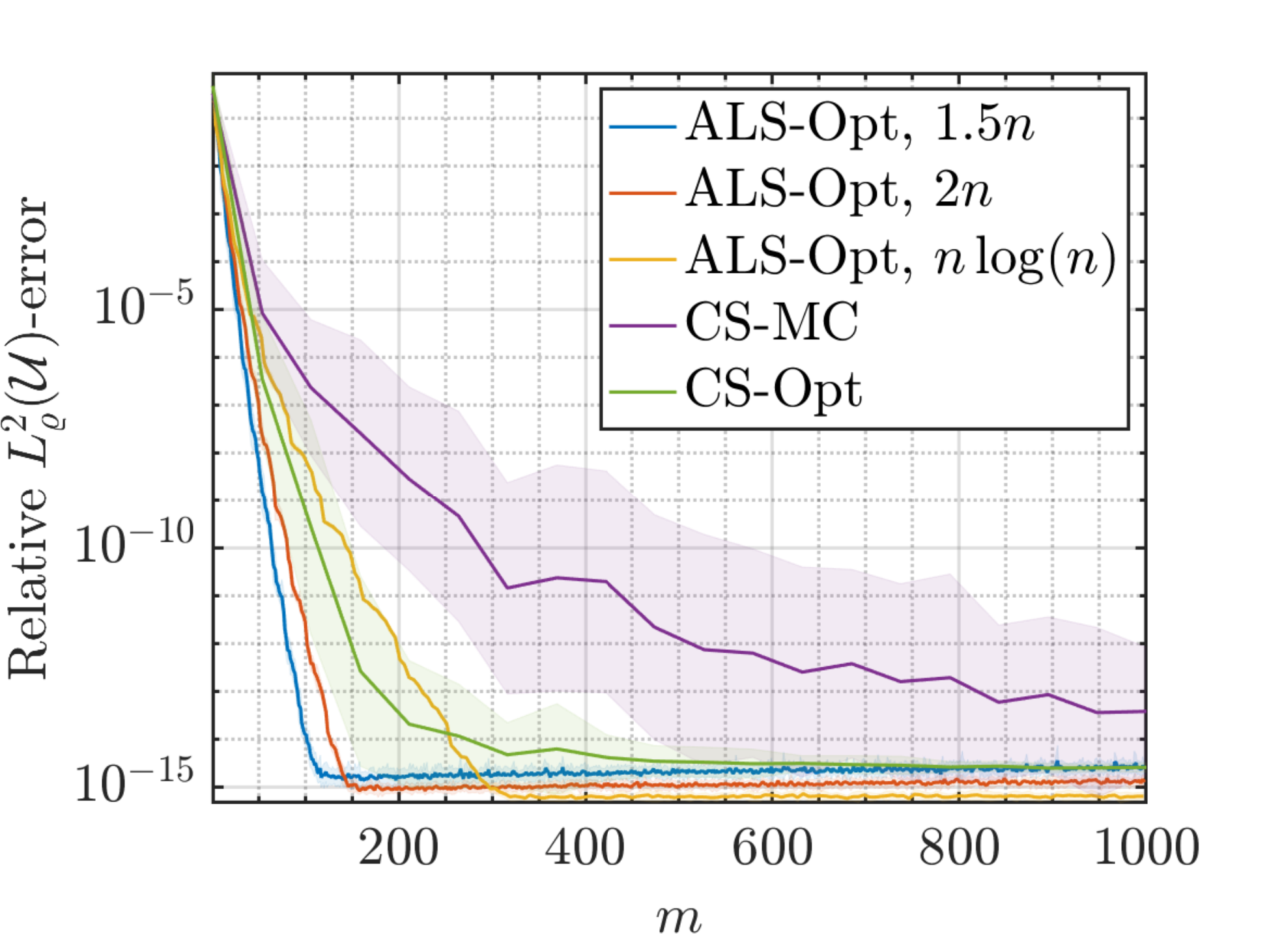}&
\includegraphics[width = \errplotimg]{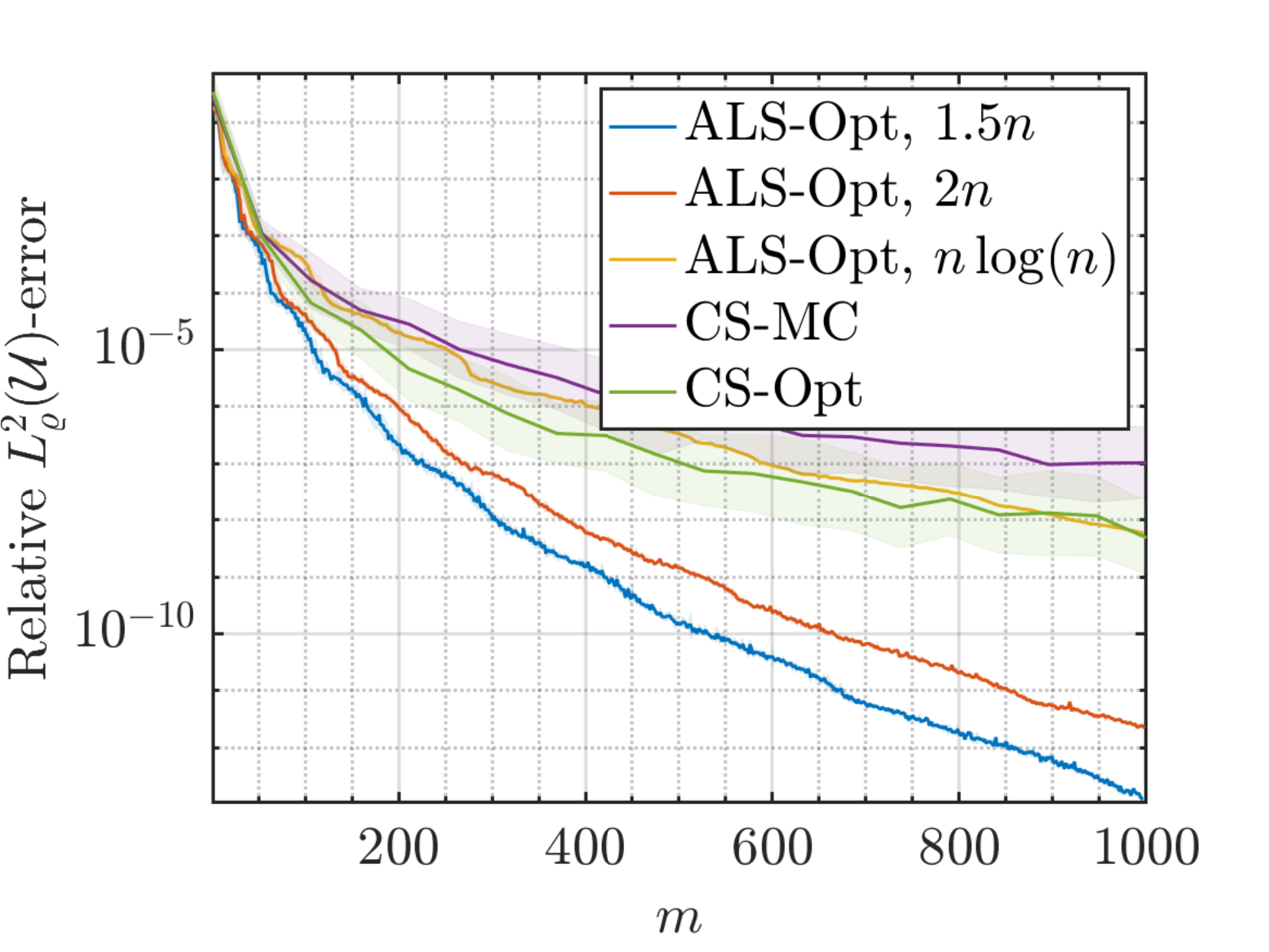}
\\[\errplotgraphsp]
\includegraphics[width = \errplotimg]{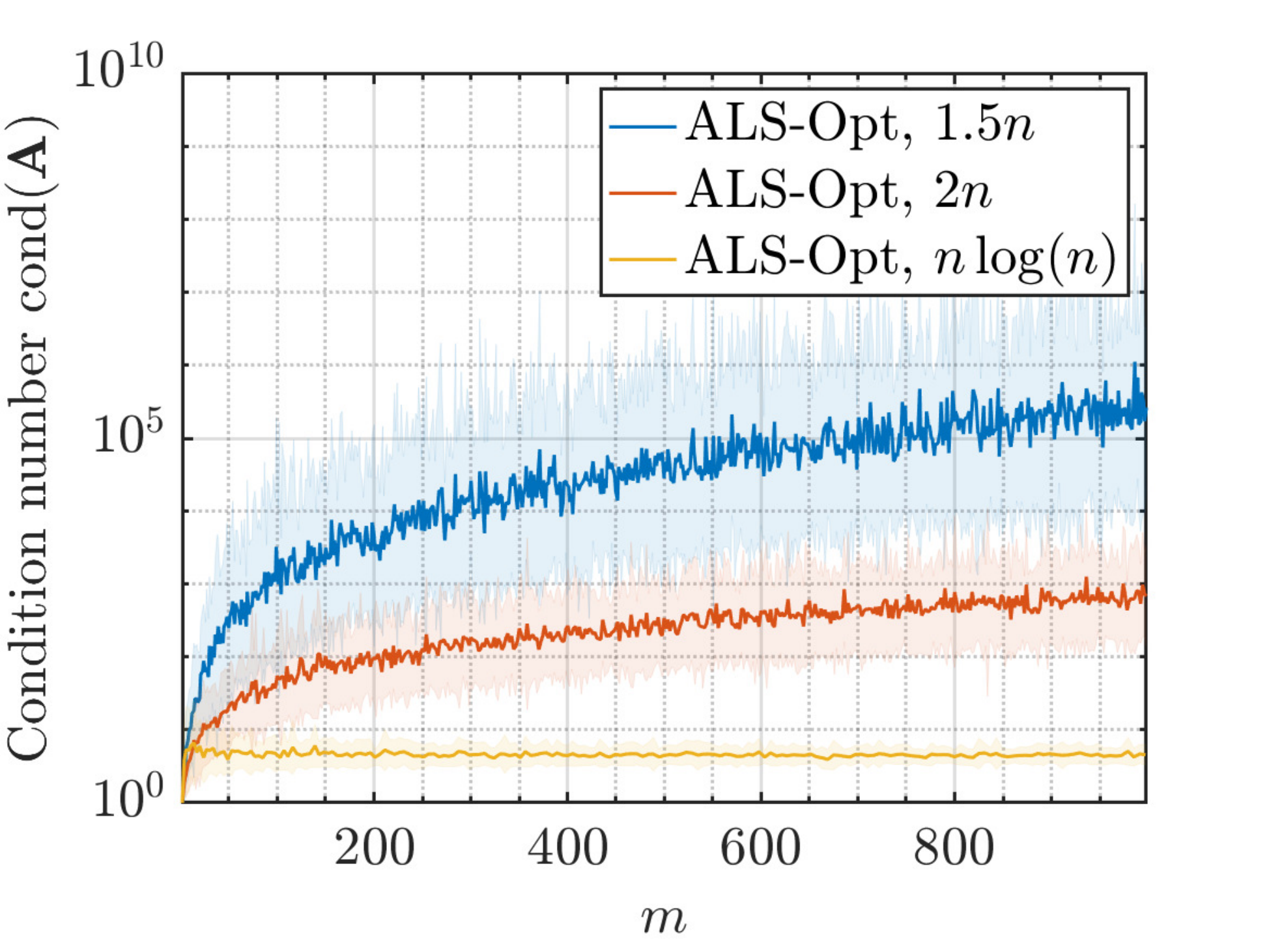}&
\includegraphics[width = \errplotimg]{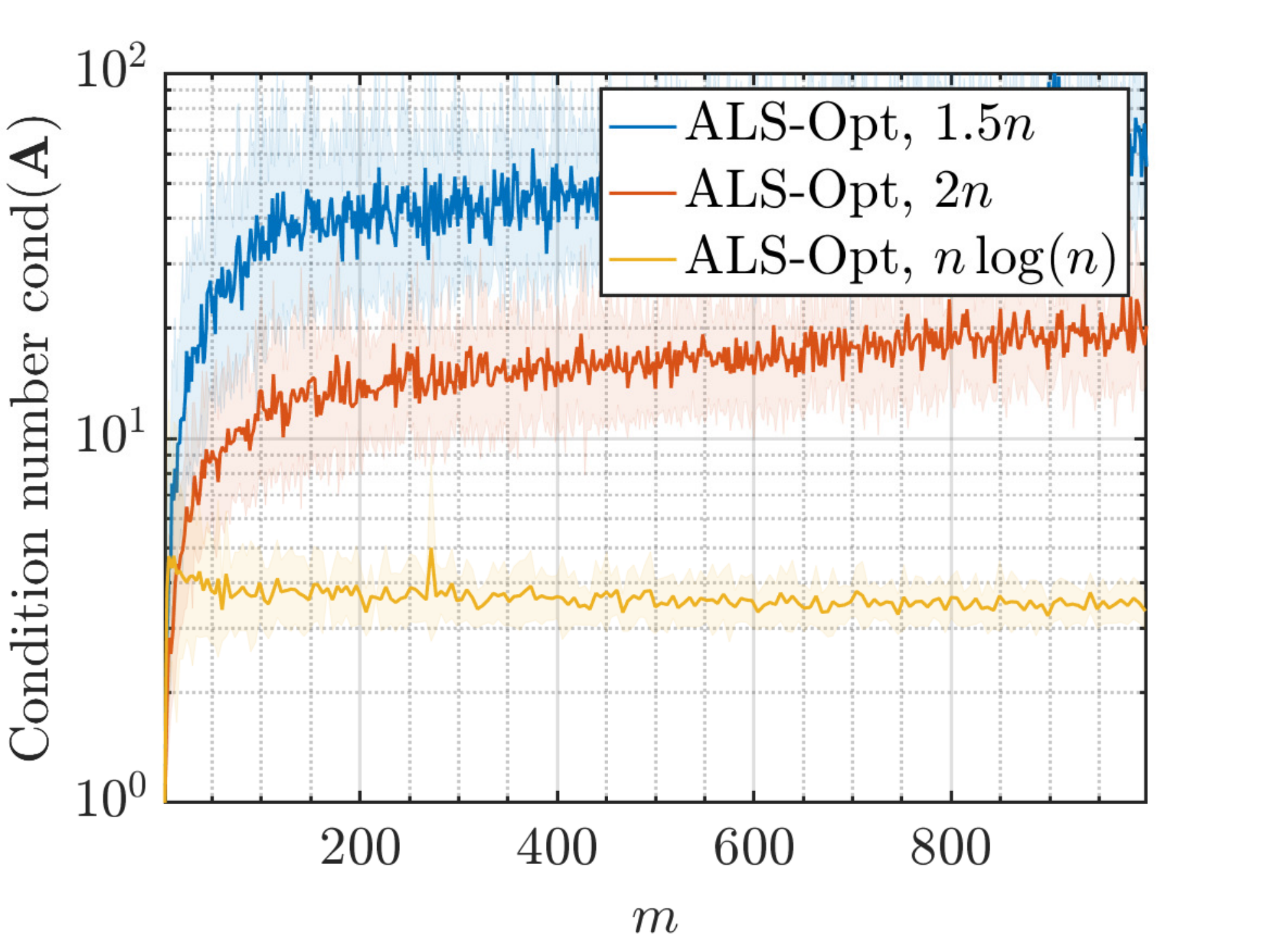}&
\includegraphics[width = \errplotimg]{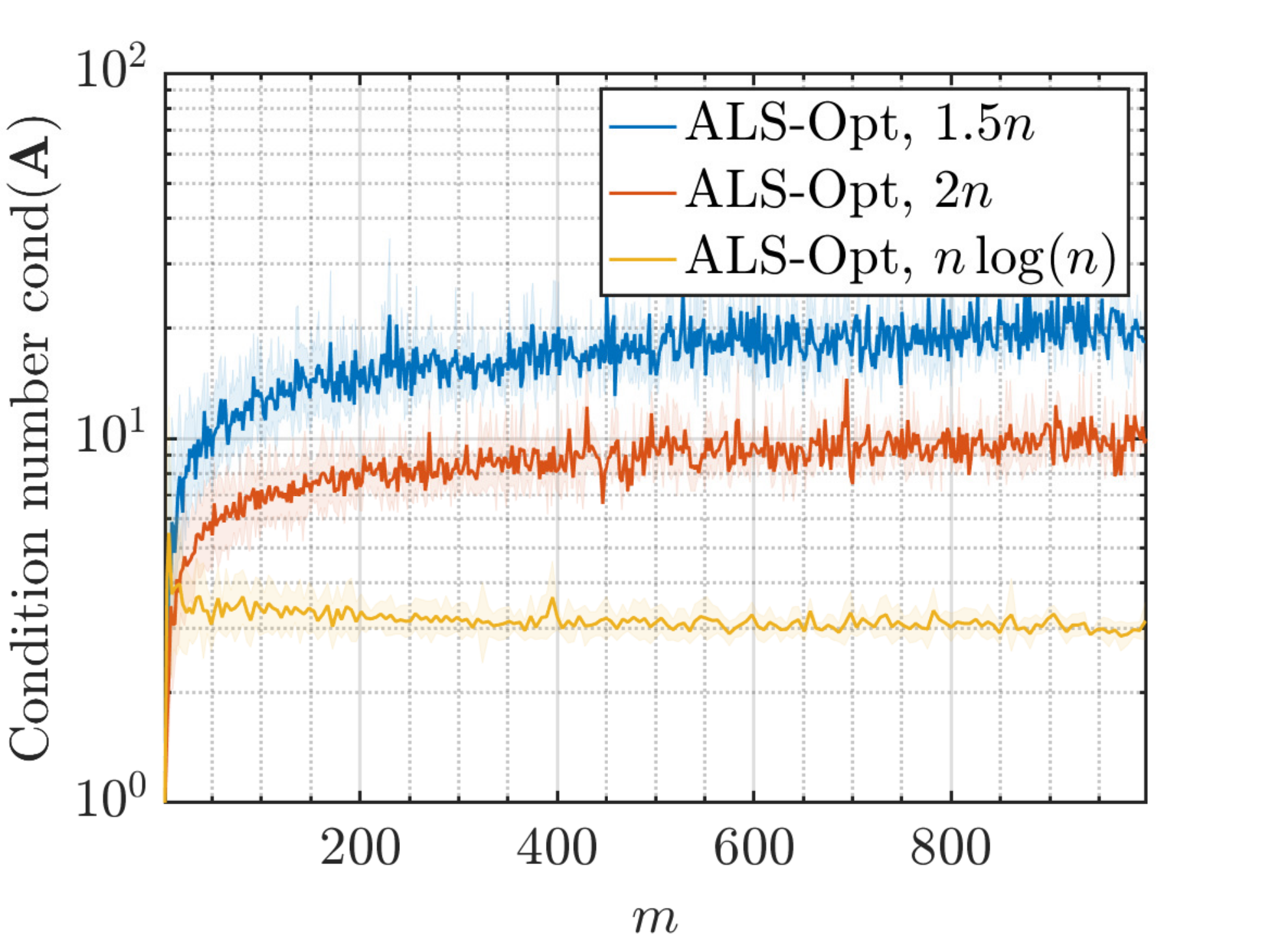}
\\[\errplottextsp]
$d = 1$ & $d = 2$ & $d = 4$
\\[\errplottextsp]
\includegraphics[width = \errplotimg]{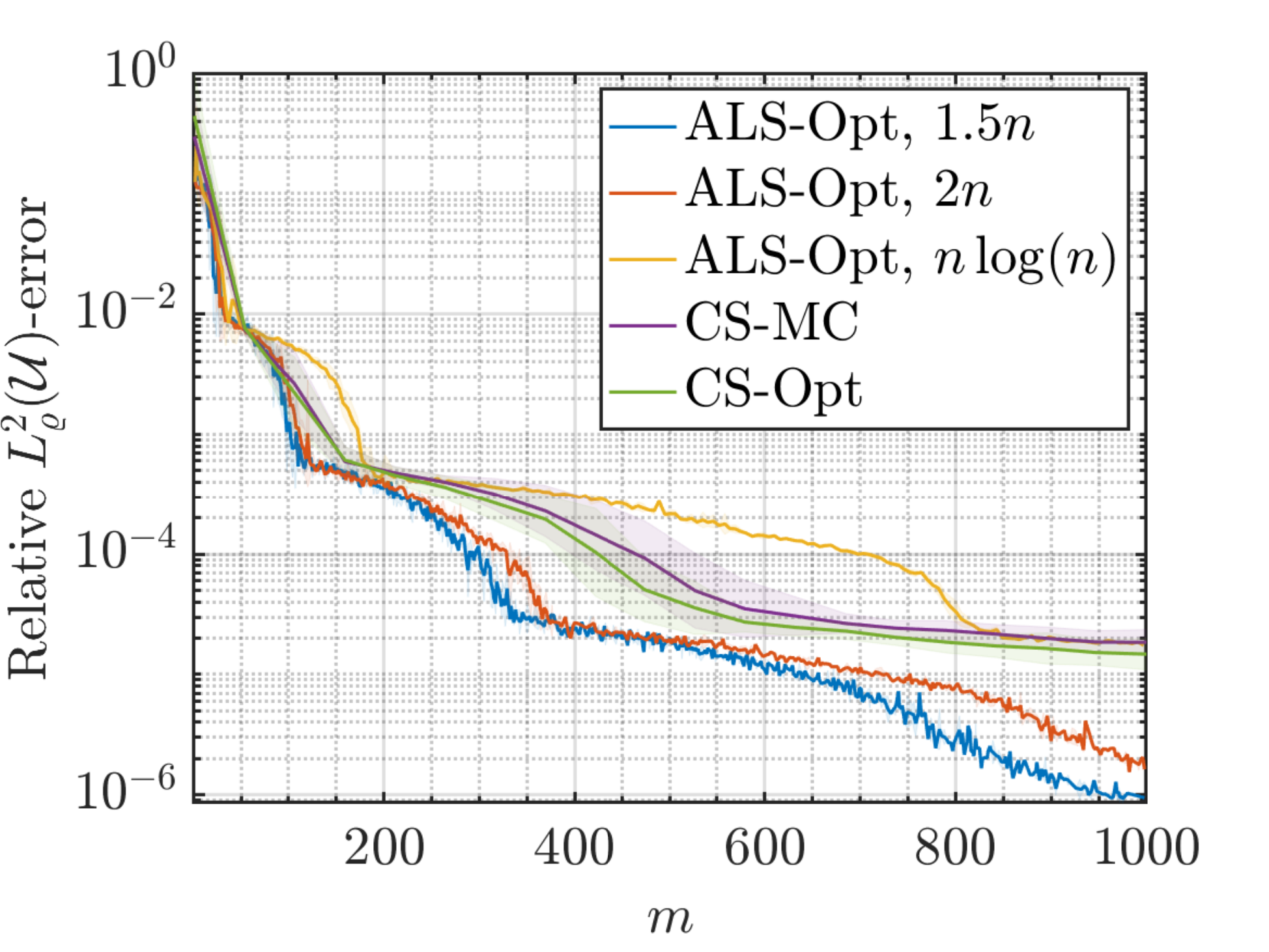}&
\includegraphics[width = \errplotimg]{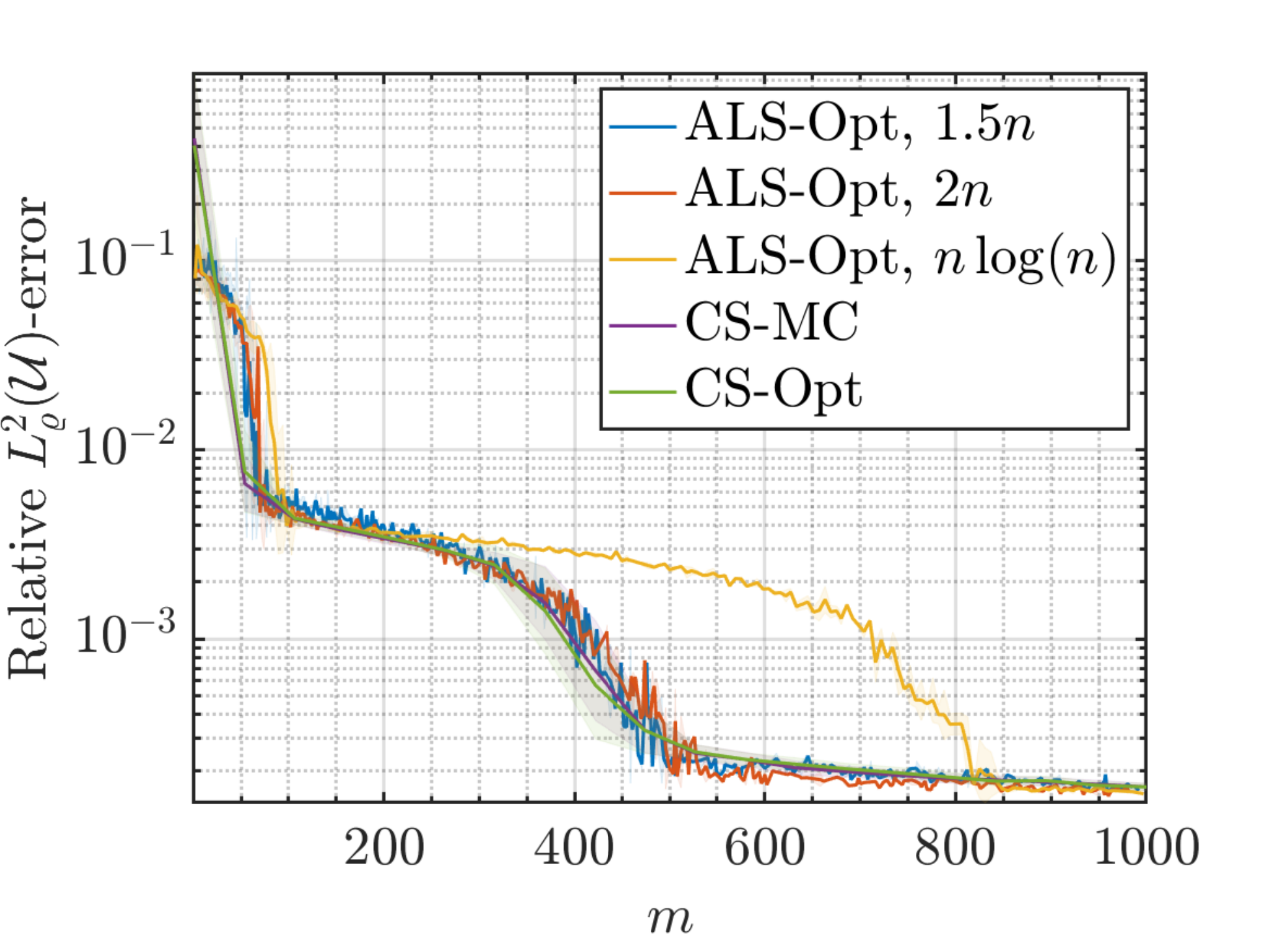}&
\includegraphics[width = \errplotimg]{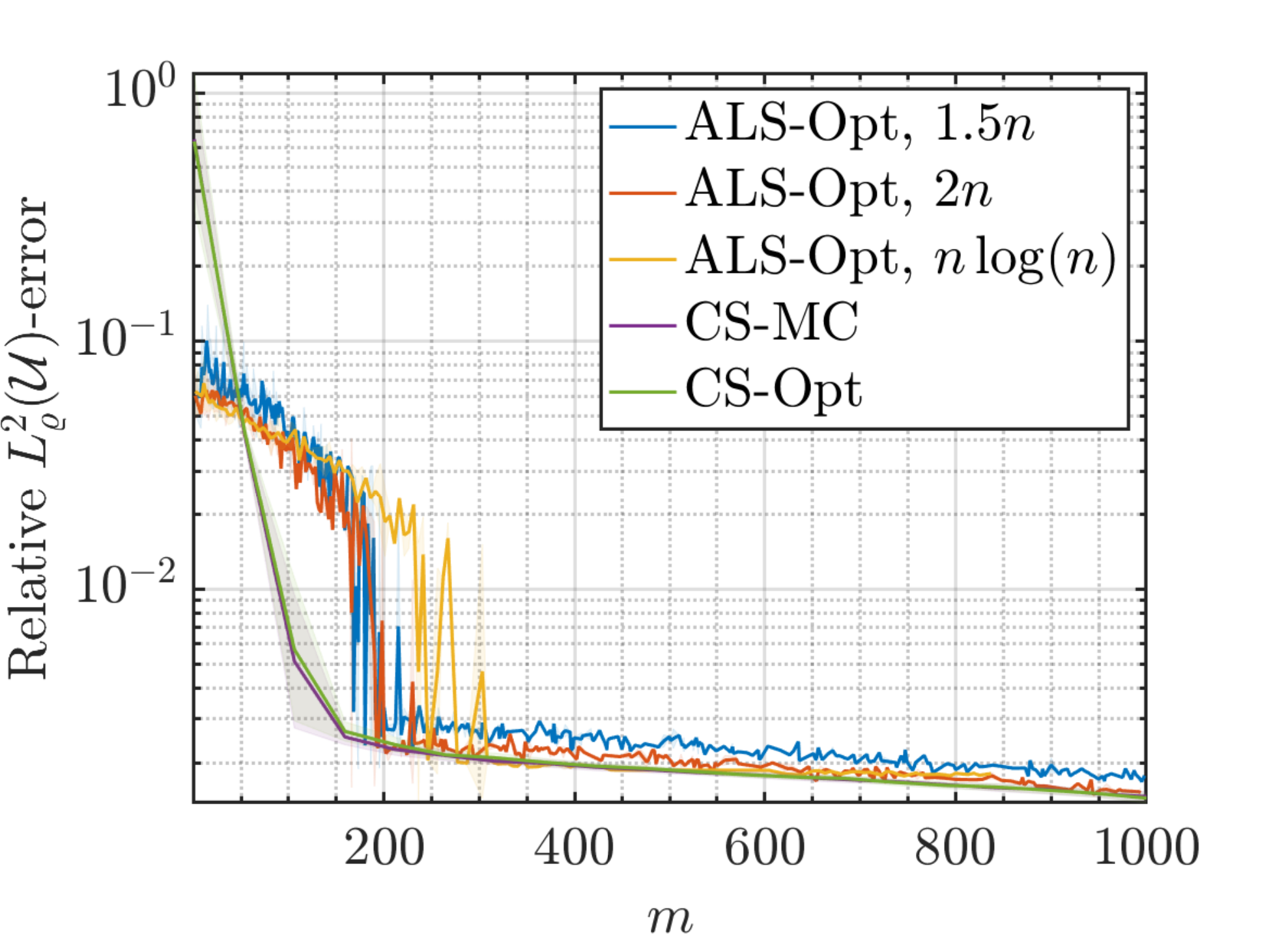}
\\[\errplotgraphsp]
\includegraphics[width = \errplotimg]{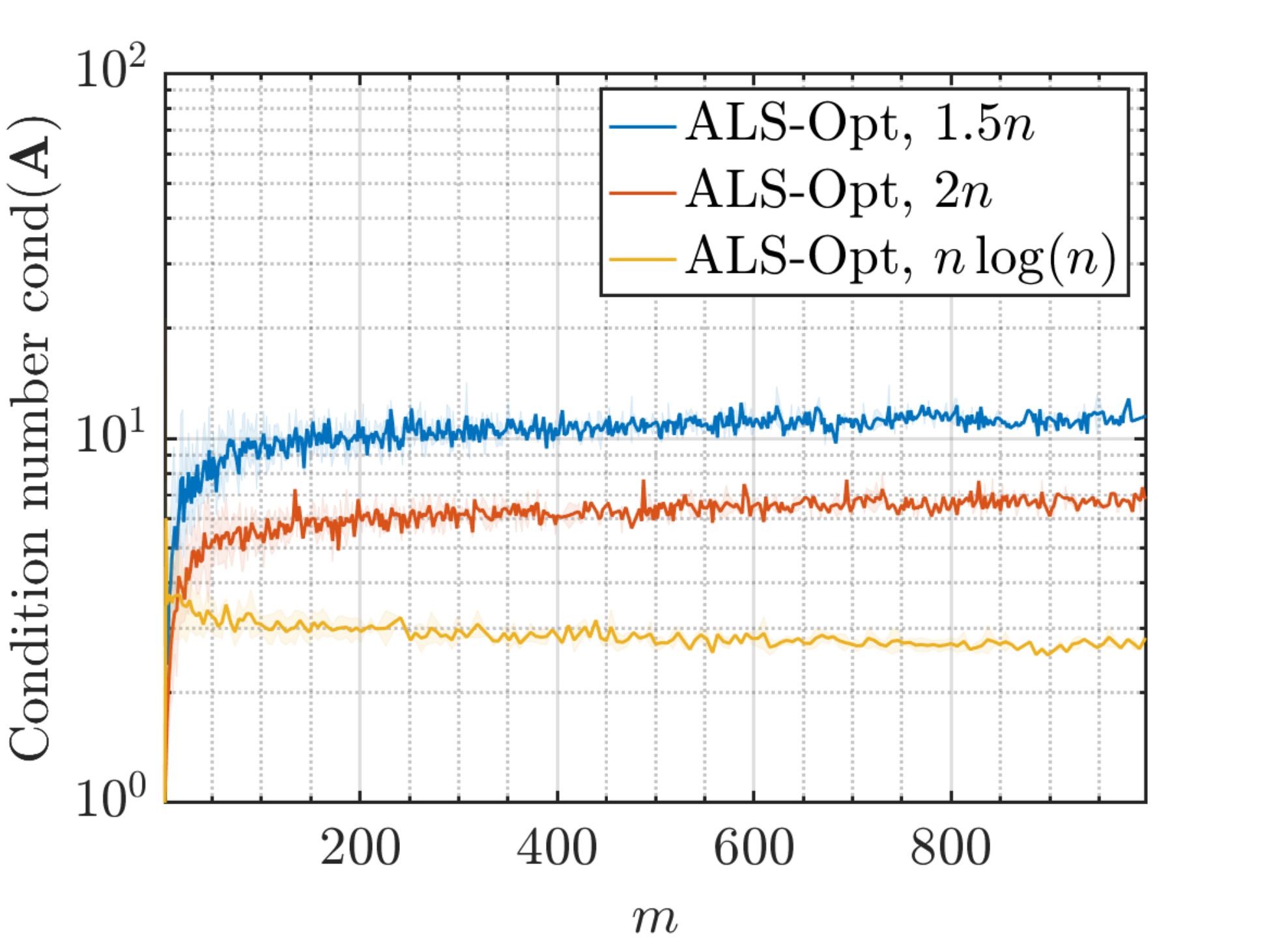}&
\includegraphics[width = \errplotimg]{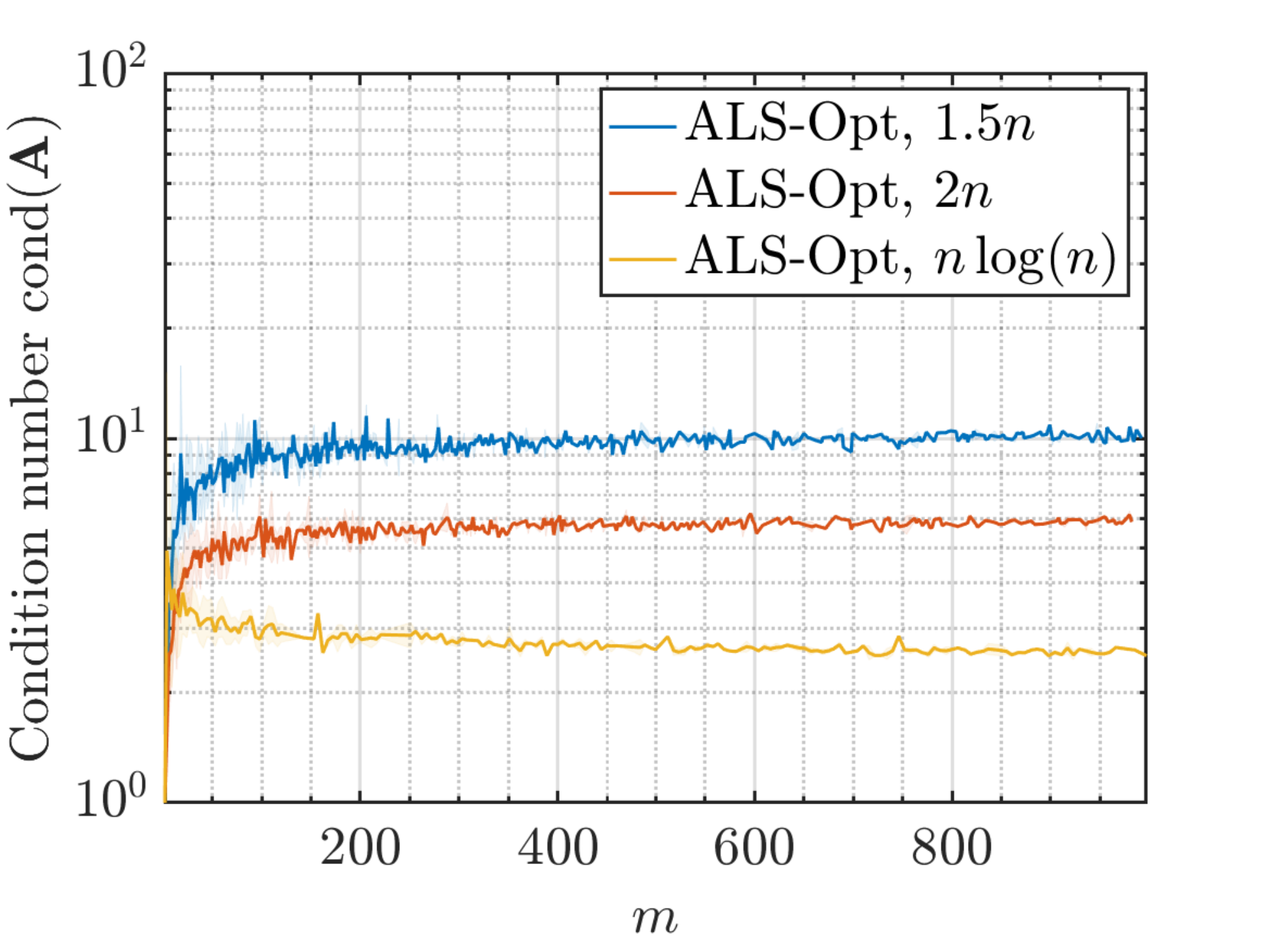}&
\includegraphics[width = \errplotimg]{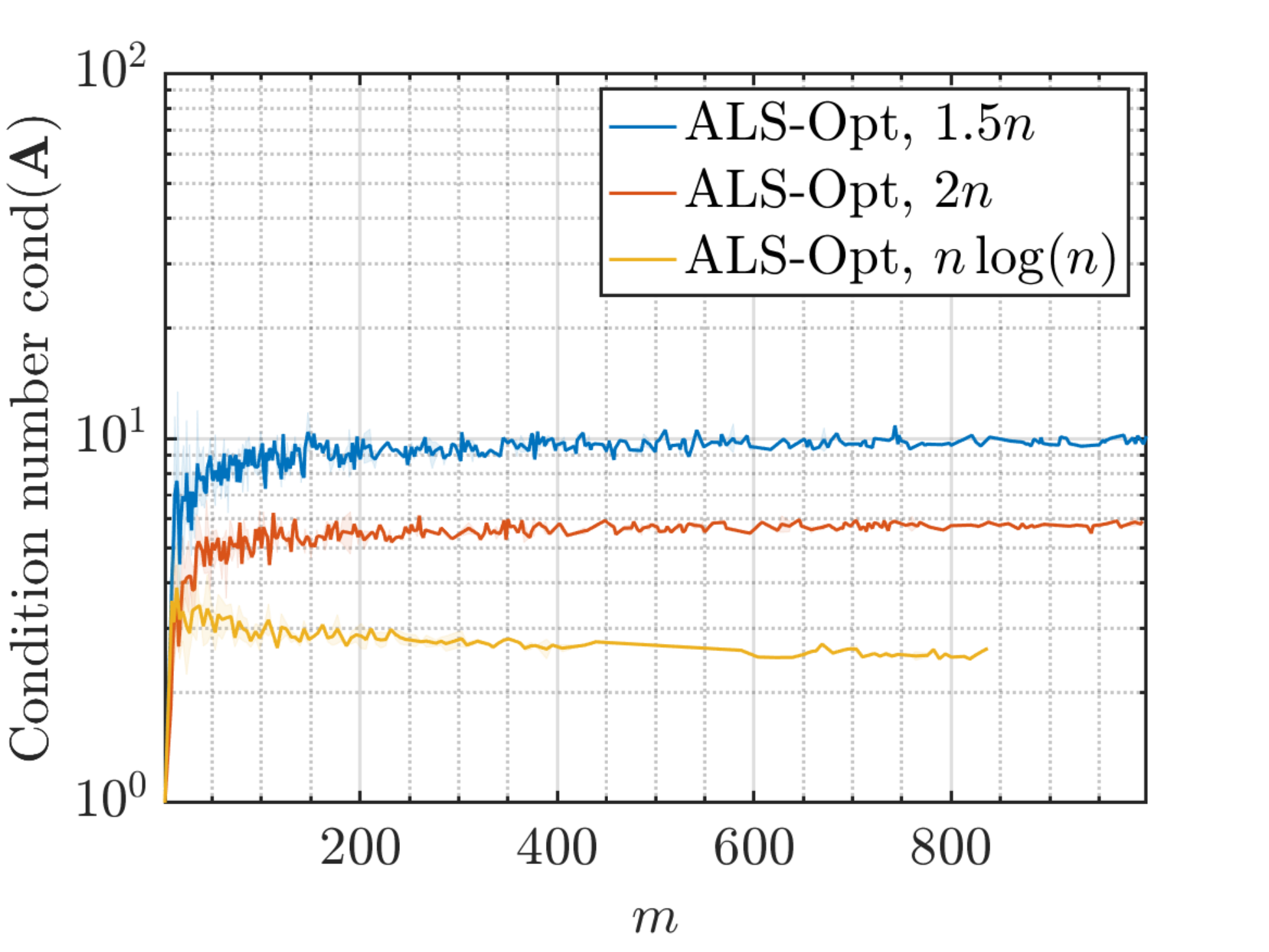}
\\[\errplottextsp]
$d = 8$ & $d = 16$ & $d = 32$
\end{tabular}
\end{small}
\end{center}
\caption{The same as Fig.\ \ref{fig:fig10}, except with additional graphs corresponding to ALS with the linear scalings \eqref{m-scaling-LS-linear}. The first and third rows show the approximation errors. The second and fourth rows show the condition numbers for the ALS approximations.} 
\label{fig:figSM6}
\end{figure}

\begin{figure}[t!]
\begin{center}
\begin{small}
 \begin{tabular}{@{\hspace{0pt}}c@{\hspace{\errplotsp}}c@{\hspace{\errplotsp}}c@{\hspace{0pt}}}
\includegraphics[width = \errplotimg]{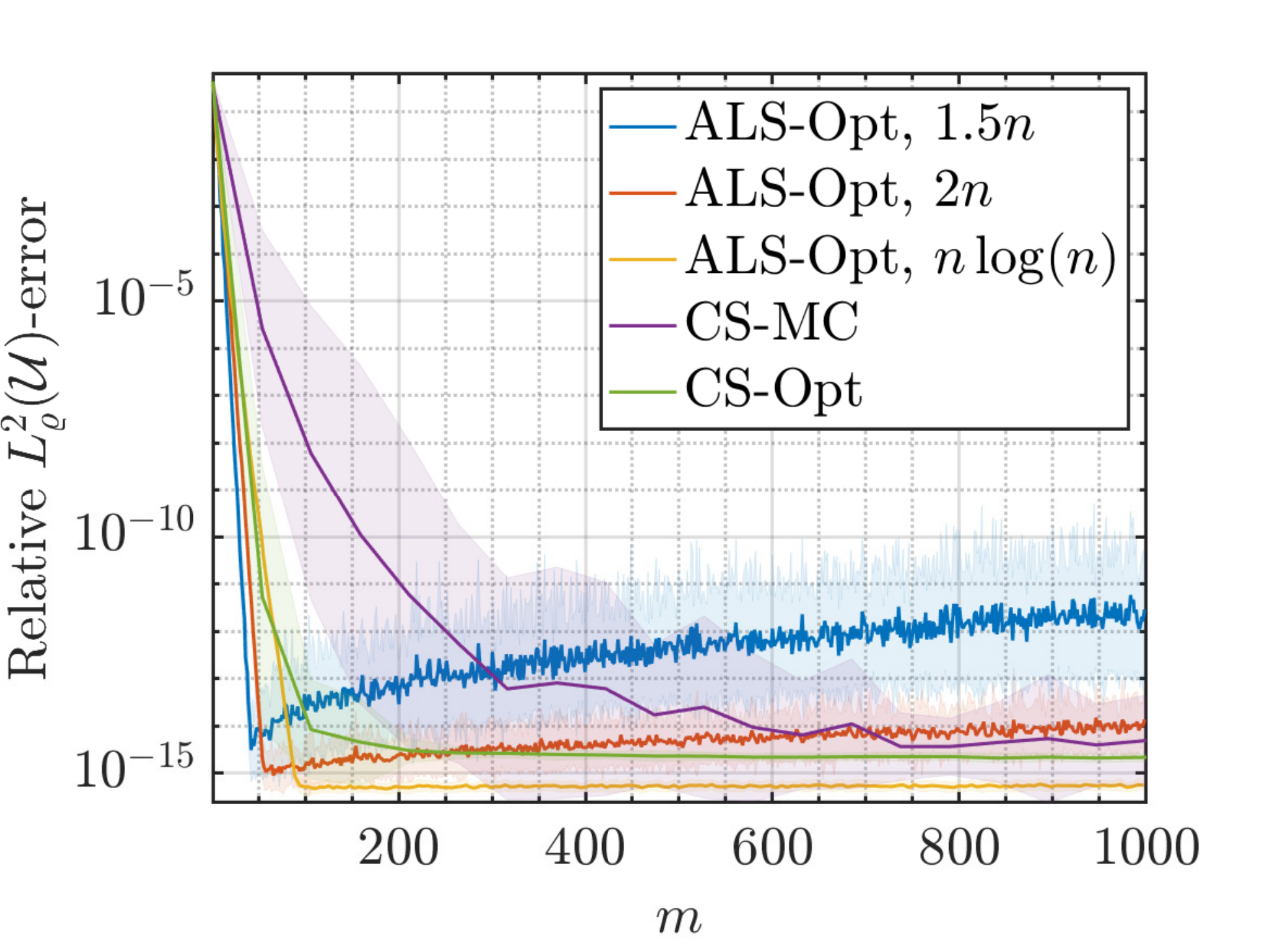}&
\includegraphics[width = \errplotimg]{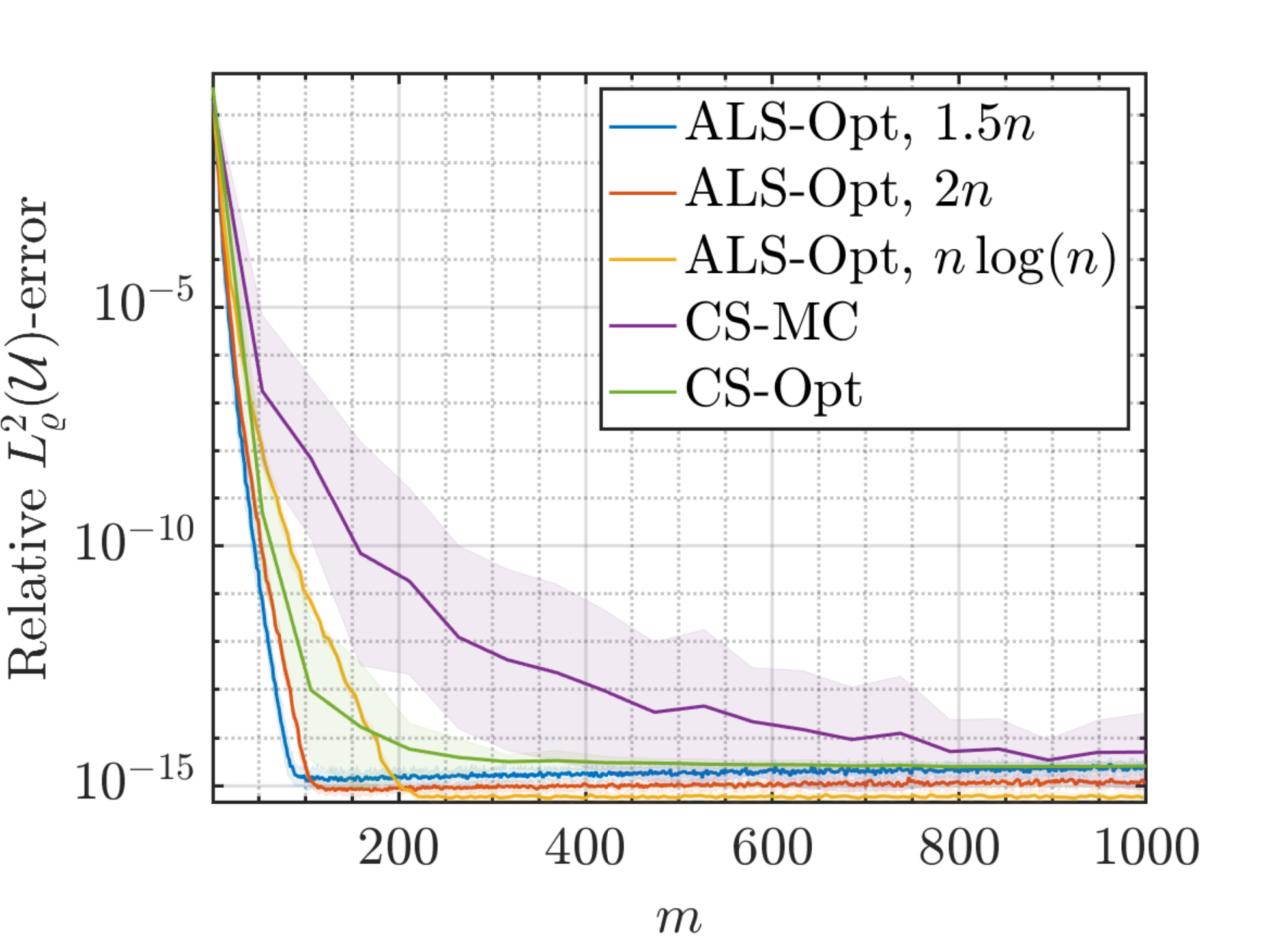}&
\includegraphics[width = \errplotimg]{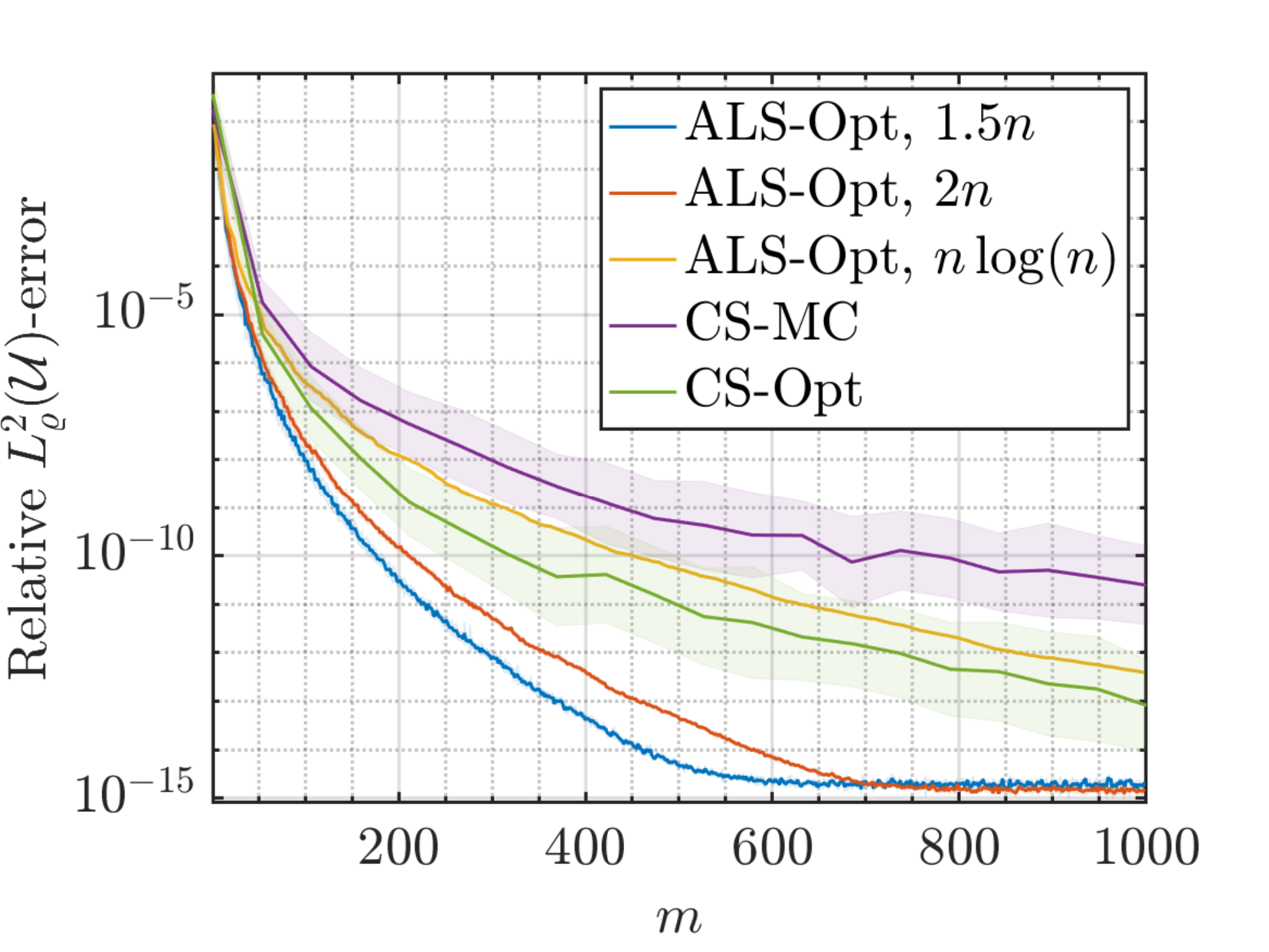}
\\[\errplotgraphsp]
\includegraphics[width = \errplotimg]{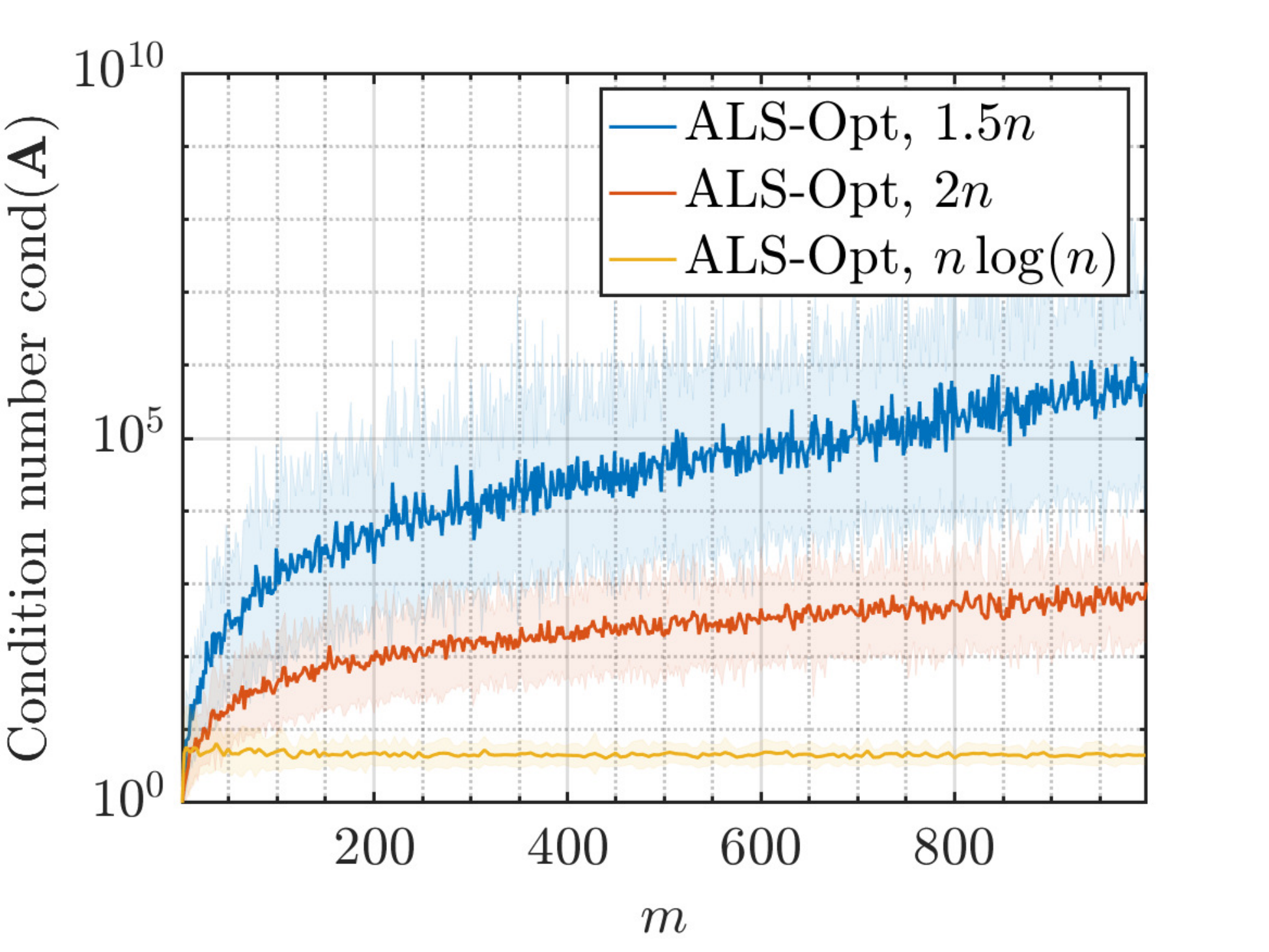}&
\includegraphics[width = \errplotimg]{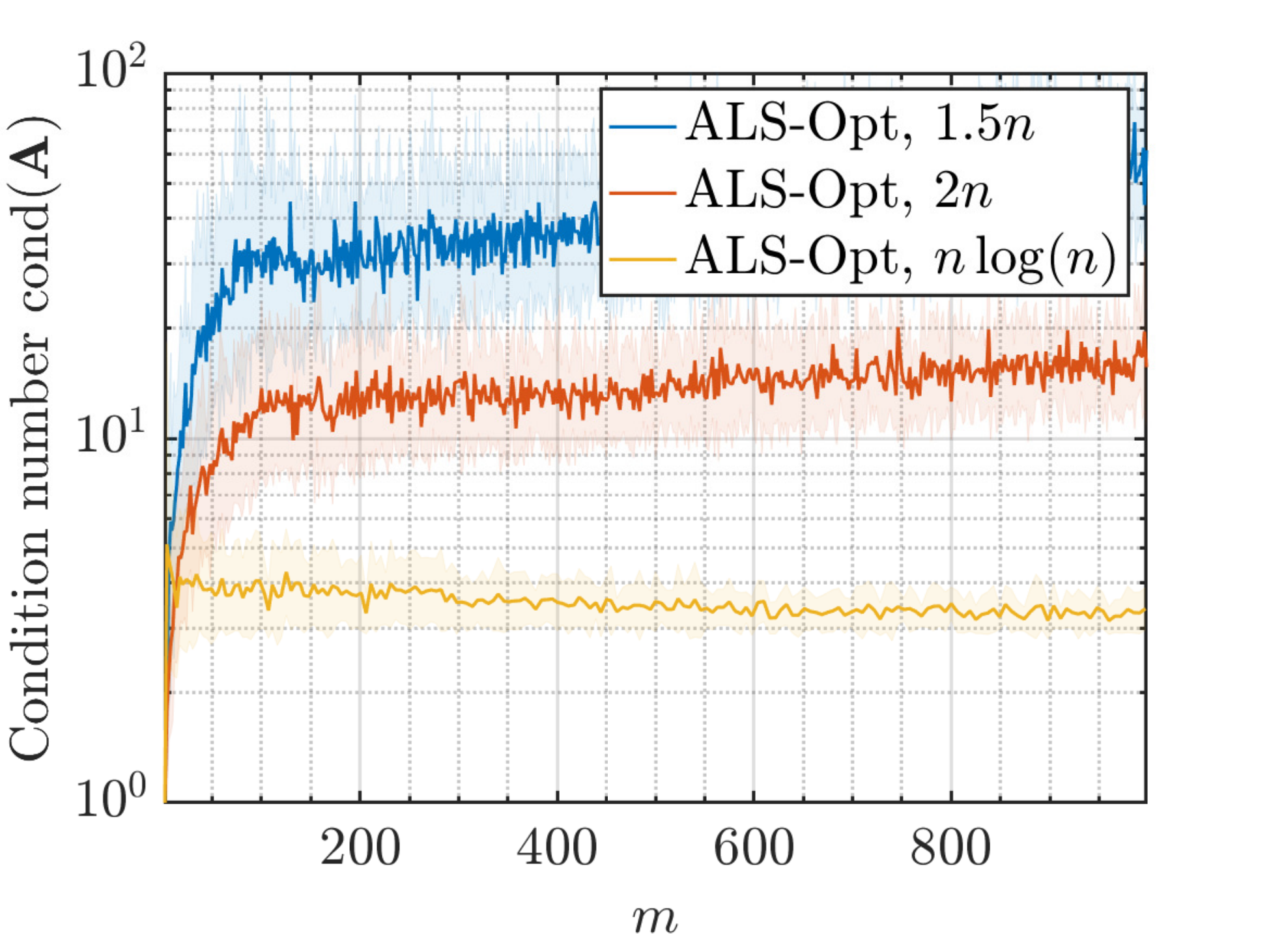}&
\includegraphics[width = \errplotimg]{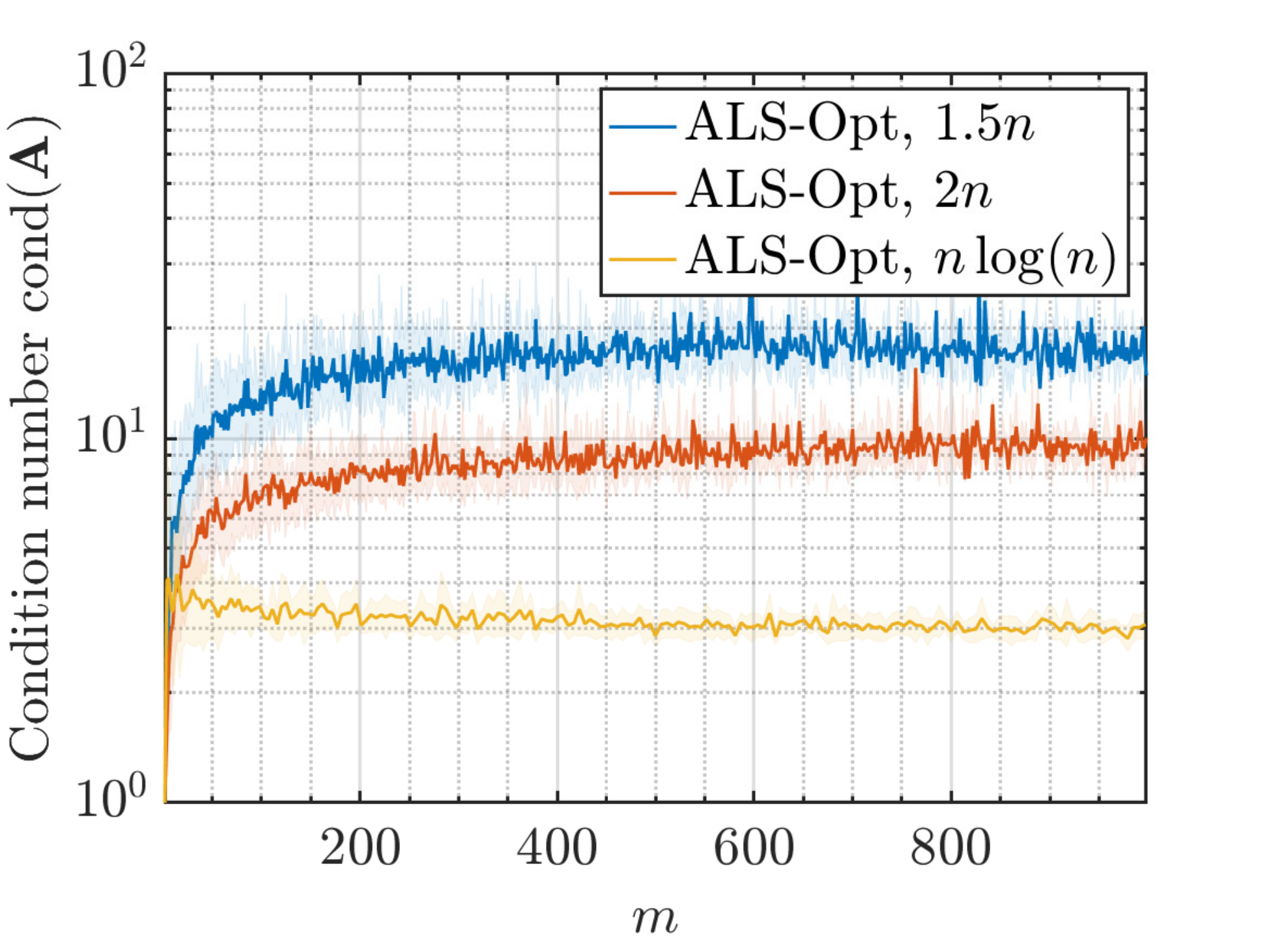}
\\[\errplottextsp]
$d = 1$ & $d = 2$ & $d = 4$
\\[\errplottextsp]
\includegraphics[width = \errplotimg]{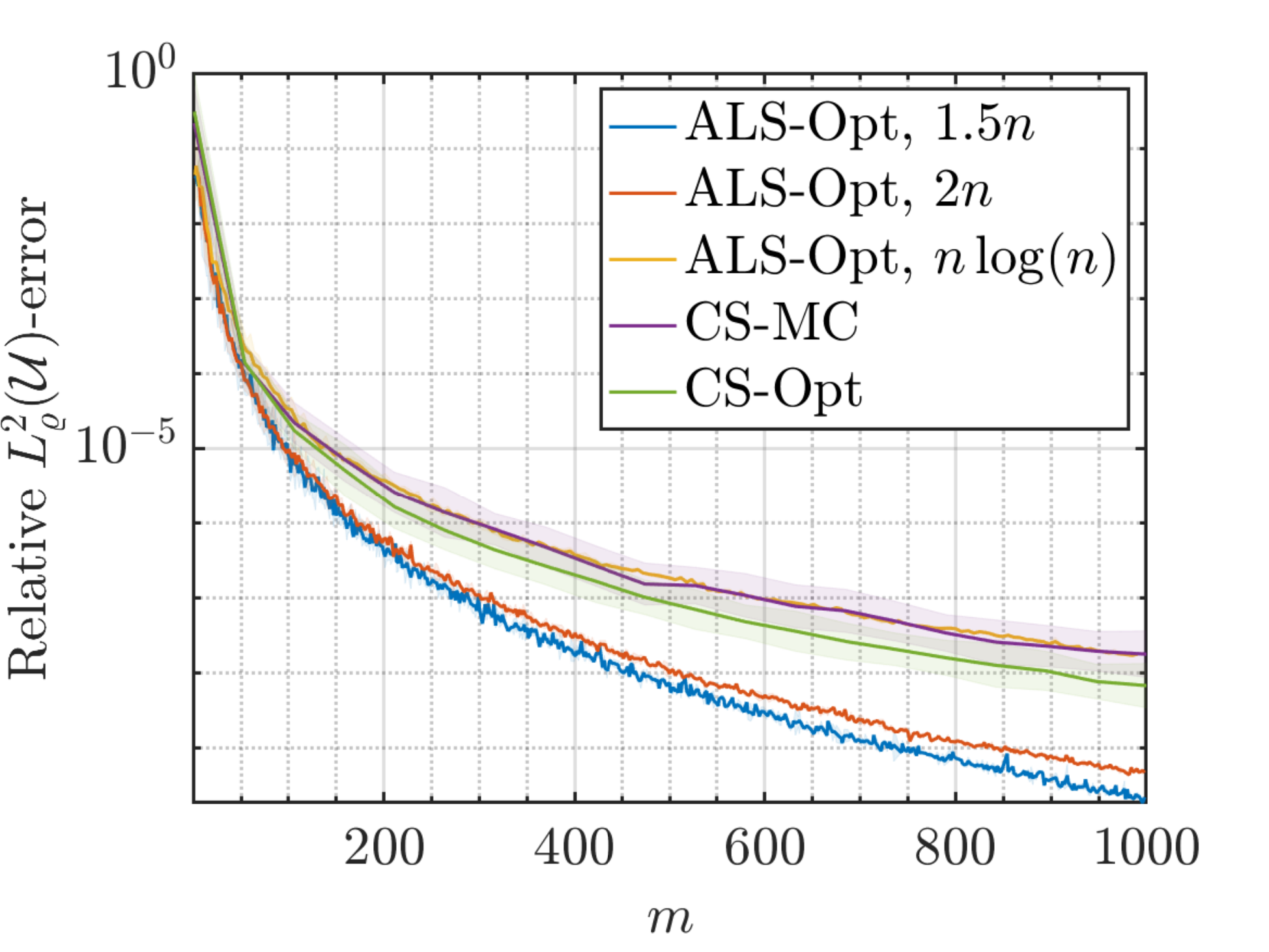}&
\includegraphics[width = \errplotimg]{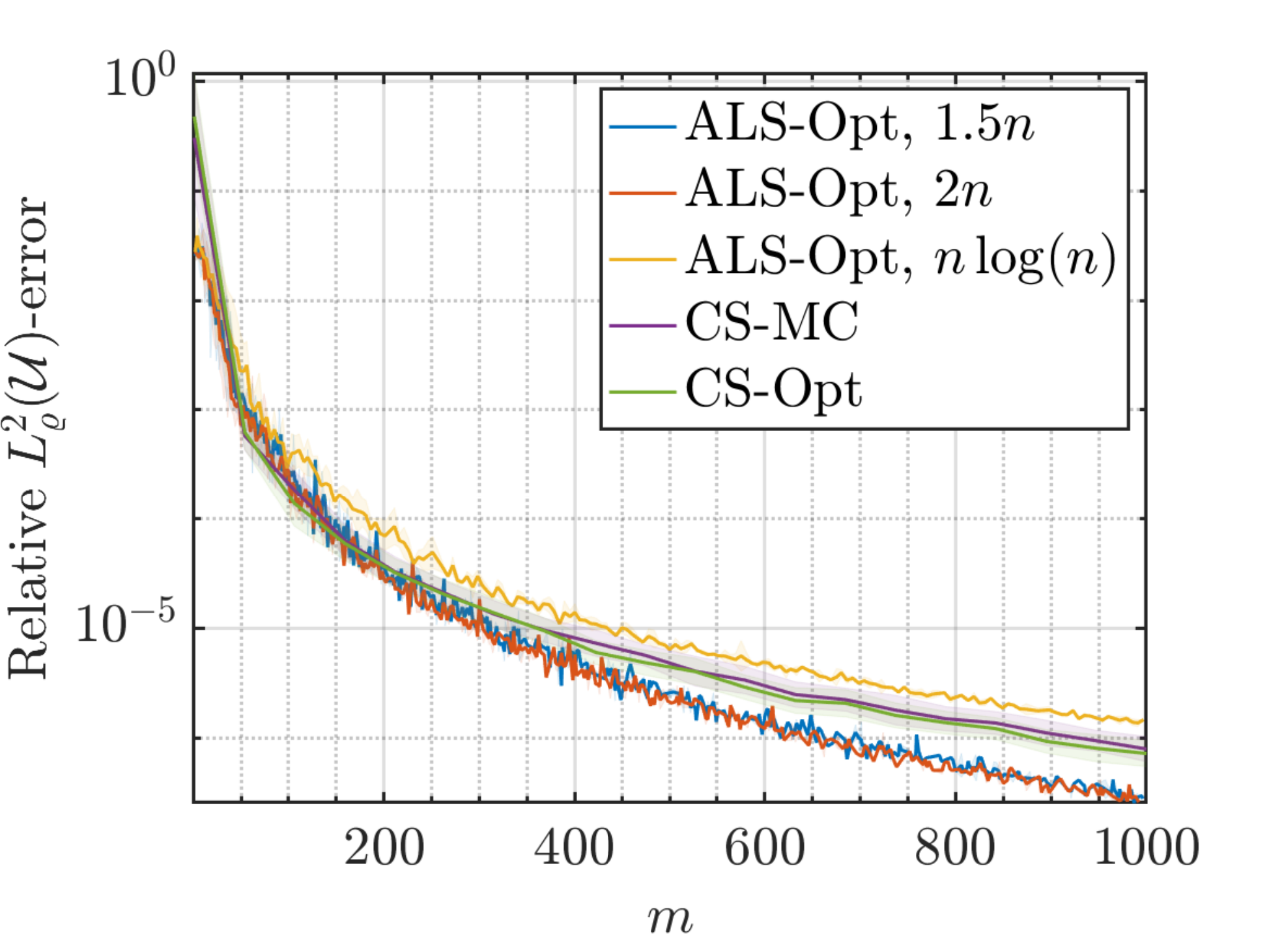}&
\includegraphics[width = \errplotimg]{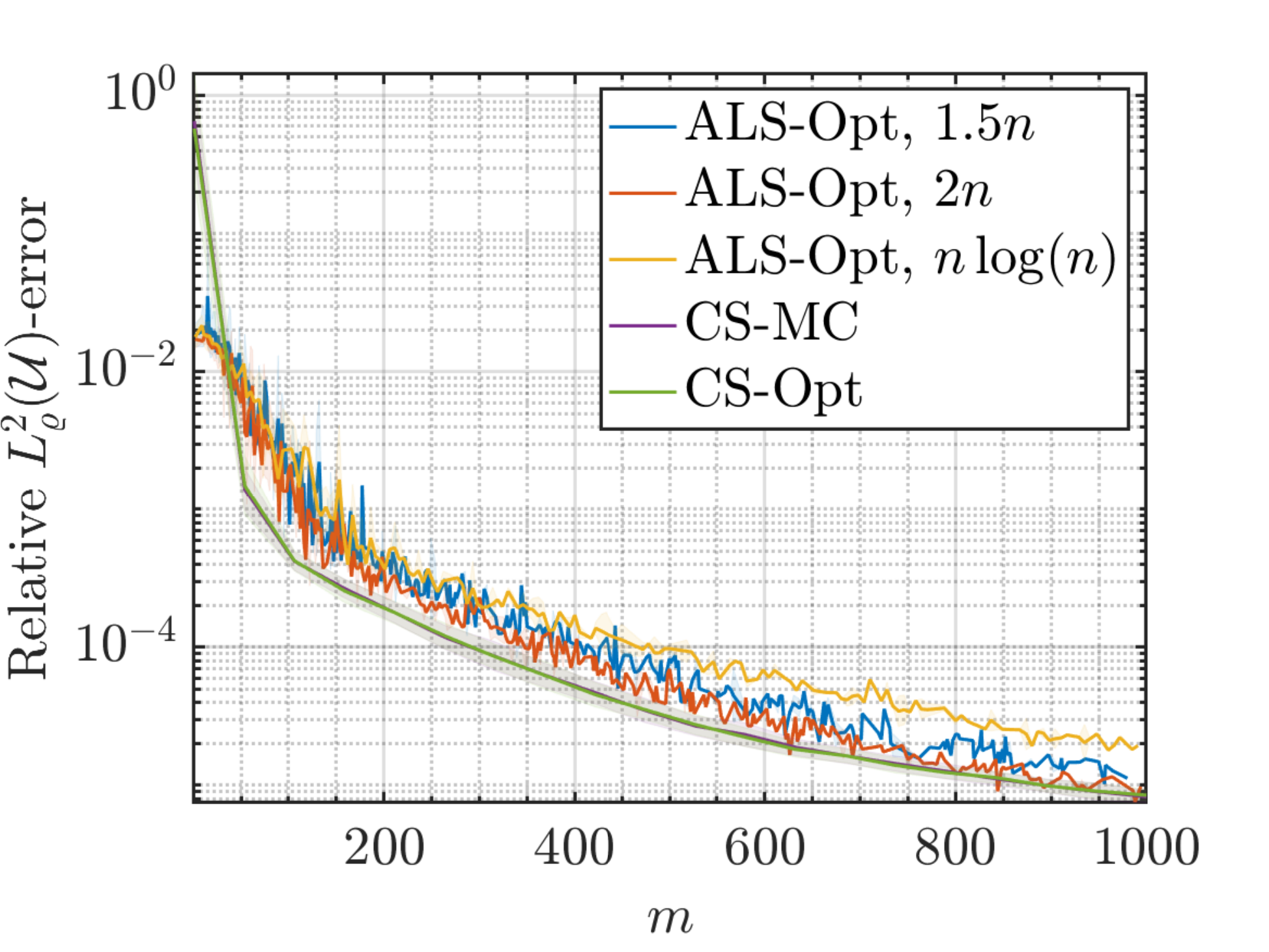}
\\[\errplotgraphsp]
\includegraphics[width = \errplotimg]{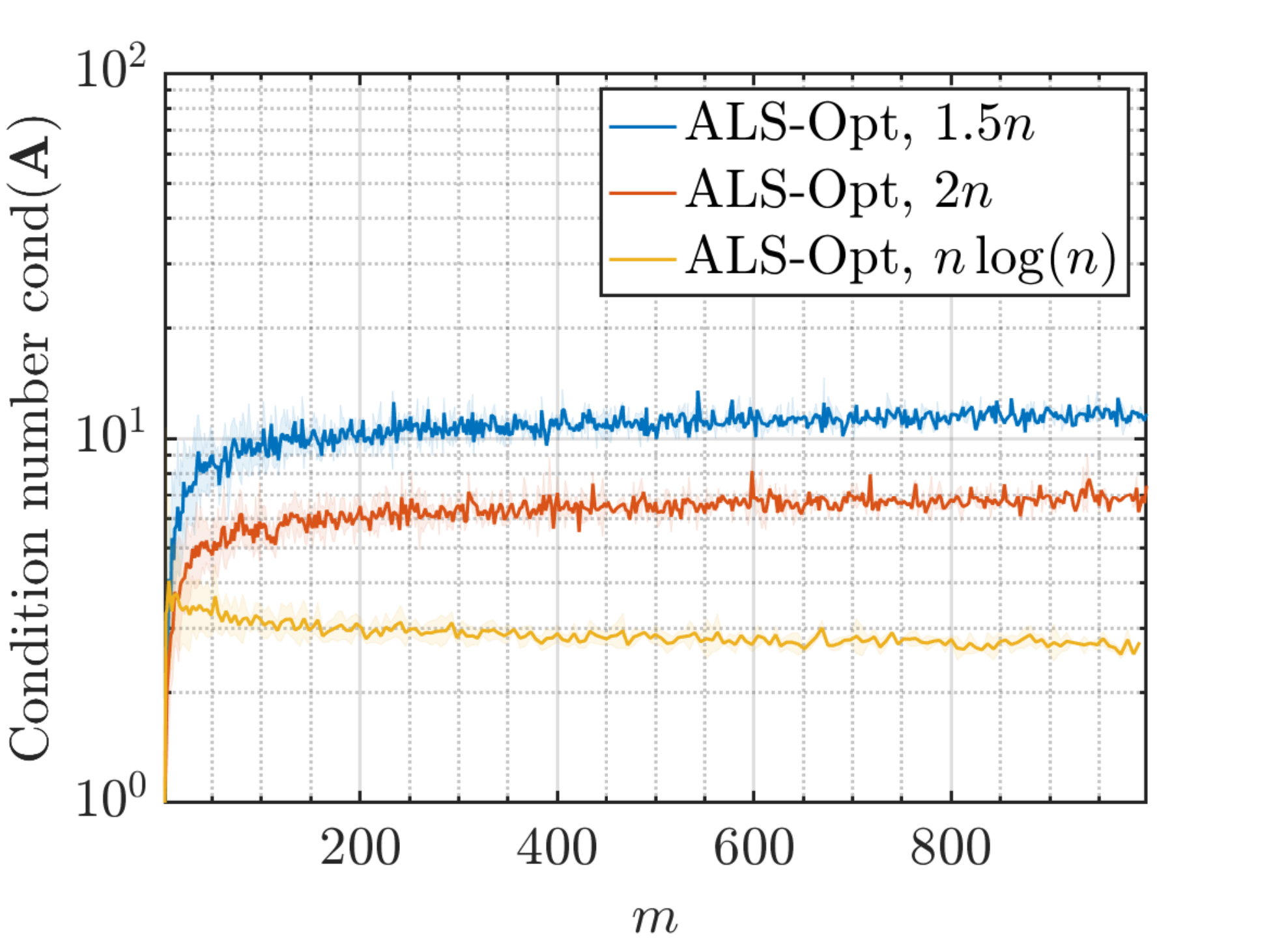}&
\includegraphics[width = \errplotimg]{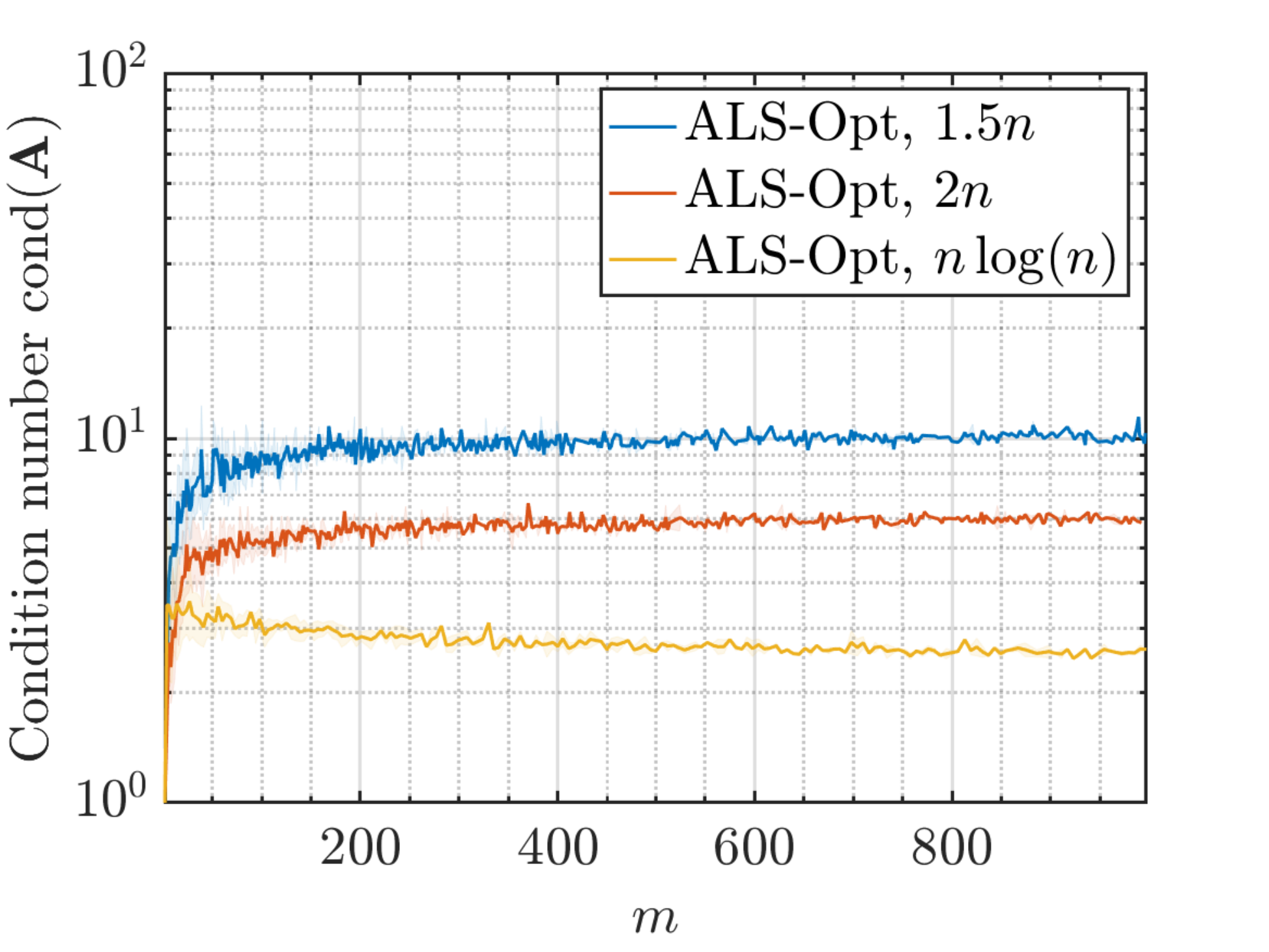}&
\includegraphics[width = \errplotimg]{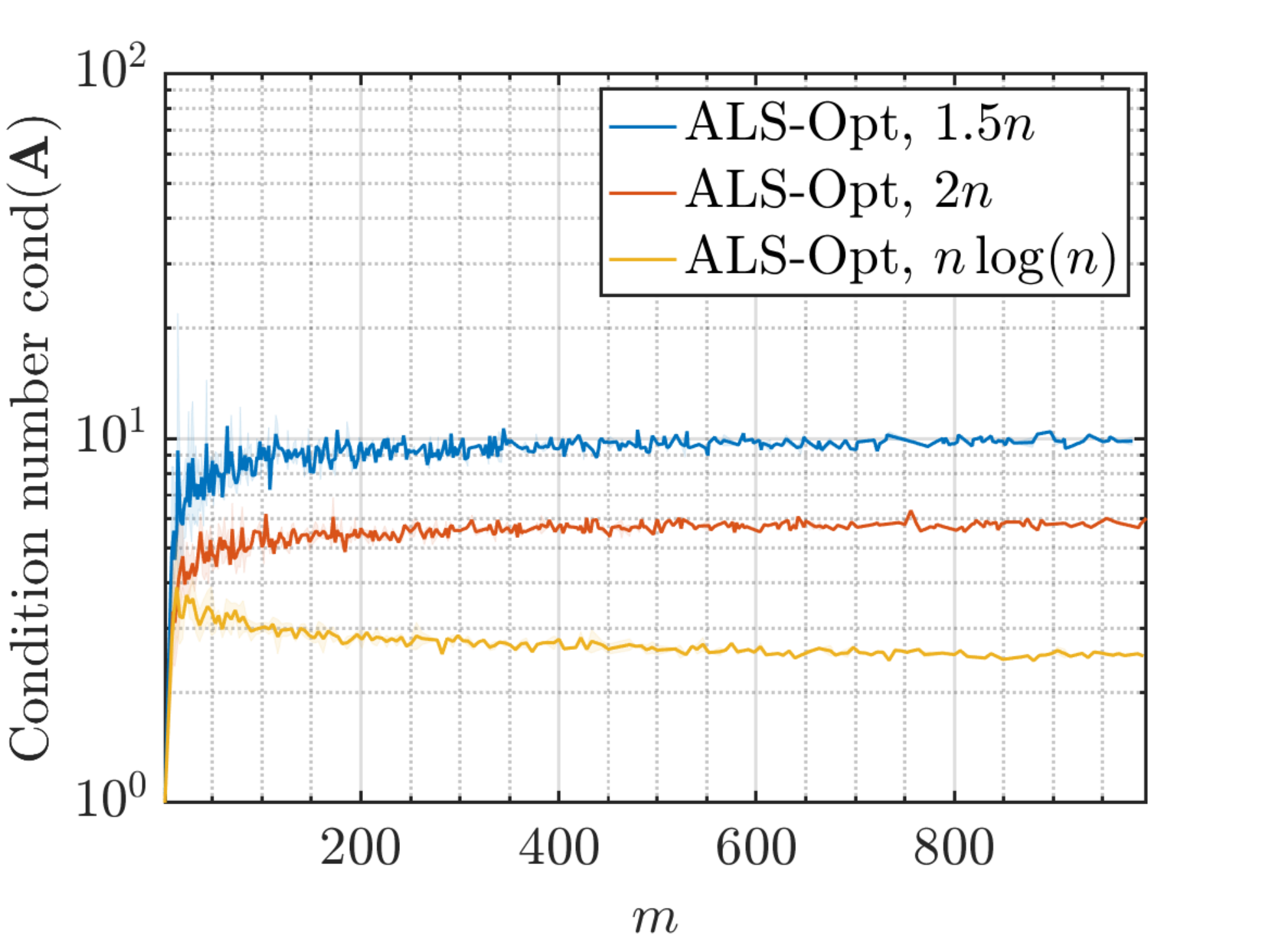}
\\[\errplottextsp]
$d = 8$ & $d = 16$ & $d = 32$
\end{tabular}
\end{small}
\end{center}
\caption{The same as Fig.\ \ref{fig:fig12}, except with the modifications made in Fig.\ \ref{fig:figSM6}.} 
\label{fig:figSM7}
\end{figure}

\begin{figure}[t!]
\begin{center}
\begin{small}
 \begin{tabular}{@{\hspace{0pt}}c@{\hspace{\errplotsp}}c@{\hspace{\errplotsp}}c@{\hspace{0pt}}}
\includegraphics[width = \errplotimg]{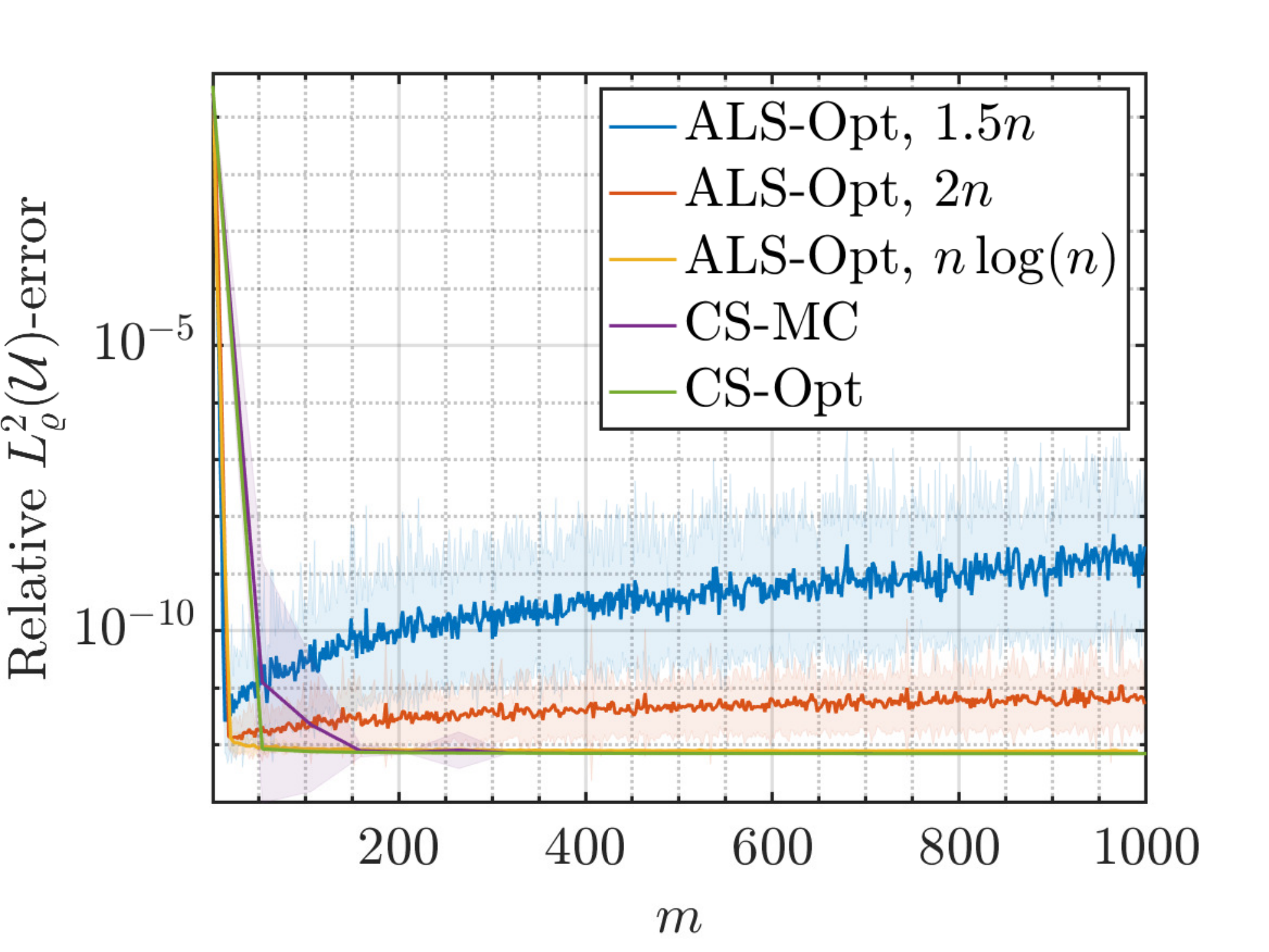}&
\includegraphics[width = \errplotimg]{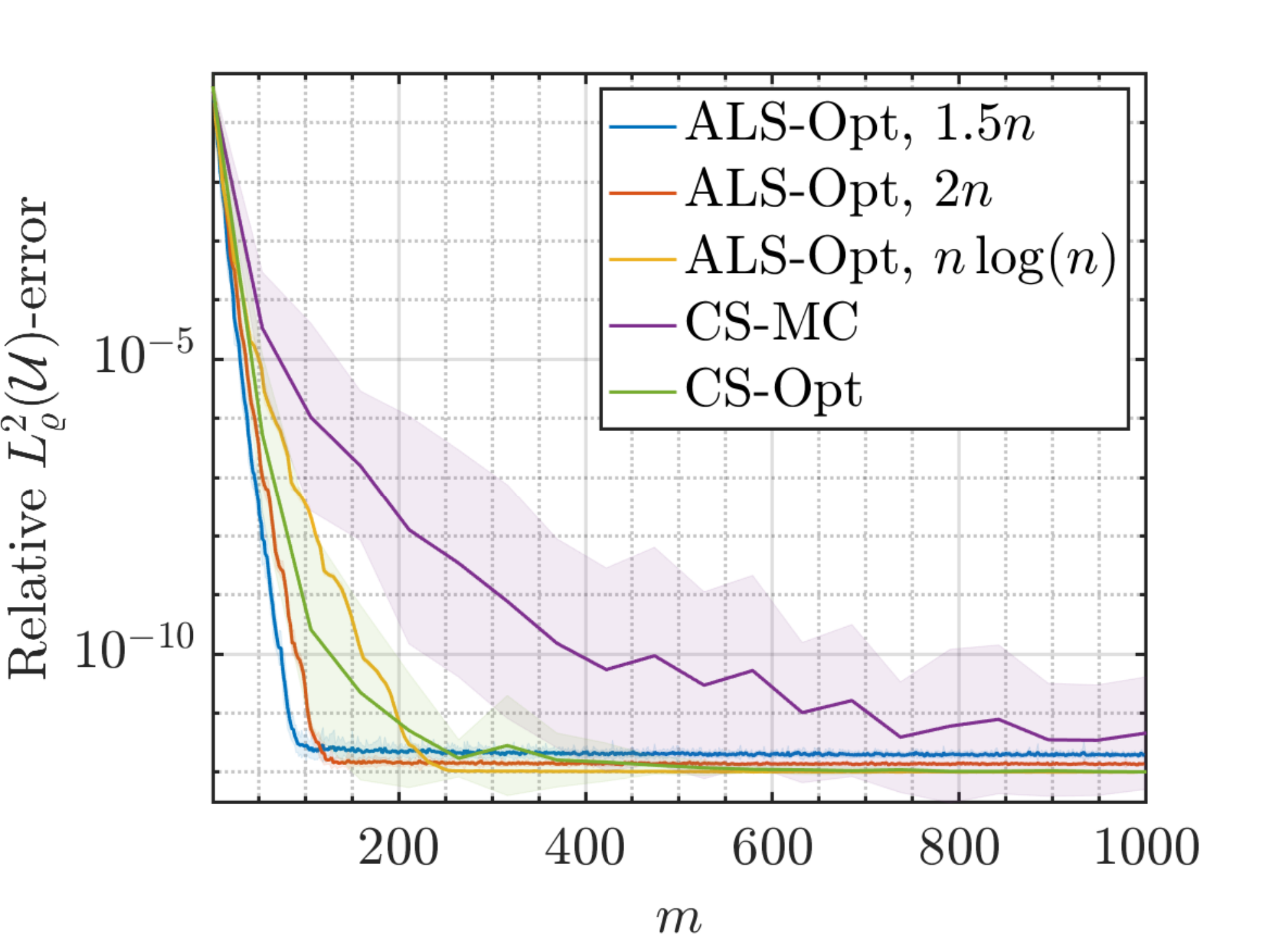}&
\includegraphics[width = \errplotimg]{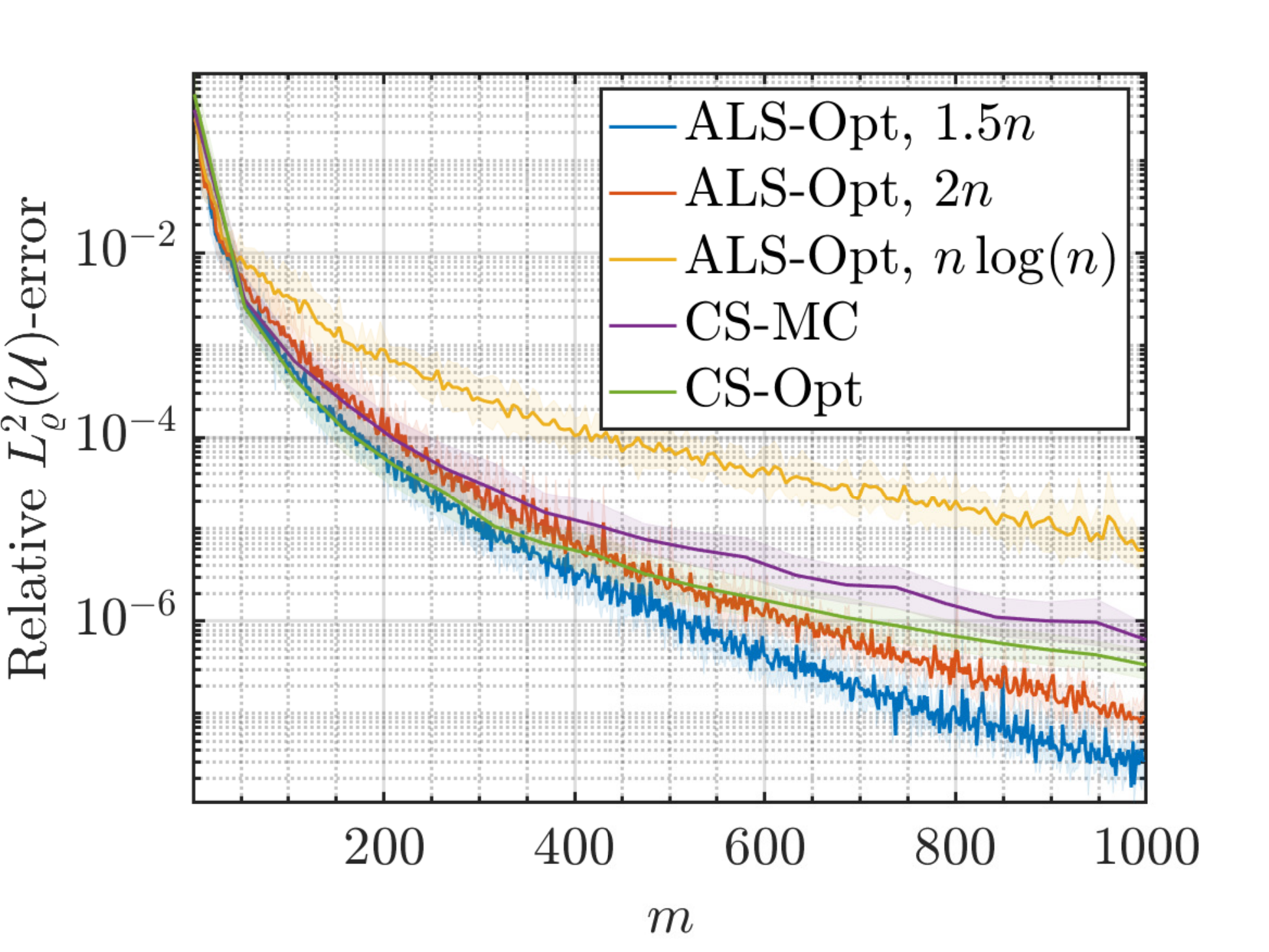}
\\[\errplotgraphsp]
\includegraphics[width = \errplotimg]{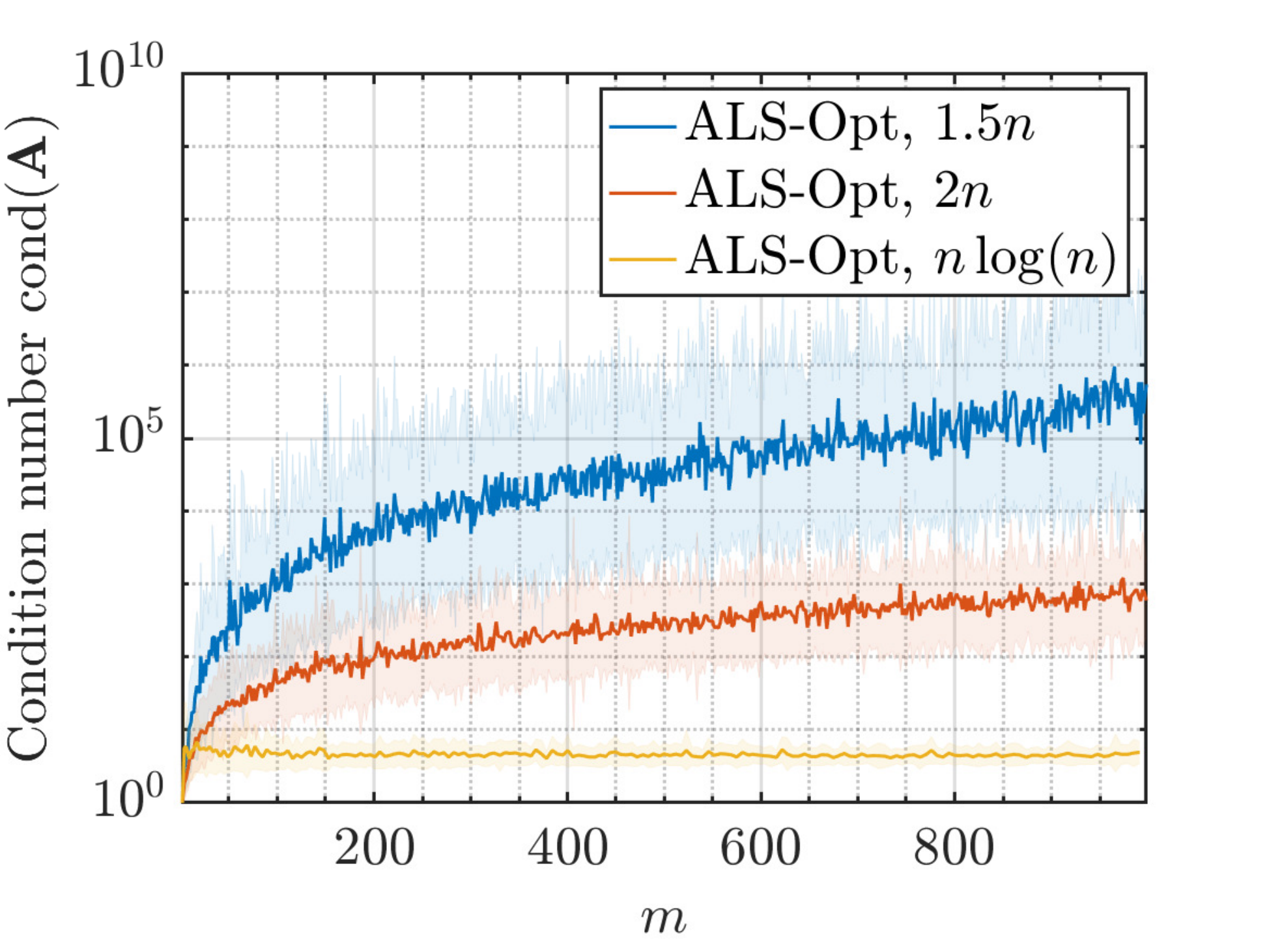}&
\includegraphics[width = \errplotimg]{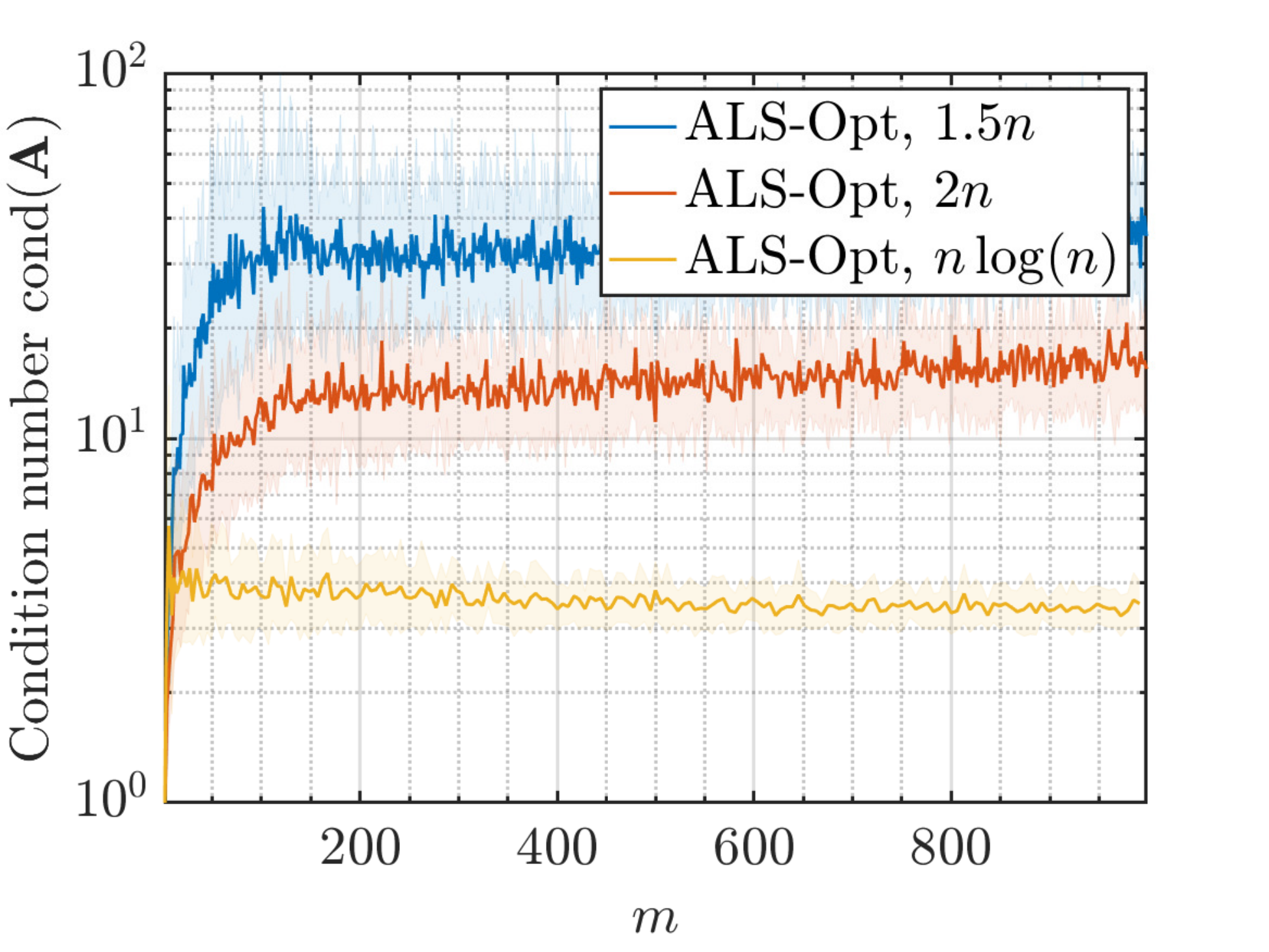}&
\includegraphics[width = \errplotimg]{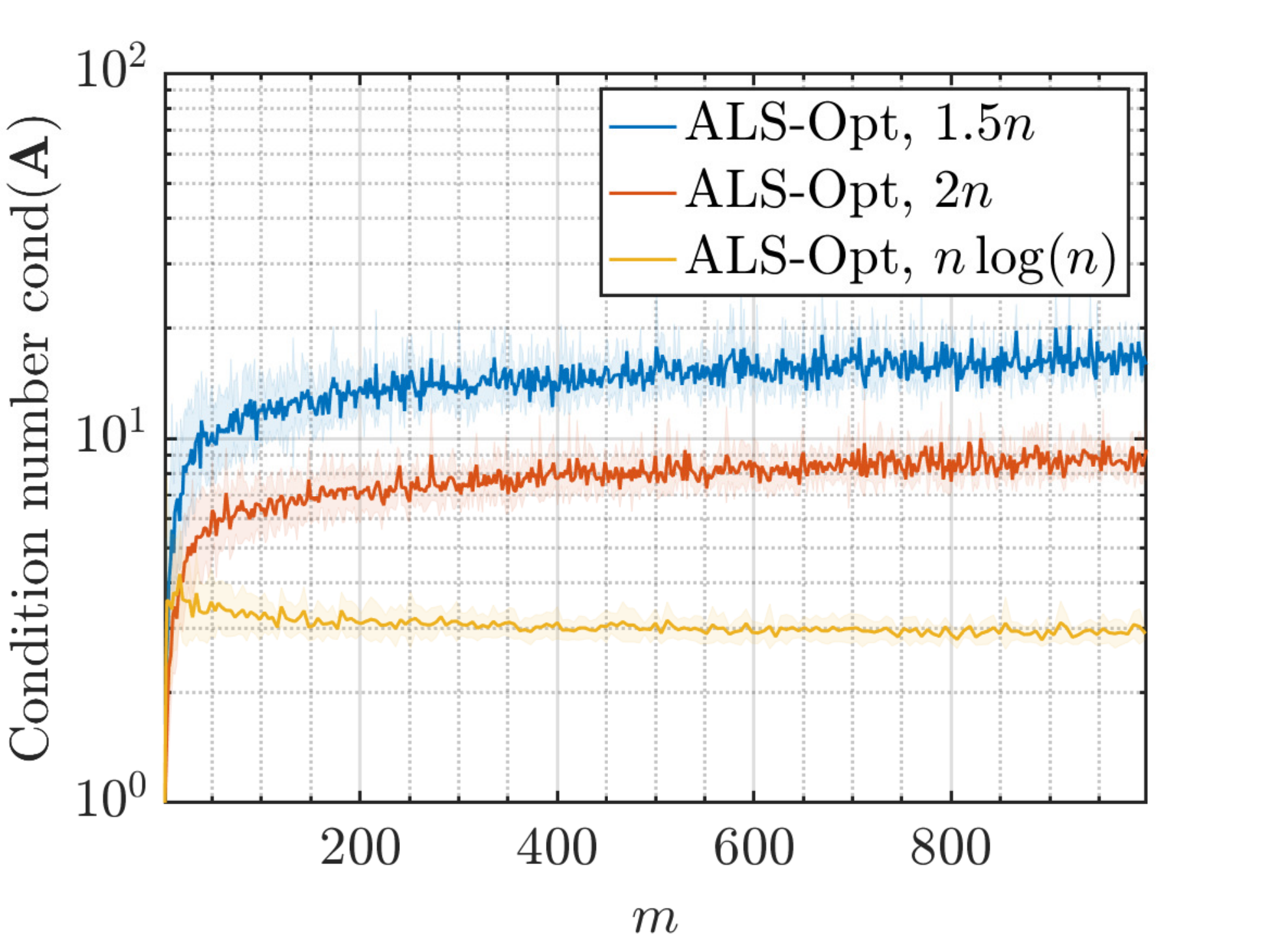}
\\[\errplottextsp]
$d = 1$ & $d = 2$ & $d = 4$
\\[\errplottextsp]
\includegraphics[width = \errplotimg]{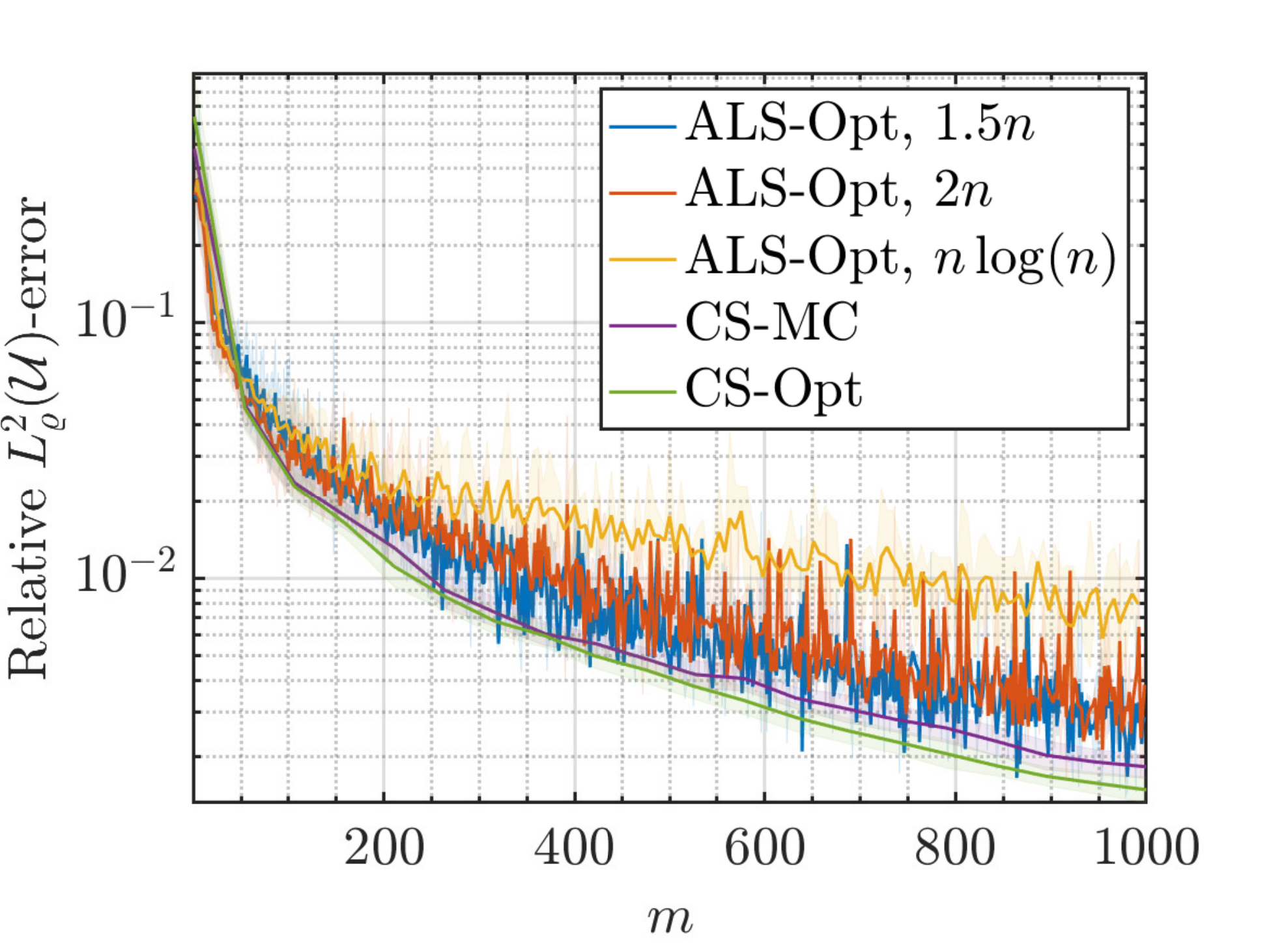}&
\includegraphics[width = \errplotimg]{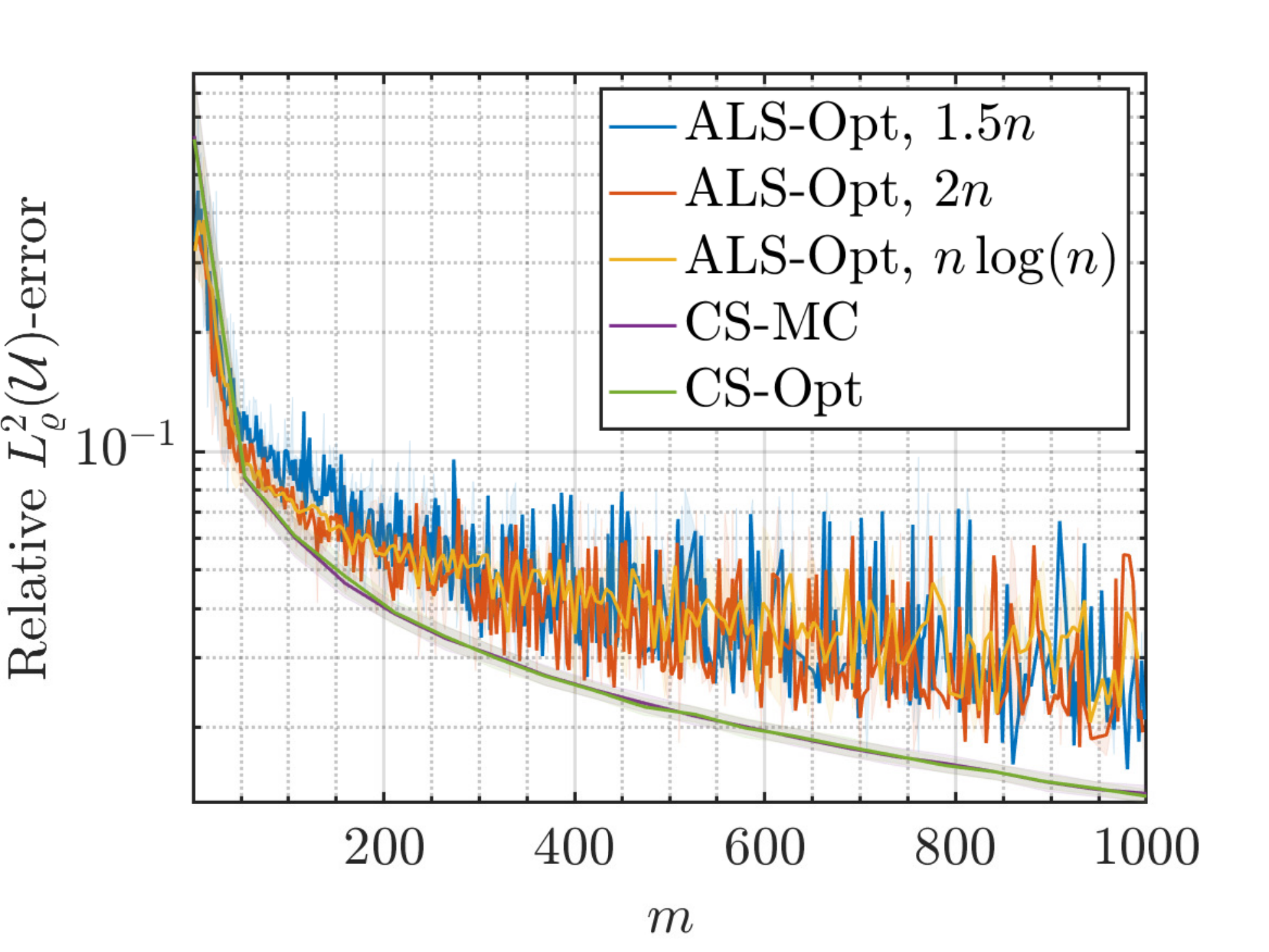}&
\includegraphics[width = \errplotimg]{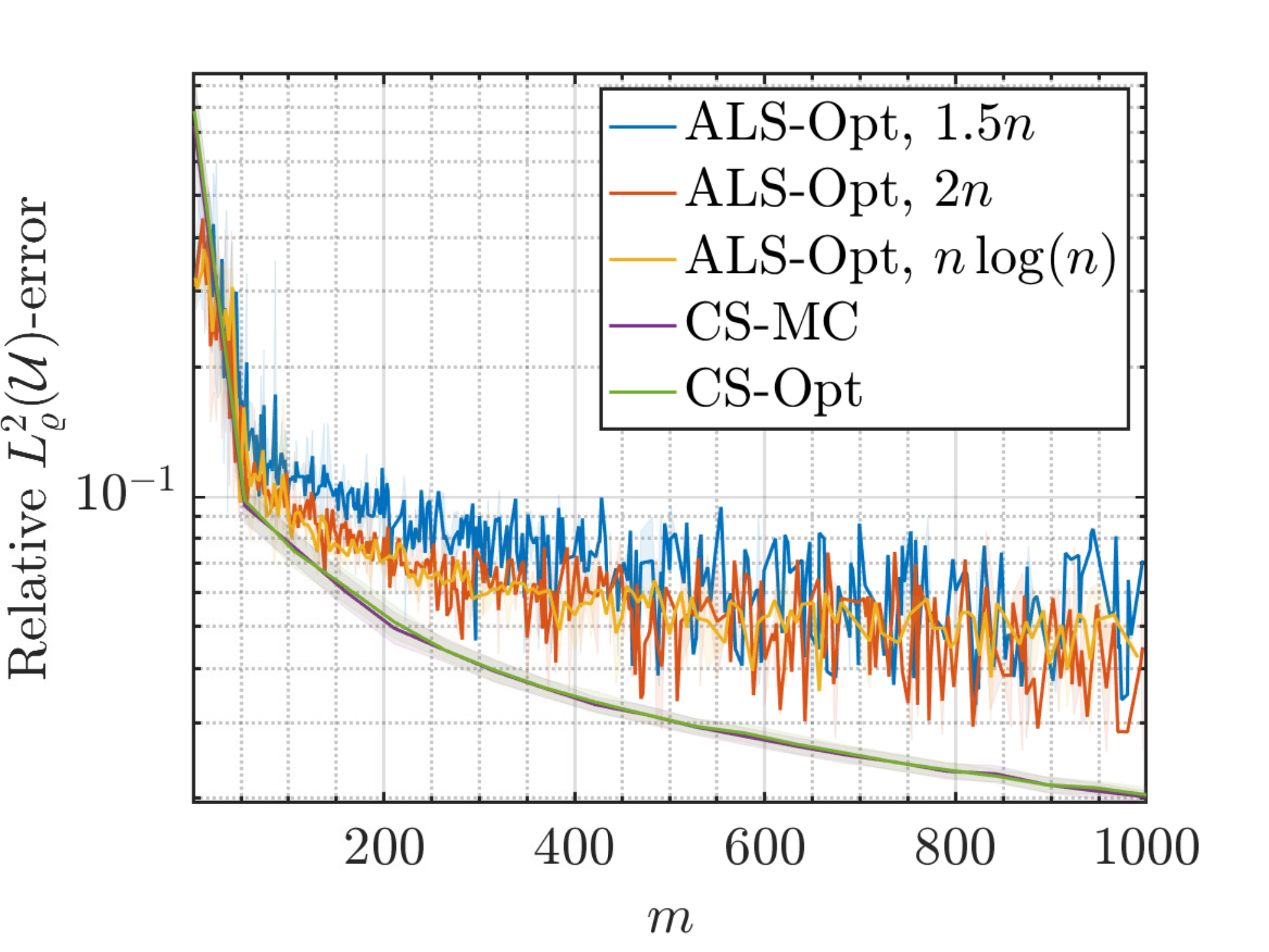}
\\[\errplotgraphsp]
\includegraphics[width = \errplotimg]{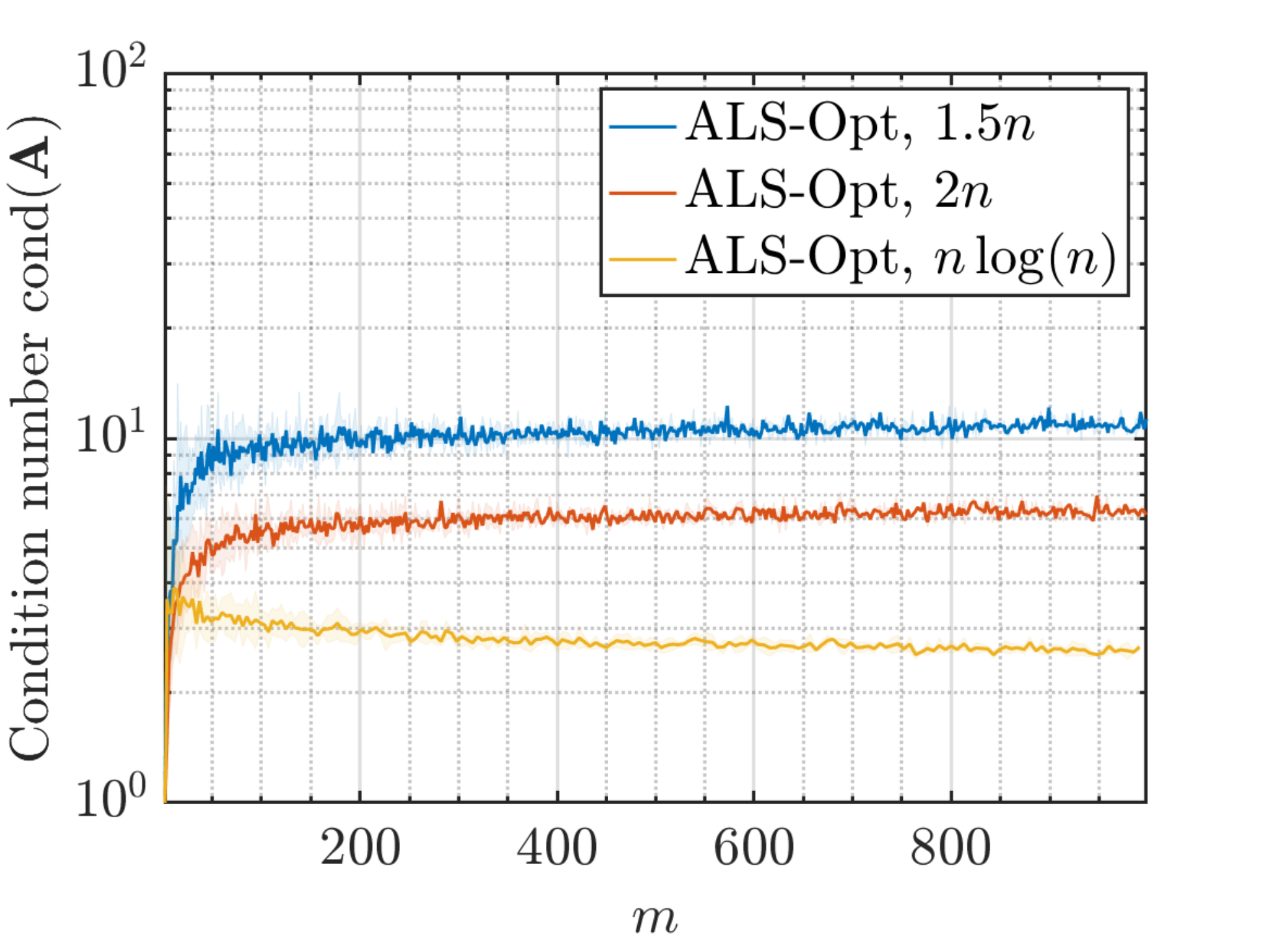}&
\includegraphics[width = \errplotimg]{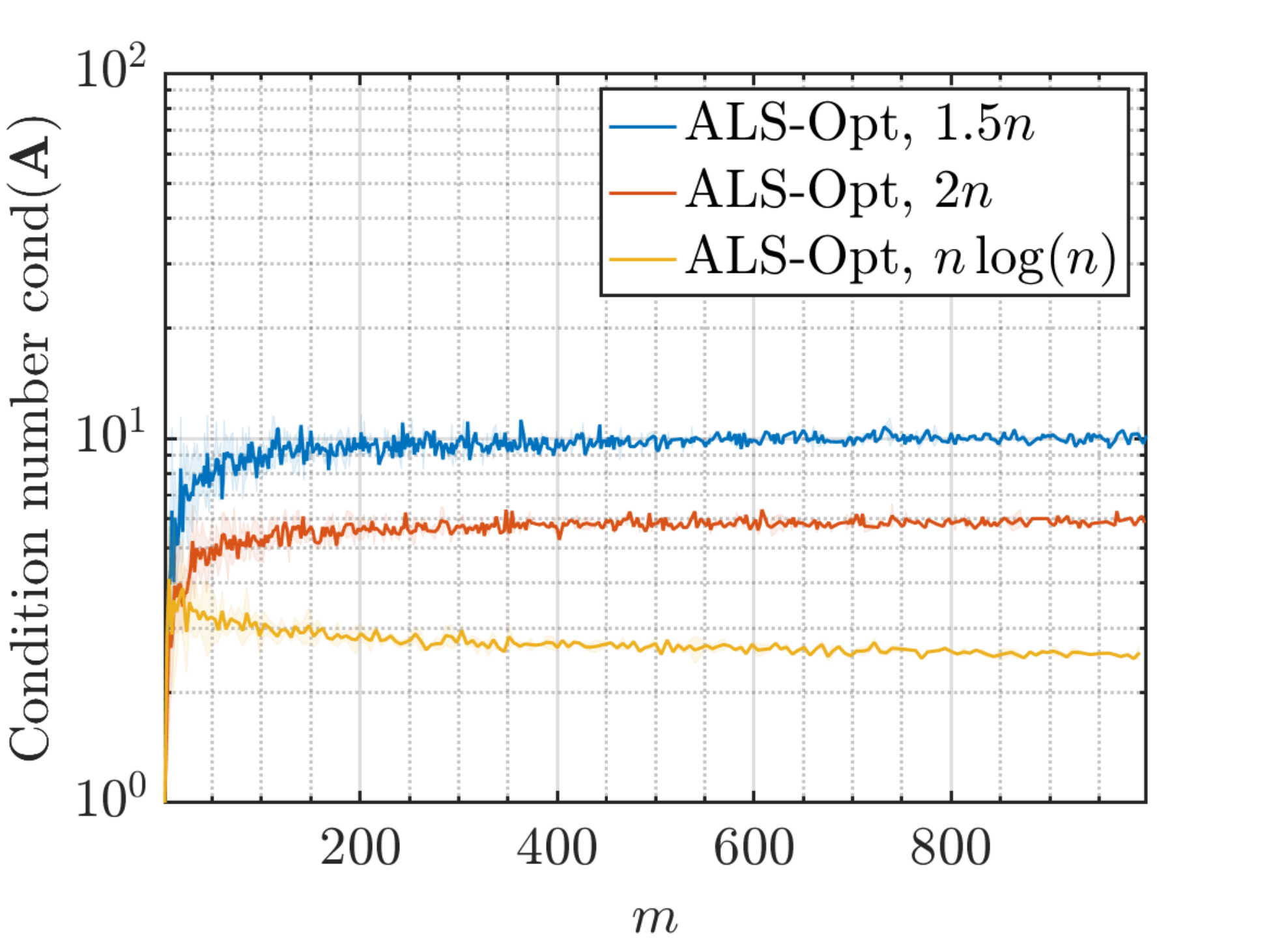}&
\includegraphics[width = \errplotimg]{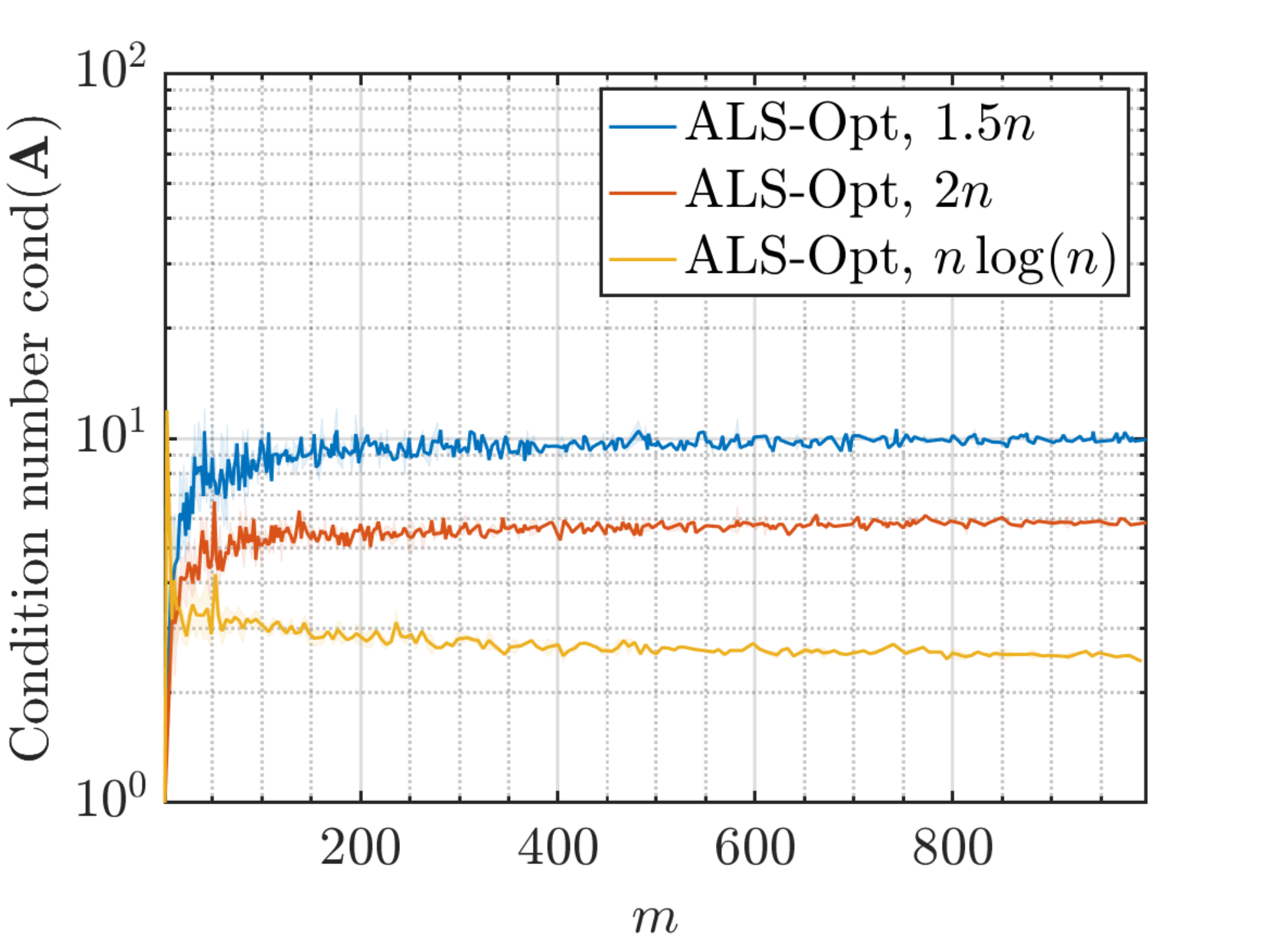}
\\[\errplottextsp]
$d = 8$ & $d = 16$ & $d = 32$
\end{tabular}
\end{small}
\end{center}
\caption{The same as Fig.\ \ref{fig:fig13}, except with the modifications made in Fig.\ \ref{fig:figSM6}. }
\label{fig:figSM8}
\end{figure}

Specifically, Figs.\ \ref{fig:figSM6}--\ref{fig:figSM8} we show the same experiments as in Figs.\ \ref{fig:fig10}--\ref{fig:fig14}, but we now also include the two linear oversampling factors
\be{
\label{m-scaling-LS-linear}
m = \lceil 1.5 n^{(l)} \rceil,\qquad m = 2 n^{(l)}.
}
As is evident, in lower dimensions these scalings can lead to somewhat better performance than CS. However, this comes at the price of a much worse condition number. Indeed, in $d = 1$ dimensions the condition numbers can grow to roughly $10^7$, while the condition number of log-linear oversampling remains less than $10^1$. Furthermore, the benefits of these scalings in terms of error reduction diminishes as the dimension increases.

\subsection{Uniform-norm experiments}\label{ss:uniform-norm-exp}

\begin{figure}[t!]
\begin{center}
\begin{small}
 \begin{tabular}{@{\hspace{0pt}}c@{\hspace{\errplotsp}}c@{\hspace{\errplotsp}}c@{\hspace{0pt}}}
\includegraphics[width = \errplotimg]{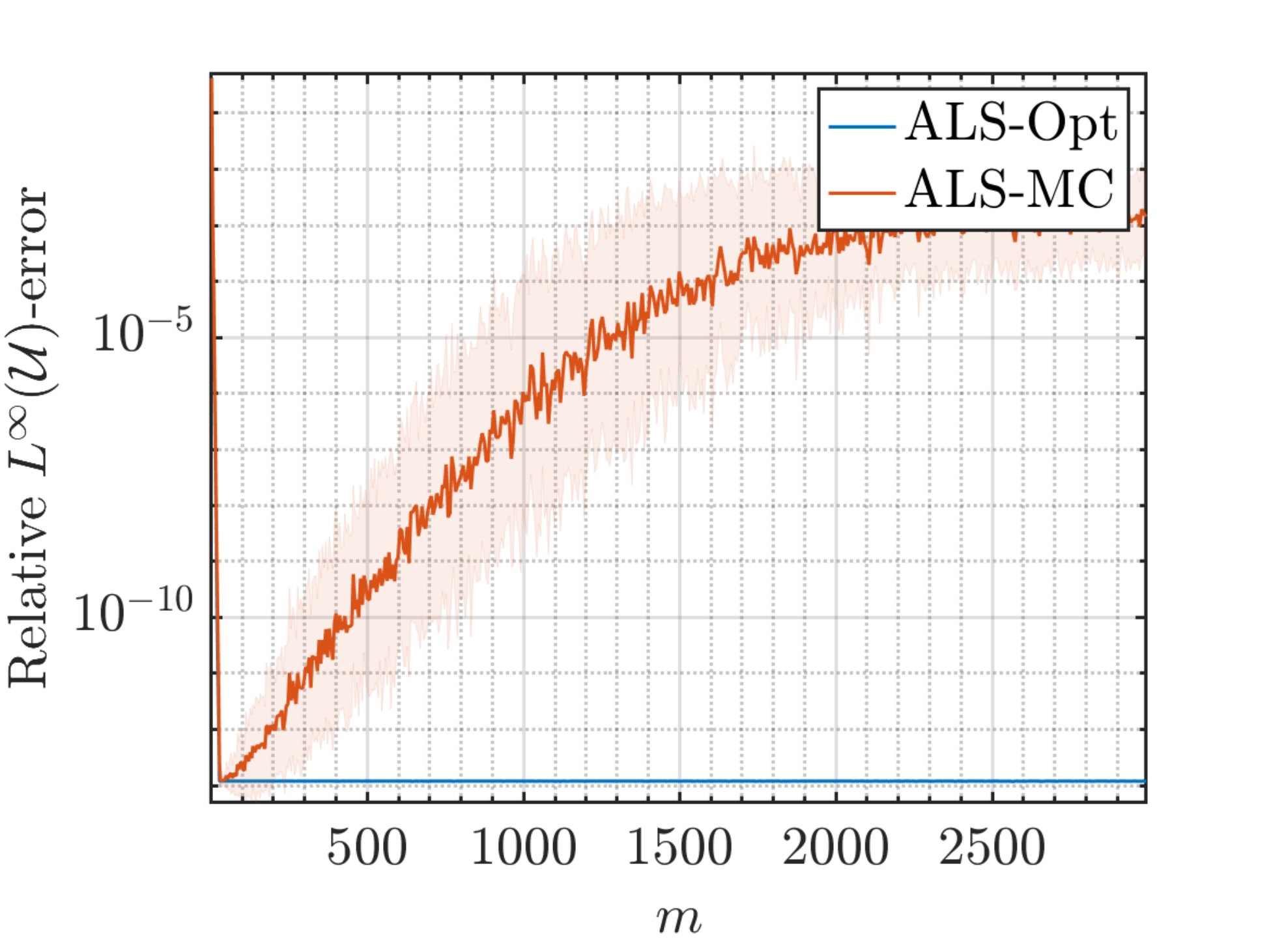}&
\includegraphics[width = \errplotimg]{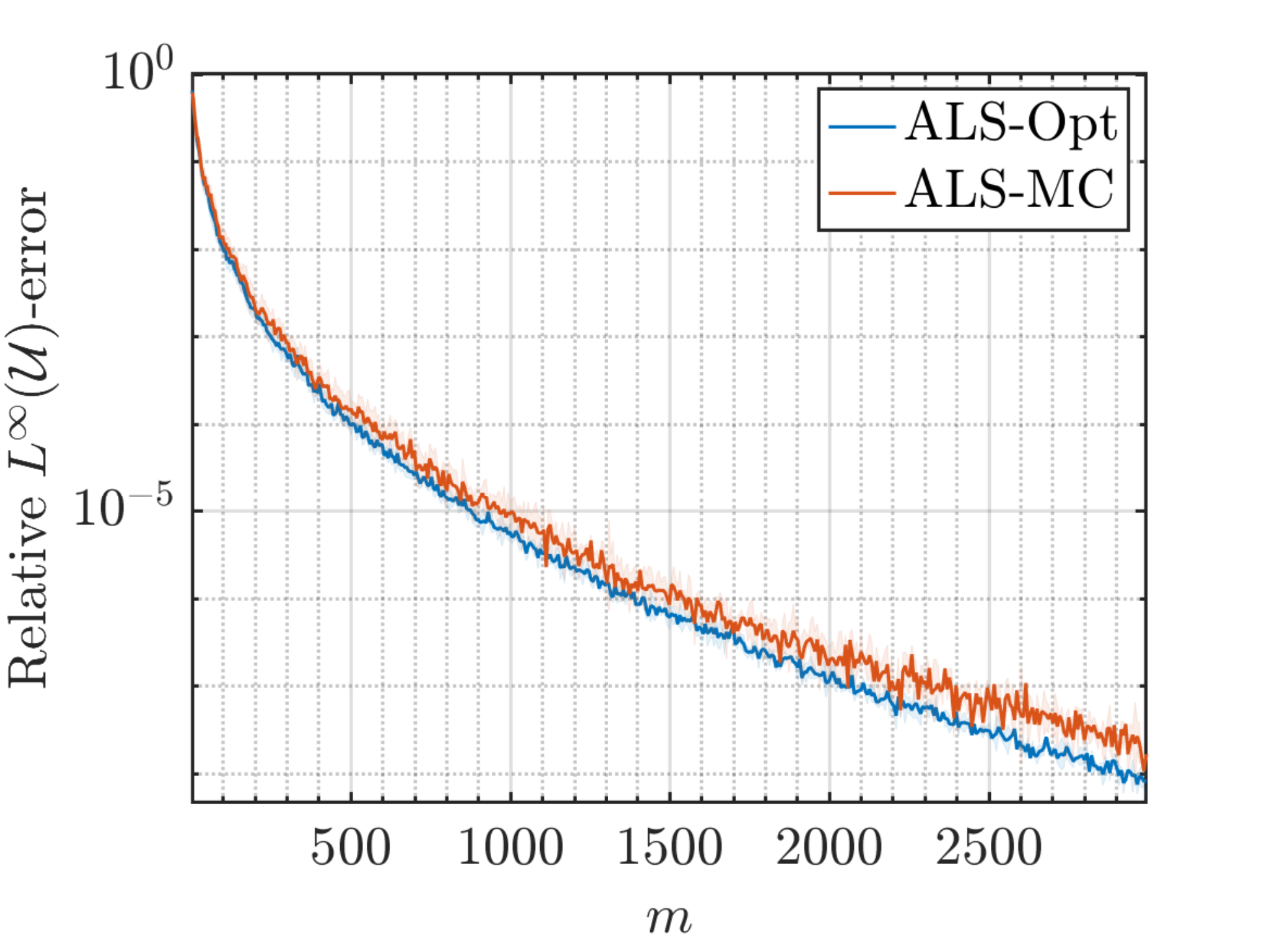}&
\includegraphics[width = \errplotimg]{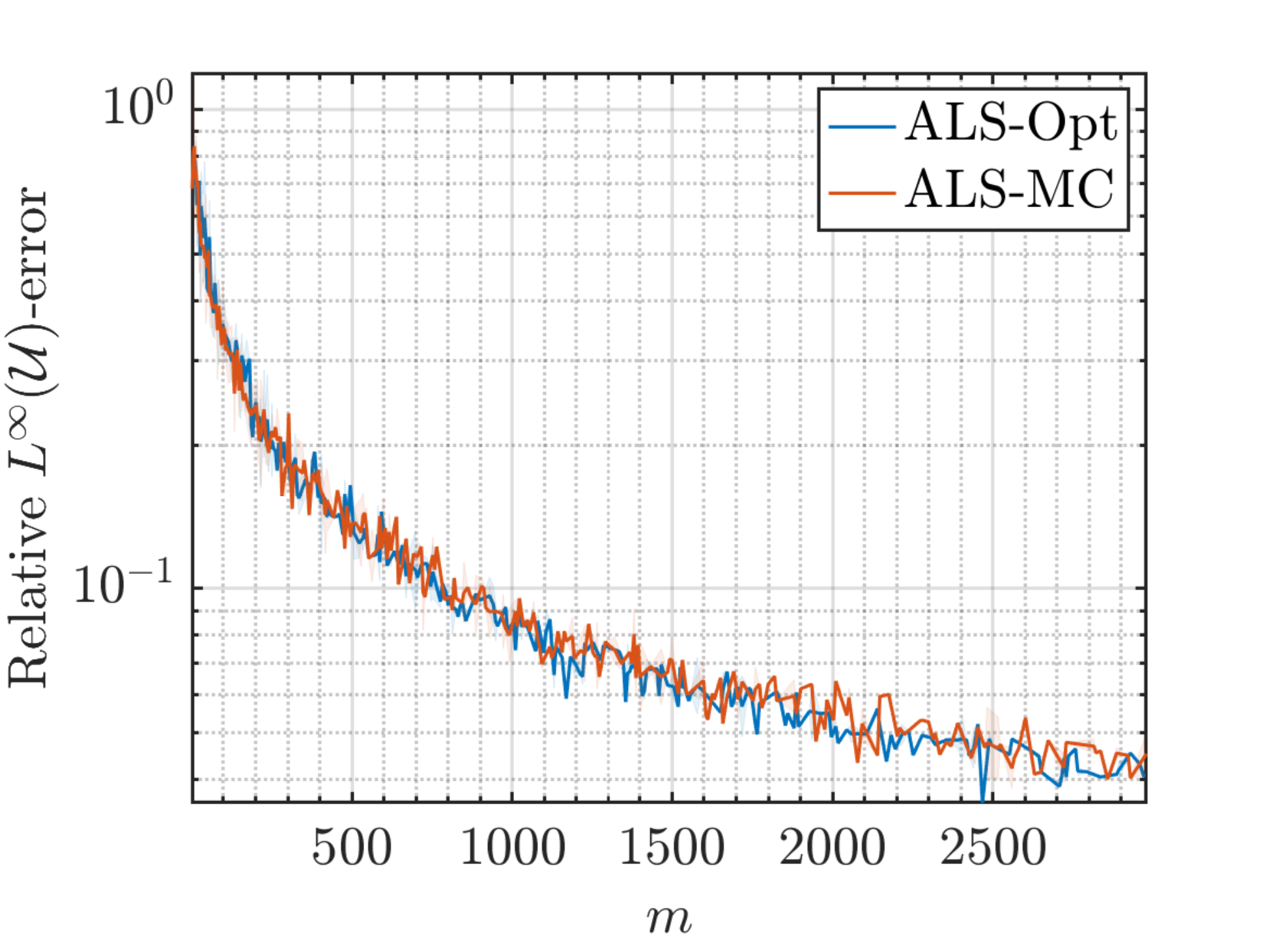}
\\[\errplotgraphsp]
\includegraphics[width = \errplotimg]{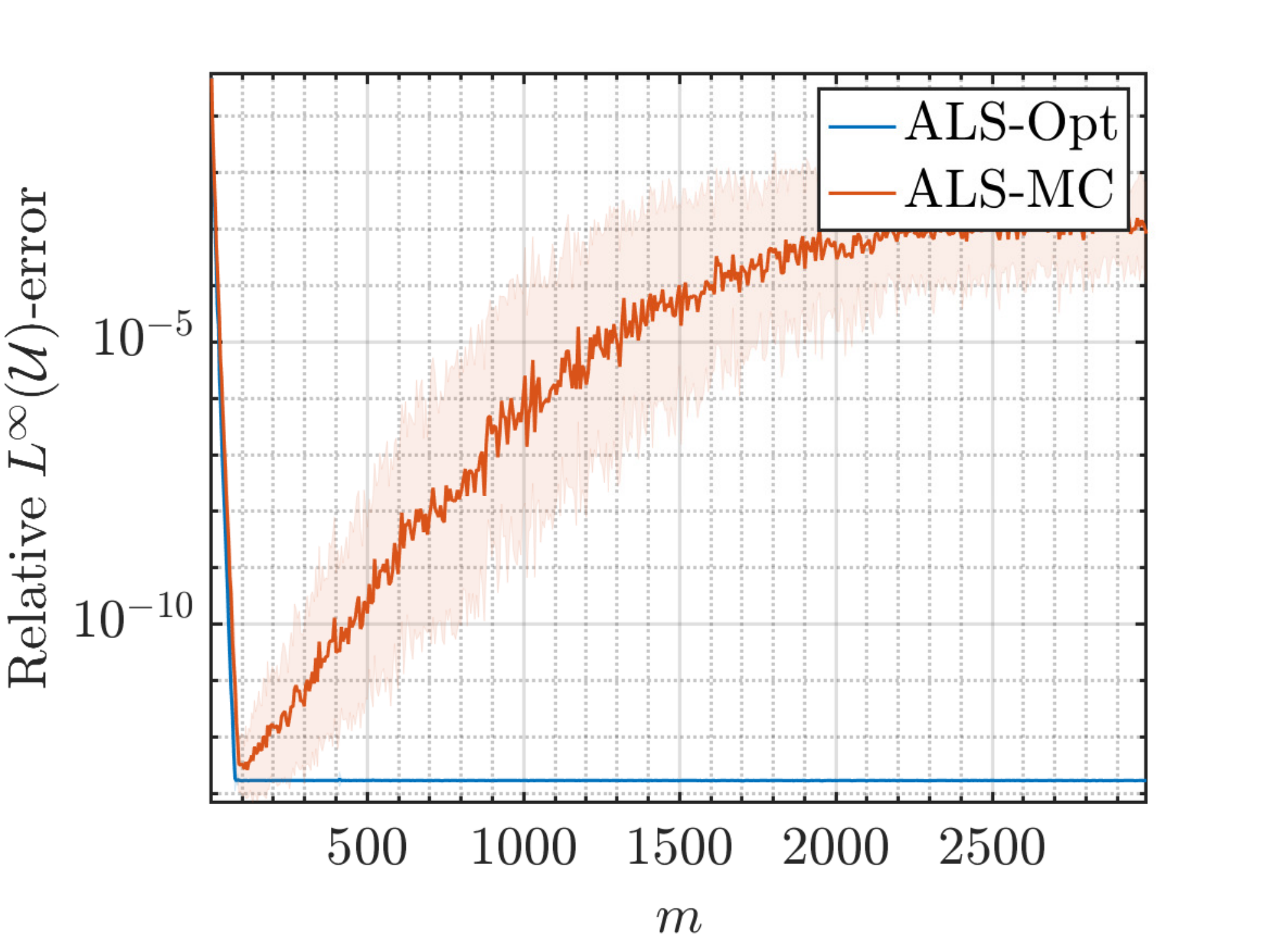}&
\includegraphics[width = \errplotimg]{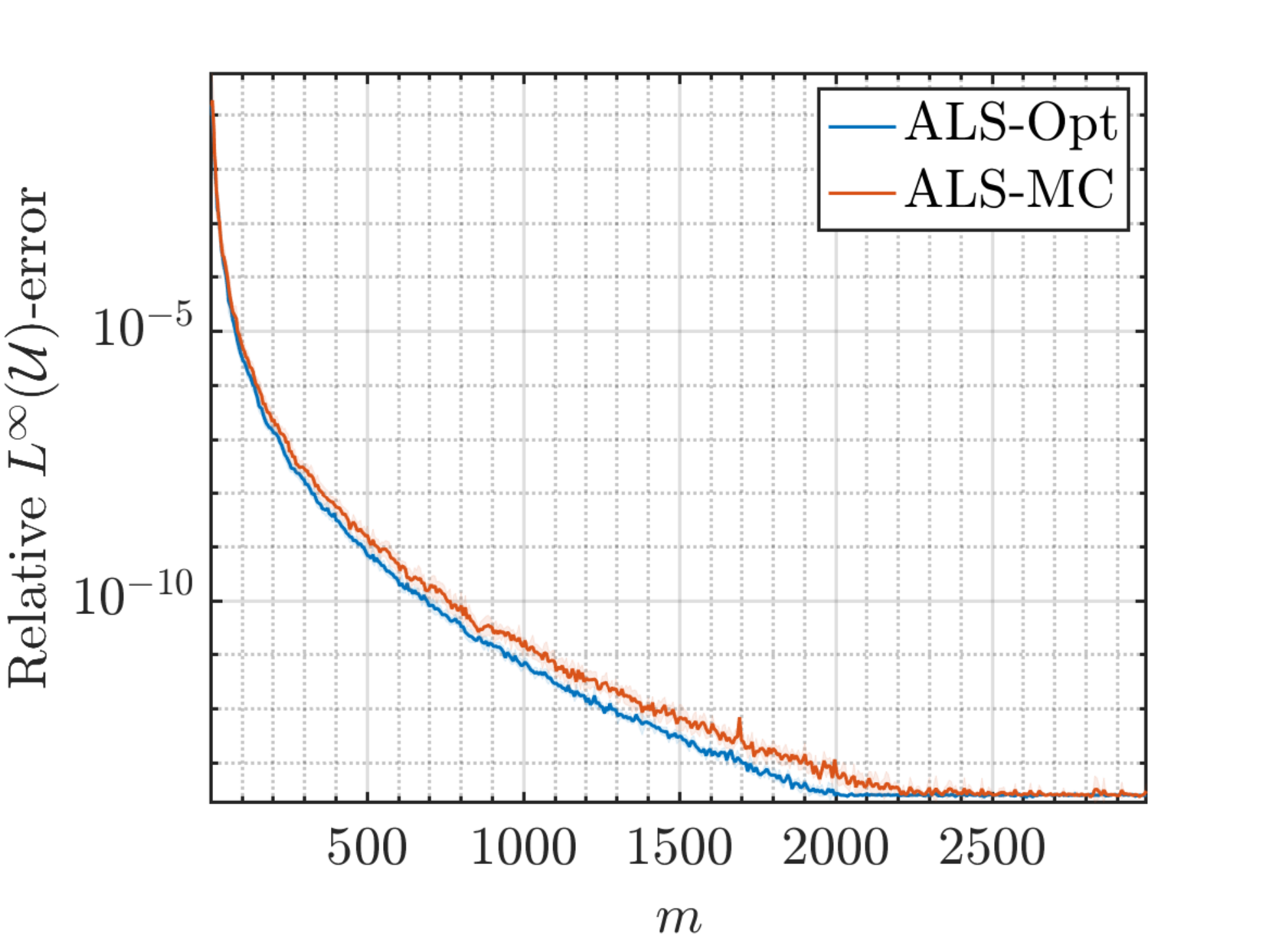}&
\includegraphics[width = \errplotimg]{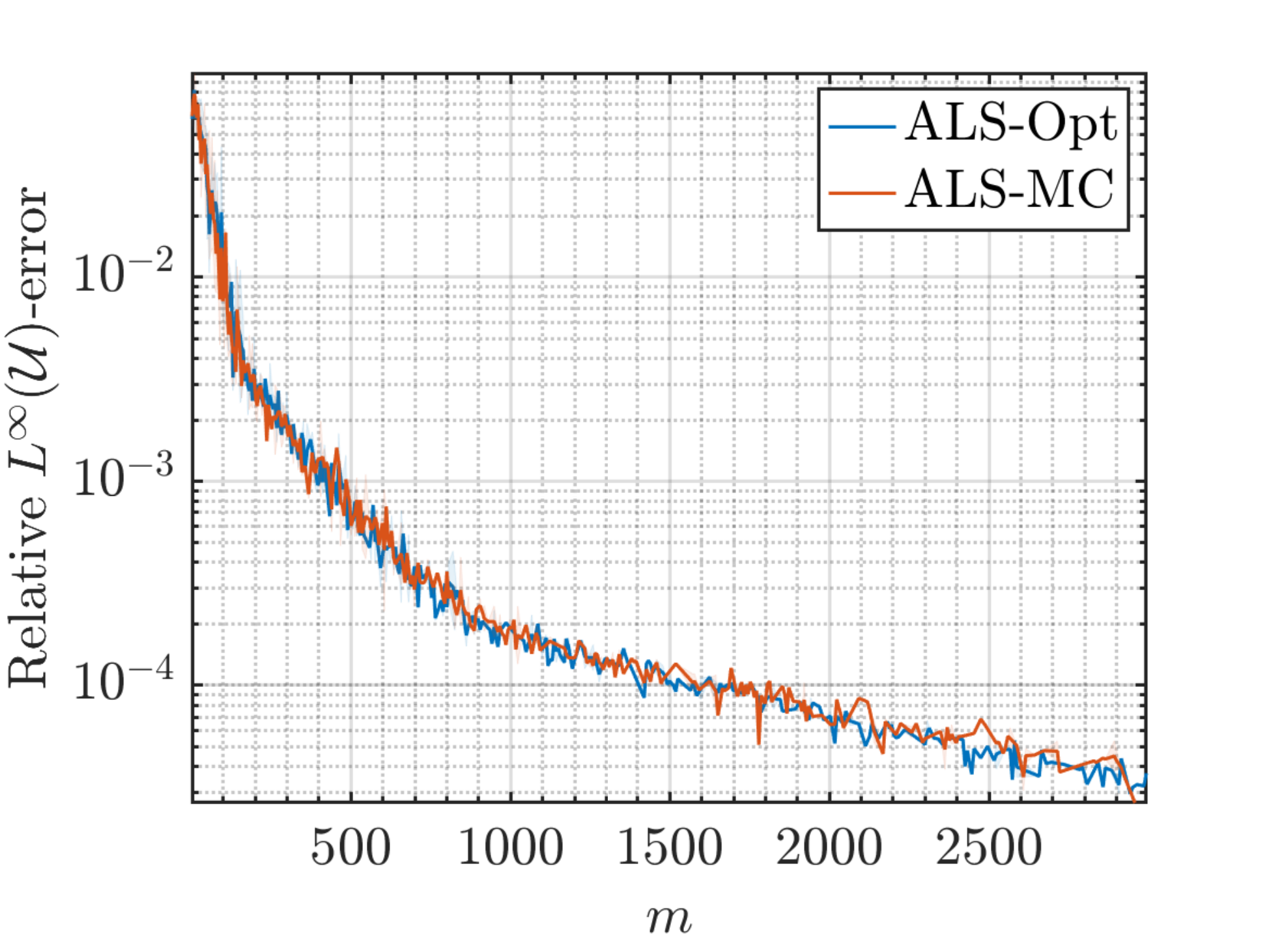}
\\[\errplotgraphsp]
\includegraphics[width = \errplotimg]{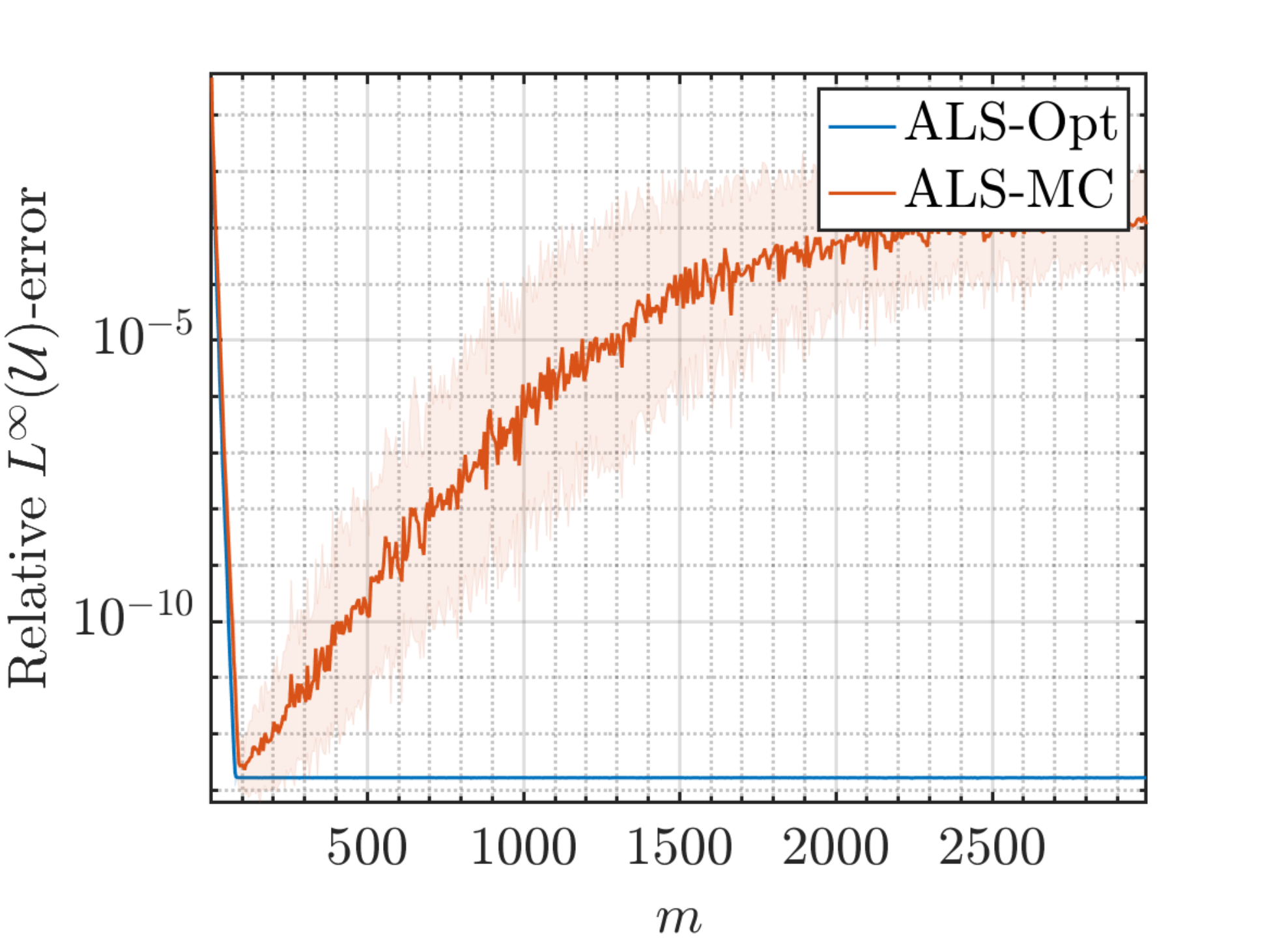}&
\includegraphics[width = \errplotimg]{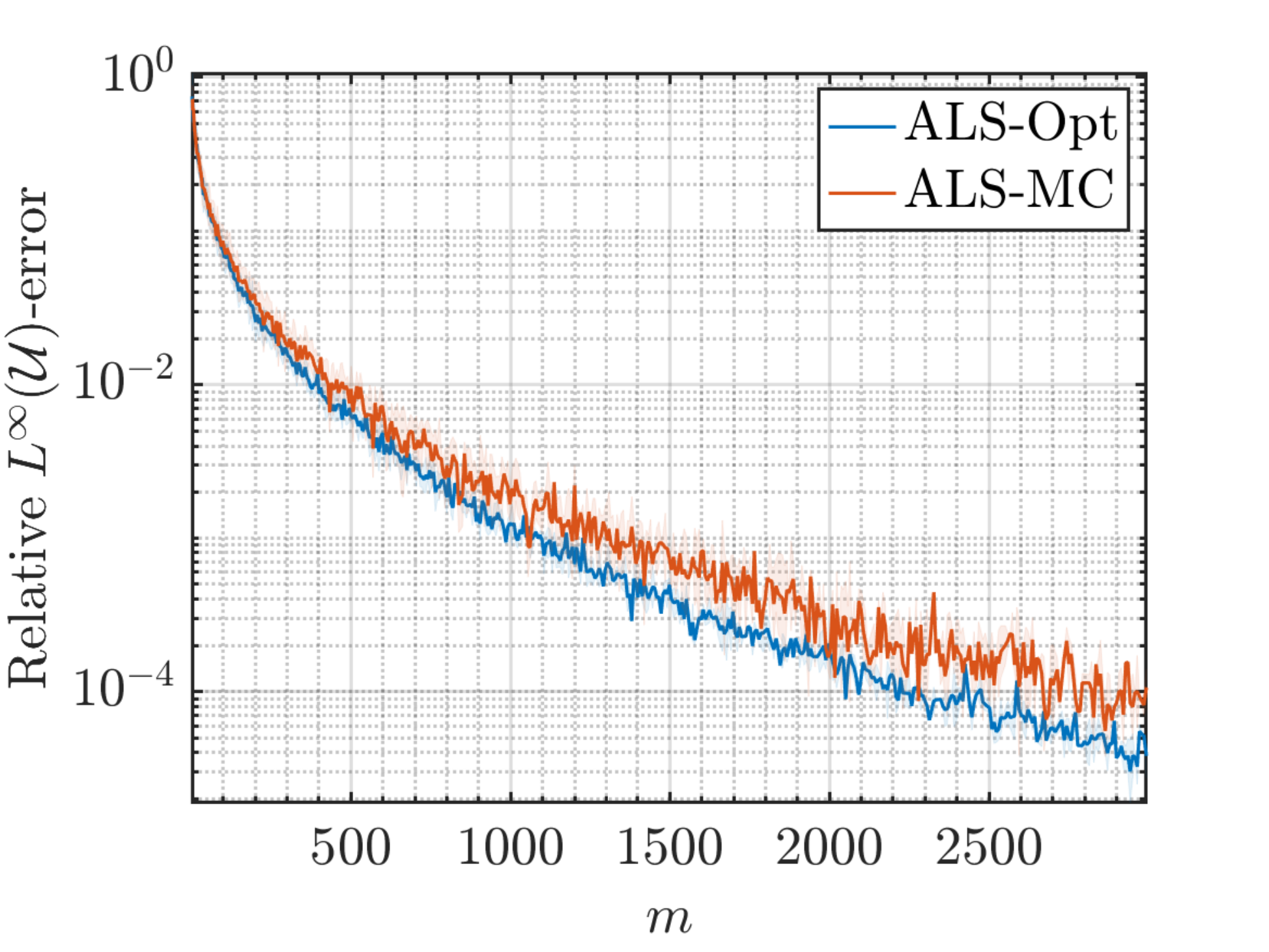}&
\includegraphics[width = \errplotimg]{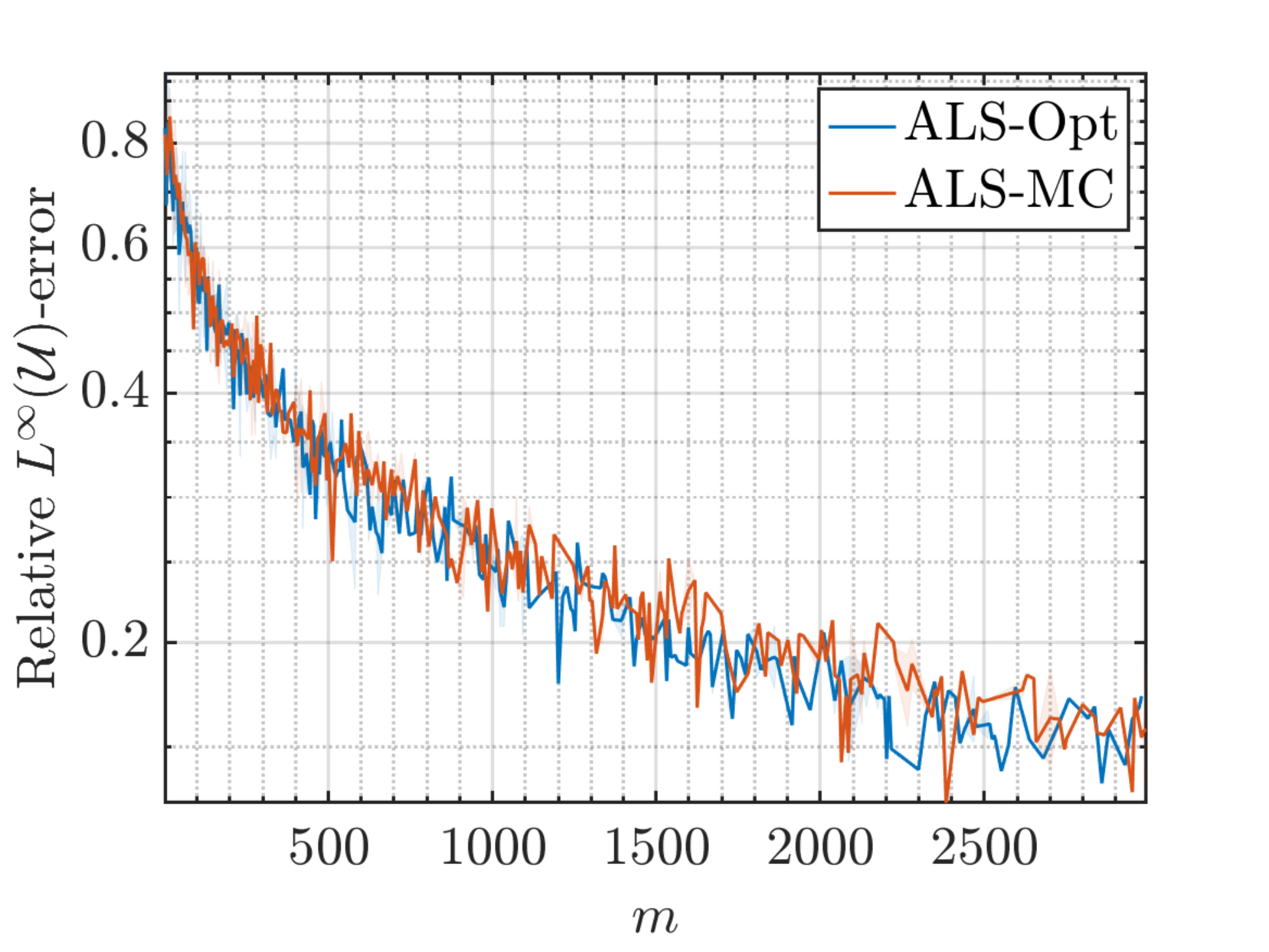}
\\[\errplotgraphsp]
\includegraphics[width = \errplotimg]{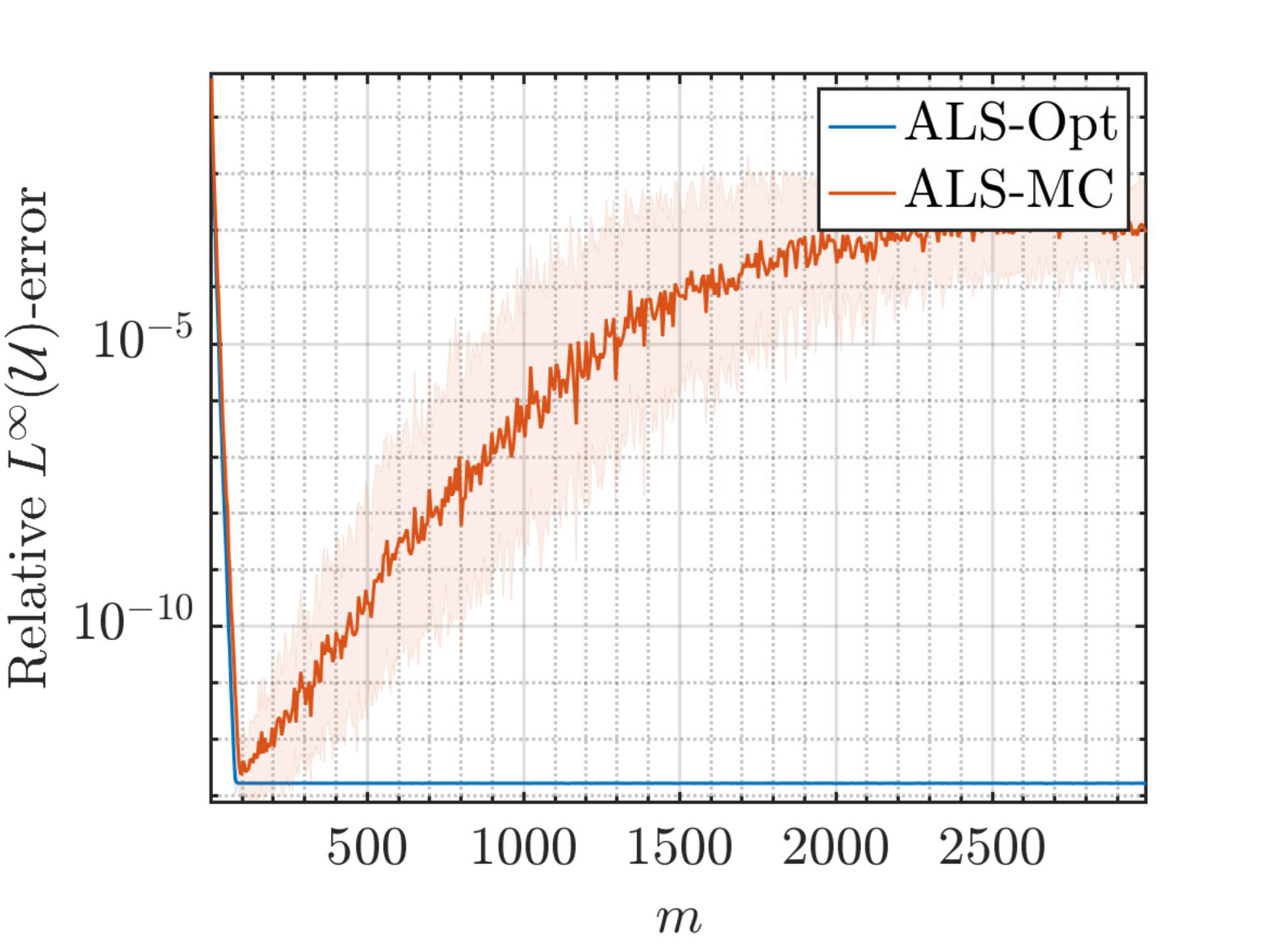}&
\includegraphics[width = \errplotimg]{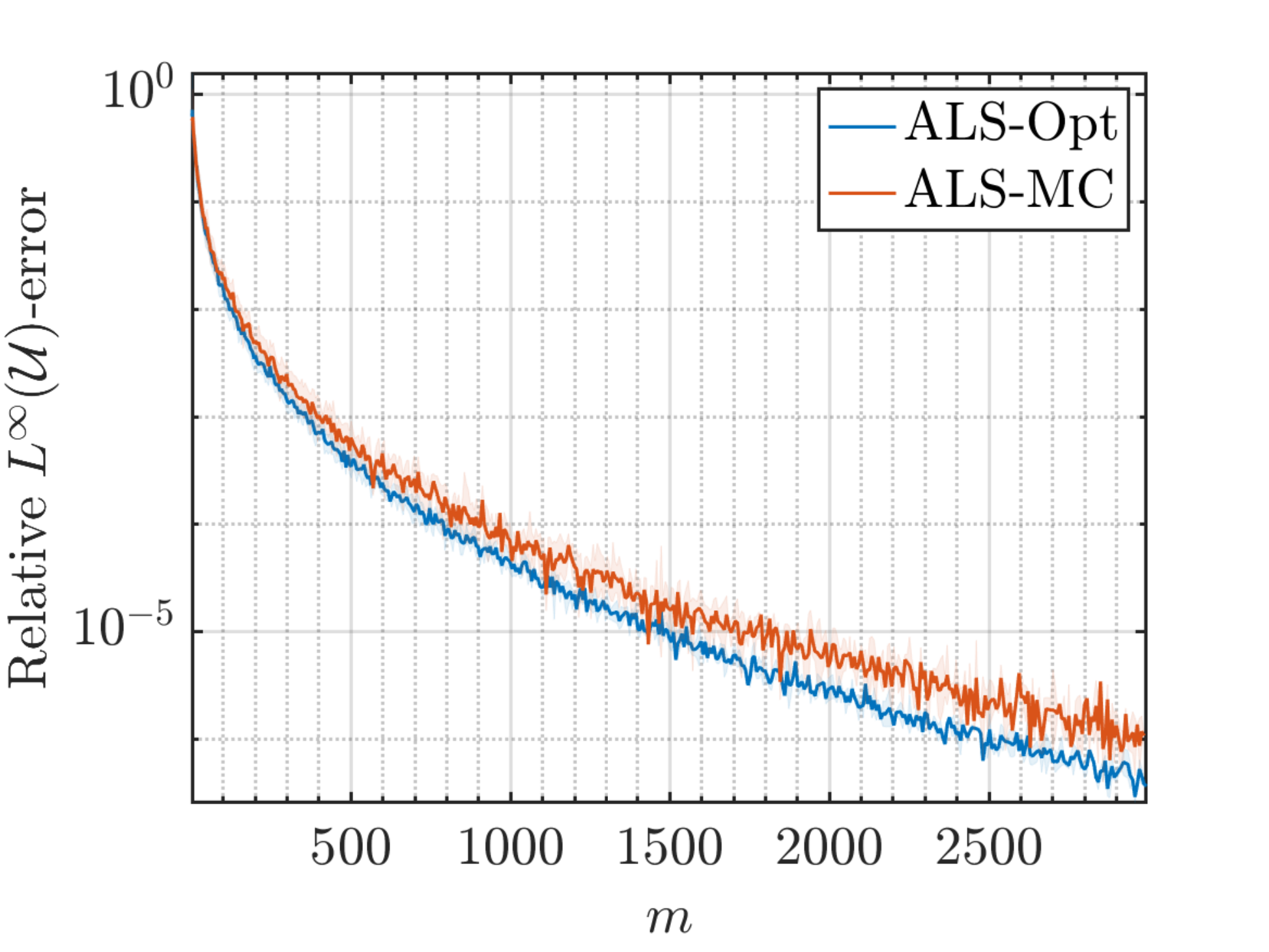}&
\includegraphics[width = \errplotimg]{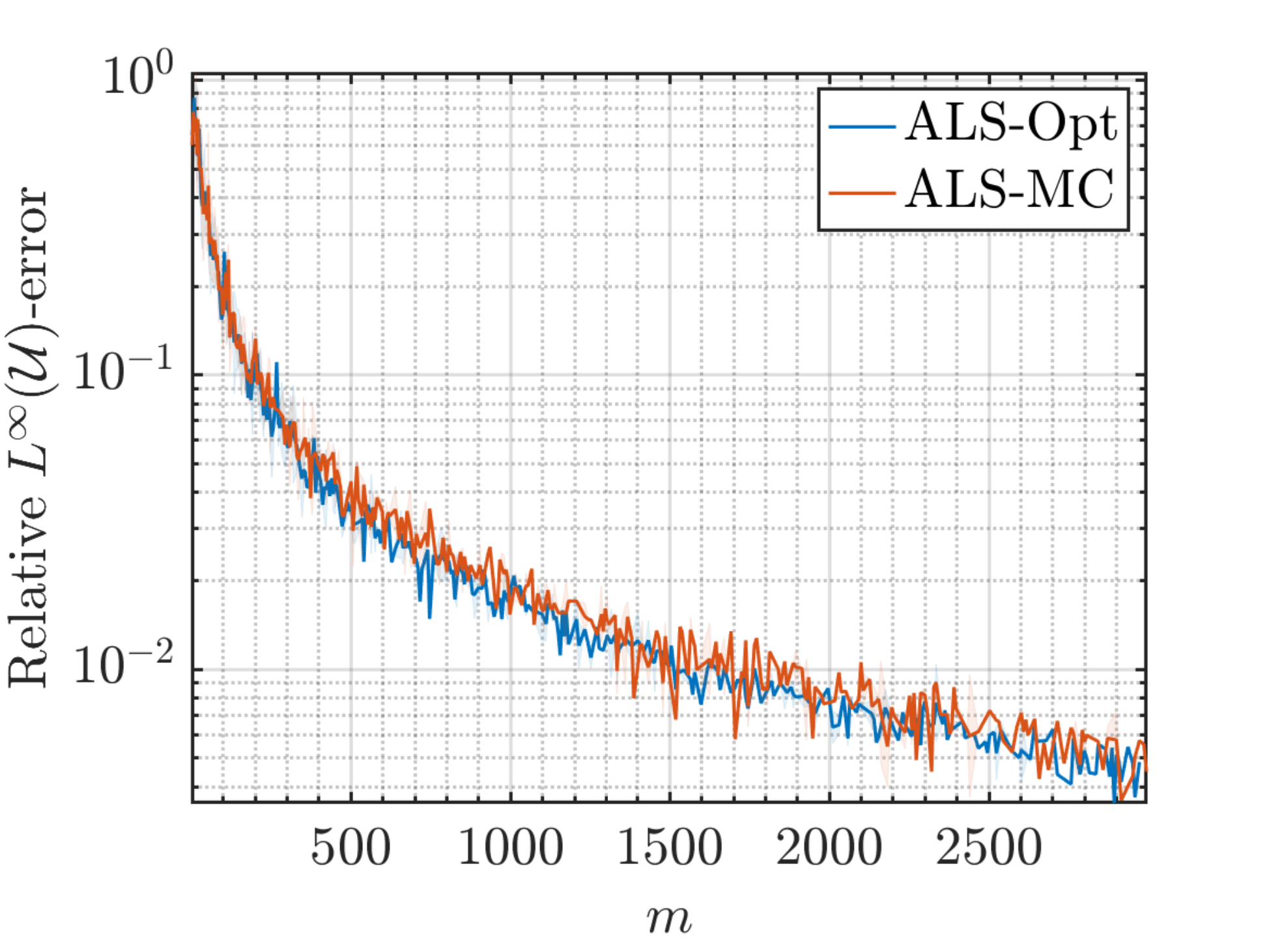}
\\[\errplottextsp]
$d = 1$ & $d = 4$ & $d = 32$
\end{tabular}
\end{small}
\end{center}
\caption{{ALS approximation of the functions $f = f_1$ (first row), $f = f_2$ (second row), $f = f_3$ with $\delta_i = i$ (third row) and $\delta_i = i^2$ (fourth row) and $f = f_4$ (fifth row) using Legendre polynomials and either near-optimal sampling (\S \ref{ss:opt-samp-scheme}) or MC sampling (\S \ref{ss:MC-samp-scheme}). This figure shows the relative $L^{\infty}(\cU)$-norm error versus $m$.}}
\label{fig:figSM14}
\end{figure}

{
Finally, in Fig.\ \ref{fig:figSM14} we repeat several experiments from the main paper but using the relative (discrete) $L^{\infty}$-norm error, rather than the $L^2_{\varrho}$-norm error. Specifically, we show the error versus $m$ for the ALS scheme for the functions considered in Figs.\ \ref{fig:fig2}--\ref{fig:fig5}. Note that this error is computed as in \eqref{f-error-fun}, after replacing the summation by a maximum.   
}

{
The main conclusion from this experiment is that the same conclusion holds in this stronger norm. Namely, in low dimensions MC sampling performs very badly, yet as the dimension increases, it offers very similar performance to the near-optimal sampling scheme. 
}

\bibliographystyle{siamplain}
\bibliography{MCgoodrefs}

%% file: MCgood_shared_Revision.tex
\usepackage{lipsum}
\usepackage{amsfonts}
\usepackage{graphicx}
\usepackage{epstopdf}

\usepackage{algorithmic}
\ifpdf
  \DeclareGraphicsExtensions{.eps,.pdf,.png,.jpg}
\else
  \DeclareGraphicsExtensions{.eps}
\fi

\usepackage{benstyle_SIAM} % Ben's macros
\usepackage{cite}

\newcommand\errplotimg{0.301\textwidth}
\newcommand\errplotsp{-4pt}
\newcommand\errplottextsp{-3.25pt}
\newcommand\errplotgraphsp{-5.5pt}

\newsiamremark{remark}{Remark}
\newsiamremark{hypothesis}{Hypothesis}
\crefname{hypothesis}{Hypothesis}{Hypotheses}
\newsiamthm{claim}{Claim}

\headers{MC is a good sampling strategy}{B. Adcock and S. Brugiapaglia}

\title{Monte Carlo is a good sampling strategy for polynomial approximation in high dimensions\thanks{Submitted to the editors DATE.
\funding{BA acknowledges the support of NSERC through grant RGPIN-2021-611675.SB acknowledges the support of NSERC through grant RGPIN-2020-06766, the Faculty of Arts and Science of Concordia University and the CRM Applied Math Lab.}}}

\author{Ben Adcock\thanks{Department of Mathematics, Simon Fraser University, 8888 University Drive, Burnaby BC, Canada, V5A 1S6
  (\email{ben\_adcock@sfu.ca}, \url{https://www.benadcock.ca}).}
\and Simone Brugiapaglia\thanks{Department of Mathematics and Statistics, Concordia University, J.W. McConnell Building, 1400 De Maisonneuve Blvd. W., Montr\'eal, QC, Canada, H3G 1M8 
  (\email{simone.brugiapaglia@concordia.ca}, \url{https://www.simonebrugiapaglia.ca}).}}

\usepackage{amsopn}